\pdfoutput=1
\documentclass[10pt]{article}
\usepackage[a4paper, top=1in, bottom=1in, left=1in, right=1in]{geometry}
\usepackage[numbers]{natbib}
\usepackage{amssymb}
\usepackage{comment}
\usepackage{amsmath}
\usepackage{amsthm}
\usepackage{mathrsfs}
\usepackage{stmaryrd}
\usepackage{setspace}
\usepackage{graphicx}
\usepackage{xcolor}
\usepackage{stmaryrd}
 \usepackage{algorithm} 
\usepackage{algpseudocode}
\algnewcommand\algorithmicinput{\textbf{INPUT:}}
\algnewcommand\INPUT{\item[\algorithmicinput]}
\algnewcommand\algorithmicoutput{\textbf{OUTPUT:}}
\algnewcommand\OUTPUT{\item[\algorithmicoutput]}
\usepackage{hyperref}[]
\hypersetup{
    colorlinks=true,
    linkcolor=blue,
    filecolor=magenta,      
    urlcolor=cyan,
      citecolor=blue
}
\usepackage{subcaption} 
\usepackage{algorithm}
\usepackage{romannum}

\usepackage{titletoc}

\usepackage{graphicx}
\usepackage{lineno}
\allowdisplaybreaks
\usepackage{authblk}

\newtheorem{assumption}{Assumption}
\newtheorem{lemma}{Lemma}[section]
\newtheorem{proposition}{Proposition}[section]
\newtheorem{corollary}{Corollary}[section]

\newtheorem{remark}{Remark}
\newtheorem{theorem}{Theorem}[section]

\newcommand\AppendixToC{%
  \startcontents
  \setcounter{tocdepth}{2} 
  \printcontents{}{1}{\section*{Appendix Contents}\vskip3pt\hrule\vskip5pt} 
  \vskip3pt
}

\newcommand{\cA}{{\cal A}}

\newcommand{\cG}{{\cal G}}

\newcommand{\cT}{{\cal T}}

\newcommand{\cX}{{\cal X}}
\newcommand{\cY}{{\cal Y}}

\newcommand{\mcA}{\mathcal{A}}
\newcommand{\mcB}{\mathcal{B}}

\newcommand{\mcE}{\mathcal{E}}

\newcommand{\mcG}{\mathcal{G}}

\newcommand{\mcN}{\mathcal{N}}
\newcommand{\mcO}{\mathcal{O}}
\newcommand{\mcP}{\mathcal{P}}

\newcommand{\mcT}{\mathcal{T}}

\newcommand{\mcX}{\mathcal{X}}
\newcommand{\mcY}{\mathcal{Y}}
\newcommand{\mcZ}{\mathcal{Z}}

\newcommand{\ulambda}{\underline{\lambda}}  
\newcommand{\olambda}{\overline{\lambda}}   
\newcommand{\oR}{\overline{R}}              
\newcommand{\oor}{\overline{r}}              
\newcommand{\op}{\overline{p}}              

\newcommand{\RN}[1]{\Romannum{#1}}

\newcommand{\whZ}{\widehat{Z}}
\newcommand{\whE}{\widehat{E}}
\newcommand{\whT}{\widehat{T}}
\newcommand{\whU}{\widehat{U}}
\newcommand{\whDelta}{\widehat{\Delta}}

\newcommand{\wtU}{\widetilde{U}}

\newcommand{\bel}{\begin{eqnarray}\label}
\newcommand{\eel}{\end{eqnarray}}
\newcommand{\bes}{\begin{eqnarray*}}
\newcommand{\ees}{\end{eqnarray*}}

\def\Mat{\hbox{\rm Mat}}
\def\Vec{\hbox{\rm Vec}}

\def\P{\mathbb{P}}

\def\R{\mathbb{R}}

\title{Statistical Inference for Low-Rank Tensor Models}
\author[1]{Ke Xu\thanks{kxu6@nd.edu}}
\author[2]{Elynn Chen\thanks{elynn.chen@stern.nyu.edu}}
\author[1]{Yuefeng Han\thanks{yuefeng.han@nd.edu}}
\affil[1]{Department of Applied and Computational Mathematics and Statistics, University of Notre Dame}
\affil[2]{Department of Technology, Operations, and Statistics,
New York University}

\date{}

\begin{document}

\maketitle
\pagenumbering{arabic}

\begin{abstract}
Statistical inference for tensors has emerged as a critical challenge in analyzing high-dimensional data in modern data science. This paper introduces a unified framework for inferring general and low-Tucker-rank linear functionals of low-Tucker-rank signal tensors for several low-rank tensor models. 
Our methodology tackles two primary goals: achieving asymptotic normality and constructing minimax-optimal confidence intervals. By leveraging a debiasing strategy and projecting onto the tangent space of the low-Tucker-rank manifold, we enable inference for general and structured linear functionals, extending far beyond the scope of traditional entrywise inference. Specifically, in the low-Tucker-rank tensor regression or PCA model, we establish the computational and statistical efficiency of our approach, achieving near-optimal sample size requirements (in regression model) and signal-to-noise ratio (SNR) conditions (in PCA model) for general linear functionals without requiring sparsity in the loading tensor. Our framework also attains both computationally and statistically optimal sample size and SNR thresholds for low-Tucker-rank linear functionals. Numerical experiments validate our theoretical results, showcasing the framework's utility in diverse applications. This work addresses significant methodological gaps in statistical inference, advancing tensor analysis for complex and high-dimensional data environments.
\end{abstract}



\section{Introduction} 



In recent years, the study of tensors and high-dimensional arrays has gained significant attention across fields such as statistics, applied mathematics, machine learning, and data science. Tensors frequently appear in scientific domains such as compressed sensing \citep{caiafa2013multidimensional, friedland2014compressive}, neuroimaging \citep{li2018tucker, bi2021tensors}, recommendation systems \citep{frolov2017tensor, zhang2022tensor}, and econometrics and finance \citep{li2015tensor, wang2024high}. These tensors are often high-dimensional, with their ambient dimensions far exceeding the sample size. Yet, many practical scenarios reveal that tensors possess low-dimensional structures, such as low rank or sparsity \citep{kolda2009tensor, udell2019big}, driving advancements in tensor estimation and structural recovery techniques.

Despite substantial progress in tensor estimation \citep{chen2019non, han2022optimal, kressner2014low, zhang2019cross, xia2019polynomial}, statistical inference for tensors remains relatively underexplored. Most existing studies focus on matrix/tensor completion and recovery from missing data \citep{chen2019inference, ma2024statistical,xia2021normal}, with less attention paid to fully observed tensor data. Furthermore, much of the existing work has emphasized entrywise inference \citep{chen2019inference, agterberg2024statistical}. In contrast, the inference of general linear functionals---expressed as $\langle \mcT, \mcA \rangle$, where $\mcA$ is a loading tensor and $\mcT \in \mathbb{R}^{p_1 \times p_2 \times p_3}$ represents the signal tensor of interest---has been largely overlooked. These linear functionals are crucial in various applications, offering flexibility in capturing both localized and aggregated features. To demonstrate this, we present three example tasks.

\paragraph{\it Task 1: Inference of Specific Entries.} In many real-world applications, researchers often focus on inferring specific entries of the signal tensor $\mcT$. In these cases, the loading tensor $\mcA$ is typically sparse, containing only a few non-zero elements corresponding to the targeted entries. For example, tensor-based methods have been used to analyze spatiotemporal gene expression data, capturing interactions among genes, regions, and time points using a three-mode gene expression tensor \citep{liu2022characterizing}. If the goal is to compare the expression levels of gene $g$ in the region $s$ at two different time points $t_1$ and $t_2$, the difference of interest can be expressed as $\langle \mcT, \mcA \rangle = \mcT_{g,s,t_1} - \mcT_{g,s,t_2}$.

\paragraph{\it Task 2: Inference of a Subgroup of Entries.} In some scenarios, the linear functional $\langle \mcT, \mcA \rangle$ involves a subgroup of entries across one or more modes of the signal tensor $\mcT$. For instance, in recommender systems, the $(i, j, k)$-th entry of a third-order tensor may represent interactions among user $i$, item $j$, and context $k$. To decide whether to recommend item either $j_1$ or $j_2$ to a group of users, indexed by $\mathcal{G} \subseteq [p_1]$, this decision-making can be formalized by testing the condition
$
\langle \mcT, \mcA \rangle = \sum_{i \in \mathcal{G}} \langle e_i \otimes (e_{j_1} - e_{j_2})^\top \otimes e_k, \mcT \rangle > 0,
$
as demonstrated in \citet{zhang2022tensor, xia2021statistical}. Similarly, in network traffic analysis, traffic data indexed by source $s$, destination $g$, and time $t$ can be represented as a tensor \citep{zhou2016robust}. The total traffic at time $t$ is captured by the linear functional $\langle \mcT, \mcA \rangle = \sum_{s, g} \mcT_{s, g, t}$, which aggregates entries across specific modes.

\paragraph{\it Task 3: Inference of General Linear Functionals.} 
In more general cases, the loading tensor $\mcA$ exhibits a more complex structure. For example, in neuroimaging studies, the relationship between brain imaging data and continuous clinical or cognitive assessment scores is modeled using a parameter tensor $\mcT$ \citep{li2016sparse}. In this context, given an input image $\mcA$, the fitted clinical score is $\langle \mcT, \mcA \rangle$. Unlike simpler cases, the input image $\mcA$ is not necessarily of low Tucker rank. Inference on this general linear functional is essential to quantify the uncertainty of the predicted clinical outcomes.

\vspace{4pt}
\indent These examples underscore the need for statistical inference methods capable of handling linear functionals of signal tensors, with loading tensors $\mcA$ exhibiting diverse sparsity and structural complexity. While entrywise inference (Task 1) has been studied in tensor completion, Tasks 2 and 3 involve linear functionals with potentially many nonzero entries and more complex structures, which are largely underexplored. To accommodate this diversity, we propose to model the complexity of the loading tensor through its Tucker rank. For entrywise inference (Task 1), the Tucker rank of $\mcA$ is near minimal, typically $\operatorname{rank}(\mcA) \approx (1,1,1)$. For subgroup aggregate inference (Task 2), the Tucker rank may range from $(1,1,1)$ to a small tuple, such as $(R_1, R_2, R_3)$. In more general applications (Task 3), $\mcA$ may reach full rank, i.e., $(p_1, p_2, p_3)$. Using Tucker rank to model $\mcA$'s complexity provides a unified framework for handling different structural requirements.

Tensor regression and tensor principal component analysis (PCA) are two fundamental problems that motivate the exploration of signal tensors. Both approaches leverage low-rank tensor structures to address the challenges of high-dimensional data. 
For a comprehensive discussion on low-rank tensor models, 
we refer readers to \citet{bi2021tensors, liu2022tensor}. This naturally raises a critical question:
\begin{itemize}
\item \textbf{What are the sample size or signal-to-noise ratio (SNR) requirements for inferring general linear functionals of the parameters in low-rank tensor models? }
\end{itemize}

The field of tensor inference builds on advancements in statistical inference for low-rank matrices, as matrices can be seen as mode-2 tensors.
\citet{cai2016geometric} investigated the inference of general linear functionals under the low-rank matrix trace regression framework with Gaussian design. However, their results required a sample size proportional to $p_1 p_2$ to construct a valid confidence interval for a low-rank parameter matrix $M \in \mathbb{R}^{p_1 \times p_2}$ with rank $r$. 
This prompts the question: can valid inferences for general linear functionals be achieved with fewer samples than the order of $p_1 p_2 p_3$ for a low-rank signal tensor $\cT\in \mathbb{R}^{p_1 \times p_2 \times p_3}$ in tensor regression?

In this paper, we focus on statistical inference for both full-Tucker-rank (referred to as {\bf general linear functionals} throughout this article) and {\bf low-Tucker-rank linear functionals} of low-Tucker-rank signal tensors.
Our investigation centers on two fundamental settings: tensor regression and tensor PCA. In both cases, the least squares estimator in tensor regression and the estimator in tensor PCA can be viewed as the true signal tensor corrupted by either complex or simple noise terms.
We develop a unified inference framework that debiases the initial estimate and projects it onto the singular spaces of the signal tensor. This framework constructs confidence intervals by identifying the pivot quantity on the tangent space of the low-Tucker-rank manifold. It accommodates a wide variety of loading tensor structures, extending beyond entrywise inference to facilitate applications involving general linear functionals of the signal tensor.
With certain modifications, our framework could extend beyond tensor regression and PCA to accommodate a range of low-rank tensor models, including tensor completion \citep{zhang2019cross,xia2021statistically}, generalized tensor estimation \citep{han2022optimal}, and high-order tensor clustering \citep{han2022exact,luo2022tensor}, enabling inference for diverse low-rank structures.

\subsection{Main results}


Table \ref{tab: main results} provides an overview of our key findings on the inference of low-rank and general linear functionals in both tensor regression and tensor PCA settings. To simplify notation and emphasize the main contributions, in this section, we assume that the condition number $\kappa$ of the signal tensor $\mcT = \mcG \times_1 U_1 \times_2 U_2 \times_3 U_3 \in \mathbb{R}^{p_1 \times p_2 \times p_3}$ is fixed. 
For tensor regression, we further assume that $\ulambda \geq c$, where $c$ is a positive constant and $\ulambda$ is the smallest singular value across all mode-$j$ matricizations of $\mathcal{T}$. 

\begin{table}[h!] 
    \centering
    \setlength{\tabcolsep}{2pt} 
    \renewcommand{\arraystretch}{2} 
        \centering
         \begin{tabular}{c|c|c|c}
            \hline
            \textbf{Loading Tensor $\mcA$} & \textbf{Sample Size} ($n$) & 
            \textbf{Incoherence Condition} & 
            \textbf{Sample Splitting}\\
             \hline
            low-rank $(R_1, R_2, R_3)$ 
            & $\max\big\{\kappa^2 \op / \ulambda^2, \op\oor^2\big\}$ &  $\frac{\left\|\mcP_{U_j}A_j\left(\mcP_{U_{j+2}}\otimes \mcP_{U_{j+1}}\right)\right\|_{\mathrm{F}}}{\left\|\mcP_{U_{j\perp}}A_j\left(\mcP_{U_{j+2}}\otimes \mcP_{U_{j+1}}\right)]\right\|_{\mathrm{F}}}\lesssim \ulambda\oor^{\frac{1}{2}}\op^{-\frac{1}{2}}$& 
            Required \\
            \hline
            low-rank $(R_1, R_2, R_3)$ 
            & $\max \big\{\kappa^2 \op / \ulambda^2, \op^\frac{3}{2}\oor \big\}$ & $\frac{\left\|\mcP_{U_j}A_j\mcP_{ (U_{j+2} \otimes U_{j+1} )G_j^{\top}} \right\|_{\mathrm{F}}}{\left\|\mcP_{U_{j\perp}}A_j\mcP_{(U_{j+2} \otimes U_{j+1})G_j^{\top}} \right\|_{\mathrm{F}}} \lesssim \ulambda\oor^{\frac{1}{2}}\op^{-\frac{1}{4}}$& 
            Required \\
            \hline
            full rank $(p_1, p_2, p_3)$ & $\max\{\kappa^2 \op / \ulambda^2, \op^2\oor\}$ & Not required & Not required\\
            \hline
        \end{tabular} \\
    \vspace{0.3cm}
      (a) Tensor Regression Setting \\
    \vspace{0.3cm} 
        \centering
        \begin{tabular}{c|c|c}
            \hline
            \textbf{Loading Tensor $\mcA$} & \textbf{SNR Requirement} ($\ulambda$) & \textbf{Incoherence Condition} \\
            \hline
            low-rank $(R_1, R_2, R_3)$ & $\max\big\{ \kappa\op^\frac{1}{2}, \op^{\frac{3}{4}}\oor^{\frac{1}{2}}\big\}$ & $\frac{\left\|\mcP_{U_j}A_j\mcP_{\left(U_{j+2} \otimes U_{j+1}\right)G_j^{\top}}\right\|_{\mathrm{F}}}{\left\|\mcP_{U_{j\perp}}A_j\mcP_{\left(U_{j+2} \otimes U_{j+1}\right)G_j^{\top}}\right\|_{\mathrm{F}}} \lesssim \oor^{\frac{1}{2}}\op^{-\frac{1}{4}}$\\
            \hline
            full rank $(p_1, p_2, p_3)$ & $\max\big\{ \kappa\op^\frac{1}{2}, \op\oor^{\frac{1}{2}}\big\}$ & Not required \\
            \hline
        \end{tabular} \\
    \vspace{0.3cm}
      (b) Tensor PCA Setting
      \vspace{0.3cm}
    \caption{Summary of main results on the inference of low-Tucker-rank and general linear functionals under tensor regression and PCA settings. Here, $\kappa$ denotes the condition number of $\mcT=\mcG \times_1 U_1 \times_2 U_2 \times_3 U_3$, $U_j \in \mathbb{O}^{p_j \times r_j}$ are left singular spaces of $\Mat_j(\mcT)$, $\mcG \in \mathbb{R}^{r_1 \times r_2 \times r_3}$ is the core tensor of the signal tensor $\mcT$, and $\op = \max\{p_1, p_2, p_3\},\oor = \max\{r_1, r_2, r_3\}$. The projection matrix onto the left singular space of matrix $A$ is denoted by $\mcP_A$. In addition to the sample size/ SNR requirement (up to a logarithmic factor of $\op$) and incoherence conditions, alignment conditions, which establish a lower bound on variance components, are necessary for valid inference in all scenarios. A detailed discussion on the interplay between the incoherence and alignment conditions is presented in Section \ref{sec: The Interplay between the Incoherence Condition and Alignment Condition}. } \label{tab: main results}
\end{table}

\paragraph*{Inference for General Linear Functionals} We establish central limit theorems that facilitate the statistical inference for general linear functionals in both tensor regression and PCA frameworks. To the best of our knowledge, this is the first work to address the inference of general linear functionals without imposing incoherence conditions on the loading tensor $\mathcal{A}$ in tensor models, filling a significant gap in the literature. Previous studies have primarily concentrated on tensor completion, recovery from missing data, and entrywise inference, leaving more complex inference problems underexplored. We introduce a unified framework for inferring general linear functionals of signal tensors, thereby expanding the applicability of tensor analysis to a broader range of practical scenarios. While our focus is on third-order tensors for clarity, the framework is readily extensible to higher-order tensors.

Our results achieve nearly computationally optimal sample sizes in tensor regression and optimal signal-to-noise ratio (SNR) requirements in tensor PCA, all without additional incoherence conditions. Specifically, for a signal tensor $\mcT \in \mathbb{R}^{p_1 \times p_2 \times p_3}$ with $\op = \max_{j=1,2,3}{p_j}$, we demonstrate that a sample size of $n\gtrsim \op^2\oor$ suffices to infer any linear functional in tensor regression. This improves upon the matrix trace regression results by \citet{cai2016geometric}, which require $n\gtrsim p_1 p_2$ for a low-rank parameter matrix $M \in \mathbb{R}^{p_1 \times p_2}$. For tensor PCA, we show that an SNR ratio of $\ulambda \gtrsim \op\oor^{1/2}$ is sufficient for accurate inference of any linear functional of the signal tensor, where the SNR $\ulambda$ is defined as the smallest singular value across all mode-$j$ matricizations of $\mathcal{T}$, denoted by
$
\ulambda = \min_{j=1,2,3} \sigma_{\text{min}}(\operatorname{Mat}_j(\mcT)).
$

\paragraph*{Inference for Low-Tucker-Rank Linear Functionals} In this context, we achieve computationally optimal results w.r.t. $\op$: a sample size requirement of $n \gtrsim \op^{3/2}\oor$ in tensor regression with sample splitting and an SNR requirement of $\ulambda\gtrsim \op^{3/4}\oor^{1/2}$ in tensor PCA. These results are derived under relaxed incoherence conditions of $\mathcal{O}(\oor^{1/2}\op^{-1/4})$, significantly less restrictive than the conventional incoherence conditions of $\mathcal{O}(\oor^{1/2}\op^{-1/2})$ commonly assumed in the matrix and tensor inference literature  \citep{chen2019inference, xia2021statistical, agterberg2024statistical, ma2024statistical}. Furthermore, under stricter incoherence conditions of $\mathcal{O}(\oor^{1/2}\op^{-1/2})$ and assuming homogeneous singular values, we establish that a statistically optimal sample size of $n \gtrsim \op\oor^2 $ w.r.t. $\op$ in tensor regression is sufficient to ensure valid inference.

\paragraph*{Methodological Innovations} Our framework differs from leave-one-out approaches that establish entrywise confidence intervals using eigenvector distributions and $\ell_{2,\infty}$ perturbation analysis \citep{chen2019inference, xia2021statistical, ma2024statistical, agterberg2024statistical}. Instead, we employ a refined $\ell_2$ perturbation analysis of the spectral representation of perturbation terms. Building upon the spectral representation of singular space estimation proposed by \citet{xia2021normal}, our approach leverages the first-order perturbation term for inference while demonstrating that higher-order terms are negligible. The assumption of a low Tucker rank for the loading tensor $\mathcal{A}$ simplifies the complexity of perturbations in the $\ell_2$ norm, thereby achieving computational optimality.

Additionally, we develop novel concentration inequalities for quadratic and higher-order sub-Gaussian polynomials, essential for bounding negligible terms in our analysis. The random components in the estimated linear forms are expressed as polynomials of sub-Gaussian entries from noise or design tensors. 
The sub-Gaussian design differentiates our approach from tensor completion settings, which typically assume uniform sampling. Moreover, we relax the dependence on the condition number, ensuring that our analysis remains computationally optimal even as the condition number $\kappa$ of the signal tensor $\mathcal{T}$ grows at a rate of $\mathcal{O}(\op^{1/4})$ and becomes arbitrarily large.



\subsubsection{Tensor Regression}
Statistical inference in the tensor regression framework, despite its importance, has been relatively unexplored. Using double-sample splitting, we demonstrate that asymptotic normality can be achieved with a computationally optimal sample size of $n \asymp \op^{3/2}\oor$, w.r.t. $\op$,
\begin{align*}
\sqrt{n} \big( \langle \widehat{\cT} , \mathcal{A} \rangle- \langle \mcT, \mathcal{A}\rangle \big)/ (\sigma_\xi / \sigma \cdot s_{\mcA}) \stackrel{d}{\rightarrow} \mathcal{N}(0,1) ,
\end{align*}
provided that the following {\it incoherence condition} holds,
\begin{align}
{\big\|\mcP_{U_j}A_j\mcP_{(U_{j+2}\otimes U_{j+1})G_j^{\top}}\big\|_{\mathrm{F}}} / {\big\|\mcP_{U_{j\perp}}A_j\mcP_{(U_{j+2}\otimes U_{j+1})G_j^{\top}}\big\|_{\mathrm{F}}} \lesssim \ulambda \oor^{1/2}\op^{-1/4}, \label{eq: computationally optimal incoherence condition in tensor regression with sample splitting}
\end{align}
for $j=1,2,3$, where $\sigma,\sigma_\xi$ are parameters of the design tensor and the noise (See Assumption~\ref{assump: subgaussian design tensor in tensor regression} and \ref{assump: subgaussian noise in tensor regression} below), $A_j=\Mat_j(\mcA)$ is the mode-$j$ matricization of the loading tensor $A$, $\otimes$ denotes Kronecker product, $\times_j$ represent mode-$j$ product, and the variance component of the estimated linear functional is
\begin{align}
s_{\mcA}^2=\sum_{j=1}^3\big\|\mcP_{U_{j \perp}} A_j \mcP_{(U_{j+2} \otimes U_{j+1}) G_j^{\top}}\big\|_{\mathrm{F}}^2+\|\mathcal{A} \times_1 U_1 \times_2 U_2 \times_3 U_3 \|_{\mathrm{F}}^2 .  \label{eq: variance component sA2}
\end{align}
Here, $\mcP_{U_j}$ denotes the projection matrix onto the left singular space of $\operatorname{Mat}_j(\mcT)$, while $\mcP_{(U_{j+2}\otimes U_{j+1})G_j^{\top}}$ represents the projection matrix onto its right singular space. 

Asymptotic normality requires an initial estimate $\widehat{\mcT}^{\text{init}}$ satisfying $\|\widehat{\mcT}^{\text{init}} - \mcT\|_{\mathrm{F}} = o[(r\log(\op))^{-1}]$, alongside an alignment condition to ensure a lower bound on the variance component $s_\mcA^2$. This prevents the cases where $\mcA$ is nearly orthogonal to the singular space of the signal tensor $\mcT$. The alignment conditions are elaborated upon in the following sections. Under stronger incoherence and alignment conditions, specifically
\begin{align}
& \max_{j=1,2,3}{\big\|\mcP_{U_j}A_j(\mcP_{U_{j+2}}\otimes \mcP_{U_{j+1}})\big\|_{\mathrm{F}}} / {\big\|\mcP_{U_{j\perp}}A_j\mcP_{(U_{j+2}\otimes U_{j+1})G_j^{\top}}\big\|_{\mathrm{F}}} \lesssim \ulambda \oor^{1/2}\op^{-1/2}, \label{eq: statistically optimal incoherence condition in tensor regression with sample splitting}
\end{align}
we achieve a statistically optimal sample size of $n \asymp \op\oor^2$, w.r.t. $\op$.

 
In the absence of sample splitting, a larger sample size of $n \asymp \op^2\oor$ is required to overcome the dependence between the initial estimate and the bias-correction step due to repeated use of the dataset. In this case, under the minimax optimal initial estimation error $\|\widehat{\mcT}^{\text{init}} - \mcT\|_{\mathrm{F}} \asymp (\op\oor/n)^{1/2}$ and the alignment condition, the asymptotic normality of any linear functional $\langle \widehat{\mcT}, \mcA\rangle$ can be established without any incoherence condition.

\subsubsection{Tensor PCA}

In the tensor PCA framework, a computational optimal signal-to-noise ratio (SNR) of $\ulambda\gtrsim \op^{3/4}\oor^{1/2}$, w.r.t. $\op$, suffices for valid inference of low-Tucker-rank linear functionals $\langle \mcT, \mcA \rangle$. This is achieved by establishing the asymptotic normality
\begin{align*}
\big( \langle \widehat{\cT} , \mathcal{A} \rangle- \langle \mcT, \mathcal{A} \rangle \big)/ (\sigma \cdot s_{\mcA} ) \stackrel{d}{\rightarrow} \mathcal{N}(0,1)  ,  
\end{align*}
where the variance component $s_{\mcA}$ is defined in \eqref{eq: variance component sA2}, $\sigma$ is the parameter of the noise tensor (See Assumption~\ref{assump: sigma-sub-Gaussian noise in tensor PCA} below), provided that the following {\it incoherence condition} 
\begin{align}
{\big\|\mcP_{U_j}A_j\mcP_{(U_{j+2}\otimes U_{j+1})G_j^{\top}} \big\|_{\mathrm{F}}} / {\big\|\mcP_{U_{j\perp}}A_j\mcP_{(U_{j+2}\otimes U_{j+1})G_j^{\top}}\big\|_{\mathrm{F}}} \lesssim \oor^{1/2}\op^{-1/4}, \label{eq: computationally optimal incoherence condition in tensor PCA}
\end{align} 
along with the corresponding alignment condition. For the special case of entrywise inference, where the loading tensor $\mcA = e_{j_1} \otimes e_{j_2} \otimes e_{j_3}$, the incoherence condition simplifies to $\max_j \|U_j\|_{2,\infty} \lesssim \oor^{1/2}\op^{-1/4}$. This is significantly weaker than the stricter conditions commonly assumed in matrix/tensor completion literature \citep{chen2019inference,xia2021statistical,ma2024statistical} and in recent tensor PCA results \citep{agterberg2024statistical}, such as $\max_j\|U_j\|_{2,\infty} \asymp \oor^{1/2}\op^{-1/2}$. For rank-one signal tensors, our results are consistent with recent findings on entrywise inference for rank-one tensors \citep{xia2022inference}. For more general linear functionals $\langle \mcT, \mcA \rangle$, an SNR of $\ulambda \gtrsim \overline{p}\oor^{1/2}$ is required to ensure asymptotic normality. 

\subsection{Related work}

Statistical inference for tensor regression remains a largely unexplored area of research. However, insights can be drawn from methodologies developed for matrix and tensor completion under uniform sampling. Most work in low-rank matrix completion focuses on entrywise inference, typically requiring incoherence conditions. For example, \citet{xia2021statistical} employed spectral perturbation techniques, while \citet{chen2019inference} used a leave-one-out approach for entrywise statistical inference in matrix completion. Both methods rely on incoherence conditions for the singular spaces $U_j$, imposing an upper bound on the $\ell_{2,\infty}$ norm of $U_j$ to ensure uniform noise magnitudes across all rows. A comprehensive discussion of $\ell_{2,\infty}$ geometry can be found in \citet{cape2019two}. However, verifying incoherence conditions in practical applications is often non-trivial. The standard assumption $\|U_j  \|_{2,\infty} \asymp \sqrt{\oor/ \op}$ represents the most stringent scenario, achieving the lower bound of the $\ell_{2,\infty}$ norm.

In the realm of tensor completion, \citet{ma2024statistical} recently tackled the inference problem for general linear functionals of the signal tensor, extending the matrix completion framework of \citet{xia2021statistical} to tensors. However, their approach depends on the incoherence conditions of \citet{agterberg2024statistical} and requires that $\|\mcA \|_{\ell_1} / \|\mcA \|_{\mathrm{F}}$ remains bounded, thereby enforcing sparsity on the loading tensor $\mcA$. This sparsity aligns with Task 1 scenarios,
where the linear functional $\langle \mcA, \mcT \rangle$ involves only a few tensor entries.

Significant progress in tensor inference has been made by \citet{xia2022inference}, who developed methods for singular space inference for low-Tucker-rank tensors and entrywise inference for rank-one tensors under Gaussian noise in tensor PCA.
Building on this work, \citet{agterberg2024statistical} extended entrywise inference to low-Tucker-rank signal tensors with heteroskedastic sub-Gaussian noise in the tensor PCA setting, achieving computationally optimal SNR under the incoherence condition $\|U_j \|_{2,\infty} \asymp \sqrt{\oor/ \op}$ for each mode $j$.

While entrywise inference has been extensively explored in tensor completion and PCA models (\citealp{xia2021statistical, ma2024statistical, agterberg2024statistical}) as well as in matrix completion models (\citealp{chen2019inference, xia2021statistical}), inference involving linear functionals with many non-zero entries, as described in Tasks 2 and 3, remains underdeveloped.
In such cases, the $\ell_1$ norm of $\mcA$ becomes substantial, rendering existing methods unsuitable.

A crucial component of our framework is obtaining an accurate initial estimate of the signal tensor in tensor regression. Prior research has proposed various methods for low-rank tensor estimation in regression settings \citep{chen2019non, raskutti2019convex, zhou2013tensor, tomioka2013convex}. More recently, \citet{zhang2020islet} introduced an importance sketching framework for tensor estimation under Gaussian designs, achieving minimax optimal error rates. \citet{han2022optimal} further extended non-convex optimization approaches to more general settings, demonstrating minimax optimal rates for low-Tucker-rank tensors under both Gaussian and sub-Gaussian noise. These advancements in tensor estimation provide the foundation for our inference methodology.

For both tensor regression and tensor PCA, our framework relies on accurately estimating the singular spaces of each tensor mode. In tensor PCA, \citet{zhang2018tensor} established the minimax optimal rate for estimating the singular spaces of a parameter tensor $\mcT \in \mathbb{R}^{p_1 \times p_2 \times p_3}$ in tensor PCA, under the assumption of i.i.d. sub-Gaussian noise. These rates serve as critical benchmarks for our approach.

\subsection{Organization}

The remainder of this paper is organized as follows. In Section \ref{sec: notation}, we define the notation and terminology used throughout. Section \ref{sec: inference for tensor regression} introduces a comprehensive inference framework for tensor regression, presenting algorithms for constructing debiased estimators both without sample splitting and with sample splitting, along with their theoretical guarantees for asymptotic normality. This section also explores the relationship between incoherence and alignment conditions and constructs confidence intervals with proven theoretical properties, demonstrating their minimax optimality.
We then extend our framework to tensor PCA in Section \ref{sec: tensor_pca}, detailing the corresponding algorithms, theoretical guarantees, and the construction of confidence intervals, while also establishing their minimax optimality in this context. A proof sketch of the main theorems is presented in Section \ref{sec: proof sketch of main theorems}.
Numerical experiments are provided in Section \ref{sec: Numerical Experiments} to validate our methods, followed by a discussion of future work in Section \ref{sec: Discussion}. Technical proofs are provided in the Supplementary Material.

\section{Notation} \label{sec: notation}

Throughout this paper, we use the following notation. Tensors are denoted by calligraphic letters, such as $\mcT$, $\mcG$, and $\mcA$. Matrices are represented by uppercase letters like $A$ and $B$, while vectors are indicated by lowercase letters such as $u$ and $v$. For two sequences of real numbers $\{a_n\}$ and $\{b_n\}$, write $a_n\lesssim b_n$ (resp. $a_n\gtrsim b_n$) if there exist a constant $C$ independent of $n$ such that $a_n\le Cb_n$ (resp. $a_n\ge Cb_n$), and write $a_n \asymp b_n$ if there are positive constants $c$
and $C$ such that $c \le a_n/b_n \le C$ for all $n$.

For a matrix $A \in \mathbb{R}^{p_1 \times p_2}$, we denote its projection matrix onto the left singular space by $\mcP_A$ and its projection matrix onto the right singular space by $\mcP_{A^{\top}}$. The singular value decomposition (SVD) of $A$ is expressed as $A = U_1 \Lambda U_2^\top$, where $\Lambda \in \mathbb{R}^{r \times r}$ is a diagonal matrix containing the singular values of $A$, and $U_1 \in \mathbb{R}^{p_1 \times r}, U_2 \in \mathbb{R}^{p_2 \times r}$ are unitary matrices. From this, $\mcP_A = U_1 U_1^\top$ and $\mcP_{A^\top} = U_2 U_2^\top$. Let $\Lambda=\text{diag}(\sigma_1(A), \sigma_2(A), ..., \sigma_{r}(A))$, where the singular values are arranged in descending order $\sigma_{\max}(A) := \sigma_1(A)\ge\sigma_2(A)\ge \cdots\ge \sigma_{\min}(A) := \sigma_{r}(A)\ge 0$. For any orthonormal matrix $U \in \mathbb{R}^{p \times r}$, we denote its orthogonal complement by $U_{\perp} \in \mathbb{R}^{p \times (p - r)}$. 

This paper focuses on the analysis of 3-mode tensors. For such tensors, the Tucker rank of a tensor $\mcT$, denoted as $\operatorname{rank}(\mathcal{T})$, is defined as a tuple $(r_1, r_2, r_3)$ representing the ranks along each mode. Let ${\rm{Vec}}(\cdot)$ be the vectorization of matrices and tensors. The Frobenius norm of a tensor or matrix is given by $\|\cT\|_{\mathrm{F}}=\| \Vec(\cT) \|_2$, while the spectral norm of a matrix is denoted by $\|A\|=\sigma_1(A)$. The inner product of two tensor, $\cA,\cT$, is defined as $\langle \cT, \cA \rangle = \Vec(\cT)^\top \Vec(\cA)$. The vectorized $\ell_p$-norm of a tensor or matrix is written as $\|\cT\|_{\ell_p}=\|\Vec(\cT)\|_p$, where, in particular, $\|\cT\|_{\ell_\infty}$ represents the maximum absolute value among all entries.

For a tensor $\mcT \in \mathbb{R}^{p_1 \times p_2 \times p_3}$ with Tucker rank $(r_1, r_2, r_3)$, we define $\overline{p} = \max\{p_1, p_2, p_3\}$ and $\overline{r} = \max\{r_1, r_2, r_3\}$ to simplify expressions involving the largest dimension or rank. The operator $\times_j$ denotes the mode-$j$ product between a tensor and a matrix. Specifically, for a tensor $\mathcal{G} \in \mathbb{R}^{r_1 \times r_2 \times r_3}$ and a matrix $V_1 \in \mathbb{R}^{p_1 \times r_1}$, the mode-1 product is defined as
$ [\cG \times_1 V_1]_{i_1,i_2,i_3} = \sum_{j_1=1}^{r_1} [\mathcal{G} ]_{j_1, i_2, i_3} [V_1]_{i_1, j_1}$ for $i_1\le p_1, i_2\le r_2, i_3\le r_3$.  

The set $\mathcal{M}_{(r_1, r_2, r_3)}$ of tensors with a fixed Tucker rank $\mathbf{r}=(r_1, r_2, r_3)$ forms a smooth embedded submanifold of $\mathbb{R}^{n_1 \times n_2 \times n_3}$, with dimension of $r_1r_2r_3+\sum_{j=1}^3 r_j (p_j-r_j)$. For matrices in $\mathbb{R}^{p \times r}$, the Stiefel manifold, denoted by $\mathbb{O}^{p\times r}$ is the set of matrices with orthonormal columns, i.e. $\mathbb{O}^{p\times r} = \{M \in \mathbb{R}^{p \times r} \mid M^{\top}M = \mathcal{I}_{r} \}$.

We use $\operatorname{Mat}_j(\cdot)$ to denote the mode-$j$ matricization (unfolding) of a tensor, which rearranges the tensor into a matrix by stacking its mode-$j$ fibers as columns. For a tensor $\mcT \in \mathbb{R}^{p_1 \times p_2 \times p_3}$, the mode-$j$ matricization $\operatorname{\text{Mat}}_j(\mcT) \in \mathbb{R}^{p_j \times p_{j+2}p_{j+1}}$ is defined as $[\operatorname{\text{Mat}}_j(\mcT)]_{k_j, k_{j+1} + p_{j+1} (k_{j+2} - 1)} = \mcT_{k_1, k_2, k_3}$, where $j+1$ and $j+2$ are computed modulo 3. For simplicity, we denote the mode-$j$ matricization of a tensor $\mcT$ by $T_j := \operatorname{Mat}_j(\mcT)$. This notation is used consistently throughout the paper.

For a tensor $\mcT$ with Tucker ranks $(r_1, \cdots, r_m)$, we define its signal strength of $\mcT$ as
\begin{align}\label{def:ulambda}
\ulambda := \lambda_{\min}(\mcT) = \min_j \{ \sigma_{r_j} (\operatorname{\text{Mat}}_j(\mcT) ) \}, 
\end{align}
which represents the smallest positive singular value among all tensor matricizations. Similarly, the maximum signal strength is defined as 
\begin{align}\label{def:mlambda}
\overline{\lambda} := \lambda_{\max}(\mcT) = \max_j \{\sigma_1(\operatorname{\text{Mat}}_j(\mcT) ) \}.
\end{align}
The condition number of $\mcT$, reflecting the ratio of the maximum to minimum signal strength, is given by $\kappa:= \kappa(\mcT) = \overline{\lambda}/\ulambda$.

\section{Inference for Tensor Regression}\label{sec: inference for tensor regression}

\subsection{Problem Setting}

Consider a collection of i.i.d. random samples $(\mcX_i, Y_i)_{i=1}^n$ modeled as
\begin{align}\label{eq: tensor regression model}
Y_i = \langle \mcT, \mcX_i \rangle + \xi_i,  
\end{align}
where $\mcT$ is a low-rank tensor with Tucker rank $(r_1, r_2, r_3)$, capturing the relationship between the scalar responses $\{Y_i\}_{i=1}^n$ and the tensor covariates $\{\mcX_i\}_{i=1}^n$. The terms $\{\xi_i\}_{i=1}^n$ represent independent noise.
Our objective is to perform statistical inference on the linear functional $\langle \mcT, \mcA \rangle$, where $\mcA \in \mathbb{R}^{p_1 \times p_2 \times p_3}$ is a prespecified loading tensor. In many practical scenarios, the loading tensor $\mcA$ may also have a low Tucker rank, denoted by $(R_1, R_2, R_3)$. For instance, entrywise inference corresponds to a special case where $\mcA$ has a Tucker rank of $(1,1,1)$. The low-rank structures of the signal tensor $\mcT$ and the loading tensor $\mcA$ will be leveraged to construct efficient estimators.

Our model, referred to as scalar-on-tensor regression, relates a multiway predictor to a scalar response. This framework has broad applications across various fields, such as predicting clinical outcomes or attributes from medical images \citep{fang2019image, zhou2013tensor, spencer2022bayesian}. While parameter estimation in tensor regression has been extensively studied \citep{zhou2013tensor,chen2019non,raskutti2019convex}, statistical inference in this context remains largely underexplored.

We formalize the assumptions on the signal tensor $\mcT$ and the loading tensor $\mcA$ as follows.

\begin{assumption}[Structures]\label{assump: structural assumption on the signal tensor and loading tensor}

(i). The signal tensor $\mcT$ follows a Tucker low-rank structure and can be expressed as
$$
\mcT = \mcG \times_1 U_1 \times_2 U_2 \times_3 U_3 \in \mathbb{R}^{p_1 \times p_2 \times p_3}, 
$$
where $\mcG \in \mathbb{R}^{r_1 \times r_2 \times r_3}$ is the core tensor, and $U_j \in \mathbb{O}^{p_j \times r_j}$ is the factor matrix for $j=1,2,3$-th mode of $\mcT$. 

(ii). The loading tensor $\mcA$ can be represented as
$$
\mcA = \mathcal{B} \times_1 V_1 \times_2 V_2 \times_3 V_3 \in \mathbb{R}^{p_1 \times p_2 \times p_3},
$$
where $\mathcal{B} \in \mathbb{R}^{R_1 \times R_2 \times R_3}$ is the core tensor, and $V_j \in \mathbb{O}^{p_j \times R_j}$ is the factor matrix for $j=1,2,3$-th mode of $\mcA$. 
\end{assumption}

\begin{remark}
In this paper, we assume that the signal tensor $\mcT$ has a fixed or slowly growing Tucker rank $(r_1, r_2, r_3 )$. In contrast, the Tucker rank of the loading tensor $\mcA$, denoted by $(R_1, R_2, R_3 )$, may vary across different scenarios. Notably, we consider the case where $\mcA$ has a full Tucker rank, i.e., $(R_1, R_2, R_3 ) = (p_1, p_2, p_3 )$.
\end{remark}

Throughout our analysis, we assume that both the design tensors and the noise terms are sub-Gaussian. Specifically, we impose the following assumptions.
\begin{assumption}[Subgaussian design tensor]\label{assump: subgaussian design tensor in tensor regression}

In the tensor regression model \eqref{eq: tensor regression model}, the design tensors $\{\mcX_i \}_{i=1}^n$ are i.i.d. copies of a random tensor $\mcX$, whose entries $[\mcX]_{j_1,j_2,j_3}$ are i.i.d., mean-zero, and $\sigma$-subgaussian. Specifically, for each entry, we have $\|[\mcX]_{j_1,j_2,j_3}\|_{\psi_2} \leq \sigma$, where $\sigma$ is a positive constant, $\|\cdot \|_{\psi_2}$ denotes the Orcliz $\psi_2$-norm, and there exist positive constants $c$ and $C$ such that $c\sigma^2 \leq \operatorname{Var}(\mcX_{j_1,j_2,j_3}) \leq C\sigma^2$.
\end{assumption}

\begin{assumption}[Subgaussian noise]\label{assump: subgaussian noise in tensor regression}

In the tensor regression model \eqref{eq: tensor regression model}, the noise terms $\{\xi_i \}_{i=1}^n$ are i.i.d., mean-zero, and $\sigma_\xi$-subgaussian. Specifically, we have $\|\xi_i\|_{\psi_2} \leq \sigma_\xi$, where $\sigma_\xi$ is a positive constant.
\end{assumption}

We begin by addressing the inference of general linear functionals.

\subsection{Debiased Estimator of Linear Functionals without Sample Splitting} \label{sec:debias_nonsplit} 

To efficiently estimate $\langle \mcT, \mcA \rangle$, we propose a multi-step algorithm leveraging an initial estimator.

\subsubsection*{Step 1: Initialization} Obtain initial estimates of the signal tensor $\widehat{\mcT}^{\text{init}} \in \mathbb{R}^{p_1 \times p_2 \times p_3}$ and the factor matrices $\whU_1^{\text{init}}$, $\whU_2^{\text{init}}$, and $\whU_3^{\text{init}}$ from the observed data $\{Y_i, \mcX_i\}_{i=1}^n$. Set $\whU_j^{(0)} := \whU_j^{\text{init}}$ for $j = 1, 2, 3$.

To ensure the effectiveness of our inference procedures, we impose the following assumptions on the initial estimates.

\begin{assumption}[Error bound for the initial signal tensor estimate]\label{assump: initial signal tensor estimate error bound in tensor regression} 

The initial estimate $\widehat{\mcT}^{\text{init}}$ is assumed to have the same Tucker rank $(r_1, r_2, r_3)$ as the true signal tensor $\mcT$. Additionally, the estimation error satisfies $\| \widehat{\mcT}^{\text{init}} - \mcT \|_{\mathrm{F}} \leq \Delta$, where $\Delta = o(1)$ with probability at least $1 - \mathbb{P}(\mcE_{\Delta})$, and the event $\mcE_{\Delta}$ is defined as
$$
\mcE_{\Delta} = \{\| \widehat{\mcT}^{\text{init}} - \mcT \|_{\mathrm{F}} > \Delta\}.
$$

\end{assumption}

\begin{remark}
The error bound $\|\widehat{\mcT} - \mcT\|_{\mathrm{F}} \asymp (\sigma_\xi/\sigma) \sqrt{ \op \oor/n }$ is minimax optimal for low-Tucker-rank tensor estimation, as shown in \citet{han2022optimal, zhang2020islet}. Here, $\op = \max\{p_1, p_2, p_3\}$, $\oor = \max\{r_1, r_2, r_3\}$. In Theorem~\ref{thm: main theorem in tensor regression without sample splitting}, we analyze how this initial error affects convergence rates for asymptotic normality. This optimal bound for the initial estimate ensures the nearly computational optimal sample size $n \gtrsim \op^2\oor$ for general linear functional inference.
Using sample splitting, as shown in Theorem~\ref{thm: main theorem in tensor regression with sample splitting}, can remove the dependence between the initial estimator and the bias correction, and allow $\Delta = o[(\oor\log(\op))^{-1}]$ to suffice for valid statistical inference.
\end{remark}

\begin{assumption}[Error bound for the initial singular space estimate]\label{assump: initial singular space estimate error bound in tensor regression}

The initial estimates of the singular spaces are assumed to achieve minimax optimal error rates. Specifically, for each mode $j = 1, 2, 3$, the estimation of projection matrices satisfy $\|\mcP_{\whU_j^{(0)}} - \mcP_{U_j}\| \leq c_0 (\sigma_\xi/\sigma) \sqrt{\op /n }$, with probability at least $1 - \mathbb{P}(\mcE_{U}^{\text{reg}})$, where the event $\mcE_{U}^{\text{reg}}$ is defined as
$$
\mcE_{U}^{\text{reg}}= \cup_{j=1,2,3}\{\|\mcP_{\whU_j^{(0)}} - \mcP_{U_j} \| > c_0 (\sigma_\xi/\sigma) \sqrt{\op /n }  \}.
$$
\end{assumption}

\begin{remark}
The minimax rate $\|\mcP_{\whU_j^{(0)}} - \mcP_{U_j}\| \asymp (\sigma_\xi/\sigma) \sqrt{\op /n }$ ensures that initial singular space estimates are statistically optimal. If this rate is not achieved, then applying Higher-order Orthogonal Iteration (\texttt{HOOI}) \citep{de2000best} to the initial estimate of singular space can refine the estimates as long as $\|\mcP_{\whU_j^{(0)}} - \mcP_{U_j} \| \leq 1/2$. Without loss of generality, we assume the initial estimates already satisfy this optimal rate.


\end{remark}

\subsubsection*{Step 2: Debiasing} 

The initial estimate $\widehat{\mcT}^{\text{init}}$ is typically biased. To address this, we compute a debiased estimator using residuals from the initial model fit
$$
\widehat{\mcT}^{\text{unbs}} = \widehat{\mcT}^{\text{init}} + \frac{1}{n \sigma^2} \sum\nolimits_{i=1}^{n} \big( Y_i - \langle \widehat{\mcT}^{\text{init}}, \mcX_i \rangle \big) \mcX_i.
$$
Let $\widehat{\Delta} = \mcT - \widehat{\mcT}^{\text{init}}$ be the estimation error. The debiased estimator can be decomposed as 
\begin{align}
\widehat{\mcT}^{\text{unbs}} &=\mcT+\frac{1}{n \sigma^2} \sum\nolimits_{i=1}^n \xi_i \mathcal{X}_i+ \frac{1}{n \sigma^2} \sum\nolimits_{i=1}^n \big[\langle\widehat{\Delta}, \mathcal{X}_i\rangle \mathcal{X}_i-\sigma^2\cdot \widehat{\Delta} \big]  := \mcT+\widehat{\mcZ}^{(1)} + \widehat{\mcZ}^{(2)}, \label{eq: debiased estimator and its decomposition in tensor regression}
\end{align}
with $\widehat{\mcZ} = \widehat{\mcZ}^{(1)} + \widehat{\mcZ}^{(2)}$. Here, $\widehat{\mcZ}^{(1)}$ serves as the candidate pivot quantity, and $\widehat{\mcZ}^{(2)}$ accounts for the bias correction due to the initial estimation error.

\subsubsection*{Step 3: Two-step Power Iteration} 

Using the initial singular space estimates $\widehat{U}_j^{(0)}$ for $j = 1, 2, 3$, we refine the singular space estimates through a two-step power iteration. Specifically, $\widehat{U}_j^{(1)}$ and $\widehat{U}_j^{(2)}$ represent the estimates after the first and second iterations, respectively. For each iteration $k = 1, 2$ and mode $j = 1, 2, 3$, the power iteration is performed as follows:

For each mode $j=1,2,3$, $\whU_j^{(k)}$ is obtained as the leading $r_1$ left singular vectors of
$$
\Mat_j\big( \widehat{\mcT}^{\text{unbs}} \times_{j+1} \whU_{j+1}^{(k-1)\top} \times_{j+2} \whU_{j+2}^{(k-1)\top} \big) = \Mat_j\big( \widehat{\mcT}^{\text{unbs}} \big) \big( \whU_{j+2}^{(k-1)} \otimes \whU_{j+1}^{(k-1)} \big).
$$

After completing the iterations, the final estimates are set as $\widehat{U}_j := \widehat{U}_j^{(2)}$ for each mode $j$. 

\subsubsection*{Step 4: Projection and Plug-in Estimator }

With the refined singular space estimates $\widehat{U}_j = \widehat{U}_j^{(2)}$ for $j = 1, 2, 3$ and their corresponding projection matrices $\mcP_{\widehat{U}_j} = \widehat{U}_j \widehat{U}_j^{\top}$, we compute the projected tensor as
\begin{align}
\widehat{\mcT} = \widehat{\mcT}^{\text{unbs}} \times_1 \mcP_{\whU_1} \times_2 \mcP_{\whU_2} \times_3 \mcP_{\whU_3}. \label{eq: definition of whT}
\end{align}
The linear functional $\langle \mcA, \mcT \rangle$ is then estimated by $\langle \mcA, \widehat{\mcT} \rangle$.

\begin{remark}\label{remark: role of projection}
The projection step is critical for both algorithmic performance and theoretical guarantees. A candidate pivot quantity,
$\widehat{\mcZ}^{(1)} = (n\sigma^2)^{-1} \sum_{i=1}^n \xi_i \mcX_{i}$, is introduced after bias correction.
However, directly using $\langle \mcA, \widehat{\mcZ}^{(1)} \rangle$ results in suboptimal confidence interval lengths, as
$\widehat{\mcZ}^{(1)}$ resides in the high-dimensional space $\mathbb{R}^{p_1 \times p_2 \times p_3}$ rather than the low-Tucker-rank manifold $\mathcal{M}_{(r_1, r_2, r_3)}$.
To address this, we project $\widehat{\mcZ}^{(1)}$ onto the estimated singular spaces via $\mcP_{\widehat{U}_j}$, restricting it to the low-Tucker-rank manifold. Since low-Tucker-rank manifolds lack certain structural properties, the tangent space $\mathbb{T}_{\mcT}\mathcal{M}_{(r_1, r_2, r_3)}$ at the signal tensor $\mcT$ on $\mathcal{M}_{(r_1, r_2, r_3)}$ serves as a first-order approximation. Confidence intervals are constructed by 
projecting the candidate pivot quantity $\widehat{\mcZ}^{(1)}$ onto the tangent space at the true parameter tensor $\mcT$ within the manifold. Then, the projection in the tangent space $\mcP_{\mathbb{T}_{\mcT}\mathcal{M}_{(r_1, r_2, r_3)}}(\widehat{\mcZ}^{(1)})$ is the pivot quantity for constructing the confidence interval.

\end{remark}

\subsection{Asymptotic Normality of Estimated Linear Functionals without Sample Splitting} \label{sec: Asymptotic Normality of Estimated Linear Functional without Sample Splitting}

In this section, we establish the asymptotic normality of the estimator $\langle \mcA, \widehat{\mcT} \rangle$ obtained from Section \ref{sec:debias_nonsplit}.

\begin{assumption}[Sample size requirement]\label{assump: sample size requirement in tensor regression}

The sample size $n$ for tensor regression satisfies $n \geq C\max\{\kappa^2 \op / \ulambda^2, \op\oor \}$, where $\op = \max\{p_1, p_2, p_3\}, \oor = \max\{r_1, r_2, r_3\}$, $C$ is a constant depending on the noise scales $\sigma_\xi$ and $\sigma$ as defined in Assumption~\ref{assump: subgaussian noise in tensor regression}, and $\kappa = \olambda / \ulambda$ is the condition number of the signal tensor $\mcT$, $\ulambda,\olambda$ are defined in \eqref{def:ulambda} and \eqref{def:mlambda}, respectively. To simplify the presentation of our results, we assume $\op \gtrsim \oor^2$ throughout this paper. 
\end{assumption}

\begin{remark} 
This assumption ensures two key requirements. (i) Spectral Representation Validity: The condition $n \gtrsim \kappa^2 \op/ \ulambda^2$, supports the spectral representation necessary for asymptotic normality, as in Theorem 1 of \citet{xia2021normal}. (ii) Degrees of Freedom: The sample size must scale with the tensor's degrees of freedom, $r_1 r_2 r_3 + \sum_{j=1}^3 p_j r_j$, which is statistically optimal for inference. Notably, larger singular values ($\ulambda$) reduce the sample size requirement.

\end{remark}

With these assumptions, we proceed to establish the asymptotic normality of $\langle \mcA, \widehat{\mcT} \rangle$.

\begin{theorem}[Main Theorem: Asymptotic Normality in Tensor Regression]\label{thm: main theorem in tensor regression without sample splitting}

Consider the low-Tucker-rank tensor regression model \eqref{eq: tensor regression model}. Suppose that Assumptions~\ref{assump: structural assumption on the signal tensor and loading tensor}-\ref{assump: sample size requirement in tensor regression} hold, and assume that $\|\mcT\|_{\mathrm{F}}^2 \geq C_1 \sigma_\xi^2/\sigma^2 $ and $\ulambda \geq c_1$ for some constants $C_1, c_1 > 0$. Let $\widehat{\mcT}$ and $\widehat{U}_j$ denote the outputs of the debiasing procedure in Section \ref{sec:debias_nonsplit}. Then, the estimator satisfies
\begin{align*}
& \sup _{x \in \mathbb{R}} \left\lvert\, \mathbb{P}\left(\frac{\sqrt{n} \big( \langle \widehat{\mcT}, \mcA \rangle-\langle \mcT, \mcA\rangle \big) }{(\sigma_{\xi}/\sigma) s_{\mcA} } \leq x\right)-\Phi(x)\right| \\
\lesssim & \underbrace{\sqrt{\frac{1}{n}}}_{\substack{\text{rate of asymptotic}\\ \text{normal terms}}} + \underbrace{\frac{ \Omega_1 + \Omega_2 + \Omega_3 + \Omega_4 }{(\sigma_{\xi}/\sigma)  s_{\mcA} \sqrt{\frac{1}{n}}}}_{\text{rate of negligible terms}} + \underbrace{\left[ \op^{-c} + e^{-cn} + \mcP_{\mcE_{U}^{\text{reg}}} + \mcP_{\mcE_{\Delta}}\right]}_{\text{rate of initial estimates}},
\end{align*}
where $c$ is a positive constant, the variance component $s_{\mcA}$ is defined in \eqref{eq: variance component sA2}, and $\Omega_1, \Omega_2, \Omega_3, \Omega_4$ are upper bounds for various negligible error terms:
{\small \begin{align*}
\Omega_1
= & \|\mcA\times_1 U_1 \times_2 U_2 \times_3 U_3 \|_{\mathrm{F}} \cdot \frac{\sigma_\xi^2}{\sigma^2} \left(\frac{\sqrt{\op \oor\log(\op)}}{n \ulambda} +\Delta\cdot \frac{\op \oor^{1/2}}{n \ulambda}\right), \\
\Omega_2
= & \sum_{j=1}^3 \big\|\mcP_{U_{j\perp}} A_j \mcP_{(U_{j+2} \otimes U_{j+1} ) G_j^{\top}} \big\|_{\mathrm{F}} \cdot  \frac{\sigma_\xi^2}{\sigma^2} \left(\frac{\sqrt{\op \oor^2 \log (\op)}}{n \ulambda }+\Delta \cdot \frac{\op\oor^{1/2}}{n \ulambda}\right) \\
+ & \sum_{j=1}^3\big\|\mcA \times_j U_j\big\|_{\mathrm{F}} \cdot  \frac{\sigma_\xi^2}{\sigma^2} \left(\frac{\oor^{3/2} \log (\op)}{n \ulambda }+\Delta \cdot \frac{\sqrt{ \op \oR \oor \log (\op)}}{n \ulambda}+\Delta^2 \cdot \frac{\op \oor^{1/2}}{n \ulambda }\right) \\
+ & \big\|\mcA\big\|_{\mathrm{F}} \cdot  \frac{\sigma_\xi^3}{\sigma^3} \left(\frac{\oor^{3/2} \log (\op)^{3 / 2}}{n^{3 / 2} \ulambda^2 }+\Delta \cdot \frac{\op^{1 / 2} \oR \oor^{1/2}\log (\op)}{n^{3 / 2} \ulambda^2}+\Delta^2 \cdot \frac{\op \oR \oor^{1/2}\log (\op)}{n^{3 / 2} \ulambda^2 }+\Delta^3 \cdot \frac{\op^{3 / 2} \oor^{1/2} }{n^{3 / 2} \ulambda^2}\right), \\
\Omega_3
= & \sum_{j=1}^3 \big\|\mcP_{U_j} A_j \mcP_{(U_{j+2} \otimes U_{j+1}) G_j^{\top}}\big\|_{\mathrm{F}} \cdot  \frac{\sigma_\xi^2}{\sigma^2}\cdot \frac{\op \oor^{1/2}}{n\ulambda} \\
+ & \sum_{j=1}^3\big\|\mcA \times_{j+1} U_{j+1} \times_{j+2} U_{j+2}\big\|_{\mathrm{F}} \cdot \frac{\sigma_\xi^3}{\sigma^3} \left(\frac{\op\oor^{1/2}\log(\op)^{1/2}}{n^{3/2} \ulambda^2 }+\Delta\cdot \frac{\op^{3/2} \oor^{1/2}}{n^{3/2} \ulambda^2 }\right), \\
\Omega_4 =&  \sum_{j=1}^3\big\|\mcP_{U_{j\perp}}A_j\mcP_{(U_{j+2}\otimes U_{j+1} )G_j^{\top}}\big\|_{\mathrm{F}}\cdot \Delta \cdot \frac{\sigma_\xi}{\sigma} \sqrt{\frac{\op \oor }{n}} + \big\|\mcA\times_1 U_1\times_2 U_2\times_3 U_3\big\|_{\mathrm{F}}\cdot \Delta \cdot \frac{\sigma_\xi}{\sigma} \sqrt{\frac{\op\oor }{n}}.
\end{align*}}%
Here, $\oR=\max\{ R_1,R_2,R_3\}$ is allowed to divergent, and $\Delta$ denotes the initial error.
\end{theorem}

The theorem asserts that the estimator $\langle \widehat{\mcT}, \mcA \rangle$ is asymptotically normal, centered at the true linear functional $\langle \mcT, \mcA \rangle$, and scaled by the variance term $(\sigma_{\xi}/\sigma)  s_{\mcA} $. The variance component $s_{\mcA}$ reflects the variability introduced by projecting both the loading tensor $\mcA$ and the noise tensor $\widehat{\mcZ}^{(1)} = (n\sigma^2)^{-1} \sum_{i=1}^n \xi_i \mcX_i$ onto the tangent space of the low-Tucker-rank manifold $\mathcal{M}_{(r_1, r_2, r_3)}$ at the true tensor $\mcT$. The asymptotic normal term $\langle \mcP_{\mathbb{T}_{\mcT} \mathcal{M}_{\mathrm{r}}} (\mathcal{A}), \mcP_{\mathbb{T}_{\mcT} \mathcal{M}_{\mathrm{r}}}(\widehat{\mathcal{Z}}^{(1)}) \rangle$, as introduced in Remark~\ref{remark: role of projection}, represents the Riemann metric on the tangent space $\mathbb{T}_{\mcT}$ between the loading tensor $\mcA$ and the candidate pivot quantity $\widehat{\mathcal{Z}}^{(1)}$, defined in \eqref{eq: debiased estimator and its decomposition in tensor regression}. This projection leverages a first-order approximation of $\widehat{\mcT}$ around $\mcT$ within the manifold, facilitating the normal approximation.

\begin{remark} \label{rmk:regression_errors}
The error terms $\Omega_1$, $\Omega_2$, $\Omega_3$, and $\Omega_4$ represent various sources of negligible errors. Specifically, $\Omega_1$ accounts for negligible errors introduced by the artificial noise $\widehat{\mcZ}^{(1)}$ during bias correction, $\Omega_2$ captures errors arising from projection onto the estimated singular spaces, $\Omega_3$ provides a common upper bound for errors from both bias correction and projection onto the estimated singular spaces. Among these, $\Omega_4$, introduced by the initialization error in bias correction, typically dominates the negligible terms. The dependence between the initial estimate and the bias-correction complicates the analysis and requires a stricter sample size to ensure asymptotic normality.

\end{remark}


Theorem~\ref{thm: main theorem in tensor regression without sample splitting} establishes conditions for valid inference of the linear functional $\langle \mcA, \mcT \rangle$, even when the loading tensor $\mcA$ has full Tucker rank. These results are summarized in the following corollary.

\begin{corollary}[Asymptotic normality of estimated general linear functionals] \label{corollary: asymptotic normality of estimated general linear functional in tensor regression without sample splitting}

Under the conditions of Theorem~\ref{thm: main theorem in tensor regression without sample splitting}, assume the initial estimate is minimax optimal, i.e., $\Delta \asymp \sqrt{ \op\oor/n }$. Further, suppose the sample size satisfies $n \geq C_1\max \{\kappa^2 \op / \ulambda^2, \op^2\oor \}$. Additionally, assume the following alignment condition holds
\begin{align}
s_{\mcA} \ge C_2\max\left\{\oor\op^{-1/2}\ulambda^{-1} \|\mcA \times_j U_j \|_{\mathrm{F}}, \oor\op^{-1}\ulambda^{-2} \|\mcA\|_{\mathrm{F}}\right\}, \label{eq: general alignment condition in tensor regression without sample splitting}
\end{align}
where $C_1$ and $C_2$ are two constants depending only on the noise scales $\sigma_\xi$ and $\sigma$. Then, for any loading tensor $\mcA \in \mathbb{R}^{p_1\times p_2\times p_3}$, the estimator satisfies
$
\sqrt{n} ( \langle \widehat{\cT} , \mathcal{A} \rangle- \langle \mcT, \mathcal{A} \rangle )/ (\sigma_\xi/\sigma \cdot s_{\mcA}) \stackrel{d}{\rightarrow} \mathcal{N}(0,1)  .  
$

\end{corollary}

The minimax optimal initial estimation error $\Delta$ and the sample size requirement $n \gtrsim \max\{\op^2 \oor, \kappa^2 \op / \ulambda^2\}$ are crucial for mitigating the leading error term caused by the dependence between the initial estimate and the bias-correction step. Minimax optimal initial estimators can be achieved using methods such as projected gradient descent \citep{han2022optimal} or sketching \citep{zhang2020islet} within the tensor regression framework.

For general linear functionals where the loading tensor $\mcA$ has full Tucker rank, no additional incoherence conditions are needed if the sample size is sufficiently large. The alignment condition ensures that the variance component $s_{\mcA}$ is sufficiently large, which is facilitated by a strong signal tensor $\mcT$ with a larger minimum singular value $\ulambda$. A stronger signal leads to a more favorable alignment condition. Further discussion on the interplay between incoherence and alignment conditions is provided in Section \ref{sec: The Interplay between the Incoherence Condition and Alignment Condition}.

The sub-Gaussian design tensor assumption (Assumption~\ref{assump: subgaussian design tensor in tensor regression}) aligns with the sub-Gaussian sampling framework often used in compressed sensing \citep{carpentier2019uncertainty}. Notably, our results demonstrate that under sub-Gaussian designs, valid inference for general linear functionals is achievable without requiring the sample size to scale with the total number of entries ($p_1 p_2 p_3$) in the parameter tensor $\mcT$, provided the initial estimates are accurate. This highlights the data compression benefits of low-Tucker-rank modeling, in contrast to low-rank matrix trace regression where $n\gtrsim p_1 p_2$ is typically required for inference with a coefficient matrix $M\in\R^{p_1\times p_2}$ \citep{cai2016geometric}. 

\subsection{Debiased Estimator of Linear Functionals with Sample Splitting} \label{sec: Debiased Estimator of Linear Functionals with Sample Splitting}

In Theorem~\ref{thm: main theorem in tensor regression without sample splitting} and Corollary~\ref{corollary: asymptotic normality of estimated general linear functional in tensor regression without sample splitting}, we established that the general linear functional $\langle \mcA, \mcT \rangle$ can be inferred with a sample size of $n\gtrsim \op^2$, assuming $\oor\asymp 1$. While this guarantees valid inference, 
the required sample size falls short of the computationally optimal rate $n\gtrsim \op^{3/2}$. This limitation arises from reusing data for both bias correction and initial estimation, which constrains the convergence rate. 

In tensor regression, there exists a gap between the computationally optimal sample size $n\gtrsim \op^{3/2}$, as highlighted in Remark 4.3 of \citet{han2022optimal}, and the lower bound of the sample size determined by the degrees of freedom for a tensor with Tucker rank-$(r_1, r_2, r_3)$, given by $n \gtrsim r_1 r_2 r_3 + \sum_{j=1}^3 p_j r_j$. This discrepancy raises the natural question: {Can the statistically optimal sample size be achieved while ensuring valid inference?}

In the following section, we address this affirmatively by employing a sample-splitting strategy. This approach eliminates the dependence between the initial estimation and the debiasing step.
Specifically, the observed data is divided into two disjoint subsets: dataset $\mathrm{\RN{1}}$ and dataset $\mathrm{\RN{2}}$. One subset is used for initial estimation, and the other for debiasing, ensuring that the debiasing process is independent of the data used for obtaining initial estimates. The detailed algorithm for tensor regression with sample splitting is provided below.


\subsubsection*{Step 1: Initialization}
Using the first dataset $\mathrm{\RN{1}}:=\{Y_{i_1}^{(\mathrm{\RN{1}})}, X_{i_1}^{(\mathrm{\RN{1}})}\}_{i_1=1}^{n_1}$, apply the initial estimation procedure to obtain the initial (typically biased) tensor estimate $\widehat{\mcT}^{\text{init},(\mathrm{\RN{1}})} \in \mathbb{R}^{p_1 \times p_2 \times p_3}$ and initial estimates of the factor matrices $\whU_1^{\text{init},(\mathrm{\RN{1}})}$, $\whU_2^{\text{init},(\mathrm{\RN{1}})}$, $\whU_3^{\text{init},(\mathrm{\RN{1}})}$. Similarly, using the second dataset, $\mathrm{\RN{2}}:=\{Y_{i_2}^{(\mathrm{\RN{2}})}, X_{i_2}^{(\mathrm{\RN{2}})}\}_{i_2=1}^{n_2}$, apply the same procedure to obtain the tensor estimate $\widehat{\mcT}^{\text{init},(\mathrm{\RN{2}})} \in \mathbb{R}^{p_1 \times p_2 \times p_3}$ and factor matrices $\whU_1^{\text{init},(\mathrm{\RN{2}})}$, $\whU_2^{\text{init},(\mathrm{\RN{2}})}$, $\whU_3^{\text{init},(\mathrm{\RN{2}})}$.
Set $\whU_j^{(0),(\mathrm{\RN{1}})}:=\whU_j^{\text{init},(\mathrm{\RN{1}})}, \whU_j^{(0),(\mathrm{\RN{2}})}:=\whU_j^{\text{init},(\mathrm{\RN{2}})}$ for $j=1,2,3$.

\subsubsection*{Step 2: Debiasing} 
The initial estimates are debiased using the complementary dataset as follows,
\begin{align*}
\widehat{\mcT}^{\text{unbs}, (\mathrm{\RN{1}})}
= & \widehat{\mcT}^{\text{init},(\mathrm{\RN{2}})}+\frac{1}{n_1\sigma^2} \sum\nolimits_{i_1=1}^{n_1}\big(Y_{i_1}^{(\mathrm{\RN{1}})}-\langle\widehat{\mcT}^{\text{init},(\mathrm{\RN{2}})}, \mathcal{X}_{i_1}^{(\mathrm{\RN{1}})}\rangle \big) \mathcal{X}_{i_1}^{(\mathrm{\RN{1}})} ,\\
\widehat{\mcT}^{\text{unbs}, (\mathrm{\RN{2}})}
= & \widehat{\mcT}^{\text{init},(\mathrm{\RN{1}})}+\frac{1}{n_2\sigma^2} \sum\nolimits_{i_2=1}^{n_2}\big(Y_{i_2}^{(\mathrm{\RN{2}})}-\langle\widehat{\mcT}^{\text{init},(\mathrm{\RN{1}})}, \mathcal{X}_{i_2}^{(\mathrm{\RN{2}})}\rangle \big) \mathcal{X}_{i_2}^{(\mathrm{\RN{2}})} .
\end{align*}
Here, $\widehat{\Delta}^{(\mathrm{\RN{1}})} = \mcT - \widehat{\mcT}^{\text{init},(\mathrm{\RN{1}})}$ and $\widehat{\Delta}^{(\mathrm{\RN{2}})} = \mcT - \widehat{\mcT}^{\text{init},(\mathrm{\RN{2}})}$ represent the estimation errors from Dataset \RN{1} and Dataset \RN{2}, respectively. The debiasing step uses one dataset to correct the bias in the estimates obtained from the other dataset, effectively eliminating the dependence between the initial estimates and the bias correction.

\subsubsection*{Step 3: One-step Power Iteration}
Using the initial estimates $\whU_j^{(0)}$, $j=1,2,3$, we perform a one-step power iteration to refine the estimates of the singular space. Specifically,
$\whU_j^{(1),(\mathrm{\RN{1}})}$ and $\whU_j^{(1), (\mathrm{\RN{2}})}$ are obtained as the leading $r_j$ left singular vectors of
$$
\Mat_j\big(\widehat{\mcT}^{\text{unbs}, (\mathrm{\RN{1}})} \big) \big( \whU_{j+2}^{(0), (\mathrm{\RN{2}})} \otimes \whU_{j+1}^{(0), (\mathrm{\RN{2}})} \big), \quad \text{ and} \quad \Mat_j\big( \widehat{\mcT}^{\text{unbs}, (\mathrm{\RN{2}})} \big) \big( \whU_{j+2}^{(0), (\mathrm{\RN{1}})} \otimes \whU_{j+1}^{(0), (\mathrm{\RN{1}})} \big),
$$
for mode $j=1,2,3$, respectively. 

Unlike the algorithm without sample splitting in Section \ref{sec:debias_nonsplit}, a single iteration suffices due to the independence introduced by sample splitting, which simplifies the perturbation analysis.

\subsubsection*{Step 4: Projection and plug-in Estimator}
The final projected estimator is computed by averaging the contributions from both datasets,
\begin{align*}
&\widehat{\mcT} = \frac{n_1}{n}\widehat{\mcT}^{\text{unbs},(\mathrm{\RN{1}})} \times_1 \mcP_{\whU_{1}^{(\mathrm{\RN{1}})}} \times_2 \mcP_{\whU_{2}^{(\mathrm{\RN{1}})}} \times_3 \mcP_{\whU_{3}^{(\mathrm{\RN{1}})}}+ \frac{n_2}{n}\widehat{\mcT}^{\text{unbs},(\mathrm{\RN{2}})} \times_1 \mcP_{\whU_{1}^{(\mathrm{\RN{2}})}} \times_2 \mcP_{\whU_{2}^{(\mathrm{\RN{2}})}} \times_3 \mcP_{\whU_{3}^{(\mathrm{\RN{2}})}},
\end{align*}
where $n=n_1+n_2$. Finally, the linear functional $\langle\mcA, \mcT\rangle$ is estimated by $\langle\mcA, \widehat{\mcT}\rangle$.

\subsection{Asymptotic Normality of Estimated Linear Functionals with Sample Splitting} \label{sec: Asymptotic Normality of Estimated Linear Functional with Sample Splitting}

In this section, we establish the asymptotic normality of the estimator $\langle \mcA, \widehat{\mcT} \rangle$ derived from the bias-correction procedure with sample-splitting.

\begin{theorem}\label{thm: main theorem in tensor regression with sample splitting}
Consider the low-Tucker-rank tensor regression model \eqref{eq: tensor regression model}. Suppose that Assumptions~\ref{assump: structural assumption on the signal tensor and loading tensor}-\ref{assump: sample size requirement in tensor regression} hold for each sub-dataset and its corresponding initial estimators. Assume that $\|\mcT\|_{\mathrm{F}}^2 \geq C_1 \sigma_\xi^2/\sigma^2 $ and $\ulambda \geq c_1$ for some constants $C_1, c_1 > 0$, 
and let $n_1, n_2 \asymp n$ with $n_1 + n_2 = n$. Let $\widehat{\mcT}$ and $\widehat{U}_j$ be the outputs of the debiasing procedure in Section \ref{sec: Debiased Estimator of Linear Functionals with Sample Splitting}. Then, the estimator satisfies
\begin{align*}
& \sup _{x \in \mathbb{R}} \left\lvert\, \mathbb{P}\left(\frac{\sqrt{n} \big( \langle \widehat{\mcT}, \mcA \rangle-\langle \mcT, \mcA\rangle \big) }{(\sigma_{\xi}/\sigma) s_{\mcA} } \leq x\right)-\Phi(x)\right| \\
\lesssim & \underbrace{\sqrt{\frac{1}{n}}}_{\substack{\text{rate of asymptotic}\\ \text{normal terms}}} + \underbrace{\frac{ \Omega_1 + \Omega_2 + \Omega_3 }{(\sigma_{\xi}/\sigma)  s_{\mcA}  \sqrt{\frac{1}{n}}}}_{\text{rate of negligible terms}} + \underbrace{\left[\op^{-c} + e^{-cn}  + \mcP_{\mcE_{U}^{\text{reg}}} + \mcP_{\mcE_{\Delta}}\right]}_{\text{rate of initial estimates}},
\end{align*}
where $c > 0$ is a constant, the variance component $s_{\mcA}$ is defined in \eqref{eq: variance component sA2}, and $\Omega_1, \Omega_2, \Omega_3$ are upper bounds for various negligible error terms:
{\small \begin{align*}
\Omega_1
= & \big\|\mcA\times_1 U_1 \times_2 U_2 \times_3 U_3 \big\|_{\mathrm{F}} \left(\frac{\sigma_\xi}{\sigma} \cdot \Delta \cdot \sqrt{\frac{\oor^2\log(\op)}{n}} + \frac{\sigma_\xi^2}{\sigma^2} \cdot \frac{\op^{1/2}\oor^{1/2}\log(\op)^{1/2}}{n \ulambda} + \frac{\sigma_\xi^4}{\sigma^4} \cdot \frac{\op^2\oor^{1/2}}{n^2\ulambda^3} \right), \\
\Omega_2
= & \sum_{j=1}^3 \big\|\mcP_{U_{j\perp}} A_j \mcP_{(U_{j+2} \otimes U_{j+1}) G_j^{\top}} \big\|_{\mathrm{F}} \cdot \left(\frac{\sigma_\xi}{\sigma} \cdot \Delta  \sqrt{\frac{\oor^2\log(\op)}{n}} + \frac{\sigma_{\xi}^3}{\sigma^3} \cdot \frac{\op \oor \log (\op)^{1/2}}{n^{3 / 2}\ulambda^2}\right) \\
+ & \sum_{j=1}^3 \big\|\mcA \times_j U_j \big\|_{\mathrm{F}} \cdot \frac{\sigma_{\xi}^2}{\sigma^2}  \left(\frac{\oor^{3/2} \log (\op)}{n \ulambda} + \Delta \cdot \frac{\oR \oor^{1/2} \log (\op)}{n \ulambda}\right) \\
+ & \|\mcA\|_{\mathrm{F}} \cdot \frac{\sigma_{\xi}^3}{\sigma^3}   \left(\frac{\oor^2 \log (\op)^{3 / 2}}{n^{3 / 2} \ulambda^2} + \Delta \cdot \frac{\oR^{3 / 2}\oor^{1/2} \log (\op)^{3 / 2}}{n^{3 / 2} \ulambda^2 }\right), \\
\Omega_3
= & \sum_{j=1}^3\big\|\mcP_{U_j}A_j\mcP_{(U_{j+2}\otimes U_{j+1})G_j^{\top}}\big\|_{\mathrm{F}} \cdot \frac{\sigma_\xi^2}{\sigma^2}\cdot \frac{\op\oor^{1/2}}{n \ulambda} \\
+ & \sum_{j=1}^3\big\|\mcA \times_{j+1} U_{j+1} \times_{j+2} U_{j+2}\big\|_{\mathrm{F}} \cdot \frac{\sigma_\xi^3}{\sigma^3} \cdot \frac{\op\oor\log(\op)^{1/2}}{n^{3/2} \ulambda^2} .  
\end{align*}}
Here, $\oR=\max\{ R_1,R_2,R_3\}$ is allowed to divergent, and $\Delta$ denotes the initial error.
\end{theorem}

The double-sample-splitting debiasing procedure outlined in Section \ref{sec: Asymptotic Normality of Estimated Linear Functional with Sample Splitting} partitions the data into independent subsets for initial estimation and bias correction. This separation eliminates the dependence between these steps, reducing the leading error term $\Omega_4$ identified in Theorem~\ref{thm: main theorem in tensor regression without sample splitting}, which stems from repeated data use. As a result, sample splitting enables both computationally and statistically optimal sample sizes under appropriate conditions.

The following corollary establishes the detailed conditions required for achieving the asymptotic normality of the estimated low-Tucker-rank linear form, with computationally and statistically optimal sample sizes. Notably, it shows that the initial estimate need not attain the minimax optimal rate for valid low-rank linear functional inference. 

\begin{corollary}[Asymptotic normality of estimated low-Tucker-rank linear functionals] \label{corollary: asymptotic normality of estimated low-Tucker-rank linear functional in tensor regression with sample splitting}

Under the conditions of Theorem~\ref{thm: main theorem in tensor regression with sample splitting}, assume the Tucker rank of the loading tensor $\operatorname{rank}(\mcA)=(R_1, R_2, R_3)$ is fixed and independent of $\op$. Given the sample size requirement $n \geq C\max \{\op\oor, \kappa^2 \op/\ulambda^2 \}$, where $C>0$ is a constant depending only on $R$ and the noise scales $\sigma_\xi$ and $\sigma$, as specified in Assumption~\ref{assump: sample size requirement in tensor regression}, the following holds
\begin{align*}
& \sup _{x \in \mathbb{R}} \left\lvert\, \mathbb{P}\left(\frac{\sqrt{n} \big( \langle \widehat{\mcT}, \mcA \rangle-\langle \mcT, \mcA\rangle \big) }{(\sigma_{\xi}/\sigma) s_{\mcA} } \leq x\right)-\Phi(x)\right| \\
\lesssim & \Delta \cdot \oor \sqrt{\log(\op)} + [\op^{-c} + e^{-cn} + \mcP_{\mcE_{U}^{\text{reg}}} + \mcP_{\mcE_{\Delta}}] \\
& +  \frac{1}{(\sigma_{\xi}/\sigma) s_{\mcA} \cdot \sqrt{\frac{1}{n}}} \cdot \Bigg\{\sum_{j=1}^3 \big\|\mcP_{U_j} A_j \mcP_{(U_{j+2} \otimes U_{j+1}) G_j^{\top}}\big\|_{\mathrm{F}} \cdot \frac{\sigma_{\xi}^2}{ \sigma^2} \cdot \frac{\op\oor^{1/2}}{n \ulambda} \\
& +  \big\|\mcA\times_1 U_1 \times_2 U_2 \times_3 U_3\big\|_{\mathrm{F}} \cdot \left(\frac{\sigma_\xi^2}{\sigma^2} \cdot \frac{\oor^{1/2}\op^{1/2}\log(\op)^{1/2}}{n \ulambda} + \frac{\sigma_\xi^4}{\sigma^4} \cdot \frac{\op^2\oor^{1/2}}{n^2\ulambda^3} \right) \\
& +  \sum_{j=1}^3\big\|\mcA \times_{j+1} U_{j+1} \times_{j+2} U_{j+2}\big\|_{\mathrm{F}} \cdot \frac{\sigma_\xi^3}{\sigma^3} \left(\frac{\op\oor\log(\op)^{1/2}}{n^{3/2} \ulambda^2} \right)  \\
& +  \sum_{j=1}^3 \big\|\mcA \times_j U_j \big\|_{\mathrm{F}} \cdot \frac{\sigma_{\xi}^2}{ \sigma^2} \cdot \frac{\oor^{3/2}\log (\op)}{n \ulambda}  + \|\mcA \|_{\mathrm{F}} \cdot \frac{\sigma_{\xi}^3}{ \sigma^3} \cdot \frac{\oor^2 \log (\op)^{3 / 2}}{n^{3 / 2} \ulambda^2} \Bigg\} ,
\end{align*}
where $c > 0$ is a constant, the variance component $s_{\mcA}$ is defined in \eqref{eq: variance component sA2}.
\end{corollary}

{\it Statistical optimal sample size w.r.t. $\op$.} If the following {\it incoherence condition} \eqref{eq: statistically optimal incoherence condition in tensor regression with sample splitting} holds,
\begin{align*}
\max_{j=1,2,3}\big\|\mcP_{U_j}A_j (\mcP_{U_{j+2}}\otimes \mcP_{U_{j+1}} ) \big\|_{\mathrm{F}} / \big\|\mcP_{U_{j\perp}}A_j(\mcP_{U_{j+2}}\otimes \mcP_{U_{j+1}} ) \big\|_{\mathrm{F}} \leq c_2\ulambda \oor^{1/2} \op^{-1/2} 
\end{align*}
and the following {\it alignment condition} holds,
{\small \begin{align}
s_{\mcA} \geq C_2\max_{j=1,2,3}\big\{\ulambda^{-2}\|\mcA \times_{j+1} U_{j+1} \times_{j+2} U_{j+2}\|_{\mathrm{F}}, \oor^{1/2}\op^{-1/2}\ulambda^{-1} \|\mcA\times_j U_j \|_{\mathrm{F}}, \op^{-1}\ulambda^{-2} \|\mcA \|_{\mathrm{F}}\big\}, \label{eq: statistically optimal alignment condition in tensor regression with sample splitting} 
\end{align}}%
then a sample size of $n \geq C\max \{\kappa^2 \op / \ulambda^2, \op\oor^2\log(\op)^{3/2} \}$ is sufficient for asymptotic normality of the estimated linear functional, where $C$, $C_2$ and $c_2$ are positive constants depending only on $\oR$, the fixed rank of the loading tensor $\mcA$, and the noise scales $\sigma_\xi$ and $\sigma$. When $\Delta = o[(\oor\log(\op) )^{-1}]$, the estimator satisfies
$
\sqrt{n} ( \langle \widehat{\cT} , \mathcal{A} \rangle-  \langle \mcT, \mathcal{A} \rangle )/  (\sigma_\xi/\sigma \cdot s_{\mcA} ) \stackrel{d}{\rightarrow} \mathcal{N}(0,1)    
$
as $n,\op \rightarrow \infty$.

{\it Computational optimal sample size w.r.t. $\op$.} If the following {\it incoherence condition} \eqref{eq: computationally optimal incoherence condition in tensor regression with sample splitting} holds,
\begin{align*}
\max_{j=1,2,3}\big\|\mcP_{U_j}A_j\mcP_{(U_{j+1}\otimes U_{j+1})G_j^{\top}} \big\|_{\mathrm{F}} / \big\|\mcP_{U_{j\perp}}A_j\mcP_{(U_{j+1}\otimes U_{j+1})G_j^{\top}} \big\|_{\mathrm{F}} \leq c_2 \ulambda \oor^{1/2}\op^{-1/4} ,
\end{align*}
and the following {\it alignment condition} holds,
{\small \begin{align}
s_{\mcA} \geq C_2\max_{j} \Big\{
& \op^{-1/2}\ulambda^{-1}\|\mcA \times_{j+1} U_{j+1} \times_{j+2} U_{j+2}\|_{\mathrm{F}}, \oor \op^{-3/4}\ulambda^{-1}\|\mcA \times_j U_j\|_{\mathrm{F}}, \oor\op^{-3/2}\ulambda^{-2}\|\mcA\|_{\mathrm{F}}\Big\}, \label{eq: computationally optimal alignment condition in tensor regression with sample splitting} 
\end{align}}%
then a sample size of $n \geq C\max \{\kappa^2 \op / \ulambda^2, \op^{3/2}\oor\log(\op)^{3/2} \}$ is sufficient for valid inference, where $C$, $C_2$ and $c_2$ are positive constants depending only on $\oR$, the fixed rank of the loading tensor $\mcA$, and the noise scales $\sigma_\xi$ and $\sigma$. When $\Delta = o[(\oor\log(\op))^{-1}]$, the estimator satisfies
$
\sqrt{n} ( \langle \widehat{\cT} , \mathcal{A} \rangle-  \langle \mcT, \mathcal{A} \rangle  )/  (\sigma_\xi/\sigma \cdot s_{\mcA} ) \stackrel{d}{\rightarrow} \mathcal{N}(0,1) 
$
as $n,\op \rightarrow \infty$.

\begin{remark}[Role of $\ulambda$]
So far, our discussion of computational and statistical optimality has focused on the case $\ulambda\asymp 1$. In the tensor regression setting, our analysis reveals that increasing the signal strength $\ulambda$ relaxes the requirements for incoherence and alignment conditions necessary for valid inference. For example, if $\ulambda \gtrsim \op^{1/4}$ and $\oor \asymp 1, \kappa \asymp 1$, a sample size $n\gtrsim \op^{3/2}$ suffices for low-rank linear functional inference without requiring an incoherence condition. If $\ulambda \gtrsim \op^{1/2}$ and $\oor \asymp 1, \kappa \asymp 1$, even sample size $n\gtrsim \op$ is sufficient for valid inference, again with no incoherence condition needed. However, the role of $\ulambda$ in determining estimation quality and the computationally/statistically optimal sample size (lower bounds) in the tensor regression framework is underexplored in the existing literature. Similar results have been reported in \citep{zhang2020islet}. In practice, researchers often have limited prior knowledge of the signal strength of $\mcT$, highlighting a gap that warrants further investigation.
\end{remark}

\subsection{The Incoherence and Alignment Conditions: A Geometric Interpretation on the Low-Tucker-Rank Manifold} \label{sec: The Interplay between the Incoherence Condition and Alignment Condition}

To clarify the roles of the incoherence and alignment conditions in our framework, we analyze them from the perspective of low-Tucker-rank manifold geometry. Let $\mathcal{M}_{(r_1, r_2, r_3)}$ denote the manifold of tensors $\mcT \in \mathbb{R}^{p_1 \times p_2 \times p_3}$ with Tucker rank $(r_1, r_2, r_3)$. At a point $\mcT = \mathcal{G} \times_1 U_1 \times_2 U_2 \times_3 U_3$, the tangent space is parameterized as
$$
\mathbb{T}_{\mcT} \mathcal{M}_{(r_1, r_2, r_3)} = \Big\{\widetilde{\mathcal{G}} \times_1 U_1 \times_2 U_2 \times_3 U_3 + \sum_{j=1}^3 \mathcal{G} \times_j \widetilde{U}_j \times_{j+1} U_{j+1} \times_{j+2} U_{j+2} \,\Big|\, \widetilde{U}_j^{\top} U_j = 0 \Big\},
$$
where $\widetilde{\mathcal{G}} \in \mathbb{R}^{r_1 \times r_2 \times r_3}$ and $\widetilde{U}_j \in \mathbb{R}^{p_j \times r_j}$ are free parameters that represent perturbations to the core tensor and factor subspaces, respectively.
The orthogonal projection of a tensor $\mathcal{A} \in \mathbb{R}^{p_1 \times p_2 \times p_3}$ onto the tangent space $\mathbb{T}_{\mcT} \mathcal{M}_{(r_1, r_2, r_3)}$ is defined as
$$
\begin{aligned}
& \mcP_{\mathbb{T}_{\mcT} \mathcal{M}_{(r_1, r_2, r_3)}}: \mathbb{R}^{p_1 \times p_2 \times p_3} \rightarrow \mathbb{T}_{\mcT} \mathcal{M}_{(r_1, r_2, r_3)}, \\
& \mathcal{A} \mapsto \mathcal{A} \times_1 \mcP_{U_1} \times_2 \mcP_{U_2} \times_3 \mcP_{U_3} + \sum_{j=1}^3 \operatorname{Mat}_j^{-1}\big(\mcP_{U_{j\perp}} A_j \mcP_{(U_{j+2} \otimes U_{j+1}) G_j^{\top}} \big),
\end{aligned}
$$
where $\mcP_{U_j}$ is the projection matrix onto the subspace spanned by $U_j$, $\mcP_{U_{j\perp}} = I - \mcP_{U_j}$ projects onto the orthogonal complement of $U_j$, and $\operatorname{Mat}_j^{-1}$ is the inverse of the mode-$j$ matricization operator.

The alignment conditions are commonly assumed in matrix inference \citep{xia2021statistical, chen2019inference} and tensor inference \citep{agterberg2024statistical, ma2024statistical}. These alignment conditions (\eqref{eq: general alignment condition in tensor regression without sample splitting}, \eqref{eq: statistically optimal alignment condition in tensor regression with sample splitting}, and \eqref{eq: computationally optimal alignment condition in tensor regression with sample splitting}) ensure that the loading tensor $\mathcal{A}$ is sufficiently aligned with the tangent space of the low-Tucker-rank manifold at $\mcT$. Specifically, the alignment conditions guarantee that the magnitude of the asymptotic normal term $\langle \mcP_{\mathbb{T}_{\mcT} \mathcal{M}_{(r_1, r_2, r_3)}} (\mathcal{A}), \mcP_{\mathbb{T}_{\mcT} \mathcal{M}_{(r_1, r_2, r_3)}} (\widehat{\mathcal{Z}}^{(1)}) \rangle$ dominates the perturbation terms in the normal space. 

The tangent space $\mathbb{T}_{\mcT} \mathcal{M}{(r_1, r_2, r_3)}$ can be decomposed into a direct sum of the following subspaces:
$$
\mathbb{T}_{\mcT} \mathcal{M}_{(r_1, r_2, r_3)} = \mathbb{T}_{\mathcal{G}} \oplus \mathbb{T}_1 \oplus \mathbb{T}_2 \oplus \mathbb{T}_3,
$$
where
$$
\mathbb{T}_{\mathcal{G}}=\big\{\widetilde{\mathcal{G}} \times_1 U_1 \times_2 U_2 \times_3 U_3 \mid \widetilde{\mathcal{G}} \in \mathbb{R}^{r_1 \times r_2 \times r_3}\big\},
$$
captures perturbations to the core tensor, and
$$
\mathbb{T}_j=\big\{\mathcal{G} \times_j \widetilde{U}_j \times_{j+1} U_{j+1} \times_{j+2} U_{j+2} \mid \widetilde{U}_j^{\top} U_j=0, \widetilde{U}_j \in \mathbb{R}^{p_j \times r_j} \big\},
$$
for $j=1,2,3$, captures perturbations to the factor matrices in each mode. This decomposition facilitates the analysis of how the loading tensor $\mathcal{A}$ interacts with different components of the tangent space. A similar characterization of perturbation sources on the tangent space is also employed in the proof of the minimax lower bound in Theorem~\ref{thm: minimax lower bound in tensor regression}.

In addition to perturbations in the normal space, the error terms $\Omega_1$ and $\Omega_3$ in the theorems also include perturbations in the direction of $\mathcal{A} \times_1 \mcP_{U_1} \times_2 \mcP_{U_2} \times_3 \mcP_{U_3}$, which lie within the tangent space—specifically in the subspace $\mathbb{T}_{\mathcal{G}}$. These perturbations are not controlled by the alignment condition and constitute leading terms in higher-order perturbations. To address this, the incoherence condition (\eqref{eq: computationally optimal incoherence condition in tensor regression with sample splitting} and \eqref{eq: computationally optimal incoherence condition in tensor PCA}) requires that the projection of $\mathcal{A}$ onto $\mathbb{T}_{\mathcal{G}}$ is relatively small compared to its projection onto the entire tangent space. For each mode $j = 1, 2, 3$, the condition imposes an upper bound on 
$$
\frac{\|\mcP_{U_j}A_j\mcP_{(U_{j+2}\otimes U_{j+1})G_j^{\top}} \|_{\mathrm{F}}}{ \|\mcP_{U_{j\perp}}A_j\mcP_{(U_{j+2}\otimes U_{j+1})G_j^{\top}} \|_{\mathrm{F}}} \quad \text{or} \quad \frac{\|\mcP_{U_j}A_j(\mcP_{U_{j+2}} \otimes \mcP_{U_{j+1}}) \|_{\mathrm{F}}}{  \|\mcP_{U_{j\perp}}A_j(\mcP_{U_{j+2}} \otimes \mcP_{U_{j+1}})\|_{\mathrm{F}}},
$$
where the condition limits the influence of $\mathcal{A}$ in the direction of core tensor perturbations. 
This incoherence condition is satisfied when $\|U_j^{\top} V_j \|$, representing the angle between the mode-$j$ subspace of $\mcT = \mathcal{G} \times_1 U_1 \times_2 U_2 \times_3 U_3$ and $\mathcal{A} = \mathcal{B} \times_1 V_1 \times_2 V_2 \times_3 V_3$, is sufficiently small for $j=1,2,3$. This implies that the loading tensor $\mathcal{A}$ is not fully aligned with the subspace $\mathbb{T}_{\mathcal{G}}$ of the tangent space, which represents perturbations in the core tensor. Since the core tensor subspace is a component of the tangent space, the incoherence condition and the alignment condition are, in some sense, inherently at odds with one another.

To illustrate this more clearly, consider the entrywise inference as an example.
The combination of the incoherence and alignment conditions assumed in prior works \citep{chen2019inference, ma2024statistical, xia2021statistical} requires $\max_j\|U_j\|_{2,\infty} \asymp \sqrt{\oor/\op}$, effectively enforcing that $\|U_j \|_{2,\infty}$ reaches its lower bound. This imposes highly restrictive constraints on the factor matrices, limiting these methods to scenarios where the factor matrices exhibit a high level of incoherence. 

In contrast, our framework relaxes these stringent requirements, enhancing the flexibility and applicability of inferential procedures in low-rank tensor settings. Specifically, by not requiring the incoherence condition, our normal approximation framework for general linear functionals, presented in Section \ref{sec: Asymptotic Normality of Estimated Linear Functional without Sample Splitting}, offers greater flexibility while achieving nearly computationally optimal sample size requirements. Furthermore, when combined with sample splitting, our framework attains computationally optimal sample sizes for low-rank linear functionals under a weaker incoherence condition \eqref{eq: computationally optimal incoherence condition in tensor regression with sample splitting}, where $\max_j\|U_j\|_{2,\infty} \lesssim \oor^{1/2} \op^{-1/4}$, compared to those in the existing literature. 

Additionally, our results show that increasing the signal strength $\ulambda$ of the underlying parameter tensor allows for even weaker incoherence and alignment conditions. Consequently, our tensor regression framework, along with the inferential procedures for tensor PCA introduced subsequently, extends the applicability and flexibility of statistical inference in low-rank tensor models.

\subsection{Data-driven Inference of Estimated Linear functionals} \label{sec: Data-driven Inference of Estimated Linear functionals in tensor regression}

The asymptotic normality of the estimator $\langle \widehat{\mcT}, \mcA \rangle$, established in the previous section, provides a foundation for statistical inferences about the linear functional $\langle \mcT, \mcA \rangle$. To construct confidence intervals or perform hypothesis testing in practical applications, it is crucial to accurately estimate the variance of $\langle \widehat{\mcT}, \mcA \rangle$.

To estimate the noise variance $\sigma_\xi^2$, we define the following estimators: 
Without sample splitting,
\begin{align}
\widehat{\sigma}_{\xi}^2
= &\frac{1}{n} \sum_{i=1}^{n} \big(Y_i- \langle\widehat{\mcT}^{\text{init}}, \mcX_i \rangle \big)^2 \label{eq: estimate of sigmaxi2 in tensor regression without sample splitting} 
\end{align}
as specified in Theorem~\ref{thm: main theorem in tensor regression without sample splitting}.
With sample splitting,
\begin{align}
\widehat{\sigma}_{\xi}^2
= &\frac{1}{n} \sum_{i_1=1}^{n_1} \big(Y_i^{\rm (\RN{1})}- \langle\widehat{\mcT}^{\text{init, (\RN{2})}}, \mcX_i^{\rm (\RN{1})} \rangle \big)^2+\frac{1}{n} \sum_{i_2=1}^{n_2}\big(Y_i^{(\rm \RN{2})}- \langle\widehat{\mcT}^{\text{init, (\RN{1})}}, \mcX_i^{\rm (\RN{2})}\rangle \big)^2 \label{eq: estimate of sigmaxi2 in tensor regression with sample splitting}
\end{align}
as specified in Theorem~\ref{thm: main theorem in tensor regression with sample splitting}.
For the design variance $\sigma^2$, we use
\begin{align}
\hat\sigma^2 = & \frac{1}{np_1p_2p_3}\sum_{i=1}^n \|\mcX_i \|_{\mathrm{F}}^2. \label{eq: estimate of sigma2 in tensor regression}
\end{align} 
The low-rank property of the initial estimate $\widehat{\cT}$, imposed in Assumption~\ref{assump: initial signal tensor estimate error bound in tensor regression}, ensures the consistency of these variance estimators.

To estimate the variance components $s_{\mcA}^2$ in \eqref{eq: variance component sA2}, associated with the linear functional, we define
{\small \begin{align}
& \widehat{s}_{\mcA}^2 = \sum_{j=1}^3\big\|\big(\mathcal{I}-\mcP_{\whU_j}\big)A_j\big(\whU_{j+2}\otimes \whU_{j+1}\big)\widehat{W}_j\widehat{W}_j^{\top}\big(\whU_{j+2}\otimes \whU_{j+1}\big)^{\top}\big\|_{\mathrm{F}}^2 + \big\|\mcA \times_1 \whU_1 \times_2 \whU_2 \times_3 \whU_3 \big\|_{\mathrm{F}}^2, \label{eq: estimate of sA2}
\end{align}}%
where
\begin{align}
& \widehat{W}_j=\mathrm{QR}\big[\Mat_j\big(\widehat{\mcT} \times_1 \whU_1^{\top} \times_2 \whU_2^{\top} \times_3 \whU_3^{\top} \big)^{\top}\big] \label{eq: estimate of right singular space of T}
\end{align}
is the estimate of the right singular space of the mode-$j$ matricization of the core tensor $\mathcal{G} \in \mathbb{R}^{r_1 \times r_2 \times r_3}$, obtained via QR decomposition for each mode $j=1,2,3$. In the sample splitting case, $\whU_j$ can be either $\whU_j^{(\text{\RN{1}})}$ or  $\whU_j^{(\text{\RN{2}})}$.

The following theorem establishes that the asymptotic normality of $\langle \widehat{\mcT}, \mcA\rangle$ remains valid when the variance is replaced by plug-in estimates.

\begin{theorem}\label{thm: main theorem in tensor regression with plug-in estimates}

Under the conditions in Corollary~\ref{corollary: asymptotic normality of estimated general linear functional in tensor regression without sample splitting}, let the variance components ${\sigma}_{\xi}^2$, ${\sigma}^2$, and ${s}_{\mathcal{A}}^2$ be estimated by $\widehat{\sigma}_{\xi}^2$, $\widehat{\sigma}^2$, and $\widehat{s}_{\mathcal{A}}^2$, as defined in \eqref{eq: estimate of sigmaxi2 in tensor regression without sample splitting}, \eqref{eq: estimate of sigma2 in tensor regression}, and \eqref{eq: estimate of sA2}, respectively. Then, 
\begin{align*}
& \frac{\sqrt{n} \big( \langle\widehat{\mcT}, \mcA \rangle-\langle \mcT, \mcA\rangle \big)}{(\widehat{\sigma}_{\xi}/\widehat{\sigma}) \widehat{s}_{\mcA} } \rightarrow \mathcal{N}(0,1).
\end{align*}

In the sample-splitting case, as specified in Corollary~\ref{corollary: asymptotic normality of estimated low-Tucker-rank linear functional in tensor regression with sample splitting}, where $\widehat{\sigma}_\xi^2$ is defined in \eqref{eq: estimate of sigmaxi2 in tensor regression with sample splitting}, the same asymptotic normality result holds.
\end{theorem}

A generalized version of this theorem, including non-asymptotic results, is provided in the appendix. In particular, Theorem~\ref{thm: main theorem in tensor regression with plug-in estimates} does not impose additional assumptions beyond those stated in Theorem~\ref{thm: main theorem in tensor regression without sample splitting} and Theorem~\ref{thm: main theorem in tensor regression with sample splitting}. Under these conditions, the variance estimators ensure that the asymptotic normality of $\langle \widehat{\mcT}, \mcA \rangle$ holds, enabling the construction of confidence intervals for $\langle \mcT, \mcA \rangle$ using the plug-in variance estimates. Specifically, the $100(1 - \alpha)\%$ confidence interval is given by
$$
\widehat{\mathrm{CI}}_{\mcA, \cT}^{\alpha}=\big[\langle  \widehat{\mcT}, \mcA \rangle-z_{\alpha / 2} (\widehat{\sigma}_{\xi}/\widehat{\sigma})  \widehat{s}_{\mcA} /\sqrt{n}, \langle \widehat{\mcT}, \mcA \rangle + z_{\alpha / 2} (\widehat{\sigma}_{\xi}/\widehat{\sigma}) \widehat{s}_{\mcA} /\sqrt{n}\big]  ,
$$
where $\alpha \in (0, 1)$, and $z_\theta = \Phi^{-1}(1 - \theta)$ denotes the upper $\theta$ quantile of the standard normal distribution.

\subsection{Minimax Optimality of the Confidence Interval Length} \label{sec: Minimax Optimality of the Confidence Interval Length in tensor regression}


A natural question is how the proposed inferential procedures compare to other methods. To evaluate their performance, we analyze the length of the confidence intervals constructed in our framework and compare them with the minimax lower bound for the tensor regression setting. In this section, we establish that these confidence intervals achieve minimax rate optimality, demonstrating the efficiency of our approach.

Our analysis builds upon results from related literature \citep{cai2015confidence} but diverges from the Cram\'er--Rao lower bound traditionally studied in information geometry \citep{smith2005covariance,ma2024statistical}. Unlike the Cram\'er--Rao framework, which focuses on parameter estimation under unbiasedness constraints, our work characterizes the optimal performance of confidence intervals under general perturbations in the tangent space of the low-Tucker-rank manifold. This perspective allows us to rigorously show that the proposed procedures attain the fundamental limits of inference accuracy in tensor regression.

In what follows, the parameter space $\Theta(\underline{\lambda}, \kappa)$ is defined as
{\small \begin{align}
\Theta(\ulambda, \kappa )
:= \big\{ \mcT = \mcG \times_1 U_1 \times_2 U_2 \times_3 U_3 \; \big| \; \mcG \in \mathbb{R}^{r_1 \times r_2 \times r_3}, \; U_j \in \mathbb{O}^{p_j \times r_j}, \ulambda \leq \lambda \leq \kappa \ulambda \big\} , \label{eq: parameter sapce}
\end{align}}%
where $\lambda$ represents any nonzero singular value of $\operatorname{Mat}_j(\mcT )$ for each mode $j=1,2,3$.
We also define the class of confidence intervals with coverage level $100(1-\alpha)\%$ as
\begin{align}
\mathcal{I}_\alpha(\Theta,\mcA ):=\big\{\mathrm{CI}_{\mcA }^\alpha(\mcT,\mathcal D ) =[l(\mathcal D), u(\mathcal D)]: \inf_{\mcT \in \Theta } \mathbb{P} \big(l(\mathcal D) \leq \langle \mcT, \mcA\rangle \leq u(\mathcal D) \big) \geq 1-\alpha\big\}, \label{eq: space of confidence intervals}
\end{align}
where the observed data is denoted as $\mathcal D:=\{(Y_1,\cX_1),...,(Y_n,\cX_n) \}$.

The following theorem establishes a minimax lower bound for the expected length of any valid confidence interval in this tensor regression setting, assuming Gaussian design and noise.

\begin{theorem}\label{thm: minimax lower bound in tensor regression}
Suppose that the significance level $\alpha$ satisfies $0<\alpha<1 / 2$, and let $\overline{\lambda} = \kappa\ulambda > \ulambda \geq c$. Additionally, assume that the noise terms $\xi_i$'s are i.i.d $N(0,\sigma_\xi^2)$ and the design tensors $\{\cX_i\}$'s are i.i.d. random tensors with i.i.d. $N(0,\sigma^2)$ entries. Let $L(\cdot)$ denote the length of a confidence interval. Then, under the tensor regression model \eqref{eq: tensor regression model}, where the loading tensor $\mcA$ and signal tensor $\mcT$ satisfy Assumption~\ref{assump: structural assumption on the signal tensor and loading tensor}, there exists some constant $c_3>0$ such that
\begin{align*}
& \inf_{\mathrm{CI}_{\mcA }^\alpha(\mcT,\mathcal D) \in \mathcal{I}_\alpha(\Theta,\mcA)} \sup_{ \mcT \in \Theta(\ulambda, \kappa)} \mathbb{E} L\big(\mathrm{CI}_{\mcA }^\alpha ( \mcT, \mathcal D)\big) \\
\geq & \frac{c_3 \sigma_\xi}{\sigma\sqrt{n}} \sqrt{\sum_{j=1}^3\big\|\mcP_{U_{j\perp}}A_j\mcP_{(U_{j+2}\otimes U_{j+1})G_j^{\top}}\big\|_{\mathrm{F}}^2 +  \big\|\mcA \times_1 U_1 \times_2 U_2 \times_3 U_3\big\|_{\mathrm{F}}^2}.
\end{align*}
\end{theorem}

Theorem~\ref{thm: minimax lower bound in tensor regression} demonstrates that the confidence intervals derived from the asymptotic normality results in Sections \ref{sec: Asymptotic Normality of Estimated Linear Functional without Sample Splitting} and \ref{sec: Asymptotic Normality of Estimated Linear Functional with Sample Splitting}
achieve minimax rate optimality. The lower bound on the expected confidence interval length under the tensor regression setting is novel, 
and grounded in the characterization of perturbation along possible directions in the tangent space $\mathbb{T}_{\mcT}\mathcal{M}_{(r_1, r_2, r_3)}$ of the low-Tucker-rank manifold $\mathcal{M}_{(r_1, r_2, r_3)}$ at 
$\mcT$. This aligns with the discussion in Section \ref{sec: The Interplay between the Incoherence Condition and Alignment Condition}. 

\section{Inference for Tensor PCA} \label{sec: tensor_pca}

\subsection{Problem Setting}

In this section, we focus on the tensor Principal Component Analysis (PCA) model, defined as
\begin{align}
\mcY=\mcT+\mcZ , \label{eq: tensor PCA model}
\end{align}
where $\mcY \in \mathbb{R}^{p_1 \times p_2 \times p_3}$ is the observed tensor, $\mcT$ is the underlying signal tensor, and $\mathcal{Z}$ is a noise tensor. The signal tensor $\mcT$ admits a low Tucker-rank decomposition $\mcT=\mcG \times_1 U_1 \times_2 U_2 \times_3 U_3$, where $U_j \in \mathbb{O}^{p_j \times r_j}$. The goal is to perform valid statistical inference on the linear functional $\langle \mcT, \mcA \rangle$, where $\mcA \in \mathbb{R}^{p_1 \times p_2 \times p_3}$ is a prespecified loading tensor. We also consider cases where $\mcA$ is a potentially low-rank tensor with Tucker rank $(R_1, R_2, R_3)$, expressed as $\mcA = \mcB \times_1 V_1 \times_2 V_2 \times_3 V_3$. Entrywise inference is a special case where the loading tensor has Tucker rank $(1,1,1)$. Structural assumptions on the signal tensor $\mcT$ and the loading tensor $\mcA$ are detailed in Assumption~\ref{assump: structural assumption on the signal tensor and loading tensor}.

In our analysis, we assume that the noise tensor $\mcZ$ has i.i.d. sub-Gaussian entries, formalized as follows.
\begin{assumption}[Sub-Gaussian Noise in tensor PCA]\label{assump: sigma-sub-Gaussian noise in tensor PCA}
The noise tensor $\mcZ$ is entrywise i.i.d. with mean zero and sub-Gaussian, and its $\psi_2$ Orlicz norm is bounded by $\sigma > 0$, i.e., $\|[\mcZ]_{j_1,j_2,j_3}\|_{\psi_2} \leq \sigma$.
\end{assumption}
This implies that the variance of each entry satisfies $\operatorname{Var} ( [\mcZ ]_{j_1,j_2,j_3} ) = \mathbb{E}([\mcZ ]_{j_1,j_2,j_3}^2) \leq C\sigma^2$ for some constant $C>0$.

\subsection{Estimation of Linear Functionals}\label{sec:tpca_alg}

To estimate $\langle \mcT, \mcA \rangle$ under the tensor PCA setting, we propose the following algorithm.

\subsubsection*{Step 1: Initialization.} \label{step:step1 initialization in main algorithm}

We use the observed tensor $\mcY$ as initial estimate for the signal tensor $\mcT$, and obtain initial estimates of the loading factors $\whU_1^{\text{init}}$, $\whU_2^{\text{init}}$, and $\whU_3^{\text{init}}$ via Higher-Order SVD (\texttt{HOSVD}, \citet{de2000multilinear}). For shorthand, we denote $\whU_j^{(0)}:=\whU_j^{\text{init}}$ for $j=1,2,3$. 
Since $\mcY$ is already an unbiased estimator of $\mcT$, {\bf no debiasing step} is required in tensor PCA.

\begin{assumption}[Error Bound for Initial Estimates of Singular Spaces]\label{assump: initialization error bound of singular spaces in tensor PCA}
We assume that the initial singular space estimates 
satisfy the minimax-optimal error bound $ \|\widehat{U}_j^{(0)}\widehat{U}_j^{(0)\top} - U_jU_j^{\top}\| \le c_0 \sigma\sqrt{\op}/\ulambda $, for $j=1,2,3$, with probability at least $1 - \mathbb{P}(\mcE_{U}^{\text{PCA}})$, where the event $\mcE_{U}^{\text{PCA}}$ is defined as
$$
\mcE_{U}^{\text{PCA}} = \cup_{j=1,2,3} \{\|\whU_j^{(0)}\whU_j^{(0)\top} - U_jU_j^{\top} \| > c_0 \sigma\sqrt{\op}/\ulambda \}.
$$
\end{assumption}

\begin{remark}
The Higher-Order Orthogonal Iteration (\texttt{HOOI}\citet{ de2000best}) method achieves minimax-optimal error bound under i.i.d. sub-Gaussian noise \citep{zhang2018tensor}.
Thus, we assume minimax-optimal initialization without loss of generality.
\end{remark}

\subsubsection*{Step 2: Two-step Power Iteration.} \label{step3: Two-step power iteration in main algorithm}
Using the initial estimates $\whU_j^{(0)}$, we refine the singular space estimates via two-step power iteration. For each iteration $k = 1, 2$ and mode $j = 1, 2, 3$, the power iteration is performed as follows:

For each $j = 1, 2, 3$, $\whU_j^{(k)}$ is obtained as the leading $r_j$ left singular vectors of
$$
\Mat_j\big( \mcY \times_{j+1} \whU_{j+1}^{(k-1) \top} \times_{j+2} \whU_{j+2}^{(k-1) \top} \big) = \Mat_j( \mcY ) \big( \whU_{j+2}^{(k-1)} \otimes \whU_{j+1}^{(k-1)} \big).
$$

After completing the two iterations, the final estimates are $\whU_j := \whU_j^{(2)}$. 


\subsubsection*{Step 3: Projection and Plug-in Estimator.} 
The final signal tensor estimate is obtained by projecting $\mcY$ onto the estimated singular spaces:
\begin{align*}
& \widehat{\mcT} = \mcY \times_1 \mcP_{\whU_{1}} \times_2 \mcP_{\whU_{2}} \times_3 \mcP_{\whU_{3}},
\end{align*}
where $\mcP_{\whU_j} = \whU_j \whU_j^\top$ are projection matrices. The linear functional $\langle \mcT, \mcA \rangle$ is then estimated by $ \langle \widehat{\mcT}, \mcA \rangle$.

\subsection{Asymptotic Normality of Estimated Linear Functionals} \label{sec: Asymptotic Normality of Estimated Linear Functionals}

In this section, we establish the asymptotic normality of the estimator $\langle \widehat{\mcT} , \mcA \rangle$ obtained from Section \ref{sec:tpca_alg}. 

\begin{theorem}[Main Theorem: asymptotic normality in Tensor PCA]\label{thm: main theorem in tensor PCA}
Consider the tensor PCA model \eqref{eq: tensor PCA model}. Suppose that Assumptions~\ref{assump: structural assumption on the signal tensor and loading tensor}, \ref{assump: sigma-sub-Gaussian noise in tensor PCA}, and \ref{assump: initialization error bound of singular spaces in tensor PCA} hold, and $\ulambda \geq C\kappa \op^\frac{1}{2}$, where $C$ is a positive constant depending only on the noise scale $\sigma$. Then
{\small\begin{align*}
& \sup _{x \in \mathbb{R}} \left\lvert\, \mathbb{P}\left(\frac{\langle \widehat{\mcT}, \mcA\rangle-\langle \cT, \mcA\rangle}{\sigma \cdot s_{\mcA}} \leq x\right)-\Phi(x)\right| \lesssim \underbrace{\Psi}_{\substack{\text{rate of asymptotic}\\ \text{normal terms}}} + \underbrace{\frac{\Omega_1 + \Omega_2 + \Omega_3}{\sigma \cdot s_{\mcA}}}_{\text{rate of negligible terms}} + \underbrace{\big(\op^{-c} + \mcP_{\mcE_{U}^{\text{PCA}}} \big)}_{\text{rate of initial estimates}} ,
\end{align*}}%
where $c$ is a positive constant, and the variance component $s_{\mcA}$ is defined in \eqref{eq: variance component sA2}.
Here,
\begin{align}
\Psi = \frac{K_3 \big(\sum_{j=1}^3\big\|\mcP_{U_{j\perp}} A_j\mcP_{(U_{j+2} \otimes U_{j+1})G_j^{\top}}\big\|_{\ell_\infty}+\|\mcA \times_1 U_1 \times_2 U_2 \times_3 U_3\|_{\ell_\infty}\big)}{\sigma^3 \big(\sum_{j=1}^3\big\|\mcP_{U_{j\perp}} A_j\mcP_{(U_{j+2} \otimes U_{j+1})G_j^{\top}}\big\|_{\mathrm{F}}^2+\|\mcA \times_1 U_1 \times_2 U_2 \times_3 U_3\|_{\mathrm{F}}^2\big)^{\frac{1} {2}}} \label{eq: rate of Berry-Essen bound in tensor PCA}
\end{align}
represents the Berry-Esseen bound for asymptotic normality, with $K_3 = \mathbb{E}(|\mcZ_{i,j,k}|^3) \lesssim \sigma^3$,
{\small \begin{align*}
\Omega_1
= & \|\mcA \times_1 U_1 \times_2 U_2 \times_3 U_3\|_{\mathrm{F}} \left(\frac{\sigma^2\op^{1/2} \oor^{1/2} \log (\op)^{1/2}}{\ulambda}+\frac{\sigma^3 \op^{3/2}\oor^{1/2} }{\ulambda^2}\right), \\
\Omega_2
= & \sum_{j=1}^3 \big\|\mcP_{U_{j\perp}} A_j \mcP_{(U_{j+2} \otimes U_{j+1}) G_j^{\top}}\big\|_{\mathrm{F}} \left(\frac{\sigma^2 \op^{1/2} \oor \log(\op)^{1/2}}{\ulambda}+\frac{\sigma^3 \op^{3/2}\oor^{1/2} }{\ulambda^2}\right) \\
+ & \sum_{j=1}^3 \|\mcA \times_j U_j \|_{\mathrm{F}} \left(\frac{\sigma^2 \oor^{3/2} \log (\op)}{\ulambda}+\frac{\sigma^3 \op^{1/2} \oR \oor^{1/2} \log (\op)}{\ulambda^2}+\frac{\sigma^4 \op^{3 / 2}\oR^{1/2} \oor^{1/2} \log (\op)^{1/2}}{\ulambda^3}\right)  \\
+ & \|\mcA\|_{\mathrm{F}} 
\left(\frac{\sigma^3 \oor^2 \log (\op)^{3 / 2}}{\ulambda^2} 
+ \frac{\sigma^4  \op^{1 / 2} \oR^{3 / 2}  \oor^{1/2}\log (\op)^{3 / 2}}{\ulambda^3} 
+ \frac{\sigma^5 \op^{3 / 2} \oR \oor^{1/2}   \log (\op)}{\ulambda^4}\right), \\
\Omega_3
= & \sum_{j=1}^3 \big\|\mcP_{U_j} A_j \mcP_{(U_{j+2} \otimes U_{j+1}) G_j^{\top}}\big\|_{\mathrm{F}}  \frac{\sigma^2 \op \oor^{1/2} }{\ulambda} + \sum_{j=1}^{3} \big\|\mcA \times_{j+1} U_{j+1} \times_{j+2} U_{j+2}\big\|_{\mathrm{F}}  \frac{\sigma^3 \op \oor \log(\op)^{1/2}}{\ulambda^2} 
\end{align*}}%
are the upper bounds of negligible terms involved due to the noise $\mcZ$ and the projection error. 
\end{theorem}

Similar to Theorem~\ref{thm: main theorem in tensor regression without sample splitting}, the upper bounds $\Omega_1, \Omega_2$, and $\Omega_3$, capture different sources of error. Specifically, $\Omega_1$ arises from the noise tensor, $\Omega_2$ originates from negligible terms when using first-order perturbation for normal approximation, $\Omega_3$ reflects shared contributions between both steps. Theorem~\ref{thm: main theorem in tensor PCA} demonstrates that the convergence rate of the estimated liner form $\langle\widehat{\mcT}, \mcA \rangle$ depends on the ratio of the $\ell_\infty$ norm and the $\ell_2$ norm of the variance components.

Building on Theorem~\ref{thm: main theorem in tensor PCA}, we explore two types of inference problems: inference for low-rank linear functionals and general linear functionals.

\begin{corollary}[Asymptotic normality of estimated low-Tucker-rank linear functionals] \label{corollary: asymptotic normality of estimated low-Tucker-rank linear functional in tensor PCA}

Under the conditions of Theorem~\ref{thm: main theorem in tensor PCA}, assume that the Tucker rank of the loading tensor $\mcA$ satisfies $\operatorname{rank}(\mcA)=(R_1, R_2, R_3)$ is fixed and independent of $\op$.
If $\ulambda \geq C \max\{\kappa\op^{1/2}, \op^{3/4}\oor^{1/2}\log(\op)\}$, and the variance component $s_\mcA$ defined in \eqref{eq: variance component sA2} satisfies the {\it alignment condition}
{\small \begin{align}
s_\mcA \geq C_2\max_{j=1,2,3}\big\{\op^{-1/2}\big\|\mcA\times_{j+1} U_{j+1} \times_{j+2} U_{j+2}\big\|_{\mathrm{F}}, \oor\op^{-3/4}\big\|\mcA\times_j U_j\big\|_{\mathrm{F}}, \oor\op^{-3/2}\|\mcA\|_{\mathrm{F}}\big\}  , \label{eq: computationally optimal alignment condition in tensor PCA}
\end{align}}%
and the {\it incoherence condition} \eqref{eq: computationally optimal incoherence condition in tensor PCA}
\begin{align*}
\max_{j=1,2,3}\big\|\mcP_{U_j}A_j\mcP_{(U_{j+2}\otimes U_{j+1})G_j^{\top}}\big\|_{\mathrm{F}} / \big\|\mcP_{U_{j\perp}}A_j\mcP_{(U_{j+2}\otimes U_{j+1})G_j^{\top}}\big\|_{\mathrm{F}} \leq c_2 \oor^{1/2}\op^{-1/4}, 
\end{align*}
where $C$, $C_2$ and $c_2$ are positive constants depending only on $\oR$, the fixed rank of the loading tensor $\mcA$, and the noise scale $\sigma$, then it holds that
$
(\langle \widehat{\mcT}, \mcA \rangle-\langle \mcT, \mcA\rangle ) / (\sigma \cdot s_\mcA) \overset{d}{\rightarrow} \mathcal{N}(0,1).
$

\end{corollary}

In the special case of entrywise inference under the tensor PCA setting, our approach imposes significantly weaker incoherence and alignment conditions compared to \citet{agterberg2024statistical}, which assumes $\|U_j \|_{2, \infty} \asymp \sqrt{\oor/\op}$ and $s_\mcA \gtrsim \oor^2\op^{-5/4} \|\mcA \|_{\mathrm{F}}$ for loading tensor $\cA$ of the form $e_i\otimes e_k \otimes e_l$. Our incoherence and alignment conditions align with those in Theorem 7 of \citet{xia2022inference}, which focused on entrywise inference for rank-one signal tensors in tensor PCA. However, our results generalize this to tensors of any low Tucker rank. Additionally, unlike existing literature \citep{agterberg2024statistical, ma2024statistical, xia2021statistical}, our framework remains computationally optimal even when the condition number $\kappa$ diverges at a rate of $\mcO(\op^{1/4})$. Beyond entrywise inference, our method naturally extends to other low-Tucker-rank linear functionals, offering broad applicability.

In the tensor completion setting, \citet{ma2024statistical} demonstrated that statistically optimal sample sizes and SNRs are sufficient for accurate inference using $\ell_{2,\infty}$ perturbation analysis. The statistically optimal SNR is achievable because the incoherence condition significantly simplifies uniform sampling. However, the $\ell_{2,\infty}$ perturbation analysis in \citet{ma2024statistical} cannot be readily extended to sub-Gaussian settings, especially when the loading tensor $\mcA$ is not sparse, as required by the condition $\|\mcA\|_{\ell_1} / \|\mcA\|_{\mathrm{F}}$ being bounded in \citet{ma2024statistical}. Furthermore, a consistent initialization via \texttt{HOSVD} \citet{de2000multilinear} requires computationally optimal size (Remark 1, \citet{zhang2018tensor}). Consequently, it remains unclear whether the statistically optimal SNR is attainable for inference under tensor PCA settings. 

In addition to inference for low-Tucker-rank linear functionals, Theorem~\ref{thm: main theorem in tensor PCA} implies the following asymptotic normality for estimated general linear functionals without requiring incoherence conditions. 

\begin{corollary}[Asymptotic normality of estimated general linear functionals] \label{corollary: asymptotic normality of estimated general linear functional in tensor PCA}

Suppose that the loading tensor $\mcA$ has rank at most $(p_1, p_2, p_3)$. When $\ulambda \geq C\max \{\kappa\op^{1/2}, \op\oor^{1/2} \}$, and $s_{\mcA}$ defined in \eqref{eq: variance component sA2} satisfies the {\it alignment condition}
\begin{align}
s_\mcA \geq C_2\max_{j=1,2,3}\big\{\oor\op^{-1/2}\|A\times_j U_j\|_{\mathrm{F}}, \oor\op^{-1}\|\mcA\|_{\mathrm{F}} \big\}, \label{eq: general alignment condition in tensor PCA}
\end{align}
where $C$ and $C_2$ are positive constants depending only on the noise scale $\sigma$, then it holds that
$
(\langle \widehat{\mcT}, \mcA \rangle-\langle \mcT, \mcA\rangle )/ (\sigma \cdot s_\mcA ) \overset{d}{\rightarrow} \mathcal{N}(0,1) .
$

\end{corollary}

\subsection{Data-driven Inference of Estimated Linear Functionals} \label{sec: Data-driven Inference of Estimated Linear Functionals in tensor PCA}

Building on the asymptotic normality established in Theorem~\ref{thm: main theorem in tensor PCA}, we extend the methodology to enable data-driven inference for $ \langle \mcT, \mcA  \rangle$ under the tensor PCA setting, incorporating plug-in variance estimates. To estimate the noise variance $\sigma^2$, we define 
\begin{align}
& \widehat{\sigma}^2=\big\|\mcY-\mcY \times_1 \mcP_{\whU_1} \times_2 \mcP_{\whU_2} \times_3 \mcP_{\whU_3} \big\|_{\mathrm{F}}^2/(p_1 p_2 p_3), \label{eq: estimate of sigma2 in tensor PCA}
\end{align}
and use \eqref{eq: estimate of sA2} to estimate the variance component $s_{\mcA}^2$, where $\widehat{W}_j=\mathrm{QR} [\operatorname{Mat}_j (\widehat{\mcT} \times_1 \whU_1^{\top} \times_2 \whU_2^{\top} \times_3  \whU_3^{\top} )^{\top} ]$,
for $j=1,2,3$, estimating the right singular space of the mode-$j$ matricization of the core tensor $\mcG \in \mathbb{R}^{r_1 \times r_2 \times r_3}$. 
The noise variance estimator $\widehat{\sigma}^2$ in \eqref{eq: estimate of sigma2 in tensor PCA} follows the construction in \citet{xia2022inference} (see their Lemma 1). Its accuracy relies on the assumption that the noise tensor $\mcZ$ has i.i.d. entries.

We further extend the asymptotic normality result to demonstrate that these variance estimates are valid for practical statistical inference. Specifically, the following theorem establishes the validity of using the plug-in variance estimates $\widehat{\sigma}^2$ and $\widehat{s}_{\mcA}^2$.

\begin{theorem} \label{thm: main theorem in tensor PCA with plug-in estimates}

Under the conditions in Corollary~\ref{corollary: asymptotic normality of estimated low-Tucker-rank linear functional in tensor PCA} or Corollary~\ref{corollary: asymptotic normality of estimated general linear functional in tensor PCA}, let the variance components ${\sigma}^2$, and ${s}_{\mathcal{A}}^2$ be estimated by $\widehat{\sigma}^2$ and $\widehat{s}_{\mcA}^2$, as defined as in \eqref{eq: estimate of sigma2 in tensor PCA} and \eqref{eq: estimate of sA2}, respectively. Then
$ ( \langle\widehat{\mcT}, \mcA \rangle-\langle \mcT, \mcA\rangle ) / (\widehat{\sigma} \cdot \widehat{s}_{\mcA}  ) \rightarrow \mathcal{N}(0,1).$

\end{theorem}

A generalized version of this theorem, including non-asymptotic results, is presented in the appendix. 
This theorem indicates that, under the same conditions as Corollary~\ref{corollary: asymptotic normality of estimated low-Tucker-rank linear functional in tensor PCA} or Corollary~\ref{corollary: asymptotic normality of estimated general linear functional in tensor PCA}, the plug-in variance estimates $\widehat{\sigma}^2$ and $\widehat{s}_{\mcA}^2$ enable valid construction of confidence intervals and hypothesis tests for $\langle \mcT, \mcA  \rangle$. Specifically, a $100(1 - \alpha)\%$ confidence interval for $\langle \mcT, \mcA  \rangle$ is given by
$$
\widehat{\mathrm{CI}}_{\mcA,\mcT}^{\alpha}=\big[\langle \widehat{\mcT}, \mcA \rangle-z_{\alpha / 2}  \widehat{\sigma} \widehat{s}_{\mcA} ,\langle  \widehat{\mcT}, \mcA \rangle+z_{\alpha / 2} \cdot \widehat{\sigma} \widehat{s}_{\mcA} \big],
$$
where $z_{\alpha/2}$ is the upper $\alpha/2$ quantile of the standard normal distribution.

\subsection{Minimax Optimality of the Confidence Interval Length} \label{sec: Minimax Optimality of the Confidence Interval Length in tensor PCA}

We define the parameter space of the signal tensor for the tensor PCA problem and the set of confidence intervals as in \eqref{eq: parameter sapce} and \eqref{eq: space of confidence intervals}, respectively, consistent with the tensor regression problem. The observed data is denoted by $\mathcal{D}=\cY$. The following theorem establishes the minimax lower bound on the expected length of confidence intervals, for estimating linear functionals of the signal tensor $\mcT$ under the tensor PCA model.

\begin{theorem}\label{thm: minimax lower bound in tensor PCA}
Suppose that the significance level $\alpha$ satisfies $0<\alpha<1 / 2$, and let $\kappa >1$. Additionally, assume that the entries of $[\mcZ ]_{i,j,k}$'s are i.i.d. $N(0,\sigma^2)$. Then, under the tensor PCA model \eqref{eq: tensor PCA model}, there exists some constant $c_3>0$ such that
\begin{align*}
& \inf_{\mathrm{CI}_{\mcA }^\alpha(\mcT, \mathcal D) \in \mathcal{I}_\alpha(\Theta,\mcA)} \sup_{ \mcT \in \Theta(\ulambda, \kappa)}\mathbb{E} L\big(\mathrm{CI}_{\mcA}^\alpha(\mcT, \mathcal D )\big) \\
\geq & c_3 \sigma \sqrt{\sum_{j=1}^3\big\|\mcP_{U_{j \perp}} A_j \mcP_{(U_{j+2} \otimes U_{j+1} ) G_j^{\top}}\big\|_{\mathrm{F}}^2+\big\|\mcA \times_1 U_1 \times_2 U_2 \times_3 U_3\big\|_{\mathrm{F}}^2},
\end{align*}
where $L(\cdot)$ is the length of the confidence interval.
  
\end{theorem}

Theorem~\ref{thm: minimax lower bound in tensor PCA} demonstrates that the confidence intervals derived from our estimators, which achieve this lower bound, are minimax optimal. In the context of entrywise inference under the tensor PCA setting, \citet{agterberg2024statistical} provides a similar minimax lower bound with a differently defined parameter space. This result highlights the effectiveness of our inference procedure for constructing confidence intervals for linear functionals of low-rank tensors in high-dimensional settings. 

\section{Numerical Experiments} \label{sec: Numerical Experiments}

In this section, we conduct numerical simulations to validate the proposed central limit theorems for estimated linear functionals. We consider three settings: tensor regression without sample splitting (Theorem~\ref{thm: main theorem in tensor regression without sample splitting}), tensor regression with sample splitting (Theorem~\ref{thm: main theorem in tensor regression with sample splitting}), and tensor PCA (Theorem~\ref{thm: main theorem in tensor PCA}).

Our simulations are designed to assess the performance of our inference framework under various scenarios. We begin by generating a core tensor $\mcG \in \mathbb{R}^{3 \times 3 \times 3}$ with full Tucker rank. The diagonal entries $G_{j,j,j}$ are sampled uniformly from the interval $[\ulambda, \overline{\lambda}]$, where $\overline{\lambda} = \kappa\ulambda$ while the off-diagonal entries are set to be zero. Tensor operations are implemented in R using the \texttt{rTensor} package \citep{li2018rtensor}. To construct the singular subspaces, we consider both coherent and incoherent settings. Coherent singular subspaces $U_j$ are generated by performing SVD on matrices where the first contains only a large value equals $\sqrt{\op}$, while the remaining entries are sampled from a standard normal distribution. In contrast, incoherent singular subspaces are generated by applying SVD to random Gaussian matrices with i.i.d. standard normal entries. Using these subspaces, the signal tensor $\mcT$ is constructed as $\mcT = \mcG \times_1 U_1 \times_2 U_2 \times_3 U_3$. We set $p_1=p_2=p_3=\op$.

We consider three scenarios for the loading tensor $\mcA$. In the first scenario with a full-Tucker-rank loading tensor, $\mcA$ is generated with entries drawn independently from a standard normal distribution and then rescaled to have unit Frobenius norm. This ensures that $\mcA$ has full Tucker rank $(p_1, p_2, p_3)$ with probability 1. 
In the second case, we use a low-Tucker-rank loading tensor defined as
$$
\mcA = \frac{1}{\sqrt{\lfloor2p_1^{1/4}\rfloor\lfloor2p_2^{1/4}\rfloor\lfloor2p_3^{1/4}\rfloor} } \sum\nolimits_{j_1=1}^{\lfloor2p_1^{1/4}\rfloor} \sum\nolimits_{j_2=1}^{\lfloor2p_2^{1/4}\rfloor} \sum\nolimits_{j_3=1}^{\lfloor2p_3^{1/4}\rfloor} e_{j_1} \otimes e_{j_2} \otimes e_{j_3},
$$
which has a Tucker rank of $(1, 1, 1)$ and maintains unit Frobenius norm. This low-rank construction presents a significant challenge for existing frameworks, such as \citet{ma2024statistical}, because it leads to $\|\mcA\|_{\ell_1} / \|\mcA\|_{\mathrm{F}} \gtrsim \op^{1/4}$, thereby violating their specified condition. To ensure compliance with the incoherence conditions in \eqref{eq: computationally optimal incoherence condition in tensor regression with sample splitting} and \eqref{eq: computationally optimal incoherence condition in tensor PCA}, the singular subspaces are generated to be incoherent ($\|U_j\|_{2,\infty} \lesssim \oor^{1/2}\op^{-1/2}$) with high probability. In the third scenario, we perform an entrywise inference using these incoherent singular spaces to validate the proposed central limit theorem with a statistically optimal sample size in the context of tensor regression.

\subsection{Simulation under the Tensor Regression Setting}

In the tensor regression setting, the observation noise $\{\xi_i\}_{i=1}^{n}$ and design tensors $\{\mcX_i\}_{i=1}^{n}$ are generated with i.i.d. standard normal entries, such that $\sigma_\xi = \sigma = 1$. 
Additionally, we set $p_1=p_2=p_3=40:=\op$.

The signal strength and condition number of the signal tensor $\mcT$ are controlled by setting $\underline{\lambda} = 1$ and $\overline{\lambda} = \kappa\underline{\lambda}$. For the inference of general linear functionals and low-Tucker-rank linear functionals, we set $\kappa=\sqrt{\op}$ and $\kappa=\op^{1/4}$ respectively, reflecting a more relaxed condition number compared to previous work \citep{xia2021statistical}. For entrywise inference, we set $\kappa=1$ to achieve the statistically optimal sample size.

For initializing general linear functional inference without sample splitting, the signal tensor $\mcT$ is estimated as
$\widehat{\mcT}^{\text{init}} = \mcT + \frac{\overline{\mcX}}{\|\overline{\mcX}\|} \cdot \sqrt{\frac{\op  \oor}{n}},$
where $\overline{\mcX}$ is the average of the design tensors. This initialization provides a dependent estimate of $\mcT$ at the minimax optimal rate $\sqrt{\op \oor/n}$. The initial singular subspaces $U_j^{\text{init}} \in \mathbb{R}^{p_j \times r_j}$ are then computed using 
\texttt{HOSVD} of $\widehat\mcT^{\text{init}}$.
For inferring a low-Tucker-rank linear functional with sample splitting, the dataset
is divided into two equal subsets \RN{1} and \RN{2}. The signal tensor is initialized separately for each subset as
$
\widehat{\mcT}^{\text{init}, (\mathrm{\RN{1}})} = \mcT +  \frac{\overline{\mcX}^{(\mathrm{\RN{1}})}}{\|\overline{\mcX}^{(\mathrm{\RN{1}})}\|} \cdot \sqrt{\frac{\op  \oor}{n_1}}, \widehat{\mcT}^{\text{init}, (\mathrm{\RN{2}})} = \mcT + \frac{\overline{\mcX}^{(\mathrm{\RN{2}})}}{\|\overline{\mcX}^{(\mathrm{\RN{2}})} \|} \cdot \sqrt{\frac{\op  \oor}{n_2}},
$
where $n_1$ and $n_2$ are the sample sizes for the respective subsets.
This maintains the minimax optimal rate of $\sqrt{\op  \oor/ n}$. Initial singular subspaces are obtained via \texttt{HOSVD} for each $\widehat\mcT^{\text{init},(i)}$, and the debiased estimate is constructed using the data from the alternate subset.

\begin{figure}
        \centering
        \includegraphics[page=1, width=\textwidth]{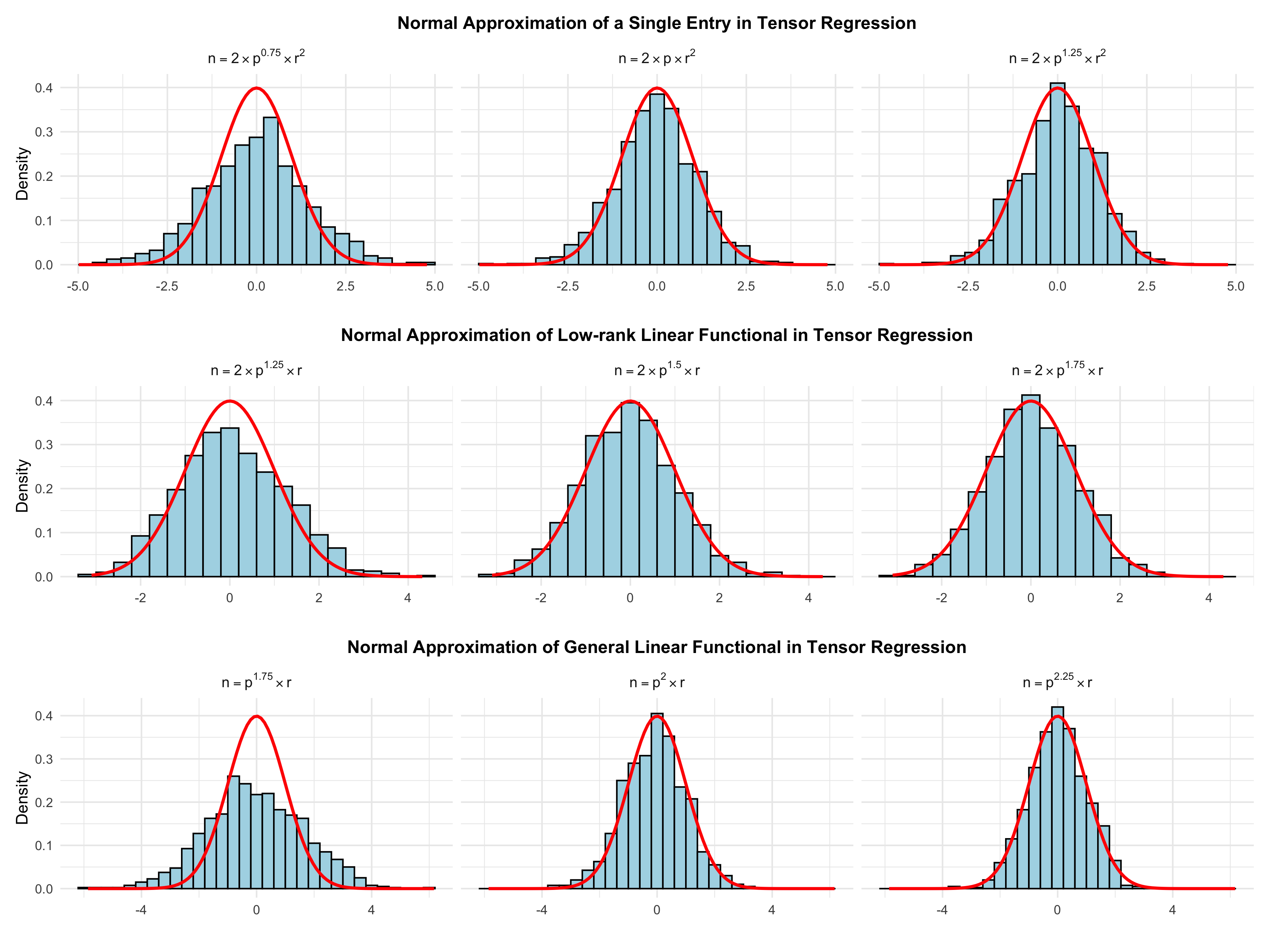}
    \caption{Histogram of normal approximation under the tensor regression setting based on 1000 independent replications, with $\op=40$ and $\oor=3$. For single-entry inference, $n \in \{2\op^{3/4}\oor, 2\op\oor, 2\op^{5/4}\oor\}$. For low-Tucker-rank linear functional inference, $n \in \{2\op^{5/4}\oor, 2\op^{3/2}\oor, 2\op^{7/4}\oor\}$. For general linear functional inference, $n \in \{\op^{7/4}\oor, \op^{2}\oor, \op^{9/4}\oor\}$. } 
    \label{fig: simulation_regression}
\end{figure}

For inferring a general linear functional, the sample size $n$ is varied around $n = \op^2\oor$, specifically using $n \in \{\op^{7/4}\oor, \op^{2}\oor$ and $\op^{9/4}\oor\}$. For inferring a low-Tucker-rank linear functional, $n$ is varied around the computationally optimal sample size $n = 2\op^{3/2}\oor$, with $n \in \{2\op^{5/4}\oor, 2\op^{3/2}\oor$ and  $2p^{7/4}\oor\}$, assuming $\kappa=\op^{1/4}$. For entrywise inference, $n$ is varied around the statistically optimal sample size $2\op\oor^2$, with $n \in \{2\op^{3/4}\oor^2, 2\op\oor^2$ and $2p^{5/4}\oor^2\}$ assuming $\kappa=1$. To account for potential power loss due to sample splitting, we set a multiplier of 2 for the sample size. Figure \ref{fig: simulation_regression} demonstrates the normal approximations achieved across all inference scenarios, provided the sample size requirements are met.

\subsection{Simulation under the Tensor PCA Setting}

In the tensor PCA setting, the noise tensor $\mcZ$ is generated with i.i.d. standard normal entries, satisfying $\sigma = 1$. 
The observed tensor is constructed as $\mcY = \mcT + \mcZ$, and the initial singular subspaces $\widehat{U}_j^{\text{init}} \in \mathbb{R}^{p_j \times r_j}$ are obtained via HOSVD on $\mcY$. Furthermore, the dimensions are set to $p_1=p_2=p_3=100:=\op$.

\begin{figure}
    \centering
        \centering
        \includegraphics[page=2, width=\textwidth]{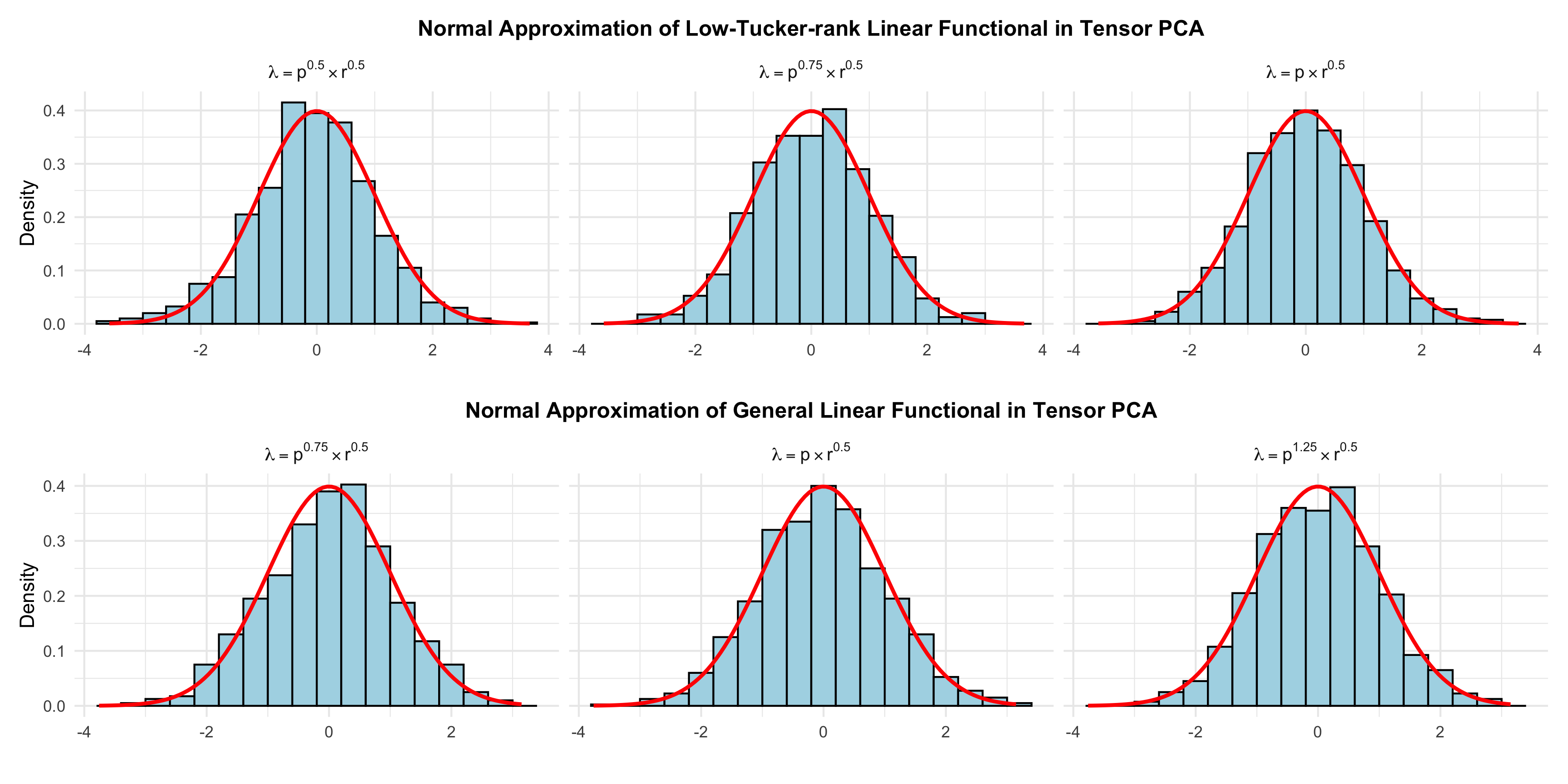}
        \caption{Histogram of normal approximation under the tensor PCA setting based on 1000 independent replications, with $\op=100$ and $\oor=3$. For low-Tucker-rank linear functional inference, $\ulambda \in \{\op^{1/2}\oor^{1/2}, \op^{3/4}\oor^{1/2}, \op\oor^{1/2}\}$. For general linear functional inference, $\ulambda \in \{\op^{3/4}\oor^{1/2}, \op\oor^{1/2}, \op^{5/4}\oor^{1/2}\}$. }
        \label{fig: simulation_PCA}
\end{figure}

For general linear functional inference, the signal strength $\ulambda$ is varied as $\underline{\lambda} \in \{\op^{1/2}\oor^{1/2}, \op^{3/4}\oor^{1/2}, \op\oor^{1/2}\}$, with the condition number set to $\kappa = \op^{1/4}$. For low-Tucker-rank linear functional inference, the signal strength $\ulambda$ is varied as $\ulambda \in \{\op\oor^{3/4}, \op\oor^{1/2}, \op^{5/4}\oor^{1/2}\}$, also with $\kappa = \op^{1/4}$. Figure \ref{fig: simulation_PCA} illustrates the normal approximations achieved under varying signal strengths and condition numbers for both inference scenarios.

To benchmark our approach against existing work \citet{agterberg2024statistical}, we consider the following setting for entrywise inference in tensor PCA. The signal tensor is defined as $\mathcal{T} = u_1 \otimes u_2 \otimes u_3$, where the singular vector for mode-$j$ is given by
$$
u_j=\frac{1}{\sqrt{4\op^{1/2} + (\op-1)}} \cdot  \big(2\cdot \op^{1/4}, \underbrace{1, \cdots, 1}_{(\op-1)\ \text{copies of 1}} \big)^{\top}.
$$ 
This construction satisfies the incoherence condition $\|u_j\|_{2, \infty} \asymp (\op^{-1/4})$. The linear functional is defined as $\mathcal{A} = e_1 \otimes e_1 \otimes e_1$, where $e_1 \in \mathbb{R}^{p}$ is a unit vector with the first entry equal to $1$. 
Our results, summarized in Table \ref{tab:coverage}, demonstrate that the coverage rates of confidence intervals align closely with the theoretical guarantees. This suggests that the characterization of asymptotic variance in our analysis is more precise than the one provided in \citet{agterberg2024statistical} under these settings.

\begin{table}
\renewcommand{\arraystretch}{1.5}
\small
\begin{tabular}{|c|c|c|c|c|c|}
\hline
& \multicolumn{2}{|c|}{$\mathcal{T}_{1,1,1}$} & \multicolumn{2}{|c|}{$\mathcal{T}_{1,1,1} - \mathcal{T}_{2,2,2}$} \\
\cline{2-5}
SNR & \multicolumn{2}{|c|}{Empirical Coverage Rate of CI} & \multicolumn{2}{|c|}{Empirical Coverage Rate of CI}\\
\cline{2-5}
& Proposed Method & Competing Method & Proposed Method & Competing Method \\
\hline
$\lambda / \sigma = \op^{1/2}\oor^{1/2}$ & 0.967 & 0.974 & 0.967 & 0.976\\
$\lambda / \sigma = \op^{3/4}\oor^{1/2}$ & 0.953 & 0.979 & 0.953 & 0.980\\
$\lambda / \sigma = \op\oor^{1/2}$ & 0.953 & 0.971 & 0.953 & 0.974\\
\hline
Average CI Length & 1.756 & 1.953 & 1.758 & 1.954\\
\hline
\end{tabular}
\caption{Empirical coverage rates of confidence intervals of a single entry $\mathcal{T}_{1,1,1}$, and the difference between two entries, $\mathcal{T}_{1,1,1} - \mathcal{T}_{2,2,2}$, using plug-in variance estimates from
\eqref{eq: estimate of sigma2 in tensor PCA} over 1000 independent replications. The results compare the performance of our proposed method with the competing method from \citet{agterberg2024statistical}.} \label{tab:coverage}
\end{table}

\section{Discussion} \label{sec: Discussion}

In this work, we present a unified framework for statistical inference on general linear functionals of signal tensors in both tensor regression and tensor PCA settings. Our approach nearly achieves computationally optimal sample size or signal-to-noise ratio (SNR) requirements for general linear functionals and precisely meets these optimal requirements for low-rank linear functionals, paving the way for further advancements in tensor-based statistical inference.

While our analysis focuses on sub-Gaussian noise, future research could explore the robustness of our framework under heavy-tailed distributions or Huber contamination models to broaden its applicability. Robust estimation techniques for matrix and tensor parameters have been studied in areas such as low-rank matrix recovery \citep{yu2024low}, matrix completion \citep{wang2024robust}, and tensor decomposition \citep{shen2023quantile}. However, inference under heavy-tailed noise remains relatively unexplored. Another promising direction for future work is extending our framework to handle structured sparse loading tensors, as discussed in \citet{zhang2019optimal}.

Practical implementation of our framework on large-scale datasets necessitates efficient computation and storage of projection matrices, alongside iterative updates in high-dimensional settings. Future studies might focus on algorithmic enhancements, such as randomized or distributed approaches, to handle large tensors while reducing computational and memory demands. For instance, sketching algorithms have been explored in tensor estimation problems, including tensor regression \citep{zhang2020islet} and tensor PCA \citep{malik2018low}. However, the application of sketching techniques for statistical inference remains an open area for exploration.


While the discussion primarily focuses on scenarios with i.i.d. sub-Gaussian noise in tensor PCA or design tensors in tensor regression, the framework can be extended to handle heteroskedastic sub-Gaussian noise, as shown in \citet{agterberg2024statistical} for tensor PCA. This extension holds if the entrywise noise of the design tensor in tensor regression or the observational noise in tensor PCA, denoted by $\sigma_{j_1, j_2, j_3}$, satisfies $\underline{\sigma} \leq \sigma_{j_1, j_2, j_3} \leq \overline{\sigma}$.

Assume the vectorized design tensor, $\operatorname{Vec}(\mcX_i) \in \mathbb{R}^{p_1 p_2 p_3}$, is a sub-Gaussian random vector with mean zero and covariance matrix $\Sigma \in \mathbb{R}^{p_1 p_2 p_3 \times p_1 p_2 p_3}$. The debiased initial estimator is then expressed as
$$
\widehat{\mcT}^{\text{unbs}}=\widehat{\mcT}^{\text{init}} + \frac{1}{n}\sum_{i=1}^{n}\big(y_i - \langle \mcX_i, \widehat{\mcT}^{\text{init}}\rangle \big) \operatorname{Vec}^{-1}\big(\Sigma^{-1}(\operatorname{Vec}(\mcX_i))\big).
$$
By following the procedures outlined in Section \ref{sec:debias_nonsplit} and Section \ref{sec: Debiased Estimator of Linear Functionals with Sample Splitting}, we can construct a similar debiased projected estimator $\widehat\cT$.
Under the conditions of Theorem~\ref{thm: main theorem in tensor regression without sample splitting} or Theorem~\ref{thm: main theorem in tensor regression with sample splitting}, the asymptotic normality follows
$$
\sqrt{n} \big( \langle\widehat{\mcT}, \mcA\rangle-\langle\mcT, \mcA\rangle \big)/(\sigma_{\xi} \widetilde{s}_{\mcA} ) \stackrel{d}{\rightarrow} N(0,1),
$$
where the variance component $\widetilde{s}_{\mcA}^2$ is adjusted compared to the i.i.d. case $s_{\mcA}^2$ in \eqref{eq: variance component sA2}, 
$$
\widetilde{s}_{\mcA}^2 = \sum_{j=1}^3\big\|\Sigma^{-\frac{1}{2}} \operatorname{Vec}\big(\mcP_{U_{j\perp}} A_j \mcP_{\left(U_{j+2} \otimes U_{j+1}\right) G_j^{\top}}\big)\big\|_{\mathrm{F}}^2+\| \Sigma^{-\frac{1}{2}}(\mcP_{U_3} \otimes \mcP_{U_2} \otimes \mcP_{U_1}) \operatorname{Vec}(\mcA) \|_2^2.
$$
Estimating the covariance matrix of the covariate tensor, $\Sigma$, typically requires additional structural assumptions, which are beyond the scope of this paper. 


\bibliographystyle{apalike}
\bibliography{citations} 

\newpage
\newpage

\appendix

\begin{center}
{\LARGE {Supplementary Material to ``Statistical Inference for Low-Rank Tensor Models''}}
\end{center}

\bigskip \centerline{\large Ke Xu, Elynn Chen, and Yuefeng Han}

\centerline{\large University of Notre Dame, and New York University}

\AppendixToC

\pagenumbering{arabic}

\newpage


This appendix is structured into two main parts, providing supplementary details and proofs supporting the main text.

The first part contains the technical proofs of the main theorems. 
We start with a proof sketch of the main theorems in Section \ref{sec: proof sketch of main theorems}. While the proofs for tensor regression and tensor PCA share similar ideas, the regression case is considerably more complex. The rest of the appendix is structured as follows: in Section \ref{sec: proof of main theorem in tensor regression without sample splitting} and Section \ref{sec: proof of main theorem in tensor regression with sample splitting}, we prove the asymptotic normality of the estimated linear functional $\langle \widehat{\mcT}, \mcA \rangle$ under the setting of tensor regression without or with sample splitting, corresponding to Theorem~\ref{thm: main theorem in tensor regression without sample splitting} and Theorem~\ref{thm: main theorem in tensor regression with sample splitting}, respectively. We also prove the asymptotic normality of the estimated linear functional $\langle \widehat{\mcT}, \mcA \rangle$ under the setting of tensor PCA, corresponding to Theorem~\ref{thm: main theorem in tensor PCA}, which closely parallels the technical results in the regression case.

We consider scenarios where variance components of the signal tensor and noise are estimated from observed data. The asymptotic normality of the estimated linear functional with these estimated variances is discussed in Section \ref{sec: Proof of Asymptotic Normality with Plug-in Estimates}, specifically for tensor regression in Section \ref{subsec: Proof of Asymptotic Normality with Plug-in Estimates in Tensor Regression} and for tensor PCA in Section \ref{subsec: Proof of Asymptotic Normality with Plug-in Estimates in Tensor PCA}.

Furthermore, Section \ref{sec: Proof of Minimax Optimal Length of Confidence Interval} establishes the minimax optimality of the confidence interval length under tensor regression (see Section \ref{subsec: Proof of Minimax Optimal Length of Confidence Interval in Tensor Regression}) and tensor PCA (see Section \ref{subsec: Proof of Minimax Optimal Length of Confidence Interval in Tensor PCA}).

The second part of the appendix adds more details to the main text. Several examples of confidence intervals for inferring low-Tucker-rank linear functionals are provided in Section \ref{sec: Specific Examples}. In addition, we summarize the algorithms of estimating the underlying linear functional $\langle \mcA, \mcT \rangle$. Specifically, the inference procedures for tensor regression (with and without sample splitting) and tensor PCA are summarized in Section \ref{sec: Inference procedure for tensor regression} and Section \ref{sec: Inference procedure for tensor PCA}.

\section{Proof Sketch of Main Theorems} \label{sec: Proof Sketch of Main Theorems} \label{sec: proof sketch of main theorems}

In this section, we outline the proof of the main theorem for tensor regression without sample splitting (Theorem~\ref{thm: main theorem in tensor regression without sample splitting}). The proof of asymptotic normality under the sample splitting setting, with adjustments for data partitioning follows a similar framework. Notably, sample splitting simplifies the analysis by eliminating the dependence between the initial estimation and the bias correction.

The proof of Theorem~\ref{thm: main theorem in tensor PCA} in Section \ref{sec: tensor_pca}, which establishes the asymptotic normality of the estimated linear functional $\langle \widehat{\mcT}, \mcA \rangle$ in the tensor PCA setting, also aligns closely with the tensor regression framework. However, a key distinction is that tensor PCA does not require an initial estimate, setting it apart from tensor regression.

Our objective is to analyze the asymptotic distribution of the linear functional $\langle \widehat{\mcT} - \mcT, \mcA \rangle$, which can be decomposed as:
\begin{align*}
\left\langle\widehat{\mcT}  -\mcT, \mcA \right\rangle
= &   \underbrace{\left\langle\widehat{\mcZ} \times_1 \mcP_{\whU_1} \times_2 \mcP_{\whU_2} \times_3 \mcP_{\whU_3}, \mcA\right\rangle}_{\text{Step 1}} + \underbrace{\left\langle \mcT \times_1 \mcP_{\whU_1} \times_2 \mcP_{\whU_2} \times_3 \mcP_{\whU_3} -\mcT, \mcA \right\rangle}_{\text{Step 2}}.
\end{align*}

We analyze each term separately.

\subsubsection*{Step 1: Analyze $\langle\widehat\mcZ \times_1 \mcP_{\whU_1} \times_2 \mcP_{\whU_2} \times_3 \mcP_{\whU_3}, \mcA\rangle$}

We further decompose $\mcP_{\whU_j} $ as $ (\mcP_{\whU_j} - \mcP_{U_j}) + \mcP_{U_j}$ for $j=1,2,3$. This decomposition introduces negligible terms, which we aim to bound
\begin{align*}
\text{Step 1.1}\quad & \left\langle \widehat{\mcZ}\times_j \left(\mcP_{\whU_j} - \mcP_{U_j}\right) \times_{j+1} \mcP_{U_{j+1}} \times_{j+2} \mcP_{U_{j+2}}, \mcA \right\rangle, \\
\text{Step 1.2}\quad & \left\langle \widehat{\mcZ} \times_j \left(\mcP_{\whU_j} - \mcP_{U_j}\right) \times_{j+1} \left(\mcP_{\whU_{j+1}} -\mcP_{U_{j+1}}\right) \times_{j+2} \mcP_{U_{j+2}} ,  \mcA \right\rangle , \\ 
\text{Step 1.3}\quad & \left\langle \widehat{\mcZ} \times_j \left(\mcP_{\whU_j}- \mcP_{U_j}\right) \times_{j+1} \left(\mcP_{\whU_{j+1}}- \mcP_{U_{j+1}}\right) \times_{j+2} \left(\mcP_{\whU_{j+2}}- \mcP_{U_{j+2}}\right), \mcA \right\rangle .
\end{align*}

The negligibility of terms in Step 1 relies on the fact that $\mcP_{\whU_j} - \mcP_{U_j}$ is small. Using the spectral representation from \citet{xia2021normal}, we have
\begin{align*}
\mcP_{\whU_j} - \mcP_{U_j} 
= &\sum_{k_j=1}^{+\infty} S_{G_j, k_j}\left(\whE_j\right) = \underbrace{S_{G_j, 1}(\whE_j)}_{\text{first-order term}}+ \underbrace{\sum_{k_j=2}^{+\infty} S_{G_j, k_j}\left(\whE_j\right)}_{\text{higher-order terms}} ,
\end{align*}
where the perturbation is
\begin{align*}
\whE_j=\whT_j\left(\mcP_{\whU_{j+1}}\otimes \mcP_{\whU_{j+2}}\right)\whT_j^{\top}-T_j\left(\mcP_{U_{j+1}}\otimes \mcP_{U_{j+2}}\right)T_j^{\top}, 
\end{align*}
and $T_j := \operatorname{Mat}_j(\mcT)$.
The expansion terms $S_{G_j, k_j}(\whE_j)$ can be found in the supplementary material, 
which involves the projection matrices $\mcP_j^{0}=\mcP_{U_{j\perp}}$, $\mcP_j^{-s}=U_j(G_jG_j^{\top})^{-s}U_j^{\top}$, and the perturbation $\whE_j$. The first-order perturbation error terms are $\|\mcP_{U_j}\whE_jU_j\mcP_j^{-1/2}\|$, $\|\mcP_j^{-1/2}\whE_jU_{j\perp}\mcP_j^{-1/2}\|$, $\|\mcP_{U_{j\perp}}\whE_j\mcP_{U_{j\perp}}\|$,
and the spectral norm of higher-order spectral projector satisfies $\|\mcP_j^{-s}\|= \|U_j\left(G_jG_j^{\top}\right)^{-s}U_j^{\top}\|\leq \ulambda^{-2s}$. 

The low-rank structure of $\mcA=\mathcal{B} \times_1 V_1 \times_2 V_2 \times_3 V_3$ reduces the perturbation error through projection onto lower-dimensional subspaces spanned by $V_jV_j^{\top}$, resulting in new first-order perturbation error terms 
$\|V_jV_j^{\top}\mcP_{U_{j\perp}}\whE_jU_j\mcP_j^{-1/2} \|$,  $\|V_jV_j^{\top}\mcP_{U_{j\perp}}\whE_jU_{j\perp} \|.$
For example, under certain conditions, $\|V_jV_j^{\top}\mcP_{U_{j\perp}}\whE_j\mcP_j^{-1/2} \|=\mcO_p[\sqrt{\oR\log (\op)/n}+\Delta \sqrt{\op/n} ]$,
while $\|\mcP_{U_{j\perp}}\whE_j\mcP_j^{-1/2} \|=\mcO_p (\sqrt{\op/n} )$.

Since the spectral representation is applied to all three modes, it is essential to analyze the leading terms introduced by the Kronecker product across multiple modes. These terms involve complex polynomials of $\whE_j$, $j=1,2,3$.
To bound these complex error terms, we derived new concentration inequalities for sub-Gaussian polynomials. Additionally, we developed novel concentration bounds for expressions like
$\operatorname{tr}[B\whZ_1^{(2)\top}C\whZ_1^{(1)}]$ and $\operatorname{tr}[B\whZ_1^{(2)\top}C\whZ_1^{(2)}]$, where $B$ and $C$ are arbitrary fixed matrices.
These results are particularly useful in the sample-splitting case, as they relax the dependency on the accuracy of the initial estimate.

\subsubsection*{Step 2: Analyze $\langle\mcT \times_1 \mcP_{\whU_1} \times_2 \mcP_{\whU_2} \times_3 \mcP_{\whU_3}-\mcT, \mcA\rangle$}

Since $\mcP_{\whU_j} = (\mcP_{\whU_j} - \mcP_{U_j}) + \mcP_{U_j}$ for $j=1,2,3$, we analyze the following terms separately
\begin{align*}
\text{Step 2.1}\quad & \left\langle\mcT \times_j \left(\mcP_{\whU_j} - \mcP_{U_j}\right) \times_{j+1} \mcP_{U_{j+1}} \times_{j+2} \mcP_{U_{j+2}}, \mcA \right\rangle, \\
\text{Step 2.2}\quad & \left\langle\mcT \times_j \left(\mcP_{\whU_j} - \mcP_{U_j}\right) \times_{j+1} \left(\mcP_{\whU_{j+1}} -\mcP_{U_{j+1}}\right) \times_{j+2} \mcP_{U_{j+2}} ,  \mcA \right\rangle , \\ 
\text{Step 2.3}\quad & \left\langle \mcT \times_j \left(\mcP_{\whU_j}- \mcP_{U_j}\right) \times_{j+1} \left(\mcP_{\whU_{j+1}}- \mcP_{U_{j+1}}\right) \times_{j+2} \left(\mcP_{\whU_{j+2}}- \mcP_{U_{j+2}}\right), \mcA \right\rangle .
\end{align*} 
Our goal is to derive upper bounds for the negligible terms in Steps 2.2 and 2.3. Step 2.1 includes both a first-order component essential for establishing asymptotic normality and higher-order negligible terms.

By combining the decomposition in Step 2.1 with the upper bounds from Steps 2.2 and 2.3, we obtain,
\begin{align*}
& \left\langle\mcT \times_1 \mcP_{\whU_1} \times_2 \mcP_{\whU_2} \times_3 \mcP_{\whU_3} - \mcT, \mcA \right\rangle = \sum_{j=1}^3 \underbrace{\left\langle \mcP_{U_{j \perp}} \whZ_j \mcP_{\left(U_{j+2} \otimes U_{j+1}\right) G_j^{\top}}, A_j\right\rangle}_{\text{approximately normal terms in Step 2.1}} + \left(\text{negligible terms}\right) .
\end{align*}
The first-order perturbation term $S_{G_j,1}(\whE_j)=\mcP_{U_{j\perp}} \whE_j \mcP_{U_{j}}+ \mcP_{U_{j}} \whE_j \mcP_{U_{j\perp}}$ is crucial for the normal approximation. The characterization of upper bounds for the negligible terms follows the same approach as in Step 1.


\subsubsection*{Step 3: Analyze the Asymptotic Normal Terms in Step 1 and Step 2}

Finally, we analyze the terms 
$
\langle \widehat{\mcZ} \times_1 \mcP_{U_1} \times_2 \mcP_{U_2} \times_3 \mcP_{U_3}, \mcA \rangle + \sum_{j=1}^{3}  \langle \mcP_{U_{j \perp}} \whZ_j \mcP_{\left(U_{j+2} \otimes U_{j+1}\right) G_j^{\top}}, A_j  \rangle
$. The noise term $
\widehat{\mcZ}
= (n\sigma^2)^{-1}\sum_{i=1}^{n} \xi_i\mathcal{X}_i+ (n\sigma^2)^{-1}\sum_{i=1}^{n} [\langle \mathcal{X}_i, \widehat{\Delta} \rangle \mathcal{X}_i -\sigma^2\cdot \widehat{\Delta}]= \widehat{\mcZ}^{(1)} + \widehat{\mcZ}^{(2)},
$
where $\whZ^{(2)}$ is negligible if the error bound of the initial estimator $\widehat{\Delta}$ is sufficiently small. The leading-order term then becomes
\begin{align*}
& \left\langle \widehat{\mcZ}^{(1)} \times_1 \mcP_{U_1} \times_2 \mcP_{U_2} \times_3 \mcP_{U_3}, \mcA\right\rangle + \sum_{j=1}^{3} \left\langle \mcP_{U_{j \perp}} \whZ_j^{(1)} \mcP_{\left(U_{j+2} \otimes U_{j+1}\right) G_j^{\top}}, A_j \right\rangle,
\end{align*}
which has variance $(n\sigma^2)^{-1}\sigma_\xi^2 \cdot \left(\|\mcA \times_1 U_1 \times_2 U_2 \times_3 U_3  \|_{\mathrm{F}}^2 + \sum_{j=1}^{3} \|U_{j \perp}^{\top} A_j \mcP_{ (U_{j+2} \otimes U_{j+1}) G_j^{\top}}\|_{\mathrm{F}}^2\right)$.
The negligible term is
$\langle \widehat{\mcZ}^{(2)} \times_1 \mcP_{U_1} \times_2 \mcP_{U_2} \times_3 \mcP_{U_3}, \mcA \rangle + \sum_{j=1}^{3}  \langle \mcP_{U_{j \perp}} \whZ_j^{(2)} \mcP_{\left(U_{j+2} \otimes U_{j+1}\right) G_j^{\top}}, A_j \rangle$ and depends on the initial estimate error bound $\Delta$.

Asymptotic normality is established by interpreting the leading-order term as a weighted sum of the entries of $\widehat{\mcZ}^{(1)}$ corresponding to the entries of $\mcA \times_1 \mcP_{U_1} \times_2 \mcP_{U_2} \times_3 \mcP_{U_3} $ and $\mcP_{U_{j\perp}}A_j\mcP_{(U_{j+2} \otimes U_{j+1}) G_j^{\top}}$ for $j=1,2,3$. When sample splitting is employed, the dependence of initial estimate and bias-correction terms is removed, and thus negligible terms decay at faster rates. 

\section{Proof of Theorem~\ref{thm: main theorem in tensor regression without sample splitting}} \label{sec: proof of main theorem in tensor regression without sample splitting} 



We begin by establishing that certain events hold with high probability. Specifically, for each mode $j=1,2,3$, assume that $\left\|\mcP_{\whU_j^{(0)}} - \mcP_{U_j}\right\| \leq \frac{\sigma_\xi}{\sigma}\sqrt{\frac{\op}{n}}$, which holds with probability at least probability at least $1- \mathbb{P}\left(\mcE_{U}^{\text{reg}}\right)$ for $j=1,2,3$, , where the event $\mcE_{U}^{\text{reg}}$ is defined as 
$\mcE_{U}^{\text{reg}} = \left\{\max_{j}\|\mcP_{\whU_j^{(0)}} - \mcP_{U_j}\| > \frac{\sigma_\xi}{\sigma}\sqrt{\frac{\op}{n}}\right\}$.

Applying Lemma~\ref{lemma: error contraction of l2 error of singular space in tensor regression}, we deduce that
\begin{align*}
\left\|\mcP_{\whU_j^{(1)}} - \mcP_{U_j}\right\| \leq \frac{\sigma_\xi}{\sigma}\sqrt{\frac{\op}{n}} \quad \text{and} \quad
\left\|\mcP_{\whU_j} - \mcP_{U_j}\right\|= \left\|\mcP_{\whU_j^{(2)}} - \mcP_{U_j}\right\| \leq \frac{\sigma_\xi}{\sigma}\sqrt{\frac{\op}{n}},
\end{align*}
which hold with probability at least $1-\exp(-c\op) - \mathbb{P}\left(\mcE_U^{\text{reg}}\right)$ for each $j=1,2,3$.

Furthermore, we assume that the initial estimation error satisfies
$
\left\|\whT^{\text{init}} - \mcT\right\|_{\mathrm{F}} \leq \Delta
$
which holds with probability at leas $1-\mathbb{P}(\mcE_{\Delta})$. Here, the event $\mcE_{\Delta}$ is defined as $
\mcE_{\Delta}= \left\{\left\|\whT^{\text{init}} - \mcT\right\|_{\mathrm{F}} > \Delta \right\}$.

With these high-probability events established, we proceed through the following steps to complete the proof of the main theorem.

\subsection*{Step 1: Upper Bound of Negligible Terms in $\langle \widehat{\mcZ} \times_1 \mcP_{\widehat{U}_1} \times_2 \mcP_{\widehat{U}_2} \times_3 \mcP_{\widehat{U}_3}, \mcA \rangle$}

In Step 3, we establish the asymptotic normality of $\left\langle \widehat{\mcZ} \times_1 \mcP_{U_1} \times_2 \mcP_{U_2} \times_3 \mcP_{U_3}, \mcA \right\rangle$, by deriving the Berry-Esseen bound for its components. To facilitate this, it is essential to quantify upper bounds for the negligible terms. Due to symmetry, it suffices to consider the upper bounds of the following terms:
\begin{align}
\text{(Step 1.1)}: & \left|\left\langle \widehat{\mcZ} \times_1 \left(\mcP_{\whU_1} - \mcP_{U_1}\right) \times_2 \mcP_{U_2} \times_3 \mcP_{U_3}, \mcA \right\rangle\right|,  \label{eq: step 1.1 in tensor regression without sample splitting}\\
\text{(Step 1.2)}: & \left|\left\langle \widehat{\mcZ} \times_1 \left(\mcP_{\whU_1} - \mcP_{U_1}\right) \times_2 \left(\mcP_{\whU_2} - \mcP_{U_2}\right) \times_3 \mcP_{U_3}, \mcA  \right\rangle\right|, \label{eq: step 1.2 in tensor regression without sample splitting} \\
\text{(Step 1.3)}: & \left|\left\langle \widehat{\mcZ} \times_1 \left(\mcP_{\whU_1} - \mcP_{U_1}\right) \times_2 \left(\mcP_{\whU_2} - \mcP_{U_2}\right) \times_3 \left(\mcP_{\whU_3} - \mcP_{U_3}\right), \mcA \right\rangle\right|. \label{eq: step 1.3 in tensor regression without sample splitting} 
\end{align}

\subsubsection*{Step 1.1: Upper Bound of Negligible Terms in $\langle \widehat{\mcZ} \times_1 \left(\mcP_{\widehat{U}_1} - \mcP_{U_1}\right) \times_2 \mcP_{U_2} \times_3 \mcP_{U_3}, \mcA \rangle$} 

First, consider the following decomposition:
\begin{align}
& \left|\left\langle \widehat{\mcZ} \times_1 \left(\mcP_{\whU_1} - \mcP_{U_1}\right) \times_2 \mcP_{U_2} \times_3 \mcP_{U_3}, \mcA\right\rangle\right| \notag \\
\leq & \underbrace{\left|\operatorname{tr}\left[\left(\mcP_{U_3}\otimes \mcP_{U_2}\right)A_1^{\top}\mcP_{U_1}\left(\mcP_{\whU_1}-\mcP_{U_1}\right)\mcP_{U_1}\whZ_1\left(\mcP_{U_3}\otimes \mcP_{U_2}\right)\right]\right|}_{\mathrm{\RN{1}}} \label{eq: term 1 in step 1.1 in tensor regression without sample splitting}\\
+ & \underbrace{\left|\operatorname{tr}\left[\left(\mcP_{U_3}\otimes \mcP_{U_2}\right)A_1^{\top}V_1V_1^{\top}\mcP_{U_{1\perp}}\left(\mcP_{\whU_1}-\mcP_{U_1}\right)\mcP_{U_1}\whZ_1\left(\mcP_{U_3}\otimes \mcP_{U_2}\right)\right]\right|}_{{\mathrm{\RN{2}}}} \label{eq: term 2 in step 1.1 in tensor regression without sample splitting}\\
+ & \underbrace{\left|\operatorname{tr}\left[\left(\mcP_{U_3}\otimes \mcP_{U_2}\right)A_1^{\top}\mcP_{U_1}\left(\mcP_{\whU_1}-\mcP_{U_1}\right)\mcP_{U_{1\perp}}\whZ_1\left(\mcP_{U_3}\otimes \mcP_{U_2}\right)\right]\right|}_{{\mathrm{\RN{3}}}} \label{eq: term 3 in step 1.1 in tensor regression without sample splitting}\\
+ & \underbrace{\left|\operatorname{tr}\left[\left(\mcP_{U_3}\otimes \mcP_{U_2}\right)A_1^{\top}V_1V_1^{\top}\mcP_{U_{1\perp}}\left(\mcP_{\whU_1}-\mcP_{U_1}\right)\mcP_{U_{1\perp}}\whZ_1\left(\mcP_{U_3}\otimes \mcP_{U_2}\right)\right]\right|}_{{\mathrm{\RN{4}}}} \label{eq: term 4 in step 1.1 in tensor regression without sample splitting}.
\end{align}

We begin with the upper bound for the first term \RN{1} \eqref{eq: term 1 in step 1.1 in tensor regression without sample splitting}: 
\begin{align}
\mathrm{\RN{1}}
\leq & \left\|\left(\mcP_{U_3}\otimes \mcP_{U_2}\right)A_1^{\top}\mcP_{U_1}\right\|_{\mathrm{F}}\cdot \underbrace{\left\|\mcP_{U_1}\left(\mcP_{\whU_1}-\mcP_{U_1}\right)\mcP_{U_1}\right\|}_{\eqref{eq: high-prob upper bound of PU1(PUhat1-PU1)PU1 in tensor regression}} \cdot \underbrace{\left\|\mcP_{U_1}\whZ_1\left(\mcP_{U_3}\otimes \mcP_{U_2}\right)\right\|_{\mathrm{F}}}_{\eqref{eq: high-prob upper bound of U1Zhat1(U3oU2) in tensor regression without sample splitting}} \notag\\
\leq & \left\|\mcA \times_1 U_1 \times_2 U_2 \times_3 U_3\right\|_{\mathrm{F}} \cdot \left[\frac{\sigma_\xi^3\oor^{1/2}}{\ulambda^2\sigma^3}\cdot \left(\frac{\op\sqrt{\oor\log(\op)}}{n^{3/2}}+\Delta\cdot \frac{\op^{3/2}}{n^{3/2}}\right)\right]. \label{eq: upper bound of term 1 in step 1.1 in tensor regression without sample splitting}
\end{align}

For the second term $\mathrm{\RN{2}}$ \eqref{eq: term 2 in step 1.1 in tensor regression without sample splitting}, we have
\begin{align}
\mathrm{\RN{2}}
\leq & \underbrace{\left\|\left(\mcP_{U_3}\otimes \mcP_{U_2}\right)A_1^{\top}\mcP_{U_{1\perp}}\left(\mcP_{\whU_1}-\mcP_{U_1}\right)\mcP_{U_1}\right\|_{\mathrm{F}}}_{\eqref{eq: high-prob upper bound of V1tPU1p(PUhat1-PU1)U1 in tensor regression without sample splitting}} \cdot \underbrace{\left\|\mcP_{U_1}\whZ_1\left(\mcP_{U_3}\otimes \mcP_{U_2}\right)\right\|_{\mathrm{F}}}_{\eqref{eq: high-prob upper bound of U1Zhat1(U3oU2) in tensor regression without sample splitting}} \notag\\
\lesssim & \left\|\mcA \times_2 U_2 \times_3 U_3\right\|_{\mathrm{F}} \cdot \left[\frac{\sigma_\xi^2\oor^{1/2}}{\ulambda\sigma^2}\left(\frac{\oor\log(\op)}{n}+\Delta^2\cdot \frac{\op}{n}\right)\right] \label{eq: upper bound of term 2 in step 1.1 in tensor regression without sample splitting}.
\end{align}

For the third term $\mathrm{\RN{3}}$ \eqref{eq: term 3 in step 1.1 in tensor regression without sample splitting}, we first have 
\begin{align}
\mathrm{\RN{3}} 
\leq & \underbrace{\left|\operatorname{tr}\left[\left(\mcP_{U_3}\otimes \mcP_{U_2}\right)A_1^{\top}\mcP_{U_1} S_{G_1,1}\left(\whE_1\right)\mcP_{U_{1\perp}}\whZ_1\mcP_{U_2}\otimes \mcP_{U_3}\right]\right|}_{\mathrm{\RN{3}}.\mathrm{\RN{1}}} \label{eq: term 3.1 in step 1.1 in tensor regression without sample splitting}\\
+ & \underbrace{\left|\operatorname{tr}\left[\left(\mcP_{U_3}\otimes \mcP_{U_2}\right)A_1^{\top}\mcP_{U_1}\sum_{k_1=2}^{+\infty} S_{G_1,k_1}\left(\whE_1\right)\mcP_{U_{1\perp}}\whZ_1\mcP_{U_2}\otimes \mcP_{U_3}\right]\right|}_{\mathrm{\RN{3}}.\mathrm{\RN{2}}}. \label{eq: term 3.2 in step 1.1 in tensor regression without sample splitting}
\end{align}

Note that $\whZ_1=\underbrace{\frac{1}{n\sigma^2}\sum_{i=1}^n\xi_i\Mat_1\left(\mcX_i\right)}_{\whZ_1^{(1)}}+ \underbrace{\frac{1}{n\sigma^2}\sum_{i=1}^n\left[\left\langle\mcX_i, \widehat{\Delta}\right\rangle\Mat_1\left(\mcX_i\right)-\sigma^2\cdot \widehat{\Delta}_1\right]}_{\whZ_1^{(2)}} $. Then by Lemma~\ref{lemma: high-prob upper bound of tr(BZ(1)tCZ(1))}, for the term $\mathrm{\RN{3}}.\mathrm{\RN{1}}$ \eqref{eq: term 3.1 in step 1.1 in tensor regression without sample splitting}, we have
\begin{align}
& \mathrm{\RN{3}}.\mathrm{\RN{1}} \notag\\
\leq & \left|\operatorname{tr}\left[\left(\mcP_{U_3}\otimes \mcP_{U_2}\right)A_1^{\top}U_1\left(G_1G_1^{\top}\right)^{-1}G_1\left(U_3\otimes U_2\right)^{\top}\whZ_1^{(1)\top}\mcP_{U_{1\perp}}\whZ_1^{(1)}\left(\mcP_{U_3}\otimes \mcP_{U_2}\right)\right]\right| \notag\notag\\ 
+ & \left|\operatorname{tr}\left[\left(\mcP_{U_3}\otimes \mcP_{U_2}\right)A_1^{\top}U_1\left(G_1G_1^{\top}\right)^{-1}G_1\left(U_3\otimes U_2\right)^{\top}\whZ_1^{(1)\top}\mcP_{U_{1\perp}}\whZ_1^{(2)}\left(\mcP_{U_3}\otimes \mcP_{U_2}\right)\right]\right| \notag\\ 
+ & \left|\operatorname{tr}\left[\left(\mcP_{U_3}\otimes \mcP_{U_2}\right)A_1^{\top}U_1\left(G_1G_1^{\top}\right)^{-1}G_1\left(U_3\otimes U_2\right)^{\top}\whZ_1^{(2)\top}\mcP_{U_{1\perp}}\whZ_1^{(1)}\left(\mcP_{U_3}\otimes \mcP_{U_2}\right)\right]\right| \notag\\ 
+ & \left|\operatorname{tr}\left[\left(\mcP_{U_3}\otimes \mcP_{U_2}\right)A_1^{\top}U_1\left(G_1G_1^{\top}\right)^{-1}G_1\left(U_3\otimes U_2\right)^{\top}\whZ_1^{(2)\top}\mcP_{U_{1\perp}}\whZ_1^{(2)}\left(\mcP_{U_3}\otimes \mcP_{U_2}\right)\right]\right| \notag\\ 
\lesssim & \underbrace{\left|\operatorname{tr}\left[\left(U_3\otimes U_2\right)A_1^{\top}U_1\left(G_1G_1^{\top}\right)^{-1}G_1\right]\operatorname{tr}\left(\mcP_{U_{1\perp}}\right)\right|\cdot \frac{\sigma_\xi^2}{\sigma^2}\frac{1}{n} + \left\|\left(\mcP_{U_3}\otimes \mcP_{U_2}\right)A_1^{\top}\mcP_{U_1}\right\|_{\mathrm{F}}\cdot \frac{\sigma_\xi^2}{\ulambda\sigma^2}\cdot \frac{\sqrt{\op\log(\op)}}{n}}_{\eqref{eq: high-prob upper bound of tr(BZ(1)tCZ(1))}} \notag
\\
+ & \left\|\left(\mcP_{U_3}\otimes \mcP_{U_2}\right)A_1^{\top}\mcP_{U_1}\right\|_{\mathrm{F}}\cdot \frac{\sigma_\xi}{\ulambda\sigma}\cdot \Delta\sqrt{\frac{\op\oor}{n}}\cdot \frac{\sigma_\xi}{\sigma}\sqrt{\frac{\op}{n}}  + \left\|\left(\mcP_{U_3}\otimes \mcP_{U_2}\right)A_1^{\top}\mcP_{U_1}\right\|_{\mathrm{F}}\cdot \frac{\sigma_\xi}{\ulambda\sigma}\cdot \Delta\sqrt{\frac{\op\oor}{n}}\cdot \frac{\sigma_\xi}{\sigma}\sqrt{\frac{\op}{n}} \notag \\
+ & \left\|\left(\mcP_{U_3}\otimes \mcP_{U_2}\right)A_1^{\top}\mcP_{U_1}\right\|_{\mathrm{F}}\cdot \frac{\sigma_\xi^2\oor^{1/2}}{\ulambda \sigma^2} \cdot\Delta^2\cdot \frac{\op}{n} \notag\\
\lesssim & \left\|\mcP_{U_1}A_1\mcP_{\left(U_3\otimes U_2\right)G_1^{\top}}\right\|_{\mathrm{F}} \cdot \frac{\sigma_\xi^2\oor^{1/2}}{\ulambda\sigma^2}\cdot \frac{\op}{n} + \left\|\mcA\times_1 U_1 \times_2 U_2 \times_3 U_3\right\|_{\mathrm{F}} \cdot \frac{\sigma_\xi^2\oor^{1/2}}{\ulambda\sigma^2} \cdot \left( \frac{\sqrt{\op\log(\op)}}{n} +\Delta\cdot \frac{\op}{n}\right). \label{eq: upper bound of the 3.1 term in step 1.1 in tensor regression without sample splitting}
\end{align}

Furthermore, for the term $\mathrm{\RN{3}}.\mathrm{\RN{2}}$ \eqref{eq: term 3.2 in step 1.1 in tensor regression without sample splitting}, we have
\begin{align}
\mathrm{\RN{3}}.\mathrm{\RN{2}} 
\leq & \left\|\left(\mcP_{U_3}\otimes \mcP_{U_2}\right)A_1^{\top}\mcP_{U_1}\right\|_{\mathrm{F}} \cdot \underbrace{\left\|\mcP_1^{-2}\whE_1\mcP_1^{0}\right\|}_{\eqref{eq: high-prob upper bound of P1(0)Ehat1P1(-1/2) in tensor regression}} \cdot \underbrace{\left\|\mcP_1^{0}\whE_1\mcP_1^{0}\right\|}_{\eqref{eq: high-prob upper bound of P1(0)Ehat1P1(0) in tensor regression}} \cdot \underbrace{\left\|\mcP_{U_{1\perp}}\whZ_1\left(\mcP_{U_3}\otimes \mcP_{U_2}\right)\right\|_{\mathrm{F}}}_{\eqref{eq: high-prob upper bound of U1pZhat1(U3oU2) in tensor regression without sample splitting}} \notag\\
+ & \left\|\left(\mcP_{U_3}\otimes \mcP_{U_2}\right)A_1^{\top}\mcP_{U_1}\right\|_{\mathrm{F}} \cdot \underbrace{\left\|\mcP_1^{-1}\whE_1\mcP_1^{-\frac{1}{2}}\right\|}_{\eqref{eq: high-prob upper bound of P1(-1/2)Ehat1P1(-1/2) in tensor regression without sample splitting}} \cdot \underbrace{\left\|\mcP_1^{-\frac{1}{2}}\whE_1\mcP_1^{0}\mcP_{U_{1\perp}}\right\|}_{\eqref{eq: high-prob upper bound of P1(0)Ehat1P1(-1/2) in tensor regression}} \cdot \underbrace{\left\|\mcP_{U_{1\perp}}\whZ_1\left(\mcP_{U_3}\otimes \mcP_{U_2}\right)\right\|_{\mathrm{F}}}_{\eqref{eq: high-prob upper bound of U1pZhat1(U3oU2) in tensor regression without sample splitting}} \notag\\
+ & \left\|\left(\mcP_{U_3}\otimes \mcP_{U_2}\right)A_1^{\top}\mcP_{U_1}\right\|_{\mathrm{F}} \cdot \left\|\mcP_{U_1}\sum_{k_1=3}^{+\infty} S_{G_1,k_1}\left(\whE_1\right)\mcP_{U_{1\perp}}\right\| \cdot \underbrace{\left\|\mcP_{U_{1\perp}}\whZ_1\left(\mcP_{U_3}\otimes \mcP_{U_2}\right)\right\|_{\mathrm{F}}}_{\eqref{eq: high-prob upper bound of U1pZhat1(U3oU2) in tensor regression without sample splitting}} \notag \\
\lesssim & \left\|\mcA\times_1 U_1 \times_2 U_2 \times_3 U_3\right\|_{\mathrm{F}} \cdot \left[\frac{\sigma_\xi^3\oor^{1/2}}{\ulambda^2\sigma^3}\cdot \left(\frac{\op\sqrt{\oor\log(\op)}}{n^{3/2}} + \Delta \cdot \frac{\op^{3/2}}{n^{3/2}}\right)\right]. \label{eq: upper bound of the 3.2 term in step 1.1 in tensor regression without sample splitting}
\end{align}

It implies that 
\begin{align}
\mathrm{\RN{3}} 
\lesssim & \left\|\mcP_{U_1}A_1\mcP_{\left(U_3\otimes U_2\right)G_1^{\top}}\right\|_{\mathrm{F}} \cdot \frac{\sigma_\xi^2\oor^{1/2}}{\ulambda\sigma^2}\cdot \frac{\op}{n} + \left\|\mcA\times_1 U_1 \times_2 U_2 \times_3 U_3\right\|_{\mathrm{F}} \cdot \left[\frac{\sigma_\xi^2\oor^{1/2}}{\ulambda\sigma^2} \left(\frac{\sqrt{\op\log(\op)}}{n} +\Delta\cdot\frac{\op}{n}\right)\right]. \label{eq: upper bound of term 3 in step 1.1 in tensor regression without sample splitting}
\end{align}

For the fourth term $\mathrm{\RN{4}}$ \eqref{eq: term 4 in step 1.1 in tensor regression without sample splitting}, we have
\begin{align}
\mathrm{\RN{4}} 
\leq & \underbrace{\left\|\left(\mcP_{U_3}\otimes \mcP_{U_2}\right)A_1^{\top}\mcP_{U_{1\perp}}\left(\mcP_{\whU_1}-\mcP_{U_1}\right)\mcP_{U_{1\perp}}\right\|_{\mathrm{F}}}_{\eqref{eq: high-prob upper bound of V1tPU1p(PUhat1-PU1)U1p in tensor regression without sample splitting}} \cdot \underbrace{\left\|\mcP_{U_{1\perp}}\whZ_1\left(\mcP_{U_3}\otimes \mcP_{U_2}\right)\right\|_{\mathrm{F}}}_{\eqref{eq: high-prob upper bound of U1pZhat1(U3oU2) in tensor regression without sample splitting}} \notag\\
\lesssim & \left\|\mcA \times_2 U_2 \times_3 U_3\right\|_{\mathrm{F}} \cdot \left[\frac{\sigma_\xi^3\oor^{1/2}}{\ulambda^2\sigma^3}\left(\frac{\op\sqrt{\oor\log(\op)}}{n^{3/2}}+\Delta\cdot \frac{\op^{3/2}}{n^{3/2}}\right)\right] \label{eq: upper bound of term 4 in step 1.1 in tensor regression without sample splitting}.
\end{align}

Therefore, we have the following upper bound for \eqref{eq: step 1.1 in tensor regression without sample splitting}
\begin{align*}
& \left| \left\langle\widehat{\mcZ} \times_1\left(\mcP_{\whU_1}-\mcP_{U_1}\right) \times_2 \mcP_{U_2} \times_3 \mcP_{U_3}, \mcA\right\rangle\right| 
\lesssim  \eqref{eq: upper bound of term 1 in step 1.1 in tensor regression without sample splitting} + \eqref{eq: upper bound of term 2 in step 1.1 in tensor regression without sample splitting} + \eqref{eq: upper bound of term 3 in step 1.1 in tensor regression without sample splitting} + \eqref{eq: upper bound of term 4 in step 1.1 in tensor regression without sample splitting} \\
\lesssim & \left\|\mcP_{U_1}A_1\mcP_{\left(U_3\otimes U_2\right)G_1^{\top}}\right\|_{\mathrm{F}} \cdot \frac{\sigma_\xi^2\oor^{1/2}}{\ulambda\sigma^2}\cdot \frac{\op}{n} + \left\|\mcA\times_1 U_1 \times_2 U_2 \times_3 U_3\right\|_{\mathrm{F}} \cdot \left[\frac{\sigma_\xi^2\oor^{1/2}}{\ulambda\sigma^2}\cdot \left(\frac{\sqrt{\op\log(\op)}}{n} +\Delta\cdot \frac{\op}{n}\right)\right] \\
+ & \left\|\mcA \times_2 U_2 \times_3 U_3\right\|_{\mathrm{F}} \cdot \left[\frac{\sigma_\xi^2\oor^{1/2}}{\ulambda\sigma^2}\left(\frac{\oor\log(\op)}{n}+\Delta^2\cdot \frac{\op}{n}\right) + \frac{\sigma_\xi^3\oor^{1/2}}{\ulambda^2\sigma^3}\left(\frac{\op\sqrt{\oor\log(\op)}}{n^{3/2}}+\Delta\cdot \frac{\op^{3/2}}{n^{3/2}}\right)\right].
\end{align*}

\subsubsection*{Step 1.2: Upper Bound of $\langle \widehat{\mcZ} \times_1 \left(\mcP_{\widehat{U}_1} - \mcP_{U_1}\right) \times_2 \left(\mcP_{\widehat{U}_2} - \mcP_{U_2}\right) \times_3 \mcP_{U_3}, \mcA \rangle$}

Then we consider the following decomposition of \eqref{eq: step 2.1 in tensor regression without sample splitting}:
\begin{align}
& \left|\left\langle \widehat{\mcZ} \times_1 \left(\mcP_{\whU_1} - \mcP_{U_1}\right) \times_2 \left(\mcP_{\whU_2} - \mcP_{U_2}\right) \times_3 \mcP_{U_3}, \mcA \right\rangle \right|\notag \\
\leq & \underbrace{\left|\left\langle \widehat{\mcZ} \times_1 \left(\mcP_{\whU_1} - \mcP_{U_1}\right) \times_2 \left(\mcP_{\whU_2} - \mcP_{U_2}\right) \times_3 \mcP_{U_3}, \mcA \times_1 \mcP_{U_1} \times_2 \mcP_{U_2} \times_3 \mcP_{U_3}\right\rangle\right|}_{\mathrm{\RN{1}}} \label{eq: term 1 in step 1.2 in tensor regression without sample splitting}\\ 
+ & \underbrace{\left|\left\langle \widehat{\mcZ} \times_1 \left(\mcP_{\whU_1} - \mcP_{U_1}\right) \times_2 \left(\mcP_{\whU_2} - \mcP_{U_2}\right) \times_3 \mcP_{U_3}, \mcA \times_1 \mcP_{U_{1\perp}} \times_2 \mcP_{U_2} \times_3 \mcP_{U_3}\right\rangle\right|}_{\mathrm{\RN{2}}} \label{eq: term 2 in step 1.2 in tensor regression without sample splitting}\\ 
+ & \underbrace{\left|\left\langle \widehat{\mcZ} \times_1 \left(\mcP_{\whU_1} - \mcP_{U_1}\right) \times_2 \left(\mcP_{\whU_2} - \mcP_{U_2}\right) \times_3 \mcP_{U_3}, \mcA \times_1 \mcP_{U_1} \times_2 \mcP_{U_{2\perp}} \times_3 \mcP_{U_3}\right\rangle\right|}_{\mathrm{\RN{3}}} \label{eq: term 3 in step 1.2 in tensor regression without sample splitting}\\
+ & \underbrace{\left|\left\langle \widehat{\mcZ} \times_1 \left(\mcP_{\whU_1} - \mcP_{U_1}\right) \times_2 \left(\mcP_{\whU_2} - \mcP_{U_2}\right) \times_3 \mcP_{U_3}, \mcA \times_1 \mcP_{U_{1\perp}} \times_2 \mcP_{U_{2\perp}} \times_3 \mcP_{U_3}\right\rangle\right|}_{\mathrm{\RN{4}}} \label{eq: term 4 in step 1.2 in tensor regression without sample splitting}.
\end{align}

Here, for the first term \RN{1} \eqref{eq: term 1 in step 1.2 in tensor regression without sample splitting}, we have
\begin{align}
\mathrm{\RN{1}} 
\leq & \left\|\left(\mcP_{U_2}\otimes \mcP_{U_1}\right)A_3^{\top}\mcP_{U_3}\right\|_{\mathrm{F}}\cdot\sup_{\substack{W_1\in \mathbb{R}^{p_1\times r_1}, \left\|W_1\right\|=1\\ W_2 \in \mathbb{R}^{p_2\times r_2}, \left\|W_2\right\|=1}} \left\|U_3^{\top}\whZ_3\left(W_2\otimes W_1\right)\right\|_{\mathrm{F}} \cdot \prod_{j=1}^2 \left\|\mcP_{\whU_j} - \mcP_{U_j}\right\| \notag\\
\lesssim & \left\|\left(\mcP_{U_2}\otimes \mcP_{U_1}\right)A_3^{\top}\mcP_{U_3}\right\|_{\mathrm{F}} \cdot \frac{\sigma_\xi}{\sigma} \sqrt{\frac{\op\oor}{n}}\cdot \left(\frac{\sigma_\xi}{\ulambda\sigma}\sqrt{\frac{\op}{n}}\right)^2 = \left\|\mcA \times_1 U_1 \times_2 U_2 \times_3 U_3 \right\|_{\mathrm{F}}\cdot \frac{\sigma_\xi^3\oor^{1/2}}{\ulambda^2\sigma^3}\cdot \frac{\op^{3/2}}{n^{3/2}}. \label {eq: upper bound of term 1 in step 1.2 in tensor regression without sample splitting}
\end{align}

Second, we decompose the second term \RN{2} \eqref{eq: term 2 in step 1.2 in tensor regression without sample splitting} as follows:
\begin{align}
\mathrm{\RN{2}}
\leq & \underbrace{\left|\operatorname{tr}\left[\mcP_{U_3}\whZ_3\left[\mcP_{U_2}\left(\mcP_{\whU_2} - \mcP_{U_2}\right)\mcP_{U_2} \otimes \mcP_{U_1}\left(\mcP_{\whU_1} - \mcP_{U_1}\right)\mcP_{U_{1\perp}}V_1V_1^{\top}\right]A_3^{\top}\mcP_{U_3}\right]\right|}_{\mathrm{\RN{2}}.\mathrm{\RN{1}}} \label{eq: term 2.1 in step 1.2 in tensor regression without sample splitting}\\
+ & \underbrace{\left|\operatorname{tr}\left[\mcP_{U_3}\whZ_3\left[\mcP_{U_2}\left(\mcP_{\whU_2} - \mcP_{U_2}\right)\mcP_{U_2} \otimes \mcP_{U_{1\perp}}\left(\mcP_{\whU_1} - \mcP_{U_1}\right)\mcP_{U_{1\perp}}V_1V_1^{\top}\right]A_3^{\top}\mcP_{U_3}\right]\right|}_{\mathrm{\RN{2}}.\mathrm{\RN{2}}} \label{eq: term 2.2 in step 1.2 in tensor regression without sample splitting}\\
+ & \underbrace{\left|\operatorname{tr}\left[\mcP_{U_3}\whZ_3\left[\mcP_{U_{2\perp}}\left(\mcP_{\whU_2} - \mcP_{U_2}\right)\mcP_{U_2} \otimes \mcP_{U_1}\left(\mcP_{\whU_1} - \mcP_{U_1}\right)\mcP_{U_{1\perp}}V_1V_1^{\top}\right]A_3^{\top}\mcP_{U_3}\right]\right|}_{\mathrm{\RN{2}}.\mathrm{\RN{3}}} \label{eq: term 2.3 in step 1.2 in tensor regression without sample splitting}\\
+ & \underbrace{\left|\operatorname{tr}\left[\mcP_{U_3}\whZ_3\left[\mcP_{U_{2\perp}}\left(\mcP_{\whU_2} - \mcP_{U_2}\right)\mcP_{U_2} \otimes \mcP_{U_{1\perp}}\left(\mcP_{\whU_1} - \mcP_{U_1}\right)\mcP_{U_{1\perp}}V_1V_1^{\top}\right]A_3^{\top}\mcP_{U_3}\right]\right|}_{\mathrm{\RN{2}}.\mathrm{\RN{4}}} \label{eq: term 2.4 in step 1.2 in tensor regression without sample splitting}. 
\end{align}
For the term $\mathrm{\RN{2}}.\mathrm{\RN{1}}$ in \eqref{eq: term 2.4 in step 1.2 in tensor regression without sample splitting}, we have
\begin{align}
\mathrm{\RN{2}}.\mathrm{\RN{1}}
\leq & \underbrace{\left\|\mcP_{U_3}\whZ_3\left(\mcP_{U_2}\otimes \mcP_{U_1}\right)\right\|_{\mathrm{F}}}_{\eqref{eq: high-prob upper bound of U1Zhat1(U3oU2) in tensor regression without sample splitting}} \cdot \underbrace{\left\|\mcP_{U_1}\big(\mcP_{\whU_1} - \mcP_{U_1}\big)\mcP_{U_{1\perp}}V_1\right\|}_{\eqref{eq: high-prob upper bound of V1tPU1p(PUhat1-PU1)U1 in tensor regression without sample splitting}} \cdot \underbrace{\left\| \mcP_{U_2}\big(\mcP_{\whU_2} - \mcP_{U_2}\big)\mcP_{U_2}\right\|}_{\eqref{eq: high-prob upper bound of PU1(PUhat1-PU1)PU1 in tensor regression}} \cdot \left\|A_1\left(\mcP_{U_3}\otimes \mcP_{U_2}\right)\right\|_{\mathrm{F}} \notag\\
\lesssim & \left\|\mcA\times_2 U_2\times_3 U_3 \right\|_{\mathrm{F}} \cdot \frac{\sigma_\xi^4\oor^{1/2}}{\ulambda^3\sigma^4} \cdot \left[\left(\frac{\op\sqrt{\oR  \oor}\log(\op)}{n^2}+\Delta^2\cdot \frac{\op^2}{n^2}\right)\right] \label{eq: upper bound of term 2.1 in step 1.2 in tensor regression without sample splitting}.
\end{align}
For the term $\mathrm{\RN{2}}.\mathrm{\RN{2}}$ in the second term \eqref{eq: term 2.4 in step 1.2 in tensor regression without sample splitting}, we have
\begin{align}
\mathrm{\RN{2}}.\mathrm{\RN{2}}
\leq & \sup_{\substack{W_1\in \mathbb{R}^{p_1\times 2r_1}, \left\|W_1\right\|=1\\ W_2\in \mathbb{R}^{p_2\times 2r_2}, \left\|W_2\right\|=1}}\left\|\mcP_{U_3}\whZ_3\left(W_2\otimes W_1\right)\right\|_{\mathrm{F}}  \underbrace{\left\|\mcP_{U_{1\perp}}\left(\mcP_{\whU_1} - \mcP_{U_1}\right)\mcP_{U_{1\perp}}V_1\right\|}_{\eqref{eq: high-prob upper bound of V1tPU1p(PUhat1-PU1)U1p in tensor regression without sample splitting}} \notag \\
& \cdot \underbrace{\left\|\mcP_{U_2}\left(\mcP_{\whU_2} - \mcP_{U_2}\right)\mcP_{U_2}\right\|}_{\eqref{eq: high-prob upper bound of PU1(PUhat1-PU1)PU1 in tensor regression}} \left\|\mcA \times_2 U_2\times_3 U_3\right\|_{\mathrm{F}} \notag\\
\lesssim & \left\|\mcA\times_2 U_2\times_3 U_3 \right\|_{\mathrm{F}} \cdot \frac{\sigma_\xi^5\oor^{1/2}}{\ulambda^4\sigma^5} \left[ \left(\frac{\op^2\sqrt{\oR\log(\op)}}{n^{5/2}} + \Delta\cdot \frac{\op^{5/2}}{n^{5/2}}\right)\right]. \label{eq: upper bound of term 2.2 in step 1.2 in tensor regression without sample splitting}
\end{align}

For the term $\mathrm{\RN{2}}.\mathrm{\RN{3}}$ in the third term \eqref{eq: term 2.3 in step 1.2 in tensor regression without sample splitting}, we have
\begin{align}
\mathrm{\RN{2}}.\mathrm{\RN{3}} 
\leq & \sup_{\substack{W_1\in \mathbb{R}^{p_1\times 2r_1}, \left\|W_1\right\|=1 \\ W_2\in \mathbb{R}^{p_2\times 2r_2}, \left\|W_2\right\|=1}} \left\|\mcP_{U_3}\whZ_3 (W_2\otimes W_1 )\right\|_{\mathrm{F}}  \underbrace{\left\|\mcP_{U_1}\big(\mcP_{\whU_1}-\mcP_{U_1}\big)\mcP_{U_{1\perp}}V_1\right\|}_{\eqref{eq: high-prob upper bound of V1tPU1p(PUhat1-PU1)U1 in tensor regression without sample splitting}} \notag \\
& \cdot \underbrace{\left\|\mcP_{U_{2\perp}}\big(\mcP_{\whU_2} - \mcP_{U_2}\big)\mcP_{U_{2\perp}}V_2\right\|}_{\eqref{eq: high-prob upper bound of V1tPU1p(PUhat1-PU1)U1p in tensor regression without sample splitting}}  \left\|\mcA \times_2 U_2\times_3 U_3\right\|_{\mathrm{F}} \notag\\
\lesssim & \left\|\mcA\times_2 U_2\times_3 U_3 \right\|_{\mathrm{F}}  \cdot \left[\frac{\sigma_\xi^4\oor^{1/2}}{\ulambda^3\sigma^4} \left(\frac{\op\oR\log(\op)}{n^2}+ \Delta^2\cdot \frac{\op^2}{n^2}\right)\right] \label{eq: upper bound of term 2.3 in step 1.2 in tensor regression without sample splitting}.
\end{align}

For the term $\mathrm{\RN{2}}.\mathrm{\RN{4}}$ in the fourth term \eqref{eq: term 2.4 in step 1.2 in tensor regression without sample splitting}, we have
\begin{align}
\mathrm{\RN{2}}.\mathrm{\RN{4}} 
\leq & \sup_{\substack{W_1\in \mathbb{R}^{p_1\times 2r_1}, \left\|W_1\right\|=1\\ W_2\in \mathbb{R}^{p_2\times 2r2}, \left\|W_2\right\|=1}}\left\|\mcP_{U_3}\whZ_3\left(W_2\otimes W_1\right)\right\|_{\mathrm{F}} \underbrace{ \left\|\mcP_{U_{1\perp}}\big(\mcP_{\whU_1} - \mcP_{U_1}\big)\mcP_{U_{1\perp}}V_1\right\|}_{\eqref{eq: high-prob upper bound of V1tPU1p(PUhat1-PU1)U1p in tensor regression without sample splitting}} 
  \left\|\mcP_{U_{2\perp}}\big(\mcP_{\whU_1} - \mcP_{U_1}\big)\mcP_{U_2}\right\|  \left\|\mcA\times_2 U_2\times_3 U_3\right\|_{\mathrm{F}} \notag\\
\lesssim & \left\|\mcA\times_2 U_2\times_3 U_3 \right\|_{\mathrm{F}} \cdot \frac{\sigma_\xi^4\oor^{1/2}}{\ulambda^3\sigma^4}\cdot \left[\left(\frac{\op^{3/2}\sqrt{\oR\log(\op)}}{n^2}+ \Delta\cdot \frac{\op^2}{n^2}\right)\right] \label{eq: upper bound of term 2.4 in step 1.2 in tensor regression without sample splitting}.
\end{align}

Therefore, we have the following upper bound of the second term \RN{2} \eqref{eq: term 2 in step 1.2 in tensor regression without sample splitting}
\begin{align}
\mathrm{\RN{2}} 
\lesssim & \left\|\mcA \times_2 U_2\times_3 U_3 \right\|_{\mathrm{F}} \cdot \frac{\sigma_\xi^4\oor^{1/2}}{\ulambda^3\sigma^4}\cdot \left[\left(\frac{\op^{3/2}\sqrt{\oR\log(\op)}}{n^2}+ \Delta\cdot \frac{\op^2}{n^2}\right)\right] . 
\label{eq: upper bound of term 2 in step 1.2 in tensor regression without sample splitting}
\end{align}

Similar to the proof of \RN{2}, we have the same upper bound for \RN{3} as for \RN{2}.


Finally, consider the fourth term $\mathrm{\RN{4}}$ \eqref{eq: term 4 in step 1.2 in tensor regression without sample splitting}. We have the following decomposition:
\begin{align}
\mathrm{\RN{4}} 
\leq & \underbrace{\left|\operatorname{tr}\left[\left(\mcP_{U_3}\otimes \mcP_{U_2}\left(\mcP_{\whU_2}-\mcP_{U_2}\right)\mcP_{U_{2\perp}}\right)A_1^{\top}\mcP_{U_{1\perp}}\left(\mcP_{\whU_1}-\mcP_{U_1}\right)\mcP_{U_1}\whZ_1\left(\mcP_{U_3}\otimes \mcP_{U_2}\right)\right]\right|}_{\mathrm{\RN{4}}.\mathrm{\RN{1}}} \label{eq: term 4.1 in step 1.2 in tensor regression without sample splitting}\\
+ & \underbrace{\left|\operatorname{tr}\left[\left(\mcP_{U_3}\otimes \mcP_{U_{2\perp}}\left(\mcP_{\whU_2}-\mcP_{U_2}\right)\mcP_{U_{2\perp}}\right)A_1^{\top}\mcP_{U_{1\perp}}\left(\mcP_{\whU_1}-\mcP_{U_1}\right)\mcP_{U_1}\whZ_1\left(\mcP_{U_3}\otimes \mcP_{U_{2\perp}}\right)\right]\right|}_{\mathrm{\RN{4}}.\mathrm{\RN{2}}} \label{eq: term 4.2 in step 1.2 in tensor regression without sample splitting}\\
+ & \underbrace{\left|\operatorname{tr}\left[\left(\mcP_{U_3}\otimes \mcP_{U_2}\left(\mcP_{\whU_2}-\mcP_{U_2}\right)\mcP_{U_{2\perp}}\right)A_1^{\top}\mcP_{U_{1\perp}}\left(\mcP_{\whU_1}-\mcP_{U_1}\right)\mcP_{U_{1\perp}}\whZ_1\left(\mcP_{U_3}\otimes \mcP_{U_2}\right)\right]\right|}_{\mathrm{\RN{4}}.\mathrm{\RN{3}}} \label{eq: term 4.3 in step 1.2 in tensor regression without sample splitting}\\
+ & \underbrace{\left|\operatorname{tr}\left[\left(\mcP_{U_3}\otimes \mcP_{U_{2\perp}}\left(\mcP_{\whU_2}-\mcP_{U_2}\right)\mcP_{U_{2\perp}}\right)A_1^{\top}\mcP_{U_{1\perp}}\left(\mcP_{\whU_1}-\mcP_{U_1}\right)\mcP_{U_{1\perp}}\whZ_1\left(\mcP_{U_3}\otimes \mcP_{U_{2\perp}}\right)\right]\right|}_{\mathrm{\RN{4}}.\mathrm{\RN{4}}}. \label{eq: term 4.4 in step 1.2 in tensor regression without sample splitting}
\end{align}

Here, for the first term $\mathrm{\RN{4}}.\mathrm{\RN{1}}$ \eqref{eq: term 4.1 in step 1.2 in tensor regression without sample splitting}, we have
\begin{align}
\mathrm{\RN{4}}.\mathrm{\RN{1}}
\leq & \underbrace{\left\|\left(\mcP_{U_3} \otimes \mcP_{U_2}\left(\mcP_{\whU_2}-\mcP_{U_2}\right) \mcP_{U_{2 \perp}}\right) A_1^{\top} \mcP_{U_{1 \perp}}\left(\mcP_{\whU_1}-\mcP_{U_1}\right) \mcP_{U_1}\right\|_{\mathrm{F}}}_{\eqref{eq: high-prob upper bound of AoPU1oPU2p(PUhat2-P_U2)PU2oPU3p(PUhat3-P_U3)PU3 in tensor regression without sample splitting}} \cdot \underbrace{\left\|\mcP_{U_1}  \whZ_1\left(\mcP_{U_3} \otimes \mcP_{U_2}\right)\right\|_{\mathrm{F}}}_{\eqref{eq: high-prob upper bound of U1Zhat1(U3oU2) in tensor regression without sample splitting}} \notag \\
\lesssim & \left\|\mcA \times_3 U_3 \right\|_{\mathrm{F}} \cdot \left[\frac{\sigma_{\xi}^3\oor^{1/2}}{\ulambda^2 \sigma^3} \cdot\left(\frac{\oor^{3/2}\log(\op)^{3/2}}{n^{3/2}}+\Delta \cdot \frac{\sqrt{\oor\oR\op\log(\op)}}{n^{3/2}}+\Delta^2\cdot \frac{\op\sqrt{\oR\log(\op)}}{n^{3/2}}+\Delta^3 \cdot \frac{\op^{3/2}}{n^{3/2}}\right)\right]. \label{eq: upper bound of term 4.1 in step 1.2 in tensor regression without sample splitting}
\end{align}

For the second term $\mathrm{\RN{4}}.\mathrm{\RN{2}}$ \eqref{eq: term 4.2 in step 1.2 in tensor regression without sample splitting}, we have
\begin{align}
\mathrm{\RN{4}}.\mathrm{\RN{2}}
\leq & \underbrace{\left\|\left(\mcP_{U_3} \otimes \mcP_{U_{2\perp}}\left(\mcP_{\whU_2}-\mcP_{U_2}\right) \mcP_{U_{2 \perp}}\right) A_1^{\top} \mcP_{U_{1 \perp}}\left(\mcP_{\whU_1}-\mcP_{U_1}\right) \mcP_{U_1}\right\|_{\mathrm{F}}}_{\eqref{eq: high-prob upper bound of AoPU1oPU2p(PUhat2-P_U2)PU2oPU3p(PUhat3-P_U3)PU3 in tensor regression without sample splitting}} \cdot \underbrace{\left\|\mcP_{U_1}  \whZ_1\left(\mcP_{U_3} \otimes \mcP_{U_{2\perp}}\right)\right\|_{\mathrm{F}}}_{\eqref{eq: high-prob upper bound of U1pZhat1(U3oU2) in tensor regression without sample splitting}} \notag \\
= & \left\|\mcA \times_3 U_3\right\|_{\mathrm{F}}  \cdot \left[\frac{\sigma_{\xi}^4\oor^{1/2}}{\ulambda^3 \sigma^4}\left(\frac{\oor\op\log (\op)}{n^2}+\Delta \cdot \frac{\oR^{1/2} \op^{3/2}}{n^2}+\Delta^2 \cdot \frac{\op^2}{n^2}\right)\right] .\label{eq: upper bound of term 4.2 in step 1.2 in tensor regression without sample splitting}
\end{align}

By symmetry, we have the following upper bound for the third term $\mathrm{\RN{4}}.\mathrm{\RN{3}}$ \eqref{eq: term 4.3 in step 1.2 in tensor regression without sample splitting} in $\mathrm{\RN{4}}$ \eqref{eq: term 4 in step 1.2 in tensor regression without sample splitting}:
\begin{align}
\mathrm{\RN{4}}.\mathrm{\RN{3}} 
\lesssim & \left\|\mcA \times_3 U_3\right\|_{\mathrm{F}}  \cdot \left[\frac{\sigma_{\xi}^4\oor^{1/2}}{\ulambda^3 \sigma^4}\left(\frac{\oor\op\log (\op)}{n^2}+\Delta \cdot \frac{\oR^{1/2} \op^{3/2}}{n^2}+\Delta^2 \cdot \frac{\op^2}{n^2}\right)\right] \label{eq: upper bound of term 4.3 in step 1.2 in tensor regression without sample splitting}.
\end{align}

For the fourth term $\mathrm{\RN{4}}.\mathrm{\RN{4}}$ \eqref{eq: term 4.4 in step 1.2 in tensor regression without sample splitting} in term $\mathrm{\RN{4}}$ \eqref{eq: term 4 in step 1.2 in tensor regression without sample splitting}, we have
\begin{align}
\mathrm{\RN{4}}.\mathrm{\RN{4}} 
\leq & \underbrace{\left\|\left(\mcP_{U_3}\otimes \mcP_{U_{2\perp}}\left(\mcP_{\whU_2}-\mcP_{U_2}\right)\mcP_{U_{2\perp}}\right)A_1^{\top}\mcP_{U_{1\perp}}\left(\mcP_{\whU_1}-\mcP_{U_1}\right)\mcP_{U_{1\perp}}\right\|_{\mathrm{F}}}_{\eqref{eq: high-prob upper bound of AoPU1oPU2p(PUhat2-P_U2)PU2oPU3p(PUhat3-P_U3)PU3 in tensor regression without sample splitting}} \cdot \sup_{\substack{W_1 \in \mathbb{R}^{p_1\times 2r_1}, \left\|W_1\right\|=1\\ W_2 \in \mathbb{R}^{p_2\times 2r_2}, \left\|W_2\right\|=1}} \left\|W_1^{\top}\widehat{Z}_1\left(U_3\otimes W_2\right)\right\|_{\mathrm{F}} \notag \\
\lesssim & \left\|\mcA \times_3 U_3\right\|_{\mathrm{F}} \cdot \left[\frac{\sigma_{\xi}^5\oor^{1/2}}{\ulambda^4 \sigma^5}\left(\frac{\oor\op^{3/2}\log(\op)}{n^{5/2}}+\Delta \cdot \frac{\oR^{1/2} \op^2}{n^{5/2}}+\Delta^2 \cdot \frac{\op^{5/2}}{n^{5/2}}\right)\right] \label{eq: upper bound of term 4.4 in step 1.2 in tensor regression without sample splitting}.
\end{align}

Therefore, we have the following upper bound of \eqref{eq: term 4 in step 1.2 in tensor regression without sample splitting}:
\begin{align}
\mathrm{\RN{4}}
\leq & \left\|\mcA \times_3 U_3 \right\|_{\mathrm{F}} \cdot \left[\frac{\sigma_{\xi}^3\oor^{1/2}}{\ulambda^2 \sigma^3} \cdot\left(\frac{\oor^{3/2}\log(\op)^{3/2}}{n^{3/2}}+\Delta \cdot \frac{\sqrt{\oor\oR\op\log(\op)}}{n^{3/2}}+\Delta^2\cdot \frac{\op\sqrt{\oR\log(\op)}}{n^{3/2}}+\Delta^3 \cdot \frac{\op^{3/2}}{n^{3/2}}\right)\right] \label{eq: upper bound of term 4 in step 1.2 in tensor regression without sample splitting}.
\end{align}

Then, combining all the results above, we finally have the following upper bound of \eqref{eq: step 1.2 in tensor regression without sample splitting}:
\begin{align*}
& \left|\left\langle\mcZ \times_1\left(\mcP_{\whU_1}-\mcP_{U_1}\right) \times_2\left(\mcP_{\whU_2}-\mcP_{U_2}\right) \times_3 \mcP_{U_3}, \mcA\right\rangle\right| 
\lesssim  \eqref{eq: upper bound of term 1 in step 1.2 in tensor regression without sample splitting} + \eqref{eq: upper bound of term 2 in step 1.2 in tensor regression without sample splitting} +  \eqref{eq: upper bound of term 4 in step 1.2 in tensor regression without sample splitting} \\
\lesssim & \left\|\mcA \times_1 U_1 \times_2 U_2 \times_3 U_3\right\|_{\mathrm{F}} \cdot \frac{\sigma_{\xi}^3 \oor^{1/2}}{\ulambda^2 \sigma^3} \cdot \frac{\op^{3/2}}{n^{3/2}} \\
+ & \left(\left\|\mcA \times_2 U_2 \times_3 U_3\right\|_{\mathrm{F}} + \left\|\mcA \times_1 U_1 \times_3 U_3\right\|_{\mathrm{F}}\right) \cdot \frac{\sigma_{\xi}^4 \oor^{1/2}}{\ulambda^3 \sigma^4} \cdot\left[\left(\frac{\op^{3/2} \sqrt{\oR \log (\op)}}{n^2}+\Delta \cdot \frac{\op^2}{n^2}\right)\right] \\
+ & \left\|\mcA \times_3 U_3 \right\|_{\mathrm{F}} \cdot \left[\frac{\sigma_{\xi}^3\oor^{1/2}}{\ulambda^2 \sigma^3} \cdot\left(\frac{\oor^{3/2}\log(\op)^{3/2}}{n^{3/2}}+\Delta \cdot \frac{\sqrt{\oor\oR\op\log(\op)}}{n^{3/2}}+\Delta^2\cdot \frac{\op\sqrt{\oR\log(\op)}}{n^{3/2}}+\Delta^3 \cdot \frac{\op^{3/2}}{n^{3/2}}\right)\right] .
\end{align*}

\subsubsection*{Step 1.3: Upper Bound of $\langle \widehat{\mcZ} \times_1 \left(\mcP_{\widehat{U}_1} - \mcP_{U_1}\right) \times_{2} \left(\mcP_{\widehat{U}_{2}} - \mcP_{U_2}\right) \times_{3} \left(\mcP_{\widehat{U}_{3}} - \mcP_{U_3}\right), \mcA \rangle$}

By symmetry, it suffices to consider
\begin{align}
\mathrm{\RN{1}}=& \left|\left\langle \widehat{\mcZ} \times_1 \left(\mcP_{\whU_1} - \mcP_{U_1}\right) \times_{2} \left(\mcP_{\whU_{2}} - \mcP_{U_2}\right) \times_{3} \left(\mcP_{\whU_{3}} - \mcP_{U_3}\right), \mcA \times_1 \mcP_{U_1} \times_2 \mcP_{U_2} \times_3 \mcP_{U_3} \right\rangle\right| , 
\label{eq: term 1 in step 1.3 in tensor regression without sample splitting}\\
\mathrm{\RN{2}}=& \left|\left\langle \widehat{\mcZ} \times_1 \left(\mcP_{\whU_1} - \mcP_{U_1}\right) \times_{2} \left(\mcP_{\whU_{2}} - \mcP_{U_2}\right) \times_{3} \left(\mcP_{\whU_{3}} - \mcP_{U_3}\right), \mcA \times_1 \mcP_{U_{1\perp}} \times_2 \mcP_{U_2} \times_3 \mcP_{U_3} \right\rangle\right| ,
\label{eq: term 2 in step 1.3 in tensor regression without sample splitting} \\
\mathrm{\RN{3}}=& \left|\left\langle \widehat{\mcZ} \times_1 \left(\mcP_{\whU_1} - \mcP_{U_1}\right) \times_{2} \left(\mcP_{\whU_{2}} - \mcP_{U_2}\right) \times_{3} \left(\mcP_{\whU_{3}} - \mcP_{U_3}\right), \mcA \times_1 \mcP_{U_1} \times_2 \mcP_{U_{2\perp}} \times_3 \mcP_{U_3} \right\rangle\right| ,
\label{eq: term 3 in step 1.3 in tensor regression without sample splitting}\\
\mathrm{\RN{4}}=& \left|\left\langle \widehat{\mcZ} \times_1 \left(\mcP_{\whU_1} - \mcP_{U_1}\right) \times_{2} \left(\mcP_{\whU_{2}} - \mcP_{U_2}\right) \times_{3} \left(\mcP_{\whU_{3}} - \mcP_{U_3}\right), \mcA \times_1 \mcP_{U_{1\perp}} \times_2 \mcP_{U_{2\perp}} \times_3 \mcP_{U_{3\perp}} \right\rangle\right| \label{eq: term 4 in step 1.3 in tensor regression without sample splitting}.
\end{align}

We first consider the first term 
$\mathrm{\RN{1}}$ \eqref{eq: term 1 in step 1.3 in tensor regression without sample splitting}
\begin{align}
\mathrm{\RN{1}} 
\leq & \sup_{\substack{W_j \in \mathbb{R}^{p_j\times 2r_j}, \left\|W_j\right\|=1\\ j=1,2,3}}\left\|W_1^{\top}\whZ_1\left(W_3\otimes W_2\right)\right\|_{\mathrm{F}}\cdot \prod_{j=1}^3 \left\|\mcP_{\whU_j} -\mcP_{U_j}\right\| \cdot \left\|\mcA\times U_1\times_2 U_2\times_3 U_3\right\|_{\mathrm{F}} \notag \\
\lesssim & \left\|\mcA \times_1 U_1 \times_2 U_2 \times_3 U_3\right\|_{\mathrm{F}} \cdot \frac{\sigma_\xi^4\oor^{1/2}}{\ulambda^3\sigma^4}\cdot \frac{\op^2}{n^2}. \label{eq: upper bound of term 1 in step 1.3 in tensor regression without sample splitting}
\end{align}

Then similarly, consider the second term \eqref{eq: term 2 in step 1.3 in tensor regression without sample splitting}
\begin{align}
\mathrm{\RN{2}} 
\leq & \underbrace{\left\|\left(\mcP_{U_3}\otimes \mcP_{U_2}\right)A_1^{\top}\mcP_{U_{1\perp}}\left(\mcP_{\whU_1} - \mcP_{U_1}\right)\right\|_{\mathrm{F}}}_{\eqref{eq: high-prob upper bound of V1tPU1p(PUhat1-PU1) in tensor regression without sample splitting}} \cdot \sup_{\substack{W_j \in \mathbb{R}^{p_j\times r_j}, \left\|W_j\right\|=1\\ j=1,2,3}}\left\|W_1^{\top}\whZ_1\left(W_3\otimes W_2\right)\right\|_{\mathrm{F}} \cdot \prod_{j=2}^3 \left\|\mcP_{\whU_j} -\mcP_{U_j}\right\| \notag \\
\lesssim & \left\|\mcA \times_2 U_2 \times_3 U_3\right\|_{\mathrm{F}} \cdot \left[\frac{\sigma_\xi^4\oor^{1/2}}{\ulambda^3\sigma^4}\cdot \left(\frac{\op^{3/2}\sqrt{\oor\log(\op)}}{n^2} + \Delta\cdot \frac{\op^2}{n^2}\right)\right] \label{eq: upper bound of term 2 in step 1.3 in tensor regression without sample splitting}.
\end{align}

Furthermore, for the third term $\mathrm{\RN{3}}$ \eqref{eq: term 3 in step 1.3 in tensor regression without sample splitting}, we have
\begin{align*}
\mathrm{\RN{3}} 
= & \left|\left\langle \mcZ \times_1 \left(\mcP_{\whU_1} - \mcP_{U_1}\right) \times_{2} \left(\mcP_{\whU_{2}} - \mcP_{U_{2}}\right) \times_{3} \left(\mcP_{\whU_{3}} - \mcP_{U_3}\right), \mcA \times_1 \mcP_{U_{1\perp}} \times_2 \mcP_{U_{2\perp}} \times_3 \mcP_{U_3} \right\rangle \right|\\
= & \underbrace{\left| \operatorname{tr}\left[\mcP_{U_3}\left(\mcP_{\whU_{3}} - \mcP_{U_3}\right)\mcP_{U_3}\whZ_3\left(\mcP_{U_2}\left(\mcP_{\whU_{2}} - \mcP_{U_{2}}\right)\mcP_{U_{2\perp}} \otimes \mcP_{U_1}\left(\mcP_{\whU_1} - \mcP_{U_1}\right)\mcP_{U_{1\perp}}\right)\right]A_3\mcP_{U_3}\right|}_{\mathrm{\RN{3}}.\mathrm{\RN{1}}} \\
+ & \underbrace{\left| \operatorname{tr}\left[\mcP_{U_3}\left(\mcP_{\whU_{3}} - \mcP_{U_3}\right)\mcP_{U_{3\perp}}\whZ_3\left(\mcP_{U_2}\left(\mcP_{\whU_{2}} - \mcP_{U_{2}}\right)\mcP_{U_{2\perp}} \otimes \mcP_{U_1}\left(\mcP_{\whU_1} - \mcP_{U_1}\right)\mcP_{U_{1\perp}}\right)\right]A_3\mcP_{U_3}\right|}_{\mathrm{\RN{3}}.\mathrm{\RN{2}}} \\
+ & \underbrace{\left| \operatorname{tr}\left[\mcP_{U_3}\left(\mcP_{\whU_{3}} - \mcP_{U_3}\right)\mcP_{U_3}\whZ_3\left(\mcP_{U_2}\left(\mcP_{\whU_{2}} - \mcP_{U_{2}}\right)\mcP_{U_{2\perp}} \otimes \mcP_{U_{1\perp}}\left(\mcP_{\whU_1} - \mcP_{U_1}\right)\mcP_{U_{1\perp}}\right)\right]A_3\mcP_{U_3}\right|}_{\mathrm{\RN{3}}.\mathrm{\RN{3}}} \\
+ & \underbrace{\left| \operatorname{tr}\left[\mcP_{U_3}\left(\mcP_{\whU_{3}} - \mcP_{U_3}\right)\mcP_{U_3}\whZ_3\left(\mcP_{U_{2\perp}}\left(\mcP_{\whU_{2}} - \mcP_{U_{2}}\right)\mcP_{U_{2\perp}} \otimes \mcP_{U_1}\left(\mcP_{\whU_1} - \mcP_{U_1}\right)\mcP_{U_{1\perp}}\right)\right]A_3\mcP_{U_3}\right|}_{\mathrm{\RN{3}}.\mathrm{\RN{4}}} \\
+ & \underbrace{\left| \operatorname{tr}\left[\mcP_{U_3}\left(\mcP_{\whU_{3}} - \mcP_{U_3}\right)\mcP_{U_{3\perp}}\whZ_3\left(\mcP_{U_2}\left(\mcP_{\whU_{2}} - \mcP_{U_{2}}\right)\mcP_{U_{2\perp}} \otimes \mcP_{U_{1\perp}}\left(\mcP_{\whU_1} - \mcP_{U_1}\right)\mcP_{U_{1\perp}}\right)\right]A_3\mcP_{U_3}\right|}_{\mathrm{\RN{3}}.\mathrm{\RN{5}}} \\
+ & \underbrace{\left| \operatorname{tr}\left[\mcP_{U_3}\left(\mcP_{\whU_{3}} - \mcP_{U_3}\right)\mcP_{U_{3\perp}}\whZ_3\left(\mcP_{U_{2\perp}}\left(\mcP_{\whU_{2}} - \mcP_{U_{2}}\right)\mcP_{U_{2\perp}} \otimes \mcP_{U_1}\left(\mcP_{\whU_1} - \mcP_{U_1}\right)\mcP_{U_{1\perp}}\right)\right]A_3\mcP_{U_3}\right|}_{\mathrm{\RN{3}}.\mathrm{\RN{6}}} \\
+ & \underbrace{\left| \operatorname{tr}\left[\mcP_{U_3}\left(\mcP_{\whU_{3}} - \mcP_{U_3}\right)\mcP_{U_3}\whZ_3\left(\mcP_{U_{2\perp}}\left(\mcP_{\whU_{2}} - \mcP_{U_{2}}\right)\mcP_{U_{2\perp}} \otimes \mcP_{U_{1\perp}}\left(\mcP_{\whU_1} - \mcP_{U_1}\right)\mcP_{U_{1\perp}}\right)\right]A_3\mcP_{U_3}\right|}_{\mathrm{\RN{3}}.\mathrm{\RN{7}}} \\
+ & \underbrace{\left| \operatorname{tr}\left[\mcP_{U_3}\left(\mcP_{\whU_{3}} - \mcP_{U_3}\right)\mcP_{U_{3\perp}}\whZ_3\left(\mcP_{U_{2\perp}}\left(\mcP_{\whU_{2}} - \mcP_{U_{2}}\right)\mcP_{U_{2\perp}} \otimes \mcP_{U_{1\perp}}\left(\mcP_{\whU_1} - \mcP_{U_1}\right)\mcP_{U_{1\perp}}\right)\right]A_3\mcP_{U_3}\right|}_{\mathrm{\RN{3}}.\mathrm{\RN{8}}}. 
\end{align*}

Here, we have
\begin{align}
\mathrm{\RN{3}}.\mathrm{\RN{1}} 
\leq & \underbrace{\left\|\mcP_{U_3} \whZ_3\left(\mcP_{U_2}\otimes \mcP_{U_1}\right)\right\|_{\mathrm{F}}}_{\eqref{eq: high-prob upper bound of U1Zhat1(U3oU2) in tensor regression without sample splitting}}  \underbrace{\left\|\mcP_{U_3}\left(\mcP_{\whU_3}-\mcP_{U_3}\right) \mcP_{U_3}\right\|}_{\eqref{eq: high-prob upper bound of PU1(PUhat1-PU1)PU1 in tensor regression}} 
 \underbrace{\left\|\left(\mcP_{U_2}\left(\mcP_{\whU_2}-\mcP_{U_2}\right) \mcP_{U_{2\perp}} \otimes \mcP_{U_1}\left(\mcP_{\whU_1}-\mcP_{U_1}\right) \mcP_{U_{1\perp}}\right) A_3 \mcP_{U_3}\right\|_{\mathrm{F}}}_{\eqref{eq: high-prob upper bound of AoPU1oPU2p(PUhat2-P_U2)PU2oPU3p(PUhat3-P_U3)PU3 in tensor regression without sample splitting}} \notag \\
\lesssim & \left\|\mcA\right\|_{\mathrm{F}} \cdot \left[\frac{\sigma_{\xi}^3\oor^{1/2}}{\ulambda^3 \sigma^4} \cdot\left(\frac{\oor^{3/2}\op\log(\op)^{3/2}}{n^{5/2}}+ \Delta\cdot \frac{\oor^{1/2}\oR^{1/2}\op^{3/2}\log(\op)}{n^{5/2}} + \Delta^2\cdot \frac{\op^{2}\sqrt{\oR\log(\op)}}{n^{5/2}} +\Delta^3\cdot \frac{\op^{5/2}}{n^{5/2}}\right)\right] \label{eq: upper bound of term 3.1 in step 1.3 in tensor regression without sample splitting}.
\end{align}

Then, consider
\begin{align}
\mathrm{\RN{3}}.\mathrm{\RN{2}} 
\leq & \underbrace{\left\|\mcP_{U_{3\perp}} \whZ_3 (\mcP_{U_2} \otimes \mcP_{U_1} )\right\|_{\mathrm{F}}}_{\eqref{eq: high-prob upper bound of U1Zhat1(U3oU2) in tensor regression without sample splitting}}  \underbrace{\left\|\mcP_{U_3}\big(\mcP_{\whU_3}-\mcP_{U_3}\big) \mcP_{U_3}\right\|}_{\eqref{eq: high-prob upper bound of PU1(PUhat1-PU1)PU1 in tensor regression}} 
 \underbrace{\left\|\left[\mcP_{U_2}\big(\mcP_{\whU_2}-\mcP_{U_2}\big) \mcP_{U_{2 \perp}} \otimes \mcP_{U_1}\big(\mcP_{\whU_1}-\mcP_{U_1}\big) \mcP_{U_{1 \perp}}\right] A_3 \mcP_{U_3}\right\|_{\mathrm{F}}}_{\eqref{eq: high-prob upper bound of AoPU1oPU2p(PUhat2-P_U2)PU2oPU3p(PUhat3-P_U3)PU3 in tensor regression without sample splitting}} \notag \\
\lesssim & \left\|\mcA \times_3 U_3\right\|_{\mathrm{F}} \cdot \left[\frac{\sigma_{\xi}^5\oor^{1/2}}{\ulambda^4 \sigma^5} \cdot\left(\frac{\oor\op^{3/2}\log(\op)}{n^{5/2}}+\Delta \cdot \frac{\op^2\sqrt{\oR\log(\op)}}{n^{5/2}}+\Delta^2 \cdot \frac{\op^{5/2}}{n^{5/2}}\right)\right] \label{eq: upper bound of term 3.2 in step 1.3 in tensor regression without sample splitting}.
\end{align}

Then, consider
\begin{align}
\mathrm{\RN{3}}.\mathrm{\RN{3}} 
\leq & \underbrace{\left\|\mcP_{U_3}\left(\mcP_{\whU_3}-\mcP_{U_3}\right) \mcP_{U_3}\right\|}_{\eqref{eq: high-prob upper bound of PU1(PUhat1-PU1)PU1 in tensor regression}} \cdot \sup_{W_2 \in \mathbb{R}^{p_2\times 2r_2}, \left\|W_2\right\|=1} \left\|\mcP_{U_3} \whZ_3\left(W_2 \otimes \mcP_{U_1} \right)\right\| \notag \\
& \cdot \underbrace{\left\| \left(\mcP_{U_2}\left(\mcP_{\whU_2}-\mcP_{U_2}\right) \mcP_{U_{2 \perp}} \otimes \mcP_{U_{1 \perp}}\left(\mcP_{\whU_1}-\mcP_{U_1}\right) \mcP_{U_{1 \perp}}\right) A_3 \mcP_{U_3}\right\|_{\mathrm{F}}}_{\eqref{eq: high-prob upper bound of AoPU1oPU2p(PUhat2-P_U2)PU2poPU3p(PUhat3-P_U3)PU3 in tensor regression without sample splitting}} \notag \\
\lesssim & \left\|\mcA \times_3 U_3\right\|_{\mathrm{F}} \cdot \left[\frac{\sigma_{\xi}^6\oor^{1/2}}{\ulambda^5 \sigma^6}\left(\frac{\oor\op^2\log (\op)}{n^3}+ \Delta \cdot \frac{\op^{5/2}\sqrt{\oR\log(\op)}}{n^3} +\Delta^2 \cdot \frac{\op^3}{n^3}\right)\right] \label{eq: upper bound of term 3.3 in step 1.3 in tensor regression without sample splitting}.
\end{align}

By symmetry, \RN{3}.\RN{4} has the same upper bound as \RN{3}.\RN{3}. 

Then, consider
\begin{align}
\mathrm{\RN{3}}.\mathrm{\RN{5}} 
\leq & \underbrace{\left\|\mcP_{U_{3 \perp}} \whZ_3 (\mcP_{U_2} \otimes \mcP_{U_1} )\right\|_{\mathrm{F}}}_{\eqref{eq: high-prob upper bound of U1pZhat1(U3oU2) in tensor regression without sample splitting}}  \left\|\mcP_{U_3}\big(\mcP_{\whU_3}-\mcP_{U_3}\big) \mcP_{U_{3 \perp}}\right\| 
  \underbrace{\left\|\left(\mcP_{U_2}\big(\mcP_{\whU_2}-\mcP_{U_2}\big) \mcP_{U_{2 \perp}} \otimes \mcP_{U_{1 \perp}}\big(\mcP_{\whU_1}-\mcP_{U_1}\big) \mcP_{U_{1 \perp}}\right) A_3 \mcP_{U_3}\right\|_{\mathrm{F}}}_{\eqref{eq: high-prob upper bound of AoPU1oPU2p(PUhat2-P_U2)PU2poPU3p(PUhat3-P_U3)PU3 in tensor regression without sample splitting}} \notag\\
\lesssim & \left\|\mcA \times_3 U_3\right\|_{\mathrm{F}} \cdot \left[\frac{\sigma_{\xi}^5\oor^{1/2}}{\ulambda^4 \sigma^5}\left(\frac{\oor\op^{3/2}\log(\op)}{n^{5/2}}+\Delta \cdot \frac{\op^{3/2}\sqrt{\oR\log(\op)}}{n^{5/2}}+\Delta^2 \cdot \frac{\op^{5/2}}{n^{5/2}}\right)\right]. \label{eq: upper bound of term 3.5 in step 1.3 in tensor regression without sample splitting}
\end{align}

By symmetry, \RN{3}.\RN{6} has the same upper bound as \RN{3}.\RN{5}. 

Furthermore, we have
\begin{align}
\mathrm{\RN{3}}.\mathrm{\RN{7}} 
\leq & \sup_{\substack{W_1\in \mathbb{R}^{p_1\times 2r_1}, \left\|W_1\right\|=1\\W_2\in \mathbb{R}^{p_2\times 2r_2}, \left\|W_2\right\|=1}}\left\|\mcP_{U_{3 \perp}} \whZ_3\left(W_2 \otimes W_1\right)\right\|_{\mathrm{F}} \cdot \underbrace{\left\|\mcP_{U_3}\left(\mcP_{\whU_3}-\mcP_{U_3}\right) \mcP_{U_3}\right\|}_{\eqref{eq: high-prob upper bound of PU1(PUhat1-PU1)PU1 in tensor regression}} \notag \\
& \cdot \underbrace{\left\|\left(\mcP_{U_{2 \perp}}\left(\mcP_{\whU_2}-\mcP_{U_2}\right) \mcP_{U_{2 \perp}} \otimes \mcP_{U_{1\perp}}\left(\mcP_{\whU_1}-\mcP_{U_1}\right) \mcP_{U_{1 \perp}}\right) A_3 \mcP_{U_3}\right\|_{\mathrm{F}}}_{\eqref{eq: high-prob upper bound of AoPU1oPU2p(PUhat2-P_U2)PU2poPU3p(PUhat3-P_U3)PU3p in tensor regression without sample splitting}} \notag \\
\lesssim & \left\|\mcA \times_3 U_3\right\|_{\mathrm{F}} \cdot \left[\frac{\sigma_{\xi}^7\oor^{1/2}}{\ulambda^6 \sigma^7}\left(\frac{\oor\op^{5/2}\log(\op)}{n^{7/2}}+\Delta\cdot \frac{\op^3\sqrt{\oR\log(\op)}}{n^{7/2}}+\Delta^2 \cdot \frac{\op^{7/2}}{n^{7/2}}\right)\right]. \label{eq: upper bound of term 3.7 in step 1.3 in tensor regression without sample splitting}
\end{align}

Finally, consider
\begin{align}
\mathrm{\RN{3}}.\mathrm{\RN{8}} 
\leq & \left\|\mcP_{U_3}\left(\mcP_{\whU_3}-\mcP_{U_3}\right) \mcP_{U_{3 \perp}}\right\|_{\mathrm{F}} \cdot \sup_{\substack{W_1 \in \mathbb{R}^{p_1\times 2r_1}, \left\|W_1\right\|=1\\W_2 \in \mathbb{R}^{p_2\times 2r_2}, \left\|W_2\right\|=1}} \left\|P_{U_{3 \perp}} \whZ_3\left(W_2\otimes W_1\right)\right\| \notag \\
& \cdot \underbrace{\left\|\left(\mcP_{U_{2 \perp}}\left(\mcP_{\whU_2}-\mcP_{U_2}\right) \mcP_{U_{2 \perp}} \otimes \mcP_{U_{1 \perp}}\left(\mcP_{\whU_1}-\mcP_{U_1}\right) \mcP_{U_{1 \perp}}\right) A_3 \mcP_{U_3}\right\|_{\mathrm{F}}}_{\eqref{eq: high-prob upper bound of AoPU1oPU2p(PUhat2-P_U2)PU2poPU3p(PUhat3-P_U3)PU3 in tensor regression without sample splitting}} \notag \\
\lesssim & \left\|\mcA \times_3 U_3\right\|_{\mathrm{F}} \cdot \left[\frac{\sigma_{\xi}^6\oor^{1/2}}{\ulambda^5 \sigma^6}\left(\frac{\oor\op^2\log (\op)}{n^3}+\Delta\cdot \frac{\op^{5/2}\sqrt{\oR\log(\op)}}{n^3}+\Delta^2 \cdot \frac{\op^3}{n^3}\right)\right] \label{eq: upper bound of term 3.8 in step 1.3 in tensor regression without sample splitting}.
\end{align}

Combining the results above, we have
\begin{align}
\mathrm{\RN{3}}
\lesssim & \left\|\mcA\right\|_{\mathrm{F}} \cdot \left[\frac{\sigma_{\xi}^3\oor^{1/2}}{\ulambda^3 \sigma^4} \cdot\left(\frac{\oor^{3/2}\op\log(\op)^{3/2}}{n^{5/2}}+ \Delta\cdot \frac{\oor^{1/2}\oR^{1/2}\op^{3/2}\log(\op)}{n^{5/2}} + \Delta^2\cdot \frac{\op^{2}\sqrt{\oR\log(\op)}}{n^{5/2}} +\Delta^3\cdot \frac{\op^{5/2}}{n^{5/2}}\right)\right] \notag \\
+ & \left\|\mcA \times_3 U_3\right\|_{\mathrm{F}} \cdot\left[\frac{\sigma_{\xi}^5 \oor^{1 / 2}}{\ulambda^4 \sigma^5} \cdot\left(\frac{\oor\op^{3/2}\log(\op)}{n^{5 / 2}}+\Delta \cdot \frac{\op^2 \sqrt{\oR \log (\op)}}{n^{5 / 2}}+\Delta^2 \cdot \frac{\op^{5 / 2}}{n^{5 / 2}}\right)\right] . \label{eq: upper bound of term 3 in step 1.3 in tensor regression without sample splitting}
\end{align}

Finally, for the fourth term $\mathrm{\RN{4}}$ \eqref{eq: term 4 in step 1.3 in tensor regression without sample splitting}, we have
\begin{align*}
\mathrm{\RN{4}} 
\leq & \underbrace{\left|\operatorname{tr}\left[A_1^{\top}\mcP_{U_{1\perp}}\left(\mcP_{\whU_1}-\mcP_{U_1}\right)\mcP_{U_1}\whZ_1\left[\mcP_{U_3}\left(\mcP_{\whU_3}-\mcP_{U_3}\right)\mcP_{U_{3\perp}} \otimes \mcP_{U_2}\left(\mcP_{\whU_2}-\mcP_{U_2}\right)\mcP_{U_{2\perp}}\right]\right]\right|}_{\mathrm{\RN{4}}.\mathrm{\RN{1}}} \\
+ & \underbrace{\left|\operatorname{tr}\left[A_1^{\top}\mcP_{U_{1\perp}}\left(\mcP_{\whU_1}-\mcP_{U_1}\right)\mcP_{U_1}\whZ_1\left[\mcP_{U_3}\left(\mcP_{\whU_3}-\mcP_{U_3}\right)\mcP_{U_{3\perp}} \otimes \mcP_{U_{2\perp}}\left(\mcP_{\whU_2}-\mcP_{U_2}\right)\mcP_{U_{2\perp}}\right]\right]\right|}_{\mathrm{\RN{4}}.\mathrm{\RN{2}}} \\
+ & \underbrace{\left|\operatorname{tr}\left[A_1^{\top}\mcP_{U_{1\perp}}\left(\mcP_{\whU_1}-\mcP_{U_1}\right)\mcP_{U_1}\whZ_1\left[\mcP_{U_{3\perp}}\left(\mcP_{\whU_3}-\mcP_{U_3}\right)\mcP_{U_{3\perp}}\otimes \mcP_{U_2}\left(\mcP_{\whU_2}-\mcP_{U_2}\right)\mcP_{U_{2\perp}}\right]\right]\right|}_{\mathrm{\RN{4}}.\mathrm{\RN{3}}} \\
+ & \underbrace{\left|\operatorname{tr}\left[A_1^{\top}\mcP_{U_{1\perp}}\left(\mcP_{\whU_1}-\mcP_{U_1}\right)\mcP_{U_1}\whZ_1\left[\mcP_{U_{3\perp}}\left(\mcP_{\whU_3}-\mcP_{U_3}\right)\mcP_{U_{3\perp}}\otimes \mcP_{U_{2\perp}}\left(\mcP_{\whU_2}-\mcP_{U_2}\right)\mcP_{U_{2\perp}}\right]\right]\right|}_{\mathrm{\RN{4}}.\mathrm{\RN{4}}} \\
+ & \underbrace{\left|\operatorname{tr}\left[A_1^{\top}\mcP_{U_{1\perp}}\left(\mcP_{\whU_1}-\mcP_{U_1}\right)\mcP_{U_{1\perp}}\whZ_1\left[\mcP_{U_3}\left(\mcP_{\whU_3}-\mcP_{U_3}\right)\mcP_{U_{3\perp}}\otimes \mcP_{U_2}\left(\mcP_{\whU_2}-\mcP_{U_2}\right)\mcP_{U_{2\perp}}\right]\right]\right|}_{\mathrm{\RN{4}}.\mathrm{\RN{5}}} \\
+ & \underbrace{\left|\operatorname{tr}\left[A_1^{\top}\mcP_{U_{1\perp}}\left(\mcP_{\whU_1}-\mcP_{U_1}\right)\mcP_{U_{1\perp}}\whZ_1\left[\mcP_{U_3}\left(\mcP_{\whU_3}-\mcP_{U_3}\right)\mcP_{U_{3\perp}}\otimes \mcP_{U_{2\perp}}\left(\mcP_{\whU_2}-\mcP_{U_2}\right)\mcP_{U_{2\perp}}\right]\right]\right|}_{\mathrm{\RN{4}}.\mathrm{\RN{6}}} \\
+ & \underbrace{\left|\operatorname{tr}\left[A_1^{\top}\mcP_{U_{1\perp}}\left(\mcP_{\whU_1}-\mcP_{U_1}\right)\mcP_{U_{1\perp}}\whZ_1\left[\mcP_{U_{3\perp}}\left(\mcP_{\whU_3}-\mcP_{U_3}\right)\mcP_{U_{3\perp}}\otimes \mcP_{U_2}\left(\mcP_{\whU_2}-\mcP_{U_2}\right)\mcP_{U_{2\perp}}\right]\right]\right|}_{\mathrm{\RN{4}}.\mathrm{\RN{7}}} \\
+ & \underbrace{\left|\operatorname{tr}\left[A_1^{\top}\mcP_{U_{1\perp}}\left(\mcP_{\whU_1}-\mcP_{U_1}\right)\mcP_{U_{1\perp}}\whZ_1\left[\mcP_{U_{3\perp}}\left(\mcP_{\whU_3}-\mcP_{U_3}\right)\mcP_{U_{3\perp}}\otimes \mcP_{U_{2\perp}}\left(\mcP_{\whU_2}-\mcP_{U_2}\right)\mcP_{U_{2\perp}}\right]\right]\right|}_{\mathrm{\RN{4}}.\mathrm{\RN{8}}}.
\end{align*}

First, consider
\begin{align}
& \mathrm{\RN{4}}.\mathrm{\RN{1}} \notag \\
\leq & \underbrace{\left\|\mcP_{U_1} \whZ_1\left(\mcP_{U_3} \otimes \mcP_{U_2}\right)\right\|_{\mathrm{F}}}_{\eqref{eq: high-prob upper bound of U1Zhat1(U3oU2) in tensor regression without sample splitting}} \cdot \underbrace{\left\|\left[\mcP_{U_3}\left(\mcP_{\whU_3}-\mcP_{U_3}\right) \mcP_{U_{3 \perp}} \otimes \mcP_{U_2}\left(\mcP_{\whU_2}-\mcP_{U_2}\right) \mcP_{U_{2 \perp}}\right]A_1^{\top} \mcP_{U_{1 \perp}}\left(\mcP_{\whU_1}-\mcP_{U_1}\right) \mcP_{U_1}\right\|_{\mathrm{F}}}_{\eqref{eq: high-prob upper bound of AoPU1p(PUhat1-P_U1)PU1oPU2p(PUhat2-P_U2)PU2oPU3p(PUhat3-P_U3)PU3 in tensor regression without sample splitting}} \notag \\
\lesssim & \left\|\mcA\right\|_{\mathrm{F}}  \left[\frac{\sigma_{\xi}^4\oor^{1/2}}{\ulambda^3 \sigma^4} \left(\frac{\oor^2\log(\op)^2}{n^2}+  \frac{\Delta\oor^{1/2}\op^{1/2} \oR \log (\op)^{3/2}}{n^2} +  \frac{\Delta^2\op\oR\log(\op)}{n^2}+\Delta^3 \cdot \frac{\op^{3/2}\oR\log (\op)}{n^2} + \Delta^4\cdot \frac{\op^2}{n^2}\right)\right]. \label{eq: upper bound of term 4.1 in step 1.3 in tensor regression without sample splitting}
\end{align}

Then, consider
\begin{align}
& \mathrm{\RN{4}}.\mathrm{\RN{2}} \notag \\
\leq & \underbrace{\left\|\mcP_{U_{2\perp}} \whZ_2\left(\mcP_{U_1} \otimes \mcP_{U_3}\right)\right\|_{\mathrm{F}}}_{\eqref{eq: high-prob upper bound of U1pZhat1(U3oU2) in tensor regression without sample splitting}} \cdot \underbrace{\left\|\left[\mcP_{U_3}\left(\mcP_{\whU_3}-\mcP_{U_3}\right)\mcP_{U_{3\perp}} \otimes \mcP_{U_{2\perp}}\left(\mcP_{\whU_2}-\mcP_{U_2}\right)\mcP_{U_{2\perp}}\right]A_1^{\top} \mcP_{U_{1 \perp}}\left(\mcP_{\whU_1}-\mcP_{U_1}\right) \mcP_{U_1}\right\|_{\mathrm{F}}}_{\eqref{eq: high-prob upper bound of AoPU1p(PUhat1-P_U1)PU1poPU2p(PUhat2-P_U2)PU2oPU3p(PUhat3-P_U3)PU3 in tensor regression without sample splitting}} \notag \\
= & \left\|\mcA\right\|_{\mathrm{F}} \cdot \left[\frac{\sigma_\xi^5\oor^{1/2}}{\ulambda^4\sigma^5}\cdot \left(\frac{\oor^{3/2}\op\log(\op)^{3/2}}{n^{5/2}} + \Delta\cdot \frac{\op^{3/2}\oR\log(\op)}{n^{5/2}} + \Delta^2\cdot \frac{\op^2\oR\log(\op)}{n^{5/2}} + \Delta^3\cdot \frac{\op^{5/2}}{n^{5/2}}\right)\right]. \label{eq: upper bound of term 4.2 in step 1.3 in tensor regression without sample splitting}
\end{align}
Similar to the proof of \RN{4}.\RN{2}, we have the same upper bound for \RN{4}.\RN{3} and \RN{4}.\RN{5} as for \RN{4}.\RN{2}.


Furthermore, we have
\begin{align}
\mathrm{\RN{4}}.\mathrm{\RN{4}}
\leq & \sup _{\substack{W_3 \in \mathbb{R}_3^{p_3 \times 2 r_3},\left\|W_3\right\|=1 \\ W_2 \in \mathbb{R}^{p_2 \times 2 r_2},\left\|W_2\right\|=1}}\left\|\mcP_{U_1} \whZ_1\left(W_3 \otimes W_2\right)\right\|_{\mathrm{F}} \notag \\
& \cdot \underbrace{\left\|\left[\mcP_{U_{3\perp}}\left(\mcP_{\whU_3}-\mcP_{U_3}\right) \mcP_{U_{3 \perp}} \otimes \mcP_{U_{2 \perp}}\left(\mcP_{\whU_2}-\mcP_{U_2}\right) \mcP_{U_{2 \perp}}\right] A_1^{\top} \mcP_{U_{1 \perp}}\left(\mcP_{\whU_1}-\mcP_{U_1}\right) \mcP_{U_1}\right\|_{\mathrm{F}}}_{\eqref{eq: high-prob upper bound of AoPU1p(PUhat1-P_U1)PU1oPU2p(PUhat2-P_U2)PU2poPU3p(PUhat3-P_U3)PU3pp in tensor regression without sample splitting}} \notag \\
\lesssim & \left\|\mcA\right\|_{\mathrm{F}} \cdot \left[\frac{\sigma_{\xi}^6\oor^{1/2}}{\ulambda^5 \sigma^6} \cdot\left(\frac{\oor^{3/2}\op^{3/2}\log (\op)^{3/2}}{n^3}+\Delta \cdot \frac{\op^2 \oR \log (\op)}{n^3}+\Delta^2 \cdot \frac{\op^{5/2} \oR \log (\op)}{n^3}+\Delta^3 \cdot \frac{\op^3}{n^3}\right)\right] \label{eq: upper bound of term 4.4 in step 1.3 in tensor regression without sample splitting}.
\end{align}
Similar to the proof of \RN{4}.\RN{4}, we have the same upper bound for \RN{4}.\RN{6} and \RN{4}.\RN{7} as for \RN{4}.\RN{4}.


Finally, consider
\begin{align}
\mathrm{\RN{4}}.\mathrm{\RN{8}} 
\leq & \sup _{\substack{W_3 \in \mathbb{R}_3^{p_3 \times 2 r_3},\left\|W_3\right\|=1 \\ W_2 \in \mathbb{R}^{p_2 \times 2 r_2},\left\|W_2\right\|=1}}\left\|W_1^{\top} \whZ_1\left(W_3 \otimes W_2\right)\right\|_{\mathrm{F}} \notag \\
& \cdot \underbrace{\left\|\left[\mcP_{U_{3\perp}}\left(\mcP_{\whU_3}-\mcP_{U_3}\right) \mcP_{U_{3 \perp}} \otimes \mcP_{U_{2 \perp}}\left(\mcP_{\whU_2}-\mcP_{U_2}\right) \mcP_{U_{2 \perp}}\right] A_1^{\top} \mcP_{U_{1 \perp}}\left(\mcP_{\whU_1}-\mcP_{U_1}\right) \mcP_{U_{1\perp}}\right\|_{\mathrm{F}}}_{\eqref{eq: high-prob upper bound of AoPU1p(PUhat1-P_U1)PU1oPU2p(PUhat2-P_U2)PU2poPU3p(PUhat3-P_U3)PU3pp in tensor regression without sample splitting}} \notag \\
\lesssim & \left\|\mcA\right\|_{\mathrm{F}} \cdot \left[\frac{\sigma_\xi^7\oor^{1/2} }{\ulambda^6\sigma^7}\cdot \left(\frac{\oor^{3/2}\op^2\log (\op)^{3/2}}{n^{7/2}} + \Delta\cdot \frac{\op^{5/2}\oR\log(\op)}{n^{7/2}} + \Delta^2 \cdot \frac{\op^3\oR\log(\op)}{n^{7/2}} + \Delta^3\cdot \frac{\op^{7/2}}{n^{7/2}}\right)\right] \label{eq: upper bound of term 4.8 in step 1.3 in tensor regression without sample splitting}.
\end{align}

Therefore, we have
\begin{align}
\mathrm{\RN{4}}
\lesssim & \left\|\mcA\right\|_{\mathrm{F}} \left[\frac{\sigma_{\xi}^4\oor^{1/2}}{\ulambda^3 \sigma^4} \left(\frac{\oor^2\log(\op)^2}{n^2}+  \frac{\Delta\oor^{1/2}\op^{1/2} \oR \log (\op)^{3/2}}{n^2} +  \frac{\Delta^2\op\oR\log(\op)}{n^2}+  \frac{\Delta^3\op^{3/2}\oR\log (\op)}{n^2} +  \frac{\Delta^4\op^2}{n^2}\right)\right] \label{eq: upper bound of term 4 in step 1.3 in tensor regression without sample splitting}.
\end{align}

Therefore, by symmetry, we have
\begin{align}
& \left|\left\langle \widehat{\mcZ} \times_1 \left(\mcP_{\whU_1} - \mcP_{U_1}\right) \times_{2} \left(\mcP_{\whU_{2}} - \mcP_{U_{2}}\right) \times_{3} \left(\mcP_{\whU_{3}} - \mcP_{U_3}\right), \mcA \right\rangle\right| 
\lesssim  \eqref{eq: upper bound of term 1 in step 1.3 in tensor regression without sample splitting} + \eqref{eq: upper bound of term 2 in step 1.3 in tensor regression without sample splitting} + \eqref{eq: upper bound of term 3 in step 1.3 in tensor regression without sample splitting} + \eqref{eq: upper bound of term 4 in step 1.3 in tensor regression without sample splitting} \notag \\
\lesssim & \left\|\mcA \times_1 U_1 \times_2 U_2 \times_3 U_3\right\|_{\mathrm{F}} \cdot \frac{\sigma_{\xi}^4 \oor^{1 / 2}}{\ulambda^3 \sigma^4} \cdot \frac{\op^2}{n^2} 
+  \sum_{j=1}^3 \left\|\mcA \times_{j+1} U_{j+1} \times_{j+2} U_{j+2}\right\|_{\mathrm{F}}  \left[\frac{\sigma_{\xi}^4 \oor^{1 / 2}}{\ulambda^3 \sigma^4} \left(\frac{\oor \op^{3 / 2} \sqrt{\log (\op)}}{n^2}+\Delta \cdot \frac{\op^2}{n^2}\right)\right] \notag \\
+ & \sum_{j=1}^3 \left\|\mcA \times_j U_j\right\|_{\mathrm{F}} \cdot \left[\frac{\sigma_{\xi}^5 \oor^{1 / 2}}{\ulambda^4 \sigma^5} \cdot\left(\frac{\oor\op^{3/2}\log(\op)}{n^{5 / 2}}+\Delta \cdot \frac{\op^2 \sqrt{\oR \log (\op)}}{n^{5 / 2}}+\Delta^2 \cdot \frac{\op^{5 / 2}}{n^{5 / 2}}\right)\right] \notag \\
+ & \left\|\mcA\right\|_{\mathrm{F}} \cdot \left[\frac{\sigma_{\xi}^4\oor^{1/2}}{\ulambda^3 \sigma^4} \cdot\left(\frac{\oor^2\log(\op)^2}{n^2}+\Delta \cdot \frac{\oor^{1/2}\op^{1/2} \oR \log (\op)^{3/2}}{n^2} + \Delta^2\cdot \frac{\op\oR\log(\op)}{n^2}+\Delta^3 \cdot \frac{\op^{3/2}\oR\log (\op)}{n^2} + \Delta^4\cdot \frac{\op^2}{n^2}\right)\right] \notag .
\end{align}

\subsection*{Step 2: Upper Bound of Negligible Terms in $\langle {\mcT} \times_1 \mcP_{\widehat{U}_{1}} \times_2 \mcP_{\widehat{U}_{2}} \times_3 \mcP_{\widehat{U}_{3}}-\mcT, \mcA \rangle$} 

By symmetry, it remains to consider
\begin{align}
\text{(Step 2.1)}: & \left\langle \mcT \times_1 \left(\mcP_{\whU_1} - \mcP_{U_1}\right) \times_2 \mcP_{U_2} \times_3 \mcP_{U_3}, \mcA \right\rangle , \label{eq: step 2.1 in tensor regression without sample splitting}\\
\text{(Step 2.2)}: & \left\langle \mcT \times_1 \left(\mcP_{\whU_1} - \mcP_{U_1}\right) \times_2 \left(\mcP_{\whU_2} - \mcP_{U_2}\right) \times_3 \mcP_{U_3}, \mcA \right\rangle , \label{eq: step 2.2 in tensor regression without sample splitting}\\
\text{(Step 2.3)}: & \left\langle \mcT \times_1 \left(\mcP_{\whU_1} - \mcP_{U_1}\right) \times_2 \left(\mcP_{\whU_2} - \mcP_{U_2}\right) \times_3 \left(\mcP_{\whU_3} - \mcP_{U_3}\right), \mcA \right\rangle . \label{eq: step 2.3 in tensor regression without sample splitting}
\end{align}

\subsubsection*{Step 2.1: Upper Bound of Negligible Terms in $\langle \mcT \times_1 \left(\mcP_{\widehat{U}_1} - \mcP_{U_1}\right) \times_2 \mcP_{U_2} \times_3 \mcP_{U_3}, \mcA \rangle$} 

Note that
\begin{align}
\left\langle \mcT \times_1 \left(\mcP_{\whU_1} - \mcP_{U_1}\right) \times_2 \mcP_{U_2} \times_3 \mcP_{U_3}, \mcA \right\rangle 
= &  \underbrace{\left\langle \mcT \times_1 S_{G_1,1}\left(\whE_1\right) \times_2 \mcP_{U_2} \times_3 \mcP_{U_3}, \mcA \right\rangle}_{\mathrm{\RN{1}}} \label{eq: term 1 in step 2.1 in tensor regression without sample splitting}\\
& +  \underbrace{\left\langle \mcT \times_1 \sum_{k_1=2}^{+\infty}S_{G_1,k_1}\left(\whE_1\right) \times_2 \mcP_{U_2} \times_3 \mcP_{U_3}, \mcA \right\rangle}_{\mathrm{\RN{2}}}. \label{eq: term 2 in step 2.1 in tensor regression without sample splitting}
\end{align}
We first consider \RN{1} \eqref{eq: term 1 in step 2.1 in tensor regression without sample splitting}:
\begin{align*}
\mathrm{\RN{1}}
= & \underbrace{\left\langle \mcP_{U_{1\perp}} \whZ_1\left(\mcP_{U_3} \otimes \mcP_{U_2}\right) \mcP_{\left(U_3 \otimes U_2\right)G_1^{\top}}, A_1\right\rangle}_{\mathrm{\RN{1}}.\mathrm{\RN{1}}, \text{asymptotically normal}} \\
+ & \underbrace{\left\langle \mcP_{U_{1\perp}} \whZ_1\left[\mcP_{U_3} \otimes \left(\mcP_{\whU_2^{(1)}}-\mcP_{U_2}\right)\right] \mcP_{\left(U_3 \otimes U_2\right)G_1^{\top}}, A_1\right\rangle}_{\mathrm{\RN{1}}.\mathrm{\RN{2}}, \text{negligible}} \\
+ & \underbrace{\left\langle \mcP_{U_{1\perp}} \whZ_1\left[\left(\mcP_{\whU_3^{(1)}}-\mcP_{U_3}\right) \otimes \mcP_{U_2}\right] \mcP_{\left(U_3 \otimes U_2\right)G_1^{\top}} , A_1\right\rangle}_{\mathrm{\RN{1}}.\mathrm{\RN{3}}, \text{negligible}} \\
+ & \underbrace{\left\langle \mcP_{U_{1\perp}} \whZ_1\left[\left(\mcP_{\whU_3^{(1)}}-\mcP_{U_3}\right) \otimes \left(\mcP_{\whU_2^{(1)}}-\mcP_{U_2}\right)\right] \mcP_{\left(U_3 \otimes U_2\right)G_1^{\top}}, A_1\right\rangle}_{\mathrm{\RN{1}}.\mathrm{\RN{4}}, \text{negligible}}\\
+ & \underbrace{\left\langle \mcP_{U_{1\perp}} \whZ_1\left(\mcP_{\whU_3^{(1)}} \otimes \mcP_{\whU_2^{(1)}}\right) \whZ_1^{\top} U_1 \left(G_1G_1^{\top}\right)^{-1}  G_1 \left(U_3 \otimes U_2\right)^{\top} , A_1\right\rangle}_{\mathrm{\RN{1}}.\mathrm{\RN{5}}, \text{negligible}}.
\end{align*}
Intuitively, note that $\left\|\mcP_{\whU_j^{(1)}}- \mcP_{U_j}\right\|=\left\|\whU_j^{(1)}\whU_j^{(1)\top}- U_j U_j^{\top}\right\|$ should be sufficiently small with good initialization. We leave the proof of the asymptotic normality $\left\langle \mcP_{U_{1\perp}} \whZ_1 \mcP_{\left(U_3 \otimes U_2\right)G_1^{\top}} , A_1\right\rangle$ to Step 3. We then focus on finding upper bound of negligible terms above.

First, consider the upper bound for $\mathrm{\RN{1}}.\mathrm{\RN{2}} 
= \left\langle \mcP_{U_{1\perp}} \whZ_1\left[\mcP_{U_3} \otimes \left(\mcP_{\whU_2^{(1)}}-\mcP_{U_2}\right)\right] \mcP_{\left(U_3 \otimes U_2\right)G_1^{\top}}, A_1\right\rangle$.
\begin{align*}
\mathrm{\RN{1}}.\mathrm{\RN{2}}
= & \underbrace{\left|\operatorname{tr}\left[\Mat_2\left(\widetilde{A}_1\right)^{\top}\mcP_{U_2}S_{G_2,1}\left(\whE_2^{(0)}\right)\whZ_2\left(\mcP_{U_{1\perp}}\otimes \mcP_{U_3}\right)\right]\right|}_{\mathrm{\RN{1}}.\mathrm{\RN{2}}.\mathrm{\RN{1}}} 
+ & \underbrace{\left|\operatorname{tr}\left[\Mat_2\left(\widetilde{A}_1\right)^{\top}\mcP_{U_2}\sum_{k_2=1}^{+\infty}S_{G_2,k_2}\left(\whE_2^{(0)}\right)\whZ_2\left(\mcP_{U_{1\perp}}\otimes \mcP_{U_3}\right)\right]\right|}_{\mathrm{\RN{1}}.\mathrm{\RN{2}}.\mathrm{\RN{2}}}.
\end{align*}
where $\widetilde{A}_1=  \mcP_{U_{1\perp}}A_1 \mcP_{\left(U_3\otimes U_2\right)G_1^{\top}}$. It then follows that
\begin{align*}
\mathrm{\RN{1}}.\mathrm{\RN{2}}.\mathrm{\RN{1}} 
\leq & \underbrace{\left|\operatorname{tr}\left[\left(\mcP_{U_{1\perp}}\otimes \mcP_{U_3}\right)\Mat_2\left(\widetilde{A}_1\right)^{\top}U_2\left(G_2G_2^{\top}\right)^{-1}G_2\left(U_1\otimes U_3\right)^{\top}\left(\mcP_{\whU_1^{(0)}}\otimes \mcP_{\whU_3^{(0)}}\right)\whZ_2^{\top}\mcP_{U_{2\perp}}\whZ_2\right]\right|}_{\mathrm{\RN{1}}.\mathrm{\RN{2}}.\mathrm{\RN{1}}.\mathrm{\RN{1}}} \\
+ &  \underbrace{\left|\operatorname{tr}\left[\left(\mcP_{U_{1\perp}}\otimes \mcP_{U_3}\right)\Mat_2\left(\widetilde{A}_1\right)^{\top}U_2\left(G_2G_2^{\top}\right)^{-1}U_2^{\top}\whZ_2\left(\mcP_{\whU_1^{(0)}}\otimes \mcP_{\whU_3^{(0)}}\right)\whZ_2^{\top}\mcP_{U_{2\perp}}\whZ_2\right]\right|}_{\mathrm{\RN{1}}.\mathrm{\RN{2}}.\mathrm{\RN{1}}.\mathrm{\RN{2}}}.
\end{align*}
For $\mathrm{\RN{1}}.\mathrm{\RN{2}}.\mathrm{\RN{1}}.\mathrm{\RN{1}}$, we have
\begin{align}
& \mathrm{\RN{1}}.\mathrm{\RN{2}}.\mathrm{\RN{1}}.\mathrm{\RN{1}} \notag \\
\leq & \underbrace{\left| \operatorname{tr}\left[\left(\mcP_{U_{1\perp}}\otimes \mcP_{U_3}\right)\Mat_2\left(\widetilde{A}_1\right)^{\top}U_2\left(G_2G_2^{\top}\right)^{-1}G_2\left(U_1\otimes U_3\right)^{\top}\whZ_2^{(1)\top}\mcP_{U_{2\perp}}\whZ_2^{(1)}\right] \right|}_{\mathrm{\RN{1}}.\mathrm{\RN{2}}.\mathrm{\RN{1}}.\mathrm{\RN{1}}.\mathrm{\RN{1}}} \label{eq: term 1.2.1.1.1 in step 2.1 in tensor regression without sample splitting} \\
+ & \underbrace{\left| \operatorname{tr}\left[\left(\mcP_{U_{1\perp}}\otimes \mcP_{U_3}\right)\Mat_2\left(\widetilde{A}_1\right)^{\top}U_2\left(G_2G_2^{\top}\right)^{-1}G_2\left(U_1\otimes U_3\right)^{\top}\whZ_2^{(1)\top}\mcP_{U_{2\perp}}\whZ_2^{(2)}\right] \right|}_{\mathrm{\RN{1}}.\mathrm{\RN{2}}.\mathrm{\RN{1}}.\mathrm{\RN{1}}.\mathrm{\RN{2}}} \label{eq: term 1.2.1.1.2 in step 2.1 in tensor regression without sample splitting} \\
+ & \underbrace{\left| \operatorname{tr}\left[\left(\mcP_{U_{1\perp}}\otimes \mcP_{U_3}\right)\Mat_2\left(\widetilde{A}_1\right)^{\top}U_2\left(G_2G_2^{\top}\right)^{-1}G_2\left(U_1\otimes U_3\right)^{\top}\whZ_2^{(2)\top}\mcP_{U_{2\perp}}\whZ_2^{(1)}\right] \right|}_{\mathrm{\RN{1}}.\mathrm{\RN{2}}.\mathrm{\RN{1}}.\mathrm{\RN{1}}.\mathrm{\RN{3}}} \label{eq: term 1.2.1.1.3 in step 2.1 in tensor regression without sample splitting} \\
+ & \underbrace{\left| \operatorname{tr}\left[\left(\mcP_{U_{1\perp}}\otimes \mcP_{U_3}\right)\Mat_2\left(\widetilde{A}_1\right)^{\top}U_2\left(G_2G_2^{\top}\right)^{-1}G_2\left(U_1\otimes U_3\right)^{\top}\whZ_2^{(2)\top}\mcP_{U_{2\perp}}\whZ_2^{(2)}\right] \right|}_{\mathrm{\RN{1}}.\mathrm{\RN{2}}.\mathrm{\RN{1}}.\mathrm{\RN{1}}.\mathrm{\RN{4}}} \label{eq: term 1.2.1.1.4 in step 2.1 in tensor regression without sample splitting} \\
+ & \underbrace{\left\|\mcP_{U_{2\perp}}\whZ_2 (\mcP_{U_{1\perp}}\otimes \mcP_{U_3} )\Mat_2 (\widetilde{A}_1 )^{\top}U_2\right\|_{\mathrm{F}} \left\| (G_2G_2^{\top} )^{-1}G_2 (U_1\otimes U_3 )^{\top}\left[\big(\mcP_{\whU_1^{(0)}}-\mcP_{U_1}\big)\otimes \mcP_{U_3}\right]\whZ_2^{\top}\mcP_{U_{2\perp}}\right\|_{\mathrm{F}}}_{\mathrm{\RN{1}}.\mathrm{\RN{2}}.\mathrm{\RN{1}}.\mathrm{\RN{1}}.\mathrm{\RN{5}}} \label{eq: term 1.2.1.1.5 in step 2.1 in tensor regression without sample splitting} \\
+ & \underbrace{\left\|\mcP_{U_{2\perp}}\whZ_2 (\mcP_{U_{1\perp}}\otimes \mcP_{U_3} )\Mat_2 (\widetilde{A}_1 )^{\top}U_2\right\|_{\mathrm{F}} \left\| (G_2G_2^{\top} )^{-1}G_2 (U_1\otimes U_3 )^{\top}\left[\mcP_{U_1}\otimes \big(\mcP_{\whU_3^{(0)}} - \mcP_{U_3}\big)\right]\whZ_2^{\top}\mcP_{U_{2\perp}}\right\|_{\mathrm{F}}}_{\mathrm{\RN{1}}.\mathrm{\RN{2}}.\mathrm{\RN{1}}.\mathrm{\RN{1}}.\mathrm{\RN{6}}} \label{eq: term 1.2.1.1.6 in step 2.1 in tensor regression without sample splitting} \\
+ & \underbrace{\left| \operatorname{tr}\left[ (\mcP_{U_{1\perp}}\otimes \mcP_{U_3} )\Mat_2 (\widetilde{A}_1 )^{\top}U_2 (G_2G_2^{\top} )^{-1}G_2 (U_1\otimes U_3 )^{\top}\left[\big(\mcP_{\whU_1^{(0)}}-\mcP_{U_1}\big)\otimes  \big(\mcP_{\whU_3^{(0)}} - \mcP_{U_3} \big)\right]\whZ_2^{\top}\mcP_{U_{2\perp}}\whZ_2\right] \right|}_{\mathrm{\RN{1}}.\mathrm{\RN{2}}.\mathrm{\RN{1}}.\mathrm{\RN{1}}.\mathrm{\RN{7}}} \label{eq: term 1.2.1.1.7 in step 2.1 in tensor regression without sample splitting},
\end{align}
where $\whZ_2^{(1)}=\frac{1}{n\sigma^2}\sum_{i=1}^n\xi_i\operatorname{Mat}_2\left(\mcX_i\right)$ and $\whZ_2^{(1)}=\frac{1}{n\sigma^2}\sum_{i=1}^n\left[\left\langle\mcX_i, \widehat{\Delta}\right\rangle\operatorname{Mat}_2\left(\mcX_i\right)-\sigma^2\cdot \widehat{\Delta}\right]$.

Here, we have
\begin{align}
\mathrm{\RN{1}}.\mathrm{\RN{2}}.\mathrm{\RN{1}}.\mathrm{\RN{1}}.\mathrm{\RN{5}} + \mathrm{\RN{1}}.\mathrm{\RN{2}}.\mathrm{\RN{1}}.\mathrm{\RN{1}}.\mathrm{\RN{6}} + \mathrm{\RN{1}}.\mathrm{\RN{2}}.\mathrm{\RN{1}}.\mathrm{\RN{1}}.\mathrm{\RN{7}} 
\lesssim & \left\|\mcP_{U_{1\perp}}A_1 \mcP_{\left(U_3 \otimes U_2\right) G_1^{\top}}\right\|_{\mathrm{F}} \cdot \frac{\sigma_\xi^3\oor^{1/2}}{\ulambda^2\sigma^3}\cdot \frac{\op^{3/2}}{n^{3/2}}. \label{eq: upper bound of term 1.2.1.1.5+1.2.1.1.6+1.2.1.1.7 in step 2.1 in tensor regression without sample splitting}
\end{align}
It remains to find upper bounds for $\mathrm{\RN{1}}.\mathrm{\RN{2}}.\mathrm{\RN{1}}.\mathrm{\RN{1}}.\mathrm{\RN{1}}$ \eqref{eq: term 1.2.1.1.1 in step 2.1 in tensor regression without sample splitting}, $\mathrm{\RN{1}}.\mathrm{\RN{2}}.\mathrm{\RN{1}}.\mathrm{\RN{1}}.\mathrm{\RN{2}}$ \eqref{eq: term 1.2.1.1.2 in step 2.1 in tensor regression without sample splitting}, $\mathrm{\RN{1}}.\mathrm{\RN{2}}.\mathrm{\RN{1}}.\mathrm{\RN{1}}.\mathrm{\RN{3}}$ \eqref{eq: term 1.2.1.1.3 in step 2.1 in tensor regression without sample splitting} and $\mathrm{\RN{1}}.\mathrm{\RN{2}}.\mathrm{\RN{1}}.\mathrm{\RN{1}}.\mathrm{\RN{4}}$ \eqref{eq: term 1.2.1.1.4 in step 2.1 in tensor regression without sample splitting}. By Lemma~\ref{lemma: high-prob upper bound of tr(BZ(1)tCZ(1))}, it follows that
\begin{align}
\mathrm{\RN{1}}.\mathrm{\RN{2}}.\mathrm{\RN{1}}.\mathrm{\RN{1}}.\mathrm{\RN{1}} 
\lesssim & \left\|\mcP_{U_{1\perp}}A_1\mcP_{\left(U_3\otimes U_2\right)G_1^{\top}}\right\|_{\mathrm{F}} \cdot \left(\frac{\sigma_\xi^2}{\ulambda\sigma^2}\cdot \frac{\sqrt{\op\log(\op)}}{n}\right). \label{eq: upper bound of term 1.2.1.1.1 in step 2.1 in tensor regression without sample splitting}
\end{align}
where the first inequality follows from
\begin{align*}
& \mathbb{E}\operatorname{tr}\left[\left(\mcP_{U_{1\perp}}\otimes \mcP_{U_3}\right)\Mat_2\left(\widetilde{A}_1\right)^{\top}U_2\left(G_2G_2^{\top}\right)^{-1}G_2\left(U_1\otimes U_3\right)^{\top}\whZ_2^{(1)\top}\mcP_{U_{2\perp}}\whZ_2^{(1)}\right]  = 0 ,
\end{align*}
and the second inequality follows from 
\begin{align*}
\left\|U_2^{\top}\Mat_2\left(\widetilde{A}_1\right)\left(\mcP_{U_{1\perp}}\otimes \mcP_{U_3}\right)\right\|_{\mathrm{F}}
= & \left\|\mcP_{U_{1\perp}}\widetilde{A}_1\left(\mcP_{U_3}\otimes \mcP_{U_2}\right)\right\|_{\mathrm{F}} = \left\|\mcP_{U_{1\perp}}A_1\mcP_{\left(U_3\otimes U_2\right)G_1^{\top}}\right\|_{\mathrm{F}}.
\end{align*}
In addition, we have
\begin{align}
\mathrm{\RN{1}}.\mathrm{\RN{2}}.\mathrm{\RN{1}}.\mathrm{\RN{1}}.\mathrm{\RN{2}}
\leq & \frac{1}{\ulambda} \cdot \left\|\mcP_{U_{2\perp}}\whZ_2^{(2)}\left(\mcP_{U_{1\perp}}\otimes \mcP_{U_3}\right)\Mat_2\left(\widetilde{A}_1\right)^{\top}U_2\right\|_{\mathrm{F}} \cdot \left\|\left(U_1\otimes U_3\right)^{\top}\whZ_2^{(1)\top}\mcP_{U_{2\perp}}\right\|_{\mathrm{F}} \notag \\
\lesssim & \left\|\mcP_{U_{1\perp}}A_1 \mcP_{\left(U_3 \otimes U_2\right) G_1^{\top}}\right\|_{\mathrm{F}}\cdot \frac{\sigma_\xi^2\oor^{1/2}}{\ulambda\sigma^2}\cdot \Delta\cdot \frac{\op}{n} \label{eq: upper bound of term 1.2.1.1.2 in step 2.1 in tensor regression without sample splitting}
\end{align}
and by symmetry, 
\begin{align}
\mathrm{\RN{1}}.\mathrm{\RN{2}}.\mathrm{\RN{1}}.\mathrm{\RN{1}}.\mathrm{\RN{3}} \lesssim & \left\|\mcP_{U_{1\perp}}A_1 \mcP_{\left(U_3 \otimes U_2\right) G_1^{\top}}\right\|_{\mathrm{F}}\cdot \frac{\sigma_\xi^2\oor^{1/2}}{\ulambda\sigma^2}\cdot \Delta\cdot \frac{\op}{n}. \label{eq: upper bound of term 1.2.1.1.3 in step 2.1 in tensor regression without sample splitting}
\end{align}
Moreover, we have
\begin{align}
\mathrm{\RN{1}}.\mathrm{\RN{2}}.\mathrm{\RN{1}}.\mathrm{\RN{1}}.\mathrm{\RN{4}}
\leq & \frac{1}{\ulambda} \cdot \left\|\mcP_{U_{2\perp}}\whZ_2^{(2)}\left(\mcP_{U_{1\perp}}\otimes \mcP_{U_3}\right)\Mat_2\left(\widetilde{A}_1\right)^{\top}U_2\right\|_{\mathrm{F}} \cdot \left\|\left(U_1\otimes U_3\right)^{\top}\whZ_2^{(2)\top}\mcP_{U_{2\perp}}\right\|_{\mathrm{F}} \notag \\
\lesssim & \left\|\mcP_{U_{1\perp}}A_1 \mcP_{\left(U_3 \otimes U_2\right) G_1^{\top}}\right\|_{\mathrm{F}}\cdot \frac{\sigma_\xi^2\oor^{1/2}}{\ulambda\sigma^2}\cdot \Delta^2 \cdot \frac{\op}{n}. \label{eq: upper bound of term 1.2.1.1.4 in step 2.1 in tensor regression without sample splitting}
\end{align}
Therefore, we have
\begin{align}
\mathrm{\RN{1}}.\mathrm{\RN{2}}.\mathrm{\RN{1}}.\mathrm{\RN{1}}
\lesssim & \left\|\mcP_{U_{1\perp}}A_1 \mcP_{\left(U_3 \otimes U_2\right) G_1^{\top}}\right\|_{\mathrm{F}} \cdot \left[\frac{\sigma_\xi\oor^{1/2}}{\ulambda\sigma}\cdot \left(\frac{\sqrt{\op\log(\op)}}{n}+ \Delta\cdot \frac{\op}{n}\right)\right]. \label{eq: upper bound of term 1.2.1.1 in step 2.1 in tensor regression without sample splitting}
\end{align}
For $\mathrm{\RN{1}}.\mathrm{\RN{2}}.\mathrm{\RN{1}}.\mathrm{\RN{2}}$,
\begin{align}
\mathrm{\RN{1}}.\mathrm{\RN{2}}.\mathrm{\RN{1}}.\mathrm{\RN{2}} 
\leq & \frac{1}{\ulambda^2}  \left\|\mcP_{U_{2\perp}}\whZ_2 (\mcP_{U_{1\perp}}\otimes \mcP_{U_3} )\Mat_2 (\widetilde{A}_1 )^{\top}U_2\right\| \cdot \left\|U_2^{\top}\whZ_2\left(\mcP_{U_1}\otimes \mcP_{U_3}\right)\right\|_{\mathrm{F}} \cdot \left\|\mcP_{U_{2\perp}}\whZ_2\left(\mcP_{U_1}\otimes \mcP_{U_3}\right)\right\|_{\mathrm{F}} \notag \\
+ & \frac{1}{\ulambda^2} \left\|\mcP_{U_{2\perp}}\whZ_2 (\mcP_{U_{1\perp}}\otimes \mcP_{U_3} )\Mat_2 (\widetilde{A}_1 )^{\top}U_2\right\| 
\left(\big\|\mcP_{\whU_3^{(0)}}-\mcP_{U_3}\big\| + \big\|\mcP_{\whU_1^{(0)}}-\mcP_{U_1}\big\| + \big\|\mcP_{\whU_3^{(0)}}-\mcP_{U_3}\big\|\big\|\mcP_{\whU_1^{(0)}}-\mcP_{U_1}\big\|\right) \notag \\
& \cdot \sup_{\substack{W_1\in \mathbb{R}^{p_1\times 2r_1}, \left\|W_1\right\|=1\\ W_3\in \mathbb{R}^{p_3\times 2r_3}, \left\|W_3\right\|=1}}\left\|U_2^{\top}\whZ_2\left(W_3\otimes W_1\right)\right\|_{\mathrm{F}} \cdot \sup_{\substack{W_1\in \mathbb{R}^{p_1\times 2r_1}, \left\|W_1\right\|=1\\ W_3\in \mathbb{R}^{p_3\times 2r_3}, \left\|W_3\right\|=1}}\left\|U_{2\perp}^{\top}\whZ_2\left(W_1\otimes W_3\right)\right\|_{\mathrm{F}} \notag \\
\lesssim & \left\|\mcP_{U_{1\perp}}A_1\mcP_{\left(U_3\otimes U_2\right)G_1^{\top}}\right\|_{\mathrm{F}}\cdot \frac{\sigma_\xi^3}{\ulambda^2\sigma^3}\cdot \frac{\op}{n} \cdot \left[\left(\sqrt{\frac{\oor\log(\op)}{n}}+\Delta\sqrt{\frac{\op}{n}}\right)+ \frac{\sigma_\xi}{\ulambda\sigma}\cdot \frac{\op}{n}\right]. \label{eq: upper bound of term 1.2.1.2 in step 2.1 in tensor regression without sample splitting}
\end{align}
Hence, we have
\begin{align}
\mathrm{\RN{1}}.\mathrm{\RN{2}}.\mathrm{\RN{1}}
\lesssim & \eqref{eq: upper bound of term 1.2.1.1 in step 2.1 in tensor regression without sample splitting} + \eqref{eq: upper bound of term 1.2.1.2 in step 2.1 in tensor regression without sample splitting} \lesssim \left\|\mcP_{U_{1 \perp}} A_1 \mcP_{\left(U_3 \otimes U_2\right) G_1^{\top}}\right\|_{\mathrm{F}} \cdot \left[\frac{\sigma_{\xi}\oor^{1/2}}{\ulambda \sigma} \cdot \left(\frac{\sqrt{\op \log (\op)}}{n} + \Delta \cdot \frac{\op}{n}\right)\right] \label{eq: upper bound of term 1.2.1 in step 2.1 in tensor regression without sample splitting}.
\end{align}

Furthermore, we have
\begin{align}
\mathrm{\RN{1}}.\mathrm{\RN{2}}.\mathrm{\RN{2}} 
\leq & \left\|\left(\mcP_{U_{1\perp}} \otimes \mcP_{U_3}\right) \operatorname{Mat}_2\left(\widetilde{A}_1\right)^{\top} \mcP_{U_2}\right\|_{\mathrm{F}} \cdot \sup_{W_2 \in \mathbb{R}^{p_2\times r_2}, \left\|W_2\right\|=1} \left\|W_2^{\top} \whZ_2 \left(\mcP_{U_{1\perp}} \otimes \mcP_{U_3}\right)\right\|_{\mathrm{F}} \notag \\
& \cdot \left(\left\|\mcP_{U_2} \sum_{k_2=2}^{+\infty} S_{G_2, k_2} \left(\whE_2\right) \mcP_{U_2}\right\| + \left\|\mcP_{U_2} \sum_{k_2=2}^{+\infty} S_{G_2, k_2}\left(\whE_2\right) \mcP_{U_{2 \perp}}\right\|\right) \notag \\
\lesssim & \left\|\mcP_{U_{1 \perp}} A_1 \mcP_{\left(U_3 \otimes U_2\right) G_1^{\top}}\right\|_{\mathrm{F}} \cdot \frac{\sigma_\xi^3\oor^{1/2}}{\ulambda^2 \sigma^3} \cdot \frac{\op^{3/2}}{n^{3/2}}. \label{eq: upper bound of term 1.2.2 in step 2.1 in tensor regression without sample splitting}
\end{align}
Therefore, we have
\begin{align}
\mathrm{\RN{1}}.\mathrm{\RN{2}} 
\lesssim & \eqref{eq: upper bound of term 1.2.1 in step 2.1 in tensor regression without sample splitting} + \eqref{eq: upper bound of term 1.2.2 in step 2.1 in tensor regression without sample splitting} \lesssim \left\|\mcP_{U_{1 \perp}} A_1 \mcP_{\left(U_3 \otimes U_2\right) G_1^{\top}}\right\|_{\mathrm{F}} \cdot\left[\frac{\sigma_{\xi}\oor^{1/2}}{\ulambda \sigma} \cdot\left(\frac{\sqrt{\op \log (\op)}}{n}+\Delta \cdot \frac{\op}{n}\right)\right] \label{eq: upper bound of term 1.2 in step 2.1 in tensor regression without sample splitting}.
\end{align}

Similar to the proof of \RN{1}.\RN{2}, we have the same upper bound for \RN{1}.\RN{3} as for \RN{1}.\RN{2}.

Furthermore, 
\begin{align}
\mathrm{\RN{1}}.\mathrm{\RN{4}} 
\leq & \left\|\mcP_{U_{1\perp}}A_1 \mcP_{\left(U_3 \otimes U_2\right)G_1^{\top}}\right\|_{\mathrm{F}} \cdot \sup_{\substack{W_2 \in \mathbb{R}^{p_2\times r_2}, \left\|W_2\right\|=1 \\ W_3 \in \mathbb{R}^{p_3\times r_3}, \left\|W_3\right\|=1}} \left\| \mcP_{U_{1\perp}} \whZ_1\left(W_3\otimes W_2\right) \right\|_{\mathrm{F}} \cdot \prod_{j=2}^3 \left\|\mcP_{\whU_j}-\mcP_{U_j}\right\| \notag \\
\lesssim & \left\|\mcP_{U_{1\perp}}A_1\mcP_{\left(U_3 \otimes U_2\right)G_1^{\top}}\right\|_{\mathrm{F}}\cdot \frac{\sigma_\xi^3\oor^{1/2}}{\ulambda^2\sigma^3}\cdot \frac{\op^{3/2}}{n^{3/2}}. \label{eq: upper bound of term 1.4 in step 2.1 in tensor regression without sample splitting}
\end{align}
Furthermore, 
\begin{align*}
\mathrm{\RN{1}}.\mathrm{\RN{5}} 
\leq & \underbrace{\left|\operatorname{tr}\left[\left(G_1G_1^{\top}\right)^{-1}  G_1 \left(U_3 \otimes U_2\right)^{\top}A_1^{\top}\mcP_{U_{1\perp}} \whZ_1\left(\mcP_{U_3} \otimes \mcP_{U_2}\right) \whZ_1^{\top}U_1 \right]\right|}_{\mathrm{\RN{1}}.\mathrm{\RN{5}}.\mathrm{\RN{1}}} \\
+ & \underbrace{\left|\operatorname{tr}\left[\left(G_1G_1^{\top}\right)^{-1}  G_1 \left(U_3 \otimes U_2\right)^{\top}A_1^{\top}\mcP_{U_{1\perp}} \whZ_1\left[\left(\mcP_{\whU_3^{(1)}}-\mcP_{U_3}\right) \otimes \mcP_{U_2}\right] \whZ_1^{\top}U_1 \right]\right|}_{\mathrm{\RN{1}}.\mathrm{\RN{5}}.\mathrm{\RN{2}}} \\
+ & \underbrace{\left|\operatorname{tr}\left[\left(G_1G_1^{\top}\right)^{-1}  G_1 \left(U_3 \otimes U_2\right)^{\top}A_1^{\top}\mcP_{U_{1\perp}} \whZ_1\left[\mcP_{U_3} \otimes \left(\mcP_{\whU_2}-\mcP_{U_2}\right)\right] \whZ_1^{\top}U_1 \right]\right|}_{\mathrm{\RN{1}}.\mathrm{\RN{5}}.\mathrm{\RN{3}}} \\
+ & \underbrace{\left|\operatorname{tr}\left[\left(G_1G_1^{\top}\right)^{-1}  G_1 \left(U_3 \otimes U_2\right)^{\top}A_1^{\top}\mcP_{U_{1\perp}} \whZ_1\left[\left(\mcP_{\whU_3}-\mcP_{U_3}\right) \otimes \left(\mcP_{\whU_2}-\mcP_{U_2}\right)\right] \whZ_1^{\top}U_1 \right]\right|}_{\mathrm{\RN{1}}.\mathrm{\RN{5}}.\mathrm{\RN{4}}} .
\end{align*}
Here, we first have
\begin{align}
\mathrm{\RN{1}}.\mathrm{\RN{5}}.\mathrm{\RN{1}}
\leq & \left\|\left(G_1G_1^{\top}\right)^{-1}  G_1 \left(U_3 \otimes U_2\right)^{\top}A_1^{\top}\mcP_{U_{1\perp}}\whZ_1\left(\mcP_{U_3} \otimes \mcP_{U_2}\right)\right\|_{\mathrm{F}} \cdot \underbrace{\left\|\left(\mcP_{U_3} \otimes \mcP_{U_2}\right) \whZ_1^{\top}U_1\right\|_{\mathrm{F}}}_{\eqref{eq: high-prob upper bound of U1Zhat1(U3oU2) in tensor regression without sample splitting}} \notag \\
\lesssim & \left\|\mcP_{U_{1\perp}}A_1\mcP_{\left(U_3\otimes U_2\right)G_1^{\top}}\right\|_{\mathrm{F}} \cdot \frac{\sigma_\xi^2\oor^{1/2}}{\ulambda\sigma^2}\cdot \left(\frac{\oor\log(\op)}{n}+ \Delta^2 \cdot \frac{\op}{n}\right). \label{eq: upper bound of term 1.5.1 in step 2.1 in tensor regression without sample splitting} 
\end{align}
Second, we have
\begin{align}
\mathrm{\RN{1}}.\mathrm{\RN{5}}.\mathrm{\RN{2}}
\leq & \frac{1}{\ulambda}\left\|\mcP_{U_{1\perp}}A_1\mcP_{\left(U_3\otimes U_2\right)G_1^{\top}}\right\|_{\mathrm{F}}  \sup_{\substack{W_2\in \mathbb{R}^{p_2\times 2r_2}, \left\|W_2\right\|=1\\ W_3 \in \mathbb{R}^{p_3\times 2r_3}, \left\|W_3\right\|=1}}\left\|\mcP_{U_{1\perp}}\whZ_1\left(W_3\otimes W_2\right)\right\| 
  \left\|\mcP_{\whU_2^{(1)}}-\mcP_{U_2}\right\| \notag\\
  &\cdot \sup_{\substack{W_2\in \mathbb{R}^{p_2\times 2r_2}, \left\|W_2\right\|=1\\ W_3 \in \mathbb{R}^{p_3\times 2r_3}, \left\|W_3\right\|=1}}\left\|\mcP_{U_1}\whZ_1\left(W_3\otimes W_2\right)\right\|_{\mathrm{F}} \notag \\
\lesssim & \left\|\mcP_{U_{1\perp}}A_1\mcP_{\left(U_3\otimes U_2\right)G_1^{\top}}\right\|_{\mathrm{F}}\cdot \frac{\sigma_\xi^3\oor^{1/2}}{\ulambda^2\sigma^3}\cdot \frac{\op^{3/2}}{n^{3/2}}. \label{eq: upper bound of term 1.5.2 in step 2.1 in tensor regression without sample splitting}
\end{align}
Similar to the proof of \RN{1}.\RN{5}.\RN{2}, we have the same upper bound for \RN{1}.\RN{5}.\RN{3} as for \RN{1}.\RN{5}.\RN{2}.


Finally, we have
\begin{align}
\mathrm{\RN{1}}.\mathrm{\RN{5}}.\mathrm{\RN{4}}
\leq & \frac{1}{\ulambda}\left\|\mcP_{U_{1\perp}}A_1\mcP_{\left(U_3\otimes U_2\right)G_1^{\top}}\right\|_{\mathrm{F}}\cdot \sup_{\substack{W_2\in \mathbb{R}^{p_2\times 2r_2}, \left\|W_2\right\|=1\\ W_3 \in \mathbb{R}^{p_3\times 2r_3}, \left\|W_3\right\|=1}}\left\|\mcP_{U_{1\perp}}\whZ_1\left(W_3\otimes W_2\right)\right\| \notag \\
& \cdot \left\|\mcP_{\whU_2^{(1)}}-\mcP_{U_2}\right\| \cdot \left\|\mcP_{\whU_3^{(1)}}-\mcP_{U_3}\right\| \cdot \sup_{\substack{W_2\in \mathbb{R}^{p_2\times 2r_2}, \left\|W_2\right\|=1\\ W_3 \in \mathbb{R}^{p_3\times 2r_3, \left\|W_3\right\|=1}}}\left\|\mcP_{U_1}\whZ_1\left(W_3\otimes W_2\right)\right\|_{\mathrm{F}} \notag \\
\lesssim & \left\|\mcP_{U_{1\perp}}A_1\mcP_{\left(U_3\otimes U_2\right)G_1^{\top}}\right\|_{\mathrm{F}}\cdot \frac{\sigma_\xi^4\oor^{1/2}}{\ulambda^3\sigma^4}\cdot \frac{\op^2}{n^2} . \label{eq: upper bound of term 1.5.4 in step 2.1 in tensor regression without sample splitting}.
\end{align}
It further implies that
\begin{align}
\mathrm{\RN{1}}.\mathrm{\RN{5}} 
\lesssim & \eqref{eq: upper bound of term 1.5.1 in step 2.1 in tensor regression without sample splitting} + \eqref{eq: upper bound of term 1.5.2 in step 2.1 in tensor regression without sample splitting} + \eqref{eq: upper bound of term 1.5.4 in step 2.1 in tensor regression without sample splitting} \lesssim \left\| \mcP_{U_{1 \perp}} A_1 \mcP_{\left(U_3 \otimes U_2\right) G_1^{\top}} \right\|_{\mathrm{F}} \cdot \left[\frac{\sigma_{\xi}^2\oor^{1/2}}{\ulambda \sigma^2} \cdot\left(\frac{\oor\log(\op)}{n}+\Delta^2 \cdot \frac{\op}{n}\right)\right]. \label{eq: upper bound of term 1.5 in step 2.1 in tensor regression without sample splitting}
\end{align}

Combining the results above, we have
\begin{align}
\mathrm{\RN{1}} 
\lesssim & \eqref{eq: upper bound of term 1.2 in step 2.1 in tensor regression without sample splitting} + \eqref{eq: upper bound of term 1.4 in step 2.1 in tensor regression without sample splitting} + \eqref{eq: upper bound of term 1.5 in step 2.1 in tensor regression without sample splitting} \lesssim \left\|\mcP_{U_{1 \perp}} A_1 \mcP_{\left(U_3 \otimes U_2\right) G_1^{\top}}\right\|_{\mathrm{F}} \cdot\left[\frac{\sigma_{\xi}^2 \oor^{1/2}}{\ulambda \sigma^2} \cdot\left(\frac{\sqrt{\op \log (\op)}}{n}+\Delta \cdot \frac{\op}{n}\right)\right] \label{eq: upper bound of term 1 in step 2.1 in tensor regression without sample splitting}.
\end{align}

Then, consider an upper bound for the remaining higher-order terms. It follows that
\begin{align}
& \mathrm{\RN{2}} \notag \\
\leq & \left| \operatorname{tr}\left[\left(\mcP_{U_3}\otimes \mcP_{U_2}\right)A_1^{\top}S_{G_1,2}\left(\whE_1\right)U_1G_1\left(U_3 \otimes U_2\right)^{\top}\right] \right| 
+  \left| \operatorname{tr}\left[\left(\mcP_{U_3}\otimes \mcP_{U_2}\right)A_1^{\top}\sum_{k_1=3}^{+\infty} S_{G_1,k_1}\left(\whE_1\right)U_1G_1\left(U_3 \otimes U_2\right)^{\top}\right] \right| \notag \\
\lesssim & \left\|U_1^{\top} A_1 \mcP_{\left(U_3 \otimes U_2\right) G_1^{\top}}\right\|_{\mathrm{F}} \cdot \frac{\sigma_\xi^2\oor^{1/2}}{\ulambda\sigma^2}\cdot \frac{\op}{n} \notag \\
+ & \underbrace{\left\|U_{1\perp}^{\top}A_1^{\top}\mcP_{\left(U_3\otimes U_2\right)G_1^{\top}}\right\|_{\mathrm{F}} \cdot \left[\frac{\sigma_\xi}{\sigma} \cdot \left(\sqrt{\frac{\oor\log(\op)}{n}} +\Delta\cdot \sqrt{\frac{\op}{n}}\right)\right]}_{\eqref{eq: high-prob upper bound of V1tP1(0)Ehat1P1(-1/2) in tensor regression without sample splitting}} \cdot \left[\underbrace{\frac{\sigma_\xi\oor^{1/2}}{\ulambda\sigma} \cdot \left(\sqrt{\frac{\oor \log (\op)}{n}} + \Delta\sqrt{\frac{\op}{n}}\right)}_{\eqref{eq: high-prob upper bound of P1(-1/2)Ehat1P1(-1/2) in tensor regression without sample splitting}}\right] \notag \\
+ & \underbrace{\left\|U_{1\perp}^{\top}A_1^{\top}\mcP_{\left(U_3\otimes U_2\right)G_1^{\top}}\right\|_{\mathrm{F}} \cdot \frac{\oor^{1/2}}{\ulambda^2} \cdot \left[\frac{\sigma_\xi^2}{\sigma^2} \cdot \left(\frac{\sqrt{\oor \log (\op)} \cdot \sqrt{\op}}{n} + \Delta \cdot \frac{\op}{n}\right)\right]}_{\eqref{eq: high-prob upper bound of V1tP1(0)Ehat1P1(0) in tensor regression without sample splitting}} \cdot \frac{\sigma_\xi}{\sigma} \sqrt{\frac{\op}{n}} \notag \\
+ & \left\|U_1^{\top}A_1\mcP_{\left(U_3\otimes U_2\right)G_1^\top}\right\|_{\mathrm{F}}\cdot \frac{\sigma_{\xi}^3}{\ulambda^2 \sigma^3} \cdot \frac{\op^{3 / 2}}{n^{3 / 2}} + \underbrace{\left\|U_{1\perp}^{\top}A_1\mcP_{\left(U_3\otimes U_2\right)G_1^\top}\right\|_{\mathrm{F}}\cdot \left[\frac{\sigma_{\xi}^3\oor^{1/2}}{\ulambda^2 \sigma^3} \cdot\left(\frac{\op\sqrt{\oor\log(\op)}}{n^{3 / 2}}+\Delta \cdot \frac{\op^{3 / 2}}{n^{3 / 2}}\right)\right]}_{\eqref{eq: order3 term in step 2.1 in tensor regression without sample splitting}} \notag \\
\lesssim & \left\|U_1^{\top}A_1\mcP_{\left(U_3\otimes U_2\right)G_1^\top}\right\|_{\mathrm{F}}\cdot \frac{\sigma_\xi^2\oor^{1/2}}{\ulambda\sigma^2}\cdot \frac{\op}{n} \notag \\
+ & \left\|U_{1\perp}^{\top}A_1^{\top}\mcP_{\left(U_3\otimes U_2\right)G_1^{\top}}\right\|_{\mathrm{F}}\cdot \left[\frac{\sigma_\xi^3\oor^{1/2}}{\ulambda\sigma^3} \cdot \left(\frac{\oor\log(\op)}{n} + \Delta^2\cdot \frac{\op}{n}\right)+\frac{\sigma_{\xi}^3\oor^{1/2}}{\ulambda^2 \sigma^3} \cdot\left(\frac{\op\sqrt{\oor\log(\op)}}{n^{3 / 2}}+\Delta \cdot \frac{\op^{3 / 2}}{n^{3 / 2}}\right)\right]. \label{eq: upper bound of term 2 in step 2.1 in tensor regression without sample splitting}
\end{align}

Here, we have
\begin{align}
& \left|\operatorname{tr}\left[\left(\mcP_{U_3} \otimes \mcP_{U_2}\right) A_1^{\top} \sum_{k_1=3}^{+\infty} S_{G_1, k_1}\left(\whE_1\right) U_1 G_1\left(U_3 \otimes U_2\right)^{\top}\right]\right| \notag \\
\leq & \left|\operatorname{tr}\left[\left(\mcP_{U_3} \otimes \mcP_{U_2}\right) A_1^{\top}\mcP_{U_1}\sum_{k_1=3}^{+\infty} S_{G_1, k_1}\left(\whE_1\right) U_1 G_1\left(U_3 \otimes U_2\right)^{\top}\right]\right| \notag \\
& +  \underbrace{\left|\operatorname{tr}\left[\left(\mcP_{U_3} \otimes \mcP_{U_2}\right) A_1^{\top}\mcP_{U_{1\perp}}\sum_{k_1=3}^{+\infty} S_{G_1, k_1}\left(\whE_1\right) U_1 G_1\left(U_3 \otimes U_2\right)^{\top}\right]\right|}_{\eqref{eq: high-prob upper bound of PUporder3PUpPV in tensor regression without sample splitting}} \notag \\
\lesssim & \left\|U_1^{\top}A_1\mcP_{\left(U_3\otimes U_2\right)G_1^\top}\right\|_{\mathrm{F}}\cdot \frac{\sigma_{\xi}^3\oor^{1/2}}{\ulambda^2 \sigma^3} \cdot \frac{\op^{3 / 2}}{n^{3 / 2}} + \left\|U_{1\perp}^{\top}A_1\mcP_{\left(U_3\otimes U_2\right)G_1^\top}\right\|_{\mathrm{F}} \left[\frac{\sigma_{\xi}^3\oor^{1/2}}{\ulambda^2 \sigma^3} \left(\frac{\op\sqrt{\oor\log(\op)}}{n^{3 / 2}}+\Delta \cdot \frac{\op^{3 / 2}}{n^{3 / 2}}\right)\right]. \label{eq: order3 term in step 2.1 in tensor regression without sample splitting}
\end{align} 

Therefore, we have 
\begin{align*}
& \left|\left\langle\mcT \times_1 \left(\mcP_{\whU_1} - \mcP_{U_1}\right) \times_2 \mcP_{U_2} \times_3 \mcP_{U_3}, \mcA\right\rangle-\left\langle\mcP_{U_{1\perp}} \whZ_1\mcP_{\left(U_3 \otimes U_2\right) G_1^{\top}}, A_1\right\rangle \right| 
\lesssim  \eqref{eq: upper bound of term 1 in step 2.1 in tensor regression without sample splitting} + \eqref{eq: upper bound of term 2 in step 2.1 in tensor regression without sample splitting} \\
\lesssim & \left\|U_1^{\top}A_1\mcP_{\left(U_3\otimes U_2\right)G_1^\top}\right\|_{\mathrm{F}}\cdot \frac{\sigma_\xi^2\oor^{1/2}}{\ulambda\sigma^2}\cdot \frac{\op}{n} + \left\|\mcP_{U_{1 \perp}} A_1 \mcP_{\left(U_3 \otimes U_2\right) G_1^{\top}}\right\|_{\mathrm{F}} \cdot\left[\frac{\sigma_{\xi}^2 \oor^{1/2}}{\ulambda \sigma^2} \cdot\left(\frac{\sqrt{\oor\op \log (\op)}}{n}+\Delta \cdot \frac{\op}{n}\right)\right].
\end{align*}

\subsubsection*{Step 2.2: Upper Bound of $ \langle \mcT\times_1 \left(\mcP_{\whU_1} - \mcP_{U_1}\right) \times_2 \left(\mcP_{\whU_2} - \mcP_{U_2}\right) \times_3 \mcP_3, \mcA \rangle$}

Consider the following decomposition:
\begin{align}
& \left|\left\langle \mcT\times_1 \left(\mcP_{\whU_1} - \mcP_{U_1}\right) \times_2 \left(\mcP_{\whU_2} - \mcP_{U_2}\right) \times_3 \mcP_3, \mcA\right\rangle\right| \notag \\
\leq & \underbrace{\left|\operatorname{tr}\left[\left(\mcP_{U_3}\otimes \mcP_{U_2}\left(\mcP_{\whU_2}-\mcP_{U_2}\right)\mcP_{U_{2\perp}}\right)A_1\mcP_{U_{1\perp}}\left(\mcP_{\whU_1}-\mcP_{U_1}\right)U_1G_1\left(U_3\otimes U_2\right)^{\top}\right]\right|}_{\mathrm{\RN{1}}} \label{eq: term 1 in step 2.2 in tensor regression without sample splitting}\\
+ &  \underbrace{\left|\operatorname{tr}\left[\left(\mcP_{U_3}\otimes \mcP_{U_2}\left(\mcP_{\whU_2}-\mcP_{U_2}\right)\mcP_{U_2}\right)A_1\mcP_{U_{1\perp}}\left(\mcP_{\whU_1}-\mcP_{U_1}\right)U_1G_1\left(U_3\otimes U_2\right)^{\top}\right]\right|}_{\mathrm{\RN{2}}} \label{eq: term 2 in step 2.2 in tensor regression without sample splitting}\\
+ &  \underbrace{\left|\operatorname{tr}\left[\left(\mcP_{U_3}\otimes \mcP_{U_2}\left(\mcP_{\whU_2}-\mcP_{U_2}\right)\mcP_{U_{2\perp}}\right)A_1\mcP_{U_1}\left(\mcP_{\whU_1}-\mcP_{U_1}\right)U_1G_1\left(U_3\otimes U_2\right)^{\top}\right]\right|}_{\mathrm{\RN{3}}} \label{eq: term 3 in step 2.2 in tensor regression without sample splitting}\\
+ &  \underbrace{\left|\operatorname{tr}\left[\left(\mcP_{U_3}\otimes \mcP_{U_2}\left(\mcP_{\whU_2}-\mcP_{U_2}\right)\mcP_{U_2}\right)A_1\mcP_{U_1}\left(\mcP_{\whU_1}-\mcP_{U_1}\right)U_1G_1\left(U_3\otimes U_2\right)^{\top}\right]\right|}_{\mathrm{\RN{4}}}. \label{eq: term 4 in step 2.2 in tensor regression without sample splitting}
\end{align}

First, we have
\begin{align}
\mathrm{\RN{1}}
\leq & \underbrace{\left\|\left(\mcP_{U_3} \otimes \mcP_{U_2}\left(\mcP_{\whU_2}-\mcP_{U_2}\right) \mcP_{U_{2 \perp}}\right)\left(\mcP_{U_3} \otimes \mcP_{U_{2 \perp}}\right) A_1 \mcP_{U_{1 \perp}}\left(\mcP_{\whU_1}-\mcP_{U_1}\right) \mcP_1^{\frac{1}{2}}\right\|_{\mathrm{F}}}_{\eqref{eq: high-prob upper bound of AoPU1oPU2p(PUhat2-P_U2)PU2oPU3p(PUhat3-P_U3)PU3 in tensor regression without sample splitting}} \cdot \left\|\mcP_1^{-\frac{1}{2}}U_1 G_1\left(U_3 \otimes U_2\right)^{\top}\right\|_{\mathrm{F}} \notag \\
\lesssim & \left\|\mcA \times_1 U_1\right\|_{\mathrm{F}} \cdot \left[\frac{\sigma_{\xi}^2\oor^{1/2}}{\ulambda \sigma^2} \cdot\left(\frac{\oor\log(\op)}{n}+\Delta \cdot \frac{\sqrt{\oR \op\log (\op)}}{n}+\Delta^2 \cdot \frac{\op}{n}\right)\right] \label{eq: upper bound of term 1 in step 2.2 in tensor regression without sample splitting}.
\end{align}

Then, consider
\begin{align}
\mathrm{\RN{2}}
\leq & \underbrace{\left\|\mcP_{U_2}\left(\mcP_{\whU_2}-\mcP_{U_2}\right) \mcP_{U_2}\right\|_{\mathrm{F}}}_{\eqref{eq: high-prob upper bound of PU1(PUhat1-PU1)PU1 in tensor regression}} \cdot \left\|\left(\mcP_{U_3} \otimes \mcP_{U_2}\right) A_1 \mcP_{U_{1 \perp}}\whZ_1\mcP_{\left(U_3 \otimes U_2\right)G_1^{\top}}\right\|_{\mathrm{F}} \notag \\
+ & \underbrace{\left\|\mcP_{U_2}\left(\mcP_{\whU_2}-\mcP_{U_2}\right) \mcP_{U_2}\right\|_{\mathrm{F}}}_{\eqref{eq: high-prob upper bound of PU1(PUhat1-PU1)PU1 in tensor regression}} \cdot \sup_{\substack{W_2\in \mathbb{R}^{p_2\times r_2}, \left\|W_2\right\|=1\\ W_3 \in \mathbb{R}^{p_3\times r_3}, \left\|W_3\right\|=1}}\left\|\left(\mcP_{U_3} \otimes \mcP_{U_2}\right) A_1 \mcP_{U_{1 \perp}}\whZ_1\left(W_3\otimes W_2\right)\right\|_{\mathrm{F}} \notag \\
& \cdot \left(\left\|\mcP_{\whU_3^{(1)}}- \mcP_{U_3}\right\| + \left\|\mcP_{\whU_2^{(1)}}- \mcP_{U_2}\right\| + \prod_{j=2}^3 \left\|\mcP_{\whU_j^{(1)}}- \mcP_{U_j}\right\|\right) \notag \\
\lesssim & \left\|\mcA \times_2 U_2 \times_3 U_3\right\|_{\mathrm{F}} \cdot \left[\frac{\sigma_\xi^3\oor^{1/2}}{\ulambda^2\sigma^3}\cdot \left(\frac{\op\sqrt{\oor\log(\op)}}{n^{3/2}}+ \Delta\cdot \frac{\op^{3/2}}{n^{3/2}} \right)\right] . \label{eq: upper bound of term 2 in step 2.2 in tensor regression without sample splitting}
\end{align}

By symmetry, it also implies that
\begin{align}
\mathrm{\RN{3}}
\lesssim & \left\|\mcA \times_1 U_1 \times_3 U_3\right\|_{\mathrm{F}} \cdot \left[\frac{\sigma_\xi^3\oor^{1/2}}{\ulambda^2\sigma^3}\cdot \left(\frac{\op\sqrt{\oor\log(\op)}}{n^{3/2}}+ \Delta\cdot \frac{\op^{3/2}}{n^{3/2}} \frac{}{}\right)\right] .  \label{eq: upper bound of third term in step 2.2 in tensor regression without sample splitting}
\end{align}

Finally, consider
\begin{align}
\mathrm{\RN{4}}
\leq & \underbrace{\left\|\mcP_{U_2}\left(\mcP_{\whU_2}-\mcP_{U_2}\right) \mcP_{U_2}\right\|}_{\eqref{eq: high-prob upper bound of PU1(PUhat1-PU1)PU1 in tensor regression}} \cdot \left\|\left(\mcP_{U_3} \otimes \mcP_{U_2}\right) A_1 \mcP_{U_1}\right\|_{\mathrm{F}} \cdot \underbrace{\left\|\mcP_{U_1}\left(\mcP_{\whU_1}-\mcP_{U_1}\right) U_1 G_1\left(U_3 \otimes U_2\right)^{\top}\right\|_{\mathrm{F}}}_{\eqref{eq: high-prob upper bound of PU1(PUhat1-PU1)PU1 in tensor regression}} \notag \\
\lesssim & \left\|\mcA \times_1 U_1 \times_2 U_2 \times_3 U_3\right\|_{\mathrm{F}} \cdot \frac{\sigma_\xi^4\oor^{1/2}}{\ulambda^3\sigma^4} \cdot \frac{\op^2}{n^2}.  \label{eq: upper bound of fourth term in step 2.2 in tensor regression without sample splitting}
\end{align}

\subsubsection*{Step 2.3: Upper Bound of $\langle \mcT\times_1 \left(\mcP_{\widehat{U}_1} - \mcP_{U_1}\right) \times_2 \left(\mcP_{\widehat{U}_2} - \mcP_{U_2}\right) \times_3 \left(\mcP_{\widehat{U}_3} - \mcP_{U_3}\right), \mcA\rangle$}

By similar arguments, we have
\begin{align*}
& \left|\left\langle \mcT\times_1 \left(\mcP_{\whU_1} - \mcP_{U_1}\right) \times_2 \left(\mcP_{\whU_2} - \mcP_{U_2}\right) \times_3 \left(\mcP_{\whU_3} - \mcP_{U_3}\right), \mcA\right\rangle \right| \\
\leq & \underbrace{\left|\operatorname{tr}\left[\left(\mcP_{U_3}\left(\mcP_{\whU_3}-\mcP_{U_3}\right)\mcP_{U_3} \otimes \mcP_{U_2}\left(\mcP_{\whU_2}-\mcP_{U_2}\right)\mcP_{U_2}\right)A_1\mcP_{U_1}\left(\mcP_{\whU_1}-\mcP_{U_1}\right)U_1G_1\left(U_3\otimes U_2\right)^{\top}\right]\right|}_{\mathrm{\RN{1}}} \\
+ & \underbrace{\left|\operatorname{tr}\left[\left(\mcP_{U_3}\left(\mcP_{\whU_3}-\mcP_{U_3}\right)\mcP_{U_3} \otimes \mcP_{U_2}\left(\mcP_{\whU_2}-\mcP_{U_2}\right)\mcP_{U_{2\perp}}\right)A_1\mcP_{U_1}\left(\mcP_{\whU_1}-\mcP_{U_1}\right)U_1G_1\left(U_3\otimes U_{2\perp}\right)^{\top}\right]\right|}_{\mathrm{\RN{2}}} \\
+ & \underbrace{\left|\operatorname{tr}\left[\left(\mcP_{U_3}\left(\mcP_{\whU_3}-\mcP_{U_3}\right)\mcP_{U_{3\perp}} \otimes \mcP_{U_2}\left(\mcP_{\whU_2}-\mcP_{U_2}\right)\mcP_{U_2}\right)A_1\mcP_{U_1}\left(\mcP_{\whU_1}-\mcP_{U_1}\right)U_1G_1\left(U_3\otimes U_2\right)^{\top}\right]\right|}_{\mathrm{\RN{3}}} \\
+ & \underbrace{\left|\operatorname{tr}\left[\left(\mcP_{U_3}\left(\mcP_{\whU_3}-\mcP_{U_3}\right)\mcP_{U_{3\perp}} \otimes \mcP_{U_2}\left(\mcP_{\whU_2}-\mcP_{U_2}\right)\mcP_{U_{2\perp}}\right)A_1\mcP_{U_1}\left(\mcP_{\whU_1}-\mcP_{U_1}\right)U_1G_1\left(U_3\otimes U_2\right)^{\top}\right]\right|}_{\mathrm{\RN{4}}} \\
+ & \underbrace{\left|\operatorname{tr}\left[\left(\mcP_{U_3}\left(\mcP_{\whU_3}-\mcP_{U_3}\right)\mcP_{U_3} \otimes \mcP_{U_2}\left(\mcP_{\whU_2}-\mcP_{U_2}\right)\mcP_{U_2}\right)A_1\mcP_{U_{1\perp}}\left(\mcP_{\whU_1}-\mcP_{U_1}\right)U_1G_1\left(U_3\otimes U_2\right)^{\top}\right]\right|}_{\mathrm{\RN{5}}} \\
+ & \underbrace{\left|\operatorname{tr}\left[\left(\mcP_{U_3}\left(\mcP_{\whU_3}-\mcP_{U_3}\right)\mcP_{U_3} \otimes \mcP_{U_2}\left(\mcP_{\whU_2}-\mcP_{U_2}\right)\mcP_{U_{2\perp}}\right))A_1\mcP_{U_{1\perp}}\left(\mcP_{\whU_1}-\mcP_{U_1}\right)U_1G_1\left(U_3\otimes U_2\right)^{\top}\right]\right|}_{\mathrm{\RN{6}}} \\
+ & \underbrace{\left|\operatorname{tr}\left[\left(\mcP_{U_3}\left(\mcP_{\whU_3}-\mcP_{U_3}\right)\mcP_{U_{3\perp}} \otimes \mcP_{U_2}\left(\mcP_{\whU_2}-\mcP_{U_2}\right)\mcP_{U_2}\right)A_1\mcP_{U_{1\perp}}\left(\mcP_{\whU_1}-\mcP_{U_1}\right)U_1G_1\left(U_3\otimes U_2\right)^{\top}\right]\right|}_{\mathrm{\RN{7}}} \\
+ & \underbrace{\left|\operatorname{tr}\left[\left(\mcP_{U_3}\left(\mcP_{\whU_3}-\mcP_{U_3}\right)\mcP_{U_{3\perp}} \otimes \mcP_{U_2}\left(\mcP_{\whU_2}-\mcP_{U_2}\right)\mcP_{U_{2\perp}}\right)A_1\mcP_{U_{1\perp}}\left(\mcP_{\whU_1}-\mcP_{U_1}\right)U_1G_1\left(U_3\otimes U_2\right)^{\top}\right]\right|}_{\mathrm{\RN{8}}}.
\end{align*}

Here, we have
\begin{align}
\mathrm{\RN{1}}
\leq & \underbrace{\left\|\mcP_{U_3}\left(\mcP_{\whU_3}-\mcP_{U_3}\right) \mcP_{U_3}\right\|}_{\eqref{eq: high-prob upper bound of PU1(PUhat1-PU1)PU1 in tensor regression}}  \underbrace{\left\|\mcP_{U_2}\left(\mcP_{\whU_2}-\mcP_{U_2}\right) \mcP_{U_2}\right\|}_{\eqref{eq: high-prob upper bound of PU1(PUhat1-PU1)PU1 in tensor regression}} 
\left\|\left(\mcP_{U_3} \otimes \mcP_{U_2}\right) A_1^{\top} \mcP_{U_1}\right\|_{\mathrm{F}}  \underbrace{\left\|P_{U_1}\left(\mcP_{\whU_1}-\mcP_{U_1}\right) U_1 G_1\left(U_3 \otimes U_2\right)^{\top}\right\|_{\mathrm{F}}}_{\eqref{eq: high-prob upper bound of PU1(PUhat1-PU1)PU1 in tensor regression}} \notag \\
\lesssim & \left\|\mcA \times_1 U_1 \times_2 U_2 \times_3 U_3\right\|_{\mathrm{F}} \cdot \frac{\sigma_\xi^6\oor^{1/2}}{\ulambda^5\sigma^6} \cdot \frac{\op^3}{n^3}. \label{eq: upper bound of term 1 in step 2.3 in tensor regression without sample splitting}
\end{align}

Then consider
\begin{align}
\mathrm{\RN{2}}
\leq & \underbrace{\left\|P_{U_1}\left(\mcP_{\whU_1}-\mcP_{U_1}\right) \mcP_{U_3}\right\|}_{\eqref{eq: high-prob upper bound of PU1(PUhat1-PU1)PU1 in tensor regression}}  \underbrace{\left\|P_{U_1}\left(\mcP_{\whU_3}-\mcP_{U_3}\right) \mcP_{U_3}\right\|}_{\eqref{eq: high-prob upper bound of PU1(PUhat1-PU1)PU1 in tensor regression}} 
\underbrace{\left\|\left(\mcP_{U_1} \otimes \mcP_{U_3}\right)A_2^{\top}\mcP_{U_{2\perp}}\left(\mcP_{\whU_2}-\mcP_{U_2}\right)\mcP_{U_2}\cdot U_2G_2\left(U_1\otimes U_3\right)^{\top}\right\|}_{\eqref{eq: high-prob upper bound of V1tPU1p(PUhat1-PU1)U1 in tensor regression without sample splitting}} \notag \\
\lesssim & \left\|\mcA \times_1 U_1 \times_3 U_3\right\|_{\mathrm{F}} \cdot \left[\frac{\sigma_\xi^3\oor^{1/2}}{\ulambda^2\sigma^3}\cdot \left(\frac{\op\sqrt{\oor\log(\op)}}{n^{3/2}} + \Delta\cdot \frac{\op^{3/2}}{n^{3/2}}\right)\right]. \label{eq: upper bound of term 2 in step 2.3 in tensor regression without sample splitting}
\end{align}

By symmetry, we have
\begin{align}
\mathrm{\RN{3}}
\lesssim & \left\|\mcA \times_1 U_1 \times_2 U_2\right\|_{\mathrm{F}} \cdot \left[\frac{\sigma_\xi^3\oor^{1/2}}{\ulambda^2\sigma^3}\cdot \left(\frac{\op\sqrt{\oor\log(\op)}}{n^{3/2}} + \Delta\cdot \frac{\op^{3/2}}{n^{3/2}}\right)\right], \label{eq: upper bound of term 3 in step 2.3 in tensor regression without sample splitting}
\end{align}
and
\begin{align}
\mathrm{\RN{5}}
\lesssim & \left\|\mcA \times_2 U_2 \times_3 U_3 \right\|_{\mathrm{F}} \cdot \left[\frac{\sigma_\xi^3\oor^{1/2}}{\ulambda^2\sigma^3}\cdot \left(\frac{\op\sqrt{\oor\log(\op)}}{n^{3/2}} + \Delta\cdot \frac{\op^{3/2}}{n^{3/2}}\right)\right]. \label{eq: upper bound of term 5 in step 2.3 in tensor regression without sample splitting}
\end{align}

Then, consider
\begin{align}
\mathrm{\RN{4}} 
\leq & \underbrace{\left\|\left(\mcP_{U_3}\left(\mcP_{\whU_3}-\mcP_{U_3}\right) \mcP_{U_{3 \perp}} \otimes \mcP_{U_2}\left(\mcP_{\whU_2}-\mcP_{U_2}\right) \mcP_{U_{2 \perp}}\right)A_1 \mcP_{U_1}\right\|_{\mathrm{F}}}_{\eqref{eq: high-prob upper bound of AoPU1oPU2p(PUhat2-P_U2)PU2oPU3p(PUhat3-P_U3)PU3 in tensor regression without sample splitting}} \cdot \underbrace{\left\|\mcP_{U_1}\left(\mcP_{\whU_1}-\mcP_{U_1}\right) U_1 G_1\left(U_3 \otimes U_2\right)^{\top}\right\|_{\mathrm{F}}}_{\eqref{eq: high-prob upper bound of PU1(PUhat1-PU1)PU1 in tensor regression}} \notag \\
\lesssim & \left\|\mcA \times_1 U_1\right\|_{\mathrm{F}} \cdot \left[\frac{\sigma_{\xi}^4\oor^{1/2}}{\ulambda^3 \sigma^4} \cdot\left(\frac{\oor\op\log (\op)}{n^2}+\Delta \cdot \frac{\op^{3/2}\sqrt{\oR\log(\op)}}{n^2}+\Delta^2 \cdot \frac{\op^2}{n^2}\right)\right]. \label{eq: upper bound of term 4 in step 2.3 in tensor regression without sample splitting}
\end{align}

By symmetry, we have
\begin{align}
\mathrm{\RN{6}} 
\lesssim & \left\|\mcA \times_3 U_3 \right\|_{\mathrm{F}} \cdot \left[\frac{\sigma_{\xi}^4\oor^{1/2}}{\ulambda^3 \sigma^4} \cdot\left(\frac{\oor\op\log (\op)}{n^2}+\Delta \cdot \frac{\op^{3/2}\sqrt{\oR\log(\op)}}{n^2}+\Delta^2 \cdot \frac{\op^2}{n^2}\right)\right], \label{eq: upper bound of term 6 in step 2.3 in tensor regression without sample splitting}
\end{align}
and
\begin{align}
\mathrm{\RN{7}} 
\lesssim & \left\|\mcA \times_2 U_2\right\|_{\mathrm{F}} \cdot \left[\frac{\sigma_{\xi}^4\oor^{1/2}}{\ulambda^3 \sigma^4} \cdot\left(\frac{\oor\op\log (\op)}{n^2}+\Delta \cdot \frac{\op^{3/2}\sqrt{\oR\log(\op)}}{n^2}+\Delta^2 \cdot \frac{\op^2}{n^2}\right)\right]. \label{eq: upper bound of term 7 in step 2.3 in tensor regression without sample splitting}
\end{align}

Finally, we have
\begin{align}
& \mathrm{\RN{8}} \notag \\
\leq & \underbrace{\left\|\left(\mcP_{U_3}\left(\mcP_{\whU_3}-\mcP_{U_3}\right) \mcP_{U_{3 \perp}} \otimes \mcP_{U_2}\left(\mcP_{\whU_2}-\mcP_{U_2}\right) \mcP_{U_{2 \perp}}\right) A_1 \mcP_{U_{1 \perp}}\left(\mcP_{\whU_1}-\mcP_{U_1}\right)\mcP_1^{1/2}\right\|_{\mathrm{F}}}_{\eqref{eq: high-prob upper bound of AoPU1p(PUhat1-P_U1)PU1oPU2p(PUhat2-P_U2)PU2oPU3p(PUhat3-P_U3)PU3 in tensor regression without sample splitting}}  \left\|\mcP_1^{-1/2} U_1 G_1\left(U_3 \otimes U_2\right)^{\top}\right\|_{\mathrm{F}} \notag \\
\lesssim & \left\|\mcA\right\|_{\mathrm{F}} \cdot \left[\frac{\sigma_\xi^3\oor^{1/2} }{\ulambda^2\sigma^3}\cdot \left(\frac{\oor^{3/2}\log(\op)^{3/2}}{n^{3/2}} + \Delta\cdot \frac{\op^{1/2}\oR\log(\op)}{n^{3/2}} + \Delta^2 \cdot \frac{\op \oR \log (\op)}{n^{3 / 2}} + \Delta^3\cdot \frac{\op^{3/2}}{n^{3/2}}\right)\right] \label{eq: upper bound of term 8 in step 2.3 in tensor regression without sample splitting}.
\end{align}

Then, we have
\begin{align*}
& \left|\left\langle\mcT \times_1\left(\mcP_{\whU_1}-\mcP_{U_1}\right) \times_2\left(\mcP_{\whU_2}-\mcP_{U_2}\right) \times_3\left(\mcP_{\whU_3}-\mcP_{U_3}\right), \mcA\right\rangle\right| \\
\lesssim & \eqref{eq: upper bound of term 1 in step 2.3 in tensor regression without sample splitting} + \eqref{eq: upper bound of term 2 in step 2.3 in tensor regression without sample splitting} + \eqref{eq: upper bound of term 3 in step 2.3 in tensor regression without sample splitting} + \eqref{eq: upper bound of term 4 in step 2.3 in tensor regression without sample splitting} + \eqref{eq: upper bound of term 5 in step 2.3 in tensor regression without sample splitting}+ \eqref{eq: upper bound of term 6 in step 2.3 in tensor regression without sample splitting} + \eqref{eq: upper bound of term 7 in step 2.3 in tensor regression without sample splitting} + \eqref{eq: upper bound of term 8 in step 2.3 in tensor regression without sample splitting}\\
\lesssim & \left\|\mcA \times_1 U_1 \times_2 U_2 \times_3 U_3\right\|_{\mathrm{F}} \cdot \frac{\sigma_{\xi}^6 \oor^{1 / 2}}{\ulambda^5 \sigma^6} \cdot \frac{\op^3}{n^3} 
+  \sum_{j=1}^3 \left\|\mcA \times_{j+1} U_{j+1} \times_{j+2} U_{j+2}\right\|_{\mathrm{F}} \cdot \left[\frac{\sigma_{\xi}^3 \oor^{1 / 2}}{\ulambda^2 \sigma^3} \cdot\left(\frac{\oor \log (\op)^{1 / 2}}{n^{3 / 2}}+\Delta \cdot \frac{\op^{3 / 2}}{n^{3 / 2}}\right)\right] \\
+ & \sum_{j=1}^3 \left\|\mcA \times_j U_j\right\|_{\mathrm{F}} \cdot \left[\frac{\sigma_{\xi}^4 \oor^{1 / 2}}{\ulambda^3 \sigma^4} \cdot\left(\frac{\oor\op\log (\op)}{n^2}+\Delta \cdot \frac{\op^{3 / 2} \sqrt{\oR \log (\op)}}{n^2}+\Delta^2 \cdot \frac{\op^2}{n^2}\right)\right] \\
+ & \left\|\mcA\right\|_{\mathrm{F}} \cdot \left[\frac{\sigma_\xi^3\oor^{1/2} }{\ulambda^2\sigma^3}\cdot \left(\frac{\oor^{3/2}\log(\op)^{3/2}}{n^{3/2}} + \Delta\cdot \frac{\op^{1/2}\oR\log(\op)}{n^{3/2}} + \Delta^2 \cdot \frac{\op \oR \log (\op)}{n^{3 / 2}} + \Delta^3\cdot \frac{\op^{3/2}}{n^{3/2}}\right)\right].
\end{align*}

\subsection*{Step 3: Analysis of asymptotic normal terms} 

Recall that 
$$
\whZ_j=\underbrace{\frac{1}{n\sigma^2} \sum_{i=1}^n \xi_i \operatorname{Mat}_j\left(\mcX_i\right)}_{\widehat{\mcZ}_j^{(1)}} + \underbrace{\frac{1}{n\sigma^2}\sum_{i=1}^n \left\langle \Delta, \mcX_i \right\rangle \operatorname{Mat}_j\left(\mcX_i\right) - \sigma^2\cdot \operatorname{Mat}_j\left(\Delta\right)}_{\widehat{\mcZ}_j^{(1)}}.
$$
Let $\mcP_{\mathbb{T}_{\mcT} \mathcal{M}_{(r_1, r_2, r_3)}}\left(\mcA\right) = \mcA \times_1 \mcP_{U_1} \times_2 \mcP_{U_2} \times_3 \mcP_{U_3} + \sum_{j=1}^{3} \operatorname{Mat}_j^{-1}\left(\mcP_{U_{j\perp}}A_j\mcP_{\left(U_{j+2} \otimes U_{j+1}\right)G_j}\right)$. Then, we can write
\begin{align*}
\left\langle \widehat{\mcZ}, \mcP_{\mathbb{T}_{\mcT} \mathcal{M}_{(r_1, r_2, r_3)}}\left(\mcA\right) \right\rangle = & \left\langle \widehat{\mcZ}, \mcA \times_1 \mcP_{U_1} \times_2 \mcP_{U_2} \times_3 \mcP_{U_3} \right\rangle + \sum_{j=1}^{3}\left\langle  \whZ_j, \mcP_{U_{j\perp}}A_j\mcP_{\left(U_{j+2} \otimes U_{j+1}\right)G_j} \right\rangle \\
= & \left\langle \widehat{\mcZ}^{(1)}, \mcP_{\mathbb{T}_{\mcT} \mathcal{M}_{(r_1, r_2, r_3)}}\left(\mcA\right) \right\rangle + \left\langle \widehat{\mcZ}^{(2)}, \mcP_{\mathbb{T}_{\mcT} \mathcal{M}_{(r_1, r_2, r_3)}}\left(\mcA\right)\right\rangle.
\end{align*}

\subsubsection*{Step 3.1: Asymptotic Normality of 
$\langle \widehat{\mcZ}^{(1)}, \mcP_{\mathbb{T}_{\mcT} \mathcal{M}_{(r_1, r_2, r_3)}}\left(\mcA\right)\rangle$}

Note that $\langle \widehat{\mcZ}^{(1)}, \mcP_{\mathbb{T}_{\mcT} \mathcal{M}_{(r_1, r_2, r_3)}}\left(\mcA\right)\rangle$, for any $j=1,2,3$, is a summation of i.i.d. random variables: 
\begin{align*}
\left\langle \widehat{\mcZ}^{(1)}, \mcP_{\mathbb{T}_{\mcT} \mathcal{M}_{(r_1, r_2, r_3)}}\left(\mcA\right) \right\rangle
= \frac{1}{n\sigma^2} \left\langle \sum_{i=1}^{n} \xi_i \mcX_i, \mcP_{\mathbb{T}_{\mcT} \mathcal{M}_{(r_1, r_2, r_3)}}\left(\mcA\right) \right\rangle.
\end{align*}
To apply the Berry-Essen theorem, we calculate its second and third moments. Let $\mcX_i$ be i.i.d. copies of $\mcX$. Clearly,
\begin{align*}
\mathbb{E} \xi^2 \left[\left\langle \mcX, \mcP_{\mathbb{T}_{\mcT} \mathcal{M}_{\mathbf{r}}}\left(\mcA\right) \right\rangle\right]^2 
= & \sigma_{\xi}^2\sigma^2 \cdot \sum_{j=1}^3 \left\| \mcP_{U_{j\perp}} A_j \mcP_{\left(U_{j+2}\otimes U_{j+1}\right)G_j^{\top}} \right\|_{\mathrm{F}}^2 + \sigma_{\xi}^2\sigma^2 \cdot \left\|\mcA \times_1 \mcP_{U_1} \times_2 \mcP_{U_2} \times_3 \mcP_{U_3}\right\|_{\mathrm{F}}^2 .
\end{align*}
Next, we bound the third moment. Clearly,
\begin{align*}
& \left|\left\langle \operatorname{Mat}_j\left(\mcX\right),  \mcP_{U_{j\perp}}A_j\mcP_{\left(U_{j+2} \otimes U_{j+1}\right)G_j^{\top}} \right\rangle \right|
 \leq \left\|\mcP_{U_{j\perp}}A_j\mcP_{\left(U_{j+2} \otimes U_{j+1}\right)G_j^{\top}} \right\|_{\mathrm{F}}\cdot \left|\left\langle  X , \mcA_k\right\rangle\right|,
\end{align*}
where $\mcA_j$ is a tensor with unit Frobenius norm defined by normalizing $\mcP_{U_{j\perp}}A_j\mcP_{\left(U_{j+2} \otimes U_{j+1}\right)G_j^{\top}}$. By the sub-Gaussian assumption on $\mcX$, we know that $\left\langle  \mcX, \mcA_k \right\rangle$ is a sub-Gaussian variable with parameter $\sigma^2$. Therefore, we have 
\begin{align*}
 \mathbb{E}\left|\left\langle \operatorname{Mat}_j\left(\mcX\right),  \mcP_{U_{j\perp}}A_j\mcP_{\left(U_{j+2} \otimes U_{j+1}\right)G_j^{\top}} \right\rangle \right|^{k}
\leq & C_k \sigma^k \cdot \left\|\mcP_{U_{j\perp}}A_j\mcP_{\left(U_{j+2} \otimes U_{j+1}\right)G_j^{\top}} \right\|_{\mathrm{F}}^{k} .
\end{align*}
Similarly, we have
\begin{align*}
\mathbb{E}\left|\left\langle \mcX, \mcA\times_1 \mcP_{U_1}\times_2 \mcP_{U_2}\times_3 \mcP_{U_3}\right\rangle\right|^k 
\leq & C_k \sigma^k \cdot \left\|\mcA\times_1 \mcP_{U_1}\times_2 \mcP_{U_2}\times_3 \mcP_{U_3} \right\|_{\mathrm{F}}^{k} .
\end{align*}
In addition, note that the correlation between any two terms among $\langle X, \mcA\times_1 \mcP_{U_1}\times_2 \mcP_{U_2}\times_3 \mcP_{U_3}\rangle$ and $\langle \operatorname{Mat}_j(\mcX),  \mcP_{U_{j\perp}}A_j\mcP_{(U_{j+2} \otimes U_{j+1})G_j^{\top}} \rangle$ are uncorrelated. Therefore, by the same argument, we have 
\begin{align*}
\mathbb{E}\left|\xi\right|^3 \left|\left\langle \mcX, \mcP_{\mathbb{T}_{\mcT} \mathcal{M}_{\mathbf{r}}}\left(\mcA\right)\right\rangle\right|^3
\leq & C_3 \sigma_\xi^3 \sigma^3 \cdot \left(\sum_{j=1}^3\left\|\mcP_{U_{j\perp}} A_j\mcP_{\left(U_{j+2} \otimes U_{j+1}\right)G_j^{\top}}\right\|_{\mathrm{F}}+\left\|\mcA\times_1 U_1\times_2 U_2\times_3 U_3\right\|_{\mathrm{F}}\right)^3 .
\end{align*}
By Berry-Essen theorem \citep{berry1941accuracy, esseen1956moment} and Theorem 3.7 of \citet{chen2010normal}, we get
\begin{align*}
&\sup_{x \in \mathbb{R}} \left| \mathbb{P}  \left(\frac{\left\langle \widehat{\mcZ}^{(1)}, \mcP_{\mathbb{T}_{\mcT} \mathcal{M}_{(r_1, r_2, r_3)}}\left(\mcA\right)\right\rangle}{\frac{\sigma_{\xi}}{\sigma}\left(\sum_{j=1}^{3}\left\|\mcP_{U_{j\perp}}A_j\mcP_{\left(U_{j+2} \otimes U_{j+1}\right)G_j^{\top}} \right\|_{\mathrm{F}}^2+ \left\|\mcA \times_1 \mcP_{U_1} \times_2 \mcP_{U_2} \times_3 \mcP_{U_3}\right\|_{\mathrm{F}}^2\right)^{1/2} \cdot \sqrt{\frac{1}{n}}} \leq x\right)-\Phi(x) \right| \\
\lesssim &  \sqrt{\frac{1}{n}} \cdot \frac{\left(\sum_{j=1}^3\left\|\mcP_{U_{j\perp}} A_j\mcP_{\left(U_{j+2} \otimes U_{j+1}\right)G_j^{\top}}\right\|_{\mathrm{F}}+\left\|\mcA\times_1 U_1\times_2 U_2\times_3 U_3\right\|_{\mathrm{F}}\right)^3}{\left(\sum_{j=1}^3\left\|\mcP_{U_{j\perp}} A_j\mcP_{\left(U_{j+2} \otimes U_{j+1}\right)G_j^{\top}}\right\|_{\mathrm{F}}^2+\left\|\mcA\times_1 U_1\times_2 U_2\times_3 U_3\right\|_{\mathrm{F}}^2\right)^{\frac{3} {2}}} \lesssim \sqrt{\frac{1}{n}}.
\end{align*}

\subsubsection*{Step 3.2: Upper Bound of 
$ \sum_{j=1}^{3}\langle  \widehat{\mcZ}^{(2)}_j, \mcP_{U_j\perp}A_j\mcP_{(U_{j+2} \otimes U_{j+1})G_j^{\top}} \rangle +\langle \widehat{\mcZ}^{(2)}, \mcA\times_1 \mcP_{U_1}\times_2 \mcP_{U_2}\times_3 \mcP_{U_3}\rangle$}

Note that 
\begin{align*}
\left|\left\langle \widehat{\mcZ}^{(2)}, \mcP_{\mathbb{T}_{\mcT} \mathcal{M}_{\mathbf{r}}}\left(\mcA\right)\right\rangle\right|
\leq & \left(\sum_{j=1}^3\left\|\mcP_{U_{j\perp}}A_j\mcP_{\left(U_{j+2}\otimes U_{j+1}\right)G_j^{\top}}\right\|_{\mathrm{F}} + \left\|\mcA\times_1 U_1\times_2 U_2\times_3 U_3\right\|_{\mathrm{F}}\right) \cdot \Delta \cdot \frac{\sigma_\xi\oor^{1/2}}{\sigma} \sqrt{\frac{\op}{n}} .
\end{align*}

\subsubsection*{Step 3.3: Combining Asymptotic Normal Terms and Negligible Terms}

By the Lipschitz property of normal distribution 
function $\Phi(x)$ and note that the discussion above holds under event $\left\|\widehat{\mcT}^{\text{init}} - \mcT\right\|_{\mathrm{F}} \leq \Delta$ and $\left\|\mcP_{\whU_j} - \mcP\right\| \leq \frac{\sigma_\xi}{\sigma}\sqrt{\frac{\op}{n}}$, then finally we have
\begin{align*}
& \sup _{x \in \mathbb{R}} \left\lvert\, \mathbb{P}\left(\frac{\left\langle\widehat{\mcT}, \mcA\right\rangle-\langle T, \mcA\rangle}{\frac{\sigma_{\xi}}{\sigma} \cdot \left(\sum_{j=1}^3 \left\|\mcP_{U_{j \perp}} A_j \mcP_{\left(U_{j+2} \otimes U_{j+1}\right)G_j^{\top}} \right\|_{\mathrm{F}}^2+ \left\|\mcA \times_1 \mcP_{U_1} \times_2 \mcP_{U_2} \times_3 \mcP_{U_3} \right\|_{\mathrm{F}}^2\right)^{1/2}\cdot \sqrt{\frac{1}{n}}} \leq x\right)-\Phi(x)\right| \notag\\
\leq & \sqrt{\frac{1}{n}} + \underbrace{\left[\frac{1}{\op^c} + \exp\left(-c\op\right) + \exp\left(-cn\right) + \mcP\left(\mcE_\Delta\right) + \mcP\left(\mcE_U^{\text{reg}}\right)\right]}_{\text{rate of initial estimate}} \\
+ & \frac{1}{\frac{\sigma_{\xi}}{\sigma} \cdot s_{\mcA} \cdot \sqrt{\frac{1}{n}}} \cdot \Bigg\{ \underbrace{\left\|\mcA\times_1 U_1 \times_2 U_2 \times_3 U_3\right\|_{\mathrm{F}} \cdot \left[\frac{\sigma_\xi^2\oor^{1/2}}{\ulambda\sigma^2}\cdot \left(\frac{\sqrt{\oor\op\log(\op)}}{n} +\Delta\cdot \frac{\op}{n}\right)\right]}_{\text{from Step 1}} \\
+ & \sum_{j=1}^3 \left\|U_j^{\top} A_j \mcP_{\left(U_{j+2} \otimes U_{j+1}\right) G_j^{\top}}\right\|_{\mathrm{F}} \cdot \frac{\sigma_\xi^2\oor^{1/2}}{\ulambda\sigma^2}\cdot \frac{\op}{n} \\
+ & \underbrace{\sum_{j=1}^3\left\|\mcA \times_{j+1} U_{j+1} \times_{j+2} U_{j+2}\right\|_{\mathrm{F}} \cdot \left[\frac{\sigma_\xi^2\oor^{1/2}}{\ulambda\sigma^2}\left(\frac{\oor\log(\op)}{n}+\Delta^2\cdot \frac{\op}{n}\right) + \frac{\sigma_\xi^3\oor^{1/2}}{\ulambda^2\sigma^3}\left(\frac{\op\sqrt{\oor\log(\op)}}{n^{3/2}}\right)\right] 
}_{\text{shared between Step 1 and Step 2}}  \\
+ & \sum_{j=1}^3 \left\|\mcP_{U_{j\perp}} A_j \mcP_{\left(U_{j+2} \otimes U_{j+1}\right) G_j^{\top}}\right\|_{\mathrm{F}} \cdot\left[\frac{\sigma_{\xi}^2 \oor^{1/2}}{\ulambda \sigma^2} \cdot\left(\frac{\sqrt{\oor\op \log (\op)}}{n}+\Delta \cdot \frac{\op}{n}\right)\right] \\
+ & \sum_{j=1}^3\left\|\mcA \times_j U_j\right\|_{\mathrm{F}} \cdot\left[\frac{\sigma_{\xi}^2 \oor^{1 / 2}}{\ulambda \sigma^2} \cdot\left(\frac{\oor\log(\op)}{n}+\Delta \cdot \frac{\sqrt{\oR \op \log (\op)}}{n}+\Delta^2 \cdot \frac{\op}{n}\right)\right] \\
+ & \underbrace{\left\|\mcA\right\|_{\mathrm{F}} \cdot\left[\frac{\sigma_{\xi}^3 \oor^{1 / 2}}{\ulambda^2 \sigma^3} \cdot\left(\frac{\oor^{3/2}\log(\op)^{3/2}}{n^{3 / 2}}+\Delta \cdot \frac{\op^{1 / 2} \oR \log (\op)}{n^{3 / 2}}+\Delta^2 \cdot \frac{\op \oR \log (\op)}{n^{3 / 2}}+\Delta^3 \cdot \frac{\op^{3 / 2}}{n^{3 / 2}}\right)\right]}_{\text{from Step 2}} \\
+& \underbrace{\sum_{j=1}^3\left\|\mcP_{U_{j\perp}}A_j\mcP_{\left(U_{j+2}\otimes U_{j+1}\right)G_j^{\top}}\right\|_{\mathrm{F}}\cdot \Delta \cdot \frac{\sigma_\xi\oor^{1/2}}{\sigma} \sqrt{\frac{\op}{n}} + \left\|\mcA\times_1 U_1\times_2 U_2\times_3 U_3\right\|_{\mathrm{F}}\cdot \Delta \cdot \frac{\sigma_\xi\oor^{1/2}}{\sigma} \sqrt{\frac{\op}{n}}}_{\text{from Step 3}} \Bigg\} .
\end{align*}


\section{Proof of Theorem~\ref{thm: main theorem in tensor regression with sample splitting}} \label{sec: proof of main theorem in tensor regression with sample splitting}

In this section, we present the proof of Theorem~\ref{thm: main theorem in tensor regression with sample splitting}. Since the proof of Theorem~\ref{thm: main theorem in tensor regression with sample splitting} is similar to the proof of the Theorem~\ref{thm: main theorem in tensor regression without sample splitting}, we will focus on the parts that differ. For identical or repetitive steps, such as the decomposition of certain terms, we will provide a concise description to maintain textual conciseness. 

First,  for any $j=1,2,3$, we assume that the following events hold with high probability:
$
\|\mcP_{\whU_j^{(0)}} - \mcP_{U_j}\| \leq \frac{\sigma_\xi}{\sigma}\sqrt{\frac{\op}{n}}
$
holds with probability at least $1- \mathbb{P}\left(\mcE_{U}^{\text{reg}}\right)$, where event $\mcE_{U}^{\text{reg}}$ is defined by $\mcE_{U}^{\text{reg}}=  \{\max_{j=1,2,3}\|\mcP_{\whU_j^{(0),(\mathrm{\RN{1}})}} - \mcP_{U_j} \| > \frac{\sigma_\xi}{\sigma}\sqrt{\frac{\op}{n}} \}$.
Then by Lemma~\ref{lemma: error contraction of l2 error of singular space in tensor regression}, we know that 
$
\|\mcP_{\whU_j} - \mcP_{U_j}\|= \|\mcP_{\whU_j^{(1)}} - \mcP_{U_j}\| \leq \frac{\sigma_\xi}{\sigma}\sqrt{\frac{\op}{n}}
$
holds with probability at least $1-\exp(-c\op) - \mathbb{P}\left(\mcE_{U}^{\text{reg}}\right)$ for any $j=1,2,3$.

Besides, we assume that the following initial error bound 
$
\|\whT^{\text{init}} - \mcT\|_{\mathrm{F}} \leq \Delta
$
holds with probability at least $1-\mathbb{P}(\mcE_{\Delta})$, where event $\mcE_{\Delta}$ is given by $\mcE_{\Delta}= \{\|\whT^{\text{init}} - \mcT\|_{\mathrm{F}} > \Delta \}$.

Denote $\widehat{\Delta}=\mcT-\widehat{\mcT}^{\text{init}}$ . Recall that the estimator is of the following form
\begin{align*}
&\widehat{\mcT}-\mcT\\
=&   \underbrace{\frac{1}{n \sigma^2} \sum_{i_1=1}^{n_1} \xi_{i_1}^{(\mathrm{\RN{1}})} \mcX_{i_1}^{(\mathrm{\RN{1}})}}_{\widehat{\mcZ}^{(1),(\mathrm{\RN{1}})}}+\underbrace{\left(\frac{1}{n \sigma^2} \sum_{i_1=1}^{n_1}\left\langle\widehat{\Delta}^{(\mathrm{\RN{2}})}, \mcX_{i_1}^{(\mathrm{\RN{1}})}\right\rangle \mcX_{i_1}-\widehat{\Delta}^{(\mathrm{\RN{2}})}\right)}_{\widehat{\mcZ}^{(2),(\mathrm{\RN{1}})}} 
 + \underbrace{\frac{1}{n \sigma^2} \sum_{i_2=1}^{n_2} \xi_{i_2}^{(\mathrm{\RN{2}})} \mcX_{i_2}^{(\mathrm{\RN{2}})}}_{\widehat{\mcZ}^{(1),(\mathrm{\RN{2}})}}+\underbrace{\left(\frac{1}{n \sigma^2} \sum_{i_2=1}^{n_2}\left\langle\widehat{\Delta}^{(\mathrm{\RN{1}})}, \mcX_{i_2}^{(\mathrm{\RN{2}})}\right\rangle \mcX_{i_2}-\widehat{\Delta}^{(\mathrm{\RN{1}})}\right)}_{\widehat{\mcZ}^{(2),(\mathrm{\RN{2}})}}\\
=& \mcT+\widehat{\mcZ}^{(\mathrm{\RN{1}})}+\widehat{\mcZ}^{(\mathrm{\RN{2}})} = \mcT+\widehat{\mcZ},
\end{align*}
where $\widehat{\mcZ}=\widehat{\mcZ}^{(\mathrm{\RN{1}})}+\widehat{\mcZ}^{(\mathrm{\RN{2}})}= \widehat{\mcZ}^{(1),(\mathrm{\RN{1}})}+\widehat{\mcZ}^{(2),(\mathrm{\RN{1}})} + \widehat{\mcZ}^{(1),(\mathrm{\RN{2}})}+\widehat{\mcZ}^{(2),(\mathrm{\RN{2}})}$.
Then
\begin{align*}
\left\langle\widehat{\mcT}-\mcT, \mcA \right\rangle 
= & \left\langle \frac{n_1}{n} \widehat{\mcT}^{\text {unbs,(\RN{1}) }} \times_1 \mcP_{\whU_1^{(\mathrm{\RN{1}})}} \times_2 \mcP_{\whU_2^{(\mathrm{\RN{1}})}} \times_3 \mcP_{\whU_3^{(\mathrm{\RN{1}})}} + \frac{n_2}{n} \widehat{\mcT}^{\text {unbs,(\RN{2})}} \times_1 \mcP_{\whU_1^{(\mathrm{\RN{2}})}} \times_2 \mcP_{\whU_2^{(\mathrm{\RN{2}})}} \times_3 \mcP_{\whU_3^{(\mathrm{\RN{2}})}} -\mcT, \mcA \right\rangle \\
= & \left\langle\widehat{\mcZ} , \mcA\right\rangle = \left\langle\widehat{\mcZ}^{(\mathrm{\RN{1}})} , \mcA\right\rangle + \left\langle\widehat{\mcZ}^{(\mathrm{\RN{2}})} , \mcA\right\rangle.
\end{align*}

We then fulfill the proof of our theorem in details.

\subsection*{Step 1: Upper Bound of Negligible Terms in $\langle \widehat{\mcZ}^{(\mathrm{\RN{1}})} \times_1 \mcP_{\widehat{U}_1^{(\mathrm{\RN{1}})}} \times_2 \mcP_{\widehat{U}_2^{(\mathrm{\RN{1}})}} \times_3 \mcP_{\widehat{U}_3^{(\mathrm{\RN{1}})}} \rangle$}

Similar to the arguments in Step 1 in the proof of Theorem~\ref{thm: main theorem in tensor regression without sample splitting}, by symmetry, it suffices to consider 
\begin{align}
& \left|\left\langle\whZ^{(\mathrm{\RN{1}})} \times_1 \left(\mcP_{\whU_1^{(\mathrm{\RN{1}})}}-\mcP_{U_1}\right) \times_2 \mcP_{U_2} \times_3 \mcP_{U_3}, \mcA\right\rangle\right|, \label{eq: step 1.1 in tensor regression with sample splitting} \\
& \left|\left\langle\whZ^{(\mathrm{\RN{1}})} \times_1 \left(\mcP_{\whU_1^{(\mathrm{\RN{1}})}}-\mcP_{U_1}\right) \times_2 \left(\mcP_{\whU_2^{(\mathrm{\RN{1}})}}-\mcP_{U_2}\right) \times_3 \mcP_{U_3}, \mcA\right\rangle\right|, \label{eq: step 1.2 in tensor regression with sample splitting}\\
& \left|\left\langle\whZ^{(\mathrm{\RN{1}})} \times_1 \left(\mcP_{\whU_1^{(\mathrm{\RN{1}})}} -\mcP_{U_1}\right) \times_2 \left(\mcP_{\whU_2^{(\mathrm{\RN{1}})}} -\mcP_{U_2}\right)\times_3 \left(\mcP_{\whU_3^{(\mathrm{\RN{1}})}} -\mcP_{U_3}\right), \mcA\right\rangle\right|. \label{eq: step 1.3 in tensor regression with sample splitting}
\end{align}

\subsubsection*{Step 1.1: Upper Bound of Negligible Terms in $\langle \widehat{\mcZ}^{(\mathrm{\RN{1}})} \times_1 (\mcP_{\widehat{U}_1^{(\mathrm{\RN{1}})}} - \mcP_{U_1}) \times_2 \mcP_{U_2} \times_3 \mcP_{U_3}, \mcA \rangle$} 

First, consider
\begin{align}
\left|\left\langle \widehat{\mcZ}^{(\mathrm{\RN{1}})} \times_1 \left(\mcP_{\whU_1^{(\mathrm{\RN{1}})}} - \mcP_{U_1}\right) \times_2 \mcP_{U_2} \times_3 \mcP_{U_3}, \mcA\right\rangle\right| 
\leq & \mathrm{\RN{1}} + \mathrm{\RN{2}} + \mathrm{\RN{3}} + \mathrm{\RN{4}} \label{eq: decomposition of step 1.1 in tensor regression with sample splitting}.
\end{align}
Here, $\mathrm{\RN{1}}, \mathrm{\RN{2}}, \mathrm{\RN{3}}$, and $\mathrm{\RN{4}}$ are defined in the same manner as in \eqref{eq: term 1 in step 1.1 in tensor regression without sample splitting}, \eqref{eq: term 2 in step 1.1 in tensor regression without sample splitting}, \eqref{eq: term 3 in step 1.1 in tensor regression without sample splitting}, and \eqref{eq: term 4 in step 1.1 in tensor regression without sample splitting}, respectively, with $\widehat{Z}_1$ replaced by $\widehat{Z}_1^{(\mathrm{\RN{1}})}$.

We begin with the upper bound for the first term \RN{1} in \eqref{eq: decomposition of step 1.1 in tensor regression with sample splitting}: 
\begin{align}
\mathrm{\RN{1}}
\lesssim & \left\|\mcA \times_1 U_1 \times_2 U_2 \times_3 U_3\right\|_{\mathrm{F}} \cdot \left[\frac{\sigma_\xi^3\oor^{1/2}}{\ulambda^2\sigma^3}\cdot \left(\frac{\op\sqrt{\oor\log(\op)}}{n^{3/2}}\right)\right]. \label{eq: upper bound of term 1 in step 1.1 in tensor regression with sample splitting}
\end{align}

For the second term $\mathrm{\RN{2}}$ in \eqref{eq: decomposition of step 1.1 in tensor regression with sample splitting}, we have
\begin{align}
\mathrm{\RN{2}}
\lesssim & \left\|\mcA \times_2 U_2 \times_3 U_3\right\|_{\mathrm{F}} \cdot \left[\frac{\sigma_\xi^2\oor^{1/2}}{\ulambda^2\sigma^2}\left(\frac{\oor\log(\op)}{n}\right)\right]. \label{eq: upper bound of term 2 in step 1.1 in tensor regression with sample splitting}
\end{align}

For the third term $\mathrm{\RN{3}}$ in \eqref{eq: decomposition of step 1.1 in tensor regression with sample splitting}, we first have 
\begin{align*}
\mathrm{\RN{3}} 
= & \underbrace{\left|\operatorname{tr}\left[\left(\mcP_{U_3}\otimes \mcP_{U_2}\right)A_1^{\top}\mcP_{U_1} S_{G_1,1}\left(\whE_1^{(\mathrm{\RN{1}})}\right)\mcP_{U_{1\perp}}\whZ_1^{(\mathrm{\RN{1}})}\left(\mcP_{U_3}\otimes \mcP_{U_2}\right)\right]\right|}_{\mathrm{\RN{3}}.\mathrm{\RN{1}}} \\
+ & \underbrace{\left|\operatorname{tr}\left[\left(\mcP_{U_3}\otimes \mcP_{U_2}\right)A_1^{\top}\mcP_{U_1}\sum_{k_1=2}^{+\infty} S_{G_1,k_1}\left(\whE_1^{(\mathrm{\RN{1}})}\right)\mcP_{U_{1\perp}}\whZ_1^{(\mathrm{\RN{1}})}\left(\mcP_{U_3}\otimes \mcP_{U_2}\right)\right]\right|}_{\mathrm{\RN{3}}.\mathrm{\RN{2}}}.
\end{align*}
Note that $\whZ_1^{(\mathrm{\RN{1}})}=\underbrace{\frac{1}{n_1\sigma^2}\sum_{i_1=1}^{n_1}\xi_{i_1}^{(\mathrm{\RN{1}})} \Mat_1 (\mcX_{i_1}^{(\mathrm{\RN{1}})} ) }_{\whZ_1^{(1),(\mathrm{\RN{1}})}}+ \underbrace{\frac{1}{n_1\sigma^2}\sum_{i_1=1}^{n_1} [ \langle\mcX_{i_1}^{(\mathrm{\RN{1}})}, \widehat{\Delta}^{(\mathrm{\RN{2}})} \rangle\Mat_1 (\mcX_{i_1}^{(\mathrm{\RN{1}})} )-\sigma^2\cdot \widehat{\Delta}_1^{(\mathrm{\RN{2}})} ]}_{\whZ_1^{(2),(\mathrm{\RN{1}})} } $. By Lemma~\ref{lemma: high-prob upper bound of tr(BZ(1)tCZ(1))}, Lemma~\ref{lemma: high-prob upper bound of tr(BZ(2)tCZ(1)) in tensor regression with sample splitting} and Lemma~\ref{lemma: high-prob upper bound of tr(BZ(2)tCZ(2)) in tensor regression with sample splitting}, and by the same decomposition of $\mathrm{\RN{3}}.\mathrm{\RN{1}}$ in Step 1.1 in the proof of Theorem~\ref{thm: main theorem in tensor regression without sample splitting}, we have
\begin{align}
& \mathrm{\RN{3}}.\mathrm{\RN{1}} 
\lesssim  \left\|\mcP_{U_1}A_1\mcP_{\left(U_3\otimes U_2\right)G_1^{\top}}\right\|_{\mathrm{F}} \cdot \frac{\sigma_\xi^2\oor^{1/2}}{\ulambda\sigma^2}\cdot \frac{\op}{n} + \left\|\mcA\times_1 U_1 \times_2 U_2 \times_3 U_3\right\|_{\mathrm{F}} \cdot \frac{\sigma_\xi^2\oor^{1/2}}{\ulambda\sigma^2} \left(\frac{\sqrt{\oor\op\log(\op)}}{n} + \Delta^2 \cdot \frac{\log(\op)^{3/2}}{\sqrt{n}}\right). \label{eq: upper bound of the 3.1 term in step 1.1 in tensor regression with sample splitting}
\end{align}
Furthermore, applying the same arguments in the proof of \eqref{eq: upper bound of the 3.2 term in step 1.1 in tensor regression without sample splitting}, we have
\begin{align}
\mathrm{\RN{3}}.\mathrm{\RN{2}} 
\lesssim & \left\|\mcA\times_1 U_1 \times_2 U_2 \times_3 U_3\right\|_{\mathrm{F}} \cdot \left(\frac{\sigma_\xi^3\oor^{1/2}}{\ulambda^2\sigma^3}\cdot \frac{\op\sqrt{\oor\log(\op)}}{n^{3/2}} + \frac{\sigma_\xi^4\oor^{1/2}}{\ulambda^3\sigma^4}\cdot \frac{\op^2}{n^2}\right). \label{eq: upper bound of the 3.2 term in step 1.1 in tensor regression with sample splitting}
\end{align}
It implies that
\begin{align}
\mathrm{\RN{3}} 
\lesssim & \big\|\mcP_{U_1}A_1\mcP_{\left(U_3\otimes U_2\right)G_1^{\top}}\big\|_{\mathrm{F}} \frac{\sigma_\xi^2\oor^{1/2}}{\ulambda\sigma^2}\cdot \frac{\op}{n} \notag \\
+ & \|\mcA\times_1 U_1 \times_2 U_2 \times_3 U_3 \|_{\mathrm{F}}  \left[\frac{\sigma_\xi^2\oor^{1/2}}{\ulambda\sigma^2}  \left(\frac{\sqrt{\oor\op\log(\op)}}{n} +  \frac{\Delta^2\log(\op)^{3/2}}{\sqrt{n}}\right) + \frac{\sigma_\xi^4\oor^{1/2}}{\ulambda^3\sigma^4}\cdot \frac{\op^2}{n^2}\right]. \label{eq: upper bound of term 3 in step 1.1 in tensor regression with sample splitting}
\end{align}
For the fourth term $\mathrm{\RN{4}}$ in \eqref{eq: decomposition of step 1.1 in tensor regression with sample splitting}, we have
\begin{align}
\mathrm{\RN{4}} 
\lesssim & \left\|\mcA \times_2 U_2 \times_3 U_3\right\|_{\mathrm{F}} \cdot \left[\frac{\sigma_\xi^3\oor^{1/2}}{\ulambda\sigma^2}\left(\frac{\op\sqrt{\oor\log(\op)}}{n^{3/2}}\right)\right]. \label{eq: upper bound of term 4 in step 1.1 in tensor regression with sample splitting}
\end{align}

Therefore, we have
\begin{align*}
& \left| \left\langle\widehat{\mcZ}^{(\mathrm{\RN{1}})} \times_1\left(\mcP_{\whU_1^{(\mathrm{\RN{1}})}}-\mcP_{U_1}\right) \times_2 \mcP_{U_2} \times_3 \mcP_{U_3}, \mcA\right\rangle\right| 
\lesssim  \eqref{eq: upper bound of term 1 in step 1.1 in tensor regression with sample splitting} + \eqref{eq: upper bound of term 2 in step 1.1 in tensor regression with sample splitting} + \eqref{eq: upper bound of term 3 in step 1.1 in tensor regression with sample splitting} + \eqref{eq: upper bound of term 4 in step 1.1 in tensor regression with sample splitting} \\
\lesssim & \left\|\mcP_{U_1}A_1\mcP_{\left(U_3\otimes U_2\right)G_1^{\top}}\right\|_{\mathrm{F}} \cdot \frac{\sigma_\xi^2\oor^{1/2}}{\ulambda\sigma^2}\cdot \frac{\op}{n}
+  \left\|\mcA \times_2 U_2 \times_3 U_3\right\|_{\mathrm{F}} \cdot \left[\frac{\sigma_\xi^3\oor^{1/2}}{\ulambda\sigma^2}\left(\frac{\op\sqrt{\oor\log(\op)}}{n^{3/2}} \right)\right] \\
+ & \left\|\mcA\times_1 U_1 \times_2 U_2 \times_3 U_3\right\|_{\mathrm{F}} \cdot \left[\frac{\sigma_\xi^2\oor^{1/2}}{\ulambda\sigma^2}\cdot  \left(\frac{\sqrt{\oor\op\log(\op)}}{n} + \Delta^2 \cdot \frac{\log(\op)^{3/2}}{\sqrt{n}}\right) + \frac{\sigma_\xi^4\oor^{1/2}}{\ulambda^3\sigma^4}\cdot \frac{\op^2}{n^2}\right] .
\end{align*}

\subsubsection*{Step 1.2: Upper Bound of $\langle \widehat{\mcZ}^{(\mathrm{\RN{1}})} \times_1 \left(\mcP_{\widehat{U}_1^{(\mathrm{\RN{1}})}} - \mcP_{U_1}\right) \times_2 \left(\mcP_{\widehat{U}_2^{(\mathrm{\RN{1}})}} - \mcP_{U_2}\right) \times_3 \mcP_{U_3}, \mcA \rangle$}

Then we consider a similar decomposition as in Step 1.2 of the proof of Theorem~\ref{thm: main theorem in tensor regression without sample splitting}:
\begin{align}
\left|\left\langle \widehat{\mcZ}^{(\mathrm{\RN{1}})} \times_1 \left(\mcP_{\whU_1^{(\mathrm{\RN{1}})}} - \mcP_{U_1}\right) \times_2 \left(\mcP_{\whU_2^{(\mathrm{\RN{1}})}} - \mcP_{U_2}\right) \times_3 \mcP_{U_3}, \mcA \right\rangle\right| 
\leq & \mathrm{\RN{1}} + \mathrm{\RN{2}} + \mathrm{\RN{3}} +\mathrm{\RN{4}}. \label{eq: decomposition of step 1.2 in tensor regression with sample splitting}
\end{align}
Here, $\mathrm{\RN{1}}, \mathrm{\RN{2}}, \mathrm{\RN{3}}$, and $\mathrm{\RN{4}}$ are defined in the same manner as in \eqref{eq: term 1 in step 1.2 in tensor regression without sample splitting}, \eqref{eq: term 2 in step 1.2 in tensor regression without sample splitting}, \eqref{eq: term 3 in step 1.2 in tensor regression without sample splitting}, and \eqref{eq: term 4 in step 1.2 in tensor regression without sample splitting}, respectively, with $\widehat{Z}_1,\mcP_{\whU_j}$ replaced by $\widehat{Z}_1^{(\mathrm{\RN{1}})}$ and $\mcP_{\whU_j^{(\mathrm{\RN{1}})}}$.
Then, we can obtain

\begin{align*}
& \left|\left\langle \widehat{\mcZ}^{(\mathrm{\RN{1}})} \times_1\left(\mcP_{\whU_1^{(\mathrm{\RN{1}})}}-\mcP_{U_1}\right) \times_2\left(\mcP_{\whU_2^{(\mathrm{\RN{1}})}}-\mcP_{U_2}\right) \times_3 \mcP_{U_3}, \mcA\right\rangle\right| \\
\lesssim & \left\|\mcA \times_1 U_1 \times_2 U_2 \times_3 U_3 \right\|_{\mathrm{F}}\cdot \frac{\sigma_\xi^3\oor^{1/2}}{\ulambda^2\sigma^3}\cdot \frac{\op^{3/2}}{n^{3/2}} 
+  \left(\left\|\mcA \times_2 U_2 \times_3 U_3\right\|_{\mathrm{F}} + \left\|\mcA \times_1 U_1 \times_3 U_3\right\|_{\mathrm{F}}\right)  \left(\frac{\sigma_\xi^4\oor^{1/2}}{\ulambda^3\sigma^4}\cdot \frac{\op^{3/2}\sqrt{\oR\log(\op)}}{n^2}\right) \\
+ & \left\|\mcA \times_3 U_3 \right\|_{\mathrm{F}} \cdot \left[\frac{\sigma_{\xi}^3\oor^{1/2}}{\ulambda^2 \sigma^3} \cdot\left(\frac{\oor^{3/2}\log(\op)^{3/2}}{n^{3/2}}+\Delta \cdot \frac{\oor^{1/2}\oR\log(\op)^{3/2}}{n^{3/2}}\right) + \frac{\sigma_{\xi}^4\oor^{1/2}}{\ulambda^3 \sigma^4}\left(\frac{\oor\op\log (\op)}{n^2}+\Delta \cdot \frac{\oR\op\log(\op)}{n^2}\right)\right] .
\end{align*}

\subsubsection*{Step 1.3: Upper Bound of $\langle \widehat{\mcZ}^{(\mathrm{\RN{1}})} \times_1 \left(\mcP_{\widehat{U}_1^{(\mathrm{\RN{1}})}} - \mcP_{U_1}\right) \times_{2} \left(\mcP_{\widehat{U}_{2}^{(\mathrm{\RN{1}})}} - \mcP_{U_2}\right) \times_{3} \left(\mcP_{\widehat{U}_{3}^{(\mathrm{\RN{1}})}} - \mcP_{U_3}\right), \mcA \rangle$}

By symmetry, it suffices to consider
\begin{align}
\mathrm{\RN{1}}= & \left|\left\langle \widehat{\mcZ}^{(\mathrm{\RN{1}})} \times_1 \left(\mcP_{\whU_1^{(\mathrm{\RN{1}})}} - \mcP_{U_1}\right) \times_{2} \left(\mcP_{\whU_{2}^{(\mathrm{\RN{1}})}} - \mcP_{U_2}\right) \times_{3} \left(\mcP_{\whU_{3}^{(\mathrm{\RN{1}})}} - \mcP_{U_3}\right), \mcA \times_1 \mcP_{U_1} \times_2 \mcP_{U_2} \times_3 \mcP_{U_3} \right\rangle\right|, 
\label{eq: term 1 in step 1.3 in tensor regression with sample splitting}\\
\mathrm{\RN{2}}= & \left|\left\langle \widehat{\mcZ}^{(\mathrm{\RN{1}})} \times_1 \left(\mcP_{\whU_1^{(\mathrm{\RN{1}})}} - \mcP_{U_1}\right) \times_{2} \left(\mcP_{\whU_{2}^{(\mathrm{\RN{1}})}} - \mcP_{U_2}\right) \times_{3} \left(\mcP_{\whU_{3}^{(\mathrm{\RN{1}})}} - \mcP_{U_3}\right), \mcA \times_1 \mcP_{U_{1\perp}} \times_2 \mcP_{U_2} \times_3 \mcP_{U_3} \right\rangle\right|,
\label{eq: term 2 in step 1.3 in tensor regression with sample splitting} \\
\mathrm{\RN{3}}= & \left|\left\langle \widehat{\mcZ}^{(\mathrm{\RN{1}})} \times_1 \left(\mcP_{\whU_1^{(\mathrm{\RN{1}})}} - \mcP_{U_1}\right) \times_{2} \left(\mcP_{\whU_{2}^{(\mathrm{\RN{1}})}} - \mcP_{U_2}\right) \times_{3} \left(\mcP_{\whU_{3}^{(\mathrm{\RN{1}})}} - \mcP_{U_3}\right), \mcA \times_1 \mcP_{U_1} \times_2 \mcP_{U_{2\perp}} \times_3 \mcP_{U_3} \right\rangle\right|,
\label{eq: term 3 in step 1.3 in tensor regression with sample splitting}\\
\mathrm{\RN{4}}= & \left|\left\langle \widehat{\mcZ}^{(\mathrm{\RN{1}})} \times_1 \left(\mcP_{\whU_1^{(\mathrm{\RN{1}})}} - \mcP_{U_1}\right) \times_{2} \left(\mcP_{\whU_{2}^{(\mathrm{\RN{1}})}} - \mcP_{U_2}\right) \times_{3} \left(\mcP_{\whU_{3}^{(\mathrm{\RN{1}})}} - \mcP_{U_3}\right), \mcA \times_1 \mcP_{U_{1\perp}} \times_2 \mcP_{U_{2\perp}} \times_3 \mcP_{U_{3\perp}} \right\rangle\right| \label{eq: term 4 in step 1.3 in tensor regression with sample splitting}.
\end{align}

Applying the same decomposition for Step 1.3 in the proof of Theorem~\ref{thm: main theorem in tensor regression without sample splitting}, with $\widehat{Z}_3$ replaced by $\widehat{Z}_3^{(\mathrm{\RN{1}})}$, we can show
\begin{align}
& \left|\left\langle \widehat{\mcZ}^{(\mathrm{\RN{1}})} \times_1 \left(\mcP_{\whU_1^{(\mathrm{\RN{1}})}} - \mcP_{U_1}\right) \times_{2} \left(\mcP_{\whU_{2}^{(\mathrm{\RN{1}})}} - \mcP_{U_{2}}\right) \times_{3} \left(\mcP_{\whU_{3}^{(\mathrm{\RN{1}})}} - \mcP_{U_3}\right), \mcA \right\rangle\right| \notag \\
\lesssim & \left\|\mcA \times_1 U_1 \times_2 U_2 \times_3 U_3\right\|_{\mathrm{F}} \cdot \frac{\sigma_{\xi}^4 \oor^{1 / 2}}{\ulambda^3 \sigma^4} \cdot \frac{\op^2}{n^2} + \sum_{j=1}^3 \left\|\mcA \times_{j+1} U_{j+1} \times_{j+2} U_{j+2}\right\|_{\mathrm{F}} \cdot\left[\frac{\sigma_{\xi}^4 \oor^{1 / 2}}{\ulambda^3 \sigma^4} \cdot\left(\frac{\oor \op^{3 / 2} \sqrt{\log (\op)}}{n^2}\right)\right] \notag \\
+ & \sum_{j=1}^3 \left\|\mcA \times_j U_j\right\|_{\mathrm{F}} \cdot \left[\frac{\sigma_{\xi}^3\oor^{1/2}}{\ulambda^2 \sigma^3} \cdot\left(\frac{\oor^{3/2}\log(\op)^{3/2}}{n^{3/2}}+\Delta \cdot \frac{\oor^{1/2}\oR\log(\op)^{3/2}}{n^{3/2}}\right) + \frac{\sigma_{\xi}^4\oor^{1/2}}{\ulambda^3 \sigma^4}\left(\frac{\oor\op\log (\op)}{n^2}+\Delta \cdot \frac{\oR\op\log(\op)}{n^2}\right)\right] \notag \\
+ & \left\|\mcA\right\|_{\mathrm{F}} \cdot \left[\frac{\sigma_{\xi}^4\oor^{1/2}}{\ulambda^3 \sigma^4} \cdot\left(\frac{\oor^2\log(\op)^2}{n^2}+\Delta \cdot \frac{\oor^{1/2} \oR^{3/2} \log (\op)^2}{n^2}\right) + \frac{\sigma_\xi^5\oor^{1/2}}{\ulambda^4\sigma^5}\cdot \left(\frac{\oor^{3/2}\op\log(\op)^{3/2}}{n^{5/2}} + \Delta\cdot \frac{\oR^{3/2} \op \log (\op)^{3/2}}{n^{5/2}}\right)\right] \notag .
\end{align}

\subsection*{Step 2: Upper Bound of Negligible Terms in $\langle {\mcT} \times_1 \mcP_{\widehat{U}_{1}^{(\mathrm{\RN{1}})}} \times_2 \mcP_{\widehat{U}_{2}^{(\mathrm{\RN{1}})}} \times_3 \mcP_{\widehat{U}_{3}^{(\mathrm{\RN{1}})}}-\mcT, \mcA \rangle$}

The proof of upper bound of the negligible terms remains the same with the no-sample-splitting case. By the sample splitting procedure, we have removed the dependence between the initial estimate and the artificial noise in the de-biasing procedure. Similar to the arguments in Step 2 in the proof of Theorem~\ref{thm: main theorem in tensor regression without sample splitting}, by symmetry, it remains to consider
\begin{align}
\text{(Step 2.1)}: & \left\langle \mcT \times_1 \left(\mcP_{\whU_1^{(\mathrm{\RN{1}})}} - \mcP_{U_1}\right) \times_2 \mcP_{U_2} \times_3 \mcP_{U_3}, \mcA \right\rangle \label{eq: step 2.1 in tensor regression with sample splitting} \\
\text{(Step 2.2)}: & \left\langle \mcT \times_1 \left(\mcP_{\whU_1^{(\mathrm{\RN{1}})}} - \mcP_{U_1}\right) \times_2 \left(\mcP_{\whU_2^{(\mathrm{\RN{1}})}} - \mcP_{U_2}\right) \times_3 \mcP_{U_3}, \mcA \right\rangle \label{eq: step 2.2 in tensor regression with sample splitting} \\
\text{(Step 2.3)}: & \left\langle \mcT \times_1 \left(\mcP_{\whU_1^{(\mathrm{\RN{1}})}} - \mcP_{U_1}\right) \times_2 \left(\mcP_{\whU_2^{(\mathrm{\RN{1}})}} - \mcP_{U_2}\right) \times_3 \left(\mcP_{\whU_3^{(\mathrm{\RN{1}})}} - \mcP_{U_3}\right), \mcA \right\rangle  . \label{eq: step 2.3 in tensor regression with sample splitting}
\end{align}

\subsubsection*{Step 2.1: Upper Bound of Negligible Terms in $\langle \mcT \times_1 \left(\mcP_{\whU_1^{(\mathrm{\RN{1}})}} - \mcP_{U_1}\right) \times_2 \mcP_{U_2} \times_3 \mcP_{U_3}\rangle$ }

Note that
\begin{align*}
& \left\langle \mcT \times_1 \left(\mcP_{\whU_1^{(\mathrm{\RN{1}})}} - \mcP_{U_1}\right) \times_2 \mcP_{U_2} \times_3 \mcP_{U_3}, \mcA \right\rangle \notag \\
= & \underbrace{\left\langle \mcT \times_1 S_{G_1,1}\left(\whE_1^{(\mathrm{\RN{1}})}\right) \times_2 \mcP_{U_2} \times_3 \mcP_{U_3}, \mcA \right\rangle}_{\mathrm{\RN{1}}}  + \underbrace{\left\langle \mcT \times_1 \sum_{k_1=2}^{+\infty}S_{G_1,k_1}\left(\whE_1^{(\mathrm{\RN{1}})}\right) \times_2 \mcP_{U_2} \times_3 \mcP_{U_3}, \mcA \right\rangle}_{\mathrm{\RN{2}}}. 
\end{align*}

By \eqref{eq: definition of whEj in tensor regression with sample splitting} and \eqref{eq: definition of S1(whEj) in tensor regression with sample splitting}, it then follows that
\begin{align*}
\mathrm{\RN{1}}
= & \underbrace{\left\langle \mcP_{U_{1\perp}} \whZ_1^{(\mathrm{\RN{1}})}\left(\mcP_{U_3} \otimes \mcP_{U_2}\right) \mcP_{\left(U_3 \otimes U_2\right)G_1^{\top}}, A_1\right\rangle}_{\mathrm{\RN{1}}.\mathrm{\RN{1}}, \text{asymptotically normal}} 
+  \underbrace{\left\langle \mcP_{U_{1\perp}} \whZ_1^{(\mathrm{\RN{1}})}\left[\mcP_{U_3} \otimes \left(\mcP_{\whU_2^{(0),(\mathrm{\RN{2}})}}-\mcP_{U_2}\right)\right] \mcP_{\left(U_3 \otimes U_2\right)G_1^{\top}} , A_1\right\rangle}_{\mathrm{\RN{1}}.\mathrm{\RN{2}}, \text{negligible}}\\
+ & \underbrace{\left\langle \mcP_{U_{1\perp}} \whZ_1^{(\mathrm{\RN{1}})}\left[\left(\mcP_{\whU_3^{(0),(\mathrm{\RN{2}})}}-\mcP_{U_3}\right) \otimes \mcP_{U_2}\right] \mcP_{\left(U_3 \otimes U_2\right)G_1^{\top}}, A_1\right\rangle}_{\mathrm{\RN{1}}.\mathrm{\RN{3}}, \text{negligible}} \\
+ & \underbrace{\left\langle \mcP_{U_{1\perp}} \whZ_1^{(\mathrm{\RN{1}})}\left[\left(\mcP_{\whU_3^{(0),(\mathrm{\RN{2}})}}-\mcP_{U_3}\right) \otimes \left(\mcP_{\whU_2^{(0),(\mathrm{\RN{2}})}}-\mcP_{U_2}\right)\right] \mcP_{\left(U_3 \otimes U_2\right)G_1^{\top}}, A_1\right\rangle}_{\mathrm{\RN{1}}.\mathrm{\RN{4}}, \text{negligible}}\\
+ & \underbrace{\left\langle \mcP_{U_{1\perp}} \whZ_1^{(\mathrm{\RN{1}})}\left(\mcP_{\whU_3^{(0),(\mathrm{\RN{2}})}} \otimes \mcP_{\whU_2^{(0),(\mathrm{\RN{2}})}}\right) \whZ_1^{(\mathrm{\RN{2}}),\top} U_1 \left(G_1G_1^{\top}\right)^{-1}  G_1 \left(U_3 \otimes U_2\right)^{\top} , A_1\right\rangle}_{\mathrm{\RN{1}}.\mathrm{\RN{5}}, \text{negligible}}.
\end{align*}

We leave the proof of the asymptotic normality $\mathrm{\RN{1}}.\mathrm{\RN{1}}=\left\langle \mcP_{U_{1\perp}} \whZ_1^{(\mathrm{\RN{1}})} \mcP_{\left(U_3 \otimes U_2\right)G_1^{\top}} , A_1\right\rangle$ to Step 3.  First, consider the upper bound for
\begin{align}
\mathrm{\RN{1}}.\mathrm{\RN{2}}
\leq & \left\|\mcP_{\left(U_3\otimes U_2\right)G_1^{\top}}\right\|_{\mathrm{F}} \cdot \left\|\mcP_{\left(U_3\otimes U_2\right)G_1^{\top}}A_1\mcP_{U_{1\perp}}\whZ_1^{(\mathrm{\RN{1}})} \left[\left(\mcP_{\whU_2^{(0),(\mathrm{\RN{2}})}} - \mcP_{U_2}\right) \mcP_{U_2}\otimes \mcP_{U_3}\right]\right\|_{\mathrm{F}} \notag \\
\lesssim & \left\|\mcP_{U_{1 \perp}} A_1 \mcP_{\left(U_3 \otimes U_2\right) G_1^{\top}}\right\|_{\mathrm{F}} \cdot\left(\frac{\sigma_{\xi}^2 \oor^{1 / 2}}{\ulambda \sigma^2} \cdot \frac{\sqrt{\op\oor\log(\op)}}{n}\right) ,\label{eq: upper bound of term 1.2 in step 2.1 in tensor regression with sample splitting}
\end{align} 
where the last inequality follows from the independency between $\widehat{\mcZ}^{(\mathrm{\RN{1}})}$ and $\mcP_{\whU_2^{(0),(\mathrm{\RN{2}})}}$.

By the same argument, we have 
\begin{align}
\mathrm{\RN{1}}.\mathrm{\RN{3}} 
& \lesssim \left\|\mcP_{U_{1 \perp}} A_1 \mcP_{\left(U_3 \otimes U_2\right) G_1^{\top}}\right\|_{\mathrm{F}} \cdot\left(\frac{\sigma_{\xi}^2 \oor^{1 / 2}}{\ulambda \sigma^2} \cdot \frac{\sqrt{\op\oor\log(\op)}}{n}\right), \label{eq: upper bound of term 1.3 in step 2.1 in tensor regression with sample splitting} \\
\mathrm{\RN{1}}.\mathrm{\RN{4}} 
\lesssim & \left\|\mcP_{U_{1 \perp}} A_1 \mcP_{\left(U_3 \otimes U_2\right) G_1^{\top}}\right\|_{\mathrm{F}} \cdot\left(\frac{\sigma_{\xi}^3 \oor^{1 / 2}}{\ulambda^2 \sigma^3} \cdot \frac{\op\oor\log(\op)}{n^{3 / 2}}\right), \label{eq: upper bound of term 1.4 in step 2.1 in tensor regression with sample splitting} \\
\mathrm{\RN{1}}.\mathrm{\RN{5}} 
\lesssim & \left\|\mcP_{U_{1\perp}}A_1\mcP_{\left(U_3\otimes U_2\right)G_1^{\top}}\right\|_{\mathrm{F}} \cdot \frac{\sigma_\xi^2\oor^{1/2}}{\ulambda\sigma^2}\cdot \frac{\oor\log(\op)}{n} .\label{eq: upper bound of term 1.5 in step 2.1 in tensor regression with sample splitting}
\end{align}

Therefore, combining the results above, we have
\begin{align}
\mathrm{\RN{1}} 
\lesssim & \eqref{eq: upper bound of term 1.2 in step 2.1 in tensor regression with sample splitting} + \eqref{eq: upper bound of term 1.3 in step 2.1 in tensor regression with sample splitting} + \eqref{eq: upper bound of term 1.4 in step 2.1 in tensor regression with sample splitting} + \eqref{eq: upper bound of term 1.5 in step 2.1 in tensor regression with sample splitting} \lesssim \left\|\mcP_{U_{1\perp}}A_1\mcP_{\left(U_3\otimes U_2\right)G_1^{\top}}\right\|_{\mathrm{F}} \cdot \left(\frac{\sigma_\xi^2\oor^{1/2}}{\ulambda\sigma^2}\cdot \frac{\sqrt{\oor^2\op\log(\op)}}{n} \right) \label{eq: upper bound of term 1 in step 2.1 in tensor regression with sample splitting}.
\end{align}

Then, consider an upper bound for the remaining higher-order terms. By the same arguments in the proof of \eqref{eq: upper bound of term 2 in step 2.1 in tensor regression without sample splitting}, it follows that
\begin{align}
\mathrm{\RN{2}}
\lesssim & \left\|U_1^{\top}A_1\mcP_{\left(U_3\otimes U_2\right)G_1^\top}\right\|_{\mathrm{F}}\cdot \frac{\sigma_\xi^2\oor^{1/2}}{\ulambda\sigma^2}\cdot \frac{\op}{n} + \left\|U_{1\perp}^{\top}A_1^{\top}\mcP_{\left(U_3\otimes U_2\right)G_1^{\top}}\right\|_{\mathrm{F}} \left(\frac{\sigma_\xi^3\oor^{1/2}}{\ulambda\sigma^3} \cdot \frac{\oor\log(\op)}{n} + \frac{\sigma_{\xi}^2\oor^{1/2}}{\ulambda \sigma^2} \cdot \frac{\op\sqrt{\oor\log(\op)}}{n^{3 / 2}}\right) \label{eq: upper bound of term 2 in step 2.1 in tensor regression with sample splitting}.
\end{align}

Here, we used
\begin{align*}
& \left|\operatorname{tr}\left[\left(\mcP_{U_3} \otimes \mcP_{U_2}\right) A_1^{\top} \sum_{k_1=3}^{+\infty} S_{G_1, k_1}\left(\whE_1^{(\mathrm{\RN{1}})}\right) U_1 G_1\left(U_3 \otimes U_2\right)^{\top}\right]\right| \\
\leq & \left|\operatorname{tr}\left[\left(\mcP_{U_3} \otimes \mcP_{U_2}\right) A_1^{\top}\mcP_{U_1}\sum_{k_1=3}^{+\infty} S_{G_1, k_1}\left(\whE_1^{(\mathrm{\RN{1}})}\right) U_1 G_1\left(U_3 \otimes U_2\right)^{\top}\right]\right| \\
+ & \left|\operatorname{tr}\left[\left(\mcP_{U_3} \otimes \mcP_{U_2}\right) A_1^{\top}\mcP_{U_{1\perp}}\sum_{k_1=3}^{+\infty} S_{G_1, k_1}\left(\whE_1^{(\mathrm{\RN{1}})}\right) U_1 G_1\left(U_3 \otimes U_2\right)^{\top}\right]\right| \\
\lesssim & \left\|U_1^{\top}A_1\mcP_{\left(U_3\otimes U_2\right)G_1^\top}\right\|_{\mathrm{F}}\cdot \left(\frac{\sigma_\xi}{\ulambda\sigma}\cdot \sqrt{\frac{\op}{n}}\right)^3 \cdot \ulambda\oor^{1/2} + \underbrace{\left\|U_{1\perp}^{\top}A_1\mcP_{\left(U_3\otimes U_2\right)G_1^\top}\right\|_{\mathrm{F}}\cdot \left(\frac{\sigma_{\xi}^3}{\ulambda^3 \sigma^3} \cdot \frac{\op\sqrt{\oor\log(\op)}}{n^{3 / 2}}\right)}_{\eqref{eq: high-prob upper bound of PUorder3PUpPV in tensor regression with sample splitting}} \cdot \ulambda\oor^{1/2} \\
= & \left\|U_1^{\top}A_1\mcP_{\left(U_3\otimes U_2\right)G_1^\top}\right\|_{\mathrm{F}}\cdot \frac{\sigma_{\xi}^3\oor^{1/2}}{\ulambda^2 \sigma^3} \cdot \frac{\op^{3 / 2}}{n^{3 / 2}} + \left\|U_{1\perp}^{\top}A_1\mcP_{\left(U_3\otimes U_2\right)G_1^\top}\right\|_{\mathrm{F}}\cdot \left(\frac{\sigma_{\xi}^3\oor^{1/2}}{\ulambda^2 \sigma^3} \cdot \frac{\op\sqrt{\oor\log(\op)}}{n^{3 / 2}}\right).
\end{align*} 

Therefore, we have 
\begin{align*}
& \left| \left\langle \mcT \times_1 \left( \mcP_{\whU_1^{(\mathrm{\RN{1}})}} - \mcP_{U_1} \right) \times_2 \mcP_{U_2} \times_3 \mcP_{U_3}, \mcA \right\rangle - \left\langle \mcP_{U_{1\perp}} \whZ_1^{(\mathrm{\RN{1}})} \left( \mcP_{U_3} \otimes \mcP_{U_2} \right) \mcP_{\left( U_3 \otimes U_2 \right) G_1^{\top}}, A_1 \right\rangle \right| \\
\lesssim & \underbrace{\left\|\mcP_{U_{1 \perp}} A_1 \mcP_{\left(U_3 \otimes U_2\right) G_1^{\top}}\right\|_{\mathrm{F}} \cdot \frac{\sigma_{\xi}^2 \oor^{1 / 2}}{\ulambda \sigma^2} \cdot\left(\frac{\sqrt{\oor \op\log (\op)} }{n}+\frac{\sigma_{\xi}^3}{\ulambda^2 \sigma^3} \cdot \frac{\op^{3 / 2}}{n^{3 / 2}}\right)}_{\eqref{eq: upper bound of term 1 in step 2.1 in tensor regression with sample splitting}} \\
+ & \underbrace{\left\|U_1^{\top}A_1\mcP_{\left(U_3\otimes U_2\right)G_1^\top}\right\|_{\mathrm{F}}\cdot \frac{\sigma_\xi^2\oor^{1/2}}{\ulambda\sigma^2}\cdot \frac{\op}{n} +  \left\|U_{1\perp}^{\top}A_1^{\top}\mcP_{\left(U_3\otimes U_2\right)G_1^{\top}}\right\|_{\mathrm{F}}\cdot \left(\frac{\sigma_\xi^3\oor^{1/2}}{\ulambda\sigma^3} \cdot \frac{\oor\log(\op)}{n} + \frac{\sigma_{\xi}^3\oor^{1/2}}{\ulambda^2 \sigma^3} \cdot \frac{\op\sqrt{\oor\log(\op)}}{n^{3 / 2}}\right)}_{\eqref{eq: upper bound of term 2 in step 2.1 in tensor regression with sample splitting}} \\
\lesssim & \left\|\mcP_{U_1} A_1 \mcP_{\left(U_3 \otimes U_2\right) G_1^{\top}}\right\|_{\mathrm{F}} \cdot \frac{\sigma_\xi^2\oor^{1/2}}{\ulambda\sigma^2}\cdot \frac{\op}{n} + \left\|\mcP_{U_{1 \perp}} A_1 \mcP_{\left(U_3 \otimes U_2\right) G_1^{\top}}\right\|_{\mathrm{F}} \cdot \left(\frac{\sigma_\xi^2\oor^{1/2}}{\ulambda\sigma^2} \cdot \frac{\oor\log(\op)}{n} + \frac{\sigma_{\xi}^3\oor^{1/2}}{\ulambda^2 \sigma^3} \cdot \frac{\op\sqrt{\oor\log(\op)}}{n^{3 / 2}}\right).
\end{align*}

\subsubsection*{Step 2.2: Upper Bound of $\langle \mcT\times_1 \left(\mcP_{\widehat{U}_1^{(\mathrm{\RN{1}})}} - \mcP_{U_1}\right) \times_2 \left(\mcP_{\widehat{U}_2^{(\mathrm{\RN{1}})}} - \mcP_{U_2}\right) \times_3 \mcP_{U_3}, \mcA \rangle $}

Consider the same decomposition as in Step 2.2 in the proof of Theorem~\ref{thm: main theorem in tensor regression without sample splitting}, with $\mcP_{\whU_j}$ replaced by $\mcP_{\whU_j^{(\mathrm{\RN{1}})}}$ for any $j=1,2,3$:
\begin{align*}
\left|\left\langle \mcT\times_1 \left(\mcP_{\whU_1^{(\mathrm{\RN{1}})}} - \mcP_{U_1}\right) \times_2 \left(\mcP_{\whU_2^{(\mathrm{\RN{1}})}} - \mcP_{U_2}\right) \times_3 \mcP_{U_3}, \mcA\right\rangle\right| 
\leq & \mathrm{\RN{1}} + \mathrm{\RN{2}} + \mathrm{\RN{3}} + \mathrm{\RN{4}}.
\end{align*}
Applying similar arguments, we can show
\begin{align}
\mathrm{\RN{1}}
\lesssim & \left\|\mcA \times_3 U_3\right\|_{\mathrm{F}} \cdot \left(\frac{\sigma_{\xi}^2\oor^{1/2}}{\ulambda \sigma^2} \cdot \frac{\oor\log(\op)}{n} +\frac{\sigma_{\xi}^2\oor^{1/2}}{\ulambda \sigma^2} \cdot \Delta \cdot \frac{\oR \log (\op)}{n}\right) \label{eq: upper bound of term 1 in step 2.2 in tensor regression with sample splitting}, \\
\mathrm{\RN{2}} 
\lesssim & \left\|\mcA \times_2 U_2 \times_3 U_3\right\|_{\mathrm{F}} \cdot \left(\frac{\sigma_\xi^3\oor^{1/2}}{\ulambda^2\sigma^3} \cdot \frac{\op\sqrt{\oor\log(\op)}}{n^{3/2}}\right) 
\label{eq: upper bound of term 2 in step 2.2 with sample splitting}, \\
\mathrm{\RN{3}}
\lesssim & \left\|\mcA \times_1 U_1 \times_3 U_3\right\|_{\mathrm{F}} \cdot \left(\frac{\sigma_\xi^3\oor^{1/2}}{\ulambda^2\sigma^3} \cdot \frac{\op\sqrt{\oor\log(\op)}}{n^{3/2}}\right) 
\label{eq: upper bound of term 3 in step 2.2 with sample splitting}, \\
\mathrm{\RN{4}}
\lesssim & \left\|\mcA \times_1 U_1 \times_2 U_2 \times_3 U_3\right\|_{\mathrm{F}} \cdot \left(\frac{\sigma_\xi^3\oor^{1/2}}{\ulambda^2\sigma^3} \cdot \frac{\op\sqrt{\oor\log(\op)}}{n^{3/2}} + \frac{\sigma_\xi^4\oor^{1/2}}{\ulambda^3\sigma^4} \cdot \frac{\op^2}{n^2}\right)  .
\label{eq: upper bound of term 4 in step 2.2 with sample splitting}
\end{align}

\subsubsection*{Step 2.3: Upper Bound of $\langle \mcT\times_1 \left(\mcP_{\widehat{U}_1^{(\mathrm{\RN{1}})}} - \mcP_{U_1}\right) \times_2 \left(\mcP_{\widehat{U}_2^{(\mathrm{\RN{1}})}} - \mcP_{U_2}\right) \times_3 \left(\mcP_{\widehat{U}_3^{(\mathrm{\RN{1}})}} - \mcP_{U_3}\right), \mcA \rangle$}
By similar arguments of decomposition in Step 2.3 in the proof of Theorem~\ref{thm: main theorem in tensor regression without sample splitting}, with $\mcP_{\whU_j}$ replaced by $\mcP_{\whU_j^{(\mathrm{\RN{1}})}}$ for any $j=1,2,3$, we have the same decomposition.

Then, we can obtain
\begin{align*}
& \left|\left\langle\mcT \times_1\left(\mcP_{\whU_1^{(\mathrm{\RN{1}})}}-\mcP_{U_1}\right) \times_2\left(\mcP_{\whU_2^{(\mathrm{\RN{1}})}}-\mcP_{U_2}\right) \times_3\left(\mcP_{\whU_3^{(\mathrm{\RN{1}})}}-\mcP_{U_3}\right), \mcA\right\rangle\right| \\
\lesssim & \left\|\mcA \times_1 U_1 \times_2 U_2 \times_3 U_3\right\|_{\mathrm{F}} \cdot \frac{\sigma_\xi^6\oor^{1/2}}{\ulambda^5\sigma^6} \cdot \frac{\op^3}{n^3} + \sum_{j=1}^3 \left\|\mcA \times_{j+1} U_{j+1} \times_{j+2} U_{j+2}\right\|_{\mathrm{F}} \cdot \left(\frac{\sigma_\xi^3\oor^{1/2}}{\ulambda^2\sigma^3} \cdot \frac{\op\sqrt{\oor\log(\op)}}{n^{3/2}}\right) \\
+ & \sum_{j=1}^3 \left\|\mcA \times_j U_j\right\|_{\mathrm{F}} \cdot \left(\frac{\sigma_\xi^4\oor^{1/2}}{\ulambda^3\sigma^4} \cdot \frac{\oor^2\op\log(\op)}{n^2} + \frac{\sigma_\xi^4 \oor^{1/2}}{\ulambda^3\sigma^4}\cdot \Delta \cdot \frac{\oR \op \log (\op)}{n^2}\right) \\
+ & \left\|\mcA\right\|_{\mathrm{F}} \cdot \left(\frac{\sigma_{\xi}^3\oor^{1/2}}{\ulambda^2 \sigma^3} \cdot \frac{\oor^{3/2}\log(\op)^{3/2}}{n^{3 / 2}}+\frac{\sigma_{\xi}^3\oor^{1/2}}{\ulambda^2 \sigma^3} \cdot \Delta \cdot \frac{\oR^{3 / 2} \log (\op)^{3 / 2}}{n^{3 / 2}}\right).
\end{align*}

\subsection*{Step 3: Analysis of asymptotic normal terms}

Let 
$$
\mcP_{\mathbb{T}_{\mcT} \mathcal{M}_{\mathbf{r}}}\left(\mcA\right) := \mcP_{\mathbb{T}_{\mcT} \mathcal{M}_{(r_1, r_2, r_3)}}\left(\mcA\right) = \sum_{j=1}^3 \operatorname{Mat}_j^{-1}\left(\mcP_{U_{j\perp}}A_j\mcP_{\left(U_{j+2}\otimes U_{j+1}\right)G_j^{\top}}\right) + \mcA \times_1 \mcP_{U_1} \times_2 \mcP_{U_2} \times_3 \mcP_{U_3}.
$$

\subsubsection*{Step 3.1: Asymptotic Normality of 
$\langle \widehat{\mcZ}^{(1), (\mathrm{\RN{1}}) }, \mcP_{\mathbb{T}_{\mcT} \mathcal{M}_{(r_1, r_2, r_3)}}\left(\mcA\right) \rangle + \langle \widehat{\mcZ}^{(1), (\mathrm{\RN{2}}) }, \mcP_{\mathbb{T}_{\mcT} \mathcal{M}_{(r_1, r_2, r_3)}}\left(\mcA\right) \rangle$}

First, we have
\begin{align*}
\frac{n_1}{n} \big\langle \whZ_j^{(1),(\mathrm{\RN{1}})}, \mcP_{\mathbb{T}_{\mcT} \mathcal{M}_{\mathbf{r}}}(\mcA)\big\rangle + \frac{n_2}{n} \big\langle \whZ_j^{(1),(\mathrm{\RN{2}})}, \mcP_{\mathbb{T}_{\mcT} \mathcal{M}_{\mathbf{r}}}(\mcA)\big\rangle 
= & \frac{1}{n\sigma^2}\sum_{i_1=1}^{n_1}\xi_{i_1}^{(\mathrm{\RN{1}})} \langle \mcX_{i_1}^{(\mathrm{\RN{1}})}, \mcP_{\mathbb{T}_{\mcT} \mathcal{M}_{\mathbf{r}}}(\mcA)\rangle + \frac{1}{n\sigma^2}\sum_{i_2=1}^{n_2}\xi_{i_2}^{(\mathrm{\RN{2}})} \langle  \mcX_{i_2}^{(\mathrm{\RN{2}})}, \mcP_{\mathbb{T}_{\mcT} \mathcal{M}_{\mathbf{r}}}(\mcA)\rangle.
\end{align*}

Then, by the same arguments as in the proof of Theorem~\ref{thm: main theorem in tensor regression without sample splitting}, we have
\begin{align*}
&\sup_{x \in \mathbb{R}} \left| \mathbb{P}  \left(\frac{\frac{n_1}{n}\cdot \left\langle \widehat{\mcZ}^{(1),(\mathrm{\RN{1}})},  \mcA \right\rangle + \frac{n_2}{n} \cdot \left\langle \widehat{\mcZ}^{(1),(\mathrm{\RN{2}})}, \mcA\right\rangle}{\frac{\sigma_{\xi}}{\sigma} \cdot \left(\sum_{j=1}^{3}\left\|\mcP_{U_{j\perp}}A_j\mcP_{\left(U_{j+2} \otimes U_{j+1}\right)G_j^{\top}} \right\|_{\mathrm{F}}^2 + \left\|\mcA \times_1 U_1 \times_2 U_2 \times_3 U_3\right\|_{\mathrm{F}}^2\right)^{1/2} \cdot \sqrt{\frac{1}{n}}} \leq x\right)-\Phi(x) \right| \lesssim \sqrt{\frac{1}{n}}.
\end{align*}

\subsubsection*{Step 3.2: Upper Bound of 
$\langle \widehat{\mcZ}^{(2),(\mathrm{\RN{2}})}, \mcP_{\mathbb{T}_{\mcT} \mathcal{M}_{(r_1, r_2, r_3)}}\left(\mcA\right)\rangle$}

Furthermore, we have
\begin{align*}
& \left|\frac{n_1}{n}\sum_{j=1}^{3}\left\langle  \widehat{\mcZ}^{(2),(\mathrm{\RN{1}})}, \mcP_{\mathbb{T}_{\mcT} \mathcal{M}_{(r_1, r_2, r_3)}}\left(\mcA\right)\right\rangle + \frac{n_2}{n}\left\langle \widehat{\mcZ}^{(2),(\mathrm{\RN{2}})}, \mcP_{\mathbb{T}_{\mcT} \mathcal{M}_{(r_1, r_2, r_3)}}\left(\mcA\right)\right\rangle\right|	\\
\lesssim & \left(\sum_{j=1}^3\left\|\mcP_{U_{j\perp}}A_j\mcP_{\left(U_{j+2}\otimes U_{j+1}\right)G_j^{\top}}\right\|_{\mathrm{F}} + \left\|\mcA\times_1 U_1\times_2 U_2\times_3 U_3\right\|_{\mathrm{F}}\right) \cdot \Delta \cdot \frac{\sigma_\xi\oor^{1/2}}{\sigma} \sqrt{\frac{\oor\log(\op)}{n}}.
\end{align*}

\subsubsection*{Step 3.3: Combining Asymptotic Normal Terms and Negligible Terms}

By the Lipschitz property of $\Phi(x)$ and note that the discussion above holds under event $\left\|\widehat{\mcT}^{\text{init}} - \mcT\right\|_{\mathrm{F}} \leq \Delta$ and $\left\|\mcP_{\whU_j} - \mcP\right\| \leq \frac{\sigma_\xi}{\sigma}\sqrt{\frac{\op}{n}}$, then finally we have
{\small \begin{align*}
& \sup _{x \in \mathbb{R}} \left\lvert\, \mathbb{P}\left(\frac{\left\langle\widehat{T}, \mcA\right\rangle-\langle T, \mcA\rangle}{\frac{\sigma_{\xi}}{\sigma} \cdot \left(\sum_{j=1}^3 \left\|\mcP_{U_{j \perp}} A_j \mcP_{\left(U_{j+2} \otimes U_{j+1}\right)G_j^{\top}} \right\|_{\mathrm{F}}^2+ \left\|\mcA \times_1 \mcP_{U_1} \times_2 \mcP_{U_2} \times_3 \mcP_{U_3} \right\|_{\mathrm{F}}^2\right)^{1/2}\cdot \sqrt{\frac{1}{n}}} \leq x\right)-\Phi(x)\right| \notag\\
\leq & \sqrt{\frac{1}{n}} + \underbrace{\left[\frac{1}{\op^c} + \exp\left(-c\op\right) + \exp\left(-cn\right) + \mcP\left(\mcE_\Delta\right) + \mcP\left(\mcE_U^{\text{reg}}\right)\right]}_{\text{rate of initial estimate}}\\
&+ \frac{1}{\frac{\sigma_{\xi}}{\sigma} \cdot s_{\mcA}\cdot \sqrt{\frac{1}{n}}} \cdot \Bigg\{\underbrace{\left\|\mcA\times_1 U_1 \times_2 U_2 \times_3 U_3\right\|_{\mathrm{F}} \cdot \left[\frac{\sigma_\xi^2\oor^{1/2}}{\ulambda\sigma^2}\cdot  \left(\frac{\sqrt{\oor\op\log(\op)}}{n} + \Delta^2 \cdot \frac{\log(\op)^{3/2}}{\sqrt{n}}\right) + \frac{\sigma_\xi^4\oor^{1/2}}{\ulambda^3\sigma^4}\cdot \frac{\op^2}{n^2}\right]}_{\text{from Step 1}} \\
&+  \underbrace{\sum_{j=1}^3\left\|\mcP_{U_j}A_j\mcP_{\left(U_{j+2}\otimes U_{j+1}\right)G_j^{\top}}\right\|_{\mathrm{F}} \cdot \frac{\sigma_\xi^2\oor^{1/2}}{\ulambda\sigma^2}\cdot \frac{\op}{n} + \sum_{j=1}^3\left\|\mcA \times_{j+1} U_{j+1} \times_{j+2} U_{j+2}\right\|_{\mathrm{F}} \cdot \left[\frac{\sigma_\xi^3\oor^{1/2}}{\ulambda^2\sigma^3}\left(\frac{\op\sqrt{\oor\log(\op)}}{n^{3/2}} \right)\right]}_{\text{shared between Step 1 and Step 2}} \\
&+  \sum_{j=1}^3 \left\|\mcP_{U_{j\perp}} A_j \mcP_{\left(U_{j+2} \otimes U_{j+1}\right) G_j^{\top}}\right\|_{\mathrm{F}}  \left(\frac{\sigma_{\xi}^2 \oor^{1 / 2}}{\ulambda \sigma^2} \cdot \frac{\sqrt{\op\oor\log(\op)}}{n}\right) 
+  \sum_{j=1}^3 \left\|\mcA \times_j U_j \right\|_{\mathrm{F}}  \left(\frac{\sigma_{\xi}^2\oor^{1/2}}{\ulambda \sigma^2} \cdot \frac{\oor\log(\op)}{n} +\frac{\sigma_{\xi}^2\oor^{1/2}}{\ulambda \sigma^2} \cdot \Delta \cdot \frac{\oR \log (\op)}{n}\right) \\
&+ \underbrace{\left\|\mcA\right\|_{\mathrm{F}} \cdot \left(\frac{\sigma_{\xi}^3}{\ulambda^2 \sigma^3} \cdot \frac{\oor^{3/2}\log(\op)^{3/2}}{n^{3 / 2}}+\frac{\sigma_{\xi}^3}{\ulambda^2 \sigma^3} \cdot \Delta \cdot \frac{\oR^{3 / 2} \log (\op)^{3 / 2}}{n^{3 / 2}}\right)}_{\text{from Step 2}}  \\
&+  \underbrace{\sum_{j=1}^3\left\|\mcP_{U_{j\perp}}A_j\mcP_{\left(U_{j+2}\otimes U_{j+1}\right)G_j^{\top}}\right\|_{\mathrm{F}}\cdot \Delta \cdot \frac{\sigma_\xi\oor^{1/2}}{\sigma} \sqrt{\frac{\oor\log(\op)}{n}} + \left\|\mcA\times_1 U_1\times_2 U_2\times_3 U_3\right\|_{\mathrm{F}}\cdot \Delta \cdot \frac{\sigma_\xi\oor^{1/2}}{\sigma} \sqrt{\frac{\oor\log(\op)}{n}}}_{\text{from Step 3}} \Bigg\}.
\end{align*}}%


\section{Preliminary Upper Bounds for Tensor Regression without Sample splitting}\label{sec: Preliminary Upper Bounds for Tensor Regression without Sample splitting}

In this section, we derive preliminary upper bounds for perturbation terms in tensor regression. The spectral representation of $\whU_j\whU_j^{\top} - U_jU_j^{\top}$, as developed in \citet{xia2021normal}, plays a pivotal role in analyzing negligible terms.

After the power iteration and projection in the algorithm without sample splitting in Section \ref{sec:debias_nonsplit}, for any $j=1,2,3$, we know that $\whU_j$ contains the top- $r_j$ eigenvectors of 
$$
\widehat{\mcT}^{\text{unbs}}_j\left(\mcP_{\whU_{j+2}^{(1)}} \otimes \mcP_{\whU_{j+1}^{(1)}}\right) \widehat{\mcT}^{\text{unbs}\top}_j.
$$
Consequently, $\whU_j \whU_j^{\top}$ is the spectral projector for the top-$r_j$ left eigenvectors of
\begin{align*}
\widehat{T}_j^{\text{unbs}} \left(\mcP_{\whU_{j+2}^{(1)}} \otimes \mcP_{\whU_{j+1}^{(1)}}\right) \widehat{T}_j^{\text{unbs}\top}
= T_j\left(\mcP_{U_{j+1}} \otimes \mcP_{U_{j+2}}\right) T_j^{\top}+\whE_j 
= U_jG_jG_j^{\top} U_j^{\top} + \whE_j,
\end{align*}
where
\begin{align}
\whE_j
= & T_j\left(\mcP_{\whU_{j+2}^{(1)}} \otimes \mcP_{\whU_{j+1}^{(1)}}\right)\whZ_j^{\top} + \whZ_j\left(\mcP_{\whU_{j+2}^{(1)}} \otimes \mcP_{\whU_{j+1}^{(1)}}\right)T_j^{\top} \notag \\
& +  T_j\left(\left(\mcP_{\whU_{j+2}^{(1)}}-\mcP_{U_{j+2}}\right) \otimes \mcP_{\whU_{j+1}^{(1)}}\right) T_j^{\top}+T_j\left(\mcP_{U_{j+2}} \otimes \left(\mcP_{\whU_{j+1}^{(1)}}-\mcP_{U_{j+1}}\right)\right) T_j^{\top} \label{eq: definition of whEj in tensor regression without sample splitting}\\
& +  T_j\left(\left(\mcP_{\whU_{j+2}^{(1)}}-\mcP_{U_{j+2}}\right) \otimes \left(\mcP_{\whU_{j+1}^{(1)}}-\mcP_{U_{j+1}}\right)\right) T_j^{\top}. \notag  
\end{align}

If $\left\|\whE_j\right\| \leq \frac{1}{2}\ulambda^2$, the spectral representation formula in Theorem 1 of \citet{xia2021normal} applies:
$$
\whU_j \whU_j^{\top}-U_j U_j^{\top}=\sum_{k_j= 1}^{+\infty}\mathcal{S}_{G_j, k_j}\left(\whE_j\right).
$$
Here, for any positive integer $k$
\begin{align}
\mathcal{S}_{G_j, k_j}\left(\whE_j\right)=\sum_{s_1+\cdots+s_{k_j+1}=k_j}(-1)^{1+\tau(\mathbf{s})} \cdot \mcP_j^{-s_1} \whE_j \mcP_j^{-s_2} \whE_j \mcP_j^{-s_3} \cdots \mcP_j^{-s_{k_j}} \whE_j \mcP_j^{-s_{k_j+1}} \label{eq: definition of Sj(whEj) in tensor regression without sample splitting}
\end{align}
with $s_1, \cdots, s_{k_j+1}$ being non-negative integers $\tau(\mathbf{s})=\sum_{j=1}^{k_j+1} \mathbb{I}\left(s_{k_j}>0\right)$, and
$$
\mcP_j^{-k}= U_j\left(G_jG_j^{\top}\right)^{-k}U_j^{\top}, \quad \text{for any}\quad k \geq 1\quad \text{and}\quad \mcP_j^{0}=U_{j\perp}U_{j\perp}^{\top}.
$$

For $k=1$, the leading term simplifies to:
\begin{align}
&\mathcal{S}_{G_j, 1}\left(\whE_j\right)
=  P_j^{-1} \whE_j P_j^{0}+P_j^{0} \whE_j P_j^{-1} \notag \\
= & U_j \left(G_j G_j^{\top}\right)^{-1} G_j \left({U_{j+1}} \otimes U_{j+2}\right)^{\top} \left(\mcP_{\whU_{j+1}^{(1)}}\otimes \mcP_{\whU_{j+2}^{(1)}}\right)\whZ_j^{\top} \mcP_{U_{j\perp}} \notag \\
& + \mcP_{U_{j\perp}} \whZ_j\left(\mcP_{\whU_{j+2}^{(1)}} \otimes \mcP_{\whU_{j+1}^{(1)}}\right) \left(U_{j+2} \otimes U_{j+1}\right)G_j^{\top}\left(G_j G_j^{\top}\right)^{-1} U_j^{\top}, \label{eq: definition of S1(whEj) in tensor regression without sample splitting}
\end{align}
for any $j=1,2,3$, where the second equality, the third inequality come from the definition that $P_j^{-1}=U_j\left(G_jG_j^{\top}\right)^{-1}U_j^{\top}$.

Here, note that $\left\|\whE_j\right\| \leq \kappa\ulambda \sqrt{\frac{\op}{n}}$. Then the condition, $\left\|\whE_j\right\| \leq \frac{1}{2}\ulambda^2$, for Theorem 1 in \citet{xia2021normal} is satisfied provied that $n \gtrsim \kappa^2 \op / \ulambda^2$.

In the subsequent sections, for any $j=1,2,3$, we further assume that the following events hold with high probability:
$$
\left\|\mcP_{\whU_j^{(0)}} - \mcP_{U_j}\right\| \leq \frac{\sigma_\xi}{\sigma}\sqrt{\frac{\op}{n}}
$$
holds with probability at least $1- \mathbb{P}\left(\mcE_{U}^{\text{reg}}\right)$, where event $\mcE_{U}^{\text{reg}}$ is defined by $\mcE_{U}^{\text{reg}}= \left\{\left\|\mcP_{\whU_j^{(0)}} - \mcP_{U_j}\right\| > \frac{\sigma_\xi}{\sigma}\sqrt{\frac{\op}{n}}\right\}$.
Then by Lemma~\ref{lemma: error contraction of l2 error of singular space in tensor regression}, we know that 
$\left\|\mcP_{\whU_j^{(1)}} - \mcP_{U_j}\right\| \leq \frac{\sigma_\xi}{\sigma}\sqrt{\frac{\op}{n}}$ and $\left\|\mcP_{\whU_j} - \mcP_{U_j}\right\| \leq \frac{\sigma_\xi}{\sigma}\sqrt{\frac{\op}{n}}$
holds with probability at least $1-\exp(-c\op) - \mathbb{P}(\mcE_{U}^{\text{reg}})$ for any $j=1,2,3$.

Finally, we assume the initial error bound:
$$
\left\|\whT^{\text{init}} - \mcT\right\|_{\mathrm{F}} \leq \Delta
$$
holds with probability at least $1-\mathbb{P}(\mcE_{\Delta})$, where event $\mcE_{\Delta}$ is given by $\mcE_{\Delta}= \left\{\left\|\whT^{\text{init}} - \mcT\right\|_{\mathrm{F}} > \Delta \right\}$.

In the following subsections, we established upper bounds for perturbation terms of varying orders in the spectral representation under the setting of tensor regression without sample splitting. In particular, we will show that the first-order perturbation term is the leading term. Throughout this section, we assume that $\Delta \geq (\sigma_\xi / \sigma)\sqrt{\op / n}$, which implies that the initial estimate satisfies the minimax lower bound. Therefore, the perturbation term introduced by the dependency between the initial estimate and debiasing procedure, in some sense, dominates the negligible terms.

\subsection{Preliminary Bounds in the Proof of Theorem~\ref{thm: main theorem in tensor regression without sample splitting}} \label{subsec: preliminary bounds in the proof of main theorems in tensor regression}

\begin{proposition}\label{prop: high-probability upper bound of AoPU1p(PUhat1-P_U1)oPU2p(PUhat2-P_U2)oPU3p(PUhat3-P_U3) in tensor regression without sample splitting}
Under the same setting of Theorem~\ref{thm: main theorem in tensor regression without sample splitting}, with probability at least $1-\exp(-cn) - \frac{1}{p^C} - \P\left(\mcE_\Delta\right) - \P\left(\mcE_{U}^{\text{reg}}\right)$, where $c$ and $C$ are two universal constants, the following bounds hold for any $j=1,2,3$:
\begin{align}
& \left\|\mcA \times_j \mcP_{U_{j\perp}}\left(\mcP_{\whU_j} - \mcP_{U_j}\right)\mcP_{U_j} \times_{j+1} \mcP_{U_{j+1\perp}}\left(\mcP_{\whU_{j+1}} - \mcP_{U_{j+1}}\right)\mcP_{U_{{j+1}}} \times_{j+2} \mcP_{U_{j+2\perp}}\left(\mcP_{\whU_{j+2}} - \mcP_{U_{j+2}}\right)\mcP_{U_{j+2}}  \right\|_{\mathrm{F}} \notag \\
\leq & \left\|\mcA\right\|_{\mathrm{F}} \cdot \left[\frac{\sigma_\xi^3 }{\ulambda^3\sigma^3} \left(\frac{\oor^{3/2}\log(\op)^{3/2}}{n^{3/2}} + \Delta\cdot \frac{\op^{1/2}\oR\log(\op)}{n^{3/2}} + \Delta^2 \cdot \frac{\op \oR\log (\op)}{n^{3 / 2}} + \Delta^3\cdot \frac{\op^{3/2}}{n^{3/2}}\right)\right] \label{eq: high-prob upper bound of AoPU1p(PUhat1-P_U1)PU1oPU2p(PUhat2-P_U2)PU2oPU3p(PUhat3-P_U3)PU3 in tensor regression without sample splitting} \\
& \left\|\mcA \times_j \mcP_{U_{j\perp}}\left(\mcP_{\whU_j} - \mcP_{U_j}\right)\mcP_{U_{j\perp}} \times_{j+1} \mcP_{U_{j+1\perp}}\left(\mcP_{\whU_{j+1}} - \mcP_{U_{j+1}}\right)\mcP_{U_{j+1}} \times_{j+2} \mcP_{U_{j+2\perp}}\left(\mcP_{\whU_{j+2}} - \mcP_{U_{j+2}}\right)\mcP_{U_{j+2}}  \right\|_{\mathrm{F}} \notag \\
\leq & \left\|\mcA\right\|_{\mathrm{F}} \cdot \left[\frac{\sigma_\xi^4 }{\ulambda^4\sigma^4} \left(\frac{\oor^3\op^{1/2}\log (\op)^{3/2}}{n^2} + \Delta\cdot \frac{\op\oR\log(\op)}{n^2} + \Delta^2\cdot \frac{\op^{3/2}\oR\log(\op)}{n^2} + \Delta^5\cdot \frac{\op^2}{n^2}\right)\right] \label{eq: high-prob upper bound of AoPU1p(PUhat1-P_U1)PU1poPU2p(PUhat2-P_U2)PU2oPU3p(PUhat3-P_U3)PU3 in tensor regression without sample splitting} \\
& \left\|\mcA \times_j \mcP_{U_{j\perp}}\left(\mcP_{\whU_j} - \mcP_{U_j}\right)\mcP_{U_j} \times_{j+1} \mcP_{U_{j+1\perp}}\left(\mcP_{\whU_{j+1}} - \mcP_{U_{j+1}}\right)\mcP_{U_{{j+1}\perp}} \times_{j+2} \mcP_{U_{j+2\perp}}\left(\mcP_{\whU_{j+2}} - \mcP_{U_{j+2}}\right)\mcP_{U_{j+2\perp}}  \right\|_{\mathrm{F}} \notag \\
\leq & \left\|\mcA\right\|_{\mathrm{F}} \cdot \left[\frac{\sigma_\xi^5 }{\ulambda^5\sigma^5} \left(\frac{\oor^{3/2}\op\log(\op)^{3/2}}{n^{5/2}} + \Delta\cdot \frac{\op^{3/2}\oR\log(\op)}{n^{5/2}} + \Delta^2\cdot \frac{\op^2\oR\log(\op)}{n^{5/2}} + \Delta^3\cdot \frac{\op^{5/2}}{n^{5/2}}\right)\right] \label{eq: high-prob upper bound of AoPU1p(PUhat1-P_U1)PU1oPU2p(PUhat2-P_U2)PU2poPU3p(PUhat3-P_U3)PU3p in tensor regression without sample splitting}\\
& \left\|\mcA \times_j \mcP_{U_{j\perp}}\left(\mcP_{\whU_j} - \mcP_{U_j}\right)\mcP_{U_{j\perp}} \times_{j+1} \mcP_{U_{j+1\perp}}\left(\mcP_{\whU_{j+1}} - \mcP_{U_{j+1}}\right)\mcP_{U_{{j+1}\perp}} \times_{j+2} \mcP_{U_{j+2\perp}}\left(\mcP_{\whU_{j+2}} - \mcP_{U_{j+2}}\right)\mcP_{U_{j+2\perp}}  \right\|_{\mathrm{F}} \notag \\
\leq & \left\|\mcA\right\|_{\mathrm{F}} \cdot \left[\frac{\sigma_\xi^6 }{\ulambda^6\sigma^6} \left(\frac{\oor^{3/2}\op^{3/2}\log(\op)^{3/2}}{n^3} + \Delta\cdot \frac{\op^2\oR\log(\op)}{n^3} + \Delta^2 \cdot \frac{\op^{5/2}\oR\log(\op)}{n^3} + \Delta^3\cdot \frac{\op^3}{n^3}\right)\right] \label{eq: high-prob upper bound of AoPU1p(PUhat1-P_U1)PU1oPU2p(PUhat2-P_U2)PU2poPU3p(PUhat3-P_U3)PU3pp in tensor regression without sample splitting}.
\end{align}
\end{proposition}

\begin{proof}

By symmetry, it suffices to consider
\begin{align*}
\mathrm{\RN{1}} 
= & \left\|\left(\mcP_{U_3}\left(\mcP_{\whU_3}-\mcP_{U_3}\right) \mcP_{U_{3 \perp}} \otimes \mcP_{U_2}\left(\mcP_{\whU_2}-\mcP_{U_2}\right) \mcP_{U_{2 \perp}}\right) A_1^{\top} \mcP_{U_{1 \perp}}\left(\mcP_{\whU_1}-\mcP_{U_1}\right) \mcP_{U_1}\right\|_{\mathrm{F}} \\
\mathrm{\RN{2}} 
= & \left\|\left(\mcP_{U_3}\left(\mcP_{\whU_3}-\mcP_{U_3}\right) \mcP_{U_{3 \perp}} \otimes \mcP_{U_2}\left(\mcP_{\whU_2}-\mcP_{U_2}\right) \mcP_{U_{2 \perp}}\right) A_1^{\top} \mcP_{U_{1 \perp}}\left(\mcP_{\whU_1}-\mcP_{U_1}\right) \mcP_{U_1}\right\|_{\mathrm{F}} \\
\mathrm{\RN{3}} 
= & \left\|\left(\mcP_{U_3}\left(\mcP_{\whU_3}-\mcP_{U_3}\right) \mcP_{U_{3 \perp}} \otimes \mcP_{U_2}\left(\mcP_{\whU_2}-\mcP_{U_2}\right) \mcP_{U_{2 \perp}}\right) A_1^{\top} \mcP_{U_{1 \perp}}\left(\mcP_{\whU_1}-\mcP_{U_1}\right) \mcP_{U_{1 \perp}}\right\|_{\mathrm{F}} \\
\mathrm{\RN{4}} 
= & \left\|\left(\mcP_{U_{3 \perp}}\left(\mcP_{\whU_3}-\mcP_{U_3}\right) \mcP_{U_{3 \perp}} \otimes \mcP_{U_{2 \perp}}\left(\mcP_{\whU_2}-\mcP_{U_2}\right) \mcP_{U_{2 \perp}}\right) A_1^{\top} \mcP_{U_{1 \perp}}\left(\mcP_{\whU_1}-\mcP_{U_1}\right) \mcP_{U_1}\right\|_{\mathrm{F}}
\end{align*}

Note that
\begin{align*}
\mathrm{\RN{1}} 
\leq & \underbrace{\left\|\left(\mcP_3^{-1} \whE_3 \mcP_{U_{3 \perp}}\otimes \mcP_2^{-1} \whE_2 \mcP_{U_{2 \perp}}\right)\left(\mcP_{U_3} \otimes \mcP_{U_{2 \perp}}\right) A_1^{\top} \mcP_{U_{1 \perp}} \whE_1 \mcP_1^{-1}\right\|_{\mathrm{F}}}_{\eqref{eq: high-prob upper bound of AoPU1pEhat1P1(-1)oPU2pEhat2P2(-1)oPU3pEhat3P3(-1) in tensor regression without sample splitting}} \\
+ & \sum_{j=1}^3 \underbrace{\left\|\mcP_{j+2}^{-1} \whE_{j+2} \mcP_{U_{j+2 \perp}}  V_{j+2}\right\|}_{\eqref{eq: high-prob upper bound of V1tP1(0)Ehat1P1(-1/2) in tensor regression without sample splitting}} \cdot \underbrace{\left\|\mcP_{j+1}^{-1} \whE_{j+1} \mcP_{U_{2 \perp}}  V_{j+1}\right\|}_{\eqref{eq: high-prob upper bound of V1tP1(0)Ehat1P1(-1/2) in tensor regression without sample splitting}} \cdot\left\|A_j^{\top} \right\|_{\mathrm{F}} \cdot \underbrace{\left\|\mcP_{U_j} \sum_{k_j=2}^{+\infty} S_{G_j, k_j}\left(\whE_j\right) \mcP_{U_{j \perp}} V_j\right\|}_{\eqref{eq: high-prob upper bound of PUporder2PUpPV in tensor regression without sample splitting}} \\
+ & \sum_{j=1}^3 \underbrace{\left\|\mcP_{U_{j+2}} \sum_{k_{j+2}=2}^{+\infty} S_{G_{j+2}, k_{j+2}}\left(\whE_{j+2}\right) \mcP_{U_{{j+2} \perp}}  V_{j+2}\right\|}_{\eqref{eq: high-prob upper bound of PUporder2PUpPV in tensor regression without sample splitting}} \cdot \underbrace{\left\|\mcP_{U_{j+1}} \sum_{k_{j+1}=2}^{+\infty} S_{G_{j+1}, k_{j+1}}\left(\whE_{j+1}\right) \mcP_{U_{{j+1} \perp}}  V_{j+1}\right\|}_{\eqref{eq: high-prob upper bound of PUporder2PUpPV in tensor regression without sample splitting}} \\
& \cdot \left\|\mcA\right\|_{\mathrm{F}} \cdot \underbrace{\left\|V_j^{\top}  \mcP_{U_{j \perp}} \whE_j \mcP_j^{-1}\right\|}_{\eqref{eq: high-prob upper bound of V1tP1(0)Ehat1P1(-1/2) in tensor regression without sample splitting}} \\
+ & \underbrace{\left\|\mcP_{U_{j+2}} \sum_{k_{j+2}=2}^{+\infty} S_{G_{j+2}, k_{j+2}}\left(\whE_{j+2}\right) \mcP_{U_{{j+2} \perp}}  V_{j+2}\right\|}_{\eqref{eq: high-prob upper bound of PUporder2PUpPV in tensor regression without sample splitting}} \cdot \underbrace{\left\|\mcP_{U_{j+2}} \sum_{k_{j+1}=2}^{+\infty} S_{G_{j+1}, k_{j+1}}\left(\whE_{j+1}\right) \mcP_{U_{{j+1} \perp}}  V_{j+1}\right\|}_{\eqref{eq: high-prob upper bound of PUporder2PUpPV in tensor regression without sample splitting}} \\
& \cdot\left\|\mcA\right\|_{\mathrm{F}} \cdot \underbrace{\left\|\mcP_{U_{j+1}} \sum_{k_{j+1}=2}^{+\infty} S_{G_{j+1}, k_{j+1}}\left(\whE_{j+1}\right) \mcP_{U_{{j+1} \perp}}  V_{j+1}\right\|}_{\eqref{eq: high-prob upper bound of PUporder2PUpPV in tensor regression without sample splitting}} \\
\lesssim & \left\|\mcA\right\|_{\mathrm{F}} \cdot \left[\frac{\sigma_{\xi}^3}{\ulambda^3 \sigma^3} \cdot\left(\frac{\oor^{3/2}\log(\op)^{3/2}}{n^{3/2}}+\Delta \cdot \frac{\op^{1/2} \oR \log (\op)}{n^{3/2}}+\Delta^3 \cdot \frac{\op^{3/2} \oR \log (\op)}{n^{3/2}}\right)\right].
\end{align*}

Applying similar arguments, we obtain the bounds for \RN{2}, \RN{3}, \RN{4}.

\end{proof}

\begin{proposition}\label{prop: high-probability upper bound of AoPU1oPU2p(PUhat2-P_U2)oPU3p(PUhat3-P_U3) without sample splitting}

Under the same setting of Theorem~\ref{thm: main theorem in tensor regression without sample splitting}, with probability at least $1-\exp(-cn) - \frac{1}{p^C} - \P\left(\mcE_\Delta\right) - \P\left(\mcE_{U}^{\text{reg}}\right)$, where $c$ and $C$ are two universal constants, the following bounds hold for any $j=1,2,3$:
\begin{align}
& \left\|\mcA \times_j \mcP_{U_j} \times_{j+1} \mcP_{U_{j+1\perp}}\left(\mcP_{\whU_{j+1}} - \mcP_{U_{j+1}}\right)\mcP_{U_{j+1}} \times_{j+2} \mcP_{U_{j+2\perp}}\left(\mcP_{\whU_{j+2}} - \mcP_{U_{j+2}}\right)\mcP_{U_{j+2}}\right\|_{\mathrm{F}} \notag \\
\lesssim & \left\|\mcA \times_j U_j\right\|_{\mathrm{F}} \cdot \left[\frac{\sigma_{\xi}^2}{\ulambda^2 \sigma^2} \cdot\left(\frac{\oor\log(\op)}{n}+\Delta \cdot \frac{\sqrt{\oR \op\log(\op)}}{n}+\Delta^2 \cdot \frac{\op}{n}\right)\right] \label{eq: high-prob upper bound of AoPU1oPU2p(PUhat2-P_U2)PU2oPU3p(PUhat3-P_U3)PU3 in tensor regression without sample splitting}\\
& \left\|\mcA \times_j \mcP_{U_j} \times_{j+1} \mcP_{U_{j+1\perp}}\left(\mcP_{\whU_{j+1}} - \mcP_{U_{j+1}}\right)\mcP_{U_{j+1\perp}} \times_{j+2} \mcP_{U_{j+2\perp}}\left(\mcP_{\whU_{j+2}} - \mcP_{U_{j+2}}\right)\mcP_{U_{j+2}}\right\|_{\mathrm{F}} \notag \\
\lesssim & \left\|\mcA \times_j U_j\right\|_{\mathrm{F}} \cdot\left[\frac{\sigma_\xi^3}{\ulambda^3\sigma^3} \cdot \left(\frac{\oor^2\op^{1/2} \log (\op)}{n^{3/2}}+\Delta \cdot \frac{\op\sqrt{\oR \log(\op)}}{n^{3/2}}+\Delta^2 \cdot \frac{\op^{3/2}}{n^{3/2}}\right)\right] \label{eq: high-prob upper bound of AoPU1oPU2p(PUhat2-P_U2)PU2poPU3p(PUhat3-P_U3)PU3 in tensor regression without sample splitting} \\
& \left\|\mcA \times_j \mcP_{U_j} \times_{j+1} \mcP_{U_{j+1\perp}}\left(\mcP_{\whU_{j+1}} - \mcP_{U_{j+1}}\right)\mcP_{U_{j+1\perp}} \times_{j+2} \mcP_{U_{j+2\perp}}\left(\mcP_{\whU_{j+2}} - \mcP_{U_{j+2}}\right)\mcP_{U_{j+2\perp}}\right\|_{\mathrm{F}} \notag \\
\lesssim & \left\|\mcA \times_j U_j\right\|_{\mathrm{F}} \cdot \left[\frac{\sigma_{\xi}^4}{\ulambda^4\sigma^4} \cdot \frac{\op}{n} \cdot \left(\frac{\oor\log(\op)}{n}+\Delta \cdot \frac{\sqrt{\oR \op\log(\op)}}{n}+\Delta^2 \cdot \frac{\op}{n}\right)\right]. \label{eq: high-prob upper bound of AoPU1oPU2p(PUhat2-P_U2)PU2poPU3p(PUhat3-P_U3)PU3p in tensor regression without sample splitting}
\end{align}

\end{proposition}

\begin{proof}

By symmetry, it suffices to consider
\begin{align*}
\mathrm{\RN{1}} = & \left\|\left(\mcP_{U_3} \otimes \mcP_{U_2}\left(\mcP_{\whU_2}-\mcP_{U_2}\right) \mcP_{U_{2 \perp}}\right) A_1^{\top} \mcP_{U_{1 \perp}}\left(\mcP_{\whU_1}-\mcP_{U_1}\right) \mcP_{U_1}\right\|_{\mathrm{F}} \\
\mathrm{\RN{2}} = & \left\|\left(\mcP_{U_3} \otimes \mcP_{U_2}\left(\mcP_{\whU_2}-\mcP_{U_2}\right) \mcP_{U_{2 \perp}}\right) A_1^{\top} \mcP_{U_{1 \perp}}\left(\mcP_{\whU_1}-\mcP_{U_1}\right) \mcP_{U_{1 \perp}}\right\|_{\mathrm{F}} \\
\mathrm{\RN{3}} = & \left\|\left(\mcP_{U_3} \otimes \mcP_{U_{2 \perp}}\left(\mcP_{\whU_2}-\mcP_{U_2}\right) \mcP_{U_{2 \perp}}\right) A_1^{\top} \mcP_{U_{1 \perp}}\left(\mcP_{\whU_1}-\mcP_{U_1}\right) \mcP_{U_1 \perp}\right\|_{\mathrm{F}} .
\end{align*}

Note that
\begin{align*}
\mathrm{\RN{1}} 
\leq & \underbrace{\left\|\left(\mcP_{U_3}\otimes \mcP_2^{-1}\whE_2\mcP_{U_{2\perp}}\right)A_1^{\top}\mcP_{U_{1\perp}}\whE_1\mcP_1^{-1}\right\|_{\mathrm{F}}}_{\eqref{eq: high-prob upper bound of AoPU1oPU2pEhat2P2(-1)oPU3pEhat3P3(-1) in tensor regression without sample splitting}} 
+  \underbrace{\left\|\mcP_{U_2} \sum_{k_2=2}^{+\infty} S_{G_2, k_2}\left(\whE_2\right) \mcP_{U_{2 \perp}} V_2\right\|}_{\eqref{eq: high-prob upper bound of PUorder2PUpPV in tensor regression without sample splitting}} \cdot \left\|\mcA \times_3 U_3\right\|_{\mathrm{F}} \cdot \underbrace{\left\|V_1^{\top}\mcP_{U_{1\perp}}\whE_1\mcP_1^{-1}\right\|}_{\eqref{eq: high-prob upper bound of V1tP1(0)Ehat1P1(-1/2) in tensor regression without sample splitting}} \\
+ & \underbrace{\left\|\mcP_2^{-1} \whE_2 \mcP_{U_{2 \perp}}V_2\right\|}_{\eqref{eq: high-prob upper bound of V1tP1(0)Ehat1P1(-1/2) in tensor regression without sample splitting}} \cdot \left\|\mcA \times_3 U_3\right\|_{\mathrm{F}} \cdot \underbrace{\left\|\mcP_{U_1} \sum_{k_1=2}^{+\infty} S_{G_1, k_1}\left(\whE_1\right) \mcP_{U_{1 \perp}} V_1\right\|}_{\eqref{eq: high-prob upper bound of PUorder2PUpPV in tensor regression without sample splitting}} \\
+ & \underbrace{\left\|\mcP_{U_2} \sum_{k_2=2}^{+\infty} S_{G_2, k_2}\left(\whE_2\right) \mcP_{U_{2 \perp}} V_2\right\|}_{\eqref{eq: high-prob upper bound of PUorder2PUpPV in tensor regression without sample splitting}} \cdot \left\|\mcA \times_3 U_3\right\|_{\mathrm{F}} \cdot \underbrace{\left\|\mcP_{U_1} \sum_{k_1=2}^{+\infty} S_{G_1, k_1}\left(\whE_1\right) \mcP_{U_{1 \perp}} V_1\right\|}_{\eqref{eq: high-prob upper bound of PUorder2PUpPV in tensor regression without sample splitting}} \\
\lesssim & \left\|\mcA \times_3 U_3\right\|_{\mathrm{F}} \cdot \left[\frac{\sigma_{\xi}^2}{\ulambda^2 \sigma^2} \cdot\left(\frac{\oor\log(\op)}{n}+\Delta \cdot \frac{\sqrt{\oR \op\log(\op)}}{n}+\Delta^2 \cdot \frac{\op}{n}\right)\right] .
\end{align*}

Applying similar arguments, we obtain the bounds for \RN{2}, \RN{3}.

\end{proof}

\begin{proposition} \label{prop: high-prob upper bound of perturbation after projection in tensor regression without sample splitting}

Under the same setting of Theorem~\ref{thm: main theorem in tensor regression without sample splitting}, let $V_j \in \mathbb{R}^{p_j\times R_j}$ be a fixed matrix satisfying $\left\|V_j\right\|=1$ for any $j=1,2,3$. Then, for any $j=1,2,3$, the following bounds hold with probability at least $1-\exp(-cn) - \frac{1}{p^C} - \P\left(\mcE_\Delta\right) - \P\left[\left(\mcE_{U}^{\text{reg}}\right)\right]$, where $c$ and $C$ are two universal constants:
\begin{align}
& \left\|V_j^{\top}\mcP_{U_{j\perp}}\left(\mcP_{\whU_j} - \mcP_{U_j}\right)U_j\right\| \lesssim \frac{\sigma_{\xi}}{\ulambda \sigma} \cdot\left(\sqrt{\frac{\oR \log (\op)}{n}}+\Delta \sqrt{\frac{\op}{n}}\right), \label{eq: high-prob upper bound of V1tPU1p(PUhat1-PU1)U1 in tensor regression without sample splitting}\\
& \left\|V_j^{\top}\mcP_{U_{j\perp}}\left(\mcP_{\whU_j} - \mcP_{U_1}\right)U_{j\perp}\right\| \lesssim \frac{\sigma_\xi^2}{\ulambda^2\sigma^2}\cdot \left(\frac{\sqrt{\oR\log\left(\op\right)}\cdot \sqrt{\op}}{n}+ \cdot \Delta\cdot \frac{\op}{n} \right). \label{eq: high-prob upper bound of V1tPU1p(PUhat1-PU1)U1p in tensor regression without sample splitting}
\end{align}

Furthermore, 
\begin{align}
& \left\|V_j^{\top}\mcP_{U_{j\perp}}\left(\mcP_{\whU_j} - \mcP_{U_j}\right)\right\| \lesssim \frac{\sigma_\xi}{\ulambda\sigma} \cdot \left(\sqrt{\frac{\oR\log (\op)}{n}}+\Delta \sqrt{\frac{\op}{n}}\right) + \frac{\sigma_{\xi}^2}{\ulambda^2 \sigma^2} \cdot \frac{\op^{3/2}}{n^{3/2}}. \label{eq: high-prob upper bound of V1tPU1p(PUhat1-PU1) in tensor regression without sample splitting}
\end{align}

In addition, we have
\begin{align}
\left\|U_j^\top\left(\mcP_{\whU_j} - \mcP_{U_j}\right)U_j\right\| \lesssim \frac{\sigma_\xi^2}{\ulambda^2\sigma^2}\cdot \frac{\op}{n}. \label{eq: high-prob upper bound of PU1(PUhat1-PU1)PU1 in tensor regression}
\end{align}

Similar bounds hold when $\mcP_{\whU_j}$ is replaced by $\mcP_{\whU_j^{(1)}}$.

\end{proposition}

\begin{proof}

By symmetry, it suffices to consider $\left\|V_1^{\top}\mcP_{U_{1\perp}}\left(\mcP_{\whU_1^{(1)}} - \mcP_{U_1}\right)U_1\right\|$ and $\left\|V_1^{\top}\mcP_{U_{1\perp}}\left(\mcP_{\whU_1^{(1)}} - \mcP_{U_1}\right)U_{1\perp}\right\|$.

For the first inequality,
\begin{align*}
& \left\|V_1^{\top}\mcP_{U_{1\perp}}\left(\mcP_{\whU_1^{(1)}} - \mcP_{U_1}\right)U_1\right\| 
=  \left\|V_1^{\top}\mcP_{U_{1\perp}}\sum_{k_1=1}^{+\infty}S_{G_1, k_1}\left(\whE_1\right)U_1\right\| \\
\leq & \left\|V_1^{\top}\mcP_{U_{1\perp}}S_{G_1, 1}\left(\whE_1\right) U_1\right\|+ \left\|V_1^{\top}\mcP_{U_{1\perp}}S_{G_1, 2}\left(\whE_1\right) U_1\right\| + \left\|V_1^{\top}\mcP_{U_{1\perp}}\sum_{k_1=3}^{+\infty}S_{G_1, k_1}\left(\whE_1\right) U_1\right\| \\
\leq & \underbrace{\left\|V_1^{\top}\mcP_1^{0}\whE_1\mcP_1^{-1}U_1\right\|}_{\eqref{eq: high-prob upper bound of V1tP1(0)Ehat1P1(-1/2) in tensor regression without sample splitting}} + \underbrace{\left\|V_1^{\top}\mcP_1^{0}\whE_1\mcP_1^{-\frac{1}{2}}\right\|}_{\eqref{eq: high-prob upper bound of V1tP1(0)Ehat1P1(-1/2) in tensor regression without sample splitting}}  \underbrace{\left\|\mcP_1^{-\frac{1}{2}}\whE_1\mcP_1^{-1}U_1\right\|}_{\eqref{eq: high-prob upper bound of P1(-1/2)Ehat1P1(-1/2) in tensor regression without sample splitting}} \\
+ &  \underbrace{\left\|V_1^{\top}\mcP_1^{0}\whE_1\mcP_1^{0}\right\|}_{\eqref{eq: high-prob upper bound of V1tP1(0)Ehat1P1(0) in tensor regression without sample splitting}}  \underbrace{\left\|\mcP_1^{0}\whE_1\mcP_1^{-2}U_1\right\|}_{\eqref{eq: high-prob upper bound of P1(0)Ehat1P1(-1/2) in tensor regression}} 
 + \left\|V_1^{\top}\mcP_{U_{1\perp}}\sum_{k_1=3}^{+\infty}S_{G_1, k_1}\left(\whE_1\right) U_1\right\| \\
\lesssim & \frac{\sigma_\xi}{\ulambda\sigma} \cdot \left(\sqrt{\frac{\oR\log (\op)}{n}}+\Delta \sqrt{\frac{\op}{n}}\right).
\end{align*}

Applying similar arguments, we obtain the bounds for the rest quantities.





\end{proof}

\subsection{Upper Bound of First-Order Perturbation Terms}

\begin{lemma} \label{lemma: high-prob upper bound of the first-order perturbation terms in tensor regression without sample splitting} 

Under the same setting of Theorem~\ref{thm: main theorem in tensor regression without sample splitting}, let $\whE_j$ be defined by \eqref{eq: definition of whEj in tensor regression without sample splitting} and $\mcP_j^{-s}=U_j\left(G_jG_j^{\top}\right)^{-2s}U_j^{\top}$ for any $s>0$ and $j=1,2,3$ while $\mcP_j^{0} = U_{j\perp}U_{j\perp}^{\top}$. Then for any $j=1,2,3$, with probability at least $1-\exp(-cn) - \frac{1}{p^C} - \P\left(\mcE_\Delta\right) - \P\left[\left(\mcE_{U}^{\text{reg}}\right)\right]$, it follows that:
\begin{align}
\left\|\mcP_j^{-\frac{1}{2}}\whE_j\mcP_j^{-\frac{1}{2}}\right\| & = \left\|U_j^{\top}\left(G_jG_j^{\top}\right)^{-\frac{1}{2}}U_j^{\top}\whE_jU_j\left(G_jG_j^{\top}\right)^{-\frac{1}{2}}U_j^{\top}\right\| \lesssim \frac{\sigma_{\xi}}{\ulambda\sigma}\cdot \left(\sqrt{\frac{\oor\log(\op)}{n}} +\Delta\sqrt{\frac{\op}{n}}\right), \label{eq: high-prob upper bound of P1(-1/2)Ehat1P1(-1/2) in tensor regression without sample splitting}\\
\left\|\mcP_j^{0}\whE_j\mcP_j^{-\frac{1}{2}}\right\|
& = \left\|U_{j\perp}^{\top}\whE_jU_j\left(G_jG_j^{\top}\right)^{-\frac{1}{2}}U_j^{\top}\right\| \lesssim \frac{\sigma_\xi}{\sigma}\cdot\sqrt{\frac{\op}{n}}, \label{eq: high-prob upper bound of P1(0)Ehat1P1(-1/2) in tensor regression} \\
\left\|\mcP_j^{0}\whE_j\mcP_j^{0}\right\| 
& = \left\|U_{j\perp}^{\top}\whE_jU_{j\perp}\right\|\lesssim \frac{\sigma_\xi^2}{\sigma^2}\cdot\frac{\op}{n}  \label{eq: high-prob upper bound of P1(0)Ehat1P1(0) in tensor regression}
\end{align}
where $c$ and $C$ are two universal constants . 

Furthermore, \eqref{eq: high-prob upper bound of P1(0)Ehat1P1(-1/2) in tensor regression} and \eqref{eq: high-prob upper bound of P1(0)Ehat1P1(0) in tensor regression} hold under the conditions of Theorem~\ref{thm: main theorem in tensor regression with sample splitting}, where $\mcP_{\whU_j}$'s are the output of the algorithm with sample splitting in Section \ref{sec: Debiased Estimator of Linear Functionals with Sample Splitting}, as well.

\end{lemma}

\begin{proof}[Proof of Lemma~\ref{lemma: high-prob upper bound of the first-order perturbation terms in tensor regression without sample splitting}]

By symmetry, it suffices to consider
$$
\left\|\mcP_1^{-\frac{1}{2}}\whE_1\mcP_1^{-\frac{1}{2}} \right\|, \left\|U_{1\perp}^{\top}\whE_1\mcP_1^{-\frac{1}{2}} \right\|,\left\|U_{1\perp}^{\top}\whE_1U_{1\perp} \right\|.
$$

Let $\whE_1$ be defined as in \eqref{eq: definition of whEj in tensor regression without sample splitting}.
We have,

\begin{align}
\left\|\mcP_1^{-\frac{1}{2}}\whE_1\mcP_1^{-\frac{1}{2}}\right\|
= & \underbrace{\left\|U_1\left(G_1G_1\right)^{-\frac{1}{2}}G_1\left(U_3^{\top}\mcP_{\whU_3^{(1)}} \otimes U_2^{\top}\mcP_{\whU_2^{(1)}}\right) \whZ_1^{\top}\mcP_1^{-\frac{1}{2}}\right\|}_{\mathrm{\RN{1}}} \label{eq: term 1 in P1(-1/2)Ehat1P1(-1/2) in tensor regression without sample splitting}\\
+ & \underbrace{\left\|\mcP_1^{-\frac{1}{2}} \whZ_1\left(\mcP_{\whU_3^{(1)}}U_3 \otimes \mcP_{\whU_2^{(1)}}U_2\right) G_1^{\top}\left(G_1G_1\right)^{-\frac{1}{2}}U_1^{\top}\right\|}_{\mathrm{\RN{2}}} \label{eq: term 2 in P1(-1/2)Ehat1P1(-1/2) in tensor regression without sample splitting}\\
+ & \underbrace{\left\|\mcP_1^{-\frac{1}{2}} \whZ_1\left(\mcP_{\whU_3^{(1)}} \otimes \mcP_{\whU_2^{(1)}}\right) \whZ_1^{\top}\mcP_1^{-\frac{1}{2}}\right\|}_{\mathrm{\RN{3}}} \label{eq: term 3 in P1(-1/2)Ehat1P1(-1/2) in tensor regression without sample splitting}\\
+ & \underbrace{\left\|U_1\left(G_1G_1\right)^{-\frac{1}{2}}G_1\left[U_3^{\top}\left(\mcP_{\whU_3^{(1)}}-\mcP_{U_3}\right)U_3 \otimes U_2^{\top}\mcP_{U_2}U_2\right] G_1^{\top}\left(G_1G_1\right)^{-\frac{1}{2}}U_1^{\top}\right\|}_{\mathrm{\RN{4}}} \label{eq: term 4 in P1(-1/2)Ehat1P1(-1/2) in tensor regression without sample splitting}\\
+ & \underbrace{\left\|U_1\left(G_1G_1\right)^{-\frac{1}{2}}G_1\left[U_3^{\top}\mcP_{U_3}U_3 \otimes U_2^{\top}\left(\mcP_{\whU_2^{(1)}}-\mcP_{U_2}\right)U_2\right] G_1^{\top}\left(G_1G_1\right)^{-\frac{1}{2}}U_1^{\top}\right\|}_{\mathrm{\RN{5}}} \label{eq: term 5 in P1(-1/2)Ehat1P1(-1/2) in tensor regression without sample splitting} \\
+ & \underbrace{\left\|U_1\left(G_1G_1\right)^{-\frac{1}{2}}G_1\left[U_3^{\top}\left(\mcP_{\whU_3^{(1)}}-\mcP_{U_3}\right)U_3 \otimes U_2^{\top}\left(\mcP_{\whU_2^{(1)}}-\mcP_{U_2}\right)U_2\right] G_1^{\top}\left(G_1G_1\right)^{-\frac{1}{2}}U_1^{\top}\right\|}_{\mathrm{\RN{6}}} \label{eq: term 6 in P1(-1/2)Ehat1P1(-1/2) in tensor regression without sample splitting} . 
\end{align}
For the first term \RN{1} \eqref{eq: term 1 in P1(-1/2)Ehat1P1(-1/2) in tensor regression without sample splitting}, by decomposing $\mcP_{\whU_j^{(1)}} = \mcP_{\whU_j^{(1)}} - \mcP_{U_j} + \mcP_{U_j}$, we have 
\begin{align}
\mathrm{\RN{1}} 
\leq & \left\|\left(U_3 \otimes U_2\right)^{\top} \whZ_1^{\top}U_1\right\| \cdot \left\|U_1\left(G_1G_1\right)^{-\frac{1}{2}}U_1^{\top}\right\| \cdot \left\|\left(G_1G_1\right)^{-\frac{1}{2}}U_1^{\top}\right\| \notag \\
+ & \sup_{\substack{W_2 \in \mathbb{R}^{p_2\times r_3}, \left\|W_2\right\|=1\\W_3 \in \mathbb{R}^{p_3\times r_3}, \left\|W_3\right\|=1}} \left\|U_1^{\top}\whZ_1\left(W_3\otimes W_2\right)\right\|  \left\|\left(G_1G_1\right)^{-\frac{1}{2}}U_1^{\top}\right\|  \notag\\
&\cdot \left(\left\|\mcP_{\whU_2^{(1)}}- \mcP_{U_2}\right\| + \left\|\mcP_{\whU_3^{(1)}}- \mcP_{U_3}\right\| + \left\|\mcP_{\whU_2^{(1)}}- \mcP_{U_2}\right\| \cdot \left\|\mcP_{\whU_3^{(1)}}- \mcP_{U_3}\right\|\right) \notag \\
\lesssim & \frac{\sigma_{\xi}}{\ulambda\sigma}\cdot \left(\sqrt{\frac{\oor\log(\op)}{n}} +\Delta\sqrt{\frac{\op}{n}}\right). \label{eq: upper bound of term 1 in P1(-1/2)Ehat1P1(-1/2) in tensor regression without sample splitting}
\end{align}

Note that the Frobenius norms of terms \RN{1} \eqref{eq: term 1 in P1(-1/2)Ehat1P1(-1/2) in tensor regression without sample splitting} and \RN{2} \eqref{eq: term 2 in P1(-1/2)Ehat1P1(-1/2) in tensor regression without sample splitting} are equal. Then for the second term \eqref{eq: term 2 in P1(-1/2)Ehat1P1(-1/2) in tensor regression without sample splitting}, we have
\begin{align}
\mathrm{\RN{2}} & \lesssim \frac{\sigma_{\xi}}{\sigma}\cdot \left(\sqrt{\frac{\oor\log(\op)}{n}} +\Delta\sqrt{\frac{\op}{n}}\right). \label{eq: upper bound of term 2 in P1(-1/2)Ehat1P1(-1/2) in tensor regression without sample splitting}
\end{align}
For the third term \eqref{eq: term 3 in P1(-1/2)Ehat1P1(-1/2) in tensor regression without sample splitting}, we have
\begin{align}
\mathrm{\RN{3}} 
\leq & \frac{1}{\ulambda^2} \cdot \underbrace{\left\|U_1^{\top}\whZ_1\left(\mcP_{U_3}\otimes \mcP_{U_2}\right)\right\|}_{\eqref{eq: high-prob upper bound of U1Zhat1(U3oU2) in tensor regression without sample splitting}}^2 \cdot \left\|\left(G_1 G_1\right)^{-\frac{1}{2}}U_1^{\top}\right\| \notag \\
+ & \frac{1}{\ulambda} \cdot \sup_{\substack{W_2 \in \mathbb{R}^{p_2\times 2r_2}, \left\|W_2\right\|=1\\W_3 \in \mathbb{R}^{p_3\times 2r_3}, \left\|W_3\right\|=1}}\left\|U_1\whZ_1\left(W_3\otimes W_2\right)\right\| \cdot \sup_{\substack{W_2 \in \mathbb{R}^{p_2\times 2r_2}, \left\|W_2\right\|=1\\W_3 \in \mathbb{R}^{p_3\times 2r_3}, \left\|W_3\right\|=1}}\left\|U_1\whZ_1\left(W_3\otimes W_2\right)\right\| \notag \\
& \cdot \left(\left\|\mcP_{\whU_2^{(1)}}-\mcP_{U_2}\right\| + \left\|\mcP_{\whU_3^{(1)}}-\mcP_{U_3}\right\| + \left\|\mcP_{\whU_2^{(1)}}-\mcP_{U_2}\right\| \cdot \left\|\mcP_{\whU_3^{(1)}}-\mcP_{U_3}\right\|\right) \notag \\
\lesssim & \frac{\sigma_\xi^2}{\ulambda^2\sigma^2}\cdot \left(\frac{\oor\log(\op)}{n}+ \Delta^2\cdot \frac{\op}{n}\right). \label{eq: upper bound of term 3 in P1(-1/2)Ehat1P1(-1/2) in tensor regression without sample splitting}
\end{align}

For the fourth term \eqref{eq: term 4 in P1(-1/2)Ehat1P1(-1/2) in tensor regression without sample splitting}, we have
\begin{align}
\mathrm{\RN{4}}
\leq & \left\|U_1\left(G_1G_1\right)^{-\frac{1}{2}}G_1\right\|\cdot \left\|U_2^{\top}\left(\mcP_{\whU_2^{(1)}}-\mcP_{U_2}\right)U_2\right\| \cdot \left\|G_1^{\top}\left(G_1 G_1\right)^{-\frac{1}{2}}\right\| \lesssim \frac{\sigma_\xi^2}{\ulambda^2\sigma^2}\cdot \frac{\op}{n}. \label{eq: upper bound of term 4 in P1(-1/2)Ehat1P1(-1/2) in tensor regression without sample splitting}
\end{align}

For the fifth term \eqref{eq: term 5 in P1(-1/2)Ehat1P1(-1/2) in tensor regression without sample splitting}, by symmetry, we also have
\begin{align}
\mathrm{\RN{5}} 
\lesssim & \frac{\sigma_\xi^2}{\ulambda^2\sigma^2} \cdot \frac{\op}{n}. \label{eq: upper bound of term 5 in P1(-1/2)Ehat1P1(-1/2) in tensor regression without sample splitting}
\end{align}

For the sixth term \eqref{eq: term 6 in P1(-1/2)Ehat1P1(-1/2) in tensor regression without sample splitting}, we have
\begin{align}
\mathrm{\RN{6}} 
\leq & \frac{1}{\ulambda} \cdot \left\|U_2^{\top}\left(\mcP_{\whU_2^{(1)}}-\mcP_{U_2}\right)U_2\right\| \cdot \left\|U_3^{\top}\left(\mcP_{\whU_3^{(1)}}-\mcP_{U_3}\right)U_3\right\| \lesssim \frac{\sigma_\xi^4}{\ulambda^4\sigma^4}\cdot \frac{\op^2}{n^2}. \label{eq: upper bound of term 6 in P1(-1/2)Ehat1P1(-1/2) in tensor regression without sample splitting}
\end{align}
Therefore, we have 
\begin{align*}
& \left\|U_1\left(G_1G_1\right)^{-\frac{1}{2}}U_1^{\top}\whE_1U_1\left(G_1G_1\right)^{-\frac{1}{2}}U_1^{\top} \right\| \notag \\
\leq & \left[\underbrace{\frac{\sigma_{\xi}}{\ulambda \sigma} \cdot\left(\sqrt{\frac{\oor \log (\op)}{n}}+\Delta \sqrt{\frac{\op}{n}}\right)}_{\eqref{eq: upper bound of term 1 in P1(-1/2)Ehat1P1(-1/2) in tensor regression without sample splitting},\eqref{eq: upper bound of term 2 in P1(-1/2)Ehat1P1(-1/2) in tensor regression without sample splitting}} \right] + \left[\underbrace{\frac{\sigma_\xi^2}{\ulambda^2\sigma^2}\cdot \left(\frac{\oor\log(\op)}{n}+ \Delta^2\cdot \frac{\op}{n}\right)}_{\eqref{eq: upper bound of term 3 in P1(-1/2)Ehat1P1(-1/2) in tensor regression without sample splitting}}\right] + \left(\underbrace{\frac{\sigma_\xi^2}{\ulambda^2\sigma^2}\cdot \frac{\op}{n}}_{\eqref{eq: upper bound of term 4 in P1(-1/2)Ehat1P1(-1/2) in tensor regression without sample splitting},\eqref{eq: upper bound of term 5 in P1(-1/2)Ehat1P1(-1/2) in tensor regression without sample splitting}}\right) + \left(\underbrace{\frac{\sigma_\xi^4}{\ulambda^4\sigma^4}\cdot \frac{\op^2}{n^2}}_{\eqref{eq: upper bound of term 6 in P1(-1/2)Ehat1P1(-1/2) in tensor regression without sample splitting}}\right) \\
\lesssim & \frac{\sigma_{\xi}}{\ulambda\sigma}\cdot \left( \sqrt{\frac{\oor\log(\op)}{n}} +\Delta\cdot \sqrt{\frac{\op}{n}}\right).
\end{align*}
which leads to \eqref{eq: high-prob upper bound of P1(-1/2)Ehat1P1(-1/2) in tensor regression without sample splitting}.

Applying similar arguments, we obtain the upper bounds for the rest two terms.

\end{proof}

\begin{proposition}\label{prop: high-prob upper bound of first-order perturbation terms after projection in tensor regression without sample splitting}
Under the same setting of Theorem~\ref{thm: main theorem in tensor regression without sample splitting}, let $\whE_j$ be defined by \eqref{eq: definition of whEj in tensor regression without sample splitting} and $\mcP_j^{-s}=U_j\left(G_jG_j^{\top}\right)^{-2s}U_j^{\top}$ for any $s>0$ and $j=1,2,3$ while $\mcP_j^{0} = U_{j\perp}U_{j\perp}^{\top}$. Then with probability at least $1-\exp(-cn) - \frac{1}{p^C} - \P\left(\mcE_\Delta\right) - \P\left(\mcE_{U}^{\text{reg}}\right)$, it follows that:
\begin{align}
\left\|V_j^{\top} \mcP_{U_{j \perp}} \whE_j \mcP_j^{-\frac{1}{2}}\right\| = \left\|V_j^{\top}\mcP_{U_{j\perp}}\whE_j U_j\left(G_jG_j^{\top}\right)^{-\frac{1}{2}}U_j^{\top}\right\| \lesssim \frac{\sigma_{\xi}}{\sigma} \cdot\left(\sqrt{\frac{\oR \log (\op)}{n}}+\Delta \sqrt{\frac{\op}{n}}\right), \label{eq: high-prob upper bound of V1tP1(0)Ehat1P1(-1/2) in tensor regression without sample splitting}
\end{align}
and 
\begin{align}
\left\|V_j^{\top} \mcP_{U_{j\perp}} \whE_j U_{j\perp}\right\| \lesssim \frac{\sigma_{\xi}^2}{\sigma^2} \cdot\left(\frac{\sqrt{\oR \log (\op)} \cdot \sqrt{\op}}{n}+\Delta \cdot \frac{\op}{n}\right),\label{eq: high-prob upper bound of V1tP1(0)Ehat1P1(0) in tensor regression without sample splitting}
\end{align}
where $c$ and $C$ are two universal constants . 

\end{proposition}

\begin{proof}

By symmetry, it suffice to consider $\left\|V_1^{\top}\mcP_{U_{1\perp}}\whE_1 U_1\right\|$ and $\left\|V_1^{\top}\mcP_{U_{1\perp}}\whE_1 U_{1\perp}\right\|$.
%
%
Note that
\begin{align}
& \left\|V_1^{\top}\mcP_{U_{1\perp}}\whE_1 U_1\left(G_1G_1\right)^{-\frac{1}{2}}U_1^{\top}\right\| \notag \\
\leq & \underbrace{\left\|V_1^{\top}\mcP_{U_{1\perp}}\whZ_1\left(U_3\otimes U_2\right)G_1^{\top}\left(G_1G_1\right)^{-\frac{1}{2}}U_1^{\top}\right\|}_{\mathrm{\RN{1}}} \label{eq: term 1 in V1tU1ptEhat1U1 in tensor regression without sample splitting}\\
+ & \underbrace{\left\|V_1^{\top}\mcP_{U_{1\perp}}\whZ_1\left(\mcP_{U_3}\otimes \mcP_{U_2}\right)\whZ_1^{\top}U_1\left(G_1G_1\right)^{-\frac{1}{2}}U_1^{\top}\left(G_1G_1\right)^{-\frac{1}{2}}U_1^{\top}\right\|}_{\mathrm{\RN{2}}} \label{eq: term 2 in V1tU1ptEhat1U1 in tensor regression without sample splitting}\\
+ & \underbrace{\left\|V_1^{\top}\mcP_{U_{1\perp}}\whZ_1\left[\left(\mcP_{\whU_3^{(1)}}- \mcP_{U_3}\right)U_3 \otimes U_2\right]G_1^{\top}\left(G_1G_1\right)^{-\frac{1}{2}}U_1^{\top}\right\|}_{\mathrm{\RN{3}}} \label{eq: term 3 in V1tU1ptEhat1U1 in tensor regression without sample splitting}\\
+ & \underbrace{\left\|V_1^{\top}\mcP_{U_{1\perp}}\whZ_1\left[U_3 \otimes \left(\mcP_{\whU_2^{(1)}}- \mcP_{U_2}\right)U_2\right]G_1^{\top}\left(G_1G_1\right)^{-\frac{1}{2}}U_1^{\top}\right\|}_{\mathrm{\RN{4}}} \label{eq: term 4 in V1tU1ptEhat1U1 in tensor regression without sample splitting}\\
+ & \underbrace{\left\|V_1^{\top}\mcP_{U_{1\perp}}\whZ_1\left[\left(\mcP_{\whU_3^{(1)}}- \mcP_{U_3}\right)U_3 \otimes \left(\mcP_{\whU_2^{(1)}}- \mcP_{U_2}\right)U_2\right]G_1^{\top}\left(G_1G_1\right)^{-\frac{1}{2}}U_1^{\top}\right\|}_{\mathrm{\RN{5}}} \label{eq: term 5 in V1tU1ptEhat1U1 in tensor regression without sample splitting}\\
+ & \underbrace{\left\|V_1^{\top}\mcP_{U_{1\perp}}\whZ_1\left[\left(\mcP_{\whU_3^{(1)}}- \mcP_{U_3}\right)\mcP_{U_3} \otimes \mcP_{U_2}\right]\whZ_1^{\top}U_1\left(G_1G_1\right)^{-\frac{1}{2}}U_1^{\top}\left(G_1G_1\right)^{-\frac{1}{2}}U_1^{\top}\right\|}_{\mathrm{\RN{6}}} \label{eq: term 6 in V1tU1ptEhat1U1 in tensor regression without sample splitting}\\
+ & \underbrace{\left\|V_1^{\top}\mcP_{U_{1\perp}}\whZ_1\left[\mcP_{U_3} \otimes \left(\mcP_{\whU_2^{(1)}}- \mcP_{U_2}\right)\right]\whZ_1^{\top}U_1\left(G_1G_1\right)^{-\frac{1}{2}}U_1^{\top}\right\|}_{\mathrm{\RN{7}}} \label{eq: term 7 in V1tU1ptEhat1U1 in tensor regression without sample splitting}\\
+ & \underbrace{\left\|V_1^{\top}\mcP_{U_{1\perp}}\whZ_1\left[\left(\mcP_{\whU_3^{(1)}}- \mcP_{U_3}\right) \otimes \left(\mcP_{\whU_2^{(1)}}- \mcP_{U_2}\right)\right]\whZ_1^{\top}U_1\left(G_1G_1\right)^{-\frac{1}{2}}U_1^{\top}\right\|}_{\mathrm{\RN{8}}} . \label{eq: term 8 in V1tU1ptEhat1U1 in tensor regression without sample splitting} 
\end{align}
First, for the first term $\mathrm{\RN{1}}$ \eqref{eq: term 1 in V1tU1ptEhat1U1 in tensor regression without sample splitting}, we have 
\begin{align}
\mathrm{\RN{1}} 
\leq & \left\|V_1^{\top}\mcP_{U_{1\perp}}\whZ_1\left(U_3\otimes U_2\right)\right\| \cdot \left\|G_1^{\top}\left(G_1G_1\right)^{-\frac{1}{2}}U_1^{\top} \right\| \lesssim \frac{\sigma_\xi}{\sigma}\cdot \left(\sqrt{\frac{\oR\log(\op)}{n}}+ \Delta \cdot \sqrt{\frac{\op}{n}}\right) \label{eq: upper bound of term 1 in V1tU1ptEhat1U1 in tensor regression without sample splitting}
\end{align}
and for the second term \eqref{eq: term 2 in V1tU1ptEhat1U1 in tensor regression without sample splitting}, we have
\begin{align}
\mathrm{\RN{2}} 
\leq & \left\|V_1^{\top}\mcP_{U_{1\perp}}\whZ_1\left(\mcP_{U_3}\otimes \mcP_{U_2}\right)\right\| \cdot \left\|\left(\mcP_{U_3}\otimes \mcP_{U_2}\right)\whZ_1^{\top}U_1\left(G_1G_1\right)^{-\frac{1}{2}}U_1^{\top} \right\| 
\lesssim  \frac{\sigma_\xi^2}{\ulambda\sigma^2}\cdot \left(\frac{\sqrt{\oR\oor}\log(\op)}{n}+\Delta^2\cdot \frac{\op}{n}\right). \label{eq: upper bound of term 2 in V1tU1ptEhat1U1 in tensor regression without sample splitting}
\end{align}
Then, consider the third term \eqref{eq: term 3 in V1tU1ptEhat1U1 in tensor regression without sample splitting}
\begin{align}
\mathrm{\RN{3}} 
\leq & \left\|V_1^{\top}\mcP_{U_{1 \perp}} \whZ_1 \left[S_{G_1,1}\left(\whE_1^{(0)}\right) U_3 \otimes U_2\right] G_1^{\top}\left(G_1G_1\right)^{-\frac{1}{2}}U_1^{\top}\right\| \notag \\
+ & \left\|V_1^{\top}\mcP_{U_{1 \perp}} \whZ_1\left[\sum_{k_3=2}^{+\infty}\left(S_{G_3,k_3}\left(\whE_3^{(0)}\right)  U_3 \otimes U_2\right) G_1^{\top}\left(G_1G_1\right)^{-\frac{1}{2}}U_1^{\top} \right]\right\| \notag \\
\leq & \underbrace{\left\|V_1^{\top}\mcP_{U_{1 \perp}} \whZ_1 \left[\left(\mcP_{U_{3\perp}}\whZ_3\left(\mcP_{\whU_3^{(0)}} \otimes \mcP_{\whU_1^{(0)}} \right)\left(U_1\otimes U_3\right)G_2^{\top}\left(G_2G_2^{\top}\right)^{-1}\right) \otimes U_3\right]G_1^{\top}\left(G_1G_1\right)^{-\frac{1}{2}}U_1^{\top}\right\|}_{\mathrm{\RN{3}}.\mathrm{\RN{1}}} \label{eq: term 3.1 in V1tU1ptEhat1U1 in tensor regression without sample splitting} \\
+ & \underbrace{\left\|V_1^{\top}\mcP_{U_{1 \perp}} \whZ_1 \left[\left(\mcP_{U_{3\perp}}\whZ_3\left(\mcP_{\whU_2^{(0)}} \otimes \mcP_{\whU_1^{(0)}} \right)\whZ_3^{\top}U_3\left(G_3G_3^{\top}\right)^{-1}\right) \otimes U_2\right]G_1^{\top}\left(G_1G_1\right)^{-\frac{1}{2}}U_1^{\top}\right\|}_{\mathrm{\RN{3}}.\mathrm{\RN{2}}} \label{eq: term 3.2 in V1tU1ptEhat1U1 in tensor regression without sample splitting} \\
+ & \underbrace{\left\|V_1^{\top}\mcP_{U_{1 \perp}} \whZ_1\left[\sum_{k_3=2}^{+\infty}\left(S_{G_3,k_3}\left(\whE_3^{(0)}\right)  U_3 \otimes U_2\right)\right] G_1^{\top}\left(G_1G_1\right)^{-\frac{1}{2}}U_1^{\top}\right\|}_{\mathrm{\RN{3}}.\mathrm{\RN{3}}}. \label{eq: term 3.3 in V1tU1ptEhat1U1 in tensor regression without sample splitting}
\end{align}
Here, first, we have 
\begin{align}
\mathrm{\RN{3}}.\mathrm{\RN{1}}
\leq & \underbrace{\left\|V_1^{\top}\mcP_{U_{1 \perp}} \whZ_1 \left[\left(\mcP_{U_{3\perp}}\whZ_3\left(U_2\otimes U_1\right)G_3^{\top}\left(G_3G_3^{\top}\right)^{-1}U_3^{\top}\right) \otimes U_2 \right]G_1^{\top}\left(G_1G_1\right)^{-\frac{1}{2}}U_1^{\top}\right\|}_{\mathrm{\RN{3}}.\mathrm{\RN{1}}.\mathrm{\RN{1}}} \notag \\
+ & \underbrace{\left\|V_1^{\top}\mcP_{U_{1 \perp}} \whZ_1 \left[\left(\mcP_{U_{3\perp}}\whZ_3\left(\left(\mcP_{\whU_2^{(0)}} - \mcP_{U_2}\right)U_2\otimes U_1\right)G_3^{\top}\left(G_3G_3^{\top}\right)^{-1}\right)\otimes U_2\right]G_1^{\top}\left(G_1G_1\right)^{-\frac{1}{2}}U_1^{\top}\right\|}_{\mathrm{\RN{3}}.\mathrm{\RN{1}}.\mathrm{\RN{2}}} \notag \\
+ &  \underbrace{\left\|V_1^{\top}\mcP_{U_{1 \perp}} \whZ_1 \left[\left(\mcP_{U_{3\perp}}\whZ_3\left(U_2\otimes \left(\mcP_{\whU_1^{(0)}} -\mcP_{U_1}\right)U_1\right)G_3^{\top}\left(G_3G_3^{\top}\right)^{-1}\right)\otimes U_2\right]G_1^{\top}\left(G_1G_1\right)^{-\frac{1}{2}}U_1^{\top}\right\|}_{\mathrm{\RN{3}}.\mathrm{\RN{1}}.\mathrm{\RN{3}}} \notag \\
+ & \underbrace{\left\|V_1^{\top}\mcP_{U_{1 \perp}} \whZ_1 \left[\left(\mcP_{U_{3\perp}}\whZ_3\left(\left(\mcP_{\whU_2^{(0)}}- \mcP_{U_2}\right)U_2\otimes \left(\mcP_{\whU_1^{(0)}} -\mcP_{U_1}\right)U_1\right)G_3^{\top}\left(G_3G_3^{\top}\right)^{-1}\right)\otimes U_2\right]G_1^{\top}\left(G_1G_1\right)^{-\frac{1}{2}}U_1^{\top}\right\|}_{\mathrm{\RN{3}}.\mathrm{\RN{1}}.\mathrm{\RN{4}}}. \notag 
\end{align}

Note that
\begin{align*}
& \mathrm{\RN{3}}.\mathrm{\RN{1}}.\mathrm{\RN{2}} + \mathrm{\RN{3}}.\mathrm{\RN{1}}.\mathrm{\RN{3}} + \mathrm{\RN{3}}.\mathrm{\RN{1}}.\mathrm{\RN{4}}\\
\leq & \frac{1}{\ulambda} \cdot \sup_{\substack{W_2\in \mathbb{R}^{p_2\times r_2}, \left\|W_2\right\|=1\\ W_3\in \mathbb{R}^{p_3\times r_3}, \left\|W_3\right\|=1}} \left\|V_1^{\top}\mcP_{U_{2\perp}}\whZ_1\left(W_3\otimes W_2\right)\right\| \cdot \sup_{\substack{W_3\in \mathbb{R}^{p_3\times r_3}, \left\|W_3\right\|=1\\ W_1\in \mathbb{R}^{p_1\times r_1}, \left\|W_1\right\|=1}}\left\|\mcP_{U_{2\perp}}\whZ_2\left(W_1\otimes W_3\right)\right\| \\
& \cdot \left(\left\|\mcP_{\whU_2^{(0)}}- \mcP_{U_2}\right\| + \left\|\mcP_{\whU_3^{(0)}}- \mcP_{U_3}\right\| + \left\|\mcP_{\whU_2^{(0)}}- \mcP_{U_2}\right\| \cdot \left\|\mcP_{\whU_3^{(0)}}- \mcP_{U_3}\right\|\right) \\
\lesssim & \frac{\sigma_\xi^3}{\ulambda^2\sigma^3}\cdot \frac{\op^{3/2}}{n^{3/2}}.
\end{align*}

Therefore, we have
\begin{align*}
\mathrm{\RN{3}}.\mathrm{\RN{1}} 
\leq & \underbrace{\left\|V_1^{\top}\mcP_{U_{1 \perp}} \whZ_1^{(1)} \left[\left(\mcP_{U_{3\perp}}\whZ_3^{(1)}\left(U_2\otimes U_1\right)G_3^{\top}\left(G_3G_3^{\top}\right)^{-1}U_3^{\top}\right) \otimes U_2\right]G_1^{\top}\left(G_1G_1\right)^{-\frac{1}{2}}U_1^{\top}\right\|}_{\mathrm{\RN{3}}.\mathrm{\RN{1}}.\mathrm{\RN{1}}.\mathrm{\RN{1}}} \\
+ & \underbrace{\left\|V_1^{\top}\mcP_{U_{1 \perp}} \whZ_1^{(1)} \left[\left(\mcP_{U_{3\perp}}\whZ_3^{(2)}\left(U_2\otimes U_1\right)G_3^{\top}\left(G_3G_3^{\top}\right)^{-1}U_3^{\top}\right) \otimes U_2\right]G_1^{\top}\left(G_1G_1\right)^{-\frac{1}{2}}U_1^{\top}\right\|}_{\mathrm{\RN{3}}.\mathrm{\RN{1}}.\mathrm{\RN{1}}.\mathrm{\RN{2}}} \\
+ & \underbrace{\left\|V_1^{\top}\mcP_{U_{1 \perp}} \whZ_1^{(2)} \left[\left(\mcP_{U_{3\perp}}\whZ_3^{(1)}\left(U_2\otimes U_1\right)G_3^{\top}\left(G_3G_3^{\top}\right)^{-1}U_3^{\top}\right) \otimes U_2\right]G_1^{\top}\left(G_1G_1\right)^{-\frac{1}{2}}U_1^{\top}\right\|}_{\mathrm{\RN{3}}.\mathrm{\RN{1}}.\mathrm{\RN{1}}.\mathrm{\RN{3}}} \\
+ & \underbrace{\left\|V_1^{\top}\mcP_{U_{1 \perp}} \whZ_1^{(2)} \left[\left(\mcP_{U_{3\perp}}\whZ_3^{(2)}\left(U_2\otimes U_1\right)G_3^{\top}\left(G_3G_3^{\top}\right)^{-1}U_3^{\top}\right) \otimes U_2\right]G_1^{\top}\left(G_1G_1\right)^{-\frac{1}{2}}U_1^{\top}\right\|}_{\mathrm{\RN{3}}.\mathrm{\RN{1}}.\mathrm{\RN{1}}.\mathrm{\RN{4}}} \\
+ & \frac{\sigma_\xi^3}{\ulambda^2\sigma^3}\cdot \frac{\op^{3/2}}{n^{3/2}}.
\end{align*}
where $\whZ_1^{(1)}=\frac{1}{n\sigma^2}\sum_{i=1}^n\xi_i\Mat_1\left(\mcX_i\right)$ and $\whZ_1^{(2)}=\frac{1}{n\sigma^2}\sum_{i=1}^n\left[\left\langle\mcX_i, \widehat{\Delta}\right\rangle\Mat_1\left(\mcX_i\right) -\sigma^2\cdot \widehat{\Delta}\right]$.

It then remains to find upper bound for $\mathrm{\RN{3}}.\mathrm{\RN{1}}.\mathrm{\RN{1}}.\mathrm{\RN{1}}, \mathrm{\RN{3}}.\mathrm{\RN{1}}.\mathrm{\RN{1}}.\mathrm{\RN{2}}, \mathrm{\RN{3}}.\mathrm{\RN{1}}.\mathrm{\RN{1}}.\mathrm{\RN{3}}$ and $\mathrm{\RN{3}}.\mathrm{\RN{1}}.\mathrm{\RN{1}}.\mathrm{\RN{4}}$. By Lemma~\ref{lemma: high-prob upper bound of spectral norm of W_1tZ(1)t(BZ(1)W2oW3)}, we have
\begin{align}
\mathrm{\RN{3}}.\mathrm{\RN{1}}.\mathrm{\RN{1}} 
\leq & \left\|V_1^{\top}\mcP_{U_{1 \perp}} \whZ_1^{(1)} \left[\left(\mcP_{U_{3\perp}}\whZ_3^{(1)}\left(U_2\otimes U_1\right)G_3^{\top}\left(G_3G_3^{\top}\right)^{-1}U_3^{\top}\right) \otimes U_2\right]\right\| \lesssim \frac{\sigma_\xi^2}{\sigma^2}\cdot \frac{\sqrt{\oR\op\log(\op)}}{n}. \label{eq: upper bound of term 3.1.1 in V1tU1ptEhat1U1 in tensor regression without sample splitting}
\end{align}
where the last inequality follows from that 
$
\operatorname{tr}\left[U_3\left(G_3G_3^{\top}\right)^{-1}G_3\left(U_2\otimes U_1\right)^{\top}\left[U_2\otimes \left(\mcP_{U_{1\perp}}V_1\right)\right]\right]=0.
$

In addition, we have
\begin{align}
\mathrm{\RN{3}}.\mathrm{\RN{1}}.\mathrm{\RN{1}}.\mathrm{\RN{2}} 
\leq & \frac{1}{\ulambda} \cdot \sup_{\substack{W_2\in \mathbb{R}^{p_2\times r_2}, \left\|W_2\right\|=1\\W_3\in \mathbb{R}^{p_3\times r_3}, \left\|W_3\right\|=1}}\left\|V_1^{\top}\mcP_{U_{1 \perp}} \whZ_1^{(1)}\left(W_3\otimes W_2\right)\right\|\cdot \left\|\mcP_{U_{3\perp}}\whZ_3^{(2)}\left(U_2\otimes U_1\right)\right\| \lesssim \frac{\sigma_\xi^2}{\ulambda\sigma^2}\cdot \Delta\cdot \frac{\op}{n}, \label{eq: upper bound of term 3.1.2 in V1tU1ptEhat1U1 in tensor regression without sample splitting} \\
\mathrm{\RN{3}}.\mathrm{\RN{1}}.\mathrm{\RN{1}}.\mathrm{\RN{3}} 
\leq & \frac{1}{\ulambda} \cdot \sup_{\substack{W_2\in \mathbb{R}^{p_2\times r_2}, \left\|W_2\right\|=1\\W_3\in \mathbb{R}^{p_3\times r_3}, \left\|W_3\right\|=1}}\left\|V_1^{\top}\mcP_{U_{1 \perp}} \whZ_1^{(2)}\left(W_3\otimes W_2\right)\right\|\cdot \left\|\mcP_{U_{3\perp}}\whZ_3^{(1)}\left(U_2\otimes U_1\right)\right\| \lesssim \frac{\sigma_\xi^2}{\ulambda\sigma^2}\cdot \Delta\cdot \frac{\op}{n}, \label{eq: upper bound of term 3.1.3 in V1tU1ptEhat1U1 in tensor regression without sample splitting}
\end{align}
and
\begin{align}
\mathrm{\RN{3}}.\mathrm{\RN{1}}.\mathrm{\RN{1}}.\mathrm{\RN{4}} 
\leq & \frac{1}{\ulambda} \cdot \sup_{\substack{W_2\in \mathbb{R}^{p_2\times r_2}, \left\|W_2\right\|=1\\W_3\in \mathbb{R}^{p_3\times r_3}, \left\|W_3\right\|=1}}\left\|V_1^{\top}\mcP_{U_{1 \perp}} \whZ_1^{(2)}\left(W_3\otimes W_2\right)\right\|\cdot \left\|\mcP_{U_{3\perp}}\whZ_3^{(2)}\left(U_2\otimes U_1\right)\right\| \lesssim \frac{\sigma_\xi^2}{\ulambda\sigma^2}\cdot \Delta^2\cdot \frac{\op}{n}. \label{eq: upper bound of term 3.1.4 in V1tU1ptEhat1U1 in tensor regression without sample splitting}
\end{align}
Therefore, we have
\begin{align}
\mathrm{\RN{3}}.\mathrm{\RN{1}} 
\lesssim & \left(\underbrace{\frac{\sigma_\xi^2}{\ulambda\sigma^2}\cdot \frac{\sqrt{\oR\op\log(\op)}}{n}}_{\eqref{eq: upper bound of term 3.1.1 in V1tU1ptEhat1U1 in tensor regression without sample splitting}}\right)
+ \left(\underbrace{\frac{\sigma_\xi}{\ulambda\sigma}\cdot \Delta\cdot \frac{\op}{n}}_{\eqref{eq: upper bound of term 3.1.2 in V1tU1ptEhat1U1 in tensor regression without sample splitting},\eqref{eq: upper bound of term 3.1.3 in V1tU1ptEhat1U1 in tensor regression without sample splitting}}\right)
+ \left(\underbrace{\frac{\sigma_\xi^2}{\ulambda\sigma^2}\cdot \Delta^2\cdot \frac{\op}{n}}_{\eqref{eq: upper bound of term 3.1.4 in V1tU1ptEhat1U1 in tensor regression without sample splitting}}\right)
+ \frac{\sigma_\xi^3}{\ulambda^2\sigma^3}\cdot \frac{\op^{3/2}}{n^{3/2}} \notag \\
\lesssim & \frac{\sigma_\xi^2}{\ulambda\sigma^2}\cdot \left(\frac{\sqrt{\oR\op\log(\op)}}{n}+\Delta\cdot \frac{\op}{n}\right). \label{eq: upper bound of term 3.1 in V1tU1ptEhat1U1 in tensor regression without sample splitting}
\end{align}
Second, for $\mathrm{\RN{3}}.\mathrm{\RN{2}} =V_1^{\top} \mcP_{U_{1 \perp}} \whZ_1\left[\left(\mcP_{U_{3 \perp}} \whZ_3\left(\mcP_{\whU_2^{(0)}} \otimes \mcP_{\whU_1^{(0)}}\right) \whZ_3^{\top}U_3\left(G_3G_3^{\top}\right)^{-1}\right) \otimes U_2\right] G_1^{\top}\left(G_1G_1\right)^{-\frac{1}{2}}U_1^{\top}$, we have
\begin{align}
\mathrm{\RN{3}}.\mathrm{\RN{2}} 
\leq & \frac{1}{\ulambda^2} \sup_{\substack{W_2\in \mathbb{R}^{p_2\times r_2}, \left\|W_2\right\|=1\\W_3\in \mathbb{R}^{p_3\times r_3}, \left\|W_3\right\|=1}}\left\|V_1^{\top} \mcP_{U_{1 \perp}} \whZ_1\left(W_3\otimes W_2\right)\right\|  \sup_{\substack{W_1\in \mathbb{R}^{p_1\times r_1}, \left\|W_1\right\|=1\\W_2\in \mathbb{R}^{p_2\times r_2}, \left\|W_2\right\|=1}}\left\|\mcP_{U_{3\perp}} \whZ_3\left(W_2\otimes W_3\right)\right\| \notag \\
& \cdot \sup_{\substack{W_1\in \mathbb{R}^{p_1\times r_1}, \left\|W_1\right\|=1\\W_2\in \mathbb{R}^{p_2\times r_3}, \left\|W_2\right\|=1}}\left\|U_3^{\top} \whZ_3\left(W_2\otimes W_1\right)\right\| \notag \\
\lesssim & \frac{\sigma_\xi^3}{\ulambda^2\sigma^3}\cdot \frac{\op^{3/2}}{n^{3/2}}. \label{eq: upper bound of term 3.2 in V1tU1ptEhat1U1 in tensor regression without sample splitting}
\end{align}
In addition, it follows immediately that
\begin{align}
\mathrm{\RN{3}}.\mathrm{\RN{3}} 
\leq & \sup_{\substack{W_2\in \mathbb{R}^{p_2\times r_2}, \left\|W_2\right\|=1\\W_3\in \mathbb{R}^{p_3\times r_3}, \left\|W_3\right\|=1}}\left\|V_1^{\top} \mcP_{U_{1 \perp}} \whZ_1\left(W_3\otimes W_2\right)\right\| \cdot \sum_{k_1=2}^{+\infty}\left\|S_{G_1,k_1}\left(\whE_1^{(0)}\right)\right\| \lesssim \frac{\sigma_\xi^3}{\ulambda^2\sigma^3}\cdot \frac{\op^{3/2}}{n^{3/2}}. \label{eq: upper bound of term 3.3 in V1tU1ptEhat1U1 in tensor regression without sample splitting}
\end{align}
Combining the results  above, we have
\begin{align}
\mathrm{\RN{3}}
\lesssim & \eqref{eq: upper bound of term 3.1 in V1tU1ptEhat1U1 in tensor regression without sample splitting} + \eqref{eq: upper bound of term 3.2 in V1tU1ptEhat1U1 in tensor regression without sample splitting} + \eqref{eq: upper bound of term 3.3 in V1tU1ptEhat1U1 in tensor regression without sample splitting} 
\lesssim \frac{\sigma_\xi^2}{\ulambda\sigma^2}\cdot \left(\frac{\sqrt{\oR\op\log(\op)}}{n}+\Delta\cdot \frac{\op}{n}\right). \label{eq: upper bound of term 3 in V1tU1ptEhat1U1 in tensor regression without sample splitting}.
\end{align}

Therefore, we have
\begin{align*}
\left\|V_1^{\top}\mcP_{U_{1\perp}}\whE_1 \mcP_1^{-\frac{1}{2}}\right\| 
\lesssim & \eqref{eq: upper bound of term 1 in V1tU1ptEhat1U1 in tensor regression without sample splitting} + \eqref{eq: upper bound of term 2 in V1tU1ptEhat1U1 in tensor regression without sample splitting} + \eqref{eq: upper bound of term 3 in V1tU1ptEhat1U1 in tensor regression without sample splitting} \lesssim \frac{\sigma_\xi}{\sigma} \cdot \left(\sqrt{\frac{\oR\log (\op)}{n}}+\Delta \sqrt{\frac{\op}{n}}\right),
\end{align*}
where the second inequality holds as long as $\ulambda \geq \frac{\kappa\sigma_\xi}{\sigma}\cdot \sqrt{\frac{\op}{n}}$.

Applying similar arguments, we obtain the bounds for the second term.

\end{proof}

\subsection{Upper Bound of Higher-Order Perturbation Terms}

\begin{lemma} \label{lemma: high-prob upper bound of the higher-order perturbation terms in tensor regression without sample splitting}
Under the same setting of Theorem~\ref{thm: main theorem in tensor regression without sample splitting}, let $\whE_j$ be defined by \eqref{eq: definition of whEj in tensor regression without sample splitting}, $S_{G_j, k_j}(\whE_j)$ be defined by \eqref{eq: definition of S1(whEj) in tensor regression without sample splitting}, $\whE_j$ be defined by \eqref{eq: definition of whEj in tensor regression without sample splitting} and $\mcP_j^{-s}=U_j\left(G_jG_j^{\top}\right)^{-2s}U_j^{\top}$ for any $s>0$ and $j=1,2,3$ while $\mcP_j^{0} = U_{j\perp}U_{j\perp}^{\top}$. Then with probability at least $1-\exp(-cn) - \frac{1}{p^C} - \P\left(\mcE_\Delta\right) - \P\left(\mcE_{U}^{\text{reg}}\right)$, it follows that:
\begin{align}
& \left\|\mcP_{U_j} \sum_{k_j=2}^{+\infty} S_{G_j, k_j}\left(\whE_j\right) \mcP_{U_{j \perp}} V_j\right\| \lesssim \frac{\sigma_\xi^2}{\ulambda^2\sigma^2}\left(\frac{\sqrt{\oR\oor} \log(\op)}{n} + \Delta\cdot \frac{\sqrt{\oR\op\log(\op)}}{n} + \Delta^2\cdot \frac{\op}{n}\right) \label{eq: high-prob upper bound of PUorder2PUpPV in tensor regression without sample splitting}\\
& \left\|\mcP_{U_j} \sum_{k_j=3}^{+\infty} S_{G_j, k_j}\left(\whE_j\right) \mcP_{U_{j\perp}}V_j\right\| \lesssim \frac{\sigma_{\xi}^3}{\ulambda^3 \sigma^3} \cdot\left(\frac{\op \sqrt{\oR \log (\op)}}{n^{3 / 2}}+\Delta \cdot \frac{\op^{3 / 2}}{n^{3 / 2}}\right) \label{eq: high-prob upper bound of PUorder3PUpPV in tensor regression without sample splitting} \\
& \left\|\mcP_{U_{j\perp}} \sum_{k_j=2}^{+\infty} S_{G_j, k_j}\left(\whE_j\right) \mcP_{U_{j \perp}} V_j\right\| \lesssim \frac{\sigma_{\xi}^2}{\ulambda^2 \sigma^2} \cdot\left(\frac{\sqrt{\oR \op \log (\op)}}{n}+\Delta \cdot \frac{\op}{n}\right) \label{eq: high-prob upper bound of PUporder2PUpPV in tensor regression without sample splitting}\\
& \left\|\mcP_{U_{j \perp}} \sum_{k_j=3}^{+\infty} S_{G_j, k_j}\left(\whE_j\right) \mcP_{U_{j \perp}} V_j\right\| \lesssim \frac{\sigma_\xi^4}{\ulambda^3\sigma^4}\left(\frac{\op^{3/2}  \oR^{1/2} \log (\op)}{n^2} + \Delta\cdot \frac{\op^2}{n^2}\right) \label{eq: high-prob upper bound of PUporder3PUpPV in tensor regression without sample splitting}.
\end{align}

\end{lemma}

\begin{proof}

By symmetry, it suffices to consider
\begin{align*}
& \mathrm{\RN{1}} = \left\|\mcP_{U_1} \sum_{k_1=2}^{+\infty} S_{G_1, k_1}\left(\whE_1\right) \mcP_{U_{1 \perp}}V_1\right\|, \mathrm{\RN{2}} = \left\|\mcP_{U_1} \sum_{k_1=3}^{+\infty} S_{G_1, k_1}\left(\whE_1\right) \mcP_{U_{1 \perp}}V_1\right\|\\
& \mathrm{\RN{3}} = \left\|\mcP_{U_{1\perp}} \sum_{k_1=2}^{+\infty} S_{G_1, k_1}\left(\whE_1\right) \mcP_{U_{1 \perp}}V_1\right\|, \mathrm{\RN{4}} = \left\|\mcP_{U_{1\perp}} \sum_{k_1=3}^{+\infty} S_{G_1, k_1}\left(\whE_1\right) \mcP_{U_{1 \perp}}V_1\right\|.
\end{align*}


First, we have
\begin{align*}
\mathrm{\RN{1}} 
\lesssim & \left\|\mcP_{U_1}\mcP_1^{-2}\whE_1\mcP_1^{0}\whE_1\mcP_1^{0}V_1\right\| +  \left\|\mcP_{U_1}\mcP_1^{-1}\whE_1\mcP_1^{-1}\whE_1\mcP_1^{0}V_1\right\| + \left\|\mcP_{U_1}\sum_{k_1=3}^{+\infty} S_{G_1, k_1}\left(\whE_1\right)\mcP_{U_{1 \perp}}V_1\right\| \\
\leq & \underbrace{\left\|\mcP_1^{-2}\whE_1\mcP_1^{0}\right\|}_{\eqref{eq: high-prob upper bound of P1(0)Ehat1P1(-1/2) in tensor regression}} \cdot \underbrace{\left\|\mcP_1^{0}\whE_1\mcP_1^{0}V_1\right\|}_{\eqref{eq: high-prob upper bound of V1tP1(0)Ehat1P1(0) in tensor regression without sample splitting}} + \underbrace{\left\|\mcP_{U_1}\mcP_1^{-1}\whE_1\mcP_1^{-\frac{1}{2}}\right\|}_{\eqref{eq: high-prob upper bound of P1(-1/2)Ehat1P1(-1/2) in tensor regression without sample splitting}} \cdot \underbrace{\left\|\mcP_1^{-\frac{1}{2}}\whE_1\mcP_1^{0}V_1\right\|}_{\eqref{eq: high-prob upper bound of V1tP1(0)Ehat1P1(-1/2) in tensor regression without sample splitting}} + \underbrace{\left\|\mcP_{U_1}\sum_{k_1=3}^{+\infty} S_{G_1, k_1}\left(\whE_1\right)\mcP_{U_{1 \perp}}V_1\right\|}_{\eqref{eq: high-prob upper bound of PUorder3PUpPV in tensor regression without sample splitting}}  \\
\lesssim & \frac{\sigma_\xi^2}{\ulambda^2\sigma^2}\left(\frac{\sqrt{\oR\oor} \log(\op)}{n} + \Delta\cdot \frac{\sqrt{\oR\op\log(\op)}}{n} + \Delta^2\cdot \frac{\op}{n}\right).
\end{align*}

Note that
\begin{align*}
\mathrm{\RN{2}} 
\lesssim &  \left\|\mcP_1^{-3}\whE_1\mcP_1^{0}\whE_1\mcP_1^{0}\whE_1\mcP_{U_{1\perp}}V_1\right\| + \left\|\mcP_1^{-2}\whE_1\mcP_1^{-1}\whE_1\mcP_1^{0}\whE_1\mcP_{U_{1\perp}}V_1\right\| + \left\|\mcP_1^{-2}\whE_1\mcP_1^{0}\whE_1\mcP_1^{-1}\whE_1\mcP_{U_{1\perp}}V_1\right\| \\
+ & \left\|\mcP_1^{-1}\whE_1\mcP_1^{-2}\whE_1\mcP_1^{0}\whE_1\mcP_{U_{1\perp}}V_1\right\| + \left\|\mcP_1^{-1}\whE_1\mcP_1^{0}\whE_1\mcP_1^{-2}\whE_1\mcP_{U_{1\perp}}V_1\right\| + \left\|\mcP_1^{-1}\whE_1\mcP_1^{-1}\whE_1\mcP_1^{-1}\whE_1\mcP_{U_{1\perp}}V_1\right\|   ,\\
\mathrm{\RN{3}} 
\lesssim & \left\|\mcP_{U_{1\perp}}\mcP_1^{0}\whE_1\mcP_1^{-2}\whE_1\mcP_1^{0}V_1\right\| + \left\|\mcP_{U_{1\perp}}\sum_{k_1=3}^{+\infty} S_{G_1, k_1}\left(\whE_1\right)V_1\right\| ,\\
\mathrm{\RN{4}} 
\lesssim & \left\|\mcP_1^{0}\whE_1\mcP_1^{-3}\whE_1\mcP_1^{0}\whE_1\mcP_1^{0} V_1\right\| +  \left\|\mcP_1^{0}\whE_1\mcP_1^{-2}\whE_1\mcP_1^{-1}\whE_1\mcP_1^{0} V_1\right\| \\
+ & \left\|\mcP_1^{0}\whE_1\mcP_1^{-1}\whE_1\mcP_1^{-2}\whE_1\mcP_1^{0} V_1\right\| +  \left\|\mcP_1^{0}\whE_1\mcP_1^{0}\whE_1\mcP_1^{-3}\whE_1\mcP_1^{0} V_1\right\| .
\end{align*}
Applying similar arguments, we obtain the bounds for \RN{2}, \RN{3}, \RN{4}.

\end{proof}

\begin{lemma}\label{lemma: remaining terms of the first-order perturbation term in tensor regression without sample splitting}
Under the same setting of Theorem~\ref{thm: main theorem in tensor regression without sample splitting}, let $\whE_j$ be defined by \eqref{eq: definition of whEj in tensor regression without sample splitting} For any $j=1,2,3$. Then with probability at least $1-\exp(-cn) - \frac{1}{p^C} - \P\left(\mcE_\Delta\right) - \P\left(\mcE_{U}^{\text{reg}}\right)$, where $c$ and $C$ are two universal constants ,it holds that
\begin{align}
& \left\|U_{j\perp}^{\top} \whE_j U_{j\perp}-U_{j\perp}^{\top} \whZ_j\left(\mcP_{U_{j+2}} \otimes \mcP_{U_{j+1}}\right) \whZ_j^{\top} U_{j\perp}\right\| \lesssim \frac{\sigma_{\xi}^3}{\ulambda \sigma^3} \cdot \frac{\op^{3/2}}{n^{3/2}} \label{eq: high-prob upper bound of U1ptEhat1U1p - U1ptZhat1(PU3oPU2)Zhat1tU1p in tensor regression without sample splitting}\\
& \left\|V_j^{\top} \mcP_{U_{j \perp}} \whE_j U_{j \perp}-V_j^{\top} \mcP_{U_{j \perp}} \whZ_j\left(\mcP_{U_{j+2}} \otimes \mcP_{U_{j+1}}\right) \whZ_j^{\top} U_{j \perp}\right\| \lesssim \frac{\sigma_\xi^3}{\ulambda\sigma^3} \cdot \left(\frac{\op \sqrt{\oR \log (\op)}}{n^{3/2}} + \Delta\cdot \frac{\op^{3/2}}{n^{3/2}}\right) \label{eq: high-prob upper bound of V1V1tP1(0)Ehat1PU1p - V1V1tP1(0)Zhat1(PU3oPU2)Zhat1tU1p in tensor regression without sample splitting} \\
& \left\|U_{j\perp}^{\top} \whE_j P_j^{-1}-U_{j\perp}^{\top} \whZ_j\left(U_{j+2} \otimes U_{j+1}\right) G_j^{\top}\left(G_j G_j^{\top}\right)^{-1} U_j^{\top}\right\| \lesssim \frac{\sigma_{\xi}^2}{\ulambda^2 \sigma^2} \cdot \frac{\op}{n} \label{eq: high-prob upper bound of U1ptEhat1P1^(-1) - U1ptZhat1(U3oU2)G1t(G1G1t)^(-1)U1t in tensor regression without sample splitting}\\
& \left\|V_j^{\top}  U_{j \perp}^{\top} \whE_j P_j^{-1} -V_j^{\top} U_{j \perp}^{\top} \whZ_j\left(U_{j+2} \otimes U_{j+1}\right) G_j^{\top}\left(G_j G_j^{\top}\right)^{-1} U_j^{\top} \right\| \lesssim \frac{\sigma_\xi^2}{\ulambda^2\sigma^2}\cdot \left(\frac{\sqrt{\oR \op \log (\op)}}{n} + \Delta\cdot \frac{\op}{n}\right) \label{eq: high-prob upper bound of V1V1tP1(0)Ehat1P1^(-1) - V1V1tP1(0)Zhat1(U3oU2)G1t(G1G1t)^(-1)U1t in tensor regression without sample splitting}.
\end{align}

\end{lemma}

\begin{proof}

By symmetry, it suffices to consider upper bounds of:
\begin{align*}
\mathrm{\RN{1}} = & \left\|U_{1 \perp}^{\top} \whE_1 U_{1 \perp}-U_{1 \perp}^{\top} \whZ_1\left(\mcP_{U_3} \otimes \mcP_{U_2}\right) \whZ_1^{\top} U_{1 \perp}\right\| \\
\mathrm{\RN{2}} = & \left\| V_1 V_1^{\top} \mcP_{U_{1 \perp}} \whE_1 U_{1 \perp}-V_1^{\top} \mcP_{U_{1 \perp}} \whZ_1\left(\mcP_{U_3} \otimes \mcP_{U_2}\right) \whZ_1^{\top} U_{1 \perp} \right\| \\
\mathrm{\RN{3}} = & \left\| U_{1 \perp}^{\top} \whE_1 P_1^{-1}-U_{1 \perp}^{\top} \whZ_1\left(U_3 \otimes U_2\right) G_1^{\top}\left(G_1 G_1^{\top}\right)^{-1} U_1^{\top} \right\| \\
\mathrm{\RN{4}} = & \left\|V_1 V_1^{\top} U_{1 \perp}^{\top} \whE_1 P_1^{-1}-V_1 V_1^{\top} U_{1 \perp}^{\top} \whZ_1\left(U_3 \otimes U_2\right) G_1^{\top}\left(G_1 G_1^{\top}\right)^{-1} U_1^{\top}\right\|.
\end{align*}


It follows that
\begin{align*}
\mathrm{\RN{1}} 
\leq & \left\|U_{1 \perp}^{\top} \whZ_1\left[\left(\mcP_{\whU_3^{(1)}} - \mcP_{U_3}\right) \otimes \mcP_{U_2}\right] \whZ_1^{\top} U_{1 \perp}\right\| + \left\|U_{1 \perp}^{\top} \whZ_1\left[\mcP_{U_3} \otimes \left(\mcP_{\whU_2^{(1)}}-\mcP_{U_2}\right)\right] \whZ_1^{\top} U_{1 \perp}\right\| \\
+ & \left\|U_{1 \perp}^{\top} \whZ_1\left[\left(\mcP_{\whU_3^{(1)}} -\mcP_{U_3}\right) \otimes \left(\mcP_{\whU_2^{(1)}}-\mcP_{U_2}\right)\right] \whZ_1^{\top} U_{1 \perp}\right\| \\
\leq & \sup_{\substack{W_2 \in \mathbb{R}^{p_2\times r_2}, \left\|W_2\right\|=1\\ W_3 \in \mathbb{R}^{p_3\times r_3}, \left\|W_3\right\|=3}} \left\|U_{1\perp}^{\top}\whZ_1\left(W_3\otimes W_2\right)\right\|^2 \cdot \left(\left\|\mcP_{\whU_2^{(1)}} -\mcP_{U_2}\right\| + \left\|\mcP_{\whU_3^{(1)}} -\mcP_{U_3}\right\| + \prod_{j=2}^3 \left\|\mcP_{\whU_j^{(1)}} -\mcP_{U_j}\right\| \right) \\
\lesssim & \left(\frac{\sigma_\xi}{\sigma}\cdot \sqrt{\frac{\op}{n}}\right)^2 \cdot \left[\frac{\sigma_\xi}{\ulambda\sigma} \cdot \sqrt{\frac{\op}{n}} + \left(\frac{\sigma_\xi}{\ulambda\sigma} \cdot \sqrt{\frac{\op}{n}}\right)^2\right] \lesssim \frac{\sigma_\xi^3}{\ulambda\sigma^3} \cdot \frac{\op^{3/2}}{n^{3/2}}.
\end{align*}

Applying similar arguments and Lemma~\ref{lemma: high-prob upper bound of the first-order perturbation terms in tensor regression without sample splitting}, \ref{lemma: high-prob upper bound of the higher-order perturbation terms in tensor regression without sample splitting}, we obtain the bounds for \RN{2}, \RN{3}, \RN{4}.

\end{proof}

\subsection{Upper Bound of Leading Terms in the Spectral Representation}

\begin{lemma}\label{lemma: high-prob upper bound of leading terms in the spectral representation across three modes in tensor regression without sample splitting}

Under the same setting of Theorem~\ref{thm: main theorem in tensor regression without sample splitting}, let $\whE_j$ be defined by \eqref{eq: definition of whEj in tensor regression without sample splitting} and $\mcP_j^{-s}=U_j\left(G_jG_j^{\top}\right)^{-2s}U_j^{\top}$ for any $j=1,2,3$. Then with probability at least $1-\exp(-cn) - \frac{1}{p^C} - \P\left(\mcE_\Delta\right) - \P\left(\mcE_{U}^{\text{reg}}\right)$, it follows that:
\begin{align}
& \left\|\mcA \times_j \mcP_{U_{j \perp}} \whE_j \mcP_j^{-1} \times_{j+1} \mcP_{U_{{j+1} \perp}} \whE_{j+1} \mcP_{j+1}^{-1} \times_{j+2} \mcP_{U_{{j+2} \perp}} \whE_{j+2} \mcP_{j+2}^{-1}\right\|_{\mathrm{F}} \notag \\
\lesssim & \left\| \mcA \right\|_{\mathrm{F}} \cdot \left[\frac{\sigma_\xi^3 }{\ulambda^3\sigma^3}\cdot \left(\frac{\oor^{3/2}\log(\op)^{3/2}}{n^{3/2}} + \Delta\cdot \frac{\op^{1/2}\oR\log(\op)}{n^{3/2}} + \Delta^2 \cdot \frac{\op \oor \log (\op)}{n^{3 / 2}} + \Delta^3\cdot \frac{\op^{3/2}}{n^{3/2}}\right)\right] \label{eq: high-prob upper bound of AoPU1pEhat1P1(-1)oPU2pEhat2P2(-1)oPU3pEhat3P3(-1) in tensor regression without sample splitting}
\end{align}

\end{lemma}

\begin{proof}

By symmetry, it suffices to consider the following upper bound:
\begin{align*}
& \left\|\left(\mcP_{3}^{-1} \whE_3 \mcP_{U_{3 \perp}} \otimes \mcP_{2}^{-1} \whE_2 \mcP_{U_{2 \perp}}\right) A_1^{\top} \mcP_{U_{1 \perp}} \whE_1 \mcP_1^{-1}\right\|_{\mathrm{F}} \\
\leq & \underbrace{\left\|\left(\left(G_3G_3^{\top}\right)^{-1}G_3\left(\mcP_{U_2} \otimes \mcP_{U_1}\right) \whZ_3^{\top} \mcP_{U_{3 \perp}} \otimes \left(G_2G_2^{\top}\right)^{-1}G_2\left(\mcP_{U_1} \otimes \mcP_{U_3}\right) \whZ_2^{\top} \mcP_{U_{2 \perp}}\right)A_1 \whZ_1\left(U_3 \otimes U_2\right)G_1^{\top}\left(G_1G_1^{\top}\right)^{-1}\right\|_{\mathrm{F}}}_{\eqref{eq: high-prob upper bound of AoPU1pEhat1P1(-1)oPU2pEhat2P2(-1)oPU3pEhat3P3(-1) in tensor regression with sample splitting}} \\
+ & \sum_{j=1}^{3} \underbrace{\left\| (G_{j+2} G_{j+2}^{\top} )^{-1} G_{j+2} (\mcP_{U_{j+1}} \otimes \mcP_{U_j} ) \whZ_{j+2}^{\top} \mcP_{U_{j+2\perp}}V_{j+2}\right\|}_{\eqref{eq: high-prob upper bound of V1tP1(0)Ehat1P1(-1/2) in tensor regression without sample splitting}} 
 \cdot \underbrace{\left\| (G_{j+1} G_{j+1}^{\top} )^{-1} G_{j+1} (\mcP_{U_j} \otimes \mcP_{U_{j+2}} ) \whZ_{j+1}^{\top} \mcP_{U_{j+1 \perp}}V_{j+1}\right\|}_{\eqref{eq: high-prob upper bound of V1tP1(0)Ehat1P1(-1/2) in tensor regression without sample splitting}} \\
& \cdot \left\|A_j\right\|_{\mathrm{F}} \cdot \underbrace{\left\|V_j^{\top}\mcP_{U_{j \perp}} \whE_j P_j^{-1}-V_jV_j^{\top}\mcP_{U_{j \perp}} \whZ_j\left(U_{j+2} \otimes U_{j+1}\right) G_j^{\top}\left(G_j G_j^{\top}\right)^{-1} U_j^{\top}\right\|}_{\eqref{eq: high-prob upper bound of V1V1tP1(0)Ehat1P1^(-1) - V1V1tP1(0)Zhat1(U3oU2)G1t(G1G1t)^(-1)U1t in tensor regression without sample splitting}} \\
+ & \sum_{j=1}^3 \underbrace{\left\|\left(P_{j+2}^{-1} \whE_{j+2}^{\top} \mcP_{U_{j+2 \perp}}V_{j+2}V_{j+2}^{\top} 
- U_{j+2} \left(G_{j+2} G_{j+2}^{\top}\right)^{-1} G_{j+2}^\top \left(U_{j+1} \otimes U_{j}\right)^{\top} \whZ_{j+2}^{\top} \mcP_{U_{j+2 \perp}}V_{j+2}V_{j+2}^{\top}\right)\right\|}_{\eqref{eq: high-prob upper bound of V1V1tP1(0)Ehat1P1^(-1) - V1V1tP1(0)Zhat1(U3oU2)G1t(G1G1t)^(-1)U1t in tensor regression without sample splitting}} \\
& \cdot \underbrace{\left\|\left(P_{j+1}^{-1} \whE_{j+1}^{\top} \mcP_{U_{j+1 \perp}}V_{j+1}V_{j+1}^{\top} 
- U_{j+1} \left(G_{j+1} G_{j+1}^{\top}\right)^{-1} G_{j+1}^\top \left(U_{j+2} \otimes U_{j}\right)^{\top} \whZ_{j+1}^{\top} \mcP_{U_{j+1 \perp}}^{\top}\right)\right\|}_{\eqref{eq: high-prob upper bound of V1V1tP1(0)Ehat1P1^(-1) - V1V1tP1(0)Zhat1(U3oU2)G1t(G1G1t)^(-1)U1t in tensor regression without sample splitting}} \\
& \cdot \left\|A_j\right\|_{\mathrm{F}} \cdot \underbrace{\left\|V_j^{\top}\mcP_{U_{j \perp}} \whZ_j\left(U_{j+2} \otimes U_{j+1}\right) G_j^{\top}\left(G_j G_j^{\top}\right)^{-1} U_j^{\top}\right\|}_{\eqref{eq: high-prob upper bound of V1tP1(0)Ehat1P1(-1/2) in tensor regression without sample splitting}} \\
+ & \prod_{j=1}^3 \underbrace{\left\|P_{j}^{-\top} \whE_{j}^{\top} \mcP_{U_{j \perp}}V_jV_j^{\top} 
- U_{j} \left(G_{j} G_{j}^{\top}\right)^{-\top} G_{j}^\top \left(U_{j+2} \otimes U_{j+1}\right)^{\top} \whZ_{j}^{\top} \mcP_{U_{j \perp}}V_1\right\|}_{\eqref{eq: high-prob upper bound of V1V1tP1(0)Ehat1P1^(-1) - V1V1tP1(0)Zhat1(U3oU2)G1t(G1G1t)^(-1)U1t in tensor regression without sample splitting}} \cdot \left\|A_j\right\|_{\mathrm{F}}\\
\lesssim & \left\|\mcA\right\|_{\mathrm{F}} \cdot \left[\frac{\sigma_\xi^3}{\ulambda^3\sigma^3}\cdot \left(\frac{\oor^3 \log (\op)}{n^{3 / 2}}+\Delta \cdot \frac{\op^{1 / 2} \oR \log (\op)}{n^{3 / 2}}+\Delta^2 \cdot \frac{\op \oR \log (\op)}{n^{3 / 2}}+\Delta^3 \cdot \frac{\op^{3 / 2}}{n^{3 / 2}}\right)\right] .
\end{align*}

\end{proof}

\begin{lemma}\label{lemma: high-prob upper bound of leading terms in the spectral representation across two modes in tensor regression without sample splitting}

Under the same setting of Theorem~\ref{thm: main theorem in tensor regression without sample splitting}, let $\whE_j$ be defined by \eqref{eq: definition of whEj in tensor regression without sample splitting} and $\mcP_j^{-s}=U_j\left(G_jG_j^{\top}\right)^{-2s}U_j^{\top}$ for any $j=1,2,3$.Then with probability at least $1-\exp(-cn) - \frac{1}{p^C} - \P\left(\mcE_\Delta\right) - \P\left(\mcE_{U}^{\text{reg}}\right)$, it follows that:
\begin{align}
& \left\|\mcA \times_j U_j \times_{j+1} \mcP_{U_{j+1\perp}}\whE_{j+1}\mcP_{j+1}^{-1} \times_{j+2} \mcP_{U_{j+2\perp}}\whE_{j+2}\mcP_{j+2}^{-1}\right\|_{\mathrm{F}} \notag \\
&\leq\left\|\mcA \times_j U_j\right\|_{\mathrm{F}} \cdot \left[\frac{\sigma_{\xi}^2}{\ulambda^2 \sigma^2} \cdot\left(\frac{\oor\log(\op)}{n}+\Delta \cdot \frac{\sqrt{\oR \op\log(\op)}}{n}+\Delta^2 \cdot \frac{\op}{n}\right)\right] \label{eq: high-prob upper bound of AoPU1oPU2pEhat2P2(-1)oPU3pEhat3P3(-1) in tensor regression without sample splitting}
\end{align}
for any $j=1,2,3$.

\end{lemma}

\begin{proof}

First, consider
\begin{align*}
& \left\|\left(\mcP_3^{-1} \whE_3 \mcP_{U_{3 \perp}}V_3V_3^{\top} \otimes \mcP_2^{-1} \whE_2 \mcP_{U_{2 \perp}}V_2V_2^{\top}\right) A_1^{\top} \mcP_{U_1}\right\|_{\mathrm{F}} \\
\leq & \underbrace{\left\|\left(G_3 G_3^{\top}\right)^{-1} G_3\left(\left(\mcP_{U_2} \otimes \mcP_{U_1}\right) \whZ_3^{\top} \mcP_{U_{3 \perp}} \otimes \left(G_2 G_2^{\top}\right)^{-1} G_2\left(\mcP_{U_1} \otimes \mcP_{U_3}\right) \whZ_2^{\top} \mcP_{U_{2 \perp}}\right) A_1^{\top} \mcP_{U_1}\right\|}_{\eqref{eq: high-prob upper bound of AoPU1oPU2pEhat2P2(-1)oPU3pEhat3P3(-1) in tensor regression with sample splitting}} \\
+ & \underbrace{\left\|V_3^{\top}\mcP_{U_{3 \perp}} \whZ_3 (U_1 \otimes U_2 ) G_3^{\top} (G_3 G_3^{\top} )^{-1} U_3^{\top}
\right\|}_{\eqref{eq: high-prob upper bound of V1tP1(0)Ehat1P1(-1/2) in tensor regression without sample splitting}} 
 \underbrace{\left\|V_2^{\top}\mcP_{U_{2 \perp}} \whE_2 P_2^{-1} - V_2V_2^{\top}\mcP_{U_{2 \perp}} \whZ_2 (U_3 \otimes U_1 ) G_2^{\top} (G_2 G_2^{\top} )^{-1} U_2^{\top}
\right\|}_{\eqref{eq: high-prob upper bound of V1V1tP1(0)Ehat1P1^(-1) - V1V1tP1(0)Zhat1(U3oU2)G1t(G1G1t)^(-1)U1t in tensor regression without sample splitting}}  \left\|A_1^{\top} \mcP_{U_1}\right\|_{\mathrm{F}}\\
+ & \underbrace{\left\|V_3^{\top}\mcP_{U_{3 \perp}} \whE_3 P_3^{-1} - V_3V_3^{\top}\mcP_{U_{3 \perp}} \whZ_3 (U_1 \otimes U_2 ) G_3^{\top} (G_3 G_3^{\top} )^{-1} U_3^{\top}
\right\|}_{\eqref{eq: high-prob upper bound of V1V1tP1(0)Ehat1P1^(-1) - V1V1tP1(0)Zhat1(U3oU2)G1t(G1G1t)^(-1)U1t in tensor regression without sample splitting}} 
 \underbrace{\left\|V_2^{\top}\mcP_{U_{2 \perp}} \whZ_2 (U_3 \otimes U_1 ) G_2^{\top} (G_2 G_2^{\top} )^{-1} U_2^{\top} \right\|}_{\eqref{eq: high-prob upper bound of V1tP1(0)Ehat1P1(-1/2) in tensor regression without sample splitting}}  \left\| A_1^{\top} \mcP_{U_1}\right\|_{\mathrm{F}} \\
+ & \underbrace{\left\|V_3^{\top}\mcP_{U_{3 \perp}} \whE_3 P_3^{-1} - V_3V_3^{\top}\mcP_{U_{3 \perp}} \whZ_3\left(U_1 \otimes U_2\right) G_3^{\top}\left(G_3 G_3^{\top}\right)^{-1} U_3^{\top}
\right\|}_{\eqref{eq: high-prob upper bound of V1V1tP1(0)Ehat1P1^(-1) - V1V1tP1(0)Zhat1(U3oU2)G1t(G1G1t)^(-1)U1t in tensor regression without sample splitting}} \\
& \cdot \underbrace{\left\|V_2^{\top}\mcP_{U_{2 \perp}} \whE_2 P_2^{-1} - V_2V_2^{\top}\mcP_{U_{2 \perp}} \whZ_2\left(U_3 \otimes U_1\right) G_2^{\top}\left(G_2 G_2^{\top}\right)^{-1} U_2^{\top}
\right\|}_{\eqref{eq: high-prob upper bound of V1V1tP1(0)Ehat1P1^(-1) - V1V1tP1(0)Zhat1(U3oU2)G1t(G1G1t)^(-1)U1t in tensor regression without sample splitting}} \cdot\left\| A_1^{\top} \mcP_{U_1}\right\|_{\mathrm{F}} \\
\lesssim & \left\|\mcA \times_1 U_1\right\|_{\mathrm{F}} \cdot \left[\frac{\sigma_{\xi}^2}{\ulambda^2 \sigma^2} \cdot\left(\frac{\oor\log(\op)}{n}+\Delta \cdot \frac{\sqrt{\oR \op\log(\op)}}{n}+\Delta^2 \cdot \frac{\op}{n}\right)\right].
\end{align*}

\end{proof}

In the following subsections, we established upper bounds for perturbation terms of varying orders in the spectral representation under the setting of tensor regression with sample splitting. In particular, we will show that the first-order perturbation term is the leading term. Throughout this section, we assume that $\Delta \geq (\sigma_\xi / \sigma)\sqrt{\op / n}$, which implies that the initial estimate satisfies the minimax lower bound as well. Different from the scenario without sample splitting, sample splitting removes the dependency between the initial estimate and debiasing procedure. As a result, the higher-order pertubation terms will vanish at a faster rate. 

\section{Preliminary Upper Bounds for Tensor Regression with Sample splitting}\label{sec: Preliminary Upper Bounds for Tensor Regression with Sample splitting}

This section contains the essential lemmas for the proof of Theorem~\ref{thm: main theorem in tensor regression with sample splitting}. Although the following contents are similar to the lemmas for the proof of Theorem~\ref{thm: main theorem in tensor regression without sample splitting}, the double sample-splitting largely reduces the upper bound of negligible terms by removing the dependence between the projection and bias-correction using two separate datasets. As the sample-splitting divides the original dataset to two subsets, by symmetry, we only consider one subset.

After the power iteration and projection in the algorithm with sample splitting in Section \ref{sec: Debiased Estimator of Linear Functionals with Sample Splitting}, for any $j=1,2,3$, we know that $\whU_j^{(\mathrm{\RN{1}})}$ contains the top- $r_j$ eigenvectors of 
$$
\widehat{\mcT}^{\text{unbs},(\mathrm{\RN{1}})}_j\left(\mcP_{\whU_{j+2}^{(0),(\mathrm{\RN{2}})}} \otimes \mcP_{\whU_{j+1}^{(0),(\mathrm{\RN{2}})}}\right) \widehat{\mcT}^{\text{unbs},(\mathrm{\RN{1}})\top}_j.
$$
Consequently, $\whU_j^{(\mathrm{\RN{1}})} \whU_j^{{(\mathrm{\RN{1}})} \top}$ is the spectral projector for the left top- $r_j$ left eigenvectors of
\begin{align*}
\widehat{T}_j^{\text{unbs},(\mathrm{\RN{1}})} \left(\mcP_{\whU_{j+2}^{(0),(\mathrm{\RN{2}})}} \otimes \mcP_{\whU_{j+1}^{(0),(\mathrm{\RN{2}})}}\right) \widehat{T}_j^{\text{unbs},(\mathrm{\RN{1}})\top}
= T_j\left(\mcP_{U_{j+1}} \otimes \mcP_{U_{j+2}}\right) T_j^{\top}+\whE_j^{(\mathrm{\RN{1}})} 
= U_jG_jG_j^{\top} U_j^{\top} + \whE_j^{(\mathrm{\RN{1}})},
\end{align*}
where
\begin{align}
\whE_j^{(\mathrm{\RN{1}})}
= & T_j\left(\mcP_{\whU_{j+2}^{(0),(\mathrm{\RN{2}})}} \otimes \mcP_{\whU_{j+1}^{(0),(\mathrm{\RN{2}})}}\right)\whZ_j^{(\mathrm{\RN{1}}),\top} + \whZ_j^{(\mathrm{\RN{1}})}\left(\mcP_{\whU_{j+2}^{(0),(\mathrm{\RN{2}})}} \otimes \mcP_{\whU_{j+1}^{(0),(\mathrm{\RN{2}})}}\right)T_j^{\top} \notag \\
+ & T_j\left(\left(\mcP_{\whU_{j+2}^{(0),(\mathrm{\RN{2}})}}-\mcP_{U_{j+2}}\right) \otimes \mcP_{\whU_{j+1}^{(0),(\mathrm{\RN{2}})}}\right) T_j^{\top}+T_j\left(\mcP_{U_{j+2}} \otimes \left(\mcP_{\whU_{j+1}^{(1)}}-\mcP_{U_{j+1}}\right)\right) T_j^{\top} \label{eq: definition of whEj in tensor regression with sample splitting}\\
+ & T_j\left(\left(\mcP_{\whU_{j+2}^{(0),(\mathrm{\RN{2}})}}-\mcP_{U_{j+2}}\right) \otimes \left(\mcP_{\whU_{j+1}^{(0),(\mathrm{\RN{2}})}}-\mcP_{U_{j+1}}\right)\right) T_j^{\top}. \notag  
\end{align}

If $\left\|\whE_j^{(\mathrm{\RN{1}})}\right\| \leq \frac{1}{2}\ulambda^2$, then by Theorem 1 \citep{xia2021normal}, the following equation holds
$$
\whU_j^{(\mathrm{\RN{1}})} \whU_j^{{(\mathrm{\RN{1}})} \top}-U_j U_j^{\top}=\sum_{k_j= 1}^{+\infty}\mathcal{S}_{G_j, k_j}\left(\whE_j^{(\mathrm{\RN{1}})}\right).
$$
Here, for each positive integer $k$
\begin{align}
\mathcal{S}_{G_j, k_j}\left(\whE_j^{(\mathrm{\RN{1}})}\right)=\sum_{s_1+\cdots+s_{k_j+1}=k_j}(-1)^{1+\tau(\mathbf{s})} \cdot \mcP_j^{-s_1} \whE_j^{(\mathrm{\RN{1}})} \mcP_j^{-s_2} \whE_j^{(\mathrm{\RN{1}})} \mcP_j^{-s_3} \cdots \mcP_j^{-s_{k_j}} \whE_j^{(\mathrm{\RN{1}})} \mcP_j^{-s_{k_j+1}} \label{eq: definition of Sj(whEj) in tensor regression with sample splitting}
\end{align}
where $s_1, \cdots, s_{k_j+1}$ are non-negative integers and $\tau(\mathbf{s})=\sum_{j=1}^{k_j+1} \mathbb{I}\left(s_{k_j}>0\right)$, 
$
\mcP_j^{-k}= U_j\left(G_jG_j^{\top}\right)^{-k}U_j^{\top}
$
for any $k \geq 1$ and $\mcP_j^{0}=U_{j\perp}U_{j\perp}^{\top}$. It follows that
\begin{align}
\mathcal{S}_{G_j, 1}\left(\whE_j^{(\mathrm{\RN{1}})}\right)
= & P_j^{-1} \whE_j^{(\mathrm{\RN{1}})} P_j^{0}+P_j^{0} \whE_j^{(\mathrm{\RN{1}})} P_j^{-1} \notag \\
= & U_j \left(G_j G_j^{\top}\right)^{-1} G_j \left({U_{j+1}} \otimes U_{j+2}\right)^{\top} \left(\mcP_{\whU_{j+1}^{(0),(\mathrm{\RN{2}})}}\otimes \mcP_{\whU_{j+2}^{(0),(\mathrm{\RN{2}})}}\right)\whZ_j^{(\mathrm{\RN{1}}),\top} \mcP_{U_{j\perp}} 
 \notag \\
+ & \mcP_{U_{j\perp}} \whZ_j^{(\mathrm{\RN{1}})}\left(\mcP_{\whU_{j+2}^{(0),(\mathrm{\RN{2}})}} \otimes \mcP_{\whU_{j+1}^{(0),(\mathrm{\RN{2}})}}\right) \left(U_{j+2} \otimes U_{j+1}\right)G_j^{\top}\left(G_j G_j^{\top}\right)^{-1} U_j^{\top}, \label{eq: definition of S1(whEj) in tensor regression with sample splitting}
\end{align}
for any $j=1,2,3$, where the second equality, the third inequality come from the definition that $P_j^{-1}=U_j\left(G_jG_j^{\top}\right)^{-1}U_j^{\top}$.

Here, note that $\left\|\whE_j^{(\mathrm{\RN{1}})}\right\| \leq \kappa\ulambda \sqrt{\frac{\op}{n}}$. Then the condition, $\left\|\whE_j^{(\mathrm{\RN{1}})}\right\| \leq \frac{1}{2}\ulambda^2$, for Theorem 1 in \citet{xia2021normal} is satisfied provied that $n \gtrsim \kappa^2 \op / \ulambda^2$.

In the subsequent sections, we assume that the following events
$$
\left\|\mcP_{\whU_j^{(0), (\mathrm{\RN{1}})}} - \mcP_{U_j}\right\| \leq \frac{\sigma_\xi}{\sigma}\sqrt{\frac{\op}{n}}, \left\|\mcP_{\whU_j^{(0), (\mathrm{\RN{2}})}} - \mcP_{U_j}\right\| \leq \frac{\sigma_\xi}{\sigma}\sqrt{\frac{\op}{n}},
$$
hold with probability at least $1- \mathbb{P}\left(\mcE_{U}^{\text{reg}}\right)$, where event $\mcE_{U}^{\text{reg}}$ is defined by $\mcE_{U}^{\text{reg}}= \left\{\left\|\mcP_{\whU_j^{(0)}} - \mcP_{U_j}\right\| > \frac{\sigma_\xi}{\sigma}\sqrt{\frac{\op}{n}}\right\}$, where $\mcP_{\whU_j^{(0)}}$ can be either $\mcP_{\whU_j^{(0), (\mathrm{\RN{1}})}}$ or $\mcP_{\whU_j^{(0), (\mathrm{\RN{2}})}}$.

Then by Lemma~\ref{lemma: error contraction of l2 error of singular space in tensor regression}, we know that 
$ \left\|\mcP_{\whU_j^{(\mathrm{\RN{1}})}} - \mcP_{U_j}\right\| \leq \frac{\sigma_\xi}{\sigma}\sqrt{\frac{\op}{n}}$ and $
\left\|\mcP_{\whU_j^{(\mathrm{\RN{2}})}} - \mcP_{U_j}\right\| \leq \frac{\sigma_\xi}{\sigma}\sqrt{\frac{\op}{n}}$
hold with probability at least $1-\exp(-c\op) - \mathbb{P}(\mcE_{U}^{\text{reg}})$ for any $j=1,2,3$.

Besides, we assume that the initial error bound 
$
\left\|\whT^{\text{init}, (\mathrm{\RN{1}})} - \mcT\right\|_{\mathrm{F}} \leq \Delta, \left\|\whT^{\text{init}, (\mathrm{\RN{2}})} - \mcT\right\|_{\mathrm{F}} \leq \Delta
$
hold with probability at least $1-\mathbb{P}(\mcE_{\Delta})$, where event $\mcE_{\Delta}$ is given by $\mcE_{\Delta}= \left\{\left\|\whT^{\text{init}} - \mcT\right\|_{\mathrm{F}} > \Delta \right\}$, where $\whT^{\text{init}}$ can be either $\whT^{\text{init}, (\mathrm{\RN{1}})}$ or $\whT^{\text{init}, (\mathrm{\RN{2}})}$.

\subsection{Preliminary Bounds in the Proof of Theorem~\ref{thm: main theorem in tensor regression with sample splitting}} \label{subsec: preliminary bounds in the proof of main theorems in tensor regression with sample splitting}

\begin{proposition}\label{prop: high-probability upper bound of AoPU1pEhat1P_U1oPU2pEhat2oPU3pEhat3 with sample splitting}

Under the same setting of Theorem~\ref{thm: main theorem in tensor regression with sample splitting}, with probability at least $1-\exp(-cn) - \frac{1}{p^C} - \P\left(\mcE_\Delta\right) - \P\left(\mcE_{U}^{\text{reg}}\right)$, where $c$ and $C$ are two universal constants ,it holds that
\begin{align}
& \left\|\mcA \times_j \mcP_{U_{j\perp}}\left(\mcP_{\whU_j^{(\mathrm{\RN{1}})}} - \mcP_{U_j}\right)\mcP_{U_j} \times_{j+1} \mcP_{U_{j+1\perp}}\left(\mcP_{\whU_{j+1}^{(\mathrm{\RN{1}})}} - \mcP_{U_{j+1}}\right)\mcP_{U_{{j+1}}} \times_{j+2} \mcP_{U_{j+2\perp}}\left(\mcP_{\whU_{j+2}^{(\mathrm{\RN{1}})}} - \mcP_{U_{j+2}}\right)\mcP_{U_{j+2}}  \right\|_{\mathrm{F}} \notag \\
\lesssim & \left\|\mcA\right\|_{\mathrm{F}} \cdot \left(\frac{\sigma_{\xi}^3}{\ulambda^3 \sigma^3} \cdot \frac{\oor^{3/2}\log(\op)^{3/2}}{n^{3/2}} +\frac{\sigma_{\xi}^3}{\ulambda^3 \sigma^3} \cdot \Delta \cdot \frac{\oR^{3/2} \log (\op)^{3/2}}{n^{3/2}}\right) \label{eq: high-prob upper bound of AoPU1p(PUhat1-P_U1)PU1oPU2p(PUhat2-P_U2)PU2oPU3p(PUhat3-P_U3)PU3 in tensor regression with sample splitting}\\
& \left\|\mcA \times_j \mcP_{U_{j\perp}}\left(\mcP_{\whU_j^{(\mathrm{\RN{1}})}} - \mcP_{U_j}\right)\mcP_{U_{j\perp}} \times_{j+1} \mcP_{U_{j+1\perp}}\left(\mcP_{\whU_{j+1}^{(\mathrm{\RN{1}})}} - \mcP_{U_{j+1}}\right)\mcP_{U_{{j+1}}} \times_{j+2} \mcP_{U_{j+2\perp}}\left(\mcP_{\whU_{j+2}^{(\mathrm{\RN{1}})}} - \mcP_{U_{j+2}}\right)\mcP_{U_{j+2}}  \right\|_{\mathrm{F}} \notag \\
\lesssim & \left\|\mcA\right\|_{\mathrm{F}} \cdot  \left(\frac{\sigma_{\xi}^4}{\ulambda^4 \sigma^4} \cdot \frac{\oor\op^{1/2}\log(\op)^2}{n^2} + \frac{\sigma_{\xi}^4}{\ulambda^4 \sigma^4} \cdot \Delta \cdot \frac{\oR^{3/2} \op^{1/2} \log (\op)^{3/2}}{n^2}\right) \label{eq: high-prob upper bound of AoPU1p(PUhat1-P_U1)PU1poPU2p(PUhat2-P_U2)PU2oPU3p(PUhat3-P_U3)PU3 in tensor regression with sample splitting}\\
& \left\|\mcA \times_j \mcP_{U_{j\perp}}\left(\mcP_{\whU_j^{(\mathrm{\RN{1}})}} - \mcP_{U_j}\right)\mcP_{U_{j}} \times_{j+1} \mcP_{U_{j+1\perp}}\left(\mcP_{\whU_{j+1}^{(\mathrm{\RN{1}})}} - \mcP_{U_{j+1}}\right)\mcP_{U_{{j+1\perp}}} \times_{j+2} \mcP_{U_{j+2\perp}}\left(\mcP_{\whU_{j+2}^{(\mathrm{\RN{1}})}} - \mcP_{U_{j+2}}\right)\mcP_{U_{j+2\perp}}  \right\|_{\mathrm{F}} \notag \\
\lesssim & \left\|\mcA\right\|_{\mathrm{F}} \cdot \left(\frac{\sigma_{\xi}^5}{\ulambda^5 \sigma^5} \cdot \frac{\oor^{3/2}\op\log(\op)^{3/2}}{n^{5/2}} + \frac{\sigma_{\xi}^5}{\ulambda^5 \sigma^5} \cdot \Delta \cdot \frac{\oR^{3/2} \op \log (\op)^{3/2}}{n^{5/2}}\right) \label{eq: high-prob upper bound of AoPU1p(PUhat1-P_U1)PU1oPU2p(PUhat2-P_U2)PU2poPU3p(PUhat3-P_U3)PU3p in tensor regression with sample splitting}\\
& \left\|\mcA \times_j \mcP_{U_{j\perp}}\left(\mcP_{\whU_j^{(\mathrm{\RN{1}})}} - \mcP_{U_j}\right)\mcP_{U_{j\perp}} \times_{j+1} \mcP_{U_{j+1\perp}}\left(\mcP_{\whU_{j+1}^{(\mathrm{\RN{1}})}} - \mcP_{U_{j+1}}\right)\mcP_{U_{{j+1}\perp}} \times_{j+2} \mcP_{U_{j+2\perp}}\left(\mcP_{\whU_{j+2}^{(\mathrm{\RN{1}})}} - \mcP_{U_{j+2}}\right)\mcP_{U_{j+2\perp}}  \right\|_{\mathrm{F}} \notag \\
\lesssim & \left\|\mcA\right\|_{\mathrm{F}} \cdot \left(\frac{\sigma_{\xi}^6}{\ulambda^6 \sigma^6} \cdot \frac{\oor^{3/2}\op^{3/2}\log(\op)^{3/2}}{n^3}+\frac{\sigma_{\xi}^6}{\ulambda^6 \sigma^6} \cdot \Delta \cdot \frac{\oR^{3/2} \op^{3/2}\log (\op)^{3/2}}{n^6}\right) \label{eq: high-prob upper bound of AoPU1p(PUhat1-P_U1)PU1oPU2p(PUhat2-P_U2)PU2poPU3p(PUhat3-P_U3)PU3pp in tensor regression with sample splitting}
\end{align}
for any $j=1,2,3$.

\end{proposition}
The proof of Proposition~\ref{prop: high-probability upper bound of AoPU1pEhat1P_U1oPU2pEhat2oPU3pEhat3 with sample splitting} is similar to that of Proposition~\ref{prop: high-probability upper bound of AoPU1p(PUhat1-P_U1)oPU2p(PUhat2-P_U2)oPU3p(PUhat3-P_U3) in tensor regression without sample splitting}, and is thus omitted.

\begin{proposition}\label{prop: high-probability upper bound of AoPU1oPU2p(PUhat2-P_U2)oPU3p(PUhat3-P_U3) with sample splitting}

Under the same setting of Theorem~\ref{thm: main theorem in tensor regression with sample splitting}, with probability at least $1-\exp(-cn) - \frac{1}{p^C} - \P\left(\mcE_\Delta\right) - \P\left(\mcE_{U}^{\text{reg}}\right)$, where $c$ and $C$ are two universal constants ,it holds that
\begin{align}
& \left\|\mcA \times_j \mcP_{U_j} \times_{j+1} \mcP_{U_{j+1\perp}}\left(\mcP_{\whU_{j+1}^{(\mathrm{\RN{1}})}} - \mcP_{U_{j+1}}\right)\mcP_{U_{j+1}} \times_{j+2} \mcP_{U_{j+2\perp}}\left(\mcP_{\whU_{j+2}^{(\mathrm{\RN{1}})}} - \mcP_{U_{j+2}}\right)\mcP_{U_{j+2}}\right\|_{\mathrm{F}} \notag \\
\lesssim & \left\|\mcA \times_j U_j\right\|_{\mathrm{F}} \cdot \left(\frac{\sigma_{\xi}^2}{\ulambda^2 \sigma^2} \cdot \frac{\oor\log(\op)}{n} + \frac{\sigma_{\xi}^2}{\ulambda^2 \sigma^2} \cdot \Delta \cdot \frac{\oR \log (\op)}{n}\right) \label{eq: high-prob upper bound of AoPU1oPU2p(PUhat2-P_U2)PU2oPU3p(PUhat3-P_U3)PU3 in tensor regression with sample splitting} \\
& \left\|\mcA \times_j \mcP_{U_j} \times_{j+1} \mcP_{U_{j+1\perp}}\left(\mcP_{\whU_{j+1}^{(\mathrm{\RN{1}})}} - \mcP_{U_{j+1}}\right)\mcP_{U_{j+1\perp}} \times_{j+2} \mcP_{U_{j+2\perp}}\left(\mcP_{\whU_{j+2}^{(\mathrm{\RN{1}})}} - \mcP_{U_{j+2}}\right)\mcP_{U_{j+2}}\right\|_{\mathrm{F}} \notag \\
\lesssim & \left\|\mcA \times_j U_j\right\|_{\mathrm{F}} \cdot \left(\frac{\sigma_{\xi}^3}{\ulambda^3 \sigma^3} \cdot \frac{\oor\op^{1/2}\log(\op)}{n^{3/2}} + \frac{\sigma_{\xi}^3}{\ulambda^3 \sigma^3} \cdot \Delta \cdot \frac{\oR \op^{1/2} \log (\op)}{n^{3/2}} \right) \label{eq: high-prob upper bound of AoPU1oPU2p(PUhat2-P_U2)PU2poPU3p(PUhat3-P_U3)PU3 in tensor regression with sample splitting}\\
& \left\|\mcA \times_j \mcP_{U_j} \times_{j+1} \mcP_{U_{j+1\perp}}\left(\mcP_{\whU_{j+1}^{(\mathrm{\RN{1}})}} - \mcP_{U_{j+1}}\right)\mcP_{U_{j+1\perp}} \times_{j+2} \mcP_{U_{j+2\perp}}\left(\mcP_{\whU_{j+2}^{(\mathrm{\RN{1}})}} - \mcP_{U_{j+2}}\right)\mcP_{U_{j+2\perp}}\right\|_{\mathrm{F}} \notag \\
\lesssim & \left\|\mcA \times_j U_j\right\|_{\mathrm{F}} \cdot \left(\frac{\sigma_{\xi}^4}{\ulambda^4 \sigma^4} \cdot \frac{\oor\op\log (\op)}{n^2} + \frac{\sigma_{\xi}^4}{\ulambda^4 \sigma^4} \cdot \Delta \cdot \frac{\oR \op \log (\op)}{n^2}\right).\label{eq: high-prob upper bound of AoPU1oPU2p(PUhat2-P_U2)PU2poPU3p(PUhat3-P_U3)PU3p in tensor regression with sample splitting}
\end{align}
\end{proposition}
The proof of Proposition~\ref{prop: high-probability upper bound of AoPU1oPU2p(PUhat2-P_U2)oPU3p(PUhat3-P_U3) with sample splitting} is similar to that of Proposition~\ref{prop: high-probability upper bound of AoPU1oPU2p(PUhat2-P_U2)oPU3p(PUhat3-P_U3) without sample splitting}, and is thus omitted.

\begin{proposition} \label{prop: high-prob upper bound of perturbation after projection in tensor regression with sample splitting}

Under the same setting of Theorem~\ref{thm: main theorem in tensor regression with sample splitting}, let $V_j \in \mathbb{R}^{p_j\times R_j}$ be a fixed matrix satisfying $\left\|V_j\right\|=1$ for any $j=1,2,3$. Then with probability at least $1-\exp(-cn) - \frac{1}{p^C} - \P\left(\mcE_\Delta\right) - \P\left(\mcE_{U}^{\text{reg}}\right)$, where $c$ and $C$ are two universal constants ,it holds that
\begin{align}
& \left\|V_j^{\top}\mcP_{U_{j\perp}}\left(\mcP_{\whU_j^{(\mathrm{\RN{1}})}} - \mcP_{U_j}\right)U_j\right\| \lesssim \frac{\sigma_\xi}{\ulambda\sigma} \cdot \sqrt{\frac{\oR\log (\op)}{n}}, \label{eq: high-prob upper bound of V1tPU1p(PUhat1-PU1)U1 in tensor regression with sample splitting} \\
& \left\|V_j^{\top}\mcP_{U_{j\perp}}\left(\mcP_{\whU_j^{(\mathrm{\RN{1}})}} - \mcP_{U_1}\right)U_{j\perp}\right\| \lesssim \frac{\sigma_\xi^2}{\ulambda^2\sigma^2}\cdot 
\frac{\sqrt{\oR\log\left(\op\right)}\cdot \sqrt{\op}}{n}  \label{eq: high-prob upper bound of V1tPU1p(PUhat1-PU1)U1p in tensor regression with sample splitting}
\end{align}

Furthermore, 
\begin{align}
& \left\|V_j^{\top}\mcP_{U_{j\perp}}\left(\mcP_{\whU_j^{(\mathrm{\RN{1}})}} - \mcP_{U_j}\right)\right\| \lesssim \frac{\sigma_\xi}{\ulambda\sigma} \cdot \sqrt{\frac{\oR\log (\op)}{n}}. \label{eq: high-prob upper bound of V1tPU1p(PUhat1-PU1) in tensor regression with sample splitting}
\end{align}

\end{proposition}
The proof of Proposition~\ref{prop: high-prob upper bound of perturbation after projection in tensor regression with sample splitting} is similar to that of Proposition~\ref{prop: high-prob upper bound of perturbation after projection in tensor regression without sample splitting}, and is thus omitted.

\subsection{Upper Bound of First-Order Perturbation Terms}

\begin{proposition}\label{prop: high-prob upper bound of first-order perturbation terms after projection in tensor regression with sample splitting}

Under the same setting of Theorem~\ref{thm: main theorem in tensor regression with sample splitting}, let $V_j \in \mathbb{R}^{p_j\times R_j}$ be a fixed matrix satisfying $\left\|V_j\right\|=1$. $\whE_j$ is defined as in \eqref{eq: definition of whEj in tensor regression with sample splitting} for any $j=1,2,3$. Then with probability at least $1-\exp(-cn) - \frac{1}{p^C} - \P\left(\mcE_\Delta\right) - \P\left(\mcE_{U}^{\text{reg}}\right)$, where $c$ and $C$ are two universal constants ,it holds that
\begin{align}
& \left\|\mcP_j^{-\frac{1}{2}} \whE_j^{(\mathrm{\RN{1}})} \mcP_j^{-\frac{1}{2}}\right\| \lesssim \frac{\sigma_\xi}{\ulambda\sigma}\sqrt{\frac{\oor\log(\op)}{n}} + \frac{\sigma_\xi^2}{\ulambda^2\sigma^2} \cdot \frac{\op}{n}, \label{eq: high-prob upper bound of P1(-1/2)Ehat1P1(-1/2) in tensor regression with sample splitting} \\
& \left\|V_j^{\top} \mcP_{U_{j \perp}} \whE_j^{(\mathrm{\RN{1}})} \mcP_j^{-\frac{1}{2}}\right\| = \left\|V_j^{\top}\mcP_{U_{j\perp}}\whE_j^{(\mathrm{\RN{1}})} U_j\left(G_jG_j^{\top}\right)^{-\frac{1}{2}}U_j^{\top}\right\| \lesssim \frac{\sigma_\xi}{\sigma}\cdot \sqrt{\frac{\oR\log(\op)}{n}}, \label{eq: high-prob upper bound of V1tP1(0)Ehat1P1(-1/2) in tensor regression with sample splitting} \\
&\left\|V_j^{\top} \mcP_{U_{j\perp}} \whE_j^{(\mathrm{\RN{1}})} U_{j\perp}\right\| \lesssim \frac{\sigma_\xi}{\sigma}\cdot \sqrt{\frac{\oR\log(\op)}{n}}.\label{eq: high-prob upper bound of V1tP1(0)Ehat1P1(0) in tensor regression with sample splitting}
\end{align}
\end{proposition}
The proof of Proposition~\ref{prop: high-prob upper bound of first-order perturbation terms after projection in tensor regression with sample splitting} is similar to that of Proposition~\ref{prop: high-prob upper bound of first-order perturbation terms after projection in tensor regression without sample splitting}, and is thus omitted.

\subsection{Upper Bound of Higher-Order Perturbation Terms}

\begin{lemma} \label{lemma: high-prob upper bound of the higher-order perturbation terms in tensor regression with sample splitting}

Under the same setting of Theorem~\ref{thm: main theorem in tensor regression with sample splitting}, let $V_j \in \mathbb{R}^{p_j\times R_j}$ be a fixed matrix satisfying $\left\|V_j\right\|=1$. $\whE_j^{(\mathrm{\RN{1}})}$ is defined as in \eqref{eq: definition of whEj in tensor regression with sample splitting}, and $S_{G_j, k_j}(\whE_j^{(\mathrm{\RN{1}})})$ is defined as in \eqref{eq: definition of Sj(whEj) in tensor regression with sample splitting}, for any $j=1,2,3$. Then with probability at least $1-\exp(-cn) - \frac{1}{p^C} - \P\left(\mcE_\Delta\right) - \P\left(\mcE_{U}^{\text{reg}}\right)$, where $c$ and $C$ are two universal constants ,it holds that
\begin{align}
& \left\|\mcP_{U_j} \sum_{k_j=2}^{+\infty} S_{G_j, k_j}\left(\whE_j^{(\mathrm{\RN{1}})}\right) \mcP_{U_{j \perp}} V_j\right\| \lesssim \frac{\sigma_\xi^2}{\ulambda^2\sigma^2} \cdot \frac{\sqrt{\oR\oor} \log(\op)}{n} + \frac{\sigma_\xi^3}{\ulambda^3\sigma^3} \cdot \frac{\op \sqrt{\oR \log(\op)}}{n^{3/2}}, \label{eq: high-prob upper bound of PUorder2PUpPV in tensor regression with sample splitting}\\
& \left\|\mcP_{U_j} \sum_{k_j=3}^{+\infty} S_{G_j, k_j}\left(\whE_j^{(\mathrm{\RN{1}})}\right) \mcP_{U_{j\perp}}V_j\right\| \lesssim \frac{\sigma_{\xi}^3}{\ulambda^3 \sigma^3} \cdot \frac{\op \sqrt{\oR \log (\op)}}{n^{3 / 2}}, \label{eq: high-prob upper bound of PUorder3PUpPV in tensor regression with sample splitting} \\
& \left\|\mcP_{U_{j\perp}} \sum_{k_j=2}^{+\infty} S_{G_j, k_j}\left(\whE_j^{(\mathrm{\RN{1}})}\right) \mcP_{U_{j \perp}} V_j\right\| \lesssim \frac{\sigma_\xi^2}{\ulambda^2\sigma^2} \cdot \frac{\sqrt{\oR\op \log (\op)}}{n}, \label{eq: high-prob upper bound of PUporder2PUpPV in tensor regression with sample splitting}\\
& \left\|\mcP_{U_{j \perp}} \sum_{k_j=3}^{+\infty} S_{G_j, k_j}\left(\whE_j^{(\mathrm{\RN{1}})}\right) \mcP_{U_{j \perp}} V_j\right\| \lesssim \frac{\sigma_\xi^3}{\ulambda^3\sigma^3} \cdot \frac{\sqrt{\oR\op\oor}\log(\op)}{n^{3/2}} + \frac{\sigma_\xi^4}{\ulambda^4\sigma^4} \cdot \frac{\op^{3/2}\sqrt{\oR\log(\op)}}{n^2} \label{eq: high-prob upper bound of PUporder3PUpPV in tensor regression with sample splitting}.
\end{align}

\end{lemma}
The proof of Lemma~\ref{lemma: high-prob upper bound of the higher-order perturbation terms in tensor regression with sample splitting} is similar to that of Lemma~\ref{lemma: high-prob upper bound of the higher-order perturbation terms in tensor regression without sample splitting}, and is thus omitted.

\begin{lemma}\label{lemma: high-prob upper bound of negligible terms in the first-order term in tensor regression without sample splitting}

Under the same setting of Theorem~\ref{thm: main theorem in tensor regression with sample splitting}, let $V_j \in \mathbb{R}^{p_j\times R_j}$ be a fixed matrix satisfying $\left\|V_j\right\|=1$. $\whE_j^{(\mathrm{\RN{1}})}$ is defined as in \eqref{eq: definition of whEj in tensor regression with sample splitting} for any $j=1,2,3$. Then with probability at least $1-\exp(-cn) - \frac{1}{p^C} - \P\left(\mcE_\Delta\right) - \P\left(\mcE_{U}^{\text{reg}}\right)$, where $c$ and $C$ are two universal constants ,it holds that
\begin{align}
& \left\|V_j^{\top}  \mcP_{U_{j \perp}} \whE_j^{(\mathrm{\RN{1}})} U_{j \perp}-V_j V_j^{\top} \mcP_{U_{j \perp}} \whZ_j^{(\mathrm{\RN{1}})}\left(\mcP_{U_{j+2}} \otimes \mcP_{U_{j+1}}\right) \whZ_j^{(\mathrm{\RN{1}}),\top} U_{j \perp}\right\| \leq \frac{\sigma_\xi^3}{\ulambda\sigma^3} \cdot \frac{\op \sqrt{\oR \log (\op)}}{n^{3/2}}, \label{eq: high-prob upper bound of V1V1tP1(0)Ehat1PU1p - V1V1tP1(0)Zhat1(PU3oPU2)Zhat1tU1p in tensor regression with sample splitting} \\
& \left\|V_j  U_{j \perp}^{\top} \whE_j^{(\mathrm{\RN{1}})} P_j^{-1} -V_j V_j^{\top} U_{j \perp}^{\top} \whZ_j^{(\mathrm{\RN{1}})} \left(U_{j+2} \otimes U_{j+1}\right) G_j^{\top}\left(G_j G_j^{\top}\right)^{-1} U_j^{\top} \right\| \leq \frac{\sigma_\xi^2}{\ulambda^2\sigma^2}\cdot \frac{\sqrt{\oR \op \log (\op)}}{n}  \label{eq: high-prob upper bound of V1V1tP1(0)Ehat1P1^(-1) - V1V1tP1(0)Zhat1(U3oU2)G1t(G1G1t)^(-1)U1t in tensor regression with sample splitting}.
\end{align}

\end{lemma}

\begin{proof}

By symmetry, it suffices to consider upper bounds of:
\begin{align*}
\mathrm{\RN{1}} = & \left\| V_1 V_1^{\top} \mcP_{U_{1 \perp}} \whE_1^{(\mathrm{\RN{1}})} U_{1 \perp}-V_1 V_1^{\top} \mcP_{U_{1 \perp}} \whZ_1^{(\mathrm{\RN{1}})}\left(\mcP_{U_3} \otimes \mcP_{U_2}\right) \whZ_1^{(\mathrm{\RN{1}}), \top} U_{1 \perp} \right\| , \\
\mathrm{\RN{2}} = & \left\|V_1 V_1^{\top} \mcP_{U_{1 \perp}} \whE_1^{(\mathrm{\RN{1}})} P_1^{-1}-V_1 V_1^{\top} \mcP_{U_{1 \perp}} \whZ_1^{(\mathrm{\RN{1}})}\left(U_3 \otimes U_2\right) G_1^{\top}\left(G_1 G_1^{\top}\right)^{-1} U_1^{\top}\right\|.
\end{align*}

{\bf Step 1: $\left\|V_1 V_1^{\top} \mcP_{U_{1 \perp}} \whE_1^{(\mathrm{\RN{1}})} U_{1 \perp}-V_1 V_1^{\top} \mcP_{U_{1 \perp}} \whZ_1^{(\mathrm{\RN{1}})}\left(\mcP_{U_3} \otimes \mcP_{U_2}\right) \whZ_1^{(\mathrm{\RN{1}}), \top} U_{1 \perp}\right\|$}

Consider
\begin{align*}
\mathrm{\RN{2}} 
\leq & \underbrace{\left\|V_1^{\top}\mcP_{U_{1 \perp}} \whZ_1^{(\mathrm{\RN{1}})}\left[\left(\mcP_{\whU_3^{(0), (\mathrm{\RN{2}})}} - \mcP_{U_3}\right) \otimes \mcP_{U_2}\right] \whZ_1^{(\mathrm{\RN{1}}), \top} U_{1 \perp}\right\|}_{\mathrm{\RN{2}}.\mathrm{\RN{1}}} 
+  \underbrace{\left\|V_1^{\top}\mcP_{U_{1 \perp}} \whZ_1^{(\mathrm{\RN{1}})}\left[\mcP_{U_3} \otimes \left(\mcP_{\whU_2^{(0), (\mathrm{\RN{2}})}}-\mcP_{U_2}\right)\right] \whZ_1^{(\mathrm{\RN{1}}), \top} U_{1 \perp}\right\|}_{\mathrm{\RN{2}}.\mathrm{\RN{2}}} \\
&+  \underbrace{\left\|V_1^{\top}\mcP_{U_{1 \perp}} \whZ_1^{(\mathrm{\RN{1}})}\left[\left(\mcP_{\whU_3^{(0), (\mathrm{\RN{2}})}} -\mcP_{U_3}\right) \otimes \left(\mcP_{\whU_2^{(0), (\mathrm{\RN{2}})}}-\mcP_{U_2}\right)\right] \whZ_1^{(\mathrm{\RN{1}}), \top} U_{1 \perp}\right\|}_{\mathrm{\RN{2}}.\mathrm{\RN{3}}} \\
\lesssim & \frac{\sigma_\xi^3}{\ulambda\sigma^3} \cdot \frac{\op \sqrt{\oR \log (\op)}}{n^{3/2}} .
\end{align*}
Here, we used the independence between and $\whZ_1^{(\mathrm{\RN{1}})}$ and $\mcP_{\whU_j^{(0), (\mathrm{\RN{2}})}}$.

{\bf Step 2: $\left\|V_1^{\top}\mcP_{U_{1 \perp}} \whE_1^{(\mathrm{\RN{1}})} \mcP_1^{-1} - V_1V_1^{\top}\mcP_{U_{1 \perp}} \whZ_1^{(\mathrm{\RN{1}})}\left(U_3 \otimes U_2\right)\right\|$}

Furthermore, since $\widehat{Z}_1^{(1)}$ and $\mcP_{\widehat{U}_j^{(0), \mathrm{\RN{2}}}}$ are mutually independent, it follows that
\begin{align*}
\mathrm{\RN{2}}  
\leq & \underbrace{\left\|V_1^{\top}\mcP_{U_{1 \perp}} \whZ_1^{(\mathrm{\RN{1}})}\left[U_3 \otimes \left(\mcP_{\whU_2^{(0), (\mathrm{\RN{2}})}}-\mcP_{U_2}\right) U_2\right] G_1^{\top}\left(G_1G_1^{\top}\right)^{-1}U_1^{\top}\right\|}_{\mathrm{\RN{4}}.\mathrm{\RN{1}}} \\
+ & \underbrace{\left\|V_1^{\top}\mcP_{U_{1 \perp}} \whZ_1^{(\mathrm{\RN{1}})}\left[\left(\mcP_{\whU_3^{(0), (\mathrm{\RN{2}})}}-\mcP_{U_3}\right) U_3 \otimes U_2\right] G_1^{\top}\left(G_1G_1^{\top}\right)^{-1}U_1^{\top}\right\|}_{\mathrm{\RN{4}}.\mathrm{\RN{2}}} \\
+ & \underbrace{\left\|V_1^{\top}\mcP_{U_{1 \perp}} \whZ_1^{(\mathrm{\RN{1}})}\left[\left(\mcP_{\whU_3^{(0), (\mathrm{\RN{2}})}}-\mcP_{U_3}\right) U_3 \otimes \left(\mcP_{\whU_2^{(0), (\mathrm{\RN{2}})}}-\mcP_{U_2}\right) U_2\right] G_1^{\top}\left(G_1G_1^{\top}\right)^{-1}U_1^{\top}\right\|}_{\mathrm{\RN{4}}.\mathrm{\RN{3}}} \\
+ & \underbrace{\left\|V_1^{\top}\mcP_{U_{1 \perp}} \whZ_1^{(\mathrm{\RN{1}})}\left(\mcP_{U_3} \otimes \mcP_{U_2}\right) \whZ_1^{(\mathrm{\RN{1}}), \top} U_1\left(G_1G_1^{\top}\right)^{-1}U_1^{\top}\right\|}_{\mathrm{\RN{4}}.\mathrm{\RN{4}}} \\
+ & \underbrace{\left\|V_1^{\top}\mcP_{U_{1 \perp}} \whZ_1^{(\mathrm{\RN{1}})}\left[\mcP_{U_3} \otimes \left(\mcP_{\whU_2^{(0), (\mathrm{\RN{2}})}}-\mcP_{U_2}\right)\right] \whZ_1^{(\mathrm{\RN{1}}), \top} U_1\left(G_1G_1^{\top}\right)^{-1}U_1^{\top}\right\|}_{\mathrm{\RN{4}}.\mathrm{\RN{5}}} \\
+ & \underbrace{\left\|V_1^{\top}\mcP_{U_{1 \perp}} \whZ_1^{(\mathrm{\RN{1}})}\left[\left(\mcP_{\whU_3^{(0), (\mathrm{\RN{2}})}}-\mcP_{U_3}\right) \otimes \mcP_{U_2}\right] \whZ_1^{(\mathrm{\RN{1}}), \top} U_1\left(G_1G_1^{\top}\right)^{-1}U_1^{\top}\right\|}_{\mathrm{\RN{4}}.\mathrm{\RN{6}}} \\
+ & \underbrace{\left\|V_1^{\top}\mcP_{U_{1 \perp}} \whZ_1^{(\mathrm{\RN{1}})}\left[\left(\mcP_{\whU_3^{(0), (\mathrm{\RN{2}})}}-\mcP_{U_3}\right) \otimes \left(\mcP_{\whU_2^{(0), (\mathrm{\RN{2}})}}-\mcP_{U_2}\right)\right] \whZ_1^{(\mathrm{\RN{1}}), \top} U_1\left(G_1G_1^{\top}\right)^{-1}U_1^{\top}\right\|}_{\mathrm{\RN{4}}.\mathrm{\RN{7}}}\\
\lesssim & \frac{\sigma_\xi^3}{\ulambda^2\sigma^3}\cdot \frac{\sqrt{\oR \op \log (\op)}}{n} .
\end{align*}

\end{proof}

\subsection{Upper Bound of Leading Terms in the Spectral Representation}

\begin{lemma}\label{lemma: high-prob upper bound of leading terms in the spectral representation across three modes in tensor regression with sample splitting}

Under the same setting of Theorem~\ref{thm: main theorem in tensor regression with sample splitting}, with probability at least $1-\exp(-cn) - \frac{1}{p^C} - \P\left(\mcE_\Delta\right) - \P\left(\mcE_{U}^{\text{reg}}\right)$, where $c$ and $C$ are two universal constants ,it holds that
\begin{align}
& \left\|\mcA \times_j \mcP_{U_{j \perp}} \whE_j^{(\mathrm{\RN{1}})} \mcP_j^{-1} \times_{j+1} \mcP_{U_{{j+1} \perp}} \whE_{j+1}^{(\mathrm{\RN{1}})} \mcP_{j+1}^{-1} \times_{j+2} \mcP_{U_{{j+2} \perp}} \whE_{j+2}^{(\mathrm{\RN{1}})} \mcP_{j+2}^{-1}\right\|_{\mathrm{F}} \notag \\
\lesssim & \left\| \mcA \right\|_{\mathrm{F}} \cdot \left(\frac{\sigma_\xi^3}{\ulambda^3\sigma^3}\cdot \frac{\oor^3 \log (\op)}{n^{3 / 2}} + \frac{\sigma_\xi^3}{\ulambda^3\sigma^3}\cdot \Delta \cdot \frac{\oR^{3/2}\log(\op)^{3/2}}{n^{3/2}} \right) \label{eq: high-prob upper bound of AoPU1pEhat1P1(-1)oPU2pEhat2P2(-1)oPU3pEhat3P3(-1) in tensor regression with sample splitting}
\end{align}
for any $j=1,2,3$. Here, $\whE_j^{(\mathrm{\RN{1}})}$ is defined as in \eqref{eq: definition of whEj in tensor regression with sample splitting}.

\end{lemma}
The proof of Lemma~\ref{lemma: high-prob upper bound of leading terms in the spectral representation across three modes in tensor regression with sample splitting} is similar to that of Lemma~\ref{lemma: high-prob upper bound of leading terms in the spectral representation across three modes in tensor regression without sample splitting}, and is thus omitted.

\begin{lemma}\label{lemma: high-prob upper bound of leading terms in the spectral representation across two modes in tensor regression with sample splitting}
Under the same setting of Theorem~\ref{thm: main theorem in tensor regression with sample splitting}, with probability at least $1-\exp(-cn) - \frac{1}{p^C} - \P\left(\mcE_\Delta\right) - \P\left(\mcE_{U}^{\text{reg}}\right)$, where $c$ and $C$ are two universal constants ,it holds that
\begin{align}
& \left\|\mcA \times_j U_j \times_{j+1} \mcP_{U_{j+1\perp}}\whE_{j+1}^{(\mathrm{\RN{1}})}\mcP_{j+1}^{-1} \times_{j+2} \mcP_{U_{j+2\perp}}\whE_{j+2}^{(\mathrm{\RN{1}})}\mcP_{j+2}^{-1}\right\|_{\mathrm{F}} \notag \\
\lesssim & \left\|\mcA \times_j U_j\right\|_{\mathrm{F}} \cdot \left(\frac{\sigma_{\xi}^2}{\ulambda^2 \sigma^2} \cdot \frac{\oor\log(\op)}{n} + \frac{\sigma_\xi^2}{\ulambda^2\sigma^2} \cdot \Delta \cdot \frac{\oR\log(\op)}{n}\right) \label{eq: high-prob upper bound of AoPU1oPU2pEhat2P2(-1)oPU3pEhat3P3(-1) in tensor regression with sample splitting}
\end{align}
for any $j=1,2,3$. Here, $\whE_j^{(\mathrm{\RN{1}})}$ is defined as in \eqref{eq: definition of whEj in tensor regression with sample splitting}.
\end{lemma}
The proof of Lemma~\ref{lemma: high-prob upper bound of leading terms in the spectral representation across two modes in tensor regression with sample splitting} is similar to that of Lemma~\ref{lemma: high-prob upper bound of leading terms in the spectral representation across two modes in tensor regression without sample splitting}, and is thus omitted.

\section{Concentration Inequalities for Tensor regression}

The following two lemmas (Lemma~\ref{lemma: high-prob upper bound of tr(BZ(2)tCZ(1)) in tensor regression with sample splitting} and Lemma~\ref{lemma: high-prob upper bound of tr(BZ(2)tCZ(2)) in tensor regression with sample splitting}) characterize the concentration bound for the perturbation term introduced by the bias-correction procedure. 

\begin{lemma}\label{lemma: high-prob upper bound of tr(BZ(2)tCZ(1)) in tensor regression with sample splitting}

Let $X \in \mathbb{R}^{p \times d}$ be a random matrix with mean-zero $\sigma$ sub-Gaussian entries, and let $\left\{ X_i\right\}_{i=1}^{n}$ be $n$ i.i.d. copies of $X$ and $\left\{\xi_i\right\}_{i=1}^{n}$ are i.i.d. mean-zero $\sigma_\xi$-sub-Gaussian random variables. 

Assume that $A \in \mathbb{R}^{p \times d}$, $B\in \mathbb{R}^{d \times d}$ and $C \in \mathbb{R}^{p \times d}$ are three fixed matrices. Then it follows that
\begin{align}
& \mathbb{P}\Bigg(\Bigg|\operatorname{tr}\Bigg[\frac{1}{n\sigma^2} \Big(\sum_{i=1}^n\left\langle A,  X_i\right\rangle  BX_i^{\top} C - \sigma^2 BAC \Big) \cdot \frac{1}{n\sigma^2} \sum_{j=1}^n \xi_j X_j \Bigg]\Bigg| \geq \frac{K_3\sqrt{\log(\op)}}{n^{3/2}\sigma^2} \left\|A\right\|_{\mathrm{F}}\cdot \left\|B\right\|_{\mathrm{F}} \cdot \left\|C\right\|_{\mathrm{F}} 
 \notag \\
& \qquad + \frac{\sqrt{\log(\op)}}{n} \cdot \left\|A\right\|_{\mathrm{F}}\cdot \left\|B\right\| \cdot \left\|C\right\| \cdot t^{1/2} + \frac{1}{n^{3/2}} \cdot  \left\|A\right\|_{\mathrm{F}}\cdot \left|\operatorname{tr(B)}\operatorname{tr(C)}\right| \cdot  t^{1/2} + \frac{1}{n} \cdot \left\|A\right\|_{\mathrm{F}} \cdot \left\|B\right\|_{\mathrm{F}} \cdot \left\|C\right\|_{\mathrm{F}} \cdot t^{3/2}\Bigg) \notag \\
\leq & \exp\left(-ct\right) + \frac{1}{\op^{c}} + \exp(-cn) \label{eq: high-prob upper bound of tr(BZ(2)tCZ(1)) in tensor regression with sample splitting}
\end{align}
where $c>0$ is a constant.
\end{lemma}

\begin{proof}

Consider the following decomposition:
\begin{align*}
& \operatorname{tr}\left[\frac{1}{n\sigma^2} \left(\sum_{i=1}^n\left\langle A,  X_i\right\rangle  BX_i^{\top} C - \sigma^2 BAC\right) \cdot \frac{1}{n\sigma^2} \sum_{j=1}^n \xi_j X_j \right] \\
= & \frac{1}{n^2\sigma^4} \sum_{i=1}^n \xi_i a^{\top}x_i x_i^{\top}\left(C \otimes B\right)x_i + \frac{1}{n^2\sigma^4} \sum_{i=1}^n \sum_{j=1, j\neq 1}^n \xi_j  a^{\top}x_i x_i^{\top}\left(C \otimes B\right)x_j \frac{1}{n\sigma^2}\sum_{i=1}^{n} \operatorname{tr}\xi_i\left[BACX_i\right],
\end{align*}
where $x_i=\operatorname{Vec}\left(X_i\right)$ , $a=\operatorname{Vec}\left(A\right)$, and  $D = C^{\top}\otimes B$.

Here, first we have
$$
\mathbb{P}\left(\left|\sum_{j=1}^n \xi_j \operatorname{tr}\left[BACX_j\right]\right| \geq t \right) \leq \exp\left(-\frac{c t^2}{\sigma^2\sum_{j=1}^n\xi_j^2\left\|BAC\right\|_{\ell_{\infty}}^2} \right).
$$

It implies that
$$
\mathbb{P}\left(\frac{1}{n\sigma^2}\left|\sum_{j=1}^n \xi_j \operatorname{tr}\left[BACX_j\right]\right| \geq \frac{\sigma_\xi}{\sigma}\cdot \frac{\left\|BAC\right\|_{\ell_{\infty}}\sqrt{\log(\op)}}{\sqrt{n}} \right) \leq  \frac{1}{\op^c}.
$$

Then consider the following polynomial of $\left\{x_{i,k}\right\}_{i=1,k=1}^{p,d}$:
\begin{align*}
f(x)= & \sum_{i=1}^n \xi_i a^{\top} x_i x_i^{\top}\left(C \otimes B\right)x_i = \sum_{i=1}^n \xi_i \sum_{k=1}^{p_1p_2p_3} \sum_{l=1}^{p_1p_2p_3} \sum_{m=1}^{p_1p_2p_3} a_k x_{i,k} x_{i,l} D_{l,m} x_{i,m}.
\end{align*}

It implies that
\begin{align*}
& \mathbb{E}\left[\operatorname{tr}\left[\frac{1}{n\sigma^2} \sum_{i=1}^n \left\langle A,  X_i\right\rangle B X_i^{\top} C \frac{1}{n\sigma^2} \sum_{j=1}^n \xi_j X_j \right] \Bigg| \left\{\xi_{i}\right\}_{i=1}^n \right] = \frac{K_3}{n^2\sigma^2} \sum_{i=1}^n \xi_i \sum_{k=1}^{p_1p_2p_3} a_k D_{k,k}.
\end{align*}

Then consider
\begin{align*}
f(x)
= & \underbrace{\sum_{k=1}^{p_1p_2p_3} a_k D_{k,k} x_{i,k}^3}_{k=l=m} +  \underbrace{\sum_{k=1}^{p_1p_2p_3} \sum_{l\neq k}^{p_1p_2p_3} a_k D_{l,k} x_{i,k}^2x_{i,l}}_{k=m \neq l} + \underbrace{\sum_{k=1}^{p_1p_2p_3} \sum_{m \neq k}^{p_1p_2p_3} a_k x_{i,k}^2 D_{k,m} x_{i,m}}_{k=l \neq m} \\
& +  \underbrace{\sum_{k=1}^{p_1p_2p_3} \sum_{l \neq k}^{p_1p_2p_3} \sum_{m \neq k}^{p_1p_2p_3} \sum_{m=1}^{p_1p_2p_3} a_k x_{i,k} x_{i,l} D_{l,m} x_{i,m}}_{k \neq l, k \neq m}.
\end{align*}  

It follows that
\begin{align*}
\frac{\partial}{\partial x_{i,k}} f(x) 
= & 3a_k D_{k,k} x_{i,k}^2 + 2\sum_{k=1}^{p_1p_2p_3} \sum_{l\neq k}^{p_1p_2p_3} a_k D_{l,k} x_{i,k}x_{i,l} + 2 \sum_{k=1}^{p_1p_2p_3} \sum_{m \neq k}^{p_1p_2p_3} a_k x_{i,k} D_{k,m} x_{i,m} \\
& +  \sum_{l \neq k}^{p_1p_2p_3} \sum_{m \neq k}^{p_1p_2p_3} \sum_{m=1}^{p_1p_2p_3} a_k x_{i,l} D_{l,m} x_{i,m}.
\end{align*}  

It further implies that
$
\mathbb{E}\left[\frac{\partial}{\partial x_{i,k}} f(x)\right] = 3\sigma^2 \cdot a_k D_{k,k}  + \sigma^2 \cdot \sum_{l \neq k}^{p_1p_2p_3} a_k D_{l,l}.
$

Therefore, we have
$$
\left\|\mathbb{E}\nabla f(x)\right\|_{\text{HS}}^2 \leq C\sigma^4  \sum_{k=1}^{p_1p_2p_3}a_k^2 \left(\sum_{l =1}^{p_1p_2p_3} D_{l,l}\right)^2 \leq C\sigma^4 \left\|a\right\|_{\ell_2}^2 \cdot \left[\operatorname{tr}\left(D\right)\right]^2
$$

Then, consider the second derivative, we have
$$
\frac{\partial^2}{\partial x_{i,k}^2} f(x) = 6a_kD_{k,k}x_{i,k} + 2 \sum_{l\neq k}^{p_1p_2p_3} a_k D_{l,k} x_{i,l} .
$$

It implies that
$
\mathbb{E}\left[\frac{\partial^2}{\partial x_{i,k}^2} f(x)\right] = 0.
$

Furthermore, we have
\begin{align*}
\frac{\partial^2}{\partial x_{i,k} x_{i,l}} f(x) 
= & 2 a_k D_{l,k} x_{i,k} + 2 a_k x_{i,k} D_{k,l} + 2a_k x_{i,l}D_{l,l} + \sum_{m \neq k, m\neq l}^{p_1p_2p_3} a_k D_{l,m} x_{i,m} .
\end{align*}  

It implies that
$
\mathbb{E}\left(\frac{\partial^2}{\partial x_{i,k} x_{i,l}} f(x)\right) = 0
$
Therefore, we have
$
\left\|\mathbb{E}\nabla^2 f(x)\right\|_{\text{HS}}^2 =0.
$

Finally, consider the following third-order partial derivatives
\begin{align*}
& \frac{\partial^3}{\partial x_{i,k}^3} f(x) = 6a_kD_{k,k}, \quad \frac{\partial^2}{\partial x_{i,k}^2 x_{i,l}} f(x) 
= 2 a_k D_{l,k} + 2 a_k D_{k,l}, \quad \frac{\partial^2}{\partial x_{i,k} x_{i,l} x_{i,m}} f(x) 
= a_k D_{l,m} .
\end{align*}  

Therefore, we have
\begin{align*}
\left\|\mathbb{E}\left(\nabla^3 f(x)\right)\right\|_{\text{HS}}^2  
\leq & \sum_{k=1}^{p_1p_2p_3} \left(6a_kD_{k,k}\right)^2 + \sum_{k=1}\sum_{l=1, l \neq k} \left(2 a_k D_{l,k} + 2 a_k D_{k,l}\right)^2 + \sum_{k \neq l \neq m} \left(a_k D_{l,m}\right)^2 \\
\lesssim & \sum_{k=1} a_k^2 \sum_{l, m} D_{l,m}^2 = \left\|a\right\|_{\ell_2}^2 \cdot \left\|D\right\|_{\mathrm{F}}^2.
\end{align*} 

Combining all the results above, by Theorem 1.5 of \citet{gotze2021concentration}, we have
\begin{align*}
& \mathbb{P}\left(\left|\frac{1}{n^2\sigma^2}\sum_{i=1}^n \xi_i a^{\top} x_i x_i^{\top}\left(C \otimes B\right)x_i - \frac{K_3}{n^2\sigma^2} \sum_{i=1}^n \xi_i \sum_{k=1}^{p_1p_2p_3} a_k D_{k,k}\right| \geq \frac{t}{n^2} \Bigg| \left\{\xi_i\right\}_{i=1}^n\right) \\
\leq & \exp\left[-c\min\left(\frac{t^2}{\sigma^4\left\|\xi\right\|^2\left\|a\right\|_{\ell_2}^2 \cdot \left[\operatorname{tr}\left(D\right)\right]^2}, \left(\frac{t}{\sigma^3\left\|\xi\right\|_{\ell_2}\cdot \left\|a\right\|_{\ell_2} \cdot \left\|D\right\|_{\mathrm{F}}}\right)^{\frac{2}{3}}\right)\right].
\end{align*} 

Furthermore, note that by Bernstein-type inequality, we have
$$
\mathbb{P}\left( \left|\sum_{i=1}^n \xi_i\right| \geq C_1\sigma_\xi\sqrt{n\log(\op)}\right) \leq \frac{1}{\op^{c}}.
$$

Therefore, we have
\begin{align*}
& \mathbb{P}\left(\left|\frac{1}{n^2\sigma^2}\sum_{i=1}^n \xi_i a^{\top} x_i x_i^{\top}\left(C \otimes B\right)x_i\right| \geq \frac{K_3}{n^{3/2}\sigma^2} \left| \sum_{k=1}^{p_1p_2p_3} a_k D_{k,k}\right| + \frac{1}{n^{3/2}}   \left\|a\right\|_{\ell_2} \left(\left|\operatorname{tr(D)}\right| t^{1/2}+ \left\|D\right\|_{\mathrm{F}}  t^{2/3}\right) \right) \\
\leq & \exp\left(-ct\right) + \frac{1}{\op^{c}}.
\end{align*} 

Then consider the off-diagnoal terms, 
\begin{align*}
g(x) 
= & \sum_{i=1}^n \sum_{j=1, j \neq i}^{n} \xi_j a^{\top} x_i x_i^{\top}\left(C \otimes B\right)x_j  \\
= & \sum_{i=1}^n \sum_{j=1, j \neq i}^{n} \sum_{k=1}^{p_1p_2p_3} \sum_{m=1}^{p_1p_2p_3} a_k x_{i,k}^2 D_{k,m} \xi_j x_{j, m} + \sum_{i=1}^n \sum_{j=1, j \neq i}^{n} \sum_{k=1, k\neq l}^{p_1p_2p_3} \sum_{l=1}^{p_1p_2p_3} \sum_{m=1}^{p_1p_2p_3} a_k x_{i,k} x_{i,l} D_{l,m} \xi_j x_{j,m},
\end{align*}
which has expectation 0.

Given a fixed index $(i,j)$, where $i \neq j$, we consider the following function
\begin{align*}
h(x) =  \xi_j a^{\top} x_i x_i^{\top}\left(C \otimes B\right)x_j 
= & \sum_{k=1}^{p_1p_2p_3} \sum_{m=1}^{p_1p_2p_3} a_k x_{i,k}^2 D_{k,m} \xi_j x_{j, m} + \sum_{l=1}^{p_1p_2p_3} \sum_{m=1}^{p_1p_2p_3} a_k x_{i,k} x_{i,l} D_{l,m} \xi_j x_{j,m}.
\end{align*}

Then it follows that
\begin{align*}
\frac{\partial h}{\partial x_{i,k}} 
= & 2 \sum_{m=1}^{p_1p_2p_3} a_k x_{i,k} D_{k,m} \xi_j x_{j,m} + \sum_{l=1, l \neq k}^{p_1p_2p_3} \sum_{m=1}^{p_1p_2p_3} a_k D_{l,m} \xi_j x_{j,m}.
\end{align*}
and
\begin{align*}
\frac{\partial h}{\partial x_{j,m}} 
= & \sum_{k=1}^{p_1p_2p_3} \sum_{m=1}^{p_1p_2p_3} a_k x_{i,k}^2 D_{k,m} \xi_j + \sum_{k=1, k\neq l}^{p_1p_2p_3} \sum_{l=1}^{p_1p_2p_3} \sum_{m=1}^{p_1p_2p_3} a_k x_{i,k} x_{i,l} D_{l,m} \xi_j.
\end{align*}

It implies that
$$
\mathbb{E}\left(\frac{\partial h}{\partial x_{i,k}} \right)=0,
\quad \text{and} \quad
\mathbb{E}\left(\frac{\partial g}{\partial x_{j,m}} \right)= \sigma^2 \xi_j \sum_{k=1}^{p_1p_2p_3} \sum_{m=1}^{p_1p_2p_3} a_k D_{k,m} .
$$

It follows that 
$$
\left\|\mathbb{E}\left(\nabla f(x)\right) \right\|_{\text{HS}} \lesssim n^2\sigma^4 \sum_j \xi_j^2 \sum_m \left(\sum_{k,m}a_k D_{k,m}\right)^2 = n^2\sigma^4 \left\|\xi\right\|_{\ell_2}^2 \cdot \left\|D^{\top}a\right\|^2.
$$

Then, consider the second-order partial derivatives, we have
\begin{align*}
\frac{\partial^2 h}{\partial x_{i,k}^2} 
= & 2\sum_{m=1}^{p_1p_2p_3} a_kD_{k,m}\xi_j \xi_{j,m} ,\\
\frac{\partial^2 h}{\partial x_{i,k}\partial x_{j,m}}
= & 2\sum_{i=1}^n \sum_{j=1, j \neq i}^{n} \sum_{k=1}^{p_1p_2p_3} \sum_{m=1}^{p_1p_2p_3} a_k x_{i,k} D_{k,m} \xi_j + \sum_{i=1}^n \sum_{j=1, j \neq i}^{n} \sum_{k=1, k\neq l}^{p_1p_2p_3} \sum_{l=1}^{p_1p_2p_3} \sum_{m=1}^{p_1p_2p_3} a_k x_{i,l} D_{l,m} \xi_j.
\end{align*}

It implies that
$
\mathbb{E}\left(\frac{\partial^2 h}{\partial x_{i,k}^2}\right)= 0 \quad $ and $\mathbb{E}\left(\frac{\partial^2 h}{\partial x_{i,k}\partial x_{j,m}}\right)= 0 .
$
Thus, we have
$
\mathbb{E} \left(\nabla^2 f(x)\right)=0.
$

Finally, consider the third-order derivative, we have
$$
\frac{\partial^2 g}{\partial x_{i,k}^2\partial x_{j,m}} = 2a_k D_{k,m} \xi_j,
\quad \text{and} \quad
\frac{\partial^2 g}{\partial x_{i,l}\partial x_{i,k}\partial x_{j,m}} = a_k D_{l,m} \xi_j
$$
and the other third-order derivative should be equal to zero.

Therefore, we have
\begin{align*}
\left\|\mathbb{E} \left(\nabla^3 f(x)\right)\right\|_{\text{HS}}^2 
\lesssim & \sum_{i,j,k,l,m} a_k^2 D_{l,m}^2 \xi_j^2 \lesssim n\left\|\xi\right\|_{\ell_2}^2 \left\|D\right\|_{\mathrm{F}}^2 \cdot \left\|a\right\|_{\ell_2}^2.
\end{align*}

Combining the result above, by Theorem 1.5 of \citet{gotze2021concentration}, we have
\begin{align*}
\mathbb{P}\left( \left|\sum_{i=1}^n \sum_{j=1, j \neq i}^{n} \xi_j a^{\top} x_i\right| \geq  t\right) 
\leq & \exp\left[-c\min\left(\frac{t^2}{\sigma^2n\left\|\xi\right\|_{\ell_2}^2\cdot \left\|D\right\|_{\mathrm{F}}^2 \cdot \left\|a\right\|_{\ell_2}^2}, \left(\frac{t}{\sigma^3 \cdot \sqrt{n}\cdot \left\|\xi\right\| \cdot \left\|D\right\|_{\mathrm{F}} \cdot \left\|a\right\|_{\ell_2}}\right)^{2/3}\right)\right].
\end{align*}

Therefore, combining the results above, we have
\begin{align*}
& \mathbb{P}\left(\left|\frac{1}{n^2\sigma^2}\sum_{i=1}^n\sum_{j=1}^n \xi_i a^{\top} x_i x_i^{\top}\left(C \otimes B\right)x_j\right| \geq \frac{K_3}{n^{3/2}\sigma^2} \left\|a\right\|_{\ell_2}\cdot \left\|D\right\|_{\mathrm{F}} + \frac{1}{n^{3/2}} \cdot  \left\|a\right\|_{\ell_2}\cdot \left|\operatorname{tr(D)}\right| \cdot  t + \frac{1}{n} \cdot \left\|a\right\|_{\ell_2}\cdot \left\|D\right\|_{\mathrm{F}} \cdot t\right) \\
\leq & \exp\left[-c\min\left(t^2, t^{\frac{2}{3}}\right)\right] + \frac{1}{\op^{c}}.
\end{align*}
\end{proof}

\begin{lemma}\label{lemma: high-prob upper bound of tr(BZ(2)tCZ(2)) in tensor regression with sample splitting}

Let $X \in \mathbb{R}^{p \times d}$ be a random matrix with mean-zero $\sigma$ sub-Gaussian entries, and let $\left\{\mcX_i\right\}_{i=1}^{n}$ be $n$ i.i.d. copies of $X$ and $\left\{\xi_i\right\}_{i=1}^{n}$ are i.i.d. mean-zero $\sigma_\xi$-sub-Gaussian random variables. 

Assume that $A \in \mathbb{R}^{p \times d}$, $B\in \mathbb{R}^{d \times d}$ and $C \in \mathbb{R}^{p \times d}$ are three fixed matrices. Then it follows that
\begin{align}
& \mathbb{P}\Bigg(\Bigg|\operatorname{tr}\Bigg[B \Big(\frac{1}{n\sigma^2}\sum_{i=1}^n(\langle A, X_i \rangle X_i - \sigma^2  A)\Big)^{\top} C \Big(\frac{1}{n\sigma^2}\sum_{j=1}^n(\langle A, X_j \rangle X_j - \sigma^2   A)\Big)\Bigg] - \frac{1}{n \sigma^2}\big(\sigma^4 \|A \|_{\mathrm{F}}^2 \operatorname{tr}(B C)+\sigma^4 \operatorname{tr} (B A^{\top} C A )\big) \Bigg|  \notag \\
&\qquad \geq   \frac{1}{n} (\|A\|_{\mathrm{F}}^2  [\|B\|_{\mathrm{F}}  \|C\|_{\mathrm{F}} + |\operatorname{tr}(B)\operatorname{tr}(C)\|])  t^2 + n^{-1/2} \left\|A\right\|_{\mathrm{F}}^2 \left\|B\right\|  \left\|C\right\|   t^{3/2} \Bigg) \leq \exp\left(-ct\right). \label{eq: high-prob upper bound of tr(BZ(2)tCZ(2)) in tensor regression with sample splitting}
\end{align}

\end{lemma}

\begin{proof}

Consider 
\begin{align*}
& \frac{1}{n^2\sigma^2}\sum_{i=1}^n \sum_{j=1}^n \operatorname{tr}\left[B\left(\left\langle A, X_i \right\rangle X_i - \sigma^2 \cdot  A\right)^{\top} C \left(\left\langle A, X_j \right\rangle X_j - \sigma^2 \cdot  A\right)\right] \\
= & \frac{1}{n^2\sigma^4} \sum_{i=1}^n \sum_{j=1}^n \left\langle A, X_i\right\rangle \left\langle A, X_j\right\rangle \operatorname{tr}\left(BX_i^{\top}CX_j\right) - \frac{1}{n\sigma^2} \sum_{j=1}^n \left\langle A, X_j \right\rangle \operatorname{tr}\left(BA^{\top}CX_j\right) \\
& -  \frac{1}{n\sigma^2} \sum_{i=1}^n \left\langle A, X_i \right\rangle \operatorname{tr}\left(BX_i^{\top}CA\right) + \operatorname{tr}\left(BA^{\top}CA\right).
\end{align*}

Let 
\begin{align*}
f(X_i) 
= & \frac{1}{n\sigma^2}\left\langle A, X_i\right\rangle^2 \operatorname{tr}\left[BX_i^{\top}CX_i \right] = \frac{1}{n\sigma^2} \left(a^{\top}x_i\right)^2 x_i^{\top}Dx_i \\ 
= & \sum_{j=k, l=m}x_{i,j}^2A_{j,j}x_{i,l}^2D_{l,l} + \sum_{j=l,k=m} x_{i,j}^2A_{j,k}D_{j,k}x_{i,k}^2+\sum_{j=m, k=l} x_{i,j}^2A_{j,k}x_{i,k}^2D_{j,k},
\end{align*}
where $x_i=\operatorname{Vec}\left(X_i\right)$, $a=\operatorname{Vec}\left(A\right)$, and  $D = C^{\top}\otimes B$.

Furthermore, note that
\begin{align*}
\mathbb{E}\left(\left\langle A, X_i\right\rangle^2 \operatorname{tr}\left(BX_i^{\top}CX_i\right)\right) 
= &\mathbb{E}\left(x_i^{\top}aa^{\top}x_i x_i^{\top}Dx_i\right) 
=  \sigma^4\left\|A\right\|_{\mathrm{F}}^2\operatorname{tr}\left(BC\right) + \sigma^4\operatorname{tr}\left(BA^{\top}CA\right), \\
\mathbb{E}\left(\left\langle A, X_i\right\rangle \left\langle A, X_j\right\rangle \operatorname{tr}\left(BX_i^{\top}CX_j\right)\right) 
= & \mathbb{E}\left(x_i^{\top}aa^{\top}x_j x_j^{\top}Dx_i\right) 
= \sigma^4\operatorname{tr}\left(aa^{\top}D\right) 
= \sigma^4\operatorname{tr}\left(BA^{\top}CA\right),\\
\mathbb{E}\left[\left\langle A, X_j \right\rangle \operatorname{tr}\left(BA^{\top}CX_j\right)\right] = &\sigma^4 \operatorname{tr}\left(BA^{\top}CA\right).
\end{align*}

It follows that
\begin{align*}
& \mathbb{E}\left( \frac{1}{n^2\sigma^4}\sum_{i=1}^n \sum_{j=1}^n \operatorname{tr}\left[B\left(\left\langle A, X_i \right\rangle X_i - \sigma^2  A\right)^{\top} C \left(\left\langle A, X_j \right\rangle X_j - \sigma^2  A\right)\right] \right) 
=  \frac{1}{n}\left(\left\|A\right\|_{\mathrm{F}}^2 \operatorname{tr}\left(B C\right) + \operatorname{tr}\left(B A^{\top} C A\right)\right).
\end{align*}

We first consider the diagonal terms when $i=j$. By \eqref{eq: high-prob upper bound of tr(BZtCZ)} in Lemma~\ref{lemma: high-prob upper bound of tr(BZtCZ)}, we know that 
$$
\frac{1}{\left\|B\right\|_{\mathrm{F}}\left\|C\right\|_{\mathrm{F}}}\cdot \left[tr(BX^{\top}CX)  - \sigma^2\operatorname{tr}\left(B\right)tr\left(C\right)\right]
$$
is sub-exponential with parameter $C\sigma^2$ for some constant $C$.

By Remark 5.18 of \citet{vershynin2010introduction}, it follows that
\begin{align*}
& \left\|\left(\langle A, X\rangle\right)^2 \operatorname{tr}\left(B X^{\top} C X\right) - \mathbb{E}\left[\left(\langle A, X\rangle\right)^2 \operatorname{tr}\left(B X^{\top} C X\right)\right]\right\|_{\psi_{\frac{1}{2}}} 
\leq \left\|(\langle A, X\rangle)^2\right\|_{\psi_1} \cdot \left\|\operatorname{tr}\left(B X^{\top} C X\right)\right\|_{\psi_1} \\
\lesssim & \left\|\langle A, X\rangle\right\|_{\psi_2}^2 \cdot \left(\left\|\operatorname{tr}\left(B X^{\top} C X\right) - \sigma^2 \operatorname{tr}\left(B\right)\operatorname{tr}\left(C\right)\right\|_{\psi_1} + \left|\sigma^2 \operatorname{tr}\left(B\right)\operatorname{tr}\left(C\right)\right|\right) \\
\leq & \sigma^2\left\|A\right\|_{\mathrm{F}}^2 \cdot \left[\sigma^2\left\|B\right\|_{\mathrm{F}} \cdot \left\|C\right\|_{\mathrm{F}} + \sigma^2 \left|\operatorname{tr}\left(B\right)\operatorname{tr}\left(C\right)\right\|\right] ,
\end{align*}
where we used $\mathbb{E}\left(\operatorname{tr}\left(B X^{\top} C X\right)\right)=\sigma^2 \operatorname{tr}\left(B\right)\operatorname{tr}\left(C\right)$ and $\left\|\mathbb{E}X\right\|_{\psi_1} \lesssim \left|\mathbb{E}X\right|$.

Therefore, by Bernstein-type inequality, we have
\begin{align*}
& \mathbb{P}\left(\left|\frac{1}{n\sigma^2}\sum_{i=1}^n \left(\left(\langle A, X_i\rangle\right)^2 \operatorname{tr}\left(B X_i^{\top} C X_i\right) - \sigma^2 \mathbb{E}\left[\left(\langle A, X_i\rangle\right)^2 \operatorname{tr}\left(B X_i^{\top} C X_i\right)\right]\right)\right| \geq t\right) \\
\leq & \exp\left[-\frac{ct}{n\sigma^4\left\|A\right\|_{\mathrm{F}}^2\cdot \left(\left[\sigma^2\left\|B\right\|_{\mathrm{F}} \cdot \left\|C\right\|_{\mathrm{F}}+\sigma^2 \left| \operatorname{tr}(B) \operatorname{tr}(C) \right|\right]\right)}\right],
\end{align*}

Then consider the off-diagonal terms, first we can write it as
\begin{align*}
f(x) 
:= & \sum_{i=1}^n \sum_{j=1, j \neq i}^{n} a^{\top}x_i a^{\top}x_j\left(x_i-a\right)^{\top}D\left(x_j-a\right) \\
= & \sum_{i=1}^n\sum_{j=1, j\neq i}^n \sum_{k=1}^{d}\sum_{l=1}^{d}\sum_{q=1}^{d}a_kx_{i,k}a_lx_{j,l}\left(x_{i,k}-a_k\right)D_{k,q}\left(x_{j,q} - a_q\right) \\
&+  \sum_{i=1}^n\sum_{j=1, j\neq i}^n \sum_{k=1}^{d}\sum_{l=1}^{d}\sum_{p=1, p \neq k}^{d}\sum_{q=1}^{d}a_kx_{i,k}a_lx_{j,l}\left(x_{i,p}-a_p\right)D_{p,q}\left(x_{j,q} - a_q\right).
\end{align*}

Therefore, we have the first-order derivative
\begin{align*}
\frac{\partial f}{\partial x_{i,k}} 
= & \sum_{j=1, j\neq i}^n \sum_{l=1}^{d}\sum_{q=1}^{d} a_ka_lx_{j,l} \cdot \left(2x_{i,k} - a_k\right) D_{k,q}\left(x_{j,q}-a_q\right) \\
& +  \sum_{j=1, j\neq i}^n \sum_{l=1}^{d}\sum_{p=1, p \neq k}^{d}\sum_{q=1}^{d}a_ka_lx_{j,l}\left(x_{i,p}-a_p\right)D_{p,q}\left(x_{j,q} - a_q\right).
\end{align*}

It follows that
\begin{align*}
\mathbb{E}\left(\frac{\partial f}{\partial x_{i,k}}\right) 
= & - \sigma^2 \sum_{j=1, j\neq i}^n \sum_{l=1}^{d} a_k^2 a_l D_{k,l}a_l - \sigma^2\sum_{j=1, j\neq i}^n \sum_{l=1}^{d}\sum_{p=1, p \neq k}^{d} a_ka_la_pD_{p,l} 
= (n-1)  a_k a^{\top}Da .
\end{align*}

Therefore, we have
$
\left\|\mathbb{E} \nabla f(x)\right\|_{\mathrm{F}}^2 = n(n-1)^2 \cdot \left\|a\right\|_{\ell_2}^2 \cdot a^{\top}Da^2
\leq  n^3 \cdot \left\|a\right\|_{\ell_2}^4 \cdot \left\|D\right\|^2.
$

For the second-order derivative, it follows that
\begin{align*}
\frac{\partial^2 f}{\partial x_{i,k}^2} 
= 2 \sum_{j=1, j\neq i}^n \sum_{l=1}^{d}\sum_{q=1}^{d} a_ka_lx_{j,l} D_{k,q}\left(x_{j,q}-a_q\right) , \frac{\partial^2 f}{\partial x_{i,p} \partial x_{i,k}} 
= \sum_{j=1, j\neq i}^n \sum_{l=1}^{d}\sum_{q=1}^{d}a_ka_lx_{j,l}D_{p,q}\left(x_{j,q} - a_q\right).
\end{align*}

It follows that
\begin{align*}
\mathbb{E}\left(\frac{\partial^2 f}{\partial x_{i,k}^2} \right) = 2\sigma^2 \cdot (n-1) \cdot a_k \sum_{l=1}^d a_lD_{k,l}, \quad \text{and} \quad
\mathbb{E}\left(\frac{\partial^2 f}{\partial x_{i,p} \partial x_{i,k}}\right) = 2\sigma^2 \cdot (n-1) \cdot a_k \sum_{l=1}^d a_lD_{p,l} .
\end{align*}

Furthermore, note that
\begin{align*}
\frac{\partial f}{\partial x_{i,k}} 
= & \sum_{j=1, j\neq i}^n \sum_{l=1}^{d} a_ka_lx_{j,l} \cdot \left(2x_{i,k} - a_k\right) D_{k,l}\left(x_{j,l}-a_l\right) 
 +  \sum_{j=1, j\neq i}^n \sum_{l=1}^{d}\sum_{p=1, p \neq k}^{d}a_ka_lx_{j,l}\left(x_{i,p}-a_p\right)D_{p,l}\left(x_{j,l} - a_l\right) \\
& +  \sum_{j=1, j\neq i}^n \sum_{l=1}^{d}\sum_{q=1, q\neq l}^{d} a_ka_lx_{j,l} \cdot \left(2x_{i,k} - a_k\right) D_{k,q}\left(x_{j,q}-a_q\right) \\
& +  \sum_{j=1, j\neq i}^n \sum_{l=1}^{d}\sum_{p=1, p \neq k}^{d}\sum_{q=1, q\neq l}^{d}a_ka_lx_{j,l}\left(x_{i,p}-a_p\right)D_{p,q}\left(x_{j,q} - a_q\right).
\end{align*}

Therefore, we have
\begin{align*}
\frac{\partial^2 f}{\partial x_{j,l}\partial x_{i,k}} 
= &  a_ka_l \cdot \left(2x_{j,l} -a_l\right) \cdot \left(2x_{i,k} - a_k\right) D_{k,l} + \sum_{p=1, p \neq k}^{d}a_ka_l\left(2x_{j,l} -a_l\right)\left(x_{i,p}-a_p\right)D_{p,l} \\
& +  \sum_{q=1, q\neq l}^{d} a_ka_l \cdot \left(2x_{i,k} - a_k\right) D_{k,q}\left(x_{j,q}-a_q\right) + \sum_{p=1, p \neq k}^{d}\sum_{q=1, q\neq l}^{d}a_ka_l\left(x_{i,p}-a_p\right)D_{p,q}\left(x_{j,q} - a_q\right).
\end{align*}

It implies that
$
\mathbb{E}\left(\frac{\partial^2 f}{\partial x_{j,l}\partial x_{i,k}} \right)
= \sum_{p=1}^{d}\sum_{q=1}^{d}a_ka_la_pa_qD_{p,q}.
$

Therefore, we have
\begin{align*}
\mathbb{E}\left(\nabla^2 f\right) 
\lesssim & n^2\sigma^4\sum_{i,k} a_k^2 \left(\sum_{l=1}^d a_l D_{k, l}\right)^2 + n^2\sigma^4\sum_{i,k,p} a_k^2 \left(\sum_{l=1}^d a_l D_{p, l}\right)^2 + \sum_{i,j,k,l} a_k^2 a_l^2 \left(a^{\top}Da\right)^2 
\lesssim n^3 \sigma^4 \cdot \left\|a\right\|_{\ell_2}^4 \cdot \left\|D\right\|^2.
\end{align*}

Then, consider the third-order derivative, we have
\begin{align*}
\frac{\partial^3 f}{\partial x_{j,l}\partial x_{i,k}^2} = & 2 \sum_{j=1, j\neq i}^n a_ka_l \cdot \left(2x_{j,l} -a_l\right) D_{k,l} + 2 \sum_{q=1, q\neq l}^{d}a_ka_lD_{k,q}\left(x_{j,q}-a_q\right), \\
\frac{\partial^3 f}{\partial x_{i,p} \partial x_{j,l}\partial x_{i,k}}
= & a_ka_l\left(2x_{j,l} -a_l\right)D_{p,l} + \sum_{q=1, q\neq l}^{d}a_ka_l D_{p,q}\left(x_{j,q} - a_q\right).
\end{align*}

It implies that
\begin{align*}
\mathbb{E}\left(\frac{\partial^3 f}{\partial x_{j,l}\partial x_{i,k}^2}\right) = & -2 a_ka_l^2 D_{k,l} - 2 \sum_{q=1, q\neq l}^{d}a_ka_la_qD_{k,q} = -2 \sum_{q=1}^n a_ka_la_qD_{k,q}\\
\mathbb{E} \left(\frac{\partial^3 f}{\partial x_{i,p} \partial x_{j,l}\partial x_{i,k}}\right)
= & - a_ka_l^2D_{p,l} - \sum_{q=1, q\neq l}^{d}a_ka_la_q D_{p,q} =  \sum_{q=1}^{d} a_ka_la_q D_{p,q}.
\end{align*}

Therefore, we have
\begin{align*}
\left\|\mathbb{E}\left(\nabla^3 f\right) \right\|_{\text{HS}}^2
\lesssim & \sum_{i,j,k,l}a_k^2a_l^2\left(\sum_{q=1}^na_qD_{k,q}\right)^2 + \sum_{i,j,k,l,p}\left(\sum_{q=1}^na_qD_{p,q}\right)^2 
\leq n^2 \left\|a\right\|_{\ell_2}^2 \left\|Da\right\|_{\ell_2}^2
\end{align*}

Finally, consider the following fourth-order partial derivatives
\begin{align*}
& \frac{\partial^4 f}{\partial x_{j,l}^2 \partial x_{i,k}^2} = 4 a_ka_l D_{k,l}, \quad \frac{\partial^4 f}{\partial x_{j,q} x_{j,l}\partial x_{i,k}^2} = 2 a_ka_lD_{k,q}, \\
& \frac{\partial^4 f}{\partial x_{i,p} \partial x_{j,l}^2 \partial x_{i,k}}
= 2 a_ka_lD_{p,l}, \quad \frac{\partial^3 f}{\partial x_{j,q} \partial x_{i,p} \partial x_{j,l}\partial x_{i,k}}
= a_ka_l D_{p,q}.
\end{align*}

Therefore, we have
\begin{align*}
\left\|\mathbb{E}\left(\nabla^4 f\right) \right\|_{\text{HS}}^2
\leq & \sum_{i,j,k,l} a_k^2 a_l^2 D_{k,l}^2 + \sum_{i,j,k,q,l} a_k^2 a_l^2 D_{k,q}^2 + \sum_{i,j,p,q,k,l}a_k^2 a_l^2 D_{p,q}^2 \leq n^2 \left\|a\right\|_{\ell_2}^4 \cdot \left\|D\right\|_{\mathrm{F}}^2.
\end{align*}

Combining all the results above, by Theorem 1.5 of \citet{gotze2021concentration}, we have
\begin{align*}
& \mathbb{P}\left(\left|\frac{1}{n\sigma^2}\sum_{i=1}^n \left(\left(\langle A, X_i\rangle\right)^2 \operatorname{tr}\left(B X_i^{\top} C X_i\right) - \sigma^2 \mathbb{E}\left[\left(\langle A, X_i\rangle\right)^2 \operatorname{tr}\left(B X_i^{\top} C X_i\right)\right]\right)\right| \geq t\right) \\
\leq & \exp\left\{-c\min\left[\left(\frac{t}{\sigma^4\cdot n^\frac{3}{2} \cdot \left\|A\right\|_{\mathrm{F}}^2 \cdot \left\|B\right\|\cdot \left\|C\right\|}\right)^2, \left(\frac{t}{\sigma^4\cdot n^\frac{3}{2} \cdot \left\|A\right\|_{\mathrm{F}}^2 \cdot \left\|B\right\|\cdot \left\|C\right\|}\right), \right. \right. \\
& \left.\left. \left(\frac{t}{\sigma^4\cdot n \cdot \left\|A\right\|_{\mathrm{F}}^2 \cdot \left\|B\right\|\cdot \left\|C\right\|}\right)^{2/3}, \left(\frac{t}{\sigma^4\cdot n \cdot \left\|A\right\|_{\mathrm{F}}^2 \cdot \left\|B\right\|_{\mathrm{F}}\cdot \left\|C\right\|_{\mathrm{F}}}\right)\right]\right\}
\end{align*}

\end{proof}

\begin{lemma}\label{lemma: high-prob upper bound of spectral norm of B1tZhat1(1)tA1t(A2Zhat2(1)B2oA3Zhat3(1)B3) in tensor regression without sample splitting}
Let $\mcX \in \mathbb{R}^{p_1 \times p_2 \times p_3}$ be a random tensor with mean-zero $\sigma$ sub-Gaussian entries, and let $\left\{\mcX_i\right\}_{i=1}^{n}$ be $n$ i.i.d. copies of $\mcX$. Define $\widehat{\mcZ}^{(1)} = \frac{1}{n}\sum_{i=1}^{n} \xi_i \mcX_i$, where $\left\{\xi_i\right\}_{i=1}^{n}$ are i.i.d. mean-zero $\sigma_\xi$-sub-Gaussian random variables. Let $\widehat{\mcZ}^{(1)}_j = \Mat_j(\widehat{\mcZ}^{(1)})$ denote the mode-$j$ matricization of the random tensor $\widehat{\mcZ}^{(1)}$.
Assume that $U_j \in \mathbb{O}^{p_j \times r_j}$ for $j = 1, 2, 3$. Additionally, let $U_{j\perp} \in \mathbb{O}^{p_j \times (p_j - r_j)}$ such that $U_{j\perp} U_{j\perp}^{\top} \in \mathbb{R}^{p_j \times p_j}$ is a projection matrix that projects any vector onto the orthogonal complement of the subspace spanned by $U_j U_j^{\top}$.
Then it holds that
\begin{align}
& \mathbb{P}\left(\left\|\left(\mcP_{U_{j+2}} \otimes \mcP_{U_{j+1}}\right) \whZ_j^{(1)\top}\mcP_{U_{j\perp}}A_j \left(\mcP_{U_{j+1 \perp}}\whZ_{j+1}^{(1)}\left(\mcP_{U_j} \otimes \mcP_{U_{j+2}}\right) \otimes \mcP_{U_{j+2 \perp}}\whZ_{j+2}^{(1)}\left(\mcP_{U_{j+1}} \otimes \mcP_{U_j}\right)\right) \right\|  \right. \notag \\
& \quad \geq \left. C\cdot \frac{1}{n^3} \left\|\mcA\right\|_{\mathrm{F}} \cdot t \bigg| \left\{\xi_i\right\}_{i=1}^n \right)  \leq 7^{r_j+r_{j+1}r_{j+2}} \exp\left[-c\min\left(\frac{t^2}{\left\|\xi\right\|_{\ell_2}^2}, \frac{t^{\frac{2}{3}}}{\left\|\xi\right\|_{\ell_2}^2}\right)\right].
\end{align}

Furthermore, we have
\begin{align}
& \mathbb{P}\left(\left\|\left(\mcP_{U_{j+2}} \otimes \mcP_{U_{j+1}}\right) \whZ_j^{(1)\top}\mcP_{U_{j\perp}}A_j \left(\mcP_{U_{j+1 \perp}}\whZ_{j+1}^{(1)}\left(\mcP_{U_j} \otimes \mcP_{U_{j+2}}\right) \otimes \mcP_{U_{j+2 \perp}}\whZ_{j+2}^{(1)}\left(\mcP_{U_{j+1}} \otimes \mcP_{U_j}\right)\right) \right\| \right. \notag \\
& \quad \left. \geq C\cdot \frac{1}{n^{3/2}} \left\|\mcA\right\|_{\mathrm{F}} \cdot t\right)  \leq 7^{r_j+r_{j+1}r_{j+2}} \exp\left[-c\min\left(t^2, t^{\frac{2}{3}}\right)\right]. \label{eq: high-prob upper bound of spectral norm of B1tZhat1(1)tA1t(A2Zhat2(1)B2oA3Zhat3(1)B3) in tensor regression without sample splitting}
\end{align}

\end{lemma}

\begin{proof}

By symmetry, it suffices to find a high-probability upper bound for
\begin{align*}
& \left\|\left(\mcP_{U_3} \otimes \mcP_{U_2}\right) \whZ_1^{(1)\top}\mcP_{U_{1\perp}}A_1\left(\mcP_{U_{3\perp}}\whZ_3^{(1)}\left(\mcP_{U_2} \otimes \mcP_{U_1}\right) \otimes  \mcP_{U_{2 \perp}}\whZ_2^{(1)}\left(\mcP_{U_1} \otimes \mcP_{U_3}\right) \right) \right\|.
\end{align*}
For the first inequality, we only need to observe that conditioning on $\xi_i$, $\widehat{\mcZ}^{(1)}=\frac{1}{n}\sum_{i=1}^n \xi_i\mcX_i$ has i.i.d. sub-Gaussian entries with variance $\sigma^2\left\|\xi\right\|_{\ell_2}^2$. Then the first inequality follows immediately by applying Lemma~\ref{lemma: high-prob upper bound of spectral norm of B1tZ1(1)tA1t(A2Z2(1)B2oA3Z3(1)B3) in tensor PCA}. 

The second inequality follows from the following Bernstein-type inequality for random vectors with i.i.d. sub-Gaussian entries:
$$
\mathbb{P}\left(\left\|\xi\right\|_{\ell_2} \geq Ct \right) \leq \exp\left[-c\min\left(nt^2, nt\right)\right].
$$
Let $t=\frac{1}{n}$. Then it immediately follows that
$
\mathbb{P}\left(\|\xi\|_2 \geq C \sqrt{n}\right) \leq \exp\left(-c n\right).
$

\end{proof}

\begin{lemma}\label{lemma: high-prob upper bound of spectral norm of A1t(A2Z2(1)B2oA3Z3(1)B3) in tensor regression}
Let $\mcX \in \mathbb{R}^{p_1 \times p_2 \times p_3}$ be a random tensor with i.i.d. mean-zero $\sigma$-sub-Gaussian entries, and let $\left\{\mcX_i\right\}_{i=1}^{n}$ be $n$ i.i.d. copies of $\mcX$. Define $\widehat{\mcZ}^{(1)} = \frac{1}{n}\sum_{i=1}^{n} \xi_i \mcX_i$, where $\left\{\xi_i\right\}_{i=1}^{n}$ are i.i.d. mean-zero $\sigma_\xi$-sub-Gaussian random variables. Let $\widehat{\mcZ}^{(1)}_j = \Mat_j(\widehat{\mcZ}^{(1)})$ denote the mode-$j$ matricization of the random tensor $\widehat{\mcZ}^{(1)}$.
Assume that $U_j \in \mathbb{O}^{p_j \times r_j}$ for $j = 1, 2, 3$. Additionally, let $U_{j\perp} \in \mathbb{O}^{p_j \times (p_j - r_j)}$ such that $U_{j\perp} U_{j\perp}^{\top} \in \mathbb{R}^{p_j \times p_j}$ is a projection matrix that projects any vector onto the orthogonal complement of the subspace spanned by $U_j U_j^{\top}$.
Then it holds that
\begin{align}
& \mathbb{P}\left(\left\|\left(\mcP_{U_{j+2}} \otimes \left(\mcP_{U_j}\otimes \mcP_{U_{j+2}}\right)\whZ_{j+1}^{(1)}\mcP_{U_{j+1\perp}}\right)A_j^{\top}\mcP_{U_{j\perp}}\whZ_j^{(1)}\left(U_{j+2}\otimes U_{j+1}\right)\right\| \right. \notag \\
&\quad \geq \left. C \sigma^2 \cdot \frac{1}{n^2}\left\|\left(\mcP_{U_{j+1\perp}} \otimes \mcP_{U_{j+2}}\right) A_j^{\top} \mcP_{U_{j \perp}}\right\|_{\mathrm{F}}  t \bigg| \left\{\xi_i\right\}_{i=1}^n\right)
\leq 7^{r_j+r_{j+1}r_{j+2}} \exp\left[-c\min\left(\frac{t^2}{\left\|\xi\right\|_{\ell_2}^2}, \frac{t^{\frac{2}{3}}}{\left\|\xi\right\|_{\ell_2}^2}\right)\right] .
\end{align}

Furthermore, we have
\begin{align}
& \mathbb{P}\left(\left\|\left(\mcP_{U_{j+2}} \otimes \left(\mcP_{U_j}\otimes \mcP_{U_{j+2}}\right)\whZ_{j+1}^{(1)}\mcP_{U_{j+1\perp}}\right)A_j^{\top}\mcP_{U_{j\perp}}\whZ_j^{(1)}\left(U_{j+2}\otimes U_{j+1}\right)\right\| \right. \notag \\
&\quad \geq \left. C \sigma^2 \cdot \frac{1}{n}\left\|\left(\mcP_{U_{j+1\perp}} \otimes \mcP_{U_{j+2}}\right) A_j^{\top} \mcP_{U_{j \perp}}\right\|_{\mathrm{F}} \cdot t\right) \leq 7^{r_j+r_{j+1}r_{j+2}} \exp \left(-c \min \left(t^2, t\right)\right) + \exp\left(-cn\right). \label{eq: high-prob upper bound of spectral norm of A1t(A2Z2(1)B2oA3Z3(1)B3) in tensor regression}
\end{align}

\end{lemma}

The proof of Lemma~\ref{lemma: high-prob upper bound of spectral norm of A1t(A2Z2(1)B2oA3Z3(1)B3) in tensor regression} is similar to that of Lemma~\ref{lemma: high-prob upper bound of spectral norm of B1tZhat1(1)tA1t(A2Zhat2(1)B2oA3Zhat3(1)B3) in tensor regression without sample splitting}, and is thus omitted.



\begin{lemma}\label{lemma: high-prob upper bound of spectral norm of B1tZ1hatA1t(A2Z2hatB2oA3Z3hatB3) in tensor regression without sample splitting}
Under the same setting of Theorem~\ref{thm: main theorem in tensor regression without sample splitting}. Let $\widehat{\mcZ}=\frac{1}{n\sigma^2}\sum_{i=1}^{n}\xi\mcX_i + \frac{1}{n\sigma^2}\sum_{i=1}^n\left[\left\langle \mcX_i, \widehat{\Delta}\right\rangle \mcX_i - \sigma^2\cdot \widehat{\Delta}\right] \in \mathbb{R}^{p_1 \times p_2 \times p_3}$ be the debiased error. Let $\whZ_j = \Mat_j(\widehat{\mcZ})$ denote the mode-$j$ matricization of the tensor $\widehat{\mcZ}$. Additionally, let $U_{j\perp} \in \mathbb{O}^{p_j \times (p_j - r_j)}$ be orthonormal matrices such that $U_{j\perp} U_{j\perp}^{\top} \in \mathbb{R}^{p_j \times p_j}$ is a projection matrix that projects any vector onto the orthogonal complement of the space spanned by $U_j U_j^{\top}$. Then,
\begin{align}
& \left\|\mcA \times_j \mcP_{U_{j\perp}} \whZ_j \left(\mcP_{U_{j+2}} \otimes \mcP_{U_{j+1}}\right) \times_{j+1} \mcP_{U_{{j+1}\perp}} \whZ_{j+1} \left(\mcP_{U_j} \otimes \mcP_{U_{j+2}}\right) \times_{j+2} \mcP_{U_{{j+2}\perp}} \whZ_{j+2} \left(\mcP_{U_{j+1}} \otimes \mcP_{U_j}\right) \right\|_{\mathrm{F}} \notag \\
\lesssim & \left\|\mcA\right\|_{\mathrm{F}} \cdot \left(\frac{\sigma_{\xi}^3}{\sigma^3} \frac{\oor^{3/2}\log(\op)^{3/2}}{n^{3 / 2}}+\frac{\sigma_{\xi}^3}{\sigma^3} \cdot \Delta \cdot \frac{\op^{1 / 2} \oR \log (\op)}{n^{3 / 2}} + \Delta^2 \cdot \frac{\sigma_{\xi}^3}{\sigma^3} \cdot \frac{\op \oR^{1/2}\log(\op)^{1/2}}{n^{3 / 2}} + \Delta^3 \cdot \frac{\sigma_{\xi}^3}{\sigma^3} \cdot \frac{\op^{3 / 2}}{n^{3 / 2}}\right), \label{eq: high-prob upper bound of Ao(PU1pZhat1(PU3oPU2))o(PU2pZhat2(PU1oPU3))o(PU3pZhat3(PU2oPU1)) in tensor regression without sample splitting}
\end{align}
and
\begin{align}
& \left\|\mcA \times_j \mcP_{U_j} \times_{j+1} \mcP_{U_{{j+1}\perp}} \whZ_{j+1} \left(\mcP_{U_j} \otimes \mcP_{U_{j+2}}\right) \times_{j+2} \mcP_{U_{{j+2}\perp}} \whZ_{j+2} \left(\mcP_{U_{j+1}} \otimes \mcP_{U_j}\right) \right\|_{\mathrm{F}} \notag \\
\lesssim & \left\|\mcA \times_j U_j\right\|_{\mathrm{F}} \cdot \frac{\sigma_{\xi}^2}{\sigma^2} \cdot \left( \frac{\oor\log(\op)}{n} + \Delta \cdot \frac{\op^{1 / 2} \sqrt{\oR \log (\op)}}{n} + \Delta^2 \cdot \frac{\op}{n}\right) \label{eq: high-prob upper bound of AoPU1o(PU2pZhat2(PU1oPU3))o(PU3pZhat3(PU2oPU1)) in tensor regression without sample splitting}
\end{align}
hold with probability at least $1-\exp(-cn) - \frac{1}{\op^c} - \mathbb{P}\left(\mcE_U^{\text{reg}}\right) - \mathbb{P}\left(\mcE_\Delta\right)$ for any $j=1,2,3$.

\end{lemma}

\begin{proof}

{\bf Step 1: Proof of the first inequality}

\begin{align*}
& \left\|\left(\mcP_{U_3} \otimes \mcP_{U_2}\right) \whZ_1^{\top}\mcP_{U_{1\perp}}A_1 \left(\mcP_{U_{2 \perp}}\whZ_2\left(\mcP_{U_1} \otimes \mcP_{U_3}\right) \otimes \mcP_{U_{3 \perp}}\whZ_{3}\left(\mcP_{U_2} \otimes \mcP_{U_1}\right)\right) \right\|_{\mathrm{F}} \\
\leq & \underbrace{\left\|\left(\mcP_{U_3} \otimes \mcP_{U_2}\right) \whZ_1^{(1) \top} \mcP_{U_{1 \perp}} A_1\left(\mcP_{U_{3 \perp}} \whZ_3^{(1)}\left(\mcP_{U_2} \otimes \mcP_{U_1}\right) \otimes \mcP_{U_{2 \perp}} \whZ_2^{(1)}\left(\mcP_{U_1} \otimes \mcP_{U_3}\right)\right)\right\|_{\mathrm{F}}}_{\mathrm{\RN{1}}} \\
+ & \underbrace{\left\|\left(\mcP_{U_3} \otimes \mcP_{U_2}\right) \whZ_1^{(2) \top} \mcP_{U_{1 \perp}} A_1\left(\mcP_{U_{3 \perp}} \whZ_3^{(1)}\left(\mcP_{U_2} \otimes \mcP_{U_1}\right) \otimes \mcP_{U_{2 \perp}} \whZ_2^{(1)}\left(\mcP_{U_1} \otimes \mcP_{U_3}\right)\right)\right\|_{\mathrm{F}}}_{\mathrm{\RN{2}}} \\
+ & \underbrace{\left\|\left(\mcP_{U_3} \otimes \mcP_{U_2}\right) \whZ_1^{(1) \top} \mcP_{U_{1 \perp}} A_1\left(\mcP_{U_{3 \perp}} \whZ_3^{(2)}\left(\mcP_{U_2} \otimes \mcP_{U_1}\right) \otimes \mcP_{U_{2 \perp}} \whZ_2^{(1)}\left(\mcP_{U_1} \otimes \mcP_{U_3}\right)\right)\right\|_{\mathrm{F}}}_{\mathrm{\RN{3}}} \\
+ & \underbrace{\left\|\left(\mcP_{U_3} \otimes \mcP_{U_2}\right) \whZ_1^{(1) \top} \mcP_{U_{1 \perp}} A_1\left(\mcP_{U_{3 \perp}} \whZ_3^{(1)}\left(\mcP_{U_2} \otimes \mcP_{U_1}\right) \otimes \mcP_{U_{2 \perp}} \whZ_2^{(2)}\left(\mcP_{U_1} \otimes \mcP_{U_3}\right)\right)\right\|_{\mathrm{F}}}_{\mathrm{\RN{4}}} \\
+ & \underbrace{\left\|\left(\mcP_{U_3} \otimes \mcP_{U_2}\right) \whZ_1^{(2) \top} \mcP_{U_{1 \perp}} A_1\left(\mcP_{U_{3 \perp}} \whZ_3^{(2)}\left(\mcP_{U_2} \otimes \mcP_{U_1}\right) \otimes \mcP_{U_{2 \perp}} \whZ_2^{(1)}\left(\mcP_{U_1} \otimes \mcP_{U_3}\right)\right)\right\|_{\mathrm{F}}}_{\mathrm{\RN{5}}} \\
+ & \underbrace{\left\|\left(\mcP_{U_3} \otimes \mcP_{U_2}\right) \whZ_1^{(2) \top} \mcP_{U_{1 \perp}} A_1\left(\mcP_{U_{3 \perp}} \whZ_3^{(1)}\left(\mcP_{U_2} \otimes \mcP_{U_1}\right) \otimes \mcP_{U_{2 \perp}} \whZ_2^{(2)}\left(\mcP_{U_1} \otimes \mcP_{U_3}\right)\right)\right\|_{\mathrm{F}}}_{\mathrm{\RN{6}}} \\
+ & \underbrace{\left\|\left(\mcP_{U_3} \otimes \mcP_{U_2}\right) \whZ_1^{(1) \top} \mcP_{U_{1 \perp}} A_1\left(\mcP_{U_{3 \perp}} \whZ_3^{(2)}\left(\mcP_{U_2} \otimes \mcP_{U_1}\right) \otimes \mcP_{U_{2 \perp}} \whZ_2^{(2)}\left(\mcP_{U_1} \otimes \mcP_{U_3}\right)\right)\right\|_{\mathrm{F}}}_{\mathrm{\RN{7}}} \\
+ & \underbrace{\left\|\left(\mcP_{U_3} \otimes \mcP_{U_2}\right) \whZ_1^{(2) \top} \mcP_{U_{1 \perp}} A_1\left(\mcP_{U_{3 \perp}} \whZ_3^{(2)}\left(\mcP_{U_2} \otimes \mcP_{U_1}\right) \otimes \mcP_{U_{2 \perp}} \whZ_2^{(2)}\left(\mcP_{U_1} \otimes \mcP_{U_3}\right)\right)\right\|_{\mathrm{F}}}_{\mathrm{\RN{8}}}.
\end{align*}

Here, it follows from \eqref{eq: high-prob upper bound of spectral norm of B1tZhat1(1)tA1t(A2Zhat2(1)B2oA3Zhat3(1)B3) in tensor regression without sample splitting} in Lemma~\ref{lemma: high-prob upper bound of spectral norm of B1tZhat1(1)tA1t(A2Zhat2(1)B2oA3Zhat3(1)B3) in tensor regression without sample splitting} that
\begin{align*}
\mathrm{\RN{1}}
& \lesssim \underbrace{\left\|\mcA\right\|_{\mathrm{F}} \cdot \left(\frac{\sigma_\xi^3}{\sigma^3}\cdot \frac{\oor\log(\op)^{3 / 2}}{n^{3/2}}\right)}_{\eqref{eq: high-prob upper bound of spectral norm of B1tZhat1(1)tA1t(A2Zhat2(1)B2oA3Zhat3(1)B3) in tensor regression without sample splitting}}.
\end{align*}

Then consider
\begin{align*}
\mathrm{\RN{2}}
\leq & \left\|\left(\mcP_{U_3} \otimes \mcP_{U_2}\right) \whZ_1^{(2) \top} \mcP_{U_{1 \perp}}\right\| \cdot \left\|A_1\right\|_{\mathrm{F}} \cdot \left\|V_3^{\top}\mcP_{U_{3 \perp}} \whZ_3^{(1)}\left(\mcP_{U_2} \otimes \mcP_{U_1}\right) \right\| \cdot \left\|V_2^{\top}\mcP_{U_{2 \perp}} \whZ_2^{(1)}\left(\mcP_{U_2} \otimes \mcP_{U_1}\right) \right\| \\
\lesssim & \Delta\cdot \frac{\sigma_\xi}{\sigma}\sqrt{\frac{\op}{n}} \cdot 
\left\|\mcP_{U_{1 \perp}} A_1\left(\mcP_{U_{3 \perp}} \otimes \mcP_{U_{1 \perp}}\right)\right\|_{\mathrm{F}} \cdot \left(\frac{\sigma_\xi}{\sigma}\sqrt{\frac{\oR\log(\op)}{n}}\right)^2 = \Delta \cdot \frac{\sigma_\xi^3}{\sigma^3}\cdot \frac{\op^{1/2}\oR\log(\op)}{n^{3/2}} \cdot \left\|\mcA\right\|_{\mathrm{F}}.
\end{align*}

By symmetry, we have
\begin{align*}
\mathrm{\RN{3}}
 \lesssim \Delta \cdot \frac{\sigma_\xi^3}{\sigma^3}\cdot \frac{\op^{1/2}\oR\log(\op)}{n^{3/2}} \cdot \left\|\mcA\right\|_{\mathrm{F}},
\quad \text{and} \quad
\mathrm{\RN{4}}
 \lesssim \Delta \cdot \frac{\sigma_\xi^3}{\sigma^3}\cdot \frac{\op^{1/2}\oR\log(\op)}{n^{3/2}} \cdot \left\|\mcA\right\|_{\mathrm{F}}.
\end{align*}

Similarly, we have
\begin{align*}
\mathrm{\RN{5}}
\lesssim & \left(\Delta\cdot \frac{\sigma_\xi}{\sigma}\sqrt{\frac{\op}{n}}\right)^2 \cdot 
\left\|\mcP_{U_{1 \perp}} A_1\left(\mcP_{U_{3 \perp}} \otimes \mcP_{U_{1 \perp}}\right)\right\|_{\mathrm{F}} \cdot \left(\frac{\sigma_\xi}{\sigma}\sqrt{\frac{\oR\log(\op)}{n}}\right) = \Delta^2 \cdot \frac{\sigma_\xi^3}{\sigma^3}\cdot \frac{\op\oR^{1/2}\log(\op)^{1/2}}{n^{3/2}} \cdot \left\|\mcA\right\|_{\mathrm{F}}. 
\end{align*}

By symmetry, we have
$
\mathrm{\RN{6}} \lesssim \Delta^2 \cdot \frac{\sigma_\xi^3}{\sigma^3}\cdot \frac{\op\oR^{1/2}\log(\op)^{1/2}}{n^{3/2}} \cdot \left\|\mcA\right\|_{\mathrm{F}}.
$

Finally, we have
\begin{align*}
\mathrm{\RN{8}}
\lesssim & \left(\Delta \cdot \frac{\sigma_{\xi}}{\sigma} \sqrt{\frac{\op}{n}}\right)^3 \cdot\left\|\mcP_{U_{1 \perp}} A_1\left(\mcP_{U_{3 \perp}} \otimes \mcP_{U_{1 \perp}}\right)\right\|_{\mathrm{F}} = \Delta^3 \cdot \frac{\sigma_\xi^3}{\sigma^3}\cdot \frac{\op^{3/2}}{n^{3/2}} \cdot \left\|\mcA\right\|_{\mathrm{F}}.
\end{align*}

Combining the results above, we obtain the first inequality.

{\bf Step 2:} Proof of the second inequality

By symmetry, it suffices to consider 
\begin{align*}
& \left\|\mcP_{U_1} A_1 \left(\mcP_{U_{3 \perp}} \whZ_3 \left(\mcP_{U_2} \otimes \mcP_{U_1}\right) \otimes \mcP_{U_{2 \perp}} \whZ_2 \left(\mcP_{U_1} \otimes \mcP_{U_3}\right)\right) \right\|_{\mathrm{F}} \\
\leq & \underbrace{\left\|\mcP_{U_1} A_1 \left(\mcP_{U_{3 \perp}} \whZ_3^{(1)} \left(\mcP_{U_2} \otimes \mcP_{U_1}\right) \otimes \mcP_{U_{2 \perp}} \whZ_2^{(1)} \left(\mcP_{U_1} \otimes \mcP_{U_3}\right)\right)\right\|_{\mathrm{F}}}_{\mathrm{\RN{1}}} \\
+ & \underbrace{\left\|\mcP_{U_1} A_1 \left(\mcP_{U_{3 \perp}} \whZ_3^{(2)} \left(\mcP_{U_2} \otimes \mcP_{U_1}\right) \otimes \mcP_{U_{2 \perp}} \whZ_2^{(1)} \left(\mcP_{U_1} \otimes \mcP_{U_3}\right)\right)\right\|_{\mathrm{F}}}_{\mathrm{\RN{2}}} \\
+ & \underbrace{\left\|\mcP_{U_1} A_1 \left(\mcP_{U_{3 \perp}} \whZ_3^{(1)} \left(\mcP_{U_2} \otimes \mcP_{U_1}\right) \otimes \mcP_{U_{2 \perp}} \whZ_2^{(2)} \left(\mcP_{U_1} \otimes \mcP_{U_3}\right)\right)\right\|_{\mathrm{F}}}_{\mathrm{\RN{3}}} \\
+ & \underbrace{\left\|\mcP_{U_1} A_1 \left(\mcP_{U_{3 \perp}} \whZ_3^{(2)} \left(\mcP_{U_2} \otimes \mcP_{U_1}\right) \otimes \mcP_{U_{2 \perp}} \whZ_2^{(2)} \left(\mcP_{U_1} \otimes \mcP_{U_3}\right)\right)\right\|_{\mathrm{F}}}_{\mathrm{\RN{4}}}.
\end{align*}

Here, it follows from \eqref{eq: high-prob upper bound of spectral norm of A1t(A2Z2(1)B2oA3Z3(1)B3) in tensor regression} in Lemma~\ref{lemma: high-prob upper bound of spectral norm of A1t(A2Z2(1)B2oA3Z3(1)B3) in tensor regression} that
\begin{align*}
\mathrm{\RN{1}}
\lesssim \underbrace{\left\|\mcA \times_1 U_1\right\|_{\mathrm{F}} \cdot \left(\frac{\sigma_\xi^2}{\sigma^2}\cdot \frac{\oor\log(\op)}{n}\right)}_{\eqref{eq: high-prob upper bound of spectral norm of B1tZhat1(1)tA1t(A2Zhat2(1)B2oA3Zhat3(1)B3) in tensor regression without sample splitting}}.
\end{align*}
Similar to the proof of the first inequality, we can show
\begin{align*}
\mathrm{\RN{2}}
\lesssim & \Delta \cdot \frac{\sigma_\xi^2}{\sigma^2}\cdot \frac{\op^{1/2}\sqrt{\oR\log(\op)}}{n} \cdot \left\|\mcA \times_1 U_1\right\|_{\mathrm{F}}, \\
\mathrm{\RN{3}}
\lesssim & \Delta \cdot \frac{\sigma_\xi^2}{\sigma^2}\cdot \frac{\op^{1/2}\sqrt{\oR\log(\op)}}{n} \cdot \left\|\mcA \times_1 U_1\right\|_{\mathrm{F}},\\
\mathrm{\RN{4}}
\lesssim & \Delta^2 \cdot \frac{\sigma_\xi^2}{\sigma^2}\cdot \frac{\op\sqrt{\oR\log(\op)}}{n} \cdot \left\|\mcA \times_1 U_1\right\|_{\mathrm{F}}. 
\end{align*}
Combining the results above, we obtain the second inequality.

\end{proof}

\begin{lemma}
\label{lemma: high-prob upper bound of spectral norm of B1tZ1hatA1t(A2Z2hatB2oA3Z3hatB3) in tensor regression with sample splitting}
Under the same setting of Theorem~\ref{thm: main theorem in tensor regression with sample splitting}. Let $\widehat{\mcZ}=\frac{1}{n\sigma^2}\sum_{i=1}^{n}\xi\mcX_i + \frac{1}{n\sigma^2}\sum_{i=1}^n\left[\left\langle \mcX_i, \widehat{\Delta}\right\rangle \mcX_i - \sigma^2\cdot \widehat{\Delta}\right]$ be the debiased error. Let $\whZ_j = \Mat_j(\widehat{\mcZ})$ denote the mode-$j$ matricization of the tensor $\widehat{\mcZ}$. Additionally, let $U_{j\perp} \in \mathbb{O}^{p_j \times (p_j - r_j)}$ be orthonormal matrices such that $U_{j\perp} U_{j\perp}^{\top} \in \mathbb{R}^{p_j \times p_j}$ is a projection matrix that projects any vector onto the orthogonal complement of the subspace spanned by $U_j U_j^{\top}$.

Then,
\begin{align}
& \left\|\mcA \times_1 \mcP_{U_{1\perp}} \whZ_1 \left(\mcP_{U_3} \otimes \mcP_{U_2}\right) \times_2 \mcP_{U_{2\perp}} \whZ_2 \left(\mcP_{U_1} \otimes \mcP_{U_3}\right) \times_3 \mcP_{U_{3\perp}} \whZ_3 \left(\mcP_{U_2} \otimes \mcP_{U_1}\right) \right\|_{\mathrm{F}} \notag \\
\lesssim & \left\|\mcA\right\|_{\mathrm{F}} \cdot\left(\frac{\sigma_{\xi}^3}{\sigma^3} \frac{\oor^{3/2}\log(\op)^{3/2}}{n^{3 / 2}}+\frac{\sigma_{\xi}^3}{\sigma^3} \cdot \Delta \cdot \frac{\oR^{3 / 2} \log (\op)^{3 / 2}}{n^{3 / 2}} \right). \label{eq: high-prob upper bound of Ao(PU1pZhat1(PU3oPU2))o(PU2pZhat2(PU1oPU3))o(PU3pZhat3(PU2oPU1)) in tensor regression with sample splitting}
\end{align}
and
\begin{align}
& \left\|\mcA \times_1 \mcP_{U_j} \times_{j+1} \mcP_{U_{{j+1}\perp}} \whZ_{j+1} \left(\mcP_{U_j} \otimes \mcP_{U_{j+2}}\right) \times_{j+2} \mcP_{U_{{j+2}\perp}} \whZ_{j+2} \left(\mcP_{U_{j+1}} \otimes \mcP_{U_j}\right) \right\|_{\mathrm{F}} \notag \\
\lesssim & \left\|\mcA \times_j U_j\right\|_{\mathrm{F}} \cdot\left(\frac{\sigma_{\xi}^2}{\sigma^2} \cdot \frac{\oor\log(\op)}{n} + \frac{\sigma_{\xi}^2}{\sigma^2} \cdot \Delta \cdot \frac{\oR \log (\op)}{n}\right) \label{eq: high-prob upper bound of AoPU1o(PU2pZhat2(PU1oPU3))o(PU3pZhat3(PU2oPU1)) in tensor regression with sample splitting}
\end{align}
hold with probability at least $1-\exp(-cn) - \frac{1}{\op^c} - \mathbb{P}\left(\mcE_U^{\text{reg}}\right) - \mathbb{P}\left(\mcE_\Delta\right)$ for any $j=1,2,3$ .
\end{lemma}

\begin{proof}

{\bf Step 1: Proof of the first inequality}

Consider the same decomposition in the proof of \eqref{eq: high-prob upper bound of Ao(PU1pZhat1(PU3oPU2))o(PU2pZhat2(PU1oPU3))o(PU3pZhat3(PU2oPU1)) in tensor regression without sample splitting} in Lemma~\ref{lemma: high-prob upper bound of spectral norm of B1tZ1hatA1t(A2Z2hatB2oA3Z3hatB3) in tensor regression without sample splitting}. It follows from \eqref{eq: high-prob upper bound of spectral norm of B1tZhat1(1)tA1t(A2Zhat2(1)B2oA3Zhat3(1)B3) in tensor regression without sample splitting} in Lemma~\ref{lemma: high-prob upper bound of spectral norm of B1tZhat1(1)tA1t(A2Zhat2(1)B2oA3Zhat3(1)B3) in tensor regression without sample splitting} that
\begin{align}
\mathrm{\RN{1}} \lesssim \underbrace{\left\|\mcA\right\|_{\mathrm{F}} \cdot \left(\frac{\sigma_\xi^3}{\sigma^3}\cdot \frac{\oor\log(\op)^{3 / 2}}{n^{3/2}}\right)}_{\eqref{eq: high-prob upper bound of spectral norm of B1tZhat1(1)tA1t(A2Zhat2(1)B2oA3Zhat3(1)B3) in tensor regression without sample splitting}}. \label{eq: upper bound of term 1 in Ao(PU1pZhat1(PU3oPU2))o(PU2pZhat2(PU1oPU3))o(PU3pZhat3(PU2oPU1)) in tensor regression with sample splitting}
\end{align}
Similar to the proof of Lemma~\ref{lemma: high-prob upper bound of spectral norm of B1tZ1hatA1t(A2Z2hatB2oA3Z3hatB3) in tensor regression without sample splitting}, we can show
\begin{align}
\mathrm{\RN{2}} + \mathrm{\RN{3}} + \mathrm{\RN{4}}
\lesssim & \Delta \cdot \frac{\sigma_\xi^3}{\sigma^3}\cdot \frac{\oR^{3/2}\log(\op)^{3/2}}{n^{3/2}} \cdot \left\|\mcA\right\|_{\mathrm{F}},  \label{eq: upper bound of term 2, 3, 4 in Ao(PU1pZhat1(PU3oPU2))o(PU2pZhat2(PU1oPU3))o(PU3pZhat3(PU2oPU1)) in tensor regression with sample splitting}  \\
\mathrm{\RN{5}} + \mathrm{\RN{6}} + \mathrm{\RN{7}} 
\lesssim & \Delta^2 \cdot \frac{\sigma_\xi^3}{\sigma^3}\cdot \frac{\oR^{3/2}\log(\op)^{3/2}}{n^{3/2}} \cdot \left\|\mcA\right\|_{\mathrm{F}}, \label{eq: upper bound of term 5, 6, 7 in Ao(PU1pZhat1(PU3oPU2))o(PU2pZhat2(PU1oPU3))o(PU3pZhat3(PU2oPU1)) in tensor regression with sample splitting} \\
\mathrm{\RN{8}}
\lesssim & \Delta^3 \cdot \frac{\sigma_\xi^3}{\sigma^3}\cdot \frac{\oR^{3/2}\log(\op)^{3/2}}{n^{3/2}} \cdot \left\|\mcA\right\|_{\mathrm{F}}. \label{eq: upper bound of term 8 in Ao(PU1pZhat1(PU3oPU2))o(PU2pZhat2(PU1oPU3))o(PU3pZhat3(PU2oPU1)) in tensor regression with sample splitting}
\end{align}
Combining the results above, we obtain the first inequality. 

{\bf Step 2: Proof of the second inequality}

Consider the same decomposition in the proof of \eqref{eq: high-prob upper bound of AoPU1o(PU2pZhat2(PU1oPU3))o(PU3pZhat3(PU2oPU1)) in tensor regression without sample splitting} in Lemma~\ref{lemma: high-prob upper bound of spectral norm of B1tZ1hatA1t(A2Z2hatB2oA3Z3hatB3) in tensor regression without sample splitting}. It follows from \eqref{eq: high-prob upper bound of spectral norm of A1t(A2Z2(1)B2oA3Z3(1)B3) in tensor regression} in Lemma~\ref{lemma: high-prob upper bound of spectral norm of A1t(A2Z2(1)B2oA3Z3(1)B3) in tensor regression} that
\begin{align}
\mathrm{\RN{1}}
\leq & \underbrace{\left\|\mcA \times_1 U_1\right\|_{\mathrm{F}} \cdot \left(\frac{\sigma_\xi^2}{\sigma^2}\cdot \frac{\oor\log(\op)}{n}\right)}_{\eqref{eq: high-prob upper bound of spectral norm of B1tZhat1(1)tA1t(A2Zhat2(1)B2oA3Zhat3(1)B3) in tensor regression without sample splitting}}. \label{eq: upper bound of term 1 in AoPU1o(PU2pZhat2(PU1oPU3))o(PU3pZhat3(PU2oPU1)) in tensor regression with sample splitting}
\end{align}
Similar to the proof of Lemma~\ref{lemma: high-prob upper bound of spectral norm of B1tZ1hatA1t(A2Z2hatB2oA3Z3hatB3) in tensor regression without sample splitting}, we can show
\begin{align}
\mathrm{\RN{2}} + \mathrm{\RN{3}}
\lesssim & \Delta \cdot \frac{\sigma_\xi^2}{\sigma^2}\cdot \frac{\oR\log(\op)}{n} \cdot \left\|\mcA \times_1 U_1\right\|_{\mathrm{F}},  \label{eq: upper bound of term 2, 3 in AoPU1o(PU2pZhat2(PU1oPU3))o(PU3pZhat3(PU2oPU1)) in tensor regression with sample splitting} \\
\mathrm{\RN{4}}
\lesssim & \Delta^2 \cdot \frac{\sigma_\xi^2}{\sigma^2}\cdot \frac{\oR\log(\op)}{n} \cdot \left\|\mcA \times_1 U_1\right\|_{\mathrm{F}}. \label{eq: upper bound of term 4 in AoPU1o(PU2pZhat2(PU1oPU3))o(PU3pZhat3(PU2oPU1)) in tensor regression with sample splitting} 
\end{align}
Combining the results above, we obtain the second inequality. 


\end{proof}

\begin{lemma}\label{lemma: high-prob upper bound of tr(BZ(1)tCZ(1))}
Suppose that $B \in \mathbb{R}^{p_{j+2}p_{j+1} \times p_{j+2}p_{j+1}}$ and $C \in \mathbb{R}^{p_j \times p_j}$ are two fixed matrices. Let $\whZ_j^{(1)}=\frac{1}{n\sigma^2}\sum_{i=1}^n\xi_i\Mat_j\left(\mcX_i\right)$, where $\left\{\xi_i\right\}_{i=1}^n$'s are i.i.d. mean zero $\sigma_{\xi}$-sub-Gaussian variables and $\mcX_i \in \mathbb{R}^{p_1\times p_2\times p_3}$, $i=1,2,\cdots,n$ are i.i.d. random tensors with i.i.d. mean zero $\sigma$-sub-Gaussian entries. Then, it holds that
\begin{align}
& \mathbb{P}\left( \left|\operatorname{tr}\left(B\whZ_j^{(1)\top}C\whZ_j^{(1)}\right)-\frac{1}{n\sigma^2}\operatorname{tr}\left[B\right]\operatorname{tr}\left[C\right]\right| \geq \frac{c_1\sigma_\xi^2}{n\sigma^2}\cdot\left(\left\|B\right\|_{\mathrm{F}}\left\|C\right\|_{\mathrm{F}}\sqrt{t}+ \left\|B\right\|_{\ell_\infty}\left\|C\right\|_{\ell_\infty}t\right)\right) \notag  \\
\leq & \exp \left(-c t\right) + \exp\left(-cn\right), \label{eq: high-prob upper bound of tr(BZ(1)tCZ(1))}
\end{align}
where $c_1$ and $c_2$ are two universal constants.

\end{lemma}

\begin{proof}
By symmetry, it suffice to consider $\operatorname{tr}\left[B\whZ_1^{(1)\top}C\whZ_1^{(1)}\right]$. Conditioning on $\left\{\xi_i\right\}_{i=1}^n$, then by Lemma~\ref{lemma: high-prob upper bound of tr(BZtCZ)}, we have
\begin{align*}
& \mathbb{P}\left( \left|\operatorname{tr}\left(B\whZ_1^{(1)\top}C\whZ_1^{(1)}\right)-\frac{1}{n^2\sigma^2}\operatorname{tr}\left[B\right]\operatorname{tr}\left[C\right]\left\|\xi\right\|_{\ell_2}^2\right| \geq C\cdot \frac{t}{n^2\sigma^2} \middle| \left\{\xi_i\right\}_{i=1}^n\right) \\
\leq & \exp \left(-c \min \left(\frac{t^2}{\left\|\xi\right\|_{\ell_2}^4\left\|B\right\|_{\mathrm{F}}^2\left\|C\right\|_{\mathrm{F}}^2}, \frac{t}{ \left\|\xi\right\|_{\ell_\infty}^2\left\|B\right\|_{\ell_\infty}\left\|C\right\|_{\ell_\infty}}\right)\right).
\end{align*}
That is
\begin{align*}
& \mathbb{P}\left( \left|\operatorname{tr}\left(B\whZ_1^{(1)\top}C\whZ_1^{(1)}\right)-\frac{1}{n^2\sigma^2}\operatorname{tr}\left[B\right]\operatorname{tr}\left[C\right]\left\|\xi\right\|_{\ell_2}^2\right| \right. \\
& \quad \geq \left. \frac{C}{n^2\sigma^2}\cdot\left(\left\|\xi\right\|_{\ell_2}^2\left\|B\right\|_{\mathrm{F}}\left\|C\right\|_{\mathrm{F}}\sqrt{t}+ \left\|\xi\right\|_{\ell_\infty}^2\left\|B\right\|_{\ell_\infty}\left\|C\right\|_{\ell_\infty}t\right) \middle| \left\{\xi_i\right\}_{i=1}^n\right) \leq \exp \left(-ct\right).
\end{align*}
Then note that $\mathbb{P}\left(\left\|\xi\right\|_{\ell_\infty}\geq \sigma_\xi\sqrt{n}\right)\leq  \mathbb{P}\left(\left\|\xi\right\|_{\ell} \geq \sigma_\xi\sqrt{n}\right)\leq \exp(-cn)$. It implies that
\begin{align*}
& \mathbb{P}\left( \left|\operatorname{tr}\left(B\whZ_1^{(1)\top}C\whZ_1^{(1)}\right)-\frac{1}{n^2\sigma^2}\operatorname{tr}\left[B\right]\operatorname{tr}\left[C\right]\right| \leq \frac{C_1\sigma_\xi^2}{n\sigma^2}\cdot\left(\left\|B\right\|_{\mathrm{F}}\left\|C\right\|_{\mathrm{F}}\sqrt{t}+ \left\|B\right\|_{\ell_\infty}\left\|C\right\|_{\ell_\infty}t\right)\right) \\
\geq & \left(1-\exp \left(-ct\right)\right)\cdot \left(1-\exp\left(-cn\right)\right) \cdot \left(1-\exp\left(-cn\right)\right) \geq 1- \left[\exp \left(-ct\right) + \exp\left(-cn\right)\right],
\end{align*}
where the second inequality follows as long as $\exp \left(-ct\right)\leq 1, \exp\left(-cn\right)$ and $n\geq 1$.

\end{proof}

\begin{lemma}\label{lemma: high-prob upper bound of spectral norm of W_1tZ(1)t(BZ(1)W2oW3)}
Suppose that $W_{j+1} \in \mathbb{R}^{p_{j+1} \times R_{j+1}}, W_{j+2} \in \mathbb{R}^{p_{j+2} \times R_{j+2}}$, $B \in \mathbb{R}^{p_j\times p_j}$, and $\widetilde{W}_j \in \mathbb{R}^{p_{j+1}p_{j+2}\times R_j}$ are fixed matrices, where $R_j\leq p_j$ for any $j=1,2,3$. Let $\whZ_j^{(1)}=\frac{1}{n\sigma^2}\sum_{i=1}^n\xi_i\Mat_j\left(\mcX_i\right)$ are arbitrary given matrices, where $\left\{\xi_i\right\}_{i=1}^n$'s are i.i.d. mean zero $\sigma_{\xi}$-sub-Gaussian variables and $\mcX_i \in \mathbb{R}^{p_1\times p_2\times p_3}$, $i=1,2,\cdots,n$ are i.i.d. random tensors with i.i.d. mean zero $\sigma$-sub-Gaussian entries. 
Furthermore, suppose that $\widetilde{W}_2^{\top}\left(W_1\otimes W_3\right)=0$. Then, it holds that
\begin{align}
& \mathbb{P}\left(\left\|\widetilde{W}_j^{\top} \widehat{Z}_j^{(\mathrm{\RN{1}}) \top} \widetilde{B}^{\top} \widehat{Z}_j^{(\mathrm{\RN{1}})}\left(W_{j+1} \otimes W_{j+2}\right)\right\| \geq \frac{C_1\sigma_\xi^2}{n\sigma^2}\left(\|\widetilde{B}\|_{\mathrm{F}} \|\widetilde{W}_j \| \|W_{j+1} \| \|W_{j+2} \|\sqrt{t}+  \|\widetilde{B}\|_{\ell_\infty} \|\widetilde{W}_j \| \|W_{j+1} \| \|W_{j+2} \|t\right)\right) \notag \\
\leq & 9 \cdot 7^{R_j+R_{j+1}R_{j+2}} \cdot \exp \left(-ct\right) + \exp\left(-cn\right), \label{eq: high-prob upper bound of |W2tZ(1)2tBZ(1)2(W3oW1)| in spectral norm}
\end{align}
where $C$ and $c$ are two universal constants.
\end{lemma}

\begin{proof}

Conditioning on $\left\{\xi_i\right\}_{i=1}^n$, then by Lemma~\ref{lemma: high-prob upper bound of |W2tZ2tBZ2(W3oW1)|}, we have
\begin{align*}
& \mathbb{P}\left(\left\|\widetilde{W}_2^{\top} Z_2^{\top} \widetilde{B}^{\top} Z_2\left(W_3 \otimes W_1\right)\right\| \geq C \cdot \frac{t}{n^2\sigma^2} \middle| \left\{\xi\right\}_{i=1}^n\right) \\
&  \leq 9 \cdot 7^{R_2+R_3R_1} \cdot  \exp \left(-c \min \left(\frac{t^2}{\left\|\xi\right\|_{\ell_2}^4\left\|\widetilde{B}\right\|_{\mathrm{F}}^2\left\|W_1\right\|^2\left\|\widetilde{W}_2\right\|^2\left\|W_3\right\|^2}, \frac{t}{\left\|\xi\right\|_{\ell_2}^2\left\|\widetilde{B}\right\|_{\ell_{\infty}}\left\|\widetilde{W}_2\right\|\left\|W_3\right\|\left\|W_1\right\|}\right)\right).
\end{align*}
That is
\begin{align*}
& \mathbb{P}\left(\left\|\widetilde{W}_2^{\top} Z_2^{\top} \widetilde{B}^{\top} Z_2\left(W_3 \otimes W_1\right)\right\| \right. \\
& \left. \geq \frac{C}{n^2\sigma^2}\cdot\left(\left\|\xi\right\|_{\ell_2}^2\left\|\widetilde{B}\right\|_{\mathrm{F}}\left\|W_1\right\|\left\|\widetilde{W}_2\right\|\left\|W_3\right\|\sqrt{t}+ \left\|\xi\right\|_{\ell\infty}^2\left\|\widetilde{B}\right\|_{\ell_\infty}\left\|W_1\right\|\left\|\widetilde{W}_2\right\|\left\|W_3\right\|t\right) \middle| \left\{\xi_i\right\}_{i=1}^n\right) \leq 9 \cdot 7^{R_2+R_3R_1} \cdot \exp \left(-ct\right).
\end{align*}
Then note that $\mathbb{P}\left(\left\|\xi\right\|_{\ell_\infty}\geq \sigma_\xi\sqrt{n}\right)\leq  \mathbb{P}\left(\left\|\xi\right\|_{\ell} \geq \sigma_\xi\sqrt{n}\right)\leq \exp(-cn)$. It implies that
\begin{align*}
& \mathbb{P}\left( \left\|\widetilde{W}_2^{\top} Z_2^{\top} \widetilde{B}^{\top} Z_2\left(W_3 \otimes W_1\right)\right\| \leq \frac{C_1\sigma_\xi^2}{n\sigma^2}\cdot\left(\left\|\widetilde{B}\right\|_{\mathrm{F}}\left\|W_1\right\|\left\|\widetilde{W}_2\right\|\left\|W_3\right\|\sqrt{t}+ \left\|\widetilde{B}\right\|_{\ell_\infty}\left\|W_1\right\|\left\|\widetilde{W}_2\right\|\left\|W_3\right\|t\right)\right) \\
\geq & \left(1- 9 \cdot 7^{R_2+R_3R_1} \cdot \exp \left(-ct\right)\right)\cdot \left(1-\exp\left(-cn\right)\right)  \geq 1- \left[9 \cdot 7^{R_2+R_3R_1} \cdot \exp \left(-ct\right) + \exp\left(-cn\right)\right],
\end{align*}
where the second inequality follows as long as $\exp \left(-ct\right)\leq 1, \exp\left(-cn\right)$ and $n\geq 1$.


\end{proof}

The following Lemma~\ref{lemma: preliminary l2-norm perturbation bound of <X,Delta>U1X(U3oU2) in tensor regression}, \ref{lemma: preliminary l2-norm perturbation bound of sup_Delta <X,Delta>U1X(U3oU2) in tensor regression} established upper bounds for the spectral norm of matricization of projected sub-Gaussian random tensor.
By Lemma~\ref{lemma: high-prob upper bound of spectral norm of sumU1Xi1(U3oU2)} and \ref{lemma: preliminary l2-norm perturbation bound of <X,Delta>U1X(U3oU2) in tensor regression}, it follows immediately that
\begin{align}
& \left\|U_j^{\top}\widehat{Z}_j\left(U_{j+2} \otimes U_{j+1}\right)\right\| \lesssim \frac{\sigma_\xi}{\sigma}\left(\sqrt{\frac{\oor\log(\op)}{n}} + \Delta\sqrt{\frac{\op}{n}}\right) \label{eq: high-prob upper bound of U1Zhat1(U3oU2) in tensor regression without sample splitting}
\end{align}
with probability at least $1-\op^{-c} - \exp(-cn) - \mathbb{P}\left(\mcE_{\Delta}\right) - \mathbb{P}\left(\mcE_{U}^{\text{reg}}\right)$, and
\begin{align}
& \left\|U_{j\perp}^{\top}\widehat{Z}_j\left(U_{j+2} \otimes U_{j+1}\right)\right\| \lesssim \frac{\sigma_\xi}{\sigma} \sqrt{\frac{\op}{n}} \label{eq: high-prob upper bound of U1pZhat1(U3oU2) in tensor regression without sample splitting}
\end{align}
with probability at least $1-\exp(-c\op) - \exp(-cn) - \mathbb{P}\left(\mcE_{\Delta}\right) - \mathbb{P}\left(\mcE_{U}^{\text{reg}}\right)$.

\begin{lemma}\label{lemma: preliminary l2-norm perturbation bound of <X,Delta>U1X(U3oU2) in tensor regression}
Suppose that $\mcX \in \mathbb{R}^{p_1 \times p_2 \times p_3}$ has i.i.d. mean zero $\sigma$-sub-Gaussian entries, and $\mcX_1, \cdots, \mcX_n$ are i.i.d. copies of $\mcX$. Then there exist two universal constants $C, c>0$ such that for any fixed orthonormal matrix $\wtU_1 \in \mathbb{O}^{p_1 \times \widetilde{r}_1}, \wtU_2 \in \mathbb{O}^{p_2 \times \widetilde{r}_2}, \wtU_3 \in \mathbb{O}^{p_3 \times \widetilde{r}_3}$, where $\widetilde{r}_1 \leq p_1, \widetilde{r}_2 \leq p_2, \widetilde{r}_3 \leq p_3$, and a fixed tensor $\widetilde{\Delta} \in \mathbb{R}^{p_1 \times p_2 \times p_3}$ satisfying $\left\|\widetilde{\Delta}\right\|_{\mathrm{F}}\leq \Delta $ with probability at least $1-\mathbb{P}\left(\mcE_{\Delta}\right)$, it holds that
\begin{align}
& \mathbb{P}\left(\left\|\frac{1}{n \sigma^2} \sum_{i=1}^n\left\langle \mcX_i, \widetilde{\Delta}\right\rangle \wtU_j^{\top} \Mat_j (\mcX_i) \left(\wtU_{j+2} \otimes \wtU_{j+1}\right) -  \wtU_j^{\top} \widetilde{\Delta}_j \left(\wtU_{j+2}  \otimes \wtU_{j+1}\right)\right\| \geq C\sigma^2\cdot \Delta t\right) \notag \\
\leq & 2 \cdot 7^{\widetilde{r}_1+\widetilde{r}_2 \widetilde{r}_3} e^{-c \min \left\{n t^2, n t\right\}} + \mathbb{P}\left(\mcE_{\Delta}\right). \label{eq: high-prob upper bound of <X,Delta>U1X(U3oU2) in tensor regression}
\end{align}

\end{lemma}

\begin{proof}
By symmetry, it suffices to consider
$$
\left\|\frac{1}{n \sigma^2} \sum_{i=1}^n \left[\left\langle \mcX_i, \widetilde{\Delta}\right\rangle \wtU_1^{\top} \Mat_1 (\mcX_i) \left(\wtU_3 \otimes \wtU_2\right) - \sigma^2 \cdot \wtU_1^{\top} \Mat_1 \left(\widetilde{\Delta}\right) \left(\wtU_3  \otimes \wtU_2\right)\right]\right\|.
$$

For any fixed $a \in \mathbb{R}^{\widetilde{r}_1}, \left\|a\right\|=1$ and $b \in \mathbb{R}^{\widetilde{r}_2\widetilde{r}_3}, \left\|b\right\|=1$, we have
\begin{align*}
& \mathbb{E}\left[\left\langle \mcX_i, \widetilde{\Delta}\right\rangle a^{\top} \wtU_1^{\top} \Mat_1\left(\mcX_i\right) \left(\wtU_3 \otimes \wtU_2\right) b\right] = \sigma^2\cdot  a^{\top} \wtU_1^{\top} \widetilde{\Delta}_1 \left(\wtU_3 \otimes \wtU_2\right) b, \quad \forall i \in[n]  .
\end{align*}

Since $\mathbb{E}(\langle \mcX_i, \widetilde{\Delta}\rangle^2) = \mathbb{E}_{\widetilde{\Delta}} [\mathbb{E}_{\mcX}(\langle \mcX_i, \widetilde{\Delta}\rangle^2\mid \widetilde{\Delta})] = \sigma^2\cdot \mathbb{E}_{\widetilde{\Delta}}(\|\widetilde{\Delta}\|_{\mathrm{F}}^2) \leq \sigma^2\Delta^2$ 
and 
\begin{align*}
& \mathbb{E} \left[\left(a^{\top} \wtU_1^{\top} \Mat_1\left(\mcX_i\right) \left(\wtU_3 \otimes \wtU_2\right) b\right)^2 \right] = \mathbb{E} \left\{ \left[\left(b^{\top}\left(\wtU_3 \otimes \wtU_2\right)^{\top}\right)\otimes \left(a^{\top}\wtU_1^{\top}\right)\right]\operatorname{Vec}(\mcX_i)\right\}^2 =\sigma^2,
\end{align*}
where $\left[\left(b^{\top}\left(\wtU_3 \otimes \wtU_2\right)^{\top}\right)\otimes \left(a^{\top}\wtU_1^{\top}\right)\right]\left[\left(b^{\top}\left(\wtU_3 \otimes \wtU_2\right)^{\top}\right)\otimes \left(a^{\top}\wtU_1^{\top}\right)\right]^{\top}=1$, we have
\begin{align*}
& \left\|\left\langle \mcX_i, \widetilde{\Delta}\right\rangle a^{\top} \wtU_1^{\top} \Mat_1\left(\mcX_i\right) \left(\wtU_3 \otimes \wtU_2\right)  b-\sigma^2 \cdot a^{\top} \wtU_1^{\top}\widetilde{\Delta}_1 \left(\wtU_3 \otimes \wtU_2\right) b\right\|_{\psi_1} \\
\leq & 2\left\|\left\langle \mcX_i, \widetilde{\Delta}\right\rangle\right\|_{\psi_2}\left\|a^{\top} \wtU_1^{\top} \Mat_1\left(\mcX_i\right) \left(\wtU_3 \otimes \wtU_2\right) b\right\|_{\psi_2} \leq 2 \sigma^2 \cdot \|\widetilde{\Delta}\|_{\mathrm{F}} \leq 2\sigma^2 \cdot \Delta ,
\end{align*}
where $\|\cdot\|_{\psi_1}$ denotes Orlicz $\psi_1-$norm.
The first inequality follows from t Remark 5.18 of \citet{vershynin2010introduction} and the second inequality follows from the Cauchy-Schwarz inequality for the Orlicz norm. 

Therefore, by Bernstein's inequality, we have
\begin{align*}
& \mathbb{P}\left(\left|\sum_{i=1}^n \frac{1}{n\sigma^2}\left[\left\langle \mcX_i, \widetilde{\Delta}\right\rangle a^{\top} \wtU_1^{\top} \Mat_1\left(\mcX_i\right) \left(\wtU_3 \otimes \wtU_2\right) b- \sigma^2\cdot a^{\top} \wtU_1^{\top}\widetilde{\Delta}_1 \left(\wtU_3 \otimes \wtU_2\right) b\right]\right| \geq C\sigma^2 \cdot \widetilde{\Delta} t\right) \\
\leq & 2 \exp \left[-c \min \left\{n t^2, n t\right\}\right] + \mcP\left(\mcE_{\Delta}\right) .
\end{align*}
By Lemma 5.2 of \citet{vershynin2010introduction}, there exists a $\frac{1}{3}$-net $\mcN_{\widetilde{r}_1}$ for $\mathbb{S}^{\widetilde{r}_1}=\left\{x: x \in \mathbb{R}^{\widetilde{r}_1},\left\|x\right\|=1\right\}$ with cardinality at most $7^{\widetilde{r}_1}$ and a $\frac{1}{3}$-net $\mcN_{\widetilde{r}_2\widetilde{r}_3}$ for $\mathbb{S}^{\widetilde{r}_2\widetilde{r}_3}=\left\{x: x \in \mathbb{R}^{\widetilde{r}_2\widetilde{r}_3},\|x\|_2=1\right\}$ with cardinality at most $7^{\widetilde{r}_2\widetilde{r}_3}$. By the union bound, we have 
\begin{align*}
& \mathbb{P}\left(\sup _{a \in \mcN_{\widetilde{r}_1}, b \in \mcN_{\widetilde{r}_2\widetilde{r}_3}}\left| a^{\top}\sum_{i=1}^n \frac{1}{n\sigma^2}\left[\left\langle \mcX_i, \widetilde{\Delta}\right\rangle \wtU_1^{\top} \Mat_1\left(\mcX_i\right) \left(\wtU_3 \otimes \wtU_2\right)- \sigma^2 \wtU_1^{\top}\widetilde{\Delta}_1 \left(\wtU_3 \otimes \wtU_2\right) \right] b\right| \geq C\sigma^2 \Delta t\right) \\
&  \leq 2 \cdot 7^{\widetilde{r}_1+\widetilde{r}_2\widetilde{r}_3} \exp\left[{-c \min \left\{n t^2, n t\right\}}\right] .
\end{align*}
Then for any $a \in \mathbb{S}^{\widetilde{r}_1}$ and $b \in \mathbb{S}^{\widetilde{r}_2\widetilde{r}_3}$ there exist $\overline{a} \in \mcN_{\widetilde{r}_1}$ and $\overline{b} \in \mcN_{\widetilde{r}_2\widetilde{r}_3}$ such that $\left\|\widetilde{a}-\overline{a}\right\| \leq \frac{1}{3}$ and $\left\|\widetilde{b}-\overline{b}\right\| \leq \frac{1}{3}$. Therefore, by the same $\varepsilon$-net arguments as in the proof of Lemma~\ref{lemma: high-prob upper bound of spectral norm of B1tZ1(1)tA1t(A2Z2(1)B2oA3Z3(1)B3) in tensor PCA}, we have
\begin{align*}
& \mathbb{P}\left(\left\|\sum_{i=1}^{n} \frac{1}{n\sigma^2}\left(\left\langle \mcX_i, \widetilde{\Delta}\right\rangle \wtU_1^{\top} \operatorname{Mat}_1\left(\mcX_i\right) \left(\wtU_3 \otimes \wtU_2\right) - \sigma^2  \wtU_1^{\top} \widetilde{\Delta}_1 \left(\wtU_3 \otimes \wtU_2\right)\right)\right\| \geq C \sigma^2  \Delta t\right) \\
\leq & 1- \mathbb{P}\left(\sup_{a,b}\left|\sum_{i=1}^n \frac{1}{n\sigma^2}\left[\langle \mcX_i, \widetilde{\Delta}\rangle a^{\top} \wtU_1^{\top} (\Mat_1(\mcX_i) - \sigma^2 \widetilde{\Delta}_1) (\wtU_3 \otimes \wtU_2) b\right]\right| \leq C\sigma^2 \|\widetilde{\Delta}\|_{\mathrm{F}} t \middle| \widetilde{\Delta}\right)  \mathbb{P}\left(\|\widetilde{\Delta}\|_{\mathrm{F}} \leq \Delta\right) \\
\leq & 2 \cdot 7^{\widetilde{r}_1+\widetilde{r}_2\widetilde{r}_3} \exp\left(-c \min \left\{n t^2, n t\right\}\right) + \mcP\left(\mcE_{\Delta}\right).
\end{align*}
\end{proof}

\begin{lemma}\label{lemma: preliminary l2-norm perturbation bound of sup_Delta <X,Delta>U1X(U3oU2) in tensor regression}
Suppose $\mcX \in \mathbb{R}^{p_1 \times p_2 \times p_3}$ with i.i.d. mean zero sub-Gaussian entries, and $\mcX_1, \cdots, \mcX_n$ are i.i.d. copies of $X$. Then there exist two universal constants $C, c_1>0$ such that for any fixed $\wtU_1 \in \mathbb{O}^{p_1 \times \widetilde{r}_1}, \wtU_2 \in \mathbb{O}^{p_2 \times \widetilde{r}_2}, \wtU_3 \in \mathbb{O}^{p_3 \times \widetilde{r}_3}$ and $\widetilde{\widetilde{\Delta}} \in \mathbb{R}^{p_1 \times p_2 \times p_3}$ satisfying $\operatorname{rank}(\widetilde{\widetilde{\Delta}})\leq (\widetilde{\widetilde{r}}_1, \widetilde{\widetilde{r}}_2, \widetilde{\widetilde{r}}_3)$, where $4\widetilde{\widetilde{r}}_1 \leq p_1, 4\widetilde{\widetilde{r}}_2 \leq p_2, 4\widetilde{\widetilde{r}}_3 \leq p_3$,
\begin{align*}
& \mathbb{P}\left( \sup_{\substack{\widetilde{\widetilde{\Delta}} \in \mathbb{R}^{p_1\times p_2\times p_3}, \\\|\widetilde{\widetilde\Delta}\|_{\mathrm{F}} \leq \Delta \\ \operatorname{rank}(\widetilde{\widetilde{\Delta}})\leq (\widetilde{\widetilde{r}}_1, \widetilde{\widetilde{r}}_2, \widetilde{\widetilde{r}}_3)}}\left\|\frac{1}{n\sigma^2}\sum_{i=1}^n \left[\left\langle \mcX_i, \widetilde{\widetilde{\Delta}}\right\rangle \wtU_j^{\top}\Mat_j\left(\mcX_i\right)\left(\wtU_{j+2} \otimes \wtU_{j+1}\right)-  \sigma^2\cdot \wtU_j^{\top}\widetilde{\widetilde{\Delta}}_j\left(\wtU_{j+2} \otimes \wtU_{j+1}\right)\right]\right\| \geq C\sigma^2 \Delta t\right) \\
\leq & 2\cdot 17^{\widetilde{\widetilde{r}}_1\widetilde{\widetilde{r}}_2\widetilde{\widetilde{r}}_3+\sum_{j=1}^3p_j\widetilde{\widetilde{r}}_j} \cdot 7^{\widetilde{r}_j+\widetilde{r}_{j+1}\widetilde{r}_{j+2}} \exp\left[-c\min\left\{nt^2, nt\right\}\right].
\end{align*}
Here, we denote $\widetilde{\widetilde \Delta}_j=\operatorname{Mat}_j(\widetilde{\widetilde{\Delta}})$ for any $\widetilde{\widetilde{\Delta}} \in \mathbb{R}^{p_1\times p_2\times p_3}$. 

\end{lemma}

\begin{proof}

By symmetry, it suffices to consider a high-probability upper bound for 
$$
\left\|\frac{1}{n\sigma^2}\sum_{i=1}^n \left[\left\langle \mcX_i, \widetilde{\widetilde{\Delta}}\right\rangle \wtU_1^{\top}\Mat_1\left(\mcX_i\right)\left(\wtU_3 \otimes \wtU_2\right)-  \sigma^2\cdot \wtU_1^{\top}\widetilde{\widetilde{\Delta}}_1\left(\wtU_3 \otimes \wtU_2\right)\right]\right\|.
$$
By Lemma 11 in \citet{xia2022inference}, for the class of low-Tucker-rank tensors under the Frobenius norm, defined as 
$$
\mathcal{F} \left(\left(p_1,p_2,p_3\right), \left(\widetilde{\widetilde{r}}_1, \widetilde{\widetilde{r}}_2, \widetilde{\widetilde{r}}_3\right)\right) = \left\{\widetilde{\widetilde{\Delta}} \in \mathbb{R}^{p_1 \times p_2 \times p_3}: \operatorname{rank}\left(\widetilde{\widetilde{\Delta}}_j\right) \leq \widetilde{\widetilde{r}}_j, \, j=1,2,3, \, \|\widetilde{\widetilde{\Delta}}\|_{\mathrm{F}} \leq 1\right\},
$$
there exists an $\varepsilon$-net $\mathcal{N}_{(p_1,p_2,p_3), (\widetilde{\widetilde{r}}_1,\widetilde{\widetilde{r}}_2,\widetilde{\widetilde{r}}_3)}$ with elements $\widetilde{\widetilde{\Delta}}^{(1)}, \ldots, \widetilde{\widetilde{\Delta}}^{(N)}$ and cardinality satisfying
$
\left|\mathcal{N}_{(p_1,p_2,p_3), (\widetilde{\widetilde{r}}_1,\widetilde{\widetilde{r}}_2,\widetilde{\widetilde{r}}_3)}\right| \leq \left(\frac{8+\varepsilon}{\varepsilon}\right)^{\widetilde{\widetilde{r}}_1 \widetilde{\widetilde{r}}_2 \widetilde{\widetilde{r}}_3+\sum_{j=1}^3 p_j \widetilde{\widetilde{r}}_j}.
$
For each $\widetilde{\widetilde{\Delta}}^{(i)} \in \mathbb{R}^{p_1 \times p_2 \times p_3}$, we have $\|\widetilde{\widetilde{\Delta}}^{(i)}\|_{\mathrm{F}} \leq 1$. Moreover, for any $\widetilde{\widetilde{\Delta}} \in \mathcal{F} \left(\left(p_1,p_2,p_3\right), \left(\widetilde{\widetilde{r}}_1, \widetilde{\widetilde{r}}_2, \widetilde{\widetilde{r}}_3\right)\right)$, there exists an index $i \in \mathcal{N}_{(p_1,p_2,p_3), (\widetilde{\widetilde{r}}_1,\widetilde{\widetilde{r}}_2,\widetilde{\widetilde{r}}_3)}$ such that
$\|\widetilde{\widetilde{\Delta}}^{(i)} - \widetilde{\widetilde{\Delta}}\|_{\mathrm{F}} \leq \varepsilon$.  Therefore, it follows that 
\begin{align*}
& \left\|\frac{1}{n\sigma^2}\sum_{i=1}^n \left[\left\langle \mcX_i, \widetilde{\widetilde{\Delta}}\right\rangle \wtU_1^{\top}\Mat_1\left(\mcX_i\right)\left(\wtU_3 \otimes \wtU_2\right)-  \sigma^2\cdot \wtU_1^{\top}\widetilde{\widetilde{\Delta}}_1\left(\wtU_3 \otimes \wtU_2\right)\right]\right\| \\
\leq & \varepsilon \cdot \sup_{\widetilde{\widetilde{\Delta}} \in \mathcal{F} ((p_1,p_2,p_3),(2\widetilde{\widetilde{r}}_1, 2\widetilde{\widetilde{r}}_2, 2\widetilde{\widetilde{r}}_3) ) }\left\|\frac{1}{n\sigma^2}\sum_{i=1}^n \left[\left\langle \mcX_i, \widetilde{\widetilde{\Delta}}\right\rangle \wtU_1^{\top}\Mat_1\left(\mcX_i\right)\left(\wtU_3 \otimes \wtU_2\right)-  \sigma^2\cdot \wtU_1^{\top}\widetilde{\widetilde{\Delta}}_1\left(\wtU_3 \otimes \wtU_2\right)\right]\right\| \\
& +  \sup_{\widetilde{\widetilde{\Delta}} \in \mcN_{(p_1,p_2,p_3),(\widetilde{\widetilde{r}}_1, \widetilde{\widetilde{r}}_2, \widetilde{\widetilde{r}}_3)} }\left\|\frac{1}{n\sigma^2}\sum_{i=1}^n \left[\left\langle \mcX_i, \widetilde{\widetilde{\Delta}}\right\rangle \wtU_1^{\top}\Mat_1\left(\mcX_i\right)\left(\wtU_3 \otimes \wtU_2\right)-  \sigma^2\cdot \wtU_1^{\top}\widetilde{\widetilde{\Delta}}_1\left(\wtU_3 \otimes \wtU_2\right)\right]\right\| .
\end{align*}
Taking the supremum over $\mathcal{F} ((p_1,p_2,p_3),(2\widetilde{\widetilde{r}}_1, 2\widetilde{\widetilde{r}}_2, 2\widetilde{\widetilde{r}}_3) )$, it follows that 
\begin{align*}
& \mathbb{P}\left( \sup_{\Delta \in \mathcal{F} ((p_1,p_2,p_3),(2\widetilde{\widetilde{r}}_1, 2\widetilde{\widetilde{r}}_2, 2\widetilde{\widetilde{r}}_3) ) }\left\|\frac{1}{n\sigma^2}\sum_{i=1}^n \left\langle \mcX_i, \Delta\right\rangle \wtU_1^{\top}\Mat_1\left(\mcX_i\right)\left(\wtU_3 \otimes \wtU_2\right)-  \wtU_1^{\top}\Delta_1\left(\wtU_3 \otimes \wtU_2\right)\right\| \geq C^{\prime}\sigma^2\cdot \frac{ \Delta }{1-\varepsilon}t\right) \\
\leq & 2\cdot \left(\frac{8+\varepsilon}{\varepsilon}\right)^{\widetilde{\widetilde{r}}_1\widetilde{\widetilde{r}}_2\widetilde{\widetilde{r}}_3+\sum_{j=1}^3p_j\widetilde{\widetilde{r}}_j} \cdot 7^{\widetilde{r}_1+\widetilde{r}_2\widetilde{r}_3}\cdot e^{-c_1\min\left\{nt^2, nt\right\}}.
\end{align*}
Choose $\varepsilon=\frac{1}{2}$, we obtain the desired bound.
\end{proof}

\begin{lemma}\label{lemma: high-prob upper bound of spectral norm of sumU1Xi1(U3oU2)}
Suppose $\mcX \in \mathbb{R}^{p_1 \times p_2 \times p_3}$ is a tensor with independent zero-mean $ \sigma$-sub-Gaussian entries and $\mcX_1, \cdots, \mcX_n$ are i.i.d. copies of $\mcX$. $\wtU_j$'s, $j=1,2,3$ are fixed orthonormal matrix satisfying $\wtU_j \in\mathbb{O}^{p_j\times r_j}$. Let $\operatorname{Mat}_j\left(\mcX\right)$ be mode-$j$ matricization of $\mcX$. In addition, suppose $\xi_1, \cdots, \xi_n$ are independent zero-mean $ \sigma_{\xi}$-sub-Gaussian random variables. Then for any $x>0$, there exist absolute constants $C$ and $c$ such that
\begin{align}
& \mathbb{P}\left(\left\|\sum_{i=1}^n \xi_i \wtU_j^{\top} \Mat_j\left(\mcX_i\right)\left(\wtU_{j+2} \otimes \wtU_{j+1} \right)\right\| \geq C \sqrt{n} \sqrt{\widetilde{r}_j+\widetilde{r}_{j+1} \widetilde{r}_{j+2}+t} \cdot \sigma \sigma_\xi \right) \notag \\
\leq & 2\exp \left[-c\min\left(\frac{t^2}{\widetilde{r}_j}, t\right)\right] + \exp{\left(-c n\right)}. \label{eq: high-prob upper bound of spectral norm of sumU1Xi1(U3oU2)}
\end{align}
\end{lemma}

\begin{proof}

By symmetry, it suffices to consider and upper bound for
$
\left\|\sum_{i=1}^n \xi_i \wtU_1^{\top} \Mat_1\left(\mcX_i\right)\left(\wtU_3 \otimes \wtU_2\right)\right\|.
$
For any fixed $a=\left(a_1, \cdots, a_n\right) \in \mathbb{R}^n$, noting that the entries of $\sum_{i=1}^n a_i \Mat_1\left(\mcX_i\right)$ are independent $\sigma \left\|a\right\|_{\ell_2}$-sub-Gaussian random variables with mean 0 and variance $\sigma^2\left\|a\right\|_{\ell_2}^2$. By Lemma~\ref{lemma: concentration bound of |U1tildeX1(U3tildeoU2tilde)|}, we have
\begin{align*}
& \mathbb{P}\left(\left\|\sum_{i=1}^n a_i \wtU_1^{\top} \Mat_1\left(\mcX_i\right) \left(\wtU_3 \otimes \wtU_2\right)\right\| \geq C\left\|a\right\|_{\ell_2} \sqrt{\widetilde{r}_1+t} \cdot \sigma \right) \leq 2 \cdot 5^{\widetilde{r}_2\widetilde{r}_3}\exp \left[-c \min \left(\frac{t^2}{\widetilde{r}_1}, t \right)\right] . 
\end{align*}

Therefore,
\begin{align*}
& \mathbb{P}\left(\left\|\sum_{i=1}^n \xi_i \wtU_1^{\top} \Mat_1\left(\mcX_i\right) \left(\wtU_3 \otimes \wtU_2\right)\right\| \geq C\|\xi\|_2 \sqrt{\widetilde{r}_1+t} \cdot \sigma  \mid \left\{\xi_i\right\}_{i=1}^n \right) \leq 2\cdot 5^{\widetilde{r}_2\widetilde{r}_3} \exp \left[-c \min\left(\frac{t^2}{\widetilde{r}_1}, t\right)\right] .
\end{align*}
By Bernstein-type inequality for sub-Gaussian random variables, 
$$
\mathbb{P}\left(\|\xi\|_2 \geq C \sigma_\xi \sqrt{n} \right) \leq \exp {\left(-cn\right)} .
$$
Thus we have
\begin{align*}
& \mathbb{P}\left(\left\|\sum_{i=1}^n \xi_i \wtU_1^{\top} \Mat_1\left(\mcX_i\right) \left(\wtU_3 \otimes \wtU_2\right)\right\| \geq C \sigma \sqrt{\widetilde{r}_1+t}\cdot \sigma_\xi \sqrt{n}\right) \\
\leq & \mathbb{P}\left(\|\xi\|_2 \geq C \sigma_\xi \sqrt{n} \right) + \mathbb{P}\left(\left\|\sum_{i=1}^n \xi_i \wtU_1^{\top} \Mat_1\left(\mcX_i\right) \left(\wtU_3 \otimes \wtU_2\right)\right\| \geq  \sigma \sqrt{\widetilde{r}_1+t}\cdot \|\xi\|_2, \|\xi\|_2 \leq C \sigma_\xi \sqrt{n} \right) \\
\leq &  \mathbb{E}_{\xi}\left[ 2\cdot 5^{r_2r_3}\exp\left[-c\min\left(\frac{t^2}{\widetilde{r}_1},t\right)\right] \mathcal{I}\left(\left\|\xi\right\| \leq C\sigma_\xi \sqrt{n}\right) \right] + \exp {\left(-cn\right)}\\
\leq & 2\cdot 5^{r_2r_3}\exp\left[-c\min\left(\frac{t^2}{\widetilde{r}_1},t\right)\right] + \exp {\left(-cn\right)}. 
\end{align*}

\end{proof}

\section{Proof of Theorem~\ref{thm: main theorem in tensor PCA}}

In this section, we present the proof of Theorem~\ref{thm: main theorem in tensor PCA}. Since the proof of Theorem~\ref{thm: main theorem in tensor PCA} is similar to the proof of the Theorem~\ref{thm: main theorem in tensor PCA}, we will focus on the parts that differ. For identical or repetitive steps, such as the decomposition of certain terms, we will provide a concise description to maintain textual brevity. 

\begin{proof}[Proof of Theorem~\ref{thm: main theorem in tensor PCA}]

We begin by defining an event for the initial estimator:
\begin{align*}
\mathcal{E}_{U}^{\text{PCA}} = \left\{\left\|\mcP_{\whU_j^{(0)}} - \mcP_{U_j}\right\| > {\sigma \sqrt{\op}}/{\ulambda}, j=1,2,3\right\}
\end{align*}
that holds with probability $\mathbb{P}\left(\mathcal{E}_{U}^{\text{PCA}}\right)$. This analysis assumes that the event $\|\mcP_{\whU_j^{(0)}} - \mcP_{U_j}\| > {\sigma \sqrt{\op}}/{\ulambda}$ holds for any $j=1,2,3$. By Lemma~\ref{lemma: error contraction of l2 error of singular space in tensor PCA} , this assumption implies:
\begin{align*}
\left\|\mcP_{\whU_j^{(1)}} - \mcP_{U_j}\right\| \leq {\sigma \sqrt{\op}}/{\ulambda}, \quad \text{and} \quad \left\|\mcP_{\whU_j} - \mcP_{U_j}\right\|:=\left\|\mcP_{\whU_j^{(2)}} - \mcP_{U_j}\right\| \leq {\sigma \sqrt{\op}}/{\ulambda}
\end{align*}
holds with probability at least $1-\mathbb{P}\left(\mathcal{E}_{U}^{\text{PCA}}\right) - \op^{-c}$. This upper error bound will be utilized throughout the proof.

We then provide detailed proof of our main theorem.

\subsection*{Step 1: Upper bounds of negligible terms in $\langle \mcZ \times_1 \mcP_{\widehat{U}_1} \times_2 \mcP_{\widehat{U}_2} \times_3 \mcP_{\widehat{U}_3}, \mcA \rangle$}

Consider the decomposition of $\langle \mcZ, \mcA \times_1 \mcP_{\whU_1} \times_2 \mcP_{\whU_2} \times_3 \mcP_{\whU_3} \rangle$ similar to the arguments in Step 1 in the proof of Theorem~\ref{thm: main theorem in tensor regression without sample splitting}. We will prove the asymptotic normality of $\left\langle \mcZ \times_1 \mcP_{{U}_1} \times_2 \mcP_{{U}_2} \times_3 \mcP_{{U}_3}, \mcA \right\rangle$ later in Step 3 when finding the Berry-Essen bound of all asymptotic normal terms. 

Here, we focus on quantifying the upper bound of negligible terms. By symmetry, it suffices to consider  
\begin{align*}
\text{Step 1.1: } & \left\langle \mcZ \times_1 \left(\mcP_{\whU_1} - \mcP_{{U}_1}\right) \times_2 \mcP_{U_2} \times_3 \mcP_{U_3}, \mcA \right\rangle \\
\text{Step 1.2: } & \left\langle \mcZ \times_1 \left(\mcP_{\whU_1} - \mcP_{{U}_1}\right) \times_2 \left(\mcP_{\whU_2} - \mcP_{{U}_2}\right) \times_3 \mcP_{U_3}, \mcA  \right\rangle \\
\text{Step 1.3: } & \left\langle \mcZ \times_1 \left(\mcP_{\whU_1} - \mcP_{{U}_1}\right) \times_2 \left(\mcP_{\whU_2} - \mcP_{{U}_2}\right) \times_3 \left(\mcP_{\whU_3} - \mcP_{{U}_3}\right), \mcA \right\rangle. 
\end{align*}

\subsubsection*{Step 1.1: Upper Bound of  $\left\langle \mcZ \times_1 \left(\mcP_{\widehat{U}_1} - \mcP_{{U}_1}\right) \times_2 \mcP_{U_2} \times_3 \mcP_{U_3}, \mcA \right\rangle$}

First, consider
\begin{align}
& \left|\left\langle \mcZ \times_1 \left(\mcP_{\whU_1} - \mcP_{U_1}\right) \times_2 \mcP_{U_2} \times_3 \mcP_{U_3}, \mcA\right\rangle  \notag = \operatorname{tr}\left[\left(\mcP_{U_3}\otimes \mcP_{U_2}\right)A_1^{\top} \left(\mcP_{\whU_1} - \mcP_{U_1}\right)Z_1\right]\right| \notag \\
= & \underbrace{\left|\operatorname{tr}\left[\left(\mcP_{U_3}\otimes \mcP_{U_2}\right)A_1^{\top}U_1\left(G_1G_1^{\top}\right)^{-1}G_1\left(U_3\otimes U_2\right)^{\top} \left(\mcP_{\whU_3^{(1)}} \otimes \mcP_{\whU_2^{(1)}}\right)Z_1^{\top}\mcP_{U_{1\perp}}Z_1\right]\right|}_{\mathrm{\RN{1}}} \label{eq: term 1 in step 1.1 in tensor PCA}\\
& +  \underbrace{\left|\operatorname{tr}\left[\left(\mcP_{U_3}\otimes \mcP_{U_2}\right)A_1^{\top}U_1\left(G_1G_1^{\top}\right)^{-1}U_1^{\top}Z_1\left(\mcP_{\whU_3^{(1)}} \otimes \mcP_{\whU_2^{(1)}}\right)Z_1^{\top}\mcP_{U_{1\perp}}Z_1\right]\right|}_{\mathrm{\RN{2}}} \label{eq: term 2 in step 1.1 in tensor PCA}\\
& +  \underbrace{\left|\operatorname{tr}\left[\left(\mcP_{U_3}\otimes \mcP_{U_2}\right)A_1^{\top}\mcP_{U_{1\perp}}Z_1\left(\mcP_{\whU_3^{(1)}}\otimes \mcP_{\whU_2^{(1)}}\right)\left(U_3\otimes U_2\right)G_1^{\top}\left(G_1G_1^{\top}\right)^{-1}U_1^{\top}Z_1\right]\right|}_{\mathrm{\RN{3}}} \label{eq: term 3 in step 1.1 in tensor PCA}\\
& +  \underbrace{\left|\operatorname{tr}\left[\left(\mcP_{U_3}\otimes \mcP_{U_2}\right)A_1^{\top}\mcP_{U_{1\perp}}Z_1\left(\mcP_{\whU_3^{(1)}}\otimes \mcP_{\whU_2^{(1)}}\right)Z_1^{\top}U_1\left(G_1G_1^{\top}\right)^{-1}U_1^{\top}Z_1\right]\right|}_{\mathrm{\RN{4}}} \label{eq: term 4 in step 1.1 in tensor PCA}\\
& +  \underbrace{\left|\operatorname{tr}\left[\left(\mcP_{U_3}\otimes \mcP_{U_2}\right)A_1^{\top} \mcP_{U_1}\sum_{k_1=2}^{+\infty}S_{G_1,k_1}\left(E_1\right)Z_1\right]\right|}_{\mathrm{\RN{5}}} \label{eq: term 5 in step 1.1 in tensor PCA} \\
& +  \underbrace{\left|\operatorname{tr}\left[\left(\mcP_{U_3}\otimes \mcP_{U_2}\right)A_1^{\top} \mcP_{U_{1\perp}}\sum_{k_1=2}^{+\infty}S_{G_1,k_1}\left(E_1\right)Z_1\right]\right|}_{\mathrm{\RN{6}}} \label{eq: term 6 in step 1.1 in tensor PCA}.
\end{align}

We consider the upper bounds for \eqref{eq: term 1 in step 1.1 in tensor PCA}, \eqref{eq: term 2 in step 1.1 in tensor PCA}, \eqref{eq: term 3 in step 1.1 in tensor PCA}, \eqref{eq: term 4 in step 1.1 in tensor PCA} and \eqref{eq: term 5 in step 1.1 in tensor PCA} separately. We begin with the upper bound for the first term \eqref{eq: term 1 in step 1.1 in tensor PCA}. By Lemma \ref{lemma: high-prob upper bound of tr(BZtCZ)}, it follows that
\begin{align}
\mathrm{\RN{1}} 
\lesssim & \sigma^2\left|\operatorname{tr}\left[\left(\mcP_{U_3}\otimes \mcP_{U_2}\right)A_1^{\top}U_1\left(G_1G_1^{\top}\right)^{-1}G_1\left(U_3\otimes U_2\right)^{\top}\right]\operatorname{tr}\left[\mcP_{U_{1\perp}}\right]\right| 
 +  \left\|\left(\mcP_{U_3}\otimes \mcP_{U_2}\right)A_1^{\top}U_1\right\|_{\mathrm{F}}\cdot \frac{1}{\ulambda} \cdot \sigma^2 \sqrt{\op\log\left(\op\right)} \notag \\
& +  \left\|\left(\mcP_{U_3}\otimes \mcP_{U_2}\right)A_1^{\top}U_1\left(G_1G_1^{\top}\right)^{-1}G_1\left(U_3\otimes U_2\right)^{\top}\right\|_{\ell_\infty}\left\|\mcP_{U_{1\perp}}\right\|_{\ell_\infty}\cdot \sigma^2\log(\op) \notag \\
& +  \left\|\left(\mcP_{U_3}\otimes \mcP_{U_2}\right)A_1^{\top}U_1\right\|_{\mathrm{F}}\cdot \frac{1}{\ulambda} \cdot \frac{\sigma\sqrt{\op}}{{\ulambda}} \cdot \sigma\sqrt{\op\oor} \cdot \sigma\sqrt{\op} + \left\|\left(\mcP_{U_3}\otimes \mcP_{U_2}\right)A_1^{\top}U_1\right\|_{\mathrm{F}}\cdot \frac{1}{\ulambda} \cdot \frac{\sigma^2\op}{\ulambda^2} \cdot \sigma\sqrt{\op\oor} \cdot \sigma\sqrt{\op} \notag \\
\lesssim & \left\|U_1^{\top}A_1\mcP_{\left(U_3\otimes U_2\right)G_1^{\top}}\right\|_{\mathrm{F}}\cdot \frac{\sigma^2\oor^{1/2} \cdot \op}{\ulambda}
+ \left\|\mcA\times_1 U_1 \times_2 U_2 \times_3 U_3\right\|_{\mathrm{F}}\left(\frac{\sigma^2\oor^{1/2} \cdot \sqrt{\op\log(\op)}}{\ulambda}+\frac{\sigma^3\oor^{1/2} \cdot \op^{3/2}}{\ulambda^2}\right), \label{eq: upper bound of term 1 in step 1.1 in tensor PCA}
\end{align}
where the second inequality follows from $\left\|\mcP_{U_{1\perp}}\right\|_{\mathrm{F}}=\sqrt{p_1-r_2}\leq \sqrt{p_1} \leq \sqrt{\op}$,  $\left\|\mcP_{U_{1\perp}}\right\|_{\ell_\infty}\leq \left\|\mcP_{U_{1\perp}}\right\| \leq 1$, and
\begin{align*}
\left\|\left(\mcP_{U_3}\otimes \mcP_{U_2}\right)A_1^{\top}U_1\left(G_1G_1^{\top}\right)^{-1}G_1\left(U_3\otimes U_2\right)^{\top}\right\|_{\ell_\infty}
\leq & \frac{1}{\ulambda}\left\|\left(\mcP_{U_3}\otimes \mcP_{U_2}\right)A_1^{\top}U_1\right\|, \\
\left\|\left(\mcP_{U_3}\otimes \mcP_{U_2}\right)A_1^{\top}U_1\left(G_1G_1^{\top}\right)^{-1}G_1\left(U_3\otimes U_2\right)^{\top}\right\|_{\mathrm{F}} 
\leq & \frac{1}{\ulambda}\left\|\left(\mcP_{U_3}\otimes \mcP_{U_2}\right)A_1^{\top}U_1\right\|_{\mathrm{F}}.
\end{align*}

Similarly, for the second term ${\mathrm{\RN{2}}}$ \eqref{eq: term 2 in step 1.1 in tensor PCA}, we first have the following decomposition:
\begin{align}
\mathrm{\RN{2}}
= & \underbrace{\operatorname{tr}\left[\left(\mcP_{U_3}\otimes \mcP_{U_2}\right)A_1^{\top}\mcP_{U_{1\perp}}Z_1\left(U_3\otimes U_2\right)G_1^{\top}\left(G_1G_1^{\top}\right)^{-1}U_1^{\top}Z_1\right]}_{\mathrm{\RN{2}}.\mathrm{\RN{1}}} \label{eq: term 2.1 in step 1.1 in tensor PCA}\\
& +  \underbrace{\operatorname{tr}\left[\left(\mcP_{U_3}\otimes \mcP_{U_2}\right)A_1^{\top}\mcP_{U_{1\perp}}Z_1\left[\mcP_{U_3} \otimes \left(\mcP_{\whU_2^{(1)}} - \mcP_{U_2}\right)\right]\left(U_3\otimes U_2\right)G_1^{\top}\left(G_1G_1^{\top}\right)^{-1}U_1^{\top}Z_1\right]}_{\mathrm{\RN{2}}.\mathrm{\RN{2}}} \label{eq: term 2.2 in step 1.1 in tensor PCA}\\ 
& +  \underbrace{\operatorname{tr}\left[\left(\mcP_{U_3}\otimes \mcP_{U_2}\right)A_1^{\top}\mcP_{U_{1\perp}}Z_1\left[\left(\mcP_{\whU_3^{(1)}} - \mcP_{U_3}\right) \otimes \mcP_{U_2}\right]\left(U_3\otimes U_2\right)G_1^{\top}\left(G_1G_1^{\top}\right)^{-1}U_1^{\top}Z_1\right]}_{\mathrm{\RN{2}}.\mathrm{\RN{3}}} \label{eq: term 2.3 in step 1.1 in tensor PCA} \\
& +  \underbrace{\operatorname{tr}\left[\left(\mcP_{U_3}\otimes \mcP_{U_2}\right)A_1^{\top}\mcP_{U_{1\perp}}Z_1\left[\left(\mcP_{\whU_3^{(1)}} - \mcP_{U_3}\right)\otimes \left(\mcP_{\whU_2^{(1)}} - \mcP_{U_2}\right)\right]\left(U_3\otimes U_2\right)G_1^{\top}\left(G_1G_1^{\top}\right)^{-1}U_1^{\top}Z_1\right]}_{\mathrm{\RN{2}}.\mathrm{\RN{4}}} \label{eq: term 2.4 in step 1.1 in tensor PCA}.
\end{align}

Here, for the first term $\mathrm{\RN{2}}.\mathrm{\RN{1}}$ \eqref{eq: term 2.1 in step 1.1 in tensor PCA} in \eqref{eq: term 2 in step 1.1 in tensor PCA}, we have
\begin{align}
\mathrm{\RN{2}}.\mathrm{\RN{1}} 
\leq & \underbrace{\left\|\left(\mcP_{U_3}\otimes \mcP_{U_2}\right)A_1^{\top}\mcP_{U_{1\perp}}Z_1\left(U_3\otimes U_2\right)\right\|_{\mathrm{F}}}_{\eqref{eq: high-prob upper bound of U1Z1(U3oU2) in tensor PCA}} \cdot \left\|G_1^{\top}\left(G_1G_1^{\top}\right)^{-1}\right\|_{\mathrm{F}}\cdot \left\|U_1^{\top}Z_1\left(\mcP_{U_3}\otimes \mcP_{U_2}\right)\right\| \notag \\
\lesssim & \left\|\mcA \times_2 U_2 \times_3 U_3\right\|_{\mathrm{F}}\cdot \frac{\sigma^2\oor^{1/2} \cdot \oor\log(\op)}{\ulambda} \label{eq: upper bound of term 2.1 in step 1.1 in tensor PCA},
\end{align}
and for the second term $\mathrm{\RN{2}}.\mathrm{\RN{2}}$ \eqref{eq: term 2.2 in step 1.1 in tensor PCA} in \eqref{eq: term 2 in step 1.1 in tensor PCA}
\begin{align}
\mathrm{\RN{2}}.\mathrm{\RN{2}}
\leq & \left\|\left(\mcP_{U_3}\otimes \mcP_{U_2}\right)A_1^{\top}\mcP_{U_{1\perp}}\right\|_{\mathrm{F}}\sup_{\substack{W_2\in \mathbb{R}^{p_2\times r_2}, \left\|W_2\right\|=1\\ W_3\in \mathbb{R}^{p_3\times r_3}, \left\|W_3\right\|=1}}\left\|\mcP_{U_{1\perp}}Z_1\left(W_3\otimes W_2\right)\right\|\left\|\mcP_{U_3} \otimes \left(\mcP_{\whU_2^{(1)}} - \mcP_{U_2}\right)\right\| \notag \\
 & \cdot \left\|\left(U_3\otimes U_2\right)G_1^{\top}\left(G_1G_1^{\top}\right)^{-1}U_1^{\top}\right\|\left\|U_1^{\top}Z_1\left(\mcP_{U_3}\otimes \mcP_{U_2}\right)\right\|_{\mathrm{F}} \notag \\
\lesssim & \left\|\mcA \times_2 U_2 \times_3 U_3\right\|_{\mathrm{F}}\cdot \frac{\sigma^3\oor^{1/2} \cdot \op\sqrt{\oor\log(\op)}}{\ulambda^2} \label{eq: upper bound of term 2.2 in step 1.1 in tensor PCA}.
\end{align}
For the third term $\mathrm{\RN{2}}.\mathrm{\RN{3}}$ \eqref{eq: term 2.3 in step 1.1 in tensor PCA} in \eqref{eq: term 2 in step 1.1 in tensor PCA}, by symmetry, we have
\begin{align}
\mathrm{\RN{2}}.\mathrm{\RN{3}}
\lesssim & \left\|\mcA \times_2 U_2 \times_3 U_3\right\|_{\mathrm{F}} \cdot \frac{\sigma^3\oor^{1/2} \cdot \op\sqrt{\oor\log(\op)}}{\ulambda^2} \label{eq: upper bound of term 2.3 in step 1.1 in tensor PCA}, 
\end{align}
and for the fourth term $\mathrm{\RN{2}}.\mathrm{\RN{4}}$ \eqref{eq: term 2.4 in step 1.1 in tensor PCA} in \eqref{eq: term 2 in step 1.1 in tensor PCA}
\begin{align}
\mathrm{\RN{2}}.\mathrm{\RN{4}}
\leq & \frac{1}{\ulambda} \cdot \left\|\left(\mcP_{U_3}\otimes \mcP_{U_2}\right)A_1^{\top}\mcP_{U_{1\perp}}\right\|_{\mathrm{F}}\sup_{\substack{W_2\in \mathbb{R}^{p_2\times r_2}, \left\|W_2\right\|=1\\ W_3\in \mathbb{R}^{p_3\times r_3}, \left\|W_3\right\|=1}} \left\|\mcP_{U_{1\perp}}Z_1\left(W_3\otimes W_2\right)\right\| \notag \\
& \cdot \left\|\mcP_{U_{1\perp}}Z_1\left(W_3\otimes W_2\right)\right\|\left\|\left(\mcP_{\whU_2^{(1)}} - \mcP_{U_2}\right) \otimes \left(\mcP_{\whU_3^{(1)}} - \mcP_{U_3}\right)\right\| \cdot \left\|U_1^{\top}Z_1\left(\mcP_{U_3}\otimes \mcP_{U_2}\right)\right\|_{\mathrm{F}} \notag \\
\lesssim & \left\|\mcA \times_2 U_2 \times_3 U_3\right\|_{\mathrm{F}}\cdot  \frac{\sigma^4\oor^{1/2} \cdot \op^{3/2}\sqrt{\oor\log(\op)}}{\ulambda^3} \label{eq: upper bound of term 2.4 in step 1.1 in tensor PCA}.
\end{align}

It implies that 
\begin{align}
\mathrm{\RN{2}} 
\lesssim & \eqref{eq: upper bound of term 2.1 in step 1.1 in tensor PCA} + \eqref{eq: upper bound of term 2.2 in step 1.1 in tensor PCA} + \eqref{eq: upper bound of term 2.3 in step 1.1 in tensor PCA} + \eqref{eq: upper bound of term 2.4 in step 1.1 in tensor PCA} \lesssim \left\|\mcA \times_2 U_2 \times_3 U_3\right\|_{\mathrm{F}}\cdot \left(\frac{\sigma^2\oor^{1/2} \cdot \oor^{1/2}\log(\op)}{\ulambda}+ \frac{\sigma^3\oor^{1/2} \cdot  \op\sqrt{\oor\log(\op)}}{\ulambda^2}\right) \label{eq: upper bound of term 2 in step 1.1 in tensor PCA},
\end{align}
where the last inequality holds as long as $\frac{\ulambda}{\sigma} \geq \kappa\sqrt{\op}, \kappa \geq 1$.

For the third term \eqref{eq: term 3 in step 1.1 in tensor PCA}, we have 
\begin{align}
\mathrm{\RN{3}} 
\leq & \frac{1}{\ulambda} \cdot \left\|\left(\mcP_{U_3} \otimes \mcP_{U_2}\right) A_1^{\top} \mcP_{U_{1\perp}}Z_1\left(U_3\otimes U_2\right) \right\| \cdot \left\|U_1^{\top}Z_1\left(U_3\otimes U_2\right)\right\|_{\mathrm{F}} \notag \\ 
+ & \frac{1}{\ulambda} \cdot \left\|\left(\mcP_{U_3} \otimes \mcP_{U_2}\right) A_1^{\top}\right\|_{\mathrm{F}} \cdot \sup_{\substack{W_2\in \mathbb{R}^{p_2\times r_2}, \left\|W_2\right\|=1\\W_3\in \mathbb{R}^{p_3\times r_3},  \left\|W_3\right\|=1}}\left\|\mcP_{U_{1\perp}}Z_1\left(W_3\otimes W_2\right) \right\| \notag \\
& \cdot \left(\left\|\mcP_{\whU_2^{(1)}}-\mcP_{U_2}\right\| + \left\|\mcP_{\whU_3^{(1)}}-\mcP_{U_3}\right\| + \prod_{j=2}^3 \left\|\mcP_{\whU_j^{(1)}}-\mcP_{U_j}\right\|\right) \cdot  \left\|U_1^{\top}Z_1\left(U_3\otimes U_2\right)\right\|_{\mathrm{F}} \notag \\ 
\lesssim & \left\|\mcA \times_2 U_2 \times_3 U_3\right\|_{\mathrm{F}} \cdot \left(\frac{\sigma^2\oor^{1/2}\cdot \oor \log(\op)}{\ulambda}+ \frac{\sigma^3\oor^{1/2} \cdot \oor^{1/2}\op\sqrt{\log(\op)}}{\ulambda^2}\right) \label{eq: upper bound of term 3 in step 1.1 in tensor PCA}.
\end{align}

For the fourth term $\mathrm{\RN{4}}$ \eqref{eq: term 4 in step 1.1 in tensor PCA}, we have
\begin{align}
\mathrm{\RN{4}}
\leq & \frac{1}{\ulambda^2} \cdot \left\|\left(\mcP_{U_3} \otimes \mcP_{U_2}\right) A_1^{\top}\mcP_{U_{1\perp}} Z_1\left(U_3 \otimes U_2\right)\right\|_{\mathrm{F}} \cdot \left\|\left(U_3 \otimes U_2\right) Z_1^{\top}U_1\right\| \cdot \left\|U_1^{\top} Z_1 \left(\mcP_{U_3} \otimes \mcP_{U_2}\right)\right\|_{\mathrm{F}} \notag \\
& +  \left\|\left(\mcP_{U_3} \otimes \mcP_{U_2}\right) A_1^{\top}\mcP_{U_{1\perp}} \right\|_{\mathrm{F}} \cdot \sup_{\substack{W_2\in \mathbb{R}^{p_2\times r_2}, \left\|W_2\right\|=1\\W_3\in \mathbb{R}^{p_3\times 2r_3}, \left\|W_3\right\|=1}}\left\|\mcP_{U_{1\perp}}Z_1\left(W_3 \otimes W_2\right)\right\|_{\mathrm{F}} \notag \\
&\quad \cdot \left(\left\|\mcP_{\whU_2^{(1)}}-\mcP_{U_2}\right\| + \left\|\mcP_{\whU_3^{(1)}}-\mcP_{U_3}\right\| + \prod_{j=2}^3 \left\|\mcP_{\whU_j^{(1)}}-\mcP_{U_j}\right\|\right) \cdot \left\|\left(U_3 \otimes U_2\right) Z_1^{\top} U_1\right\|^2 \notag \\
\lesssim & \left\|\mcA \times_1 U_1 \times_2 U_2 \times_3 U_3\right\|_{\mathrm{F}}\cdot \left(\frac{\sigma^3\oor^{1/2}\cdot \oor^{3/2}\log(\op)^{3/2}}{\ulambda^2}+\frac{\sigma^3\oor^{1/2}\cdot \op\oor\log(\op)}{\ulambda^3}\right) \label{eq: upper bound of term 4 in step 1.1 in tensor PCA}.
\end{align}

Then, consider the fifth term (higher-order terms) \eqref{eq: term 5 in step 1.1 in tensor PCA}, 
\begin{align}
\mathrm{\RN{5}} 
\lesssim & \frac{1}{\ulambda^3}\cdot \left\|\left(\mcP_{U_3}\otimes \mcP_{U_2}\right)A_1^{\top} \mcP_{U_1}\right\|_{\mathrm{F}} \cdot \left\|\mcP_1^{-\frac{1}{2}}E_1U_{1\perp}\right\| \left\|U_{1\perp}^{\top}E_1U_{1\perp}\right\| \left\|U_{1\perp}^{\top}Z_1\left(\mcP_{U_3}\otimes \mcP_{U_2}\right)\right\|_{\mathrm{F}} \notag \\
& +  \frac{1}{\ulambda}\cdot \left\|\left(\mcP_{U_3}\otimes \mcP_{U_2}\right)A_1^{\top} \mcP_{U_1}\right\|_{\mathrm{F}} \cdot \left\|\mcP_1^{-\frac{1}{2}}E_1U_1\mcP_1^{-\frac{1}{2}}\right\|\cdot \left\|\mcP_1^{-\frac{1}{2}}U_1^{\top}E_1U_{1\perp}\right\|\left\|U_{1\perp}^{\top}Z_1\left(\mcP_{U_3}\otimes \mcP_{U_2}\right)\right\|_{\mathrm{F}} 
 \notag \\
& +  \frac{1}{\ulambda^2}\cdot \left\|\left(\mcP_{U_3}\otimes \mcP_{U_2}\right)A_1^{\top} \mcP_{U_1}\right\|_{\mathrm{F}} \cdot \left\|\mcP_1^{-\frac{1}{2}}E_1U_{1\perp}\right\| \cdot \left\|U_{1\perp}^{\top}E_1U_1\mcP_1^{-\frac{1}{2}}\right\| \cdot  \left\|\mcP_{U_1}Z_1\left(\mcP_{U_3}\otimes \mcP_{U_2}\right)\right\|_{\mathrm{F}} \notag \\
& +  \left|\operatorname{tr}\left[\left(\mcP_{U_3}\otimes \mcP_{U_2}\right)A_1^{\top} \mcP_{U_1}\sum_{k_1=3}^{+\infty}S_{G_1,k_1}\left(E_1\right)Z_1\right] \right|_{\mathrm{F}} \notag \\
\lesssim & \left\|\mcA \times_1 U_1 \times_2 U_2 \times_3 U_3\right\|_{\mathrm{F}}\cdot 
\left(\frac{\sigma^3\oor^{1/2}\cdot\op\sqrt{\log(\op)}}{\ulambda^2}+\frac{\sigma^4\oor^{1/2} \cdot \op^2}{\ulambda^3}\right) \label{eq: upper bound of term 5 in step 1.1 in tensor PCA}.
\end{align}

In addition, for $\mathrm{\RN{6}}$  \eqref{eq: term 6 in step 1.1 in tensor PCA}, by similar arguments, we have
\begin{align}
\mathrm{\RN{6}} 
\leq & \underbrace{\left\|\left(\mcP_{U_3}\otimes \mcP_{U_2}\right)A_1^{\top} \mcP_{U_{1\perp}}\sum_{k_1=2}^{+\infty}S_{G_1,k_1}\mcP_{U_{1\perp}}\right\|_{\mathrm{F}}}_{\eqref{eq: high-prob upper bound of PUporder2PUpPV in tensor PCA}} \cdot \left\|\mcP_{U_{1\perp}}Z_1\left(\mcP_{U_3}\otimes \mcP_{U_2}\right)\right\|_{\mathrm{F}} \notag \\
& +  \underbrace{\left\|\left(\mcP_{U_3}\otimes \mcP_{U_2}\right)A_1^{\top} \mcP_{U_{1\perp}}\sum_{k_1=2}^{+\infty}S_{G_1,k_1}\mcP_{U_1}\right\|_{\mathrm{F}}}_{\eqref{eq: high-prob upper bound of PUorder2PUpPV in tensor PCA}} \cdot \left\|\mcP_{U_1}Z_1\left(\mcP_{U_3}\otimes \mcP_{U_2}\right)\right\|_{\mathrm{F}} \notag \\
\lesssim & \left\|\mcA \times_2 U_2 \times_3 U_3\right\|_{\mathrm{F}} \cdot \left(\frac{\sigma^3\oor^{1/2} \cdot \op\sqrt{\oor\log(\op)}}{\ulambda^2} + \frac{\sigma^5\oor^{1/2} \cdot \op^{5/2}}{\ulambda^4}\right)
\label{eq: upper bound of term 6 in step 1.1 in tensor PCA}. 
\end{align}

Therefore, we have
\begin{align*}
& \left| \left\langle\mcZ \times_1\left(\mcP_{\whU_1}-\mcP_{U_1}\right) \times_2 \mcP_{U_2} \times_3 \mcP_{U_3}, \mcA\right\rangle\right| 
\lesssim  \eqref{eq: upper bound of term 1 in step 1.1 in tensor PCA} + \eqref{eq: upper bound of term 2 in step 1.1 in tensor PCA} + \eqref{eq: upper bound of term 3 in step 1.1 in tensor PCA} + \eqref{eq: upper bound of term 4 in step 1.1 in tensor PCA} + \eqref{eq: upper bound of term 5 in step 1.1 in tensor PCA} + \eqref{eq: upper bound of term 6 in step 1.1 in tensor PCA} \\
\lesssim & \left\|U_1^{\top}A_1\mcP_{\left(U_3\otimes U_2\right)G_1^{\top}}\right\|_{\mathrm{F}}\cdot \frac{\sigma^2\oor^{1/2} \cdot \op}{\ulambda}
+ \left\|\mcA\times_1 U_1 \times_2 U_2 \times_3 U_3\right\|_{\mathrm{F}}\cdot \left(\frac{\sigma^2\oor^{1/2} \cdot \sqrt{\op\log(\op)}}{\ulambda}+\frac{\sigma^3\oor^{1/2} \cdot \op^{3/2}}{\ulambda^2}\right) \\
& +  \left\|\mcA \times_2 U_2 \times_3 U_3\right\|_{\mathrm{F}} \cdot \left(\frac{\sigma^2\oor^{1/2}\cdot \oor \log(\op)}{\ulambda} + \frac{\sigma^3\oor^{1/2} \cdot \op\sqrt{\oor\log(\op)}}{\ulambda^2} + \frac{\sigma^4\oor^{1/2} \cdot \op}{\ulambda^3}\right).
\end{align*}

\subsubsection*{Step 1.2: Upper Bound of Negligible Terms in $\left\langle \mcZ \times_1 \left(\mcP_{\widehat{U}_1} - \mcP_{U_1}\right) \times_2 \left(\mcP_{\widehat{U}_2} - \mcP_{U_2}\right) \times_3 \mcP_{U_3}, \mcA \right\rangle$}

By the same decomposition in Step 1.2 in the proof of Theorem~\ref{thm: main theorem in tensor regression without sample splitting}, with $\widehat{\mcZ}$ replaced by $\mcZ$, we consider finding an upper bound for the following term:
\begin{align}
\left|\left\langle \mcZ \times_1 \left(\mcP_{\whU_1} - \mcP_{U_1}\right) \times_2 \left(\mcP_{\whU_2} - \mcP_{U_2}\right) \times_3 \mcP_{U_3}, \mcA \right\rangle\right| 
\leq & \mathrm{\RN{1}} + \mathrm{\RN{2}} + \mathrm{\RN{3}} + \mathrm{\RN{4}}. \label{eq: decomposition of step 1.2 in tensor PCA}
\end{align}

Applying similar arguments in the proof of Step 1.1, we can show
\begin{align*}
& \left|\left\langle\mcZ \times_1\left(\mcP_{\whU_1}-\mcP_{U_1}\right) \times_2\left(\mcP_{\whU_2}-\mcP_{U_2}\right) \times_3 \mcP_{U_3}, \mcA\right\rangle\right| \\
+ & {\left\|\mcA \times_3 U_3 \right\|_{\mathrm{F}} \cdot \left(\frac{\sigma^3 \oor^{1/2} \cdot \oor^{5/2}\log (\op)^{3/2}}{\ulambda^2} + \frac{\sigma^4 \oor^{1/2} \cdot \op\oor\log(\op)}{\ulambda^3} + \frac{\sigma^5 \oor^{1/2} \cdot \oR \op^{3/2} \sqrt{\log (\op)}}{\ulambda^4}+\frac{\sigma^6 \oor^{1/2} \cdot \op^{5/2} \sqrt{\oR \log (\op)}}{\ulambda^5}\right)}. 
\end{align*}

\subsubsection*{Step 1.3: Upper Bound of Negligible Terms in $\left\langle \mcZ \times_1 \left(\mcP_{\widehat{U}_1} - \mcP_{U_1}\right) \times_{2} \left(\mcP_{\widehat{U}_{2}} - \mcP_{U_{2}}\right) \times_{3} \left(\mcP_{\widehat{U}_{3}} - \mcP_{U_3}\right), \mcA \right\rangle $}

By symmetry, it suffices to consider
\begin{align}
\mathrm{\RN{1}}= & \left|\left\langle \mcZ \times_1 \left(\mcP_{\whU_1} - \mcP_{U_1}\right) \times_{2} \left(\mcP_{\whU_{2}} - \mcP_{U_{2}}\right) \times_{3} \left(\mcP_{\whU_{3}} - \mcP_{U_3}\right), \mcA \times_1 \mcP_{U_1} \times_2 \mcP_{U_2} \times_3 \mcP_{U_3} \right\rangle\right| , \label{eq: term 1 in step 1.3 in tensor PCA}\\
\mathrm{\RN{2}}= & \left|\left\langle \mcZ \times_1 \left(\mcP_{\whU_1} - \mcP_{U_1}\right) \times_{2} \left(\mcP_{\whU_{2}} - \mcP_{U_{2}}\right) \times_{3} \left(\mcP_{\whU_{3}} - \mcP_{U_3}\right), \mcA \times_1 \mcP_{U_{1\perp}} \times_2 \mcP_{U_2} \times_3 \mcP_{U_3} \right\rangle\right| , \label{eq: term 2 in step 1.3 in tensor PCA} \\
\mathrm{\RN{3}}= & \left|\left\langle \mcZ \times_1 \left(\mcP_{\whU_1} - \mcP_{U_1}\right) \times_{2} \left(\mcP_{\whU_{2}} - \mcP_{U_{2}}\right) \times_{3} \left(\mcP_{\whU_{3}} - \mcP_{U_3}\right), \mcA \times_1 \mcP_{U_{1\perp}} \times_2 \mcP_{U_{2\perp}} \times_3 \mcP_{U_3} \right\rangle\right| , \label{eq: term 3 in step 1.3 in tensor PCA}\\
\mathrm{\RN{4}}= & \left|\left\langle \mcZ \times_1 \left(\mcP_{\whU_1} - \mcP_{U_1}\right) \times_{2} \left(\mcP_{\whU_{2}} - \mcP_{U_{2}}\right) \times_{3} \left(\mcP_{\whU_{3}} - \mcP_{U_3}\right), \mcA \times_1 \mcP_{U_{1\perp}} \times_2 \mcP_{U_{2\perp}} \times_3 \mcP_{U_{3\perp}} \right\rangle\right| \label{eq: term 4 in step 1.3 in tensor PCA}.
\end{align}

Applying similar arguments in the proof of Step 1.1 and by the same decomposition in Step 1.3 in the proof of Theorem~\ref{thm: main theorem in tensor regression without sample splitting}, we can show
\begin{align*}
& \left|\left\langle \mcZ \times_1 \left(\mcP_{\whU_1} - \mcP_{U_1}\right) \times_{2} \left(\mcP_{\whU_{2}} - \mcP_{U_{2}}\right) \times_{3} \left(\mcP_{\whU_{3}} - \mcP_{U_3}\right), \mcA \right\rangle\right| \\
\lesssim & \underbrace{\left\|\mcA \times_1 U_1 \times_2 U_2 \times_3 U_3\right\|_{\mathrm{F}} \cdot \frac{\sigma^4 \oor^{1/2} \cdot \op^2}{\ulambda^3}}_{\mathrm{\RN{1}} }\\ 
+ & \underbrace{\sum_{j=1}^3 \left\|\mcA \times_{j+1} U_{j+1} \times_{j+2} U_{j+2}\right\|_{\mathrm{F}} \cdot \left(\frac{\sigma^2\oor^{1/2}\cdot \oor \log(\op)}{\ulambda} + \frac{\sigma^3\oor^{1/2} \cdot \op\sqrt{\oor\log(\op)}}{\ulambda^2} + \frac{\sigma^4\oor^{1/2} \cdot \op}{\ulambda^3}\right)}_{\mathrm{\RN{2}} }\\
+ & \underbrace{\sum_{j=1}^3 \left\|\mcA \times_j U_j\right\|_{\mathrm{F}} \cdot \left(\frac{\sigma^4 \oor^{1/2} \cdot \op\oor\log(\op)^2}{\ulambda^3}+\frac{\sigma^6 \oor^{1/2} \cdot \oR \op^{3/2} \log (\op)}{\ulambda^4}+\frac{\sigma^6 \oor^{1/2} \cdot \op^{5/2} \sqrt{\oR \log (\op)}}{\ulambda^5}\right)}_{\mathrm{\RN{3}} } \\ 
+ & \left\|\mcA\right\|_{\mathrm{F}} \cdot \left(\frac{\sigma^4 \oor^{1 / 2} \cdot \oor^4 \log (\op)^2}{\ulambda^3} +\frac{\sigma^5\oor^{1/2} \cdot \oR^{3/2}\oor^{1/2} \op^{1/2} \log (\op)^2}{\ulambda^4}\right) \\
+ & \underbrace{\left\|\mcA\right\|_{\mathrm{F}} \cdot \left( \frac{\sigma^5 \oor^{1/2} \cdot \op \oor^{3/2}\log(\op)^{3/2}}{\ulambda^4}+\frac{\sigma^6 \oor^{1/2} \cdot \oR^{3 / 2} \op^{3 / 2} \log (\op)^{3 / 2}}{\ulambda^5}+\frac{\sigma^7 \oor^{1/2} \cdot \oR \op^{5 / 2} \log (\op)}{\ulambda^6}\right)}_{\mathrm{\RN{4}} }. 
\end{align*}

\subsection*{Step 2: Upper bounds of negligible terms in $\left\langle\mcT \times_1 \mcP_{\widehat{U}_1} \times_2 \mcP_{\widehat{U}_2} \times_3 \mcP_{\widehat{U}_3}-\mcT, \mcA\right\rangle$} 

Consider the decomposition of $\left\langle\mcT \times_1 \mcP_{\whU_1} \times_2 \mcP_{\whU_2} \times_3 \mcP_{\whU_3}-\mcT, \mcA\right\rangle$ similar to the arguments in Step 2 in the proof of Theorem~\ref{thm: main theorem in tensor regression without sample splitting}. In Step 2, we will consider
\begin{align*}
\text{(Step 2.1)}: & \left\langle \mcT \times_1 \left(\mcP_{\whU_1} - \mcP_{U_1}\right) \times_2 \mcP_{U_2} \times_3 \mcP_{U_3}, \mcA \right\rangle \\
\text{(Step 2.2)}: & \left\langle \mcT \times_1 \left(\mcP_{\whU_1} - \mcP_{U_1}\right) \times_2 \left(\mcP_{\whU_2} - \mcP_{U_2}\right) \times_3 \mcP_{U_3}, \mcA \right\rangle \\
\text{(Step 2.3)}: & \left\langle \mcT \times_1 \left(\mcP_{\whU_1} - \mcP_{U_1}\right) \times_2 \left(\mcP_{\whU_2} - \mcP_{U_2}\right) \times_3 \left(\mcP_{\whU_3} - \mcP_{U_3}\right), \mcA \right\rangle,
\end{align*}

\subsubsection*{Step 2.1: Upper Bound of Negligible Terms in $\left\langle \mcT \times_1 \left(\mcP_{\widehat{U}_1} - \mcP_{U_1}\right) \times_2 \mcP_{U_2} \times_3 \mcP_{U_3}\right\rangle$}

Note that
\begin{align}
\left\langle \mcT \times_1 \left(\mcP_{\whU_1} - \mcP_{U_1}\right) \times_2 \mcP_{U_2} \times_3 \mcP_{U_3}, \mcA \right\rangle 
= & \underbrace{\left\langle \mcT \times_1 S_{G_1,1}\left(E_1\right) \times_2 \mcP_{U_2} \times_3 \mcP_{U_3}, \mcA \right\rangle}_{\mathrm{\RN{1}}} \label{eq: term 1 in T(PUhat1-PU1)PU2PU3A in tensor regression}\\
+ & \underbrace{\left\langle \mcT \times_1 \sum_{k_1=2}^{+\infty}S_{G_1,k_1}\left(E_1\right) \times_2 \mcP_{U_2} \times_3 \mcP_{U_3}, \mcA \right\rangle}_{\mathrm{\RN{2}}}. \label{eq: term 2 in T(PUhat1-PU1)PU2PU3A in tensor regression}
\end{align}
We first consider 
\begin{align*}
\mathrm{\RN{1}}
= & \underbrace{\left\langle \mcP_{U_{1\perp}} Z_1\left(\mcP_{U_3} \otimes \mcP_{U_2}\right) \mcP_{\left(U_3 \otimes U_2\right)G_1^{\top}}, A_1\right\rangle}_{\mathrm{\RN{1}}.\mathrm{\RN{1}}, \text{asymptotically normal}}
+  \underbrace{\left\langle \mcP_{U_{1\perp}} Z_1\left[\left(\mcP_{\whU_3^{(1)}}-\mcP_{U_3}\right) \otimes \mcP_{U_2}\right] \mcP_{\left(U_3 \otimes U_2\right)G_1^{\top}} , A_1\right\rangle}_{\mathrm{\RN{1}}.\mathrm{\RN{2}}, \text{negligible}} \\
&+  \underbrace{\left\langle \mcP_{U_{1\perp}} Z_1\left[\mcP_{U_3} \otimes \left(\mcP_{\whU_2^{(1)}}-\mcP_{U_2}\right)\right] \mcP_{\left(U_3 \otimes U_2\right)G_1^{\top}}, A_1\right\rangle}_{\mathrm{\RN{1}}.\mathrm{\RN{3}}, \text{negligible}}\\
&+  \underbrace{\left\langle \mcP_{U_{1\perp}} Z_1\left[\left(\mcP_{\whU_2^{(1)}}-\mcP_{U_2}\right) \otimes \left(\mcP_{\whU_3^{(1)}}-\mcP_{{U}_3}\right)\right] \mcP_{\left(U_3 \otimes U_2\right)G_1^{\top}}, A_1\right\rangle}_{\mathrm{\RN{1}}.\mathrm{\RN{4}}, \text{negligible}}\\
&+  \underbrace{\left\langle \mcP_{U_{1\perp}} Z_1\left(\mcP_{\whU_3^{(1)}} \otimes \mcP_{\whU_2^{(1)}}\right) Z_1^{\top} U_1 \left(G_1G_1^{\top}\right)^{-1}  G_1 \left(U_3 \otimes U_2\right)^{\top} , A_1\right\rangle}_{\mathrm{\RN{1}}.\mathrm{\RN{5}}, \text{negligible}}.
\end{align*}

Applying similar arguments in the proof of Step 1.1 and by the same decomposition in Step 2.1 in the proof of Theorem~\ref{thm: main theorem in tensor regression without sample splitting}, we can show
\begin{align*}
& \left|\left\langle\mcT \times_1 \left(\mcP_{\whU_1} - \mcP_{U_1}\right) \times_2 \mcP_{U_2} \times_3 \mcP_{U_3}, \mcA\right\rangle-\left\langle\mcP_{U_{1\perp}} Z_1\mcP_{\left(U_3 \otimes U_2\right) G_1^{\top}}, A_1\right\rangle \right| \\
\lesssim & \left\|\mcP_{U_{1 \perp}} A_1 \mcP_{\left(U_3 \otimes U_2\right) G_1^{\top}}\right\|_{\mathrm{F}} \cdot\left(\frac{\sigma^2 \oor \cdot \sqrt{\op \log (\op)}}{\ulambda}+\frac{\sigma^3 \oor \cdot \op^{3 / 2}}{\ulambda^2}\right) + \left\|\mcP_{U_1} A_1 \mcP_{\left(U_3 \otimes U_2\right) G_1^{\top}}\right\|_{\mathrm{F}} \cdot \frac{\sigma^2\oor^{1/2} \cdot \op}{\ulambda} \\
& +  \left\|\mcP_{U_{1 \perp}} A_1 \mcP_{\left(U_3 \otimes U_2\right) G_1^{\top}}\right\|_{\mathrm{F}} \cdot \left(\frac{\sigma^2\oor^{1/2} \cdot \oor\log(\op)}{\ulambda} + \frac{\sigma^3\oor^{1/2} \cdot \op\sqrt{\oor\log(\op)}}{\ulambda^2} + \frac{\sigma^4\oor^{1/2} \cdot \op^2}{\ulambda^3}\right) \\
\lesssim & \left\|U_1^{\top} A_1 \mcP_{\left(U_3 \otimes U_2\right) G_1^{\top}}\right\|_{\mathrm{F}} \cdot \frac{\sigma^2\oor^{1/2}\cdot \op}{\ulambda} + \left\|\mcP_{U_{1 \perp}} A_1 \mcP_{\left(U_3 \otimes U_2\right) G_1^{\top}}\right\|_{\mathrm{F}} \cdot\left(\frac{\sigma^2 \oor \cdot \sqrt{\op \log (\op)}}{\ulambda}+\frac{\sigma^3 \oor \cdot \op^{3/2}}{\ulambda^2}\right).
\end{align*}

\subsubsection*{Step 2.2: Upper Bound of Negligible Terms in $\left\langle \mcT\times_1 \left(\mcP_{\widehat{U}_1} - \mcP_{U_1}\right) \times_2 \left(\mcP_{\widehat{U}_2} - \mcP_{U_2}\right) \times_3 \mcP_3, \mcA\right\rangle $}
Consider the same decomposition in Step 2.2 in the proof of Theorem~\ref{thm: main theorem in tensor regression without sample splitting}:
\begin{align*}
\left|\left\langle \mcT\times_1 \left(\mcP_{\whU_1} - \mcP_{U_1}\right) \times_2 \left(\mcP_{\whU_2} - \mcP_{U_2}\right) \times_3 \mcP_3, \mcA\right\rangle\right| \leq & \mathrm{\RN{1}} + \mathrm{\RN{2}} + \mathrm{\RN{3}} + \mathrm{\RN{4}}. 
\end{align*}

First, by \eqref{eq: high-prob upper bound of AoPU1p(PUhat1-P_U1)PU1oPU2p(PUhat2-P_U2)PU2oPU3p(PUhat3-P_U3)PU3 in tensor PCA}, we have
\begin{align}
\mathrm{\RN{1}}
\lesssim & \left\|\mcA \times_1 U_1\right\|_{\mathrm{F}} \cdot \left(\frac{\sigma^2 \oor^{1/2} \cdot \oor\log(\op)}{\ulambda}+\frac{\sigma^3 \oor^{1/2} \cdot \oR \sqrt{\op \log (\op)}}{\ulambda^2}+\frac{\sigma^4 \oor^{1/2} \cdot \op^{3 / 2} \sqrt{\oR \log (\op)}}{\ulambda^3}\right). \label{eq: upper bound of term 1 in step 2.2 in tensor PCA}
\end{align}

Then, consider
\begin{align}
\mathrm{\RN{2}} 
\lesssim & \left\|\mcA \times_2 U_2 \times_3 U_3\right\|_{\mathrm{F}} \cdot \left(\frac{\sigma^3\oor^{1/2} \cdot \op\sqrt{\oor\log(\op)}}{\ulambda^2} + \frac{\sigma^4\oor^{1/2} \cdot \op^2}{\ulambda^3}\right). \label{eq: upper bound of term 2 in step 2.2 in tensor PCA}
\end{align}

By symmetry, it also implies that
\begin{align}
\mathrm{\RN{3}}
\lesssim & \left\|\mcA \times_1 U_1 \times_3 U_3\right\|_{\mathrm{F}} \cdot \left(\frac{\sigma^3\oor^{1/2} \cdot \op\sqrt{\oor\log(\op)}}{\ulambda^2} + \frac{\sigma^4\oor^{1/2} \cdot \op^2}{\ulambda^3}\right). \label{eq: upper bound of term 3 in step 2.2 in tensor PCA}
\end{align}

Finally, consider
\begin{align}
\mathrm{\RN{4}}
\lesssim & \left\|\mcA \times_1 U_1 \times_2 U_2 \times_3 U_3\right\|_{\mathrm{F}} \cdot \left(\frac{\sigma^4\oor^{1/2} \cdot \op^2}{\ulambda^3}\right). \label{eq: upper bound of term 4 in step 2.2 in tensor PCA}
\end{align}

\subsubsection*{Step 2.3: Upper Bound of Negligible Terms in $\left\langle \mcT\times_1 \left(\mcP_{\widehat{U}_1} - \mcP_{U_1}\right) \times_2 \left(\mcP_{\widehat{U}_2} - \mcP_{U_2}\right) \times_3 \left(\mcP_{\widehat{U}_3} - \mcP_{U_3}\right), \mcA\right\rangle$}

By similar arguments of decomposition in Step 2.3 in the proof of Theorem~\ref{thm: main theorem in tensor regression without sample splitting}, we have
\begin{align}
& \left|\left\langle \mcT\times_1 \left(\mcP_{\whU_1} - \mcP_{U_1}\right) \times_2 \left(\mcP_{\whU_2} - \mcP_{U_2}\right) \times_3 \left(\mcP_{\whU_3} - \mcP_{U_3}\right), \mcA\right\rangle \right| \notag \\
\leq & \mathrm{\RN{1}} + \mathrm{\RN{2}} + \mathrm{\RN{3}} + \mathrm{\RN{4}} + \mathrm{\RN{5}} + \mathrm{\RN{6}} + \mathrm{\RN{7}} + \mathrm{\RN{8}}.
\end{align}

Applying similar arguments in the proof of Step 1.1 and by the same decomposition of $\mathrm{\RN{3}}$ in Step 2.3 in the proof of Theorem~\ref{thm: main theorem in tensor regression without sample splitting}, we can show
\begin{align*}
& \left|\left\langle \mcT\times_1 \left(\mcP_{\whU_1} - \mcP_{U_1}\right) \times_2 \left(\mcP_{\whU_2} - \mcP_{U_2}\right) \times_3 \left(\mcP_{\whU_3} - \mcP_{U_3}\right), \mcA\right\rangle \right| \\
\lesssim & \underbrace{\left\|\mcA \times_1 U_1 \times_2 U_2 \times_3 U_3\right\|_{\mathrm{F}} 
\cdot\left(\frac{\sigma^6 \oor^{1 / 2} \cdot \op^3}{\ulambda^5}\right)}_{\mathrm{\RN{1}}}\\ 
+ & \underbrace{\sum_{j=1}^3 \left\|\mcA \times_{j+1} U_{j+1} \times_{j+2} U_{j+2}\right\|_{\mathrm{F}} 
\cdot \left(\frac{\sigma^3 \oor^{1/2} \cdot \op^2\oor \log (\op)}{\ulambda^4} 
+ \frac{\sigma^5 \oor^{1/2} \cdot \op^{7 / 2}}{\ulambda^6}\right)}_{\mathrm{\RN{2},\RN{3},\RN{5}}} \\ 
+ & \underbrace{\sum_{j=1}^3 \left\|\mcA \times_j U_j\right\|_{\mathrm{F}} 
\cdot \left(\frac{\sigma^4 \oor^{1/2} \cdot \op\oor^2 \log(\op)}{\ulambda^3} 
+ \frac{\sigma^6 \oor^{1/2} \cdot \oR \op^{3/2} \log (\op)}{\ulambda^4} 
+ \frac{\sigma^6 \oor^{1/2} \cdot \op^{5/2} \sqrt{\oR \log (\op)}}{\ulambda^5}\right)}_{\mathrm{\RN{4},\RN{6},\RN{7}}} \\ 
+ & \underbrace{\left\|\mcA\right\|_{\mathrm{F}} 
\cdot \left(\frac{\sigma^3 \oor^{1/2} \cdot \oor^{3/2}\log(\op)^{3/2}}{\ulambda^2} 
+ \frac{\sigma^4 \oor^{1/2} \cdot \oR^{3 / 2} \op^{1 / 2} \log (\op)^{3 / 2}}{\ulambda^3} 
+ \frac{\sigma^5 \oor^{1/2} \cdot \oR \op^{3 / 2} \log (\op)}{\ulambda^4}\right)}_{\mathrm{\RN{8}}} . 
\end{align*}

\subsection*{Step 3: Analysis of asymptotic normal terms} 

Let 
$$
\mcP_{\mathbb{T}_{\mcT} \mathcal{M}_{\mathbf{r}}}\left(\mcA\right) := \mcP_{\mathbb{T}_{\mcT} \mathcal{M}_{(r_1, r_2, r_3)}}\left(\mcA\right) = \sum_{j=1}^3 \operatorname{Mat}_j^{-1}\left(\mcP_{U_{j\perp}}A_j\mcP_{\left(U_{j+2}\otimes U_{j+1}\right)G_j^{\top}}\right) + \mcA \times_1 \mcP_{U_1} \times_2 \mcP_{U_2} \times_3 \mcP_{U_3}.
$$
We aim to show the normal approximation of
$
\left\langle \mcZ, \mcP_{\mathbb{T}_{\mcT} \mathcal{M}_{\mathbf{r}}}\left(\mcA\right)\right\rangle.
$
To apply the Berry-Essen theorem, we calculate its second and third moments. 

\subsubsection*{Step 3.1: Second Moment of Asymptotic Normal Terms} 

Clearly,
\begin{align*}
\mathbb{E}_\mcZ \left[\left\langle \mcZ, \mcP_{\mathbb{T}_{\mcT} \mathcal{M}_{\mathbf{r}}}\left(\mcA\right)\right\rangle\right]^2
= &\sigma^2 \cdot \sum_{j=1}^{3} \left\| \mcP_{U_{j\perp}} A_j \mcP_{\left(U_{j+2} \otimes U_{j+1}\right)G_j^{\top}}\right\|_{\mathrm{F}}^2 +\sigma^2 \cdot \left\|\mcA \times_1 \mcP_{U_1} \times_2 \mcP_{U_2} \times_3 \mcP_{U_3}\right\|_{\mathrm{F}}^2 .
\end{align*}

\subsubsection*{Step 3.2:  Third Moment of Asymptotic Normal Terms}

Next, we bound the third moment. By the entrywise i.i.d assumption of the noise tenrsor $Z$, we have
\begin{align*}
\mathbb{E}\left[\left\langle \mcZ,  \mcP_{\mathbb{T}_{\mcT} \mathcal{M}_{\mathbf{r}}}\left(\mcA\right)\right\rangle\right]^3
\leq & K_3 \left[\left\|\mcA\times_1 \mcP_{U_1} \times_2 \mcP_{U_2} \times_3 \mcP_{U_3}\right\|_{\ell_\infty}+ \left\|\mcP_{U_{j\perp}}A_j\mcP_{\left(U_{j+2}\otimes U_{j+1}\right)G_j^{\top}}\right\|_{\ell_\infty} \right]\\
& \cdot \left [\sum_{j=1}^{3} \left\| \mcP_{U_{j\perp}} A_j \mcP_{\left(U_{j+2} \otimes U_{j+1}\right)G_j^{\top}}\right\|_{\mathrm{F}}^2 + \left\|\mcA \times_1 \mcP_{U_1} \times_2 \mcP_{U_2} \times_3 \mcP_{U_3}\right\|_{\mathrm{F}}^2\right] .
\end{align*}

By a nonasymptotic version of Berry-Essen theorem \citep{berry1941accuracy, esseen1956moment} and Theorem 3.7 in \citet{chen2010normal}, we get
\begin{align*}
&\sup_{x \in \mathbb{R}} \left| \mathbb{P}  \left(\frac{\sum_{j=1}^3 \left\langle Z_j,  \mcP_{U_{j\perp}} A_j \mcP_{\left(U_{j+2} \otimes U_{j+1}\right)G_1^{\top}} \right\rangle + \left\langle \mcZ,  \mcA \times_1 \mcP_{U_1} \times_2 \mcP_{U_2} \times_3 \mcP_{U_3}\right\rangle}{\sigma \left(\sum_{j=1}^3\left\|\mcP_{U_{j\perp}} A_j\mcP_{\left(U_{j+2} \otimes U_{j+1}\right)G_j^{\top}}\right\|_{\mathrm{F}}^2+ \left\|\mcA\times_1 U_1\times_2 U_2\times_3 U_3\right\|_{\mathrm{F}}^2\right)^{1/2}} \leq x\right)-\Phi(x) \right| \\
\lesssim & \frac{CK_3 \left(\sum_{j=1}^3\left\|\mcP_{U_{j\perp}} A_j\mcP_{\left(U_{j+2} \otimes U_{j+1}\right)G_j^{\top}}\right\|_{\ell_\infty}+\left\|\mcA\times_1 U_1\times_2 U_2\times_3 U_3\right\|_{\ell_\infty}\right)}{\sigma^3 \left(\sum_{j=1}^3\left\|\mcP_{U_{j\perp}} A_j\mcP_{\left(U_{j+2} \otimes U_{j+1}\right)G_j^{\top}}\right\|_{\mathrm{F}}^2+\left\|\mcA\times_1 U_1\times_2 U_2\times_3 U_3\right\|_{\mathrm{F}}^2\right)^{\frac{1} {2}}} =: C\Psi.
\end{align*}

\subsubsection*{Step 3.3: Combining Asymptotic Normal Terms and Negligible Terms} 

By the Lipshitz property of $\psi(x)$, we then have
\begin{align*}
& \sup _{x \in \mathbb{R}} \left\lvert\, \mathbb{P}\left(\frac{\left\langle\widetilde{T} \times_1 \mcP_{\whU_1} \times_2 \mcP_{\whU_2} \times_3 \mcP_{\whU_3}, \mcA\right\rangle-\langle T, \mcA\rangle}{\sigma \cdot \left(\sum_{j=1}^3 \left\|\mcP_{U_{j \perp}} A_j \mcP_{\left(U_{j+2} \otimes U_{j+1}\right)G_j^{\top}} \right\|_{\mathrm{F}}^2+ \left\|\mcA \times_1 \mcP_{U_1} \times_2 \mcP_{U_2} \times_3 \mcP_{U_3} \right\|_{\mathrm{F}}^2\right)^{1/2}} \leq x\right)-\Phi(x)\right| \\
\lesssim & \Psi + \left[\frac{1}{\op^c} + \exp\left(-c\op\right) + \exp\left(-cn\right) + \mcP\left(\mcE_\Delta\right) + \mcP\left(\mcE_U^{\text{reg}}\right)\right] \\
+ & \frac{1}{\sigma \cdot s_{\mcA}} \cdot \Bigg\{\underbrace{\left\|\mcA \times_1 U_1 \times_2 U_2 \times_3 U_3\right\|_{\mathrm{F}} \cdot \left(\frac{\sigma^2 \sqrt{\op \log (\op)}}{\ulambda}+\frac{\sigma^3 \oor^{1/2} \cdot \op^{3/2}}{\ulambda^2}\right)}_{\text{from Step 1}} \\
+ & \underbrace{\sum_{j=1}^3 \left\|U_j^{\top} A_j \mcP_{\left(U_{j+2} \otimes U_{j+1}\right) G_j^{\top}}\right\|_{\mathrm{F}} \cdot \frac{\sigma^2\oor^{1/2}\cdot \op}{\ulambda} + \sum_{j=1}^{3} \left\|\mcA \times_{j+1} U_{j+1} \times_{j+2} U_{j+2}\right\|_{\mathrm{F}} \cdot \left(\frac{\sigma^3\oor^{1/2} \cdot \op\sqrt{\oor\log(\op)}}{\ulambda^2} + \frac{\sigma^4\oor^{1/2} \cdot \op}{\ulambda^3}\right)}_{\text{shared between Step 1 and Step 2}} \\
+ & \sum_{j=1}^3 \left\|\mcP_{U_{j\perp}} A_j \mcP_{\left(U_{j+2} \otimes U_{j+1}\right) G_j^{\top}}\right\|_{\mathrm{F}} \cdot\left(\frac{\sigma^2 \oor^{1/2} \cdot \sqrt{\op \log (\op)}}{\ulambda}+\frac{\sigma^3 \oor^{1/2} \cdot \op^{3/2}}{\ulambda^2}\right) \\
+ & \sum_{j=1}^3 \left\|\mcA \times_j U_j \right\|_{\mathrm{F}} \cdot \left(\frac{\sigma^2 \oor^{1/2} \cdot \oor\log(\op)}{\ulambda}+\frac{\sigma^3 \oor^{1/2} \cdot \oR \sqrt{\op \log (\op)}}{\ulambda^2}+\frac{\sigma^4 \oor^{1/2} \cdot \op^{3 / 2} \sqrt{\oR \log (\op)}}{\ulambda^3}\right)  \\
+ & \underbrace{\left\|\mcA\right\|_{\mathrm{F}} 
\cdot \left(\frac{\sigma^3 \oor^{1/2} \cdot \oor^{3/2}\log(\op)^{3/2}}{\ulambda^2} 
+ \frac{\sigma^4 \oor^{1/2} \cdot \oR^{3 / 2} \op^{1 / 2} \log (\op)^{3 / 2}}{\ulambda^3} 
+ \frac{\sigma^5 \oor^{1/2} \cdot \oR \op^{3 / 2} \log (\op)}{\ulambda^4}\right)}_{\text{from Step 2}} \Bigg\}.
\end{align*}

\end{proof}

\section{Preliminary Upper Bounds for Tensor PCA}\label{sec: Preliminary Upper Bounds for Tensor PCA}

This section contains particular necessary preliminary upper bounds in the tensor PCA. 

After the power iteration and projection in the algorithm in Section \ref{step:step1 initialization in main algorithm}, for any $j=1,2,3$, we know that $\whU_j$ contains the top- $r_j$ eigenvectors of 
$
\widehat{\mcT}^{\text{unbs}}_j\left(\mcP_{\whU_{j+2}^{(1)}} \otimes \mcP_{\whU_{j+1}^{(1)}}\right) \widehat{\mcT}^{\text{unbs}\top}_j.
$
As a result, $\whU_j \whU_j^{\top}$ is the spectral projector for the left top- $r_j$ left eigenvectors of
\begin{align*}
\widehat{T}_j^{\text{unbs}} \left(\mcP_{\whU_{j+2}^{(1)}} \otimes \mcP_{\whU_{j+1}^{(1)}}\right) \widehat{T}_j^{\text{unbs}\top}
= T_j\left(\mcP_{U_{j+1}} \otimes \mcP_{U_{j+2}}\right) T_j^{\top}+E_j 
= U_jG_jG_j^{\top} U_j^{\top} + E_j,
\end{align*}
where
\begin{align}
E_j
= & T_j\left(\mcP_{\whU_{j+2}^{(1)}} \otimes \mcP_{\whU_{j+1}^{(1)}}\right)\whZ_j^{\top} + \whZ_j\left(\mcP_{\whU_{j+2}^{(1)}} \otimes \mcP_{\whU_{j+1}^{(1)}}\right)T_j^{\top} \notag \\
&+  T_j\left(\left(\mcP_{\whU_{j+2}^{(1)}}-\mcP_{U_{j+2}}\right) \otimes \mcP_{\whU_{j+1}^{(1)}}\right) T_j^{\top}+T_j\left(\mcP_{U_{j+2}} \otimes \left(\mcP_{\whU_{j+1}^{(1)}}-\mcP_{U_{j+1}}\right)\right) T_j^{\top} \label{eq: definition of Ej in tensor PCA}\\
&+  T_j\left(\left(\mcP_{\whU_{j+2}^{(1)}}-\mcP_{U_{j+2}}\right) \otimes \left(\mcP_{\whU_{j+1}^{(1)}}-\mcP_{U_{j+1}}\right)\right) T_j^{\top}. \notag  
\end{align}

If $\left\|E_j\right\| \leq \frac{1}{2}\ulambda^2$, then by Theorem 1 \citep{xia2021normal}, the following equation holds
$$
\whU_j \whU_j^{\top}-U_j U_j^{\top}=\sum_{k_j= 1}^{+\infty}\mathcal{S}_{G_j, k_j}\left(E_j\right).
$$
Here, for each positive integer $k$
$$
\mathcal{S}_{G_1, k_j}\left(E_1\right)=\sum_{s_1+\cdots+s_{k_j+1}=k_j}(-1)^{1+\tau(\mathbf{s})} \cdot \mcP_j^{-s_1} E_j \mcP_j^{-s_2} E_j \mcP_j^{-s_3} \cdots \mcP_j^{-s_{k_j}} E_j \mcP_j^{-s_{k_j+1}}
$$
where $s_1, \cdots, s_{k_j+1}$ are non-negative integers and $\tau(\mathbf{s})=\sum_{j=1}^{k_j+1} \mathbb{I}\left(s_{k_j}>0\right)$, 
$
\mcP_j^{-k}= U_j\left(G_jG_j^{\top}\right)^{-k}U_j^{\top}
$
for any $k \geq 1$ and $\mcP_j^{0}=U_{j\perp}U_{j\perp}^{\top}$. It follows that
\begin{align}
&\mathcal{S}_{G_j, 1}\left(E_j\right)
=  P_j^{-1} E_j P_j^{0}+P_j^{0} E_j P_j^{-1} \notag \\
= & U_j (G_j G_j^{\top})^{-1} G_j ({U_{j+1}} \otimes U_{j+2})^{\top} (\mcP_{\whU_{j+1}^{(1)}}\otimes \mcP_{\whU_{j+2}^{(1)}})\whZ_j^{\top} \mcP_{U_{j\perp}} 
+  \mcP_{U_{j\perp}} \whZ_j(\mcP_{\whU_{j+2}^{(1)}} \otimes \mcP_{\whU_{j+1}^{(1)}}) (U_{j+2} \otimes U_{j+1})G_j^{\top}(G_j G_j^{\top})^{-1} U_j^{\top}, \label{eq: definition of S1(whEj) in tensor PCA}
\end{align}
for any $j=1,2,3$, where the second equality, the third inequality come from the definition that $P_j^{-1}=U_j\left(G_jG_j^{\top}\right)^{-1}U_j^{\top}$. Here, note that $\left\|E_j\right\| \leq \kappa\ulambda \sqrt{\op}$. Then the condition, $\left\|E_j\right\| \leq \frac{1}{2}\ulambda^2$, for Theorem 1 in \citet{xia2021normal} is satisfied provided that $\ulambda \gtrsim \kappa\sqrt{\op}$, which leads to the signal-to-noise assumption in Theorem~\ref{thm: main theorem in tensor PCA}.

In the subsequent sections, we assume that the following events hold with high probability:
$$
\left\|\mcP_{\whU_j^{(0)}} - \mcP_{U_j}\right\| \leq \sigma\sqrt{\op}
$$
holds with probability at least $1- \mathbb{P}\left(\mcE_{U}^{\text{PCA}}\right)$, where event $\mcE_{U}^{\text{PCA}}$ is defined by $\mcE_{U}^{\text{PCA}}= \left\{\left\|\mcP_{\whU_j^{(0)}} - \mcP_{U_j}\right\| > \sigma\sqrt{\op}\right\}$.
Then by Lemma~\ref{lemma: error contraction of l2 error of singular space in tensor PCA}, we know that 
$
\left\|\mcP_{\whU_j^{(1)}} - \mcP_{U_j}\right\| \leq \sigma\sqrt{\op} $ and $ 
\left\|\mcP_{\whU_j} - \mcP_{U_j}\right\| \leq \sigma\sqrt{\op}
$
holds with probability at least $1-\exp(-c\op) - \mathbb{P}(\mcE_{U}^{\text{PCA}})$ for any $j=1,2,3$.

In the following sections, we established upper bounds for perturbation terms of varying orders in the spectral representation under the setting of tensor PCA. In particular, we will show that the first-order perturbation term is the leading term. 

\subsection{Preliminary Bounds in the Proof of Theorem~\ref{thm: main theorem in tensor PCA}} \label{subsec: preliminary bounds in the proof of main theorems}

\begin{proposition}\label{prop: high-probability upper bound of AoPU1p(PUhat1-P_U1)oPU2p(PUhat2-P_U2)oPU3p(PUhat3-P_U3) in tensor PCA}

Under the same setting of Theorem~\ref{thm: main theorem in tensor PCA}, with probability at least $1-\exp(-cn) - \frac{1}{\op^c} - \P\left(\mcE_\Delta\right) - \P\left(\mcE_{U}^{\text{PCA}}\right)$, where $c$ and $C$ are two universal constants, it holds that
\begin{align}
& \left\|\mcA \times_j \mcP_{U_{j \perp}}\left(\mcP_{\whU_j}-\mcP_{U_j}\right) \mcP_{U_j} \times_{j+1} \mcP_{U_{j+1 \perp}}\left(\mcP_{\whU_{j+1}}-\mcP_{U_{j+1}}\right) \mcP_{U_{j+1}} \times_{j+2} \mcP_{U_{j+2 \perp}}\left(\mcP_{\whU_{j+2}}-\mcP_{U_{j+2}}\right) \mcP_{U_{j+2}}\right\|_{\mathrm{F}} \notag \\
\lesssim & \left\|\mcA\right\|_{\mathrm{F}} \cdot\left(\frac{\sigma^3 \oor^{3} \cdot \log (\op)^{3 / 2}}{\ulambda^3}+\frac{\sigma^4 \cdot \oR^{3 / 2} \op^{1 / 2} \log (\op)^{3 / 2}}{\ulambda^4}+\frac{\sigma^5 \cdot \oR \op^{3 / 2} \log (\op)}{\ulambda^5}\right), \label{eq: high-prob upper bound of AoPU1p(PUhat1-P_U1)PU1oPU2p(PUhat2-P_U2)PU2oPU3p(PUhat3-P_U3)PU3 in tensor PCA} \\ 
& \left\|\mcA \times_j \mcP_{U_{j \perp}}\left(\mcP_{\whU_j}-\mcP_{U_j}\right) \mcP_{U_{j\perp}} \times_{j+1} \mcP_{U_{j+1 \perp}}\left(\mcP_{\whU_{j+1}}-\mcP_{U_{j+1}}\right) \mcP_{U_{j+1}} \times_{j+2} \mcP_{U_{j+2 \perp}}\left(\mcP_{\whU_{j+2}}-\mcP_{U_{j+2}}\right) \mcP_{U_{j+2}}\right\|_{\mathrm{F}} \notag \\
\lesssim & \left\|\mcA\right\|_{\mathrm{F}} \cdot\left(\frac{\sigma^4 \oor^{3} \cdot \op^{1 / 2} \log (\op)^{3 / 2}}{\ulambda^4}+\frac{\sigma^5 \oR^{3 / 2} \op \log (\op)^{3 / 2}}{\ulambda^5}+\frac{\sigma^6 \cdot \oR \op^2 \log (\op)}{\ulambda^6}\right), \label{eq: high-prob upper bound of AoPU1p(PUhat1-P_U1)PU1poPU2p(PUhat2-P_U2)PU2oPU3p(PUhat3-P_U3)PU3 in tensor PCA} \\
& \left\|\mcA \times_j \mcP_{U_{j \perp}}\left(\mcP_{\whU_j}-\mcP_{U_j}\right) \mcP_{U_j} \times_{j+1} \mcP_{U_{j+1 \perp}}\left(\mcP_{\whU_{j+1}}-\mcP_{U_{j+1}}\right) \mcP_{U_{j+1\perp}} \times_{j+2} \mcP_{U_{j+2 \perp}}\left(\mcP_{\whU_{j+2}}-\mcP_{U_{j+2}}\right) \mcP_{U_{j+2\perp}}\right\|_{\mathrm{F}} \notag \\
\lesssim & \left\|\mcA\right\|_{\mathrm{F}} \cdot\left(\frac{\sigma^5 \oor^{3} \cdot \op \log (\op)^{3 / 2}}{\ulambda^5}+\frac{\sigma^6 \cdot \oR^{3 / 2} \op^{3 / 2} \log (\op)^{3 / 2}}{\ulambda^6}+\frac{\sigma^7 \cdot \oR \op^{5 / 2} \log (\op)}{\ulambda^7}\right), \label{eq: high-prob upper bound of AoPU1p(PUhat1-P_U1)PU1oPU2p(PUhat2-P_U2)PU2poPU3p(PUhat3-P_U3)PU3p in tensor PCA} \\
& \left\|\mcA \times_j \mcP_{U_{j \perp}}\left(\mcP_{\whU_j}-\mcP_{U_j}\right) \mcP_{U_{j\perp}} \times_{j+1} \mcP_{U_{j+1 \perp}}\left(\mcP_{\whU_{j+1}}-\mcP_{U_{j+1}}\right) \mcP_{U_{j+1\perp}} \times_{j+2} \mcP_{U_{j+2 \perp}}\left(\mcP_{\whU_{j+2}}-\mcP_{U_{j+2}}\right) \mcP_{U_{j+2\perp}}\right\|_{\mathrm{F}} \notag \\
\lesssim & \left\|\mcA\right\|_{\mathrm{F}} \cdot \left(\frac{\sigma^6 \oor^{3} \cdot \op^{3/2}\log (\op)^{3/2}}{\ulambda^6}+\frac{\sigma^7 \cdot \oR^{3/2} \op^2 \log (\op)^{3/2}}{\ulambda^7}+\frac{\sigma^8 \cdot \oR \op^3 \log (\op)}{\ulambda^8}\right) \label{eq: high-prob upper bound of AoPU1p(PUhat1-P_U1)PU1oPU2p(PUhat2-P_U2)PU2poPU3p(PUhat3-P_U3)PU3pp in tensor PCA}.
\end{align}
\end{proposition}
The proof of Proposition~\ref{prop: high-probability upper bound of AoPU1p(PUhat1-P_U1)oPU2p(PUhat2-P_U2)oPU3p(PUhat3-P_U3) in tensor PCA} is similar to that of Proposition~\ref{prop: high-probability upper bound of AoPU1p(PUhat1-P_U1)oPU2p(PUhat2-P_U2)oPU3p(PUhat3-P_U3) in tensor regression without sample splitting}, and is thus omitted.

\begin{proposition}\label{prop: high-probability upper bound of AoPU1oPU2p(PUhat2-P_U2)oPU3p(PUhat3-P_U3) in tensor PCA}

Under the same setting of Theorem~\ref{thm: main theorem in tensor PCA}, with probability at least $1-\exp(-cn) - \frac{1}{\op^c} - \P\left(\mcE_\Delta\right) - \P\left(\mcE_{U}^{\text{PCA}}\right)$, where $c$ and $C$ are two universal constants ,it holds that
\begin{align}
& \left\|\mcA\times_j \mcP_{U_j} \times_{j+1} \mcP_{U_{j+1 \perp}}\left(\mcP_{\whU_{j+1}}-\mcP_{U_{j+1}}\right) \mcP_{U_{j+1}} \times_{j+2} \mcP_{U_{{j+2} \perp}}\left(\mcP_{\whU_{j+2}}-\mcP_{U_{j+2}}\right)\mcP_{U_{j+2}}\right\|_{\mathrm{F}} \notag \\
\lesssim & \left\|\mcA \times_j U_j\right\|_{\mathrm{F}} \cdot \left(\frac{\sigma^2 \oor^2 \cdot \log (\op)}{\ulambda^2}+\frac{\sigma^3 \cdot \oR \sqrt{\op \log (\op)}}{\ulambda^3}+\frac{\sigma^4 \cdot \op^{3 / 2} \sqrt{\oR \log (\op)}}{\ulambda^4}\right), \label{eq: high-prob upper bound of AoPU1oPU2p(PUhat2-P_U2)PU2oPU3p(PUhat3-P_U3)PU3 in tensor PCA} \\
& \left\|\mcA\times_j \mcP_{U_j} \times_{j+1} \mcP_{U_{j+1 \perp}}\left(\mcP_{\whU_{j+1}}-\mcP_{U_{j+1}}\right) \mcP_{U_{{j+1}\perp}} \times_{j+2} \mcP_{U_{{j+2} \perp}}\left(\mcP_{\whU_{j+2}}-\mcP_{U_{j+2}}\right)\mcP_{U_{j+2}}\right\|_{\mathrm{F}} \notag \\
\lesssim & \left\|\mcA \times_j U_j\right\|_{\mathrm{F}} \cdot \left(\frac{\sigma^3 \oor^2 \cdot \sqrt{\op}\log (\op)}{\ulambda^3}+\frac{\sigma^4 \cdot \oR \op \sqrt{\log (\op)}}{\ulambda^4}+\frac{\sigma^5 \cdot \op^2 \sqrt{\oR \log (\op)}}{\ulambda^5}\right), \label{eq: high-prob upper bound of AoPU1oPU2p(PUhat2-P_U2)PU2poPU3p(PUhat3-P_U3)PU3 in tensor PCA}\\
& \left\|\mcA\times_j \mcP_{U_j} \times_{j+1} \mcP_{U_{j+1 \perp}}\left(\mcP_{\whU_{j+1}}-\mcP_{U_{j+1}}\right) \mcP_{U_{{j+1}\perp}} \times_{j+2} \mcP_{U_{{j+2} \perp}}\left(\mcP_{\whU_{j+2}}-\mcP_{U_{j+2}}\right)\mcP_{U_{j+2\perp}}\right\|_{\mathrm{F}} \notag \\
\lesssim & \left\|\mcA \times_j U_j\right\|_{\mathrm{F}} \cdot \left(\frac{\sigma^4 \oor^2 \cdot \op \log (\op)}{\ulambda^4}+\frac{\sigma^5 \cdot \oR \op^{3/2} \sqrt{\log (\op)}}{\ulambda^5}+\frac{\sigma^6 \cdot \op^{5 / 2} \sqrt{\oR \log (\op)}}{\ulambda^6}\right). \label{eq: high-prob upper bound of AoPU1oPU2p(PUhat2-P_U2)PU2poPU3p(PUhat3-P_U3)PU3p in tensor PCA}
\end{align}
for any $j=1,2,3$.
\end{proposition}
The proof of Proposition~\ref{prop: high-probability upper bound of AoPU1oPU2p(PUhat2-P_U2)oPU3p(PUhat3-P_U3) in tensor PCA} is similar to that of Proposition~\ref{prop: high-probability upper bound of AoPU1oPU2p(PUhat2-P_U2)oPU3p(PUhat3-P_U3) without sample splitting}, and is thus omitted.

\begin{proposition}[Perturbation bound After Projection] \label{lemma: high-prob upper bound of first-order perturbation terms after projection in tensor PCA}

Under the same setting of Theorem~\ref{thm: main theorem in tensor PCA}, let $V_j \in \mathbb{R}^{p_j\times R_j}$ be a fixed matrix satisfy $\left\|V_j\right\|=1$. Then with probability at least $1-\exp(-cn) - \frac{1}{\op^c} - \P\left(\mcE_\Delta\right) - \P\left(\mcE_{U}^{\text{PCA}}\right)$, where $c$ and $C$ are two universal constants ,it holds that
\begin{align}
& \left\|V_j^{\top}\mcP_{U_{j\perp}}\left(\mcP_{\whU_j} - \mcP_{U_j}\right)U_j\right\| \leq \frac{\sigma\sqrt{\oR\log(\op)}}{\ulambda} +\frac{\sigma^3 \cdot \op^{3/2}}{\ulambda^3}, \label{eq: high-prob upper bound of V1tPU1p(PUhat1-PU1)U1 in tensor PCA} \\
& \left\|V_j^{\top}\mcP_{U_{j\perp}}\left(\mcP_{\whU_1} - \mcP_{U_1}\right)U_{j\perp}\right\| \leq \frac{\sigma^2\sqrt{\oR\log\left(\op\right)}\cdot \sqrt{\op}}{\ulambda^2}+ \frac{\sigma^3\cdot \op^{3/2}}{\ulambda^3} \label{eq: high-prob upper bound of V1tPU1p(PUhat1-PU1)U1p in tensor PCA}.
\end{align}

Furthermore, we have
\begin{align}
\left\|U_j^\top\left(\mcP_{\whU_j} - \mcP_{U_j}\right)U_j\right\| \leq \frac{\sigma^2\cdot \op}{\ulambda^2}. \label{eq: high-prob upper bound of PU1(PUhat1-PU1)PU1 in tensor PCA}
\end{align}

\end{proposition}
The proof of Lemma~\ref{lemma: high-prob upper bound of first-order perturbation terms after projection in tensor PCA} is similar to that of Proposition~\ref{prop: high-prob upper bound of first-order perturbation terms after projection in tensor regression without sample splitting}, and is thus omitted.

\subsection{Upper Bound of First-Order Perturbation Terms}

\begin{lemma}[High-probability upper bound of the first-order perturbation terms in the spectral representation in tensor PCA]\label{lemma: high-prob upper bound of the first-order perturbation terms in tensor PCA}

Under the same setting of Theorem~\ref{thm: main theorem in tensor PCA}, with probability at least $1-\exp(-cn) - \frac{1}{\op^c} - \P\left(\mcE_\Delta\right) - \P\left(\mcE_{U}^{\text{PCA}}\right)$, where $c$ and $C$ are two universal constants, it holds that
\begin{align}
\left\|\mcP_j^{-\frac{1}{2}}E_j\mcP_j^{-\frac{1}{2}}\right\| 
= & \left\|U_j\left(G_j G_j^{\top}\right)^{-\frac{1}{2}} U_j^{\top}E_jU_j\left(G_j G_j^{\top}\right)^{-\frac{1}{2}} U_j^{\top}\right\| \lesssim \frac{\sigma \cdot \sqrt{\oor \log (\op)}}{\ulambda}+ \frac{\sigma^2 \cdot \op}{\ulambda^2} \label{eq: high-prob upper bound of P1(-1/2)E1P1(-1/2) in tensor PCA} ,\\ 
\left\|\mcP_{j}^{0}E_j\mcP_j^{-\frac{1}{2}}\right\| 
= & \left\|U_{j\perp}^{\top}E_jU_j\left(G_j G_j^{\top}\right)^{-\frac{1}{2}} U_j^{\top}\right\| \lesssim \sigma \cdot \sqrt{\op}, \label{eq: high-prob upper bound of P1(0)E1P1(-1/2) in tensor PCA}\\
\left\|\mcP_{j}^{0}E_j\mcP_{j}^{0}\right\|
= & \left\|U_{j\perp}^{\top}E_jU_{j\perp}\right\| \lesssim \sigma^2 \cdot \op \label{eq: high-prob upper bound of P1(0)E1P1(0) in tensor PCA}.
\end{align}
for each $j=1,2,3$, where $\mcP_j^{-\frac{1}{2}}=U_j\left(G_jG_j^{\top}\right)^{-\frac{1}{2}}U_j^{\top}$ and $\mcP_j^{0}=U_{j\perp}U_{j\perp}^{\top}$, and $E_j$ is defined by \eqref{eq: definition of Ej in tensor PCA}.

\end{lemma}
The proof of Lemma~\ref{lemma: high-prob upper bound of the first-order perturbation terms in tensor PCA} is similar to that of Lemma~\ref{lemma: high-prob upper bound of the first-order perturbation terms in tensor regression without sample splitting}, and is thus omitted.

\begin{lemma}\label{lemma: high-prob upper bound of the projected first-order perturbation terms in tensor PCA}

Under the same setting of Theorem~\ref{thm: main theorem in tensor PCA}, let $V_j \in \mathbb{R}^{p_j\times R_j}$ be a fixed matrix satisfying $\left\|V_j\right\|=1$. Then with probability at least $1-\exp(-cn) - \frac{1}{\op^c} - \P\left(\mcE_\Delta\right) - \P\left(\mcE_{U}^{\text{PCA}}\right)$, where $c$ and $C$ are two universal constants ,it holds that
\begin{align}
\left\|V_j^{\top} \mcP_{U_{j\perp}} E_j \mcP_j^{-\frac{1}{2}}\right\| = \left\|V_j^{\top} \mcP_{U_{j\perp}} E_j U_j\left(G_jG_j^{\top}\right)^{-\frac{1}{2}}U_j^{\top}\right\| \leq \sigma \cdot \sqrt{\oR \log (\op)} + \frac{\sigma^3 \cdot \op^{3/2}}{\ulambda^2} \label{eq: high-prob upper bound of V1tP1(0)E1P1(-1/2) in tensor PCA},
\end{align}
and 
\begin{align}
\left\|V_j^{\top} \mcP_{U_{j\perp}} E_j U_{j\perp}\right\| \leq \sigma^2 \sqrt{\oR\log (\op)} \cdot \sqrt{\op}+ \frac{\sigma^3 \cdot \op^{3/2}}{\ulambda} \label{eq: high-prob upper bound of V1tP1(0)Ehat1P1(0) in tensor PCA},
\end{align}
where $E_j$ be defined by \eqref{eq: definition of Ej in tensor PCA} and $P_j^{-2s}=U_j\left(G_jG_j^{\top}\right)^{-s}U_j^{\top}$, for each $j=1,2,3$.

\end{lemma}

\begin{proof}
By symmetry, it suffice to consider $\left\|V_1^{\top}\mcP_{U_{1\perp}}E_1 U_1\right\|$ and $\left\|V_1^{\top}\mcP_{U_{1\perp}}E_1 U_{1\perp}\right\|$.

{\bf Part 1: Proof for the first inequality}

By the same decomposition in the proof of \eqref{eq: high-prob upper bound of V1tP1(0)Ehat1P1(-1/2) in tensor regression without sample splitting}, we write
\begin{align*}
\left\|V_1^{\top}\mcP_{U_{1\perp}}E_1 \mcP_1^{-\frac{1}{2}} \right\| 
\leq & \mathrm{\RN{1}} + \mathrm{\RN{2}} + \mathrm{\RN{3}} + \mathrm{\RN{4}} + \mathrm{\RN{5}} + \mathrm{\RN{6}} + \mathrm{\RN{7}} + \mathrm{\RN{8}}.
\end{align*}

First, we have 
$
\mathrm{\RN{1}} 
\lesssim \sigma\sqrt{\oR\log(\op)} 
$
and
\begin{align*}
\mathrm{\RN{2}} 
\leq \frac{1}{\ulambda} \cdot \left\|V_1^{\top}\mcP_{U_{1\perp}}Z_1\left(U_3\otimes U_2\right)\right\| \cdot \underbrace{\left\|\left(U_3\otimes U_2\right)Z_1^{\top}U_1\right\|}_{\eqref{eq: high-prob upper bound of U1Z1(U3oU2) in tensor PCA}} \lesssim \frac{\sigma^2\sqrt{\oR\oor}\log(\op)}{\ulambda}. 
\end{align*}

Then, consider
\begin{align*}
\mathrm{\RN{3}} 
\leq & \underbrace{\left\|V_1^{\top}\mcP_{U_{1 \perp}} Z_1 \left[\left(\mcP_{U_{3\perp}}Z_3\left(\mcP_{\whU_2^{(0)}} \otimes \mcP_{\whU_1^{(0)}} \right)\left(U_2\otimes U_1\right)G_3^{\top}\left(G_3G_3^{\top}\right)^{-1}\right) \otimes U_2\right]G_1^{\top}\left(G_1G_1^{\top}\right)^{-\frac{1}{2}}U_1^{\top}\right\|}_{\mathrm{\RN{3}}.\mathrm{\RN{1}}} \\
+ & \underbrace{\left\|V_1^{\top}\mcP_{U_{1 \perp}} Z_1 \left[\left(\mcP_{U_{3\perp}}Z_3\left(\mcP_{\whU_2^{(0)}} \otimes \mcP_{\whU_1^{(0)}} \right)Z_3^{\top}\left(G_3G_3^{\top}\right)^{-1}\right) \otimes U_2\right]G_1^{\top} \left(G_1G_1^{\top}\right)^{-\frac{1}{2}}U_1^{\top}\right\|}_{\mathrm{\RN{3}}.\mathrm{\RN{2}}} \\
+ & \underbrace{\left\|V_1^{\top}\mcP_{U_{1 \perp}} Z_1\left[\sum_{k_3=2}^{+\infty}\left(S_{G_3,k_3}\left(E_3^{(0)}\right)  U_3 \otimes U_2\right)\right] G_1^{\top}\left(G_1G_1^{\top}\right)^{-\frac{1}{2}}U_1^{\top}\right\|}_{\mathrm{\RN{3}}.\mathrm{\RN{3}}}.
\end{align*}

Here, first, similar to $\mathrm{\RN{3}}.\mathrm{\RN{1}}$ in the proof of \eqref{eq: high-prob upper bound of V1tP1(0)Ehat1P1(-1/2) in tensor regression without sample splitting}, we have 
\begin{align*}
\mathrm{\RN{3}}.\mathrm{\RN{1}} 
\leq & \left[\left\|V_1^{\top}\mcP_{U_{1 \perp}} Z_1 \left[\left(\mcP_{U_{3\perp}}Z_3\left(U_2\otimes U_1\right)G_3^{\top}\left(G_3G_3^{\top}\right)^{-1}U_3^{\top}\right) \otimes U_2 \right]\right\| + \frac{\sigma^3\sqrt{\op}^3}{\ulambda^2}+ \frac{\sigma^3\sqrt{\op}^3}{\ulambda^2} + \frac{\sigma^4\op^2}{\ulambda^3}\right]. 
\end{align*}

It then remains to find an upper bound for $\left\|V_1^{\top}\mcP_{U_{1 \perp}} Z_1 \left[\left(\mcP_{U_{3\perp}}Z_3\left(U_2\otimes U_1\right)G_3^{\top}\left(G_3G_3^{\top}\right)^{-1}U_3^{\top}\right) \otimes U_2 \right]\right\|$. By \eqref{eq: high-prob upper bound of |W2tZ2tBZ2(W3oW1)| in spectral norm} in Lemma~\ref{lemma: high-prob upper bound of |W2tZ2tBZ2(W3oW1)|},
\begin{align*}
& \left\|V_1^{\top} \mcP_{U_{1 \perp}} Z_1\left[\left(\mcP_{U_{3\perp}}Z_3\left(U_2\otimes U_1\right)G_3^{\top}\left(G_3G_3^{\top}\right)^{-1}\right)\otimes U_2\right]\right\| \lesssim \frac{\sigma^2\cdot \sqrt{\op\oR\log\left(\op\right)}}{\ulambda}.
\end{align*}

Second, similar to the proof of \eqref{eq: high-prob upper bound of V1tP1(0)Ehat1P1(-1/2) in tensor regression without sample splitting}, we have $
\mathrm{\RN{3}}.\mathrm{\RN{2}} 
\lesssim \frac{\sigma^3 \op^{3/2}}{\ulambda^2}. 
$

In addition, it follows immediately that
$$
\mathrm{\RN{3}}.\mathrm{\RN{3}} 
= \left\| V_1^{\top}\mcP_{U_{1 \perp}} Z_1\left(\sum_{k_1=2}^{+\infty}\left(S_{G_1,k_1}\left(E_1^{(0)}\right)  U_3 \otimes U_2\right)\right) G_1^{\top}\left(G_1G_1^{\top}\right)^{-\frac{1}{2}}U_1^{\top}\right\|_{\mathrm{F}} \lesssim \sigma\sqrt{\op}\cdot \frac{\sigma^2 \op}{\ulambda^2} = \frac{\sigma^3 \op^{3/2}}{\ulambda^2}. 
$$

Therefore, we have
\begin{align*}
\left\|V_1^{\top}\mcP_{U_{1\perp}}E_1 \mcP_1^{-\frac{1}{2}}\right\| 
\lesssim & \sigma\sqrt{\oR\log(\op)}+ \frac{\sigma^2\sqrt{\op\oR\log(\op)}}{\ulambda} + \frac{\sigma^3\cdot \op^{3/2}}{\ulambda^2} \\
\lesssim & \sigma\sqrt{\oR\log(\op)} + \frac{\sigma^3\cdot \op^{3/2}}{\ulambda^2}.
\end{align*}
where the second inequality follows from $\ulambda \geq \kappa\sigma\sqrt{\op}$.

{\bf Part 2: Proof for the second inequality}

By the same arguments in the proof of \eqref{eq: high-prob upper bound of V1tP1(0)Ehat1P1(0) in tensor regression without sample splitting} in Proposition~\ref{prop: high-prob upper bound of first-order perturbation terms after projection in tensor regression without sample splitting}, it follows that
\begin{align*}
\left\|V_1^{\top} \mcP_{U_{1 \perp}} E_1 U_{1\perp}\right\| 
\lesssim & \sigma^2\sqrt{\oR\log\left(\op\right)}\cdot \sqrt{\op}+ \frac{\sigma^3\cdot\op^{3/2}}{\ulambda}.
\end{align*}

\end{proof}

\subsection{Upper Bound of Higher-Order Perturbation Terms}

\begin{lemma}\label{lemma: high-prob upper bound of the higher-order perturbation terms in tensor PCA}

Under the same setting of Theorem~\ref{thm: main theorem in tensor PCA}, let $V_j \in \mathbb{R}^{p_j\times R_j}$ satisfy $\left\|V_j\right\|=1$, $E_j$ be defined by \eqref{eq: definition of Ej in tensor PCA}. Then with probability at least $1-\exp(-cn) - \frac{1}{\op^c} - \P\left(\mcE_\Delta\right) - \P\left(\mcE_{U}^{\text{PCA}}\right)$, where $c$ and $C$ are two universal constants ,it holds that
\begin{align}
& \left\|\mcP_{U_j} \sum_{k_j=2}^{+\infty} S_{G_j, k_j}\left(E_j\right) \mcP_{U_{j \perp}} V_j\right\| \lesssim \frac{\sigma^2\cdot \sqrt{\oR\oor}\log(\op)}{\ulambda^2} + \frac{\sigma^3\cdot \op\sqrt{\oR\log(\op)}}{\ulambda^3} + \frac{\sigma^4\op^2}{\ulambda^4}, \label{eq: high-prob upper bound of PUorder2PUpPV in tensor PCA}\\ 
& \left\|\mcP_{U_j} \sum_{k_j=3}^{+\infty} S_{G_j, k_j}\left(E_j\right) \mcP_{U_{j \perp}} V_j\right\| \lesssim \frac{\sigma^3 \cdot \op \sqrt{\oR \log (\op)}}{\ulambda^3}+\frac{\sigma^5 \cdot \op^{5/2}}{\ulambda^5}, \label{eq: high-prob upper bound of PUorder3PUpPV in tensor PCA} \\
& \left\|\mcP_{U_{j\perp}} \sum_{k_j=2}^{+\infty} S_{G_j, k_j}\left(E_j\right) \mcP_{U_{j \perp}} V_j\right\| \lesssim \frac{\sigma^2 \cdot \sqrt{\oR \op \log (\op)}}{\ulambda^2}+\frac{\sigma^4 \cdot \op^2}{\ulambda^4}, \label{eq: high-prob upper bound of PUporder2PUpPV in tensor PCA}\\
& \left\|\mcP_{U_{j \perp}} \sum_{k_j=3}^{+\infty} S_{G_j, k_j}\left(E_j\right) \mcP_{U_{j \perp}} V_j\right\| \lesssim \frac{\sigma^4 \cdot \op^{3/2}  \oR^{1/2} \log (\op)}{\ulambda^4}+\frac{\sigma^5 \cdot \op^{5/2}}{\ulambda^5}. \label{eq: high-prob upper bound of PUporder3PUpPV in tensor PCA}
\end{align}
for any $j=1,2,3$.

\end{lemma}
The proof of Lemma~\ref{lemma: high-prob upper bound of the higher-order perturbation terms in tensor PCA} is similar to that of Lemma~\ref{lemma: high-prob upper bound of the higher-order perturbation terms in tensor regression without sample splitting}, and is thus omitted.

\begin{lemma}\label{lemma: high-prob upper bound of negligible terms in the first-order term in tensor PCA}

Under the same setting of Theorem~\ref{thm: main theorem in tensor PCA}, let $V_j \in \mathbb{R}^{p_j\times R_j}$ be a fixed matrix satisfying $\left\|V_j\right\|=1$, and let $E_j$ be defined by \eqref{eq: definition of Ej in tensor PCA}. Then with probability at least $1-\exp(-cn) - \frac{1}{\op^c} - \P\left(\mcE_\Delta\right) - \P\left(\mcE_{U}^{\text{PCA}}\right)$, where $c$ and $C$ are two universal constants ,it holds that
\begin{align}
& \left\|U_{j\perp}^{\top} E_j U_{j \perp}-U_{j\perp}^{\top} Z_j\left(\mcP_{U_{j+2}} \otimes \mcP_{U_{j+1}}\right) Z_j^{\top} U_{j \perp}\right\| \lesssim \frac{\sigma^3 \op^{3/2}}{\ulambda}, \label{eq: high-prob upper bound of U1ptE1U1pt - U1ptZ1(PU3oPU2)Z1tU1p in tensor PCA} \\ 
& \left\|V_j^{\top}  \mcP_{U_{j \perp}} E_j U_{j \perp}-V_j V_j^{\top} \mcP_{U_{j \perp}} Z_j\left(\mcP_{U_{j+2}} \otimes \mcP_{U_{j+1}}\right) Z_j^{\top} U_{j \perp}\right\| \lesssim \frac{\sigma^3 \cdot \op \sqrt{\oR \log (\op)}}{\ulambda}+ \frac{\sigma^4 \cdot \op^2}{\ulambda^2} ,\label{eq: high-prob upper bound of V1V1tP1(0)E1U1pt - V1V1tP1(0)Z1(PU3oPU2)Z1tU1p in tensor PCA}\\
& \left\|U_{j \perp}^{\top} E_j P_j^{-1}-U_{j \perp}^{\top} Z_j\left(U_{j+2} \otimes U_{j+1}\right) G_j^{\top}\left(G_j G_j^{\top}\right)^{-1} U_j^{\top}\right\| \lesssim \frac{\sigma^2 \cdot \op}{\ulambda^2}, \label{eq: high-prob upper bound of U1ptE1U1P1(-1) - U1ptZ1(PU3oPU2)G1(G1G1t)(-1)U1t in tensor PCA}\\
& \left\|V_j^{\top}  \mcP_{U_{j \perp}} E_j P_j^{-1} -V_j V_j^{\top} U_{j \perp}^{\top} Z_j\left(U_{j+2} \otimes U_{j+1}\right) G_j^{\top}\left(G_j G_j^{\top}\right)^{-1} U_j^{\top} \right\| \lesssim \frac{\sigma^2 \sqrt{\oR \op }}{\ulambda^2}+\frac{\sigma^3 \cdot \op^{3/2}}{\ulambda^3} \label{eq: high-prob upper bound of V1V1tP1(0)E1U1P1(-1) - V1V1tP1(0)Z1(PU3oPU2)G1(G1G1t)(-1)U1t in tensor PCA}.
\end{align}

\end{lemma}
The proof of Lemma~\ref{lemma: high-prob upper bound of negligible terms in the first-order term in tensor PCA} is similar to that of Lemma~\ref{lemma: high-prob upper bound of negligible terms in the first-order term in tensor regression without sample splitting} and Lemma~\ref{lemma: remaining terms of the first-order perturbation term in tensor regression without sample splitting}, and is thus omitted.

\subsection{Upper Bound of Leading Terms in the Spectral Representation}

\begin{lemma}\label{lemma: high-prob upper bound of leading terms in the spectral representation across three modes in tensor PCA}

Under the same setting of Theorem~\ref{thm: main theorem in tensor PCA}, let $V_j \in \mathbb{R}^{p_j\times R_j}$ be a fixed matrix satisfying $\left\|V_j\right\|=1$, and let $E_j$ be defined by \eqref{eq: definition of Ej in tensor PCA}. Then with probability at least $1-\exp(-cn) - \frac{1}{\op^c} - \P\left(\mcE_\Delta\right) - \P\left(\mcE_{U}^{\text{PCA}}\right)$, where $c$ and $C$ are two universal constants ,it holds that
\begin{align}
& \left\|\mcA \times_j \mcP_{U_{j \perp}} E_j \mcP_j^{-1} \times_{j+1} \mcP_{U_{{j+1} \perp}} E_{j+1} \mcP_{j+1}^{-1} \times_{j+2} \mcP_{U_{{j+2} \perp}} E_{j+2} \mcP_{j+2}^{-1}\right\|_{\mathrm{F}} \notag \\
\lesssim & \left\|\mcA\right\|_{\mathrm{F}} \cdot\left(\frac{\sigma^3 \oor^{3} \cdot \log (\op)^{3/2}}{\ulambda^3}+\frac{\sigma^4 \cdot \oR^{3/2} \op^{1/2} \log (\op)}{\ulambda^4}+\frac{\sigma^5 \cdot \oR \op^{3/2} \log (\op)}{\ulambda^5}\right). \label{eq: high-prob upper bound of AoPU1pE1P1(-1)oPU2pE2P2(-1)oPU3pE3P3(-1) in tensor PCA}
\end{align}

\end{lemma}
The proof of Lemma~\ref{lemma: high-prob upper bound of leading terms in the spectral representation across three modes in tensor PCA} is similar to that of Lemma~\ref{lemma: high-prob upper bound of leading terms in the spectral representation across three modes in tensor regression without sample splitting}, and is thus omitted.



\begin{lemma}\label{lemma: high-prob upper bound of leading terms in the spectral representation across two modes in tensor PCA}

Under the same setting of Theorem~\ref{thm: main theorem in tensor PCA}, let $V_j \in \mathbb{R}^{p_j\times R_j}$ be a fixed matrix satisfying $\left\|V_j\right\|=1$, and let $E_j$ be defined by \eqref{eq: definition of Ej in tensor PCA}. Then with probability at least $1-\exp(-cn) - \frac{1}{\op^c} - \P\left(\mcE_\Delta\right) - \P\left(\mcE_{U}^{\text{PCA}}\right)$, where $c$ and $C$ are two universal constants ,it holds that
\begin{align}
& \left\|\mcA \times_j U_j \times_{j+1} \mcP_{U_{j+1\perp}}E_{j+1}\mcP_{j+1}^{-1} \times_{j+2} \mcP_{U_{j+2\perp}}E_{j+2}\mcP_{j+2}^{-1}\right\|_{\mathrm{F}} \notag \\
\lesssim & \left\|\mcA \times_j U_j\right\|_{\mathrm{F}} \cdot \left(\frac{\sigma^2 \oor^2\cdot \log (\op)}{\ulambda^2}+\frac{\sigma^3 \cdot \oR \sqrt{\op \log (\op)}}{\ulambda^3}+\frac{\sigma^4 \cdot \op^{3 / 2} \sqrt{\oR \log (\op)}}{\ulambda^4}\right) \label{eq: high-prob upper bound of AoPU1oPU2pE2P2(-1)oPU3pE3P3(-1) in tensor PCA}
\end{align}

\end{lemma}
The proof of Lemma~\ref{lemma: high-prob upper bound of leading terms in the spectral representation across two modes in tensor PCA} is similar to that of Lemma~\ref{lemma: high-prob upper bound of leading terms in the spectral representation across two modes in tensor regression without sample splitting}, and is thus omitted.



\section{Concentration Inequalities for Tensor PCA}

\begin{lemma}\label{lemma: high-prob upper bound of spectral norm of B1tZ1(1)tA1t(A2Z2(1)B2oA3Z3(1)B3) in tensor PCA}

Let $\mcZ \in \mathbb{R}^{p_1 \times p_2 \times p_3}$ be a random tensor with i.i.d. mean-zero $\sigma$-sub-Gaussian entries. Let $\mcG \in \mathbb{R}^{r_1 \times r_2 \times r_3}$ be a given tensor. Let $Z_j = \operatorname{Mat}_j\left(\mcZ\right)$ and Let $G_j = \operatorname{Mat}_j\left(\mcG\right)$ denote the mode-$j$ matricization of the tensor $\mcZ$ and $\mcG$. Suppose that $U_j \in \mathbb{O}^{p_j \times r_j}$. Additionally, let $U_{j\perp} \in \mathbb{O}^{p_j \times (p_j - r_j)}$ be orthonormal matrices, such that $U_{j\perp}^{\top} U_{j\perp} = I_{p_j - r_j}$. The matrix $U_{j\perp} U_{j\perp}^{\top} \in \mathbb{R}^{p_j \times p_j}$ is a projection matrix that projects any vector onto the orthogonal complement of the space spanned by $U_j U_j^{\top}$.

Then it holds that
\begin{align}
& \mathbb{P}\left(\left\|\left(\mcP_{U_{j+2}} \otimes \mcP_{U_{j+1}}\right) Z_j^{\top}\mcP_{U_{j\perp}}A_j \left(\mcP_{U_{j+1 \perp}}Z_{j+1}\left(\mcP_{U_j} \otimes \mcP_{U_{j+2}}\right) \otimes \mcP_{U_{j+2 \perp}}Z_{j+2}\left(\mcP_{U_{j+1}} \otimes \mcP_{U_j}\right)\right) \right\| \right. \notag \\
& \quad \left. \geq C \left\|\mcP_{U_{1\perp}}A_1\left(\mcP_{U_{3\perp}}\otimes \mcP_{U_{1\perp}}\right)\right\|_{\mathrm{F}} \cdot t\right) \leq 7^{r_j + r_{j+1}r_{j+2}}\exp\left[-c\min\left(t^2,t^\frac{2}{3}\right)\right] . \label{eq: high-prob upper bound of Ao(PU1pZ1(PU3oPU2))o(PU2pZ2(PU1oPU3))o(PU3pZ3(PU2oPU1)) in tensor PCA}
\end{align}
\end{lemma}

\begin{proof}

By symmetry, it suffices to find a high-probability upper bound for
\begin{align*}
& \left\|\left(\mcP_{U_3} \otimes \mcP_{U_2}\right) Z_1^{\top}\mcP_{U_{1\perp}}A_1\left(\mcP_{U_{3\perp}}Z_3\left(\mcP_{U_2} \otimes \mcP_{U_1}\right) \otimes  \mcP_{U_{2 \perp}}Z_2\left(\mcP_{U_1} \otimes \mcP_{U_3}\right) \right) \right\|.
\end{align*}

Let
\begin{align*}
& C_1 = \left(\mcP_{U_{3\perp}}\otimes \mcP_{U_{2\perp}}\right)A_1^{\top}\mcP_{U_{1\perp}} \in \mathbb{R}^{p_2p_3\times p_1}, A_2 = \mcP_{U_{2\perp}} \in \mathbb{R}^{p_2\times p_2}, A_3 = \mcP_{U_{3\perp}} \in \mathbb{R}^{p_3\times p_3} \\
& B_1 = \mcP_{U_3} \otimes \mcP_{U_2} \in \mathbb{R}^{p_2p_3\times p_2p_3}, B_2 = \mcP_{U_1} \otimes \mcP_{U_3} \in \mathbb{R}^{p_1p_3\times p_1p_3}, B_3 = \mcP_{U_2} \otimes \mcP_{U_1} \in \mathbb{R}^{p_1p_2\times p_1p_2}.
\end{align*}

Then, it suffices to find an upper bound for
\begin{align*}
& \left\|B_1^{\top}Z_1^{\top}C_1^{\top}\left(A_3Z_3B_3\otimes A_2Z_2B_2\right)\right\|.
\end{align*}

Note that
$$
\left\|B_1^{\top}Z_1^{\top}C_1^{\top}\left(A_3Z_3B_3\otimes A_2Z_2B_2\right)\right\|=\sup_{\substack{u \in \mathbb{R}^{p_2p_3}, \left\|u\right\|_{\ell_2}=1\\ v \in \mathbb{R}^{p_1p_2\cdot p_2p_3}, \left\|v\right\|_{\ell_2}=1}}u^{\top} B_1^{\top} Z_1^{\top} C_1^{\top}\left(A_3Z_3B_3 \otimes A_2Z_2B_2\right) v.
$$

Therefore, we first find a high-probability upper bound for 
$$
u^{\top} B_1^{\top} Z_1^{\top} C_1^{\top}\left(A_3Z_3B_3 \otimes A_2Z_2B_2\right) v
$$
with two given $u\in \mathbb{R}^{p_2p_3}, v\in \mathbb{R}^{p_1p_2\cdot p_2p_3}$ and then apply an $\varepsilon$-net argument to derive the high-probability upper bound for $\left\|B_1^{\top} Z_1^{\top} C_1^{\top}\left(A_3Z_3B_3 \otimes A_2Z_2B_2\right)\right\|$.

Since 
\begin{align*}
\left[\mcZ\right]_{i_1,i_2,i_3}=\left[Z_1\right]_{i_1,i_2+p_2\left(i_3-1\right)}, \left[\mcZ\right]_{j_1,j_2,j_3}=\left[Z_2\right]_{j_2,j_3+p_3\left(j_1-1\right)},
\left[\mcZ\right]_{k_1,k_2,k_3}=\left[Z_3\right]_{k_3,k_1+p_1\left(k_2-1\right)}
\end{align*}
and
\begin{align*}
\left[Z_1\otimes Z_2\right]_{p_2(i_1-1)+j_2, p_3p_1\left[i_2+p_2(i_3-1)-1\right]+\left[j_3+p_3\left(j_1-1\right)\right]}=\left[Z_1\right]_{i_1,i_2+p_2(i_3-1)}\left[Z_2\right]_{j_2,j_3+p_3(j_1-1)},
\end{align*}
where $i_1,j_1,k_1=1,2,\cdots,p_1, i_2,j_2,k_2=1,2,\cdots,p_2, i_3,j_3,k_3=1,2,\cdots,p_3$,
it follows that
\begin{align*}
& \left[C_1Z_1B_1\right]_{a_1, b_1} =  \sum_{i_1,i_2,i_3}\left[A_1\right]_{a_1,i_1}\left[Z_1\right]_{i_1,i_2+p_2\left(i_3-1\right)}\left[B_1\right]_{i_2+p_2\left(i_3-1\right),b_1}, \\
& \left[A_2Z_2B_2\right]_{a_2,b_2}=  \sum_{j_1,j_2,j_3}\left[A_2\right]_{a_2,j_2}\left[Z_2\right]_{j_2,j_3+p_3\left(j_1-1\right)}\left[B_2\right]_{j_3+p_{3}\left(i_1-1\right), b_2}, \\
& \left[A_3Z_3B_3\right]_{a_3,b_3} =  \sum_{k_1,k_2,k_3}\left[A_3\right]_{a_3,k_3}\left[Z_3\right]_{k_3,k_1+p_1\left(k_2-1\right)}\left[B_3\right]_{k_1+p_1\left(k_3-1\right), b_3}.
\end{align*}

Thus, we have
\begin{align*}
& u^{\top}B_1^{\top}Z_1^{\top}C_1^{\top}\left(A_3Z_3B_3\otimes A_2Z_2B_2\right)v \\
= & \sum_{\substack{l,a_2,a_3,b_2,b_3\\i_1,i_2,i_3\\j_1,j_2,j_3\\k_1,k_2,k_3}} \left[C_1\right]_{p_2\left(a_3-1\right)+a_2, i_1}\left[B_1\right]_{p_2\left(i_3-1\right)+i_2,l} \left[A_2\right]_{a_2, j_2}\left[B_2\right]_{j_3+p_3\left(j_1-1\right), b_2}
 \left[A_3\right]_{a_3, k_3}\left[B_3\right]_{k_1+p_1\left(k_2-1\right), b_3} \\
 & \cdot u_l v_{p_2\left(b_3-1\right)+b_2} \cdot \left[\mcZ\right]_{i_1,i_2,i_3} \left[\mcZ\right]_{j_1,j_2,j_3} \left[\mcZ\right]_{k_1,k_2,k_3}.
\end{align*}

By symmetry, we consider dividing the summation above into the following cases by index:
\begin{align*}
\text{Case \RN{1}}= & \left(i_1,i_2,i_3\right) = \left(j_1,j_2,j_3\right) = \left(k_1,k_2,k_3\right) \\
\text{Case \RN{2}}:=& \left(i_1,i_2,i_3\right) = \left(j_1,j_2,j_3\right) \neq \left(k_1,k_2,k_3\right) \\
\text{Case \RN{3}}:=& \left(i_1,i_2,i_3\right) = \left(k_1,k_2,k_3\right) \neq \left(j_1,j_2,j_3\right)\\
\text{Case \RN{4}}:=& \left(j_1,j_2,j_3\right) = \left(k_1,k_2,k_3\right) \neq \left(i_1,i_2,i_3\right)\\
\text{Case \RN{5}}:=& \left(i_1,i_2,i_3\right) \neq \left(j_1,j_2,j_3\right) \neq \left(k_1,k_2,k_3\right).
\end{align*}
Correspondingly, write 
$$
u^{\top}B_1^{\top}Z_1^{\top}C_1^{\top}\left(A_3Z_3B_3\otimes A_2Z_2B_2\right)v= f_1\left(\mcZ\right) + f_2\left(\mcZ\right) + f_3\left(\mcZ\right) + f_4\left(\mcZ\right) + f_5\left(\mcZ\right).
$$
Here, $f_i(\mcZ), i=1,2,\cdots,5$ are polynomials of entries $\left\{\left[\mcZ\right]_{j_1,j_2,j_3}\right\}_{j_1,j_2,j_3}^{p_1,p_2,p_3}$ of random tensor $\mcZ$. We then apply Theorem 1.5 in \citet{gotze2021concentration} to find high-probability upper bounds of $f_i(\mcZ)$'s.

{\bf Case \RN{1}:} $\left(i_1,i_2,i_3\right) = \left(j_1,j_2,j_3\right) = \left(k_1,k_2,k_3\right)$

For Case \RN{1}, we have
\begin{align*}
f_1\left(\mcZ\right)= \sum_{\substack{l,a_2,a_3,b_2,b_3\\i_1,i_2,i_3}}
&  
\left[C_1\right]_{p_2\left(a_3-1\right)+a_2, i_1} 
\left[B_1\right]_{p_2\left(i_3-1\right)+i_2,l} 
\left[A_2\right]_{a_2, i_2} \left[B_2\right]_{i_3+p_3\left(i_1-1\right), b_2} \\
& \left[A_3\right]_{a_3, i_3} 
\left[B_3\right]_{i_1+p_1\left(i_2-1\right), b_3} \cdot u_l v_{p_2\left(b_3-1\right)+b_2} 
\cdot \left[\mcZ\right]_{i_1,i_2,i_3}^3.
\end{align*}
Then, similar to the proof of Lemma~\ref{lemma: high-prob upper bound of tr(BZ(2)tCZ(1)) in tensor regression with sample splitting} and Lemma~\ref{lemma: high-prob upper bound of tr(BZ(2)tCZ(2)) in tensor regression with sample splitting}, we can show, $\mathbb{E}\left(f_1\left(\mcZ\right)\right) =0$, $\left\| \mathbb{E}\nabla_Z^2f_1\left(\mcZ\right)\right\|_{\mathrm{HS}}^2=0$, and 
\begin{align*}
\left\| \mathbb{E}\nabla_{\mcZ}f_1\left(\mcZ\right)\right\|_{\mathrm{HS}}^2
\leq & 9\sigma^4 \left\|C_1\right\|_{\mathrm{F}}^2 \left\|B_1u\right\|_{\ell_\infty}^2 \left\|\left(B_3\otimes B_2\right)v\right\|_{\ell_\infty}^2 \leq 9\sigma^4 \left\|C_1\right\|_{\mathrm{F}}^2 \left\|B_1\right\|_{2,\infty}^2 \left\|B_3\otimes B_2\right\|_{2,\infty}^2, \\
\left\|  \mathbb{E}\nabla_Z^3f_1\left(\mcZ\right)\right\|_{\mathrm{HS}}^2 
\leq & 36 \left\|C_1\right\|_{\mathrm{F}}^2 \left\|B_1\right\|_{2,\infty}^2 \left\|B_3 \otimes B_2\right\|_{2,\infty}^2.
\end{align*}

Combining all the results above, by Theorem 1.5 in \citet{gotze2021concentration}, it then follows that
\begin{align*}
& \mathbb{P}\left(\left|f_1\left(\mcZ\right)\right| \geq C t\right) \\
\leq & \exp \left\{-c\min\left[\left(\frac{t}{\sigma^3\left\|C_1\right\|_{\mathrm{F}}\left\|B_1\right\|_{2,\infty}\left\|B_3\otimes B_2\right\|_{2,\infty}}\right)^2, \left(\frac{t}{\sigma^3\left\|C_1\right\|_{\mathrm{F}}\left\|B_1\right\|_{2,\infty}\left\|B_3\otimes B_2\right\|_{2,\infty}}\right)^{\frac{2}{3}}\right]\right\}.
\end{align*}

{\bf Case \RN{2}: $(i_1,i_2,i_3) = \left(j_1,j_2,j_3\right) \neq \left(k_1,k_2,k_3\right)$}

For Case \RN{2}, let 
\begin{align*}
f_2\left(\mcZ\right) =
\sum_{\substack{l,a_2,a_3,b_2,b_3\\i_1,i_2,i_3\\k_1,k_2,k_3\\i_1,i_2,i_3 \neq \left(k_1,k_2,k_3\right)}} 
& \left[A_2\right]_{a_2, i_2} \left[B_2\right]_{i_3+p_3\left(i_1-1\right), b_2}
\left[A_3\right]_{a_3, k_3} 
\left[B_3\right]_{k_1+p_1\left(k_2-1\right), b_3} \\
& \cdot \left[C_1\right]_{p_2\left(a_3-1\right)+a_2, i_1} 
\left[B_1\right]_{p_2\left(i_3-1\right)+i_2,l} \cdot u_l v_{p_2\left(b_3-1\right)+b_2} 
\cdot \left[\mcZ\right]_{i_1,i_2,i_3}^2
\left[\mcZ\right]_{k_1,k_2,k_3}.
\end{align*}
Then, by the independence between $\left[\mcZ\right]_{i_1,i_2,i_3}$ and $\left[\mcZ\right]_{k_1,k_2,k_3}$ and the mean-zero property, it readily follows that $ \mathbb{E}f_2\left(\mcZ\right)=0.$ Similar to the proof of Case I, we can show, $\left\|\mathbb{E}\nabla_{\mcZ}\left(f_2\left(Z\right)\right) \right\|_{\mathrm{HS}}^2=0$, $\left\|\mathbb{E}\nabla_{\mcZ}^2\left(f_2\left(Z\right)\right) \right\|_{\mathrm{HS}}^2=0$, and
\begin{align*}
& \left\|\mathbb{E}\nabla_{\mcZ}^3\left(f_2\left(Z\right)\right) \right\|_{\mathrm{HS}}^2 \le
 4 \left\|C_1\right\|_{\mathrm{F}}^2 \left\|B_1u\right\|_{\ell_\infty}^2 \left\|\left(B_3\otimes B_2\right)v\right\|_{\ell_2}^2 \leq 4 \left\|C_1\right\|_{\mathrm{F}}^2 \left\|B_1\right\|_{2,\infty}^2 \left\|B_2\right\|^2 \left\|B_3\right\|^2.
\end{align*}

Combining all the results above and apply Theorem 1.5 in \citet{gotze2021concentration}, we have
\begin{align*}
& \mathbb{P}\left(\left|f_2(\mcZ)\right| \geq Ct\right) \leq \exp\left(-c \left(\frac{t}{\sigma^3 \left\|C_1\right\|_{\mathrm{F}}  \left\|B_1\right\|_{2,\infty}\left\|B_2\right\|\left\|B_3\right\|}\right)^{\frac{2}{3}}\right).
\end{align*}

{\bf Case \RN{3}: $(i_1,i_2,i_3) = \left(k_1,k_2,k_3\right) \neq \left(j_1,j_2,j_3\right)$}

By symmetry, we have
\begin{align*}
& \mathbb{P}\left(\left|f_3(\mcZ)\right| \geq Ct\right) \leq \exp\left(-c \left(\frac{t}{\sigma^3 \left\|C_1\right\|_{\mathrm{F}}  \left\|B_1\right\|_{2,\infty}\left\|B_2\right\|\left\|B_3\right\|}\right)^{\frac{2}{3}}\right).
\end{align*}

{\bf Case \RN{4}: $(j_1,j_2,j_3) = \left(k_1,k_2,k_3\right) \neq \left(i_1,i_2,i_3\right)$}

For the Case $\mathrm{\RN{4}}$, we have
\begin{align*}
f_4\left(\mcZ\right)
= \sum_{\substack{l,a_2,a_3,b_2,b_3\\i_1,i_2,i_3\\j_1,j_2,j_3\\i_1,i_2,i_3 \neq \left(j_1,j_2,j_3\right)}} 
& \left[C_1\right]_{p_2\left(a_3-1\right)+a_2, i_1}
\left[B_1\right]_{p_2\left(i_3-1\right)+i_2,l}\cdot \left[A_2\right]_{a_2, j_2}\left[B_2\right]_{j_3+p_3\left(j_1-1\right), b_2}
\left[A_3\right]_{a_3, j_3} 
\left[B_3\right]_{j_1+p_1\left(j_2-1\right), b_3} \\
& \cdot u_l v_{p_2\left(b_3-1\right)+b_2} \cdot \left[\mcZ\right]_{i_1,i_2,i_3} 
\left[\mcZ\right]_{j_1,j_2,j_3}^2.
\end{align*}
Obviously, $\mathbb{E}\left[f_4\left(\mcZ\right)\right]=0.$ Similar to the proof of Case I, we can show $\left\|\mathbb{E}\left[\nabla_{\mcZ} f_4(\mcZ)\right]\right\|_{\mathrm{F}}^2 = 0$, $\left\|\mathbb{E}[\nabla_\mcZ^2(f_4(\mcZ)) ]\right\|_{\mathrm{HS}}=0$, and
\begin{align*}
& \left\|\mathbb{E} \nabla_{\mcZ}^3\left(f_4(Z)\right)\right\|_{\mathrm{HS}}^2 
\leq 4\left\|C_1\right\|_{\mathrm{F}}^2\left\|B_1\right\|^2\left\|B_2\right\|^2\left\|B_3\right\|^2.
\end{align*}

It follows that
$$
\mathbb{P}\left(\left|f_4(\mcZ)\right| \geq Ct\right) \leq \exp\left(-c\left(\frac{t^2}{\left\| C_1\right\|_{\mathrm{F}} \left\| B_1\right\| \left\| B_2\right\| \left\| B_3 \right\|}\right)^{\frac{2}{3}}\right).
$$

{\bf Case \RN{5}}. Again consider
\begin{align*}
f_5\left(\mcZ\right)
= \sum_{\substack{l,a_2,a_3,b_2,b_3\\i_1,i_2,i_3\\ \neq j_1,j_2,j_3\\ \neq k_1,k_2,k_3}}
&  
\left[C_1\right]_{p_2\left(a_3-1\right)+a_2, i_1}
\left[B_1\right]_{p_2\left(i_3-1\right)+i_2,l} \cdot  \left[A_2\right]_{a_2, j_2} \left[B_2\right]_{j_3+p_3\left(j_1-1\right), b_2}
\left[A_3\right]_{a_3, k_3} 
\left[B_3\right]_{k_1+p_1\left(k_2-1\right), b_3} \\
& \cdot u_l v_{p_2\left(b_3-1\right)+b_2} \cdot \left[\mcZ\right]_{i_1,i_2,i_3} 
\left[\mcZ\right]_{j_1,j_2,j_3} 
\left[\mcZ\right]_{k_1,k_2,k_3}.
\end{align*}
It follows immediately that
$
\left\|\mathbb{E}\left(\nabla_{\mcZ}\left(f_4(\mcZ)\right)\right)\right\|_{\mathrm{HS}} = \left\|\mathbb{E}\left(\nabla_{\mcZ}^2\left(f_4(\mcZ)\right)\right)\right\|_{\mathrm{HS}} = 0.
$
Similar to the previous cases, we can show
\begin{align*}
& \left\|\mathbb{E}\left(\nabla_{\mcZ}^{3} f_5\left(\mcZ\right)\right)\right\|_{\mathrm{HS}}^2 
\leq 4 \left\|C_1\right\|_{\mathrm{F}}^2 \left\|B_1\right\|^2 \left\|B_2\right\|^2 \left\|B_3\right\|^2 .
\end{align*}

Therefore, we have
$$
\mathbb{P}\left(\left|f_5(\mcZ)\right| \geq C t\right) \leq \exp \left(-c\left(\frac{t^2}{\sigma^3 \left\|C_1\right\|_{\mathrm{F}}\left\|B_1\right\|\left\|B_2\right\|\left\|B_3\right\|}\right)^{\frac{2}{3}}\right).
$$

Combining all the results in Case \RN{1}, \RN{2}, \RN{3}, \RN{4} and \RN{5}, above, for any given $u,v \in \mathbb{R}$ we have
\begin{align*}
& \mathbb{P}\left(\left|u^{\top} B_1^{\top} Z_1^{\top} C_1^{\top}\left(A_3Z_3B_3 \otimes A_2Z_2B_2\right) v\right| \geq K_3 \left\|C_1\right\|_{\mathrm{F}} \cdot \left\|B_1\right\| \cdot \left\|B_2\right\| \cdot \left\|B_3\right\| + C t\right) \\
\leq & \exp \left\{-c \min \left[\left(\frac{t}{\sigma^3\left\|C_1\right\|_{\mathrm{F}}\left\|B_1\right\|_{2, \infty}\left\|B_3 \otimes B_2\right\|_{2, \infty}}\right)^2,\left(\frac{t}{\sigma^3\left\|C_1\right\|_{\mathrm{F}}\left\|B_1\right\|\left\|B_2\right\|\left\|B_3\right\|}\right)^{\frac{2}{3}}\right]\right\}.
\end{align*}

Here, note that
$
B_1=\mcP_{\left(U_3 \otimes U_2\right) G_1^{\top}}
$
is of rank $r_1$, 
$
B_2=\left(U_3 \otimes U_1\right) G_2^{\top}\left(G_2 G_2^{\top}\right)^{-1} U_2^{\top}
$
is of rank $r_2$, and similarly $B_3$ is of rank $r_3$. Then there exists $U_1 \in \mathbb{O}^{p_2p_3\times r_1}, V_2 \in \mathbb{O}^{p_2\times r_2}, V_3 \in \mathbb{O}^{p_3\times r_3}$, such that $B_1= B_1U_1U_1^{\top}, B_2= B_2V_2V_2^{\top}$, and $B_1= B_3V_3V_3^{\top}$. Let $\widetilde{a}=U_1^{\top}a$ and $\widetilde{b}=\left(V_2\otimes V_3\right)^{\top}b$. It then follows that
\begin{align*}
a^{\top} B_1^{\top} Z_1^{\top} C_1^{\top}\left(A_3Z_3B_3 \otimes A_2Z_2B_2\right) v 
= &  \widetilde{a}U_1^{\top} B_1^{\top} Z_1^{\top} C_1^{\top}\left(A_3Z_3B_3 \otimes A_2Z_2B_2\right) \left(V_2\otimes V_3\right)\widetilde{b}.
\end{align*}

By Lemma 5.2 of \citet{vershynin2010introduction} , there exists $\mathbb{N}^{r_1} $, a $\frac{1}{3}$-net of $\left\{a \in \mathbb{R}^{r_1}:\left\|A\right\|=1\right\}$, such that $\left|\mathbb{N}^{r_1} \right| \leq 7^{r_1} $ and $\mathbb{N}^{r_2r_3} $, a $\frac{1}{3}$-net of $\left\{b \in \mathbb{R}^{r_2r_3}:\|b\|_2=1\right\}$, such that $\left|\mathbb{N}^{r_2r_3} \right| \leq 7^{r_2r_3} $. Then, applying a union bound, we have
\begin{align*}
& \mathbb{P}\left(\sup _{\substack{a, \in N_{R_2}\\b \in N_{R_3R_1}}} \left|\widetilde{a}U_1^{\top} B_1^{\top} Z_1^{\top} C_1^{\top}\left(A_3Z_3B_3 \otimes A_2Z_2B_2\right) \left(V_2\otimes V_3\right)\widetilde{b}\right| \geq Ct\right) \\
\leq & 2 \cdot 7^{r_1+r_2r_3} \exp \left\{-c \min \left[\left(\frac{t}{\sigma^3\left\|C_1\right\|_{\mathrm{F}}\left\|B_1\right\|\left\|B_2\right\|\left\|B_3\right\|}\right)^2,\left(\frac{t}{\sigma^3\left\|C_1\right\|_{\mathrm{F}}\left\|B_1\right\|\left\|B_2\right\|\left\|B_3\right\|}\right)^{\frac{2}{3}}\right]\right\}.
\end{align*}

Then for any $a \in \mathbb{S}^{R_2}$, $b \in \mathbb{S}^{R_3R_1}$ there exist $\overline{a} \in \mathbb{N}_{R_2}$ and $\overline{b} \in \mathbb{N}_{R_3R_1}$ such that $\|\overline{a}-a\| \leq \frac{1}{3}, \|\overline{b}-b\| \leq \frac{1}{3}$. Therefore,
\begin{align*}
\left\|B_1^{\top} Z_1^{\top} C_1^{\top}\left(A_3Z_3B_3 \otimes A_2Z_2B_2\right)\right\|
\le & \frac{9}{2} \sup_{\substack{\overline{a} \in \mathbb{N}^{r_1}, \left\|\overline{a}\right\|=1\\ \overline{b} \in \mathbb{N}^{r_2r_3}, \left\|\overline{b}\right\|=1}} \left|\overline{a}U_1^{\top} B_1^{\top} Z_1^{\top} C_1^{\top}\left(A_3Z_3B_3 \otimes A_2Z_2B_2\right) \left(V_2\otimes V_3\right)\overline{b}\right|,
\end{align*}
which leads to the following desired result.

\end{proof}

\begin{lemma}
Let $\mcZ \in \mathbb{R}^{p_1 \times p_2 \times p_3}$ be a random tensor with i.i.d. mean-zero $\sigma$-sub-Gaussian entries. Let $Z_j = \operatorname{Mat}_j(\mcZ)$ denote the mode-$j$ matricization of the tensor $\mcZ$. Suppose that $U_j \in \mathbb{O}^{p_j \times r_j}$. Additionally, let $U_{j\perp} \in \mathbb{O}^{p_j \times (p_j - r_j)}$ be orthonormal matrices, such that $U_{j\perp}^{\top} U_{j\perp} = I_{p_j - r_j}$. The matrix $U_{j\perp} U_{j\perp}^{\top} \in \mathbb{R}^{p_j \times p_j}$ is a projection matrix that projects any vector onto the orthogonal complement of the space spanned by $U_j U_j^{\top}$.

Then it holds that
\begin{align}
& \mathbb{P}\left(\left\|\left(\mcP_{U_{j+2}} \otimes \left(\mcP_{U_j}\otimes \mcP_{U_{j+2}}\right)Z_{j+1}\mcP_{U_{j+1\perp}}\right)A_j^{\top}\mcP_{U_{j\perp}}Z_j\left(U_{j+2}\otimes U_{j+1}\right)\right\| \right. \notag \\
& \geq \left. C \sigma^2 \cdot\left\|\left(\mcP_{U_{j+1\perp}} \otimes \mcP_{U_{j+2}}\right) A_j^{\top} \mcP_{U_{j \perp}}\right\|_{\mathrm{F}} \cdot t\right) \leq 7^{r_j + r_{j+1}r_{j+2}} \exp\left(-c \min \left(t^2, t\right)\right) . \label{eq: high-prob upper bound of AoPU1o(PU2pZ2(PU1oPU3))o(PU3pZ3(PU2oPU1)) in tensor PCA}
\end{align}
\end{lemma}

\begin{proof}

By symmetry, it suffices to consider
$$
\left\|\left[\mcP_{U_3} \otimes \left(\left(\mcP_{U_1}\otimes \mcP_{U_3}\right)Z_2^{\top}\mcP_{U_{2\perp}}\right)\right]C_1^{\top}\mcP_{U_{1\perp}}Z_1\left(\mcP_{U_3}\otimes \mcP_{U_2}\right)\right\|.
$$
Let
\begin{align*}
& C_1 = \left(\mcP_{U_3} \otimes \mcP_{U_{2\perp}}\right)C_1^{\top}\mcP_{U_{1\perp}} \in \mathbb{R}^{p_2p_3\times p_1}, A_2 = \mcP_{U_{2\perp}} \in \mathbb{R}^{p_2\times p_2}, \\
& B_1 = \mcP_{U_3}\otimes \mcP_{U_2} \in \mathbb{R}^{p_3p_2\times p_3p_2}, B_2 = \mcP_{U_1}\otimes \mcP_{U_3} \in \mathbb{R}^{p_1p_3\times p_1p_3}.
\end{align*}
Then 
\begin{align*}
& \left[\mcP_{U_3} \otimes \left(\left(\mcP_{U_1}\otimes \mcP_{U_3}\right)Z_2^{\top}\mcP_{U_{2\perp}}\right)\right]C_1^{\top}\mcP_{U_{1\perp}}Z_1\left(\mcP_{U_3}\otimes \mcP_{U_2}\right) = \left(\mcP_{U_3} \otimes B_2^{\top}Z_2^{\top}A_2^{\top}\right)C_1Z_1B_1 .
\end{align*}


We then first find an upper bound for $v\left(\mcP_{U_3} \otimes B_2^{\top}Z_2^{\top}A_2^{\top}\right)C_1Z_1B_1u$, 
where $u$ and $v$ are two given vectors, then apply the $\varepsilon$-net argument.

Since $\left[\mcZ\right]_{i_1,i_2,i_3}=\left[Z_1\right]_{i_1,p_2\left(i_3-1\right)+i_2}$  
and
$\left[\mcZ\right]_{j_1,j_2,j_3}=\left[Z_2\right]_{j_2,p_3\left(j_1-1\right)+j_3}$,
it follows that
\begin{align*}
\left[C_1Z_1B_1\right]_{a_1, b_1} = & \sum_{i_1,i_2,i_3}\left[C_1\right]_{a_1,i_1}\left[Z_1\right]_{i_1,i_2+p_2\left(i_3-1\right)}\left[B_1\right]_{i_2+p_2\left(i_3-1\right),b_1}, \\
\left[A_2Z_2B_2\right]_{a_2,b_2} = & \sum_{j_1,j_2,j_3}\left[A_2\right]_{a_2,j_2}\left[Z_2\right]_{j_2,j_3+p_3\left(j_1-1\right)}\left[B_2\right]_{j_3+p_{3}\left(i_1-1\right), b_2} .
\end{align*}

Thus, we have
\begin{align*}
& v\left(\mcP_{U_3} \otimes B_2^{\top}Z_2^{\top}A_2^{\top}\right)C_1Z_1B_1u \\
= & \sum_{\substack{a_2,b_2,a_3,b_3\\ i_1,i_2,i_3\\ j_1,j_2,j_3}} \left[C_1\right]_{p_2\left(a_3-1\right)+a_2, i_1}\left[B_1\right]_{p_2\left(i_3-1\right)+i_2, l} 
 \left[A_2\right]_{a_2,j_2}\left[B_2\right]_{j_3+p_3\left(j_1-1\right),b_2}\left[\mcP_{U_3}\right]_{b_3,a_3} u_l v_{p_2\left(b_3-1\right)+b_2} \cdot \left[Z\right]_{i_1,i_2,i_3}\left[Z\right]_{j_1,j_2,j_3}.
\end{align*}
Note that $ \left[\mcZ\right]_{i_1,i_2,i_3}$ and $ \left[\mcZ\right]_{j_1,j_2,j_3}$ are dependent if and only if $i_1=j_1, i_2=j_2, i_3=j_3$. We then write
\begin{align*}
& v\left(\mcP_{U_3} \otimes B_2^{\top}Z_2^{\top}A_2^{\top}\right)C_1Z_1B_1u \\
= & \sum_{\substack{a_2,b_2,a_3,b_3\\ i_1,i_2,i_3}} \left[C_1\right]_{p_2\left(a_3-1\right)+a_2, i_1}\left[B_1\right]_{p_2\left(i_3-1\right)+i_2, l}   \left[A_2\right]_{a_2,i_2}\left[B_2\right]_{i_3+p_3\left(i_1-1\right),b_2} 
 \left[\mcP_{U_3}\right]_{b_3,a_3} u_l v_{p_2\left(b_3-1\right)+b_2}  \left[Z\right]_{i_1,i_2,i_3}^2 \\
+ & \sum_{\substack{a_2,b_2,a_3,b_3\\ i_1,i_2,i_3\\ j_1,j_2,j_3\\ i_1,i_2,i_3 \neq \left(j_1,j_2,j_3\right)}} \left[C_1\right]_{p_2\left(a_3-1\right)+a_2, i_1}\left[B_1\right]_{p_2\left(i_3-1\right)+i_2, l}   \left[A_2\right]_{a_2,j_2}\left[B_2\right]_{j_3+p_3\left(j_1-1\right),b_2} 
 \left[\mcP_{U_3}\right]_{b_3,a_3} u_l v_{p_2\left(b_3-1\right)+b_2}  \left[Z\right]_{i_1,i_2,i_3}\left[Z\right]_{j_1,j_2,j_3}.
\end{align*}

Then,
\begin{align*}
& \mathbb{E}\left(\operatorname{tr}\left[C_1Z_1B_1\left(B_2^{\top}Z_2^{\top}A_2^{\top}\right) \otimes \mcP_{U_3}\right]\right) \\
= & \sigma^2\sum_{\substack{a_2,b_2,a_3,b_3\\ i_1,i_2,i_3}} \left[C_1\right]_{p_2\left(a_3-1\right)+a_2, i_1}\left[B_1\right]_{p_2\left(i_3-1\right)+i_2, l} \left[A_2\right]_{a_2,i_2}\left[B_2\right]_{i_3+p_3\left(i_1-1\right),b_2}\left[\mcP_{U_3}\right]_{b_3,a_3} u_l v_{p_2\left(b_3-1\right)+b_2}.
\end{align*}
The rest part of the proof is similar to that of Lemma~\ref{lemma: high-prob upper bound of spectral norm of B1tZ1(1)tA1t(A2Z2(1)B2oA3Z3(1)B3) in tensor PCA}, and is thus omitted.

\end{proof}

\begin{lemma}\label{lemma: high-prob upper bound of tr(BZtCZ)}
Suppose that $Z \in \mathbb{R}^{p_1\times p_2 \times p_3}$ is a random tensor with i.i.d. sub-Gaussian entries of variance $\sigma^2$ and $Z_j=\operatorname{\text{Mat}}_j\left(Z\right)$ is the mode-$j$ matricization of the tensor $Z$. 

For two matrices $\ddot{B} \in \mathbb{R}^{p_{j+1}p_{j+2}\times p_j}, \ddot{C} \in \mathbb{R}^{p_j\times p_j}$, it holds that
\begin{align}
\mathbb{P}\left(\left|\operatorname{tr}\left[\ddot{B} Z_j^{\top} \ddot{C} Z_j\right]-\sigma^2 \operatorname{tr}\left[\ddot{B}\right] \operatorname{tr}\left[\ddot{C}\right]\right| \geq C t\right) \leq \exp \left(-c \min \left(\frac{t^2}{\sigma^4\left\|\ddot{B}\right\|_{\mathrm{F}}^2\left\|\ddot{C}\right\|_{\mathrm{F}}^2}, \frac{t}{\sigma^2\left\|\ddot{B}\right\|_{\ell_{\infty}}\left\|\ddot{C}\right\|_{\ell_{\infty}}}\right)\right) \label{eq: high-prob upper bound of tr(BZtCZ)}.
\end{align}

\end{lemma}

\begin{proof}

Without loss of generality, we start by considering $Z_j=Z_1$. Note that
\begin{align}
\operatorname{tr}\left[\ddot{B}Z_1^{\top}\ddot{C}Z_1\right] 
= & \sum_{\substack{i=1,2,\cdots, p_2p_3\\ k=1,2,\cdots, p_1}} \left[\ddot{B}\right]_{i,i} \left[\ddot{C}\right]_{k,k}\left[Z_1\right]_{k,i}^2 + \sum_{\substack{(i, j) \neq (k, l) \\ i,j=1,2,\cdots, p_2p_3\\ k,l=1,2,\cdots, p_1}}\left[\ddot{B}\right]_{i,j} \left[Z_1\right]_{k,j} \left[\ddot{C}\right]_{k,l}\left[Z_1\right]_{l,i} .
\end{align}

For the quadratic terms, first note that
\begin{align*}
\mathbb{E}\left(\sum_{\substack{i=1,2,\cdots, p_2p_3\\ k=1,2,\cdots, p_1}} \left[\ddot{B}\right]_{i,i} \left[\ddot{C}\right]_{k,k}\left[Z_1\right]_{k,i}^2\right)
= \sigma^2 \sum_{\substack{i=1,2,\cdots, p_2p_3\\j=1,2,\cdots, p_1}} \left[\widetilde{B}\right]_{j,i}\left[\widetilde{C}\right]_{j,i}
= \sigma^2 \operatorname{tr}\left[\widetilde{B}\right] \operatorname{tr}\left[\widetilde{C}\right].
\end{align*}
Furthermore, since 
$
\left\|[Z_1]_{j,i}^2\right\|_{\psi_1} \leq \left\|[Z_1]_{j,i}\right\|_{\psi_2}^2 \leq \sigma^2,
$
by the Bernstein's inequality, it follows that
\begin{align*}
& \mathbb{P}\left(\left|\sum_{\substack{i=1,2,\cdots, p_2p_3\\ k=1,2,\cdots, p_1}} \left[\ddot{B}\right]_{i,i} \left[\ddot{C}\right]_{k,k}\left[Z_1\right]_{k,i}^2 - \sigma^2 \operatorname{tr}\left[\widetilde{B}\right] \operatorname{tr}\left[\widetilde{C}\right]\right| \geq Ct\right) \\
\le & \exp\left[c\min\left(\frac{t^2}{\sigma^4 \left\|\ddot{B}\right\|_{\mathrm{F}}^2\left\|\ddot{C}\right\|_{\mathrm{F}}^2}, \frac{t}{\sigma^2\left\|\ddot{B}\right\|_{\ell_\infty}\left\|\ddot{C}\right\|_{\ell_\infty}}\right)\right].
\end{align*}

Then, consider the summation of independent terms. Write 
$$
S_{\text{indep}}= \sum_{\substack{(i, j) \neq (k, l) \\ i,j=1,2,\cdots, p_2p_3\\ k,l=1,2,\cdots, p_1}}\left[\ddot{B}\right]_{i,j} \left[Z_1\right]_{k,j} \left[\ddot{C}\right]_{k,l}\left[Z_1\right]_{l,i}.
$$
By the decoupling method (Remark 6.1.3, \citet{vershynin2018high}) and the same arguments as in the proof of the first inequality, we have
\begin{align*}
\mathbb{E}_{\mcZ}\exp{\left[\lambda S_{\text{indep}} \right]}
\leq &  \mathbb{E}_{\mcZ}\exp{C\left[\lambda \left[\ddot{B}\right]_{i,j} \left[Z_1\right]_{k,j} \left[\ddot{C}\right]_{k,l}\left[Z_1^{\prime}\right]_{l,i}\right]}.
\end{align*}

Then by the comparison lemma (Lemma 6.2.3, \citet{vershynin2018high}), when 
$$
\lambda^2 \leq \frac{c^2}{\sigma^4\max_{i,k}\left[\ddot{B}\right]^2\max_{k,l}\left[\ddot{C}\right]^2}= \frac{c^2}{\sigma^4\left\|\ddot{B}\right\|_{\ell_\infty}^2\left\|\ddot{C}\right\|_{\ell_\infty}^2},
$$ 
it holds that
\begin{align*}
\mathbb{E}_{\mcZ}\exp{\left[\lambda S_{\text{indep}} \right]}
\leq & \mathbb{E}_{\mcZ}\exp{\left[C\lambda \sum_{\substack{i,j=1,2,\cdots, p_2p_3,\\ k,l=1,2,\cdots, p_1 }} \left[G_1\right]_{k,j}\left[\ddot{B}\right]_{i,j}\left[G_1^{\prime}\right]_{l,i}\left[\ddot{C}\right]_{k,l}\right]} 
\le \exp\left[C\lambda^2 \sigma^4 \left\|\ddot{B}\right\|_{\mathrm{F}}^2 \left\|\ddot{C}\right\|_{\mathrm{F}}^2 \right] ,
\end{align*}
where $G_1 \in \mathbb{R}^{p_1\times p_2p_3}$ and $G_1^{\prime} \in \mathbb{R}^{p_1\times p_2p_3}$ are two independent matrices with i.i.d. Gaussian entries. Furthermore, since
$
\mathbb{P}\left\{ S_{\text{indep}}  \geq \frac{t}{2}\right\} \leq \exp (-\lambda \frac{t}{2}) \mathbb{E} \exp (\lambda  S_{\text{indep}}),
$
it follows that
$$
\mathbb{P}\left\{ S_{\text{indep}}  \geq \frac{t}{2}\right\} \leq \exp\left[-\frac{\lambda t}{2}+C\lambda^2 \sigma^4 \left\|\ddot{B}\right\|_{\mathrm{F}}^2 \left\|\ddot{C}\right\|_{\mathrm{F}}^2 \right] . 
$$
Optimizing when $0\leq \lambda \leq \frac{c^2}{\sigma^2\left\|\ddot{B}\right\|_{\ell_\infty}\left\|\ddot{C}\right\|_{\ell_\infty}}$, we have
$$
\mathbb{P}\left\{ S_{\text{indep}}  \geq \frac{t}{2} \right\} \leq \exp \left(-c \min \left(\frac{t^2}{\sigma^4\left\|\ddot{B}\right\|_{\mathrm{F}}^2\left\|\ddot{C}\right\|_{\mathrm{F}}^2}, \frac{t}{\sigma^2 \left\|\ddot{B}\right\|_{\ell_\infty}\left\|\ddot{C}\right\|_{\ell_\infty}}\right)\right)
$$

Combining the bounds for the diagonal and off-diagonal terms, it holds that \eqref{eq: high-prob upper bound of tr(BZtCZ)}.

\end{proof}

\begin{lemma}\label{lemma: high-prob upper bound of |W2tZ2tBZ2(W3oW1)|}

Let $\mcZ \in \mathbb{R}^{p_1 \times p_2 \times p_3}$ be a random tensor with i.i.d. mean-zero $\sigma$-sub-Gaussian entries. Let $Z_j = \operatorname{Mat}_j(\mcZ)$ denote the mode-$j$ matricization of the tensor $\mcZ$. Suppose that $W_{j+1} \in \mathbb{R}^{p_{j+1}\times R_{j+1}}, W_{j+2} \in \mathbb{R}^{p_{j+2}\times R_{j+2}}$, $B \in \mathbb{R}^{p_j\times p_j}, \widetilde{W}_j \in \mathbb{R}^{p_{j+1}p_{j+2}\times R_j}$. 

Furthermore, suppose that $\widetilde{W}_2^{\top}\left(W_1\otimes W_3\right)=0$. Then, it holds that
\begin{align}
& \mathbb{P}\left(\left\|\widetilde{W}_j^{\top} Z_j^{\top} \widetilde{B}^{\top} Z_j\left(W_{j+2} \otimes W_{j+1}\right)\right\| \geq Ct\right) \notag \\
\leq& 2 \cdot 7^{R_j+R_{j+1} R_{j+2}} \exp \left(-c \min \left(\frac{t^2}{\sigma^4\left\|\widetilde{B}\right\|_{\mathrm{F}}^2\left\|\widetilde{W}_j\right\|^2\left\|W_{j+1}\right\|^2\left\|W_{j+2}\right\|^2}, \frac{t}{\sigma^2\left\|\widetilde{B}\right\|_{\ell_{\infty}}\left\|\widetilde{W}_j\right\|\left\|W_{j+1}\right\|\left\|W_{j+2}\right\|}\right)\right). \label{eq: high-prob upper bound of |W2tZ2tBZ2(W3oW1)| in spectral norm}
\end{align}

\end{lemma}

\begin{proof}

Let $\widetilde{W}_2=\left(\mcP_{U_1} \otimes \mcP_{U_3}\right) \in \mathbb{R}^{p_1p_3\times p_1p_3}$, $W_3=\mcP_{U_{3\perp}} \in \mathbb{R}^{p_3\times p_3}$ and $W_1 = \mcP_{U_{1\perp}} \in \mathbb{R}^{p_1\times p_1}$. Note that
\begin{align*}
a^{\top}\widetilde{W}^{\top}_2Z_2^{\top}\widetilde{B}^{\top}Z_2\left(W_1\otimes W_3\right)b 
= & \underbrace{\sum_{i,j}\widetilde{a}_i\widetilde{b}_i\left[\widetilde{B}^{\top}\right]_{j,j}\left[Z_2\right]_{j,i}^2}_{\text{diagonal terms}} + \underbrace{\sum_{(i,j)\neq (l,k)}\widetilde{a}_i\left[Z_2\right]_{j,i}\left[\widetilde{B}^{\top}\right]_{j,k}[Z_2]_{k,l}\widetilde{b}_l}_{\text{off-diagonal terms}} .
\end{align*}
By applying the same arguments used in the proof of Lemma~\ref{lemma: high-prob upper bound of tr(BZtCZ)} to the diagonal and off-diagonal terms separately, and utilizing the $\varepsilon$-net argument, we obtain the desired bounds. The details are omitted for brevity.

\end{proof}

The following Lemma~Lemma~\ref{lemma: concentration bound of |U1tildeX1(U3tildeoU2tilde)|} implies that
\begin{align}
\left\|U_j^{\top} Z_j\left(U_{j+2} \otimes U_{j+1}\right)\right\| \leq \sqrt{\oor\log(\op)} \label{eq: high-prob upper bound of U1Z1(U3oU2) in tensor PCA}
\end{align}
hold with probability at least $1-\op^{-c}$, and
\begin{align}
\left\|U_{j\perp}^{\top}Z_j\left(U_{j+2} \otimes U_{j+1}\right)\right\| \leq \sqrt{\op} \label{eq: high-prob upper bound of U1pZ1(U3oU2) in tensor PCA}
\end{align} 
hold with probability at least $1-\exp\left(-c\op\right)$  for a constant $c>0$. In addition, we have
\begin{align}
\left\|V_j^{\top} \mcP_{U_j\perp} Z_j\left(U_{j+2} \otimes U_{j+1}\right)\right\| \leq \sqrt{\oR\log(\op)} \label{eq: high-prob upper bound of V1PU1pZ1(U3oU2) in tensor PCA}
\end{align}
hold with probability at least $1-\op^{-c}$.

\begin{lemma}\label{lemma: concentration bound of |U1tildeX1(U3tildeoU2tilde)|}
Suppose $\mcZ \in \mathbb{R}^{p_1 \times p_2 \times p_3}$ is a tensor with independent, zero-mean, $\sigma$-sub-Gaussian entries. For $j = 1, 2, 3$, let $\wtU_j \in \mathbb{O}^{p_j \times \widetilde{r}_j}$, where $\widetilde{r}_j \leq p_j$ for all $j = 1, 2, 3$. Let $Z_j$ denote the mode-$j$ matricization of the random tensor $\mcZ$. Then, it holds that
\begin{align}
&\mathbb{P}\left(\left\|\wtU_j^{\top} Z_j\left(\wtU_{j+2} \otimes \wtU_{j+1}\right)\right\| \geq 2 \sigma \sqrt{\widetilde{r}_j + t}\right) \leq 2 \cdot 5^{\widetilde{r}_{j+1} \widetilde{r}_{j+2}} \exp \left[-c \min\left(\frac{t^2}{\widetilde{r}_j}, \frac{t}{\prod_{j=1}^3\left\|\wtU_j\right\|_{2,\infty}}\right)\right].
\end{align}

\end{lemma}

\begin{proof}
By symmetry, it suffices to consider an upper bound for
$
\left\|\wtU_1^{\top} Z_1\left(\wtU_3 \otimes \wtU_2\right)\right\|.
$
Note that
\begin{align}
& \left\|\wtU_1^{\top}Z_1\left(\wtU_3 \otimes \wtU_2\right)b\right\|_{\ell_2}^2 =  b^{\top}\left(\wtU_3 \otimes \wtU_2\right)^{\top}Z_1^{\top}\wtU_1\wtU_1^{\top}Z_1\left(\wtU_3 \otimes \wtU_2\right)b .
\end{align}
We then apply the decoupling method for the quadratic form to prove the Hansen-Wright type bound. 
Write $\widetilde{b}:=\left(\wtU_3 \otimes \wtU_2\right)b \in \mathbb{R}^{\widetilde{p}_2\widetilde{p}_3}$. It follows that
\begin{align*}
\left\|\wtU_1^{\top}Z_1\left(\wtU_3 \otimes \wtU_2\right)b\right\|^2 
= & \underbrace{\sum_{\substack{i=1,2, \cdots, p_2p_3,\\ j=1,2, \cdots, p_1}} \widetilde{b}_i^2\left[\mcP_{\wtU_1}\right]_{j,j} \left[Z_1\right]_{j,i}^2}_{\text{quadratic terms}} + \underbrace{\sum_{\substack{(i,j)\neq (l,k)\\ i,l=1,2,\cdots, p_2p_3\\ j,k=1,2, \cdots, p_1}}\widetilde{b}_i\widetilde{b}_l\left[\mcP_{\wtU_1}\right]_{j,k} \left[Z_1\right]_{j,i}\left[Z_1\right]_{k,l}}_{\text{summation of independent terms}}.
\end{align*}
By applying the same arguments used in the proof of Lemma~\ref{lemma: high-prob upper bound of tr(BZtCZ)} to the two terms separately, and utilizing the $\varepsilon$-net argument, we obtain the desired bounds. The details are omitted for brevity.

\end{proof}

\begin{lemma}\label{lemma: concentration bound of sup |U1tildeX1(U3tildeoU2tilde)|}
Suppose $\mcZ \in \mathbb{R}^{p_1 \times p_2 \times p_3}$ is a tensor with independent, zero-mean, $\sigma$-sub-Gaussian entries. For $j = 1, 2, 3$, let $Z_j= \Mat_j\left(\mcZ\right)$ be the mode-$j$ matricization of $\mcZ$ and $\wtU_j \in \mathbb{R}^{p_j \times \widetilde{r}_j}$ such that $\|\wtU_j\| = 1$. Then, the following inequalities hold:
\begin{align}
\mathbb{P} \left(\sup_{\substack{\wtU_j \in \mathbb{R}^{p_j\times \widetilde{r}_j}, \\ \left\|\wtU_j\right\| \leq 1, \\j=1,2,3}}\left\|\wtU_j^{\top} Z_j\left(\wtU_{j+2} \otimes \wtU_{j+1}\right)\right\| \geq 4 \sigma  \sqrt{\widetilde{r}_1+t} \right)
& \leq 2 \cdot 33^{2\sum_{j=1}^3p_j\widetilde{r}_j}5^{\widetilde{r}_{j+1}\widetilde{r}_{j+2}} \exp \left[-c \min \left(\frac{t^2}{\widetilde{r}_j}, \frac{t}{\prod_{j=1}^3\left\|\wtU_j\right\|_{2,\infty}}\right)\right]. \label{eq: high-prob upper bound of sup_U1,U2,U3 |U1tildeX1(U3tildeoU2tilde)|}
\end{align}
\end{lemma}
The proof of Lemma~\ref{lemma: concentration bound of sup |U1tildeX1(U3tildeoU2tilde)|} is similar to that of Lemma~\ref{lemma: concentration bound of |U1tildeX1(U3tildeoU2tilde)|}, and is thus omitted.

\section{Other Technical Lemmas}

\begin{lemma}[Error contraction of singular space estimation in Tensor regression]\label{lemma: error contraction of l2 error of singular space in tensor regression}

Let $\whU_j^{(0)}\whU_j^{(0)\top}$, $j=1,2,3$, denote the initial estimate of the singular space of the $j$-th mode of the signal tensor $\mcT$. Conditioning on the following event
\[
\left\| \whU^{(0)}_j\whU^{(0)\top}_j - U_jU_j^{\top} \right\| \leq \frac{1}{2}
\]
holds. Then, with probability at least $1-\exp\left(-cn\right) - \exp\left(-cp\right)$, the $k$-th iteration ($k=1,2$) of the singular space estimation, $\whU_j^{(1)}\whU_j^{(1)\top}$, produced by the algorithm without sample splitting in Section \ref{sec:debias_nonsplit} or the algorithm with sample splitting in Section \ref{sec: Debiased Estimator of Linear Functionals with Sample Splitting}, satisfies the following bounds for $j=1,2,3$:
\begin{align}
\left\|\whU^{(k)}_j\whU_j^{(k)\top} - U_jU_j^{\top} \right\| &\leq \frac{\sigma_\xi}{\ulambda\sigma} \cdot \sqrt{\frac{\op}{n}} \label{eq: error contraction of l2 error of singular space in Tensor regression} .
\end{align}

\end{lemma}


The proof is essentially the same as the proof of Theorem 1 in \citet{zhang2018tensor}, with the noise tensor under the tensor PCA setting is replaced by $\widehat{\mcZ}$.

\begin{lemma}[Theorem 1, \citet{zhang2018tensor}, Error contraction of singular space estimation in Tensor PCA]\label{lemma: error contraction of l2 error of singular space in tensor PCA}

Let $\whU_j^{(0)}\whU_j^{(0)\top}$, $j=1,2,3$, denote the initial estimate of the singular space of the $j$-th mode of the signal tensor $\mcT$. Conditioning on the following event
\[
\left\| \whU^{(0)}_j\whU^{(0)\top}_j - U_jU_j^{\top} \right\| \leq \frac{1}{2}.
\]
holds. Then, with probability at least $1 - \exp\left(-cp\right)$, the $k$-th iteration ($k=1,2$) of the singular space estimation, $\whU_j^{(1)}\whU_j^{(1)\top}$, produced by the algorithm without sample splitting in Section \ref{sec:debias_nonsplit} or the algorithm with sample splitting in Section \ref{sec: Debiased Estimator of Linear Functionals with Sample Splitting}, satisfies the following bounds for $j=1,2,3$:
\begin{align}
& \left\|\whU^{(1)}_j\whU_j^{(1)\top} - U_jU_j^{\top} \right\| \leq \frac{\sigma\cdot \sqrt{\op}}{\ulambda} \label{eq: error contraction of l2 error of singular space in tensor PCA}.
\end{align}

\end{lemma}

Finally, we prove the following sub-multiplicative property of $\ell_{2,\infty}$-norm.

\begin{lemma}\label{lemma: submultiplicative of l2-infty norm in kronecker product}

For $A \in \mathbb{R}^{p_{\text{A}} \times r_{\text{A}}}$ and $B \in \mathbb{R}^{p_{\text{B}} \times r_{\text{B}}}$, it holds that
\begin{align}
\left\|A \otimes B\right\|_{2,\infty} \leq & \left\|A\right\|_{2,\infty} \cdot \left\|B\right\|_{2,\infty}. \label{eq: submultiplicative of l2,infty norm in kronecker product}
\end{align}

\end{lemma}

\begin{proof}

Let $\mcA \in \mathbb{R}^{p_{\text{A}} \times r_{\text{A}}}$ and $\mcB \in \mathbb{R}^{p_{\text{B}} \times r_{\text{B}}}$. Then 
\begin{align*}
& \left\|A\otimes B\right\|_{2,\infty} 
= \max_{i \in \left[p_{\text{A}}p_{\text{B}}\right]} \left\|e_{i}^{\top}\left(A \otimes B\right)\right\|_{\ell_2} = \max_{\substack{j \in \left[p_{\text{A}}\right]\\ k \in \left[p_{\text{B}}\right]}} \left\|\left(e_j \otimes e_k\right) \left(A \otimes B\right)\right\|_{\ell_2} 
\leq \left\|A\right\|_{2,\infty} \cdot \left\|B\right\|_{2,\infty}.
\end{align*}
\end{proof}

\section{Proof of Asymptotic Normality with Plug-in Estimates}\label{sec: Proof of Asymptotic Normality with Plug-in Estimates}

In this section, we present proofs of the following two theorems Theorem~\ref{thm: generalized main theorem in tensor regression with plug-in estimates} and Theorem~\ref{thm: generalized main theorem in tensor PCA with plug-in estimates}, which are the non-asymptotic versions of asymptotic normality for tensor regression (Theorem~\ref{thm: main theorem in tensor regression with plug-in estimates}) and tensor PCA (Theorem~\ref{thm: main theorem in tensor PCA with plug-in estimates}) with plug-in estimates, respectively.

\begin{theorem}[Non-asymptotic version of Theorem~\ref{thm: main theorem in tensor regression with plug-in estimates}]\label{thm: generalized main theorem in tensor regression with plug-in estimates}

Under the same setting of Theorem~\ref{thm: main theorem in tensor regression without sample splitting}, let the estimate of variance components $\widehat{\sigma}_\xi^2$, $\widehat\sigma^2$ and $\widehat{s}_{\mcA}^2$ be defined as in \eqref{eq: estimate of sigmaxi2 in tensor regression without sample splitting},  \eqref{eq: estimate of sigma2 in tensor regression} and \eqref{eq: estimate of sA2}, respectively. Then
\begin{align*}
& \sup _{x \in \mathbb{R}} \left\lvert\, \mathbb{P}\left(\frac{ \sqrt{n} \big( \langle\widehat{\mcT}, \mcA\rangle-\langle \mcT, \mcA\rangle \big) }{ (\widehat{\sigma}_{\xi}/\widehat{\sigma}) \widehat{s}_{\mcA} } \leq x\right)-\Phi(x)\right| \\
\lesssim & \underbrace{\sqrt{\frac{1}{n}}}_{\text{rate of asymptotic normal terms}} + \underbrace{\frac{\widetilde{\Omega}_1 + \widetilde{\Omega}_2 }{(\sigma_{\xi}/\sigma) s_{\mcA} \sqrt{\frac{1}{n}}}}_{\text{rate of negligible terms}} + \underbrace{\big[\op^{-c} + e^{-cn} + \mcP_{\mcE_{U}^{\text{reg}}} + \mcP_{\mcE_{\Delta}}\big]}_{\text{rate of initial estimates}} ,
\end{align*}
where $c>0$ is a constant.
Here,
\begin{align*}
\widetilde{\Omega}_1
= & \sum_{j=1}^3 \left\|\mcP_{U_j} A_j \mcP_{\left(U_{j+2} \otimes U_{j+1}\right) G_j^{\top}}\right\|_{\mathrm{F}} \cdot  \frac{\sigma_\xi^2}{\sigma^2}\cdot  \frac{\oor^{1/2}}{\ulambda}\cdot \frac{\op}{n} \\
+ & \sum_{j=1}^3\left\|\mcA \times_{j+1} U_{j+1} \times_{j+2} U_{j+2}\right\|_{\mathrm{F}} \cdot \left[\frac{\sigma_{\xi}^2}{\sigma^2} \cdot \frac{\oor^{1/2}}{\ulambda} \frac{\sqrt{\op\log(\op)}}{n} + \frac{\sigma_\xi^3}{\sigma^3} \cdot \frac{\oor^{1/2}}{\ulambda^2} \cdot \Delta\cdot \frac{\op^{3/2}}{n^{3/2}}\right] \\
+ & \sum_{j=1}^3\left\|\mcA \times_j U_j\right\|_{\mathrm{F}} \cdot  \frac{\sigma_\xi^2}{\sigma^2} \left[\frac{ \oor^{1 / 2}}{\ulambda } \cdot\left(\frac{\sqrt{\oR} \log (\op)}{n}+\Delta \cdot \frac{\sqrt{\oR \op \log (\op)}}{n}+\Delta^2 \cdot \frac{\op}{n}\right)\right] \\
+ & \left\|\mcA\right\|_{\mathrm{F}} \cdot  \frac{\sigma_\xi^3}{\sigma^3} \left[\frac{\oor^{1 / 2}}{\ulambda^2 } \cdot\left(\frac{\oor^{3/2}\log(\op)^{3/2}}{n^{3 / 2}}+\Delta \cdot \frac{\op^{1 / 2} \oR \log (\op)}{n^{3 / 2}}+\Delta^2 \cdot \frac{\op \oR \log (\op)}{n^{3 / 2}}+\Delta^3 \cdot \frac{\op^{3 / 2}}{n^{3 / 2}}\right)\right], \\
\widetilde{\Omega}_2 =&  \left(\sum_{j=1}^3\left\|\mcP_{U_{j\perp}}A_j\mcP_{\left(U_{j+2}\otimes U_{j+1}\right)G_j^{\top}}\right\|_{\mathrm{F}} + \left\|\mcA\times_1 U_1\times_2 U_2\times_3 U_3\right\|_{\mathrm{F}}\right) \cdot \Delta \cdot \frac{\sigma_\xi\oor^{1/2}}{\sigma} \sqrt{\frac{\op}{n}}.
\end{align*}

Furthermore, in the sample-splitting case where the setting is specified in Theorem~\ref{thm: main theorem in tensor regression with sample splitting}, let the estimate of variance components $\widehat{\sigma}_\xi^2$, $\widehat\sigma^2$ and $\widehat{s}_{\mcA}^2$ be defined as in \eqref{eq: estimate of sigmaxi2 in tensor regression with sample splitting},  \eqref{eq: estimate of sigma2 in tensor regression} and \eqref{eq: estimate of sA2}, respectively. Then we have
\begin{align*}
\widetilde{\Omega}_1
= & \sum_{j=1}^3 \left\|U_j^{\top} A_j \mcP_{\left(U_{j+2} \otimes U_{j+1}\right) G_j^{\top}}\right\|_{\mathrm{F}} \cdot \frac{\sigma_\xi^2\oor^{1/2}}{\ulambda\sigma^2}\cdot \frac{\op}{n} 
+  \sum_{j=1}^{3} \left\|\mcA \times_{j+1} U_{j+1} \times_{j+2} U_{j+2} \right\|_{\mathrm{F}} \cdot \frac{\sigma_{\xi}^2\oor^{1/2}}{\ulambda\sigma^2}  \left(\sqrt{\frac{\op\log(\op)}{n}} + \Delta \cdot \frac{\op}{n}\right) \\
+ & \sum_{j=1}^3\left\|\mcA \times_j U_j\right\|_{\mathrm{F}} \cdot\left[\frac{\sigma_{\xi}^2 \oor^{1 / 2}}{\ulambda \sigma^2} \cdot\left(\frac{\oR^{1/2}\log(\op)}{n}+\Delta \cdot \frac{\sqrt{\oR \op \log (\op)}}{n}+\Delta^2 \cdot \frac{\op}{n}\right)\right] \\
+ & \left\|\mcA\right\|_{\mathrm{F}} \cdot\left[\frac{\sigma_{\xi}^3 \oor^{1 / 2}}{\ulambda^2 \sigma^3} \cdot\left(\frac{\oR\log(\op)^{3/2}}{n^{3/2}} + \Delta \cdot \frac{\op^{1 / 2} \oR \log (\op)}{n^{3 / 2}}+\Delta^2 \cdot \frac{\op \oR \log (\op)}{n^{3 / 2}}+\Delta^3 \cdot \frac{\op^{3 / 2}}{n^{3 / 2}}\right)\right] ,\\
\widetilde{\Omega}_2 
= & \left(\sum_{j=1}^3\left\|\mcP_{U_{j\perp}}A_j\mcP_{\left(U_{j+2}\otimes U_{j+1}\right)G_j^{\top}}\right\|_{\mathrm{F}} + \left\|\mcA\times_1 U_1\times_2 U_2\times_3 U_3\right\|_{\mathrm{F}}\right) \left(\Delta \cdot \frac{\sigma_\xi\oor^{1/2}}{\sigma} \sqrt{\frac{\oor\log(\op)}{n}}\right) .
\end{align*}

\end{theorem}

\begin{theorem}[Non-asymptotic version of Theorem~\ref{thm: main theorem in tensor PCA with plug-in estimates}]\label{thm: generalized main theorem in tensor PCA with plug-in estimates}

Under the same setting of Theorem~\ref{thm: main theorem in tensor PCA}, let the estimate of variance components $\widehat\sigma^2$ and $\widehat{s}_{\mcA}^2$ be defined as in \eqref{eq: estimate of sigma2 in tensor PCA} and \eqref{eq: estimate of sA2}, respectively. Then
\begin{align*}
& \sup _{x \in \mathbb{R}} \left\lvert\, \mathbb{P}\left(\frac{\langle \widehat{\mcT}, \mcA \rangle-\langle \cT, \mcA\rangle}{\sigma \cdot s_{\mcA}} \leq x\right)-\Phi(x)\right| \lesssim \underbrace{\Psi}_{\text{rate of asymptotic normal terms}} + \underbrace{\frac{\widetilde{\Omega}}{\sigma \cdot s_{\mcA}}}_{\text{rate of negligible terms}} + \underbrace{\big[ \op^{-c} + \mcP_{\mcE_{U}^{\text{PCA}}} \big]}_{\text{rate of initial estimates}} ,
\end{align*}
where $c>0$ is a constant, and the variance component $s_{\mcA}$ is defined in \eqref{eq: variance component sA2}.
Here, $\Psi$ is the same as in Theorem~\ref{thm: main theorem in tensor PCA}, and
\begin{align*}
\widetilde{\Omega}
= & \sum_{j=1}^3 \left\|\mcP_{U_j} A_j \mcP_{\left(U_{j+2} \otimes U_{j+1}\right) G_j^{\top}}\right\|_{\mathrm{F}} \cdot \frac{\sigma^2\oor^{1/2}\cdot \op}{\ulambda} + \sum_{j=1}^{3} \left\|\mcA \times_{j+1} U_{j+1} \times_{j+2} U_{j+2}\right\|_{\mathrm{F}} \cdot \frac{\sigma^2\oor^{1/2} \cdot \sqrt{\op \log (\op)}}{\ulambda} \\
+ & \sum_{j=1}^3 \left\|\mcA \times_j U_j \right\|_{\mathrm{F}} \cdot \left(\frac{\sigma \oor^{1/2} \cdot \sqrt{\oR}\log (\op)}{\ulambda}+\frac{\sigma^3 \oor^{1/2} \cdot \oR \sqrt{\op \log (\op)}}{\ulambda^2}+\frac{\sigma^4 \oor^{1/2} \cdot \op^{3 / 2} \sqrt{\oR \log (\op)}}{\ulambda^3}\right)  \\
+ & \left\|\mcA\right\|_{\mathrm{F}} 
\cdot \left(\frac{\sigma^2 \oor^{1/2} \cdot \oR\log (\op)^{3 / 2}}{\ulambda^2} 
+ \frac{\sigma^4 \oor^{1/2} \cdot \oR^{3 / 2} \op^{1 / 2} \log (\op)}{\ulambda^3} 
+ \frac{\sigma^5 \oor^{1/2} \cdot \oR \op^{3 / 2} \log (\op)}{\ulambda^4}\right) .
\end{align*}

\end{theorem}

\subsection{Proof of Asymptotic Normality with Plug-in Estimates in Tensor Regression (Theorem~\ref{thm: main theorem in tensor regression with plug-in estimates})} \label{subsec: Proof of Asymptotic Normality with Plug-in Estimates in Tensor Regression}

\subsubsection*{Part 1: Without Sample Splitting}

Note that
\begin{align*}
\frac{\langle \mcA, \widehat{\mcT} \rangle - \left\langle \mcA, \mcT \right\rangle}{(\widehat{\sigma}_\xi/\widehat{\sigma})  \widehat{s}_\mcA \cdot \sqrt{\frac{1}{n}}} 
&= \frac{\langle \mcA, \widehat{\mcT} \rangle - \left\langle \mcA, \mcT \right\rangle}{(\sigma_\xi/\sigma) s_\mcA \cdot \sqrt{\frac{1}{n}}} + \left(\frac{(\sigma_\xi/\sigma) s_\mcA}{(\widehat{\sigma}_\xi/\widehat{\sigma}) \widehat{s}_\mcA} - 1\right) 
\cdot 
\frac{\langle \mcA, \widehat{\mcT} \rangle - \left\langle \mcA, \mcT \right\rangle}{(\sigma_\xi/\sigma) s_\mcA \cdot \sqrt{\frac{1}{n}}}.
\end{align*}

\subsubsection*{Step 1: Upper Bound of $|\widehat{\sigma}_\xi^2 - \sigma_\xi^2 |$}

Note that
\begin{align*}
\widehat{\sigma}_\xi^2 = & \frac{1}{n} \sum_{i=1}^{n}\left(y_i - \left\langle \widehat{\mathcal{T}}^{\text{init}}, \mcX_i\right\rangle\right)^2 = \frac{1}{n} \sum_{i=1}^{n} \xi_i^2 + \frac{2}{n} \sum_{i=1}^{n} \xi_i\left\langle \widehat{\Delta}, \mcX_i\right\rangle + \frac{1}{n}\sum_{i=1}^{n}\left(\left\langle \widehat{\Delta}, \mcX_i\right\rangle\right)^2,
\end{align*}
where $\widehat{\Delta} = \widehat{\mcT}^{\text{init}} - \mcT$. As $\mathbb{E}\left(\frac{1}{n} \sum_{i=1}^n \xi_i^2\right) =\sigma_\xi^2$, it holds that
$$
\left|\frac{1}{n}\sum_{i=1}^n \xi_i^2 - \sigma_\xi^2\right| \leq \sigma_\xi^2\sqrt{\log(\op)}
$$
with probability at least $1-\op^{-c}$, where $C$ is a constant.

Then, it remains to find a high-probability upper bound of
$\frac{1}{n} \sum_{i=1}^{n}\langle\whDelta, \mcX_{i}\rangle^2$. By Hanson-Wright inequality and $\varepsilon$-net argument on low-Tucker-rank manifold, it follows that
\begin{align*}
\frac{1}{n}\sum_{i=1}^n \left\langle \whDelta, \mcX_i \right\rangle^2
= & \frac{1}{n}\sum_{i=1}^n \operatorname{Vec}\left(\mcX_i\right)^{\top}\operatorname{Vec}\left(\whDelta\right)\operatorname{Vec}\left(\whDelta\right)^{\top}\operatorname{Vec}\left(\mcX_i\right) \lesssim \Delta ^2 \cdot \left(\sigma_\xi^2\sigma^2 \cdot \sqrt{\frac{\op\oor}{n}}\right) ,
\end{align*}
which holds with probability at least $1-\mcP\left(\mcE_{\Delta}\right)- \exp(-cn)- \exp(-c\op)$.

In addition, we have
\begin{align*}
\frac{1}{n}\sum_{i=1}^{n}\xi_i\left\langle \whDelta, \mcX_i\right\rangle 
\leq & \|\whDelta \|_{\mathrm{F}} \cdot \sup_{\substack{W_j \in \mathbb{R}^{p_j \times 2r_j}, j=1,2,3}}\Big\|W_1^{\top}\Big(\frac{1}{n}\sum_{i=1}^{n}\xi_i\mcX_i\Big)\left(W_3\otimes W_2\right)\Big\|_{\mathrm{F}} \lesssim \sigma_\xi\sigma \cdot \Delta \cdot \sqrt{\frac{\op\oor}{n}}
\end{align*}
holds with probability at least $1-\mcP\left(\mcE_\Delta\right)- \exp(-cn)- \exp(-c\op)$.

Therefore, we have
$$
\left|\widehat{\sigma}_\xi^2 -\sigma_\xi^2\right| \lesssim 2\cdot \sigma_\xi\sigma \cdot \Delta \cdot \sqrt{\frac{\op\oor}{n}} + \sigma_\xi^2 \cdot \Delta^2 \cdot \sqrt{\frac{\op\oor}{n}} \lesssim \sigma_\xi\sigma \cdot \Delta \cdot \sqrt{\frac{\op\oor}{n}}.
$$

Since $\Delta = o(1)$ and $n \gtrsim \op$, it implies that $\widehat{\sigma}_\xi^2 -\sigma_\xi^2 = o(1)$. Therefore, we have
$$
\left|\frac{\sigma_\xi}{\widehat{\sigma}_\xi} -1 \right| \lesssim \left|\frac{\sigma_\xi^2}{\widehat{\sigma}_\xi^2} -1 \right| \lesssim \frac{\sigma}{\sigma_\xi} \cdot \Delta \cdot \sqrt{\frac{\op\oor}{n}},
$$
which holds with probability at least $1-\mcP\left(\mcE_\Delta\right)- \exp(-cn)- \exp(-c\op)$.

\subsubsection*{Step 2: Upper Bound of $|\widehat{\sigma}^2 -  \sigma^2 |$}

Then consider 
\begin{align*}
\widehat{\sigma}^2 
= & \frac{1}{np_1p_2p_3}\sum_{i=1}^n\left\|\mcX_i\right\|_{\mathrm{F}}^2 
= \frac{1}{np_1p_2p_3}\sum_{i=1}^n \sum_{j_1=1}^{p_1} \sum_{j_2=1}^{p_2} \sum_{j_3=1}^{p_3}\left[\mcX_i\right]_{j_1,j_2,j_3}^2.
\end{align*}
Note that $\mathbb{E}\left(\widehat\sigma^2\right) = \sigma^2$. Applying the Hanson-Wright inequality for the quadratic form once again, $\left|\widehat{\sigma}^2 -  \sigma^2\right|\leq  \sigma^2\sqrt{\log(\op)/(np_1p_3p_3)}$
with probability at least $1-\exp(-cn) - \op^{-c}$.

\subsubsection*{Step 3: Upper Bound of $|\widehat{s}_{\mcA}^2 - s_{\mcA}^2 |$}

Here, consider
\begin{align*}
& \left| \left\|\mcP_{\whU_{j \perp}}A_j\mcP_{\widehat{\left(U_{j+2} \otimes U_{j+1}\right) G_j^{\top}}}\right\|_{\mathrm{F}}^2 - \left\|\mcP_{U_{j \perp}} A_j \mcP_{\left(U_{j+2} \otimes U_{j+1}\right) G_j^{\top}}\right\|_{\mathrm{F}}^2\right| \\
\leq & \underbrace{\left\|\left(\mcP_{\whU_{j \perp}}-\mcP_{U_{j \perp}}\right) A_j \mcP_{\left(U_{j+2} \otimes U_{j+1}\right) G_j^{\top}}\right\|^2_{\mathrm{F}}}_{\mathrm{\RN{1}}^2} \\
+ & \underbrace{\left\|\mcP_{U_{j \perp}} A_j\left[\left(\whU_3\otimes \whU_2\right)\widehat{W}_1\widehat{W}_1^{\top}\left(\whU_3\otimes \whU_2\right)^{\top} - \left(U_3 \otimes U_2\right)W_1W_1^{\top}\left(U_3 \otimes U_2\right)^{\top}\right]\right\|^2_{\mathrm{F}}}_{\mathrm{\RN{2}}^2} \\
+ & \underbrace{\left\|\left(\mcP_{\whU_{j \perp}}-\mcP_{U_{j \perp}}\right) A_j\left[\left(\whU_3\otimes \whU_2\right)\widehat{W}_1\widehat{W}_1^{\top}\left(\whU_3\otimes \whU_2\right)^{\top} - \left(U_3 \otimes U_2\right)W_1W_1^{\top}\left(U_3 \otimes U_2\right)^{\top}\right]\right\|^2_{\mathrm{F}}}_{\mathrm{\RN{3}}^2},
\end{align*}
where $\mcP_{\widehat{\left(U_{j+2} \otimes U_{j+1}\right) G_j^{\top}}} = \left(\whU_{j+2} \otimes \whU_{j+1}\right)\widehat{W}_j\widehat{W}_j^{\top}\left(\whU_{j+2} \otimes \whU_{j+1}\right)^{\top}$.

Here, we have 
\begin{align*}
\mathrm{\RN{1}} \leq & \left\|\left(\widehat{\mcP}_{U_{1 \perp}}-\mcP_{U_{1\perp}}\right) \mcP_{U_1} A_1 \mcP_{\left(U_3 \otimes U_2\right) G_1^{\top}}\right\|_{\mathrm{F}} + \left\|\left(\widehat{\mcP}_{U_{1 \perp}}-\mcP_{U_{1\perp}}\right) \mcP_{U_{1\perp}} A_1 \mcP_{\left(U_3 \otimes U_2\right) G_1^{\top}}\right\|_{\mathrm{F}} \\
\lesssim & \frac{\sigma_\xi}{\ulambda\sigma}\sqrt{\frac{\op}{n}}\cdot \left\|\mcP_{U_1}A_1 \mcP_{\left(U_3 \otimes U_2\right) G_1^{\top}}\right\|_{\mathrm{F}} + \underbrace{\left[\frac{\sigma_{\xi}}{\ulambda \sigma} \cdot\left(\sqrt{\frac{\oor\log(\op)}{n}}+\Delta \sqrt{\frac{\op}{n}}\right)\right] \cdot \left\|\mcP_{U_{1\perp}}A_1 \mcP_{\left(U_3 \otimes U_2\right) G_1^{\top}}\right\|_{\mathrm{F}}}_{\eqref{eq: high-prob upper bound of V1tPU1p(PUhat1-PU1)U1 in tensor regression without sample splitting},\eqref{eq: high-prob upper bound of V1tPU1p(PUhat1-PU1)U1p in tensor regression without sample splitting}}.
\end{align*}

On the other hand, we have
\begin{align*}
\mathrm{\RN{2}}= & \underbrace{\left\|\mcP_{U_{1\perp}}A_1\left(\mcP_{U_3}\otimes \mcP_{U_2}\right)\left[\left(\whU_3\otimes \whU_2\right)\widehat{W}_1\widehat{W}_1^{\top}\left(\whU_3\otimes \whU_2\right)^{\top} - \left(U_3 \otimes U_2\right)W_1W_1^{\top}\left(U_3 \otimes U_2\right)^{\top}\right]\right\|_{\mathrm{F}}^2}_{\mathrm{\RN{2}}.\mathrm{\RN{1}}} \\
+ & \underbrace{\left\|\mcP_{U_{1\perp}}A_1\left(\mcP_{U_{3\perp}}\otimes \mcP_{U_2}\right)\left[\left(\whU_3\otimes \whU_2\right)\widehat{W}_1\widehat{W}_1^{\top}\left(\whU_3\otimes \whU_2\right)^{\top} - \left(U_3 \otimes U_2\right)W_1W_1^{\top}\left(U_3 \otimes U_2\right)^{\top}\right]\right\|_{\mathrm{F}}^2}_{\mathrm{\RN{2}}.\mathrm{\RN{2}}} \\
+ & \underbrace{\left\|\mcP_{U_{1\perp}}A_1\left(\mcP_{U_3}\otimes \mcP_{U_{2\perp}}\right)\left[\left(\whU_3\otimes \whU_2\right)\widehat{W}_1\widehat{W}_1^{\top}\left(\whU_3\otimes \whU_2\right)^{\top} - \left(U_3 \otimes U_2\right)W_1W_1^{\top}\left(U_3 \otimes U_2\right)^{\top}\right]\right\|_{\mathrm{F}}^2}_{\mathrm{\RN{2}}.\mathrm{\RN{3}}} \\
+ & \underbrace{\left\|\mcP_{U_{1\perp}}A_1\left(\mcP_{U_{3\perp}}\otimes \mcP_{U_{2\perp}}\right)\left[\left(\whU_3\otimes \whU_2\right)\widehat{W}_1\widehat{W}_1^{\top}\left(\whU_3\otimes \whU_2\right)^{\top} - \left(U_3 \otimes U_2\right)W_1W_1^{\top}\left(U_3 \otimes U_2\right)^{\top}\right]\right\|_{\mathrm{F}}^2}_{\mathrm{\RN{2}}.\mathrm{\RN{4}}} .
\end{align*}

Here, first, by Lemma 7 in \citet{zhang2020islet}, we have
\begin{align*}
& \left\| (\whU_3\otimes \whU_2 )\widehat{W}_1\widehat{W}_1^{\top} (\whU_3\otimes \whU_2 )^{\top} -  (U_3 \otimes U_2 )W_1W_1^{\top} (U_3 \otimes U_2 )^{\top} \right\| 
\leq  \frac{2\left\|\mcP_{\whU_1}\widehat{T}_1(\mcP_{\whU_3} \otimes \mcP_{\whU_2}) - T_1\right\|}{\ulambda} \lesssim \frac{\sigma_\xi}{\ulambda\sigma}\sqrt{\frac{\op}{n}}.
\end{align*}
where the second inequality follows from that
\begin{align*}
\left\|\mcP_{\whU_1}\widehat{T}_1 (\mcP_{\whU_3} \otimes \mcP_{\whU_2} ) - T_1\right\| 
\leq & \left\|\big[\mcP_{U_1} + (\mcP_{\whU_1} - \mcP_{U_1})\big] \big[T_1 + (\widehat{T}_1 - T_1 )\big]\big[(\mcP_{U_3} \otimes \mcP_{U_2}) +((\mcP_{\whU_3} \otimes \mcP_{\whU_2}) - (\mcP_{U_3} \otimes \mcP_{U_2}))\big]\right\| \\
\lesssim & \frac{\sigma_\xi}{\ulambda\sigma}\sqrt{\frac{\op}{n}} \cdot \ulambda = \frac{\sigma_\xi}{\sigma}\sqrt{\frac{\op}{n}}.
\end{align*}

Therefore, we have
\begin{align*}
\mathrm{\RN{2}}.\mathrm{\RN{1}}
\lesssim& \frac{\sigma_\xi^2}{\ulambda^2\sigma^2} \cdot \frac{\op}{n} \cdot \left\|\mcP_{U_{1\perp}}A_1\left(\mcP_{U_3}\otimes \mcP_{U_2}\right)\right\|_{\mathrm{F}}^2.
\end{align*}
Similarly, we have
\begin{align*}
\mathrm{\RN{2}}.\mathrm{\RN{2}}
\lesssim & \left\|\mcA \times_2 U_2\right\|_{\mathrm{F}}^2 \cdot \left[\frac{\sigma_{\xi}}{\ulambda \sigma} \cdot\left(\sqrt{\frac{\oR \log (\op)}{n}}+\Delta \sqrt{\frac{\op}{n}}\right)\right]^2,\\
\mathrm{\RN{2}}.\mathrm{\RN{3}} \lesssim & \left\|\mcA \times_3 U_3\right\|_{\mathrm{F}}^2 \cdot \left[\frac{\sigma_{\xi}}{\ulambda \sigma} \cdot\left(\sqrt{\frac{\oR \log (\op)}{n}}+\Delta \sqrt{\frac{\op}{n}}\right)\right]^2, \\
\mathrm{\RN{2}}.\mathrm{\RN{4}}
\lesssim & \left\|\mcA\right\|_{\mathrm{F}}^2 \cdot \left[\frac{\sigma_{\xi}}{\ulambda \sigma} \cdot\left(\sqrt{\frac{\oR \log (\op)}{n}}+\Delta \sqrt{\frac{\op}{n}}\right)\right]^4.
\end{align*}

Combining the results above, we have
\begin{align*}
& \left| \left\|\mcP_{\whU_{j \perp}}A_j\mcP_{\widehat{\left(U_{j+2} \otimes U_{j+1}\right) G_j^{\top}}}\right\|_{\mathrm{F}}^2 - \left\|\mcP_{U_{j \perp}} A_j \mcP_{\left(U_{j+2} \otimes U_{j+1}\right) G_j^{\top}}\right\|_{\mathrm{F}}^2\right| \\
\lesssim & \sum_{j=1}^{3} \left\|\mcA \times_{j+1} U_{j+1} \times_{j+2} U_{j+2} \right\|_{\mathrm{F}}^2 \cdot \frac{\sigma_{\xi}^2}{\ulambda^2\sigma^2} \cdot \frac{\op}{n} 
+  \sum_{j=1}^{3} \left\|\mcA \times_j U_j\right\|_{\mathrm{F}}^2 \cdot \left[\frac{\sigma_\xi^2}{\ulambda^2\sigma^2}\cdot \left(\frac{\oR\log(\op)}{n} + \Delta^2\cdot \frac{\op}{n}\right)\right]  \\
& +  \left\|\mcA\right\|_{\mathrm{F}}^2 \cdot \left[\frac{\sigma_\xi^4}{\ulambda^4\sigma^4}\cdot \left(\frac{\oR^2\log(\op)^2}{n^2} + \Delta^4\cdot \frac{\op^2}{n^2}\right)\right]
\end{align*}

Furthermore, we have
\begin{align}
& \left|\left\|\mcA \times_1 \mcP_{\whU_1} \times_2 \mcP_{\whU_2} \times_3 \mcP_{\whU_3}\right\|_{\mathrm{F}}^2 - \left\|\mcA \times_1 \mcP_{U_1} \times_2 \mcP_{U_2} \times_3 \mcP_{U_3}\right\|_{\mathrm{F}}^2 \right| \notag \\
\leq & \sum_{j=1}^{3} \left\|\mcA \times_{j+1} U_{j+1} \times_{j+2} U_{j+2} \right\|_{\mathrm{F}}^2 \cdot \left\|V_j^{\top}\left(\mcP_{\whU_j} - \mcP_{U_j}\right)\right\|^2 \notag \\
+ & \sum_{j=1}^{3} \left\|\mcA \times_j U_j\right\|_{\mathrm{F}}^2 \cdot \left\|V_{j+1}^{\top}\left(\mcP_{\whU_{j+1}} - \mcP_{U_{j+1}}\right)\right\|^2 \cdot \left\|V_{j+2}^{\top}\left(\mcP_{\whU_{j+2}} - \mcP_{U_{j+2}}\right)\right\|^2 \notag \\
+ & \sum_{j=1}^{3} \left\|\mcA\right\|_{\mathrm{F}}^2 \cdot \left\|V_j^{\top}\left(\mcP_{\whU_j} - \mcP_{U_j}\right)\right\|^2 \cdot \left\|V_{j+1}^{\top}\left(\mcP_{\whU_{j+1}} - \mcP_{U_{j+1}}\right)\right\|^2 \cdot \left\|V_{j+2}^{\top}\left(\mcP_{\whU_{j+2}} - \mcP_{U_{j+2}}\right)\right\|^2 \notag \\
\lesssim & \sum_{j=1}^{3} \left\|\mcA \times_{j+1} U_{j+1} \times_{j+2} U_{j+2} \right\|_{\mathrm{F}}^2 \cdot \left[\frac{\sigma_\xi^2}{\ulambda^2\sigma^2}\cdot \left(\frac{\oR\log(\op)}{n} + \Delta^2\cdot \frac{\op}{n}\right)\right] \notag \\
+ & \sum_{j=1}^{3} \left\|\mcA \times_j U_j\right\|_{\mathrm{F}}^2 \cdot \left[\frac{\sigma_\xi^4}{\ulambda^4\sigma^4}\cdot \left(\frac{\oR^2\log(\op)^2}{n^2} + \Delta^4\cdot \frac{\op}{n}\right)\right] + \left\|\mcA\right\|_{\mathrm{F}}^2 \cdot \left[\frac{\sigma_\xi^6}{\ulambda^6\sigma^6}\cdot \left(\frac{\oR^3\log(\op)^3}{n^3} + \Delta^6\cdot \frac{\op^3}{n^3}\right)\right]. \notag 
\end{align}

Combining all the results above, it follows that
\begin{align}
\widehat{s}_{\mcA}^2 - s_{\mcA}^2
\lesssim & \sum_{j=1}^{3} \left\|\mcA \times_{j+1} U_{j+1} \times_{j+2} U_{j+2} \right\|_{\mathrm{F}}^2 \cdot \frac{\sigma_{\xi}^2}{\ulambda^2\sigma^2} \cdot \frac{\op}{n} + \sum_{j=1}^{3} \left\|\mcA \times_j U_j\right\|_{\mathrm{F}}^2 \cdot \left[\frac{\sigma_\xi^2}{\ulambda^2\sigma^2}\cdot \left(\frac{\oR\log(\op)}{n} + \Delta^2\cdot \frac{\op}{n}\right)\right] \notag \\
& + \left\|\mcA\right\|_{\mathrm{F}}^2 \cdot \left[\frac{\sigma_\xi^4}{\ulambda^4\sigma^4}\cdot \left(\frac{\oR^2\log(\op)^2}{n^2} + \Delta^4\cdot \frac{\op^2}{n^2}\right)\right]. 
\end{align}

\subsubsection*{Step 4: Summary of results}

On the other hand, we have
\begin{align*}
& \left|\langle \mcA, \widehat{\mcT} \rangle -\left\langle \mcA, \mcT\right\rangle\right| \lesssim  \left|\left\langle  \widehat{\mcZ}^{(1)}, \mcP_{\mathbb{T}_{\mcT} \mathcal{M}_{(r_1, r_2, r_3)}}\left(\mcA\right)\right\rangle\right| + \left(\Omega_1 + \Omega_2 + \Omega_3 + \Omega_4\right) \\
\lesssim & s_\mcA \frac{\sigma_\xi}{\sigma}\cdot \sqrt{\frac{\log(\op)}{n}} + \left(\Omega_1 + \Omega_2 + \Omega_3 + \Omega_4\right),
\end{align*}
where $\Omega_j, j=1,2,3,4$ is defined in the statement of Theorem~\ref{thm: main theorem in tensor regression without sample splitting}.
It further implies that
\begin{align*}
&\left| \left(\frac{(\sigma_\xi/\sigma) s_{\mcA}}{ (\widehat{\sigma}_\xi/\widehat{\sigma}) \widehat{s}_{\mcA}} -1\right)  \cdot \frac{\sqrt{n} \big(\langle \mcA, \widehat{\mcT} \rangle - \langle \mcA, \mcT \rangle \big)}{(\sigma_{\xi}/\sigma) s_{\mcA} } \right| \\
\lesssim & \left(\frac{\sigma}{\sigma_\xi}\cdot \Delta \sqrt{\op\oor\log(\op)}\right) + \sqrt{\frac{\log(\op)^2}{np_1p_2p_3}}   
+  \frac{1}{(\sigma_{\xi}/\sigma) s_{\mcA} \sqrt{\frac{1}{n}}} \Bigg\{\sum_{j=1}^{3} \left\|\mcA \times_{j+1} U_{j+1} \times_{j+2} U_{j+2} \right\|_{\mathrm{F}} \cdot \frac{\sigma_{\xi}^2}{\ulambda\sigma^2} \cdot \frac{\sqrt{\op\log(\op)}}{n} \notag \\
+ & \sum_{j=1}^{3} \left\|\mcA \times_j U_j\right\|_{\mathrm{F}}  \left[\frac{\sigma_\xi^2}{\ulambda\sigma^2} \left(\frac{\sqrt{\oR}\log(\op)}{n} + \Delta  \frac{\sqrt{\op\log(\op)}}{n}\right)\right] 
+  \left\|\mcA\right\|_{\mathrm{F}}  \left[\frac{\sigma_\xi^3}{\ulambda^2\sigma^3}\cdot \left(\frac{\oR\log(\op)^{3/2}}{n^{3/2}} + \Delta^2\cdot \frac{\sqrt{\op\log(\op)}}{n}\right)\right] \Bigg\} .
\end{align*}

Compared with the terms $\Omega_1, \Omega_2, \Omega_3$, and $\Omega_4$ in Theorem~\ref{thm: main theorem in tensor regression without sample splitting}, we have the desired results.

\subsubsection*{Part 2: With Sample Splitting}

In the sample-splitting case, we consider
\begin{align*}
\widehat{\sigma}_{\xi}^2
= & \frac{1}{n} \sum_{i=1}^n \xi_i^2+ \frac{1}{n} \sum_{i_1=1}^{n_1} \left\langle\whDelta^{(\mathrm{\RN{2})}}, \mcX_{i_1}\right\rangle^2+\frac{1}{n} \sum_{i_2=1}^{n_2}\left\langle\whDelta^{(\mathrm{\RN{1}})}, \mcX_{i_2}\right\rangle^2 + \frac{2}{n} \sum_{i_1=1}^{n_1} \xi_i\left\langle\whDelta^{(\mathrm{\RN{2}}}), \mcX_{i_1}\right\rangle + \frac{2}{n} \sum_{i_2=1}^{n_2} \xi_i\left\langle\whDelta^{(\mathrm{\RN{1}})}, \mcX_{i_2}\right\rangle  .
\end{align*}

\subsubsection*{Step 1: Upper Bound of $|\widehat{\sigma}_\xi^2 - \sigma_\xi^2 |$}

By similar arguments, it follows that $\mathbb{E}\left(n^{-1} \sum_{i=1}^n \xi_i^2\right)=\sigma_{\xi}^2$. Furthermore, due to sample splitting, note that $\whDelta^{(\mathrm{II})}$ and $\mcX_{i_1}^{(\mathrm{I})}$ are independent. By Hansen-Wright inequality, it follows that
$$
\mathbb{P}\left(\left|\frac{1}{n} \sum_{i_1=1}^{n_1}\left\langle\whDelta^{(\mathrm{\RN{2}})}, \mcX_{i_1}\right\rangle^2\right|\geq \sigma^2\cdot \|\whDelta^{(\mathrm{\RN{2}})}\|_{\mathrm{F}}^2+ Ct \right) \leq \exp\left[-c\min\left(\frac{n t^2}{\sigma^2\|\whDelta^{\text {(\RN{2}) }}\|_{\mathrm{F}}^4}, \frac{n t}{\sigma\|\whDelta^{\text {(\RN{2}) }}\|^2}\right)\right]
$$

Therefore, it follows that
$$
\mathbb{P}\left(\left|\frac{1}{n} \sum_{i=1}^{n_1}\left\langle\whDelta^{(\mathrm{II})}, \mcX_{i_1}\right\rangle^2\right|\geq  \sigma^2 \cdot \Delta^2\sqrt{\frac{\log(\op)}{n}} \right) \leq \op^{-c}+ \mathbb{P}\left(\mcE_{\Delta}\right).
$$

In addition, with probability at least $1- \op^{-c} -\mcP_{\mcE_{\Delta}}$, it holds that
$$
\frac{1}{n_1}\sum_{i_1=1}^{n_1}\xi_i\left\langle \widehat{\Delta}^{(\mathrm{\RN{2}})}, X_{i_1}\right\rangle \leq \Delta\cdot \sigma_\xi\sigma\sqrt{\frac{\oor\log(\op)}{n}}.
$$

Combining the results above, with probability at least $1- \op^{-c} -\mcP_{\mcE_{\Delta}}$, it holds that
$
\widehat{\sigma}_\xi^2 - \sigma_\xi^2 \leq \Delta  \sigma_\xi\sigma\sqrt{ \oor\log(\op)/n }.
$

\subsubsection*{Step 2: Upper Bound of $|\widehat{\sigma}^2 -  \sigma^2 |$}
Similar to part 1, $\left|\widehat{\sigma}^2 -  \sigma^2\right| \leq \sqrt{\log(\op)/(np_1p_3p_3)}$ with probability at least $1-e^{-cn} - \op^{-c}$.

\subsubsection*{Step 3: Upper Bound of $\left|\widehat{s}_{\mcA}^2 - s_{\mcA}^2\right|$}

By the same argument in Part 1, we have
\begin{align*}
\widehat{s}_{\mcA}^2 - s_{\mcA}^2 
\lesssim & \sum_{j=1}^{3} \left\|\mcA \times_{j+1} U_{j+1} \times_{j+2} U_{j+2} \right\|_{\mathrm{F}}^2 \cdot \frac{\sigma_{\xi}^2}{\ulambda^2\sigma^2} \cdot \frac{\op}{n} + \sum_{j=1}^{3} \left\|\mcA \times_j U_j\right\|_{\mathrm{F}}^2 \cdot \left[\frac{\sigma_\xi^2}{\ulambda^2\sigma^2}\cdot \left(\frac{\oR\log(\op)}{n} \right)\right]\\
& +  \left\|\mcA\right\|_{\mathrm{F}}^2 \cdot \left[\frac{\sigma_\xi^4}{\ulambda^4\sigma^4}\cdot \left(\frac{\oR^2\log(\op)^2}{n^2}\right)\right].
\end{align*}

\subsubsection*{Step 4: Summary of results}

On the other hand, similar to part 1, we have
\begin{align*}
& \left| \langle \mcA, \widehat{\mcT} \rangle - \langle \mcA, \mcT \rangle\right| 
\lesssim  s_\mcA \cdot \frac{\sigma_\xi}{\sigma} \cdot \sqrt{\frac{\log(\op)}{n}}  + \left(\Omega_1 + \Omega_2 + \Omega_3 \right),
\end{align*}
where $\Omega_j, j=1,2,3$ are defined in Theorem~\ref{thm: main theorem in tensor regression with sample splitting}.
It implies the desired results in part 2.

\subsection{Proof of Asymptotic Normality with Plug-in Estimates in Tensor PCA (Theorem~\ref{thm: main theorem in tensor PCA with plug-in estimates})}\label{subsec: Proof of Asymptotic Normality with Plug-in Estimates in Tensor PCA}

Note that
\begin{align*}
& \frac{ \langle \mcA, \widehat{\mcT}  \rangle-  \langle\mcA, \mcT \rangle}{\sigma  \widehat{s}_\mcA } = \frac{ \langle \mcA, \widehat{\mcT}  \rangle - \langle \mcA, \mcT  \rangle}{\sigma  s_{\mcA} } + \frac{ \langle \mcA, \widehat{\mcT} \rangle - \langle \mcA, \mcT  \rangle}{\sigma  s_{\mcA} } \cdot\left(\frac{\widehat{\sigma} \widehat{s}_{\mcA}}{\sigma  s_{\mcA}} - 1\right) ,
\end{align*}
where $s_\mcA^2=\sum_{j=1}^3\left\|\mcP_{U_j} A_j \mcP_{\left(U_{j+2} \otimes U_{j+1}\right) G_j}\right\|_{\mathrm{F}}^2+\left\|\mcA \times_1 U_1 \times_2 U_2 \times_3 U_3\right\|_{\mathrm{F}}^2$.

\subsubsection*{Step 1: Upper Bound of $|\widehat{\sigma}^2 - \sigma^2 |$ }

First, since $\mcY = \mcT + \mcZ$, it then follows that
\begin{align*}
& \big|\|\mcY-\mcY \times_1 \mcP_{\whU_1} \times_2 \mcP_{\whU_2} \times_3 \mcP_{\whU_3}\|_{\mathrm{F}}-\|\mcZ\|_{\mathrm{F}}\big| \\
\leq & \big\|\mcZ \times_1 \mcP_{\whU_1} \times_2 \mcP_{\whU_2} \times_3 \mcP_{\whU_3}\big\|_{\mathrm{F}} + \big\|\mcT-\mcT \times_1 \mcP_{\whU_1} \times_2 \mcP_{\whU_2} \times_3 \mcP_{\whU_3}\big\|_{\mathrm{F}}.
\end{align*}

First, we have
$
|\mcZ \times_1 \mcP_{\whU_1} \times_2 \mcP_{\whU_2} \times_3 \mcP_{\whU_3}\|_{\mathrm{F}}=\|\whU_1^{\top} Z_1\left(\whU_3 \otimes \whU_2\right)\|_{\mathrm{F}} \leq \sigma \sqrt{\op \oor}.
$
In addition, 
\begin{align*}
& \left\|\mcT-\mcT \times_1 \mcP_{\whU_1} \times_2 \mcP_{\whU_2} \times_3 \mcP_{\whU_3}\right\|_{\mathrm{F}} 
\leq  \sum_{j=1}^3 \left\|\left(\mcP_{U_j}-\mcP_{\whU_j}\right)U_j G_j \left(U_{j+2}\otimes U_{j+1}\right)^{\top}\right\| 
\lesssim  \sigma \sqrt{\op \oor},
\end{align*}
with probability at least $1-e^{-c\op}$. Therefore, with probability at least $1-e^{-c \op}$, we have
$$
\left|\left\|\widehat{\mcT}- \widehat{\mcT} \times_1 \mcP_{\whU_1} \times_2 \mcP_{\whU_2} \times_3 \mcP_{\whU_3}\right\|_{\mathrm{F}}-\|\mcZ\|_{\mathrm{F}}\right| \leq \sigma \sqrt{\op \oor}.
$$
In addition, we have
$$
\mathbb{P}\left(\left| \|\mcZ\|_{\mathrm{F}}^2/\sigma^2-p_1 p_2 p_3\right| \geq C_2\left(\sqrt{p_1p_2p_3\log (\op)}+\log (\op)\right)\right) \leq \left(p_1p_2p_3\right)^{-1} .
$$
As a consequence, with probability at least $1-\op^{-3}$,
$$
\left|\|\mcZ\|_{\mathrm{F}}-\sqrt{p_1 p_2 p_3} \sigma\right| \leq \sigma\sqrt{\log (\op)}.
$$
Therefore, we know that with probability at least $1- \op^{-c}$,
\begin{align*}
\left|\widehat{\sigma}^2 / \sigma^2-1\right|
= & \left|\widehat{\sigma} / \sigma-1\right| \cdot \left| \widehat{\sigma} / \sigma+1\right| \leq 2\left| \widehat{\sigma}/\sigma-1\right|+ \left| \widehat{\sigma}/\sigma -1\right|^2 \lesssim \sqrt{\oor\op/(p_1p_2p_3)} .
\end{align*}

\subsubsection*{Step 2: Upper Bound of $\left|\widehat{s}_{\mcA}^2 -s_{\mcA}^2\right|$}

Then it remains to consider
\begin{align*}
& \left| \left\|\mcP_{\whU_{j \perp}}A_j\mcP_{\widehat{\left(U_{j+2} \otimes U_{j+1}\right) G_j^{\top}}}\right\|_{\mathrm{F}}^2 - \left\|\mcP_{U_{j \perp}} A_j \mcP_{\left(U_{j+2} \otimes U_{j+1}\right) G_j^{\top}}\right\|_{\mathrm{F}}^2\right| \\
\leq & \underbrace{\left\|\left(\mcP_{\whU_{j \perp}}-\mcP_{U_{j \perp}}\right) A_j \mcP_{\left(U_{j+2} \otimes U_{j+1}\right) G_j^{\top}}\right\|^2_{\mathrm{F}}}_{\mathrm{\RN{1}}^2} \\
+ & \underbrace{\left\|\mcP_{U_{j \perp}} A_j\left[\left(\whU_3\otimes \whU_2\right)\widehat{W}_1\widehat{W}_1^{\top}\left(\whU_3\otimes \whU_2\right)^{\top} - \left(U_3 \otimes U_2\right)W_1W_1^{\top}\left(U_3 \otimes U_2\right)^{\top}\right]\right\|^2_{\mathrm{F}}}_{\mathrm{\RN{2}}^2} \\
+ & \underbrace{\left\|\left(\mcP_{\whU_{j \perp}}-\mcP_{U_{j \perp}}\right) A_j\left[\left(\whU_3\otimes \whU_2\right)\widehat{W}_1\widehat{W}_1^{\top}\left(\whU_3\otimes \whU_2\right)^{\top} - \left(U_3 \otimes U_2\right)W_1W_1^{\top}\left(U_3 \otimes U_2\right)^{\top}\right]\right\|^2_{\mathrm{F}}}_{\mathrm{\RN{3}}^2}.
\end{align*}

Similar to the proof of step 2 in Section \ref{subsec: Proof of Asymptotic Normality with Plug-in Estimates in Tensor Regression}, we can show
\begin{align*}
& \sum_{j=1}^3 \left| \left\|\mcP_{\whU_{j \perp}}A_j\mcP_{\widehat{\left(U_{j+2} \otimes U_{j+1}\right) G_j^{\top}}}\right\|_{\mathrm{F}}^2 - \left\|\mcP_{U_{j \perp}} A_j \mcP_{\left(U_{j+2} \otimes U_{j+1}\right) G_j^{\top}}\right\|_{\mathrm{F}}^2\right| \\
\lesssim & \sum_{j=1}^3 \left\|\mcA \times_{j+1} U_{j+1} \times_{j+2} U_{j+2}\right\|_{\mathrm{F}}^2 \frac{\sigma^2 \op}{\ulambda^2} + \sum_{j=1}^3 \left\|\mcA \times_j U_j\right\|_{\mathrm{F}}^2  \left(\frac{\sigma^2 \oR\log(\op)}{\ulambda^2} + \frac{\sigma^6  \op^3}{\ulambda^6}\right)  
+  \left\|\mcA\right\|_{\mathrm{F}}^2  \left(\frac{\sigma^4 \oR^2\log(\op)^2}{\ulambda^4} + \frac{\sigma^{12} \op^6}{\ulambda^{12}}\right).
\end{align*}

\subsubsection*{Step 3: Summary of results}

On the other hand, we have
\begin{align*}
\left|\langle \mcA, \widehat{\mcT} \rangle -\left\langle \mcA, \mcT\right\rangle\right| 
\leq & \left|\left\langle \mcP_{\mathbb{T}_{\mcT} \mathcal{M}_{(r_1, r_2, r_3)}}\left(\mcA\right), \mcZ \right\rangle\right| + \left(\Omega_1 + \Omega_2 + \Omega_3 \right) \lesssim s_{\mcA}  \sigma\sqrt{\log\op} + \left(\Omega_1 + \Omega_2 + \Omega_3 \right) .
\end{align*}
where $\Omega_j, j=1,2,3$ are defined in the statement of Theorem~\ref{thm: main theorem in tensor PCA}.
Then,
\begin{align*}
& \Big|\frac{ \langle \mcA, \widehat{\mcT}  \rangle - \langle \mcA, \mcT \rangle}{\sigma\cdot s_{\mcA}} \cdot \big(\frac{\widehat{\sigma}\cdot \widehat{s}_{\mcA}}{\sigma \cdot s_{\mcA}} - 1\big)\Big| 
\leq  2 \left|\left\langle \mcA, \widehat{\mcT} \right\rangle -\left\langle \mcA, \mcT\right\rangle\right| \cdot \left( \left|\frac{\widehat{\sigma}}{\sigma}-1\right| + \left|\frac{\widehat{s}_\mcA}{s_\mcA}-1\right|\right) \\
\lesssim & \frac{s_{\mcA} \cdot \sigma\sqrt{\log(\op)} + \left(\Omega_1 + \Omega_2 + \Omega_3 + \Omega_4\right)}{\sigma\cdot s_{\mcA}} \cdot \Bigg\{\left(\kappa \sqrt{\frac{\oor \op}{p_1 p_2 p_3}}\right)  + \frac{1}{s_{\mcA}^2} \cdot \Bigg[\sum_{j=1}^3 \left\|\mcA \times_{j+1} U_{j+1} \times_{j+2} U_{j+2}\right\|_{\mathrm{F}}^2\cdot \frac{\sigma\cdot \sqrt{\op}}{\ulambda} \\
& + \sum_{j=1}^{3} \left\|\mcA \times_j U_j\right\|_{\mathrm{F}} \cdot \left(\frac{\sigma\cdot \sqrt{\oR\log(\op)}}{\ulambda} + \frac{\sigma^3 \cdot \op^{3/2}}{\ulambda^3}\right) + \left\|\mcA\right\|_{\mathrm{F}} \cdot \left(\frac{\sigma^2\cdot \oR\log(\op)}{\ulambda^2} + \frac{\sigma^{6}\cdot \op^3}{\ulambda^6}\right)\Bigg]\Bigg\} \\
\lesssim & \frac{\Omega_1 + \Omega_2 + \Omega_3}{\sigma\cdot s_{\mcA}} + \sqrt{\frac{\oor\op}{p_1p_2p_3}} + \frac{1}{\sigma\cdot s_{\mcA}} \cdot \Bigg\{ \sum_{j=1}^3 \left\|\mcA \times_{j+1} U_{j+1} \times_{j+2} U_{j+2}\right\|_{\mathrm{F}} \cdot \frac{\sigma^2\cdot \sqrt{\op\log(\op)}}{\ulambda} \\
+ & \sum_{j=1}^{3} \left\|\mcA \times_j U_j\right\|_{\mathrm{F}} \cdot \left(\frac{\sigma^2\cdot \sqrt{\oR}\log(\op)}{\ulambda} + \frac{\sigma^4 \cdot \op^{3/2}\sqrt{\log(\op)}}{\ulambda^3}\right) + \left\|\mcA\right\|_{\mathrm{F}} \cdot \left(\frac{\sigma^3\cdot \oR\log(\op)^{3/2}}{\ulambda^2} + \frac{\sigma^{7}\cdot \op^3\sqrt{\log(\op)}}{\ulambda^6}\right)\Bigg\}.
\end{align*}
Compared with $\Omega_j$ in Theorem~\ref{thm: main theorem in tensor PCA}, we have the desired results.

\section{Proof of Minimax Optimal Length of Confidence Interval}\label{sec: Proof of Minimax Optimal Length of Confidence Interval}

\subsection{Proof of Minimax Optimal Length of Confidence Interval in Tensor Regression (Theorem~\ref{thm: minimax lower bound in tensor regression})}\label{subsec: Proof of Minimax Optimal Length of Confidence Interval in Tensor Regression}


Define
$$
\overline{\mcT}:=\left(\mcG+\delta \mcA\times_1 U_1^{\top}\times_2 U_2^{\top}\times_3 U_3^{\top}\right) \times_1 \left(U_1+ \varepsilon\mcP_{U_{1\perp}}\Delta_1\right) \times_2 U_2 \times_3 U_3.
$$
It follows that
$$
\overline{\mcT}-\mcT= \varepsilon\mcG\times_1 \mcP_{U_{1\perp}}\Delta_1 \times_2 U_2 \times_3 U_3+ \delta \mcA\times_1 \mcP_{U_1} \times_2 \mcP_{U_2} \times_3 \mcP_{U_3}.
$$
Then, we have
$$
\left\|\overline{\mcT}-\mcT\right\|_{\mathrm{F}}^2= \varepsilon^2\left\|\mcG\times_1 \mcP_{U_{1\perp}}\Delta_1 \times_2 U_2 \times_3 U_3\right\|_{\mathrm{F}}^2 + \delta^2 \left\|\mcA\times_1 \mcP_{U_1} \times_2 \mcP_{U_2} \times_3 \mcP_{U_3}\right\|_{\mathrm{F}}^2 .
$$

Let $\Delta_1=A_1\left(U_3\otimes U_2\right)G_1^{\top}\left(G_1G_1^{\top}\right)^{-1}$. It follows immediately that
$$
\left\|\overline{\mcT}-\mcT\right\|_{\mathrm{F}}^2= \varepsilon^2\left\|\mcP_{U_{1\perp}}A_1\mcP_{\left(U_3\otimes U_2\right)G_1^{\top}}\right\|_{\mathrm{F}}^2 + \delta^2 \left\|\mcA\times_1 \mcP_{U_1} \times_2 \mcP_{U_2} \times_3 \mcP_{U_3}\right\|_{\mathrm{F}}^2 .
$$

By Lemma 1 of \citet{cai2015confidence}, we have that
\begin{align*}
& \inf_{\mathrm{CI}_{\mcA }^\alpha(\mcT, \mathcal{D}) \in \mathcal{I}_\alpha(\Theta,\mcA)} \sup_{ \mcT \in \Theta(\ulambda, \kappa)}\mathbb{E} L\left(\mathrm{CI}_{\mcA}^\alpha(\mcT, \mathcal{D} )\right) \geq \left|\left\langle \mcT-\overline{\mcT}, \mcA\right\rangle\right| \left(1-2 \alpha-\sqrt{\chi^2\left(f_{\pi_{\mathcal{H}_1}}, f_{\pi_{\mathcal{H}_0}}\right)}\right).
\end{align*}


Let $V_jV_j^{\top}=\mcP_{\left(U_{j+2}\otimes U_{j+1}\right)G_j^{\top}}$ and let
\begin{align*}
& a_{i}=\operatorname{Vec}\left[\mcX_i \times_1 U_1 \times_2 U_2 \times_3 U_3\right] \in \mathbb{R}^{r_1r_2r_3}, \\
& b_{i}= \operatorname{Vec}\left[\mcX_i \times_1 U_{1\perp} \times_2 U_2 \times_3 U_3\right] \in \mathbb{R}^{\left(p_1-r_1\right)r_2r_3}, \\
& c_{i}=
\left(\begin{array}{l}
\operatorname{Vec}\left[\mcX_i \times_1 U_1 \times_2 U_{2\perp} \times_3 U_3\right] \\
\operatorname{Vec}\left[\mcX_i \times_1 U_1 \times_2 U_2 \times_3 U_{3\perp}\right] \\
\operatorname{Vec}\left[\mcX_i \times_1 U_{1\perp} \times_2 U_{2\perp} \times_3 U_3\right] \\
\operatorname{Vec}\left[\mcX_i \times_1 U_{1\perp} \times_2 U_2 \times_3 U_{3\perp}\right] \\
\operatorname{Vec}\left[\mcX_i \times_1 U_1 \times_2 U_{2\perp} \times_3 U_{3\perp}\right] \\
\operatorname{Vec}\left[\mcX_i \times_1 U_{1\perp} \times_2 U_{2\perp} \times_3 U_{3\perp}\right] \\
\end{array}\right) \in \mathbb{R}^{p_1p_2p_3 - p_1r_2r_3}.
\end{align*}

It then follows that under $H_0$, we have
\begin{align*}
\left(\begin{array}{l}
Y_i \\
a_i \\
b_i \\
c_i \\
\end{array}\right) 
& \sim \mathcal{N}
\left(\begin{array}{llll}
\sigma^2\left\|\mcT\right\|_{\mathrm{F}}+\sigma_\xi^2 & \sigma^2\operatorname{Vec}\left(\mcG\right)^{\top}& & \\
\sigma^2\operatorname{Vec}\left(\mcG\right) & \sigma^2\mathcal{I}_{r_1r_2r_3} & & \\
& & \sigma^2\mathcal{I}_{\substack{\left(p_1-r_1\right)r_2r_3}} &\\
& & & \sigma^2\mathcal{I}_{p_1(p_2p_3-r_2r_3)} \\
\end{array}\right) \sim \mathcal{N}(0, \Sigma_{H_0}).
\end{align*}
and under $H_1$, we have
\begin{align*}
\left(\begin{array}{l}
Y_i \\
a_i \\
b_i \\
c_i \\
\end{array}\right) 
& \sim
\resizebox{\textwidth}{!}{$
\begin{pmatrix}
\sigma^2\left\|\overline{\mcT}\right\|_{\mathrm{F}}+\sigma_\xi^2 & \sigma^2\operatorname{Vec}\left(\mcG\right)^{\top}& \varepsilon\operatorname{Vec}\left(U_{1\perp}^{\top}A_1\mcP_{\left(U_3\otimes U_2\right)G_1^{\top}}\left(U_3\otimes U_2\right)\right)^{\top} & \\
\sigma^2\operatorname{Vec}\left(\mcG\right) & \sigma^2\mathcal{I}_{r_1r_2r_3} & & \\
\varepsilon\operatorname{Vec}\left(U_{1\perp}^{\top}A_1\mcP_{\left(U_3\otimes U_2\right)G_1^{\top}}\left(U_3\otimes U_2\right)\right) & & \sigma^2\mathcal{I}_{\substack{\left(p_1-r_1\right)r_2r_3}} &\\
& & & \sigma^2\mathcal{I}_{p_1(p_2p_3-r_2r_3)} \\
\end{pmatrix}
$}\\
& \sim \mathcal{N}(0, \Sigma_{H_1}).
\end{align*}

Let $\widetilde{a}_i=a_i-\frac{\sigma^2}{\sigma^2\left\|\mcT\right\|_{\mathrm{F}}^2+\sigma_\xi^2}Y_i\operatorname{Vec}\left(\mcG\right)$, then we have
\begin{align*}
\widetilde{\Sigma}_{H_0}
= & Q_{H_0}\Sigma_{H_0}Q_{H_0}^{-1} \\
= & \left(\begin{array}{llll}
\sigma^2\left\|\mcT\right\|_{\mathrm{F}}+\sigma_\xi^2 & & & \\
& \sigma^2\mathcal{I}_{r_1r_2r_3} - \frac{\sigma^4}{\sigma^2\left\|\mcT\right\|_{\mathrm{F}}^2+\sigma_\xi^2}\operatorname{Vec}\left(\mcG\right)\operatorname{Vec}\left(\mcG\right)^{\top} & & \\
& & \sigma^2\mathcal{I}_{\substack{\left(p_1-r_1\right)r_2r_3}} &\\
& & & \sigma^2\mathcal{I}_{p_1p_2p_3-p_1r_2r_3}
\end{array}\right).
\end{align*}

In addition, we have
\begin{align*}
\widetilde{\Sigma}_{H_1} 
= & Q_{H_0}\Sigma_{H_1}Q_{H_0}^{-1} \\
= & \resizebox{\textwidth}{!}{$
\begin{pmatrix}
\sigma^2\left\|\overline{\mcT}\right\|_{\mathrm{F}}+\sigma_\xi^2 & \sigma^2\left(1-\frac{\sigma^2\left\|\overline{\mcT}\right\|_{\mathrm{F}}^2+\sigma_\xi^2}{\sigma^2\left\|\mcT\right\|_{\mathrm{F}}^2+\sigma_\xi^2}\right)\operatorname{Vec}\left(\mcG\right)^{\top} & \sigma^2 \varepsilon \operatorname{Vec}\left(U_1^{\top} A_1 \mcP_{\left(U_3 \otimes U_2\right) G_1^{\top}}\left(U_3 \otimes U_2\right)\right)^{\top} & \\
\sigma^2\left(1-\frac{\sigma^2\left\|\overline{\mcT}\right\|_{\mathrm{F}}^2+\sigma_\xi^2}{\sigma^2\left\|\mcT\right\|_{\mathrm{F}}^2+\sigma_\xi^2}\right)\operatorname{Vec}\left(\mcG\right) & \sigma^2{I}_{r_1r_2r_3} - \frac{\sigma^4}{\sigma^2\left\|\mcT\right\|_{\mathrm{F}}^2+\sigma_\xi^2}\operatorname{Vec}\left(\mcG\right)\operatorname{Vec}\left(\mcG\right)^{\top} & & \\
\sigma^2\varepsilon \operatorname{Vec}\left(U_{1\perp}^{\top} A_1 \mcP_{\left(U_3 \otimes U_2\right) G_1^{\top}}\left(U_3 \otimes U_2\right)\right) & & \sigma^2\mathcal{I}_{\substack{\left(p_1-r_1\right)r_2r_3}} & \\
& & & \sigma^2\mathcal{I}_{p_1(p_2p_3-r_2r_3)} \\
\end{pmatrix}
$}.
\end{align*}

Therefore, 
\begin{align*}
& \widetilde{\Sigma}_{H_0}^{-\frac{1}{2}}\widetilde{\Sigma}_{H_1}\widetilde{\Sigma}_{H_0}^{-\frac{1}{2}}\\
= & 
\resizebox{\textwidth}{!}{$
\begin{pmatrix}
\frac{\sigma^2\left\|\overline{\mcT}\right\|_{\mathrm{F}}+\sigma_\xi^2}{\sigma^2\left\|\mcT\right\|_{\mathrm{F}}+\sigma_\xi^2} & \sigma^2\left(1-\frac{\sigma^2\left\|\overline{\mcT}\right\|_{\mathrm{F}}^2+\sigma_\xi^2}{\sigma^2\left\|\mcT\right\|_{\mathrm{F}}^2+\sigma_\xi^2}\right)\operatorname{Vec}\left(\mcG\right)^{\top} & & \\
\sigma^2\left(1-\frac{\sigma^2\left\|\overline{\mcT}\right\|_{\mathrm{F}}^2+\sigma_\xi^2}{\sigma^2\left\|\mcT\right\|_{\mathrm{F}}^2+\sigma_\xi^2}\right)\operatorname{Vec}\left(\mcG\right) & \sigma^2\mathcal{I}_{r_1r_2r_3} - \frac{\sigma^4}{\sigma^2\left\|\mcT\right\|_{\mathrm{F}}^2+\sigma_\xi^2}\operatorname{Vec}\left(\mcG\right)\operatorname{Vec}\left(\mcG\right)^{\top} & & \\
& & \sigma^2\mathcal{I}_{\substack{\left(p_1-r_1\right)r_2r_3}} &\\
& & & \sigma^2\mathcal{I}_{p_1(p_2p_3-r_2r_3)} \\
\end{pmatrix}
$}.
\end{align*}
It follows that
\begin{align*}
& \left\|\Sigma_{H_0}^{-\frac{1}{2}}\Sigma_{H_1}\Sigma_{H_0}^{-\frac{1}{2}}-\mathcal{I}\right\|_{\mathrm{F}}^2 = \left\|Q_{H_0}^{-1}\widetilde{\Sigma}_{H_0}^{-\frac{1}{2}}Q_{H_0} \cdot Q_{H_0}^{-1}\widetilde{\Sigma}_{H_1}Q_{H_0} \cdot Q_{H_0}^{-1}\widetilde{\Sigma}_{H_0}^{-\frac{1}{2}}Q_{H_0}-\mathcal{I}\right\|_{\mathrm{F}}^2 = \left\|\widetilde{\Sigma}_{H_0}^{-\frac{1}{2}}\widetilde{\Sigma}_{H_1}\widetilde{\Sigma}_{H_0}^{-\frac{1}{2}}-\mathcal{I}\right\|_{\mathrm{F}}^2\\
= &  \frac{\varepsilon^4\sigma^4\left\|U_{1\perp}^{\top}A_1\mcP_{\left(U_3 \otimes U_2\right) G_1^{\top}}\left(U_3 \otimes U_2\right)\right\|_{\mathrm{F}}^4}{\left(\sigma^2\|\mcT\|_{\mathrm{F}}^2
+\sigma_{\xi}^2\right)^2} 
+ \frac{2\sigma^4\cdot \sigma^4\varepsilon^4\left\|U_{1\perp}^{\top}A_1\mcP_{\left(U_3 \otimes U_2\right) G_1^{\top}}\left(U_3 \otimes U_2\right)\right\|_{\mathrm{F}}^4}{\left(\sigma^2\|\mcT\|_{\mathrm{F}}^2
+\sigma_{\xi}^2\right)^3} \left(\frac{\left\|\mcT\right\|_{\mathrm{F}}^2}{\sigma^2}+\frac{\left\|\mcT\right\|_{\mathrm{F}}^4}{\sigma_\xi^2}\right) \\
& + \frac{2\varepsilon^2\sigma^2\left\|U_{1\perp}^{\top}A_1\mcP_{\left(U_3 \otimes U_2\right) G_1^{\top}}\left(U_3 \otimes U_2\right)\right\|_{\mathrm{F}}^2}{\sigma^2\|\mcT\|_{\mathrm{F}}^2+\sigma_{\xi}^2} \\
\lesssim & 2 \cdot \frac{\varepsilon^2\sigma^4\left\|U_{1\perp}^{\top}A_1\mcP_{\left(U_3 \otimes U_2\right) G_1^{\top}}\right\|_{\mathrm{F}}^2}{\sigma_{\xi}^2}+ 2\cdot \frac{\sigma^6\varepsilon^4\left\|U_{1\perp}^{\top}A_1\mcP_{\left(U_3 \otimes U_2\right) G_1^{\top}}\right\|_{\mathrm{F}}^4}{\sigma_\xi^4}.
\end{align*}

By Theorem 1.1 of \citet{devroye2018total}, we have
$$
D_\text{TV}\left(f_{\pi_{\mathcal{H}_1}}, f_{\pi_{\mathcal{H}_0}}\right) \leq 2\sqrt{n}\left\|\Sigma_{\mathcal{H}_0}^{-1 / 2} \Sigma_{\mathcal{H}_1} \Sigma_{\mathcal{H}_0}^{-1 / 2}-I_d\right\|_\mathrm{F}.
$$

Therefore, let 
$$
\varepsilon=a\cdot \frac{\sigma_\xi}{\sigma\sqrt{n}}\left\|\mcG\times_1 \mcP_{U_{1\perp}}\Delta_1 \times_2 U_2 \times_3 U_3\right\|_{\mathrm{F}}^{-1}=a\cdot \frac{\sigma_\xi}{\sigma\sqrt{n}}\left\|\mcP_{U_{1\perp}}A_1\mcP_{\left(U_3\otimes U_2\right)G_1^{\top}}\right\|_{\mathrm{F}}^{-1}.
$$
It then follows that
$$
D_\text{TV}\left(f_{\pi_{\mathcal{H}_1}}, f_{\pi_{\mathcal{H}_0}}\right) \leq 2\sqrt{n}\left\|\Sigma_{\mathcal{H}_0}^{-1 / 2} \Sigma_{\mathcal{H}_1} \Sigma_{\mathcal{H}_0}^{-1 / 2}-I_d\right\|_\mathrm{F} \leq 2\sqrt{2n}\cdot \sigma\cdot \frac{\sigma_\xi}{\sigma}\sqrt{a^2+a^4} 
\leq 4a\sigma
$$
as long as $a$ is sufficiently small such that $a\leq 1$.

By the assumption of $\mcT$ that $\ulambda \leq \lambda_{\text{max}}(\mcT) \leq \frac{\kappa+1}{2}\ulambda$, we have
$$
\ulambda+ \left|b\right| -\left|a\right| \leq \lambda_{\text{max}}(\mcT) \leq \frac{\kappa+1}{2}\ulambda+\left|b\right| +\left|a\right|.
$$

Therefore, we have
\begin{align*}
\inf_{\mathrm{CI}_{\mcA }^\alpha(\mcT, \mathcal{D}) \in \mathcal{I}_\alpha(\Theta,\mcA)} \sup_{ \mcT \in \Theta(\ulambda, \kappa)}\mathbb{E} L\left(\mathrm{CI}_{\mcA}^\alpha(\mcT, \mathcal{D} )\right) 
& \geq \frac{\sigma_\xi}{\sigma\sqrt{n}}\left|a\right|\cdot \left\|\mcP_{U_{1\perp}}A_1\mcP_{\left(U_3\otimes U_2\right)G_1^{\top}}\right\|_{\mathrm{F}} \cdot \left(1-2\alpha-2a\sigma\right) \\
& \geq \frac{c_1\sigma}{\sigma\sqrt{n}} \cdot \left\|\mcP_{U_{1\perp}}A_1\mcP_{\left(U_3\otimes U_2\right)G_1^{\top}}\right\|_{\mathrm{F}}.
\end{align*}

By the same argument, for any $j=1,2,3$, we have
\begin{align*}
\inf_{\mathrm{CI}_{\mcA }^\alpha(\mcT, \mathcal{D}) \in \mathcal{I}_\alpha(\Theta,\mcA)} \sup_{ \mcT \in \Theta(\ulambda, \kappa)}\mathbb{E} L\left(\mathrm{CI}_{\mcA}^\alpha(\mcT, \mathcal{D} )\right) 
& \geq \frac{c_1\sigma}{\sigma\sqrt{n}} \cdot \left\|\mcP_{U_{j\perp}}A_j\mcP_{\left(U_{j+2}\otimes U_{j+1}\right)G_j^{\top}}\right\|_{\mathrm{F}}.
\end{align*}

By the same argument, we have
\begin{align*}
& \left\|\Sigma_{H_0}^{-\frac{1}{2}}\Sigma_{H_1}\Sigma_{H_0}^{-\frac{1}{2}}-{I}\right\|_{\mathrm{F}}^2 \\
= & \left(\frac{\sigma^2\|\overline{\mcT}\|_{\mathrm{F}}^2+\sigma_{\xi}^2}{\sigma^2\|\mcT\|_{\mathrm{F}}^2
+\sigma_{\xi}^2}-1\right)^2 + 2\cdot \frac{1}{\sigma^2\|\mcT\|_{\mathrm{F}}^2+\sigma_{\xi}^2}\cdot \left\|\sigma^2\left(1-\frac{\sigma^2\|\overline{\mcT}\|_{\mathrm{F}}^2+\sigma_{\varepsilon}^2}{\sigma^2\|\mcT\|_{\mathrm{F}}^2+\sigma_{\xi}^2}\right) \operatorname{Vec}(\mcG)\left(\frac{1}{\sigma^2} \mathcal{I}_{r_1 r_2 r_3}+\frac{1}{\sigma_\xi^2} \operatorname{Vec}(\mcG) \operatorname{Vec}(\mcG)^{\top}\right)^{1/2} \right\|_{\ell_2}^2\\
= &  \frac{\varepsilon^4\sigma^4\left\|\mcA \times_1 U_1 \times_2 U_2 \times_3 U_3\right\|_{\mathrm{F}}^4}{\left(\sigma^2\|\mcT\|_{\mathrm{F}}^2
+\sigma_{\xi}^2\right)^2} + 2\cdot \frac{\varepsilon^2\sigma^2\left\|\mcA \times_1 U_1 \times_2 U_2 \times_3 U_3\right\|_{\mathrm{F}}^2}{\sigma^2\|\mcT\|_{\mathrm{F}}^2+\sigma_{\xi}^2} \\
\lesssim & 2 \cdot \frac{\varepsilon^2\sigma^4\left\|\mcA \times_1 U_1 \times_2 U_2 \times_3 U_3\right\|_{\mathrm{F}}^2}{\sigma_{\xi}^2}+ 2\cdot \frac{\sigma^6\varepsilon^4\left\|\mcA \times_1 U_1 \times_2 U_2 \times_3 U_3\right\|_{\mathrm{F}}^4}{\sigma_\xi^4} .
\end{align*}

Therefore, let 
$
\delta=b\cdot \frac{\sigma_\xi}{\sigma\sqrt{n}}\left\|\mcA \times_1 U_1 \times_2 U_2 \times_3 U_3\right\|_{\mathrm{F}}^{-1}.
$
It then follows that
$$
D_\text{TV}\left(f_{\pi_{\mathcal{H}_1}}, f_{\pi_{\mathcal{H}_0}}\right) \leq 2\sqrt{n}\left\|\Sigma_{\mathcal{H}_0}^{-1 / 2} \Sigma_{\mathcal{H}_1} \Sigma_{\mathcal{H}_0}^{-1 / 2}-I_d\right\|_\mathrm{F} \leq 2\sqrt{2n}\cdot \sigma\cdot \frac{\sigma_\xi}{\sigma}\sqrt{b^2+b^4} \\
\leq 4b\sigma
$$
as long as a is sufficiently small such that $a\leq 1$

Therefore, we have
\begin{align*}
\inf_{\mathrm{CI}_{\mcA }^\alpha(\mcT, \mathcal{D}) \in \mathcal{I}_\alpha(\Theta,\mcA)} \sup_{ \mcT \in \Theta(\ulambda, \kappa)}\mathbb{E} L\left(\mathrm{CI}_{\mcA}^\alpha(\mcT, \mathcal{D} )\right)
& \geq \frac{\sigma_\xi}{\sigma\sqrt{n}}\left|a\right|\cdot \left\|\mcA \times_1 U_1 \times_2 U_2 \times_3 U_3\right\|_{\mathrm{F}} \cdot \left(1-2\alpha-2b\sigma\right) \\
& \geq \frac{c_4\sigma}{\sigma\sqrt{n}} \cdot \left\|\mcA \times_1 U_1 \times_2 U_2 \times_3 U_3\right\|_{\mathrm{F}}.
\end{align*}

Combining all the results above, we then have
\begin{align*}
& \inf_{\mathrm{CI}_{\mcA }^\alpha(\mcT, \mathcal{D}) \in \mathcal{I}_\alpha(\Theta,\mcA)} \sup_{ \mcT \in \Theta(\ulambda, \kappa)}\mathbb{E} L\left(\mathrm{CI}_{\mcA}^\alpha(\mcT, \mathcal{D} )\right) \\
& \geq c^{\prime} \frac{\sigma_\xi}{\sigma\sqrt{n}} \left(\sum_{j=1}^3\left\|\mcP_{U_{j\perp}}A_j\mcP_{\left(U_{j+2}\otimes U_{j+1}\right)G_j^{\top}}\right\|_{\mathrm{F}} +  \left\|\mcA\times_1 \mcP_{U_1} \times_2 \mcP_{U_2} \times_3 \mcP_{U_3}\right\|_{\mathrm{F}}\right) . 
\end{align*}

\subsection{Proof of Minimax Optimal Length of Confidence Interval in Tensor PCA (Theorem~\ref{thm: minimax lower bound in tensor PCA})}\label{subsec: Proof of Minimax Optimal Length of Confidence Interval in Tensor PCA}


Define
$$
\overline{\mcT}:=\left(\mcG+\delta \mcA\times_1 U_1^{\top}\times_2 U_2^{\top}\times_3 U_3^{\top}\right) \times_1 \left(U_1+ \varepsilon\mcP_{U_{1\perp}}\Delta_1\right) \times_2 U_2 \times_3 U_3.
$$
Then it follows that
$$
\overline{\mcT}-\mcT= \varepsilon\mcG\times_1 \mcP_{U_{1\perp}}\Delta_1 \times_2 U_2 \times_3 U_3+ \delta \mcA\times_1 \mcP_{U_1} \times_2 \mcP_{U_2} \times_3 \mcP_{U_3}.
$$

Therefore, first, we have
$$
\left\|\overline{\mcT}-\mcT\right\|_{\mathrm{F}}^2= \varepsilon^2\left\|\mcG\times_1 \mcP_{U_{1\perp}}\Delta_1 \times_2 U_2 \times_3 U_3\right\|_{\mathrm{F}}^2 + \delta^2 \left\|\mcA\times_1 \mcP_{U_1} \times_2 \mcP_{U_2} \times_3 \mcP_{U_3}\right\|_{\mathrm{F}}^2.
$$

Let $\Delta_1=A_1\left(U_3\otimes U_2\right)G_1^{\top}\left(G_1G_1^{\top}\right)^{-1}$. It follows immediately that
$$
\left\|\overline{\mcT}-\mcT\right\|_{\mathrm{F}}^2= \varepsilon^2\left\|\mcP_{U_{1\perp}}A_1\mcP_{\left(U_3\otimes U_2\right)G_1^{\top}}\right\|_{\mathrm{F}}^2 + \delta^2 \left\|\mcA\times_1 \mcP_{U_1} \times_2 \mcP_{U_2} \times_3 \mcP_{U_3}\right\|_{\mathrm{F}}^2.
$$

By Lemma 1 of \citet{cai2015confidence}, we have that
\begin{align*}
& \inf_{\mathrm{CI}_{\mcA }^\alpha(\mcT, \mathcal{D}) \in \mathcal{I}_\alpha(\Theta,\mcA)} \sup_{ \mcT \in \Theta(\ulambda, \kappa)}\mathbb{E} L\left(\mathrm{CI}_{\mcA}^\alpha(\mcT, \mathcal{D} )\right)  \geq \left|\left\langle \mcT-\overline{\mcT} , \mcA\right\rangle\right| \left(1-2 \alpha-\sqrt{\chi^2\left(f_{\pi_{\mathcal{H}_1}}, f_{\pi_{\mathcal{H}_0}}\right)}\right).
\end{align*}

Note that under the entry-wise i.i.d. Gaussian assumption
$$
\chi^2\left(f_{\pi_{\mathcal{H}_1}}, f_{\pi_{\mathcal{H}_0}}\right) \leq \exp\left(\frac{1}{\sigma^2} \cdot \left\|\overline{\mcT}-\mcT\right\|_{\mathrm{F}}^2\right) - 1.
$$

Assume that $\mcG$ satisfies that $\ulambda \leq \lambda_{\text{max}}(\mcG) \leq \frac{\kappa+1}{2}\ulambda$. Then, by Weyl's inequality for singular values, we have
\begin{align*}
\left|\sigma_j\left(\operatorname{Mat}_1\left(\overline{T}\right)\right)-\sigma_j\left(\operatorname{Mat}_1\left(T\right)\right)\right|
\leq & \sigma_1\left(\operatorname{Mat}_1\left(\varepsilon\mcG\times_1 \mcP_{U_{1\perp}}\Delta_1 \times_2 U_2 \times_3 U_3+ \delta \mcA\times_1 \mcP_{U_1} \times_2 \mcP_{U_2} \times_3 \mcP_{U_3}\right)\right) \\
\leq & \varepsilon\left\|\mcP_{U_{1\perp}}A_1\mcP_{\left(U_3 \otimes U_2\right)G_1^{\top}}\right\|+ \delta\left\|\mcP_{U_1}A_1\left(\mcP_{U_3} \otimes \mcP_{U_2}\right)\right\|.
\end{align*}
Therefore, let 
$
\varepsilon=a\left\|\mcG\times_1 \mcP_{U_{1\perp}}\Delta_1 \times_2 U_2 \times_3 U_3\right\|_{\mathrm{F}}^{-1}=a\left\|\mcP_{U_{1\perp}}A_1\mcP_{\left(U_3\otimes U_2\right)G_1^{\top}}\right\|_{\mathrm{F}}^{-1}
$ and
$
\delta=b\left\|\mcA\times_1 \mcP_{U_1} \times_2 \mcP_{U_2} \times_3 \mcP_{U_3}\right\|_{\mathrm{F}}^{-1}.
$
It then follows that
\begin{align*}
\left|\sigma_j\left(\operatorname{Mat}_1\left(\overline{T}\right)\right)-\sigma_j\left(\operatorname{Mat}_1\left(T\right)\right)\right|
\leq & \sigma_1\left(\varepsilon\mcG\times_1 \mcP_{U_{1\perp}}\Delta_1 \times_2 U_2 \times_3 U_3+ \delta \mcA\times_1 \mcP_{U_1} \times_2 \mcP_{U_2} \times_3 \mcP_{U_3}\right) \\
\leq & \left|a\right|+\left|b\right|.
\end{align*}
By the assumption that $\mcT$ satisfies that $\ulambda \leq \lambda_{\text{max}}(\mcT) \leq \frac{\kappa+1}{2}\ulambda$, we have
$$
\ulambda+ \left|b\right| -\left|a\right| \leq \lambda_{\text{max}}(\mcT) \leq \frac{\kappa+1}{2}\ulambda+\left|b\right| +\left|a\right|.
$$
Let $a=b$, then it follows that
$
\ulambda \leq \lambda_{\text{max}}(\mcT) \leq \frac{\kappa+1}{2}\ulambda+2\left|a\right|.
$
When $\left|a\right|$ is sufficiently small, then we have
$
\ulambda \leq \lambda_{\text{max}}(\mcT) \leq \frac{\kappa+1}{2}\ulambda+2\left|a\right| \leq \kappa\ulambda.
$

It only remains to find an upper bound for $\left\langle \overline{\mcT} - \mcT, \mcA\right\rangle$. Note that now
\begin{align*}
\left|\left\langle \overline{\mcT} - \mcT, \mcA\right\rangle\right|
= & \left|a\right|\cdot \left(\left\|\mcG \times \times_1 P_{U_{1 \perp}} \Delta_1 \times_2 U_2 \times_3 U_3\right\|_{\mathrm{F}}+\left\|\mcA \times_1 \mcP_{U_1} \times_2 \mcP_{U_2} \times_3 \mcP_{U_3}\right\|_{\mathrm{F}}\right) \\
= & \left|a\right|\cdot \left(\left\|\mcP_{U_{1\perp}}A_1\mcP_{\left(U_3\otimes U_2\right)G_1^{\top}}\right\|_{\mathrm{F}} +  \left\|\mcA\times_1 \mcP_{U_1} \times_2 \mcP_{U_2} \times_3 \mcP_{U_3}\right\|_{\mathrm{F}}\right).
\end{align*}

Thefore, we have
\begin{align*}
& \inf_{\mathrm{CI}_{\mcA }^\alpha(\mcT, \mathcal{D}) \in \mathcal{I}_\alpha(\Theta,\mcA)} \sup_{ \mcT \in \Theta(\ulambda, \kappa)}\mathbb{E} L\left(\mathrm{CI}_{\mcA}^\alpha(\mcT, \mathcal{D} )\right) \\
& \geq \left|a\right|\cdot \left(\left\|\mcP_{U_{1\perp}}A_1\mcP_{\left(U_3\otimes U_2\right)G_1^{\top}}\right\|_{\mathrm{F}} +  \left\|\mcA\times_1 \mcP_{U_1} \times_2 \mcP_{U_2} \times_3 \mcP_{U_3}\right\|_{\mathrm{F}}\right)\left(1-2\alpha-\sqrt{\sigma^2\exp \left(2a^2\right)-1}\right) \\
& \geq c_1 \sigma \left(\left\|\mcP_{U_{1\perp}}A_1\mcP_{\left(U_3\otimes U_2\right)G_1^{\top}}\right\|_{\mathrm{F}} +  \left\|\mcA\times_1 \mcP_{U_1} \times_2 \mcP_{U_2} \times_3 \mcP_{U_3}\right\|_{\mathrm{F}}\right).
\end{align*}

By the same arguments, for any $j=1,2,3$, we have
\begin{align*}
& \inf_{\mathrm{CI}_{\mcA }^\alpha(\mcT, \mathcal{D}) \in \mathcal{I}_\alpha(\Theta,\mcA)} \sup_{ \mcT \in \Theta(\ulambda, \kappa)}\mathbb{E} L\left(\mathrm{CI}_{\mcA}^\alpha(\mcT, \mathcal{D} )\right) \\
\gtrsim & \sigma \left\|\mcP_{U_{j\perp}}A_j\mcP_{\left(U_{j+2}\otimes U_{j+1}\right)G_h^{\top}}\right\|_{\mathrm{F}} +  \sigma\left\|\mcA\times_1 \mcP_{U_1} \times_2 \mcP_{U_2} \times_3 \mcP_{U_3}\right\|_{\mathrm{F}}.
\end{align*}

Combining all results above, we then have the desired lower bound.

\section{Specific Examples} \label{sec: Specific Examples}

Our framework is versatile and applicable to inferring any linear functional of the signal tensor $\mcT$. In this section, we demonstrate its utility through three specific scenarios: entrywise inference, comparing tensor entries, and inference of the average of entries along a mode. All confidence intervals discussed below are valid under the conditions specified in Corollary~\ref{corollary: asymptotic normality of estimated low-Tucker-rank linear functional in tensor regression with sample splitting} for tensor regression or Corollary~\ref{corollary: asymptotic normality of estimated low-Tucker-rank linear functional in tensor PCA} for tensor PCA. For the inference of general (full-rank) linear functionals, the conditions for valid inference are provided in Corollary~\ref{corollary: asymptotic normality of estimated general linear functional in tensor regression without sample splitting} for tensor regression and Corollary~\ref{corollary: asymptotic normality of estimated general linear functional in tensor PCA} for tensor PCA.

\subsection{Entrywise inference}

Entrywise inference aims to provide statistical statements about individual entries of the signal tensor $\mcT$. This task has been extensively studied in matrix PCA \citep{chen2019inference, ma2024statistical, farias2022uncertainty} and recently extended to tensor PCA \citep{xia2022inference, agterberg2024statistical}.

To perform entrywise inference for a specific tensor entry $\mcT_{k_1,k_2,k_3}$, we consider the linear functional $\langle \mcA, \mcT \rangle = \langle e_{k_1} \otimes e_{k_2} \otimes e_{k_3} , \mcT \rangle$. Applying Theorem~\ref{thm: main theorem in tensor regression with plug-in estimates} and Theorem~\ref{thm: main theorem in tensor PCA with plug-in estimates}, we construct the following confidence intervals in tensor regression,
\begin{align*}
\widehat{\mathrm{CI}}_{\mcA, T}^{\alpha}
:= & \left[\widehat{\mcT}_{k_1, k_2, k_3} - z_{\alpha / 2} \cdot \frac{\widehat{\sigma}_{\xi}}{\widehat{\sigma}}\cdot \widehat{s}_{\mcA} \cdot \sqrt{\frac{1}{n}}, \widehat{\mcT}_{k_1, k_2, k_3} + z_{\alpha / 2} \cdot \frac{\widehat{\sigma}_{\xi}}{\widehat{\sigma}}\cdot \widehat{s}_{\mcA} \cdot \sqrt{\frac{1}{n}}\right]
\end{align*}
where $\widehat{\mcT}$ is obtained via the debiasing procedure in Section \ref{sec:debias_nonsplit} and Section \ref{sec: Debiased Estimator of Linear Functionals with Sample Splitting}. 
In the tensor PCA setting, the entrywise confidence interval is of the same form with $(\widehat{\sigma}_\xi /\widehat{\sigma})\cdot \widehat{s}_{\mcA} \cdot \sqrt{1/n}$ replaced by $\widehat{\sigma} \cdot \widehat{s}_{\mcA}$, where $\widehat{\mcT}$ is the output from the procedure in Section \ref{sec:tpca_alg}. 
The variance component is estimated as
\begin{align*}
\widehat{s}_{\mcA}^2 
= & \sum_{j=1}^3 \left\|\left(\mathcal{I} - \mcP_{\whU_j}\right)e_{k_j}\right\|\left\|\whU_{j+1}^{\top} e_{k_{j+1}}\right\|_{\ell_2}^2\left\|\whU_{j+2}^{\top} e_{k_{j+2}}\right\|_{\ell_2}^2 + \prod_{j=1}^3  \left\|\whU_j^{\top} e_{k_j}\right\|_{\ell_2}^2.
\end{align*}

\subsection{Inference for the Row/Column Mean}

Researchers often aim to understand the average effect across one mode of the signal tensor. For instance, in collaborative filtering, this could involve inferring the average rating an item receives among all users \citep{frolov2017tensor, marin2022tensor}. In addition, spatial patterns of gene regulation can be identified by averaging gene expression measurements over all time points \citep{liu2022characterizing}.

To infer the average effect across the first mode, consider the loading tensor $\mcA=\frac{1}{p_1}\sum_{k_1=1}^{p_1}e_{k_1}\otimes e_{k_2} \otimes e_{k_3} \in \mathbb{R}^{p_1\times p_2\times p_3}$, which has a Tucker rank of $(1, 1, 1)$. The confidence interval for the average effect across the first mode of $\mcT$ in tensor regression is given by 
\begin{align*}
\widehat{\mathrm{CI}}_{\mcA, T}^{\alpha}
:=\left[\frac{1}{p_1}\sum_{k_1=1}^{p_1}\widehat{\mcT}_{k_1, k_2, k_3}-z_{\alpha / 2} \cdot \frac{\widehat{\sigma}_{\xi}}{\widehat{\sigma}} \cdot \widehat{s}_{\mcA} \cdot \sqrt{\frac{1}{n}}, \frac{1}{p_1}\sum_{k_1=1}^{p_1}\widehat{\mcT}_{k_1, k_2, k_3} + z_{\alpha / 2} \cdot \frac{\widehat{\sigma}_{\xi}}{\widehat{\sigma}} \cdot \widehat{s}_{\mcA} \cdot \sqrt{\frac{1}{n}}\right],
\end{align*}
In the tensor PCA setting, the confidence interval is constructed similarly with $(\widehat{\sigma}_\xi /\widehat{\sigma})\cdot \widehat{s}_{\mcA} \cdot \sqrt{1/n}$ replaced by $\widehat{\sigma} \cdot \widehat{s}_{\mcA}$. 
Here, the variance component $s_{\mcA}^2$ is estimated as
\begin{align*}
\widehat{s}_{\mcA}^2 
= & \frac{1}{p_1}\sum_{k_1=1}^{p_1}\sum_{j=1}^{3}\left\|\mcP_{\whU_{j\perp}}e_{k_j}\right\|_{\ell_2}^2\left\|\widehat{W}_{k_j}\left(\whU_{k_{j+2}}^{\top}e_{k_{j+2}} \otimes \whU_{k_{j+1}}^{\top}e_{k_{j+1}} \right)\right\|_{\ell_2}^2 + \frac{1}{p_1}\sum_{k_1=1}^{p_1}\left\|\whU_1^{\top}e_{k_1}\right\|_{\ell_2}^2 \prod_{j=2}^3 \left\|\whU_j^{\top}e_{k_j}\right\|_{\ell_2}.
\end{align*}

\section{Summary of Inference Procedure} \label{sec: Inference procedure for tensor regression}

\subsection{Inference procedure for tensor regression } \label{sec: Inference procedure for tensor regression without sample splitting}

We first outlines the steps to estimate the linear functional $\left\langle \mcA, \mcT \right\rangle$ for tensor regression when the entire dataset is used jointly without splitting. The summarized procedure is detailed in Algorithm \ref{alg: tensor regression without sample splitting}. Furthermore, with sample splitting, the procedure summarized in Algorithm~\ref{alg: tensor regression with sample splitting}, addresses the dependency between debiasing and estimation by partitioning the dataset into two independent subsets. Each subset is used to complement the estimation in the other.

\begin{algorithm}[hp]
\caption{Inference Procedure for Tensor Regression without Sample splitting}
\label{alg: tensor regression without sample splitting}
\begin{algorithmic}[1]
\State \textbf{Initialization}
\State \quad Obtain initial estimates $\widehat{\mcT}^{\text{init}} \in \mathbb{R}^{p_1 \times p_2 \times p_3}$ and factor matrices $\whU_j^{(0)} := \whU_j^{\text{init}}$ for $j = 1, 2, 3$.

\State \textbf{Debiasing}
\State \quad Compute the debiased tensor:
$
\widehat{\mcT}^{\text{unbs}} = \widehat{\mcT}^{\text{init}} + \frac{1}{n \sigma^2} \sum_{i=1}^{n} \left( Y_i - \langle \widehat{\mcT}^{\text{init}}, \mcX_i \rangle \right) \mcX_i.
$

\State \textbf{Two-step Power Iteration}
\For{$k = 1, 2$}
    \For{$j = 1, 2, 3$}
        \State Compute $\whU_j^{(k)}$ as the leading $r_j$ left singular vectors of
        $$
        \operatorname{Mat}_j\left( \widehat{\mcT}^{\text{unbs}} \times_{i \neq j} \whU_i^{(k-1)\top} \right) = \operatorname{Mat}_j\left( \widehat{\mcT}^{\text{unbs}} \right) \left( \otimes_{i \neq j} \whU_i^{(k-1)} \right).
        $$
    \EndFor
\EndFor
\State \textbf{Projection and Plug-in Estimator}
\State \quad Define projection matrices $\mcP_{\whU_j} = \whU_j^{(2)} \whU_j^{(2)\top}$ for $j = 1, 2, 3$. Compute the projected tensor:
$$
\widehat{\mcT} = \widehat{\mcT}^{\text{unbs}} \times_1 \mcP_{\whU_1} \times_2 \mcP_{\whU_2} \times_3 \mcP_{\whU_3}.
$$
Then, estimate the linear functional $\langle \mcA, \mcT \rangle$ by $\langle \mcA, \widehat{\mcT}\rangle.$

\State \textbf{Variance Estimation}

Let
$$
\widehat{\sigma}_{\xi}^2
= \frac{1}{n} \sum_{i=1}^{n} \left(Y_i-\left\langle\widehat{\mcT}^{\text{init}}, \mcX_i\right\rangle\right)^2
$$
be the estimate of the observational noise variance, and
$$
\widehat\sigma^2 = \frac{1}{np_1p_2p_3}\sum_{i=1}^n \left\|\mcX_i\right\|_{\mathrm{F}}^2
$$
be the estimate of design noise.

In addition, let 
\begin{align}
& \widehat{s}_{\mcA}^2 = \sum_{j=1}^3\left\|\left(I-\mcP_{\whU_j}\right)A_j\left(\whU_{j+2}\otimes \whU_{j+1}\right)\widehat{W}_j\widehat{W}_j^{\top}\left(\whU_{j+2}\otimes \whU_{j+1}\right)^{\top}\right\|_{\mathrm{F}}^2 + \left\|\mcA \times_1 \whU_1 \times_2 \whU_2 \times_3 \whU_3\right\|_{\mathrm{F}}^2, 
\end{align}
be the estimate of variance component ${s}_{\mcA}^2$, where
$ \widehat{W}_j=\mathrm{QR}\left[\Mat_j\left(\widehat{\mcT} \times_1 \whU_1^{\top} \times_2 \whU_2^{\top} \times_3 \whU_3^{\top} \right)^{\top}\right].
$

\State \textbf{Confidence Interval}

The $100(1 - \alpha)\%$ confidence interval is given by
$$
\widehat{\mathrm{CI}}_{\mcA, \cT}^{\alpha}=\left[\left\langle  \widehat{\mcT}, \mcA \right\rangle-z_{\alpha / 2} \cdot \frac{\widehat{\sigma}_{\xi}}{\widehat{\sigma}} \cdot \widehat{s}_{\mcA} \sqrt{\frac{1}{n}}, \left\langle \widehat{\mcT}, \mcA \right\rangle + z_{\alpha / 2} \cdot \frac{\widehat{\sigma}_{\xi}}{\widehat{\sigma}} \cdot \widehat{s}_{\mcA} \sqrt{\frac{1}{n}}\right]  ,
$$
where $\alpha \in (0, 1)$, and $z_\theta = \Phi^{-1}(1 - \theta)$ denotes the upper $\theta$ quantile of the standard normal distribution.

\end{algorithmic}
\end{algorithm}

\begin{algorithm}[hp]
\caption{Inference Procedure for Tensor Regression with Sample Splitting}
\label{alg: tensor regression with sample splitting}
\begin{algorithmic}[1]
\State \textbf{Initialization}
\State \quad Use dataset $\text{\RN{1}} = \left\{ Y_{i_1}^{(\text{\RN{1}})}, X_{i_1}^{(\text{\RN{1}})} \right\}_{i_1=1}^{n_1}$ to obtain initial estimates $\widehat{\mcT}^{\text{init},(\text{\RN{1}})}$ and $\whU_j^{\text{init}, (\text{\RN{1}})}$ for $j = 1, 2, 3$.
\State \quad Use dataset $\text{\RN{2}} = \left\{ Y_{i_2}^{(\text{\RN{2}})}, X_{i_2}^{(\text{\RN{2}})} \right\}_{i_2=1}^{n_2}$ to obtain initial estimates $\widehat{\mcT}^{\text{init},(\text{\RN{2}})}$ and $\whU_j^{\text{init}, (\text{\RN{2}})}$ for $j = 1, 2, 3$.

\State \textbf{Debiasing}
\State \quad Compute the debiased tensors:
\[
\widehat{\mcT}^{\text{unbs}, (\text{\RN{1}})} = \widehat{\mcT}^{\text{init}, (\text{\RN{2}})} + \frac{1}{n_1 \sigma^2} \sum_{i_1=1}^{n_1} \left( Y_{i_1}^{(\text{\RN{1}})} - \langle \widehat{\mcT}^{\text{init}, (\text{\RN{2}})}, \mcX_{i_1}^{(\text{\RN{1}})} \rangle \right) \mcX_{i_1}^{(\text{\RN{1}})},
\]
\[
\widehat{\mcT}^{\text{unbs}, (\text{\RN{2}})} = \widehat{\mcT}^{\text{init}, (\text{\RN{1}})} + \frac{1}{n_2 \sigma^2} \sum_{i_2=1}^{n_2} \left( Y_{i_2}^{(\text{\RN{2}})} - \langle \widehat{\mcT}^{\text{init}, (\text{\RN{1}})}, \mcX_{i_2}^{(\text{\RN{2}})} \rangle \right) \mcX_{i_2}^{(\text{\RN{2}})}.
\]

\State \textbf{One-step Power Iteration}
\For{$j = 1, 2, 3$}
    \State Compute $\whU_j^{(\text{\RN{1}})}$ and $\whU_j^{(\text{\RN{2}})}$ as the leading $r_j$ left singular vectors of
    \begin{align*}
    & \operatorname{Mat}_j\left( \widehat{\mcT}^{\text{unbs}, (\text{\RN{1}})} \times_{i \neq j} \whU_i^{\text{init},(\text{\RN{2}})\top} \right) = \operatorname{Mat}_j\left( \widehat{\mcT}^{\text{unbs}, (\text{\RN{1}})} \right) \left( \otimes_{i \neq j} \whU_i^{\text{init}, (\text{\RN{2}})} \right) \\
    & \operatorname{Mat}_j\left( \widehat{\mcT}^{\text{unbs}, (\text{\RN{2}})} \times_{i \neq j} \whU_i^{\text{init},(\text{\RN{1}})\top} \right) = \operatorname{Mat}_j\left( \widehat{\mcT}^{\text{unbs}, (\text{\RN{2}})} \right) \left( \otimes_{i \neq j} \whU_i^{\text{init}, (\text{\RN{1}})} \right).
    \end{align*}
\EndFor

\State \textbf{Projection and Plug-in Estimator}
\State \quad Define projection matrices $\mcP_{\whU_j^{(\text{\RN{1}})}}$ and $\mcP_{\whU_j^{(\text{\RN{2}})}}$ for $j = 1, 2, 3$. Compute the averaged projected tensor:
\[
\widehat{\mcT} = \frac{n_1}{n} \widehat{\mcT}^{\text{unbs}, (\text{\RN{1}})} \times_1 \mcP_{\whU_1^{(\text{\RN{1}})}} \times_2 \mcP_{\whU_2^{(\text{\RN{1}})}} \times_3 \mcP_{\whU_3^{(\text{\RN{1}})}} + \frac{n_2}{n} \widehat{\mcT}^{\text{unbs}, (\text{\RN{2}})} \times_1 \mcP_{\whU_1^{(\text{\RN{2}})}} \times_2 \mcP_{\whU_2^{(\text{\RN{2}})}} \times_3 \mcP_{\whU_3^{(\text{\RN{2}})}}.
\]

\State Estimate the linear functional $\langle \mcA, \mcT \rangle$ by
$
\langle \mcA, \widehat{\mcT} \rangle.
$

\State \textbf{Variance Estimation}

Let
\begin{align*}
\widehat{\sigma}_{\xi}^2
= &\frac{1}{n} \sum_{i_1=1}^{n_1} \left(Y_i^{\rm (\RN{1})}-\left\langle\widehat{\mcT}^{\text{init, (\RN{2})}}, \mcX_i^{\rm (\RN{1})}\right\rangle\right)^2+\frac{1}{n} \sum_{i_2=1}^{n_2}\left(Y_i^{(\rm \RN{2})}-\left\langle\widehat{\mcT}^{\text{init, (\RN{1})}}, \mcX_i^{\rm (\RN{2})}\right\rangle\right)^2 
\end{align*}
be the estimate of the observational noise variance, and
$
\widehat\sigma^2 = \frac{1}{np_1p_2p_3}\sum_{i=1}^n \left\|\mcX_i\right\|_{\mathrm{F}}^2. 
$
be the estimate of design noise. In addition, let 
\begin{align*}
& \widehat{s}_{\mcA}^2 = \sum_{j=1}^3\left\|\left(I-\mcP_{\whU_j}\right)A_j\left(\whU_{j+2}\otimes \whU_{j+1}\right)\widehat{W}_j\widehat{W}_j^{\top}\left(\whU_{j+2}\otimes \whU_{j+1}\right)^{\top}\right\|_{\mathrm{F}}^2 + \left\|\mcA \times_1 \whU_1 \times_2 \whU_2 \times_3 \whU_3\right\|_{\mathrm{F}}^2, 
\end{align*}
be the estimate of estimate of variance component ${s}_{\mcA}^2$, where
$U_j$ is either $\whU_j^{(\text{\RN{1}})}$ or $\whU_j^{(\text{\RN{2}})}$ and
\begin{align*}
& \widehat{W}_j=\mathrm{QR}\left[\Mat_j\left(\widehat{\mcT} \times_1 \whU_1^{\top} \times_2 \whU_2^{\top} \times_3 \whU_3^{\top} \right)^{\top}\right].
\end{align*}

\State \textbf{Confidence Interval}

The $100(1 - \alpha)\%$ confidence interval is given by
$$
\widehat{\mathrm{CI}}_{\mcA, \cT}^{\alpha}=\left[\left\langle  \widehat{\mcT}, \mcA \right\rangle-z_{\alpha / 2} \cdot \frac{\widehat{\sigma}_{\xi}}{\widehat{\sigma}} \cdot \widehat{s}_{\mcA} \sqrt{\frac{1}{n}}, \left\langle \widehat{\mcT}, \mcA \right\rangle + z_{\alpha / 2} \cdot \frac{\widehat{\sigma}_{\xi}}{\widehat{\sigma}} \cdot \widehat{s}_{\mcA} \sqrt{\frac{1}{n}}\right]  ,
$$
where $\alpha \in (0, 1)$, and $z_\theta = \Phi^{-1}(1 - \theta)$ denotes the upper $\theta$ quantile of the standard normal distribution.

\end{algorithmic}
\end{algorithm}

\subsection{Inference procedure for tensor PCA} \label{sec: Inference procedure for tensor PCA}

The inference procedure for tensor PCA is designed to estimate the linear functional $\left\langle \mcA, \mcT \right\rangle$ when the signal tensor is observed with additive noise. Algorithm \ref{alg: tensor PCA} summarizes the key steps.

\begin{algorithm}[hp]
\caption{Inference Procedure for Tensor PCA}
\label{alg: tensor PCA}
\begin{algorithmic}[1]

\State \textbf{Initialization}
\State \quad Take the observed tensor $\mcY$ as the initial estimate of the signal tensor $\mcT$. Obtain initial estimates $\whU_j^{\text{init}}$ for $j = 1, 2, 3$. Set $\whU_j^{(0)} := \whU_j^{\text{init}}$.

\State \textbf{Two-step Power Iteration}
\For{$k = 1, 2$}
    \For{$j = 1, 2, 3$}
        \State Compute $\whU_j^{(k)}$ as the leading $r_j$ left singular vectors of
        \[
        \Mat_j\left( \mcY \times_{j+1} \whU_{j+1}^{(k-1)\top} \times_{j+2} \whU_{j+2}^{(k-1)\top} \right).
        \]
    \EndFor
\EndFor
\State Set $\whU_j := \whU_j^{(2)}$ for $j = 1, 2, 3$.

\State \textbf{Projection and Plug-in Estimator}
\State \quad Define projection matrices $\mcP_{\whU_j} = \whU_j \whU_j^\top$ for $j = 1, 2, 3$. Compute the projected tensor:
$$
\widehat{\mcT} = \mcY \times_1 \mcP_{\whU_1} \times_2 \mcP_{\whU_2} \times_3 \mcP_{\whU_3}.
$$

\State \quad Estimate the linear functional $\langle \mcA, \mcT \rangle$ by $\langle \mcA, \widehat{\mcT} \rangle.$

\State \textbf{Variance Estimation}

Let
$$
\widehat{\sigma}^2=\frac{\left\|\mcY-\mcY \times_1 \mcP_{\whU_1} \times_2 \mcP_{\whU_2} \times_3 \mcP_{\whU_3} \right\|_{\mathrm{F}}^2}{p_1 p_2 p_3}
$$
be the estimate of the observational noise variance.

In addition, let 
\begin{align*}
& \widehat{s}_{\mcA}^2 = \sum_{j=1}^3\left\|\left(I-\mcP_{\whU_j}\right)A_j\left(\whU_{j+2}\otimes \whU_{j+1}\right)\widehat{W}_j\widehat{W}_j^{\top}\left(\whU_{j+2}\otimes \whU_{j+1}\right)^{\top}\right\|_{\mathrm{F}}^2 + \left\|\mcA \times_1 \whU_1 \times_2 \whU_2 \times_3 \whU_3\right\|_{\mathrm{F}}^2, 
\end{align*}
be the estimate of estimate of variance component ${s}_{\mcA}^2$, where
$ \widehat{W}_j=\mathrm{QR}\left[\Mat_j\left(\widehat{\mcT} \times_1 \whU_1^{\top} \times_2 \whU_2^{\top} \times_3 \whU_3^{\top} \right)^{\top}\right].
$

\State \textbf{Confidence Interval}

The $100(1 - \alpha)\%$ confidence interval is given by
$$
\widehat{\mathrm{CI}}_{\mcA,\mcT}^{\alpha}=\left[\langle \widehat{\mcT}, \mcA \rangle-z_{\alpha / 2} \cdot \widehat{\sigma} \widehat{s}_{\mcA} ,\langle  \widehat{\mcT}, \mcA \rangle+z_{\alpha / 2} \cdot \widehat{\sigma} \widehat{s}_{\mcA} \right],
$$
where $\alpha \in (0, 1)$, and $z_\theta = \Phi^{-1}(1 - \theta)$ denotes the upper $\theta$ quantile of the standard normal distribution.

\end{algorithmic}
\end{algorithm}


\end{document}